\documentclass[12pt,a4paper]{report}

\usepackage{amssymb,amscd,amsthm,mathrsfs}
\usepackage[ansinew]{inputenc}
\usepackage[centertags]{amsmath}
\usepackage[toc,page]{appendix}
\usepackage[dvips,final]{graphicx}
\usepackage[dvips]{geometry}
\usepackage{color} 
\usepackage{epsfig}
\usepackage{latexsym}
\usepackage{pstricks}

\newtheorem*{teo*}{Theorem}
\newtheorem*{teoe*}{Teorema}
\newtheorem*{teof*}{Th\'eor\`eme}
\newtheorem*{cor*}{Corollary}
\newtheorem{teo}{Theorem}[section]
\newtheorem{conj}[teo]{Conjecture}
\newtheorem*{conj*}{Conjecture}
\newtheorem*{conje*}{Conjetura}
\newtheorem*{conjf*}{Conjecture}
\newtheorem{quest}[teo]{Question}

\newtheorem{prop}[teo]{Proposition}
\newtheorem{cor}[teo]{Corollary}
\newtheorem*{claim}{Claim}
\newtheorem{lema}[teo]{Lemma}
\newtheorem{ad}[teo]{Addendum}
\newtheorem*{af}{Claim}
 \theoremstyle{definition}
\newtheorem*{notacion}{Notation}
\newtheorem*{convention}{Convention}
\newtheorem{defi}{Definition}[section]
\theoremstyle{remark}
\newtheorem{obs}[teo]{Remark}
\newtheorem*{obs*}{Remark}
\newtheorem*{obse*}{Observaci\'on}
\newtheorem*{obsf*}{Remarque}

\newcommand{\D}{\mathbb{D}}
\newcommand{\R}{\mathbb{R}}
\newcommand{\Z}{\mathbb Z}
\newcommand{\N}{\mathbb{N}}

\newcommand{\T}{\mathbb{T}}
\newcommand{\U}{\mathcal{U}}

\newcommand{\Fol}{\mathcal{F}}
\newcommand{\Res}{\mathcal{G}}

\newcommand {\eps}{\varepsilon}

\newcommand{\en}{\subset}

\newcommand{\bi}{\begin{itemize}}
\newcommand{\ei}{\end{itemize}}

  \def\CC{{\mathbb C}} \def\DD{{\mathbb D}}

 \def\NN{{\mathbb N}}  
\def\QQ{{\mathbb Q}} \def\RR{{\mathbb R}} \def\SS{{\mathbb S}} \def\TT{{\mathbb T}}
   
 \def\ZZ{{\mathbb Z}}

\def\cA{\mathcal{A}}  \def\cG{\mathcal{G}} \def\cM{\mathcal{M}} \def\cS{\mathcal{S}}
\def\cB{\mathcal{B}}  \def\cH{\mathcal{H}} \def\cN{\mathcal{N}} 
\def\cC{\mathcal{C}}  \def\cI{\mathcal{I}} \def\cO{\mathcal{O}} \def\cU{\mathcal{U}}
\def\cD{\mathcal{D}}  \def\cJ{\mathcal{J}} \def\cP{\mathcal{P}} \def\cV{\mathcal{V}}
\def\cE{\mathcal{E}}  \def\cK{\mathcal{K}} \def\cQ{\mathcal{Q}} \def\cW{\mathcal{W}}
\def\cF{\mathcal{F}}  \def\cL{\mathcal{L}} \def\cR{\mathcal{R}} \def\cX{\mathcal{X}}

\newcommand{\esbozo}[1]{\vspace{.05in}{\sc\noindent Sketch#1. }}
\newcommand{\dem}{\vspace{.05in}{\sc\noindent Proof. }}
\newcommand{\demo}[1]{\vspace{.05in}{\sc\noindent Proof#1. }}
\newcommand{\lqqd}{\par\hfill {$\Box$} \vspace*{.1in}}
\newcommand{\finobs}{\par\hfill{$\diamondsuit$} \vspace*{.1in}}
\newcommand{\trans}{\mbox{$\,{ \top} \;\!\!\!\!\!\!\raisebox{-.3ex}{$\cap$}\,$}}
\newcommand{\dashrv}{\vdash \!\! \dashv}

\DeclareMathOperator{\expon}{exp} \DeclareMathOperator{\Per}{Per}
\DeclareMathOperator{\Fix}{Fix} \DeclareMathOperator{\Lim}{Lim}

\DeclareMathOperator{\Diff}{Diff}
\DeclareMathOperator{\Homeo}{Homeo}
\DeclareMathOperator{\emb}{Emb} \DeclareMathOperator{\modulo}{mod}
\DeclareMathOperator{\diam}{diam} \DeclareMathOperator{\Leb}{Leb}
\DeclareMathOperator{\Int}{Int} \DeclareMathOperator{\Vol}{Vol}
\DeclareMathOperator{\length}{length}
\DeclareMathOperator{\Bas}{Bas} \DeclareMathOperator{\supp}{supp}
\DeclareMathOperator{\Reg}{Reg} \DeclareMathOperator{\Tang}{Tang}
\DeclareMathOperator{\diametro}{diam}
\DeclareMathOperator{\Hol}{Hol} \DeclareMathOperator{\Germ}{Germ}
\DeclareMathOperator{\rango}{Rank}
\DeclareMathOperator{\interior}{Int}

\begin{document}

\titlepage

\thispagestyle{empty}

\vspace*{2cm}

\begin{center}
\LARGE {\textsc{Thesis}} \\
\vspace{.50in}

{\Huge\textbf {Partial hyperbolicity and attracting regions in $3$-dimensional manifolds}}\\
\vspace{.50in}  By: Rafael Potrie Altieri\\ Advisors: Sylvain
Crovisier and Mart\'in Sambarino \vspace{0.3in}

\vspace{.50in} \LARGE{\textsc Doctorado en Matem\'atica\\
PEDECIBA\\ Universidad de la Rep\'ublica, Uruguay} \vspace{0.3in}

\vspace{.50in} \LARGE{\textsc Th\`ese \\} {Pour obtenir le grade de}\\
\LARGE{ DOCTEUR DE L'UNIVERSIT\'E PARIS 13\\} {Discipline:
Math\'ematiques} \vspace{0.3in}

\end{center}

\newpage

\thispagestyle{empty}







\newpage

\newpage
\vspace*{1cm}

\begin{flushright}
\emph{\LARGE{ A Nati}}
\end{flushright}

\newpage
\newpage
\begin{center}
\textbf{\LARGE{Acknowledgements}}
\end{center}

\medskip
\medskip
\medskip
The thesis could not have been done without the financial help by
the following institutions: ANR, CNRS, IFUM, ANII, PEDECIBA, CMAT,
CSIC, French Embassy in Uruguay (thanks in particular to Graciela
Vigo for her constant support and help). I am grateful to all
these institutions.

This has been a long trip and many people have of course been
there sharing the journey. All of them should know my deep
gratitude. However, I wish to express explicit acknowledgement to
some of the persons which were specially important for the project
itself.

Both Sylvain Crovisier and Mart\'in Sambarino, my advisors, were
of course very important and even essential in this project. I
want to thank Sylvain for his infinite patience and dedication, he
has clearly dedicated more time to this project than it is
expected and this has had an important effect on its results, I
want to profit to express my great admiration to him as a
brilliant mathematician and as a person: He is incredibly
generous, humble and kind, and even if sometimes it is hard to get
strength to speak to someone which is so good, he makes you feel
always comfortable and he is always encouraging, I could never
even try to list all what I learned from him. Mart\'in has been my
advisor and friend since the very start of my studies, he has been
my window to mathematics and he has always been open to pursue
with me in the directions which I wanted, as an advisor he has
always put priority in my formation over any other subject. As a
member of the community of mathematics in Uruguay, he has always
been and will always be a role model.

I greatly appreciate Christian Bonatti and Amie Wilkinson for
accepting the task of reading this (long) thesis, it is really an
honor to me. Christian has been part of this project all along
participating in several ways and being always very kind and
generous. Amie has accepted to referee this thesis in a difficult
moment, I am very happy to thank her also for devoting some of her
time to read and listen to me. Thanks also to Julien Barral,
Pierre Berger, \'Alvaro Rovella and Ra\'ul Ures for accepting
being part of the jury.

My stay in Paris with Natalia was an important stage for this
project. The fact that Andr\'es Sambarino was there was very
important. I thank Andr\'es and Juliana for their openness and
help in every stage of that stay and for sharing many good
moments. We met a lot of people in Paris which made our stay
great, but I would like to mention specially Freddy and Milena.

Mathematically, when staying in Paris, I was fortunate to meet and
discuss with many great mathematicians (some of them frequent
visitors of France) which were always open and they were important
for the development of this project. Thanks to: Pierre Berger,
Jerome Buzzi, Lorenzo D\'iaz, Tobias Jager, Patrice Le Calvez,
Katsutoshi Shinohara, Dawei Yang and specially to Christian
Bonatti and Nicolas Gourmelon. Also, I would like to thank the
whole Dynamical Systems Group of Villetaneuse for sharing the
seminar during the whole stay.

In my visits to Rio, I was fortunate to exchange with other
mathematicians which were also important in the development of
this project. Jairo Bochi, Pablo D\'avalos, Lorenzo D\'iaz, Andy
Hammerlindl and Enrique Pujals. Special thanks to Pablo Guarino,
Yuri Ki, Freddy Hernandez and Milena Cort\'ez not only for
discussions but also for their hospitality.

In a visit to Cuernavaca (M\'exico) before I started this project (around 2007), I was
 fortunate to meet Alberto Verjovski who in just an hour made an informal account of
 more results about topology and foliations that the ones I am able to understand even at
 this moment. Unexpectedly, those results and my further study of parts of what he taught me
 resulted very important in this project.

I wish to thank the dynamical systems group of Montevideo. My math
discussions with people in Montevideo from the start of
undergraduate studies has shaped my mathematical taste as well as
it has made a pleasure to study mathematics. I wish to thank
specially those which I consider have been essential, though for
very different reasons: Diego Armentano, Alfonso Artigue,
Joaqu\'in Brum, Matias Carrasco, Marcelo Cerminara, Mariana
Haim\footnote{Who in addition devoted a great deal of time in
helping me translate the introduction to french.}, Jorge Groisman,
Pablo Guarino, Nancy Guelman, Pablo Lessa, Jorge Lewowicz, Roberto
Markarian, Alejandro Passeggi, Miguel Paternain, Mariana Pereira,
\'Alvaro Rovella, Andr\'es Sambarino, Mart\'in Sambarino, Armando
Treibich, Jos\'e Vieitez, Mario Wschebor and Juliana Xavier. Among
these, I wish to highlight conversations with Andr\'es Sambarino
and Alejandro Passeggi which have always transcended math,
discussions with Joaqu\'in Brum and Pablo Lessa whose insight and
understanding of different mathematical topics never ceases to
surprise me and the courses with Miguel Paternain, \'Alvaro
Rovella and Mart\'in Sambarino which have been by far those which
made me grow the most as a mathematician.

I've done my best to write this thesis, but as in other tasks, my
best in writing is usually very poor, so, my apologizes go along
my gratefulness to anyone who dares to read this (or part of this)
thesis.

I wish to thank my friends, both in and outside mathematics, both
in Montevideo and in Paris, specially those which have been
present all along these years.

I wish to thank my family, specially my parents for their support
and initial education which is of course essential for pursuing
this kind of goals. I wish also to remember in this section of the
thesis my cousin Laura, which passed away during my stay in Paris
and has marked me tremendously.

Finally, I wish to thank Natalia, my wife, my partner and my
friend. Nothing would have sense if you were not there.

\newpage
\newpage
\begin{center}
{\bf Abstract}
\end{center}

\medskip

This thesis attempts to contribute to the study of differentiable
dynamics both from a semi-local and global point of view. The
center of study is differentiable dynamics in manifolds of
dimension $3$ where we are interested in the understanding of the
existence and structure of attractors as well as dynamical and
topological implications of the existence of a global partially
hyperbolic splitting. The main contributions are new examples of
dynamics without attractors where we get a quite complete
description of the dynamics around some wild homoclinic classes
(see Section \ref{SubSection-LocalizacionDeClasses} and subsection
\ref{SubSection-EjemploDA}) and two results on dynamical coherence
of partially hyperbolic diffeomorphisms of $\TT^3$ (see Chapter
\ref{Capitulo-ParcialmenteHiperbolicos}).
\medskip
\medskip

\begin{center}
{\bf Resumen}
\end{center}

\medskip
Esta tesis pretende contribuir al estudio de la din\'amica
diferenciable tanto desde sus aspectos semilocales como globales.
El estudio se centra en din\'amicas diferenciables en variedades
de dimensi\'on 3. Se busca comprender por un lado la existencia y
estructura de los atractores as\'i como propiedades topol\'ogicas
y din\'amicas implicadas por la existencia de una descomposici\'on
parcialmente hiperb\'olica global. Las contribuciones principales
son la construcci\'on de nuevos ejemplos de din\'amicas sin
atractores donde se da una descripci\'on bastante completa de la
din\'amica alrededor de una clase homocl\'inica salvaje (ver
Secci\'on \ref{SubSection-LocalizacionDeClasses} y la subsecci\'on
\ref{SubSection-EjemploDA}) y dos resultados sobre la coherencia
din\'amica de difeomorfismos parcialmente hiperb\'olicos en
$\TT^3$ (ver Cap\'itulo \ref{Capitulo-ParcialmenteHiperbolicos}).

\medskip
\medskip

\begin{center}
{\bf Resum\`e}
\end{center}

\medskip

Le but de cette th\`ese est de contribuer \`a la compr\'ehension
des dynamiques diff\'erentiables aussi bien d\'un point de vue
semilocal que global. L'etude se concentre sur les
diffeomorphismes des vari\'et\'es de dimension 3. On cherche \`a
comprendre l'existence et la structure de leurs attracteurs, mais
aussi \`a decrire les propri\'et\'es topologiques et dynamiques
des diff\'eomorphismes partiellement hyperboliques globaux. Les
contributions pricipales sont la construction de nouvelle
dynamiques sauvages (voir Section
\ref{SubSection-LocalizacionDeClasses} et subsection
\ref{SubSection-EjemploDA}) et deux r\'esultats sur la coh\'erence
dynamique des diff\'eomorphismes partiellement hyperboliques dans
$\TT^3$ (voir Chapitre \ref{Capitulo-ParcialmenteHiperbolicos}).


\tableofcontents

   \chapter*{Notations}

\bi

\item[] For a compact metric space $X$ we denote $d(\cdot,\cdot)$
to the metric of $X$.

\item[] $M^d$ will denote a compact connected Riemannian manifold
without boundary of dimension $d \in \NN$. It is a metric space
whose metric is induced by the Riemannian metric $\langle \cdot
,\cdot \rangle$.

\item[]$\Leb(\cdot)$ will denote the measure induced by any volume
form on $M$ of total measure $1$. For our purposes it will make no
difference which one is it (since we shall not assume the maps to
preserve it) and we shall only care about sets having positive,
total or zero measure.

\item[] For $X\en M$, we denote $T_X M = \bigcup_{x \in X} T_xM$
with the topology induced by the inclusion $T_XM \en TM$ into the
tangent bundle of $M$.

\item[] $\Diff^r(M)$ ($r\geq 0$) denotes the set of
$C^r$-diffeomorphisms (homeomorphisms in the case $r=0$) with the
$C^r$ topology (see \cite{HirschLibro}). We shall denote the
distance in $\Diff^r(M)$ as $d_{C^r}(\cdot,\cdot)$. It is a Baire
space. Similarly, $C^r(M,N)$ denotes the space of $C^r$-maps from
$M$ to $N$ and $\emb^r(M,N)$ the space of $C^r$-embeddings.

\item[] For $f \in \Diff^1(M)$ we denote as $D_x f : T_xM \to
T_{f(x)}M$ the derivative of $f$ over $x$. Sometimes, we shall not
make reference to the point $x$ when it is understood.

\item[] For $V, W$ submanifolds of $M$ we say that they intersect
\emph{transversally} at $x \in V\cap W$ if we have that $T_x V +
T_xW = T_xM$. The set of points in $V\cap W$ where the
intersection is transversal is denoted by $V \trans W$. When $V
\cap W = V \trans W$ we say that $V$ and $W$ intersect
\emph{transversally}.

\item[] For $V, W$ compact embedded submanifolds of $M$ (possibly
with boundary) which are diffeomorphic to a certain manifold $D$,
we define the $C^r$-distance between them as the infimum of the
$C^r$-distance between the pairs of embeddings of $D$ in $M$ whose
image is respectively $V$ and $W$.

\item[] For Baire spaces (in particular, sets which are metric and
complete or open subsets of these) we say that a set $\cG$ is
\emph{residual} (or $G_\delta-$dense) if it is a countable
intersection of open and dense subsets.

\item[] We shall say that a property verified by diffeomorphisms
in $\Diff^r(M)$ is $C^r$-\emph{generic} if it is verified by
diffeomorphisms in a residual subset of $\Diff^r(M)$. Sometimes,
hoping it makes no confusion, we will say that a diffeomorphism
$f$ is a $C^r$-\emph{generic diffeomorphism} to mean that $f$
verifies properties in a residual subset $\Diff^r(M)$ (which will
be clear from the context).

\item[] $\TT^d$ will denote the (flat) $d$-dimensional torus
$\RR^d \slash_{\ZZ^d}$ with the metric induced by the canonical
covering map $p: \RR^d \to \TT^d$ and the Euclidean metric.

\item[] $B_\eps(x)$ denotes the (open) $\eps$-neighborhood of the
point $x$, i.e. the (open) set of points at distance smaller than
$\eps$ of $x$. 

\item[] $B_\eps(K)$ denotes the (open) $\eps$-neighborhood of the
set $K$.

\item[] Given a subset $A$ of a metric space $X$ we denote
$\Int(A), \overline{A}, \partial A, A^c$ to the interior, closure,
frontier and complement of $A$ respectively.

\item[] Given a point $x\in A$ we will denote $cc_x(A)$ to the
connected component of $A$ containing $x$.

\item[] Given a compact metric space $X$, we denote $\cK(X)$ to be
the set of compact subsets of $X$ endowed with the Hausdorff
distance: $$ d_H(A,B) = \max\{ \sup_{x\in A} \inf_{y\in B} \{
d(x,y)\}, \sup_{y\in B} \inf_{x\in A} \{ d(x,y)\}\} $$ \noindent
which is compact.

\item[] Given a sequence of sets $A_n \en X$ a topological space.
We define $\limsup A_n = \bigcap_{n\geq 0}
\overline{\bigcup_{k\geq n} A_n} $.

\item[] We use the symbol $\Box$ to denote the end of a proof of a
Theorem, Lemma, Proposition or Corollary. We use the symbol
$\diamondsuit$ to denote the end of a Remark, Definition or the
proof of some Claim (inside the proof of something else).

\ei

   \chapter{Introduction and presentation of results}

\section{Introduction (English)}

\subsection{Historical account and context}

One may\footnote{We warn the reader that the historic context we
will present is completely subjective and not necessarily reflects
the true historical facts. It must be thought as a plausible
context in which the work of this thesis fits.} claim that the
main goal in dynamical systems is to understand the asymptotic
behavior of orbits for a given evolution law. Originally, the
subject began with the study of ordinary differential equations of
the form

$$ \dot x = X(x)  \qquad  X: \RR^n \to \RR^n $$

\noindent and the goal was to solve these equations analytically
and obtaining, for each initial value $x_0 \in \RR^n$ an explicit
solution  $\varphi_t(x_0)$ to the differential equation.

It was soon realized that even extremely simple equations gave
rise to complicated analytic solutions. Moreover, it was realized
that the integrated equations did not supply enough understanding
of the laws of evolutions.

By studying the famous $3$-body problem, Poincar\'e
(\cite{Poincare}) was probably the first to propose that there
should be a qualitative study of evolution rather than a
quantitative one and he proposed to study ``the behavior of
\emph{most} orbits for the \emph{majority} of systems".

At the start, the study focused on stability. Lyapunov studied
stable orbits, this means, orbits which contain a neighborhood of
points having essentially the same asymptotic behavior. Andronov
and Pontryagin, followed by Peixoto, studied stable systems, this
means, those whose dynamical properties are robust under
perturbations. But it was probably Smale (\cite{SmaleBulletin})
the first to revitalize Poincar\'e's suggestion by giving to it a
more precise formulation:

The goal is to fix a closed manifold $M$ of dimension $d$ and to
understand the dynamics of a \emph{large} subset of $\Diff^r(M)$,
the space of diffeomorphisms of $M$ endowed with the
$C^r$-topology.

Moreover, he proposed that a subset of diffeomorphisms should be
considered large if it was open and dense with this topology, or
at least, residual or dense (in a way that by understanding large
sets of diffeomorphisms one could not neglect behavior happening
in a robust fashion). We will not discuss other possible notions
of largeness used in the literature nor the reasons for
considering this ones (we refer the reader to \cite{B} or
\cite{Crov-Hab} for an explanation of this choice).

Structural stability became the center of Smale's program which
was strongly based on the hope that even if dynamical systems
could not be stable from the point of view of their dynamics (they
could be chaotic) they could be, at least the majority of them,
stable from the point of view of their orbit structure. This would
give that their dynamics and the one of their perturbations could
be understood by symbolic or probabilistic methods. Palis and
Smale \cite{PalisSmale} conjectured that structurally stable
systems coincide with hyperbolic ones.

Hyperbolicity became the paradigm. Robbin and Robinson
(\cite{Robbin,RobinsonEstabilidad}) proved that hyperbolic systems
were stable and long afterwards Ma\~ne showed (\cite{ManheIHES})
that $C^1$-structurally stable dynamics were indeed hyperbolic
(this is still unknown in other topologies). Describing the
dynamics of hyperbolic diffeomorphisms was the center of attention
for dynamicists in the 60's and 70's.


This project was tremendously successful from this point of view
and it was not only the semi-local study (through the use of
symbolic and ergodic techniques) that was understood but also very
deep global aspects as well as some understanding of the topology
of basic pieces was achieved.

On semi-local aspects, without being exhaustive, we mention
particularly the contributions of Bowen, Newhouse, Palis, Sinai,
Ruelle and Smale. We refer the reader to \cite{KH} Part 4 for a
panoramic view of a large part of the theory.

On the other hand, the global aspects of the study were mainly
associated to the work of Anosov, Bowen, Franks, Shub, Smale,
Sullivan and Williams and good part of those can be appreciated in
the nice book of Franks \cite{Franks-Homology}. It is worth
mentioning that, for different reasons, people working in this
aspects of differentiable dynamics abandoned the subject and this
may be an explanation on why these results are less popular.

However, the program of Smale, as well as the hope that
structurally stable systems should be typical among
diffeomorphisms of a manifold fell down after some examples of
robustly non-hyperbolic dynamics started to appear. The first non
hyperbolic examples were those of Abraham-Smale (\cite{AS}) and
Newhouse (\cite{Newhouse}).

This gave rise to the theory of bifurcations, where Newhouse,
Palis and Takens (among others) were pioneers and after many work
and new examples the initial program was finally adapted in order
to contemplate these new examples and to maintain the initial
philosophy of Smale. Palis' program \cite{PalisGlobalView},
however, has only a semilocal point of view.

After the paradigm of hyperbolicity began to fall, the research
started focusing on finding alternative notions, such as
non-uniform hyperbolicty (mainly by the Russian school, of which
the principal contributors were Pesin and Katok) or the partial
hyperbolicity (independently by Hirsch-Pugh-Shub \cite{HPS} and
Brin-Pesin \cite{BrinPesin}). In this thesis, we are mainly
interested in the second generalization of hyperbolicity for its
condition of geometric structure (in contrast with the measurable
structure given by non-uniform hyperbolicity) and its strong
relationship with robust dynamical properties. See \cite{BDV} for
a panorama on dynamics beyond hyperbolicity.

In his quest for a proof of the stability conjecture, Ma\~ne
(independently also Pliss and Liao \cite{pliss, Liao}) introduced
the concept of dominated splitting and showed its close
relationship with the dynamics of the tangent map over periodic
orbits.

When one studies the space of diffeomorphisms with the
$C^1$-topology the perturbation techniques developed since the
60's by Pugh, Ma\~ne, Hayashi and more recently by Bonatti and
Crovisier imply that the periodic orbits capture in a good way
(topological and statistical) the dynamics of generic
diffeomorphisms. See \cite{Crov-Hab} for a survey on this topics.

Recently, Bonatti \cite{B} has proposed a realistic program for
the study of the dynamics of $C^1$-generic diffeomorphisms which
extends Palis' program and complements it. It is also the case
that this program has a semilocal flavour.

From the global point of view, there is much less work done, and
also less proposals on how to proceed (see \cite{PujSamHandbook}
section 5 for a short survey) although there are some ideas on how
to proceed in some cases at least in dimension 3.

In what follows, we will try to present the contributions of this
thesis and explain how our results fit in this subjective account
of the development of differentiable dynamics.

\subsection{Attractors in $C^1$-generic dynamics}

It is always possible to decompose the dynamics of a homeomorphism
of a compact metric space into its chain-recurrence classes. This
is the content of Conley's theory \cite{Conley}.

This decomposition has proven very useful in the understanding of
$C^1$-generic dynamics\footnote{We will use this expression to
refer to diffeomorphisms belonging to a residual subset of
$\Diff^1(M)$ with the $C^1$-topology.} thanks to a result by
Bonatti and Crovisier (\cite{BC}) which guaranties that it is
possible to detect chain recurrence classes of a generic
diffeomorphism by its periodic orbits. In a certain sense, the
dynamics around periodic orbits has attracted most of the
attention in the study of semi-local properties of generic
diffeomorphisms and it is hoped that by understanding their
behavior one will be able to understand $C^1$-generic dynamics
(see \cite{B}).

If we wish to understand the dynamics of most of the orbits, there
are some chain-recurrence classes which stand out from the rest.
\emph{Quasi-attractors} are chain-recurrence classes which admit a
basis of neighborhoods $U_n$ verifying that $f(\overline{U_n}) \en
U_n$. Such classes always exist, and it was proven in \cite{BC}
that, for $C^1$-generic diffeomorphisms, there is a residual
subset of points in the manifold whose forward orbit accumulates
in a quasi-attractor.

Sometimes, it is possible to show that these quasi-attractors are
isolated from the rest of the chain-recurrence classes and in this
case, we say that they are \emph{attractors}. Attractors have the
property of being accumulated by the future orbit of nearby points
and being dynamically indecomposable. To determine whether
attractors exist and the topological and statistical properties of
their basins is one of the main problems in non-conservative
dynamics. In dimension two, it is possible to show that
$C^1$-generic diffeomorphisms have attractors. This was originally
shown by Araujo \cite{Araujo} but there was a gap in the proof and
this was never published\footnote{See \cite{BrunoSantiago} for a
very recent correction of the original proof with the use of the
results of Pujals and Sambarino \cite{PujSamAnnals1}.}. This
result was in a certain way incorporated to the folklore (see for
example \cite{BLY}). In this thesis, we present a proof of the
following result which appeared in \cite{PotExistenceOfAttractors}
(see Section \ref{Section-AtractoresSuperficies}).

\begin{teo*}
There exists an open and dense subset $\cU$ of the space of
$C^1$-diffeomorphisms of a surface $M$ such that every $f\in \cU$
has a hyperbolic attractor. Moreover, if $f$ cannot be perturbed
in order to have infinitely many attracting periodic orbits, then
every quasi-attractor of $f$ is a hyperbolic attractor and there
are finitely many quasi-attractors.
\end{teo*}

When a quasi-attractor is hyperbolic, it must be an attractor. To
show that when a diffeomorphism has robustly finitely many
attracting periodic orbits, all quasi-attractors are hyperbolic it
is a key step to show that they admit what is called a
\emph{dominated splitting}. This means that there exists a
$Df$-invariant splitting of the tangent bundle over the
quasi-attractor into two subbundles which verify a uniform
condition of domination (vectors in one subbundle are uniformly
less contracted than on the other).

Dominated splitting, as well as many other $Df$-invariant
geometric structure, is an important tool for studying the
dynamical properties of a chain-recurrence class, and
fundamentally, to understand how the class is accumulated by other
classes (see Section \ref{Section-EstructurasInvariantes}).

Aiming at the understanding of the dynamics close to
quasi-attractors for $C^1$-generic diffeomorphisms, we have
obtained the following partial result about the structure of those
quasi-attractors which are homoclinic classes (see Section
\ref{Section-GenericBiLyapunov} and \cite{PotGenericBiLyapunov}):

\begin{teo*}
Let $f$ be a $C^1$-generic diffeomorphism and $H$ be a
quasi-attractor which contains a periodic orbit $p$ such that the
differential of $f$ over $p$ at the period contracts volume. Then,
$H$ admits a non-trivial dominated splitting.
\end{teo*}

From the point of view of the conclusion of this theorem, it is
possible to see by means of examples that the conclusion is in
some sense optimal (see \cite{BV}). The same happens with the
hypothesis of having a periodic orbit (see \cite{BDihes}).
However, the hypothesis on the dissipation of volume along a
periodic orbit seems to follow from the fact that $H$ is a
quasi-attractor but we were not able to prove this. Proving this
seems to require what is known as an \emph{ergodic closing lemma
inside a homoclinic class} which is a problem not well understood
for the moment (see \cite{B} Conjecture 2).

The main novelty in the proof of this result is the use of a new
perturbation result due to Gourmelon \cite{GouFranksLemma} which
allows to perturb the differential over a periodic orbit while
keeping control on its invariant manifolds. The use of this result
combined with Lyapunov stability has allowed us to solve the
problem of guaranteeing that a point remains in the class after
perturbation. The result responds affirmatively to a question
posed in \cite{ABD} (Problem 5.1).

The dream of having $C^1$-generic dynamics admitting
attractors\footnote{See the introduction of \cite{BLY} for more
historic account on the problem.} has fallen recently due to a
surprising example presented in \cite{BLY} which shows how the
recent development of the theory of $C^1$-generic dynamics has had
an important influence in the way we understand dynamics and has
simplified questions which seemed unapproachable.

The examples of \cite{BLY} posses what they have called
\emph{essential attractors} (see subsection
\ref{SubSection-AttractingSets}) and it is not yet known whether
they posses attractors in the sense of Milnor. In Section
\ref{Section-EjemplosQuasiAtractores} we review their examples and
present new examples from \cite{PotWildMilnor} on which we have a
better understanding on how other classes approach their dynamics:

\begin{teo*}
There exists an open set $\cU$ of $\Diff^1(\TT^3)$ of
diffeomorphisms such that:
\begin{itemize}
\item[-] If $f\in \cU$ then $f$ has a unique quasi-attractor and
an attractor in the sense of Milnor.  Moreover, if $f$ is $C^2$,
then it has a unique SRB measure whose basin has total Lebesgue
measure. \item[-] For $f$ in a $C^r$-residual subset of $\cU$ we
know that $f$ has no attractors.  \item[-] Chain-recurrence
classes different from the quasi-attractor are contained in
periodic surfaces.
\end{itemize}
\end{teo*}

The last point of the theorem contrasts with the new results
obtained by Bonatti and Shinohara (\cite{BS}). In Section
\ref{Section-TrappingAttractors} we speculate on how both results
can fit in the same theory.

\subsection{Partial hyperbolicity in $\TT^3$}

In the previous section we have discussed problems which are of
semilocal nature (although it is of course possible to ask
questions of global nature about attractors and their topology).
In this section we treat results of global dynamics.

A well known result in differentiable dynamics, which joins
classical results by Ma\~ne (\cite{ManheErgodic}) and Franks
(\cite{FranksAnosov}) can be stated as follows:

\begin{teo*}[Ma\~ne-Franks] Let $M$ be a compact surface.
The following three properties for $f \in \Diff^1(M)$ are
equivalent:
\begin{itemize}
\item[(i)] $f$ is $C^1$-robustly transitive. \item[(ii)] $f$ is
Anosov. \item[(iii)] $f$ is Anosov and conjugated to a linear
Anosov automorphism.
\end{itemize}
\end{teo*}

Ma\~ne proved that (i) $\Rightarrow$ (ii) while Franks had proven
(ii) $\Rightarrow$ (iii). Robustness of Anosov diffeomorphisms and
the fact that linear Anosov automorphisms are transitive gives
(iii) $\Rightarrow$ (i).

If we interpret being an Anosov diffeomorphism as having a
$Df$-invariant geometric structure, we can identify the result (i)
$\Rightarrow$ (ii) as saying: ``an robust dynamical property
forces the existence of a $Df$-invariant geometric structure''. In
fact, since $C^1$-perturbations cannot break the dynamical
behavior, it is natural to expect that this geometric structure
will posses certain rigidity properties.

On the other hand, the direction (ii) $\Rightarrow$ (i) can be
thought as a converse statement, showing that $Df$-invariant
geometric structures may imply the existence of certain robust
dynamical behavior, in this case, transitivity.

In higher dimensions, the understanding of the relationship
between robust dynamical properties and $Df$-invariant geometric
structures is quite less advanced although results in the
direction of obtaining a $Df$-invariant geometric structure from a
robust dynamical property do exist. In dimension $3$, it follows
from a result of \cite{DPU} and this was generalized to higher
dimensions by \cite{BDP} (see Section
\ref{Section-EstructurasInvariantes}):

\begin{teo*}[Diaz-Pujals-Ures]
Let $M$ be a $3$-dimensional manifold and $f$ a $C^1$-robustly
transitive diffeomorphism, then, $f$ is partially hyperbolic.
\end{teo*}

The definition of partial hyperbolicity varies throughout the
literature and time. We use the definition used in \cite{BDV}
which is the one that fits best our approach (see Section
\ref{Section-EstructurasInvariantes} for precise definitions). A
partially hyperbolic diffeomorphism for us will be one which
preserves a splitting of the tangent bundle $TM=E \oplus F$
verifying a domination property between the bundles and such that
one of them is uniform. For notational purposes, we remove the
symmetry of the definition and work with partially hyperbolic
diffeomorphisms of the form $TM = E^{cs} \oplus E^u$ where $E^u$
is uniformly expanded.

The first difficulty one encounters when trying to work in the
converse direction of the previous theorem is the fact that one
has no control on the contraction in the direction $E^{cs}$ other
that it is dominated by the expansion in $E^u$. This forbids us to
gain a complete control on the dynamics in that direction as we
have in the hyperbolic case. One notable exception is the work of
\cite{PujSamDominated} where precise dynamical consequences are
obtained from the existence of a dominated splitting in dimension
$2$.

It is now also time to mention the importance of item (iii) in
Ma\~ne-Franks' Theorem which we have neglected so far. In a
certain sense, the underlying idea is that in order to obtain a
robust dynamical property out of the existence of a $Df$-invariant
geometric structure it can be important to rely on the topological
restriction this geometric structure imposes, such as the topology
of the manifold or the isotopy class of the diffeomorphism. It is
for this reason, and the difficulties that have appeared in the
attempt to obtain results in converse direction of
Diaz-Pujals-Ures' Theorem that it seems for us a good idea to
divide the study in something in the spirit of (ii) $\Rightarrow$
(iii) and (ii)+(iii) $\Rightarrow$ (i) as in Ma\~ne-Franks'
Theorem. This means that the study of partially hyperbolic
dynamics in a fixed manifold or even in a fixed isotopy class
seems to be an important step in the understanding of these
relations.

Also in the direction of obtaining results giving topological
properties by the existence of a $Df$-invariant geometric
structure the difficulty increases considerably as we raise the
dimension. In dimension 2, the sole fact of admitting an
continuous line field imposes strong restrictions on the topology
of the manifold. In dimension 3, it is well known that every
$3$-dimensional manifold admits a non-vanishing vector field and
moreover, it also admits a codimension one foliation.

This situation may be considered as very bad from the point of
view of finding topological properties out of the existence of
$Df$-invariant geometric structures. However, very recently, a
beautiful remark by Brin-Burago-Ivanov (\cite{BBI,BI}) has renewed
the hope:

\medskip

\begin{obs*}[Brin-Burago-Ivanov]
In a $3$-dimensional manifold, if $\cF$ is a foliation transverse
to the unstable direction of a partially hyperbolic
diffeomorphism, then $\cF$ has no Reeb components.
\finobs
\end{obs*}

Reeb components, and its relationship with partially hyperbolic
dynamics had already been studied in \cite{DPU} (Theorem H) and
\cite{BonWilk} (Lemma 3.7) in more restrictive contexts (assuming
dynamical coherence and transitivity). This remark is much more
general, its strength relies in that the only dependence on the
dynamics is in the fact that the unstable foliation cannot have
closed curves.

On the one hand, it is known that many $3$-manifolds do not admit
foliations without Reeb components. On the other, there exists
quite a lot of theory regarding its classification (see for
example \cite{Plante, Roussarie}) and therefore, we can expect
that progress in the classification of partially hyperbolic
diffeomorphisms is within reach.

Unfortunately, another difficulty arises: it is not known if every
partially hyperbolic diffeomorphism in a $3$-dimensional manifold
posses a foliation transverse to the unstable direction. In this
thesis, we propose the notion of \emph{almost dynamical coherence}
which we show is an open and closed property among partially
hyperbolic diffeomorphisms and expect that under this hypothesis
more progress can be made.

One of our main results (see Chapter
\ref{Capitulo-ParcialmenteHiperbolicos} and
\cite{Pot3dimPHisotopic}) guaranties that in certain isotopy
classes of diffeomorphisms of $\TT^3$ partial hyperbolicity and
almost dynamical coherence are enough to guarantee the existence
of a $f$-invariant foliation tangent to $E^{cs}$.

\begin{teo*}
Let $f:\TT^{3} \to \TT^3$ be an almost dynamically coherent
partially hyperbolic diffeomorphism isotopic to a linear Anosov
automorphism. Then, $f$ is dynamically coherent.
\end{teo*}

Being \emph{dynamically coherent} means that the bundle $E^{cs}$
integrates into an $f$-invariant foliation.

In the \emph{strong partially hyperbolic} case (i.e. where $TM
=E^s \oplus E^c \oplus E^u$ is a $Df$-invariant splitting with
domination properties and $E^s$ and $E^u$ are uniform) we can say
more.

As a starting point, in a remarkable paper \cite{BI} it was proved
that every strong partially hyperbolic diffeomorphism is almost
dynamically coherent. This was used first by \cite{BI} (following
a very simple and elegant argument of \cite{BBI}) and then in
\cite{Parwani} to give topological conditions these must satisfy.
On the other hand, these progress makes expectable that, at least
in some special cases, the following conjecture of Pujals may be
within reach:

\begin{conj*}[Pujals \cite{BonWilk}]
Let $f:M \to M$ with $M$ a $3$-dimensional manifold a strong
partially hyperbolic diffeomorphism which is transitive. Then, one
of the following possibilities holds (modulo considering finite
lifts):
\begin{itemize}
\item[-] $f$ is leaf conjugate to a linear Anosov automorphism of
$\TT^3$. \item[-] $f$ is leaf conjugate to a skew-product over a
linear Anosov automorphism of $\TT^2$ (and so $M= \TT^3$ or a
nilmanifold). \item[-] $f$  is leaf conjugate to the time one map
of an Anosov flow.
\end{itemize}
\end{conj*}

There has been some progress in the direction of this conjecture
lately. Let us mention first the work of \cite{BonWilk} which
makes considerable progress. They work without making any
hypothesis on the topology of the manifold but they demand the
existence of a closed curve tangent to the center direction and
some other technical hypothesis. Then, the work of Hammerlindl
\cite{Hammerlindl, HNil} has given a proof of the conjecture when
the manifold is $\TT^3$ or a nilmanifold by assuming a more
restrictive notion of partial hyperbolicity (partial hyperbolicity
with absolute domination). Although this notion is verified by
many examples, it is in some sense artificial and does not fit
well with the results of \cite{DPU}.

Of course, to prove Pujals' conjecture, a previous step must be to
show dynamical coherence of such diffeomorphisms since leaf
conjugacy requires this for a start (see Section
\ref{Section-HIRSHPUGHSHUB}). The work of Hammerlindl relies
heavily on previous work by Brin-Burago-Ivanov \cite{BBI2} who
have established dynamical coherence of strong partially
hyperbolic diffeomorphisms of $\TT^3$ under this more restrictive
notion of partial hyperbolicity we mentioned above.

While one could expect the use of this restrictive notion to be
mainly technical, a recent example of Rodriguez Hertz-Rodriguez
Hertz-Ures (\cite{HHU}) of a non-dynamically coherent strong
partially hyperbolic diffeomorphism of $\TT^3$ shows that the
passage to the general definition should at least use the
transitivity hypothesis in a fundamental way and that some
difficulties must be addressed. We have completed the panorama
(\cite{Pot3dimPHisotopic}) by showing:

\begin{teo*}
Let $f: \TT^3 \to \TT^3$ be a strong partially hyperbolic
diffeomorphism which does not admit neither a periodic normally
attracting torus nor a periodic normally repelling torus, then $f$
is dynamically coherent.
\end{teo*}

This theorem responds to a conjecture made by Rodriguez
Hertz-Rodriguez Hertz-Ures in $\TT^3$ and it also allows one to
prove Pujals' conjecture for $M=\TT^3$ by further use of the
techniques in the proof (see \cite{HP}).

In the end of Chapter \ref{Capitulo-ParcialmenteHiperbolicos} we
also obtain some results in higher dimensions which are analogous
to the classical results of Franks, Newhouse and Manning for
Anosov diffeomorphisms in the context of partial hyperbolicity.

\subsection{Other contributions}

In this section we briefly describe other contributions of this
thesis.

In Section \ref{SubSection-LocalizacionDeClasses} we present a
mechanism from \cite{PotWildMilnor} for the localization of
chain-recurrence classes which we consider has intrinsic value
since it can be applied in many contexts (in this thesis, it is
applied in Section \ref{Section-DescDomT2} as well as in
subsections  \ref{SubSection-EjemploDA} and
\ref{SubSection-EjemploPlykin}).

In Section \ref{Section-EjemplosQuasiAtractores} we present
several examples of quasi-attractors and of robustly transitive
diffeomorphisms some of which are modifications of well known
examples but we consider they may contribute to the understanding
of these phenomena.

Then, in Chapter \ref{Capitulo-FoliacionesEHipParcial} we present
results on foliations which we use later in Chapter
\ref{Capitulo-ParcialmenteHiperbolicos} which we believe may have
independent interest. In particular, we mention a quantitative
result about the existence of a global product structure of
certain transverse foliations which is presented in Section
\ref{Section-GlobalProductStructure}. Also, in Section
\ref{Section-DescDomT2} we give a characterization of dynamics of
globally partially hyperbolic diffeomorphisms of $\TT^2$ to show
the techniques we use later in Chapter
\ref{Capitulo-ParcialmenteHiperbolicos}.

We also include in this thesis 4 appendices where some results
which we preferred to separate from the main line of the thesis
are presented. We make particular emphasis on Appendix
\ref{Apendice-PseudoRotaciones} based on \cite{PotRecurrence}
where we prove a result about homeomorphisms of $\TT^2$ with a
unique rotation vector and Appendix \ref{Apendice-BCGP} based on
\cite{BCGP} where we present a joint work with Bonatti, Crovisier
and Gourmelon. This last work studies the bifurcations of robustly
isolated chain-recurrence classes and gives examples of such
classes which are not robustly transitive answering to a question
posed in \cite{BC}.


\section{Introducci\'on (Espa\~nol)}

\subsection{Contexto hist\'orico}

Se puede\footnote{Vale aclarar que la introducci\'on hist\'orica
que se presentar\'a es subjetiva y no necesariamente refleja con
exactitud los hechos hist\'oricos. Se puede pensar que lo que
sigue es una historia plausible que explica algunas razones por
las cuales estudiar los temas aqu\'i presentados.} decir que el
objetivo fundamental de los sistemas din\'amicos es comprender el
comportamiento asint\'otico de un estado sujeto a una ley de
evoluci\'on. Se comenz\'o por el estudio de ecuaciones
diferenciales ordinarias del tipo

$$ \dot x = X(x)  \qquad  X: \RR^n \to \RR^n $$

\noindent donde se buscaba una soluci\'on anal\'itica a la
ecuaci\'on: buscando para cada valor posible de condici\'on
inicial $x_0 \in \RR^n$ una soluci\'on expl\'icita
$\varphi_t(x_0)$.

R\'apidamente se vi\'o que ecuaciones extremadamente simples daban
lugar a soluciones complicadas que incluso luego de ser integradas
tampoco aportaban a la comprensi\'on de la ley de evoluci\'on.

Poincar\'e, interesado en estudiar el famoso problema de los
$3$-cuerpos (\cite{Poincare}) fue quiz\'as el primero en proponer
que el estudio de la evoluci\'on deber\'ia ser cualitativo, y de
alguna manera propuso estudiar "el comportamiento de la
\emph{mayor\'ia} de las \'orbitas para la \emph{mayor\'ia} de los
sistemas".

El comienzo del estudio se centr\'o en la estabilidad. Lyapunov
estudi\'o las \'orbitas estables, \'orbitas que contienen un
entorno de puntos donde el comportamiento es escencialmente el
mismo. Por otro lado, Andronov y Pontryagin seguidos por Peixoto,
estudiaron sistemas estables, es decir, aquellos cuyas
perturbaciones verifican que su estructura din\'amica es la misma.
Fue quiz\'as Smale (\cite{SmaleBulletin}) el primero en retomar la
propuesta de Poincar\'e d\'andole una formulaci\'on m\'as precisa:

El objetivo es fijar una variedad cerrada $M$ de dimensi\'on $d$ y
entender la din\'amica de un subconjunto \emph{grande} de
$\Diff^r(M)$, el conjunto de difeomorfismos de $M$ munido de la
topolog\'ia de la convergencia uniforme hasta orden $r$.

Adem\'as, propuso que un conjunto fuese considerado grande si era
abierto y denso, o en su defecto residual o denso (de esta forma
conjunto abierto de sistemas puede ser despreciado). No se har\'a
menci\'on a otras posibles formas de considerar un subconjunto
como grande ni se discutir\'an las razones por las cuales
considerar estas nociones (referimos al lector a \cite{B} o
\cite{Crov-Hab} por m\'as fundamentaci\'on).

Smale propuso un programa que se centr\'o en el estudio de la
estabilidad estructural y en la esperanza de que si bien las
din\'amicas t\'ipicas pod\'ian no ser estables (pod\'ian ser
ca\'oticas) estas ser\'ian estables desde el punto de vista de que
sus propiedades din\'amicas persistir\'ian frente a
perturbaciones. De esta forma, se podr\'ia describir la din\'amica
a trav\'es de m\'etodos simb\'olicos o estad\'isticos. Palis y
Smale \cite{PalisSmale} conjeturaron que los sistemas
estructuralmente estables coincid\'ian con los difeomorfismos
hiperb\'olicos.

La hiperbolicidad fue entonces el paradigma. Robbin y Robinson
(\cite{Robbin,RobinsonEstabilidad}) probaron que los sistemas
hiperb\'olicos eran estables. Finalmente Ma\~ne en un c\'elebre
resultado (\cite{ManheIHES}) completo la caracterizaci\'on de la
las din\'amicas estables en topolog\'ia $C^1$ (en otras topologias
a\'un es desconocido). Describir la din\'amica de los
difeomorfismos hiperb\'olicos fue la tarea que concentr\'o la
mayor atenci\'on en esos a\~nos 60 y principios de los 70.

Desde el punto de vista de comprender la din\'amica de los
difeomorfismos hiperb\'olicos el proyecto fue sumamente exitoso.
No s\'olo se logr\'o comprender cabalmente los aspectos
semilocales de su din\'amica (a traves del estudio simb\'olico y
erg\'odico de sus propiedades), sino que se obtuvieron resultados
profundos acerca de aspectos globales y de la topolog\'ia de sus
piezas b\'asicas.

Sobre los aspectos semilocales, sin ser exhaustivos, se hace
menci\'on a las contribuciones de Bowen, Newhouse, Palis, Sinai,
Ruelle y el propio Smale. Se refiere al lector a \cite{KH} Part 4
por una visi\'on panor\'amica de gran parte de la teor\'ia.

Por otro lado, los aspectos globales fueron fundamentalmente
asociados a los trabajos de Anosov, Bowen, Franks, Shub, Smale,
Sullivan y Williams y buena parte de estos se pueden apreciar en
el libro \cite{Franks-Homology} que contiene una muy linda
presentaci\'on de los trabajos conocidos acerca de la din\'amica
global de difeomorfismos hiperb\'olicos. Por diferentes razones
estas personas abandonaron estos temas, lo que puede ser una
explicaci\'on de por qu\'e estos resultados son menos conocidos.

El programa de Smale, as\'i como la esperanza de que los sistemas
estructuralmente estables fueran t\'ipicos en el espacio de
difeomorfismos de una variedad, cay\'o cuando empezaron a aparecer
ejemplos de din\'amicas robustamente no hiperb\'olicas y no
estables como los de Abraham-Smale (\cite{AS}) y de Newhouse
(\cite{Newhouse}).

Estos ejemplos dieron lugar a la teor\'ia de bifurcaciones donde
Newhouse, Palis y Takens (entre otros) fueron pioneros. Luego de
muchos trabajos al respecto los programas iniciales fueron
adaptados para contemplar dichas bifurcaciones manteniendo la
filosof\'ia inicial de Smale. El programa de Palis
\cite{PalisGlobalView} sin embargo, tiene un enfoque
principalmente semilocal.

Ca\'ido el paradigma de la hiperbolicidad se empez\'o a buscar
nociones alternativas como la hiperbolicidad no uniforme (por
parte de la escuela rusa, fundamentalmente Katok y Pesin, ver
\cite{KH} Supplement S), o la hiperbolicidad parcial
(independientemente por Hirsch-Pugh-Shub \cite{HPS} y Brin-Pesin
\cite{BrinPesin}). Este trabajo se interesa fundamentalmente por
esta segunda generalizaci\'on dada su condici\'on de estructura
geom\'etrica (en contraposici\'on a la condici\'on de propiedad
medible de la hiperbolicidad no uniforme), y su fuerte
vinculaci\'on con las propiedades robustas de la din\'amica.
V\'ease \cite{BDV} por un panorama general de la din\'amica m\'as
all\'a de la hiperbolicidad.

En su b\'usqueda de la prueba de la conjetura de estabilidad,
Ma\~ne (independientemente lo hicieron tambi\'en Pliss y Liao
\cite{pliss,Liao}) introdujo el concepto de descomposici\'on
dominada, y lo que es m\'as importante, mostr\'o su relaci\'on con
la din\'amica de la aplicaci\'on tangente sobre las \'orbitas
peri\'odicas.

Cuando se estudia el espacio de difeomorfismos con la topolog\'ia
$C^1$, gracias a las t\'ecnicas de perturbaci\'on de \'orbitas
desarrolladas desde los a\~nos 60 por Pugh, Ma\~ne, Hayashi y
m\'as recientemente por Bonatti y Crovisier, sabemos que de alguna
manera las \'orbitas peri\'odicas de los difeomorfismos
gen\'ericos capturan la din\'amica de los difeomorfismos
(topol\'ogica y estad\'isticamente). Ver \cite{Crov-Hab} por un
survey sobre estos temas.

Recientemente Bonatti \cite{B} propuso un programa realista para
la comprensi\'on de un conjunto grande de difeomorfismos con la
topolog\'ia $C^1$. Este programa extiende y complementa el
programa general de Palis ya mencionado. Tambi\'en en este caso el
programa tiene un punto de vista semilocal.

Desde el punto de vista global hay menos trabajo hecho y menos
propuestas de trabajo (ver \cite{PujSamHandbook} section 5),
aunque en dimensi\'on tres hay algunas ideas de c\'omo proceder en
ciertos casos.

En lo que sigue, se presentan las contribuciones de esta tesis
pretendiendo mostrar como \'estas encajan en este panorama
subjetivo del desarrollo de la teor\'ia.

\subsection{Atractores en din\'amica $C^1$-gen\'erica}

Siempre es posible descomponer la din\'amica de un homeomorfismo
de un espacio m\'etrico compacto en sus clases de recurrencia.
Este es el contenido de la teor\'ia de Conley \cite{Conley}.

Esta descomposici\'on es muy \'util en la comprensi\'on de las
din\'amicas $C^1$-gen\'ericas\footnote{Utilizaremos esta
expresi\'on para referirnos a difeomorfismos pertenecientes a un
conjunto residual de $\Diff^1(M)$ con la topolog\'ia  $C^1$.},
debido a un resultado de Bonatti y Crovisier (\cite{BC}) que
garantiza que es posible detectar las clases de recurrencia de un
difeomorfismo gen\'erico a partir de sus \'orbitas peri\'odicas.
En buena medida, la comprension de la din\'amica alrededor de
dichas \'orbitas peri\'odicas se lleva toda la atenci\'on y se
espera que puedan describir la din\'amica de dichos difeomorfismos
(ver \cite{B}).

Cuando queremos entender la din\'amica de la mayor\'ia de los
puntos, existen clases de recurrencia que se destacan sobre otras.
Los \emph{quasi-atractores} son clases de recurrencia que admiten
una base de entornos $U_n$ que verifican que $f(\overline{U_n})
\en U_n$. Estas clases siempre existen. Fue probado en \cite{BC}
que existe un conjunto residual de puntos de la variedad que
convergen a dichas clases cuando se trata de un difeomorfismo
gen\'erico.

Algunas veces es posible demostrar que estos quasi-atractores
est\'an aislados del resto de las clases de recurrencia, en ese
caso, decimos que son \emph{atractores}. Los atractores tienen la
propiedad de ser acumulados por la \'orbita futura de los puntos
cercanos y ser din\'amicamente indescomponibles. Determinar la
existencia de atractores y sus propiedades topol\'ogicas y
estad\'isticas es uno de los grandes problemas en sistemas
din\'amicos. En dimensi\'on dos es posible demostrar que para la
mayor parte de los sistemas din\'amicos, vistos con la topolog\'ia
$C^1$, existen atractores. Esto fue demostrado originalmente por
Araujo \cite{Araujo} aunque la demostraci\'on conten\'ia un error
y nunca fue publicada\footnote{Ver \cite{BrunoSantiago} por una
correcci\'on a la prueba original utilizando los resultados de
Pujals y Sambarino \cite{PujSamAnnals1}.}. El resultado fue de
alguna manera incorporado al folklore (ver por ejemplo
\cite{BLY}). Esta tesis presenta una prueba del siguiente
resultado aparecida por primera vez en
\cite{PotExistenceOfAttractors} (ver secci\'on
\ref{Section-AtractoresSuperficies}).

\begin{teoe*}
Existe un conjunto abierto y denso $\cU$ del espacio de
difeomorfismos de una superficie con la topolog\'ia $C^1$ tal que
si $f\in \cU$ entonces $f$ tiene un atractor hyperb\'olico. M\'as
a\'un, si $f$ no puede ser perturbado para tener infinitos puntos
peri\'odicos atractores (pozos), entonces $f$ y sus perturbados
verifican que poseen finitos quasi-atractores y estos son
atractores hiperb\'olicos.
\end{teoe*}

La hiperbolicidad de un quasi-atractor garantiza que este debe ser
un atractor. Para mostrar que cuando hay robustamente finitos
pozos todos los quasi-atractores son hiperb\'olicos, es clave
demostrar que dichos quasi-atractores poseen lo que se llama una
\emph{descomposici\'on dominada}: Existe una descomposici\'on del
fibrado tangente sobre el quasi-atractor en dos subfibrados
$Df$-invariantes que verifican una condici\'on uniforme de
dominaci\'on de uno sobre el otro (la contracci\'on de los
vectores en uno de los fibrados es uniformemente menor que en el
otro fibrado).

La descomposici\'on dominada, as\'i como otras varias posibles
estructuras geom\'etricas $Df$-invariantes, es una herramienta
importante para el estudio de las propiedades din\'amicas de una
clase de recurrencia, y fundamentalmente, para entender como dicha
clase es acumulada por otras clases (ver Secci\'on
\ref{Section-EstructurasInvariantes}).

Buscando comprender la din\'amica cerca de los quasi-atractores
para las din\'amicas $C^1$-gen\'ericas, se obtuvo el siguiente
resultado parcial acerca de la estructura de aquellos
quasi-atractores que son clases homocl\'inicas (ver Secci\'on
\ref{Section-GenericBiLyapunov} y \cite{PotGenericBiLyapunov}):

\begin{teoe*}
Para un difeomorfismo $C^1$ gen\'erico $f$, si $H$ es un
quasi-attractor que contiene un punto peri\'odico $p$ tal que el
diferencial de $f$ sobre $p$ en el per\'iodo contrae volumen,
entonces $H$ admite una descomposici\'on dominada no trivial.
\end{teoe*}

Es posible ver mediante ejemplos (\cite{BV}) que la conclusi\'on
del teorema es en cierto sentido \'optima. La hip\'otesis de que
existan puntos peri\'odicos en la clase es tambi\'en necesaria
(ver \cite{BDihes}). Sin embargo, la hip\'otesis acerca de la
disipatividad del diferencial sobre la \'orbita peri\'odica parece
ser consecuencia de las otras hip\'otesis, pero no fue posible
eliminarla. Demostrarlo parecer\'ia necesitar de alg\'un resultado
del tipo \emph{ergodic closing lemma en la clase homocl\'inica}
que es un problema que a\'un no se logra entender correctamente
(ver \cite{B} Conjecture 2).

La mayor novedad en la prueba del teorema es que se utiliza un
nuevo resultado perturbativo debido a Gourmelon
\cite{GouFranksLemma} que permite perturbar el diferencial de una
\'orbita peri\'odica con cierto control de las variedades
invariantes luego de la perturbaci\'on. El uso de dicho resultado
combinado con la estabilidad Lyapunov nos permite resolver el
problema de garantizar que un punto pertenece a la clase luego de
la perturbaci\'on. El resultado responde positivamente a una
pregunta realizada en \cite{ABD} (Problem 5.1).

El sue\~no de que las din\'amicas $C^1$-gen\'ericas tuviesen
atractores\footnote{Ver la introducci\'on de \cite{BLY}.} cay\'o
recientemente debido a un ejemplo sorprendente debido a \cite{BLY}
que muestra como el desarrollo reciente de la teor\'ia de la
din\'amica $C^1$-gen\'erica ha tenido una gran influencia en la
forma de entender la din\'amica y ha simplificado preguntas que
resultaban a simple vista inabordables.

Los ejemplos de \cite{BLY} poseen lo que llamaron \emph{atractores
esenciales} (ver subsecci\'on \ref{SubSection-AttractingSets}) y
no es a\'un sabido si poseen atractores en el sentido de Milnor.
En la Secci\'on \ref{Section-EjemplosQuasiAtractores} revemos
estos ejemplos y presentamos algunos ejemplos nuevos de
\cite{PotWildMilnor} sobre los cuales tenemos una mejor
comprensi\'on de c\'omo las otras clases se acumulan a su
din\'amica.

\begin{teoe*}
Existe un abierto $\cU$ de $\Diff^1(\TT^3)$ de difeomorfismos
tales que:
\begin{itemize}
\item[-] Si $f\in \cU$ entonces $f$ tiene un \'unico
quasi-atractor y un atractor en el sentido de Milnor. Adem\'as, si
$f$ es de clase $C^2$ posee una \'unica medida SRB cuya cuenca es
de medida total. \item[-] Para $f$ en un residual de $\cU$ se
verifica que $f$ no tiene atractores. \item[-] Las clases de
recurrencia diferentes del quasi-atractor est\'an contenidas en
superficies peri\'odicas.
\end{itemize}
\end{teoe*}

El \'ultimo punto del teorema entra en contraste con los nuevos
resultados obtenidos por Bonatti y Shinohara (\cite{BS}) y en la
Secci\'on \ref{Section-TrappingAttractors} especulamos acerca de
c\'omo ambos resultados podr\'ian llegar a entrar en una misma
teor\'ia.

\subsection{Hiperbolicidad parcial en el toro $\TT^3$}

As\'i como la secci\'on anterior trato impl\'icitamente problemas
que son de naturaleza semilocal (a pesar de que se pueden hacer
preguntas de \'indole global acerca de la existencia de atractores
y de su topolog\'ia) esta secci\'on trata fundamentalmente acerca
de problemas de din\'amica global.

Un conocido teorema en din\'amica diferenciable, que reune
resultados cl\'asicos de Ma\~ne (\cite{ManheErgodic}) y Franks
(\cite{FranksAnosov}) dice lo siguiente:

\begin{teoe*}[Ma\~ne-Franks] Sea $M$ una superficie compacta.
Entonces, las tres propiedades siguientes para $f \in \Diff^1(M)$
son equivalentes:
\begin{itemize}
\item[(i)] $f$ es $C^1$-robustamente transitivo. \item[(ii)] $f$
es Anosov. \item[(iii)] $f$ es Anosov y conjugado a un
difeomorfismo de Anosov lineal.
\end{itemize}
\end{teoe*}

Escencialmente, Ma\~ne prob\'o la implicancia (i) $\Rightarrow$
(ii) y Franks la implicancia (ii) $\Rightarrow$ (iii). La robustez
de los Anosov y la transitividad de los Anosov lineales da (iii)
$\Rightarrow$ (i).

Si se interpreta el ser difeomorfismo de Anosov como que el
diferencial de $f$ preserve una estructura geom\'etrica, podemos
de alguna manera identificar el resultado (i) $\Rightarrow$ (ii)
como diciendo que ``una propiedad din\'amica robusta fuerza la
existencia de una estructura geom\'etrica $Df$-invariante''. De
hecho, en vista que las perturbaciones $C^1$ no pueden romper la
propiedad din\'amica, no es sorprendente que dicha estructura
geom\'etrica tenga propiedades de rigidez frente a perturbaciones.

Por otro lado, la direcci\'on (ii) $\Rightarrow$ (i) se puede ver
como diciendo que la existencia de una estructura geom\'etrica
invariante est\'a tambi\'en relacionada a la existencia de una
propiedad din\'amica robusta, en este caso, la transitividad.

En dimensiones mayores la relaci\'on entre las propiedades
din\'amicas robustas y las estructuras geom\'etricas invariantes
est\'a menos desarrollada, aunque existen resultados en la
direcci\'on de obtener una estructura geom\'etrica invariante a
partir de una propiedad din\'amica robusta. En dimensi\'on 3, esto
surge de \cite{DPU} y fue generalizado a dimensiones mayores en
\cite{BDP} (ver Secci\'on \ref{Section-EstructurasInvariantes}):

\begin{teoe*}[Diaz-Pujals-Ures]
Si $M$ es una variedad de dimensi\'on $3$ y $f$ es un
difeomorfismo $C^1$-robustamente transitivo, entonces, $f$ es
parcialmente hiperb\'olico.
\end{teoe*}

Existen diversas definiciones de hiperbolicidad parcial en la
literatura, y estas tambi\'en han variado a lo largo del tiempo.
Nosotros seguimos la definici\'on que se utiliza en \cite{BDV} que
es la que mejor se ajusta a nuestro enfoque (ver Secci\'on
\ref{Section-EstructurasInvariantes} por definiciones precisas).
Un difeomorfismo parcialmente hiperb\'olico sera uno que verifica
que el fibrado tangente se descompone en una suma $Df$-invariante
$TM=E \oplus F$ que verifica una propiedad de dominaci\'on y tal
que uno de los dos fibrados es uniforme. Por simplicidad,
eliminamos la simetr\'ia y consideramos descomposiciones del tipo
$TM =E^{cs} \oplus E^u$ con $E^u$ uniformemente expandido.

Una primera dificultad que aparece si nos interesamos en entender
las propiedades din\'amicas robustas implicadas por la existencia
de una descomposici\'on parcialmente hiperb\'olica, es el no tener
control de la contracci\'on en la direcci\'on centro estable
$E^{cs}$. Esto impide que tengamos una comprensi\'on cabal de la
din\'amica en esa direcci\'on, como si la tenemos en el caso donde
los fibrados son uniformes. Una excepci\'on notable a esto es
\cite{PujSamDominated} donde se estudian las consecuencias
din\'amicas de la descomposici\'on dominada en dimensi\'on $2$.

Es importante mencionar tambi\'en, la importancia del \'item (iii)
en el Teorema de Ma\~ne y Franks. De alguna manera, la idea que
subyace es que para obtener una propiedad din\'amica robusta a
partir de la existencia de una estructura geom\'etrica invariante,
puede ser importante apoyarse en las propiedades topol\'ogicas
impuestas por dicha estructura, tanto en la topolog\'ia de la
variedad como en la clase de isotop\'ia del difeomorfismo. Es por
ello que en vista de ese teorema, y de la dificultad que se ha
tenido para obtener resultados en la direcci\'on rec\'iproca al
Teorema de Diaz-Pujals-Ures, puede ser importante dividir el
estudio en buscar resultados del tipo (ii) $\Rightarrow$ (iii) y
del tipo (ii)+(iii) $\Rightarrow$ (i) emulando el Teorema de
Ma\~ne y Franks. En particular puede ser importante estudiar
propiedades de difeomorfismos parcialmente hiperb\'olicos en
ciertas variedades, o incluso en clases de isotop\'ia fijadas.

Otra dificultad en dimensi\'on $3$ es obtener propiedades
topol\'ogicas a partir de las estructuras invariantes. En
dimensi\'on $2$, el solo hecho de preservar un campo de vectores
impone enormes restricciones en la topolog\'ia de la variedad. En
dimensi\'on $3$, es bien sabido que toda variedad admite un campo
de vectores no nulo, e incluso, una foliaci\'on de codimensi\'on
$1$.

Esta situaci\'on podr\'ia ser considerada muy mala desde el punto
de vista de obtener resultados en el sentido de encontrar
propiedades topol\'ogicas de un difeomorfismo que preserva una
estructura geom\'etrica. Sin embargo, recientemente, una simple
pero poderosa observaci\'on de Brin-Burago-Ivanov (\cite{BBI,BI})
despert\'o nuevamente la esperanza:

\medskip

\begin{obse*}[Brin-Burago-Ivanov]
En una variedad de dimensi\'on 3, si $\cF$ es una foliaci\'on
transversal a la direcci\'on inestable $E^u$ de un difeomorfismo
parcialmente hiperb\'olico $f$, entonces $\cF$ no tiene
componentes de Reeb. \finobs
\end{obse*}

Las componentes de Reeb y su relaci\'on con los difeomorfismos
parcialmente hiperb\'olicos ya hab\'ia sido estudiada, por ejemplo
en \cite{DPU} (Theorem H) y \cite{BonWilk} (Lemma 3.7) en
contextos m\'as restrictivos (asumiendo coherencia din\'amica y
transitividad). Esta observaci\'on es m\'as general, su fuerza
radica en el hecho que depende de la din\'amica s\'olo en que la
foliaci\'on inestable no tiene curvas cerradas.

Por un lado, es conocido que diversas variedades no admiten
foliaciones sin componentes de Reeb. Por otro lado, existe mucha
teor\'ia acerca de su clasificaci\'on (ver por ejemplo
\cite{Plante, Roussarie}) y por lo tanto, de alguna manera nos
hace esperar que es posible entender los difeomorfismos
parcialmente hiperb\'olicos al menos en cierto grado.

Por otro lado, aparece una nueva dificultad ya que no es conocido
si todo difeomorfismo parcialmente hiperb\'olico posee una
foliaci\'on transversal a la direcci\'on inestable. En esta tesis
se propone la noci\'on de \emph{casi coherencia din\'amica} y se
prueba que es una propiedad abierta y cerrada en el espacio de
difeomorfismos parcialmente hiperb\'olicos.

Uno de los teoremas principales (ver Cap\'itulo
\ref{Capitulo-ParcialmenteHiperbolicos} y
\cite{Pot3dimPHisotopic}) garantiza en ciertas clases de
isotop\'ia de difeomorfismos de $\TT^3$ que dicha noci\'on es
suficiente para que el fibrado $E^{cs}$ sea tangente a una
foliaci\'on invariante:

\begin{teoe*}
Sea $f:\TT^{3} \to \TT^3$ un difeomorfismo parcialmente
hiperb\'olico que es casi dinamicamente coherente y es isot\'opico
a un difeomorfismo de Anosov lineal. Entonces $f$ es
din\'amicamente coherente.
\end{teoe*}

Ser \emph{din\'amicamente coherente} significa justamente que el
fibrado $E^{cs}$ sea tangente a una foliaci\'on $f$-invariante.

En el caso \emph{parcialmente hiperb\'olico fuerte} (es decir,
donde $TM =E^s \oplus E^{c} \oplus E^u$ es una descomposici\'on
$Df$-invariante con propiedades de dominaci\'on y donde $E^s$ y
$E^u$ son uniformes) podemos decir m\'as.

En \cite{BI} fue probado que todo difeomorfismo parcialmente
hiperb\'olico fuerte es casi-din\'amicamente coherente y esto fue
aprovechado por, primero \cite{BI} (siguiendo \cite{BBI}) y luego
\cite{Parwani} para dar condiciones topol\'ogicas que estos deben
satisfacer. Por otro lado, esta condici\'on permite esperar que
una conjetura de Pujals sea atacable:

\begin{conje*}[Pujals \cite{BonWilk}]
Sea $f:M \to M$ con $M$ variedad de dimensi\'on $3$ un
difeomorfismo parcialmente hiperb\'olico fuerte y transitivo.
Entonces, tenemos las siguientes posibilidades (m\'odulo
considerar levantamientos finitos):
\begin{itemize}
\item[-] $f$ es conjugado por hojas a un difeomorfismo de Anosov
en $\TT^3$. \item[-] $f$ es conjugado por hojas a un skew-product
sobre un difeomorfismo de Anosov de $\TT^2$ (entonces la variedad
es $\TT^3$ o una nilvariedad). \item[-] $f$ es conjugado por hojas
al tiempo $1$ de un flujo de Anosov.
\end{itemize}
\end{conje*}

\'Ultimamente ha habido progreso en la direcci\'on de esta
conjetura. Para empezar, el trabajo de \cite{BonWilk} hace un
avance considerable sin dar ninguna hip\'otesis acerca de la
topolog\'ia de la variedad asumiendo la existencia de curvas
cerradas tangentes a la direcci\'on central. Luego, los trabajos
de Hammerlindl \cite{Hammerlindl,HNil} dan una prueba a la
conjetura en caso que la variedad sea $\TT^3$ o una nilvariedad
pero con una definici\'on m\'as restrictiva de hiperbolicidad
parcial. Si bien esta definici\'on es verificada por varios
ejemplos, es de alguna manera artificial y no encaja bien con el
resultado de \cite{DPU}.

Por supuesto, para conseguir probar la conjetura de Pujals, un
paso previo es mostrar que los difeomorfismos de ese tipo son
din\'amicamente coherentes ya que la definici\'on de conjugaci\'on
por hojas (ver Secci\'on \ref{Section-HIRSHPUGHSHUB}) requiere la
existencia de foliaciones invariantes. El trabajo de Hammerlindl
se apoya en trabajos previos de Brin-Burago-Ivanov \cite{BBI2} que
muestran que en $\TT^3$, con esta definici\'on restrictiva de
hiperbolicidad parcial se tiene coherencia din\'amica.

Por otro lado, recientemente apareci\'o un ejemplo debido a
Rodriguez Hertz-Rodriguez Hertz-Ures (\cite{HHU}) de un
difeomorfismo parcialmente hiperb\'olico de $\TT^3$ que no admite
foliaciones invariantes. En esta tesis se completa el panorama
(\cite{Pot3dimPHisotopic}) mostrando que:

\begin{teoe*}
Sea $f: \TT^3 \to \TT^3$ un difeomorfismo parcialmente
hiperb\'olico fuerte que no admite un toro peri\'odico normalmente
atractor ni un toro peri\'odico normalmente repulsor, entonces,
$f$ es din\'amicamente coherente.
\end{teoe*}

Esto responde a una conjetura de Rodriguez Hertz-Rodriguez
Hertz-Ures en $\TT^3$ y tambi\'en, por lo obtenido en la prueba,
permite responder a la conjetura de Pujals en $\TT^3$ (ver
\cite{HP}).

Sobre el final del Cap\'itulo
\ref{Capitulo-ParcialmenteHiperbolicos} tambi\'en se obtienen
algunos resultados en dimensiones mayores an\'alogos a los
cl\'asicos resultados de Franks, Newhouse y Manning para
difeomorfismos de Anosov en el contexto de parcialmente
hiperb\'olicos.

\subsection{Otras contribuciones}

En esta secci\'on se describen otras contribuciones de la tesis.

Por un lado, en la Secci\'on
\ref{SubSection-LocalizacionDeClasses} se presenta un mecanismo
para la localizaci\'on de clases de recurrencia aparecido en
\cite{PotWildMilnor} que se ha considerado tiene valor en si mismo
ya que puede ser aplicado en diferentes contextos (en esta tesis
se utiliza en la Secci\'on \ref{Section-DescDomT2} as\'i como en
las subsecciones \ref{SubSection-EjemploDA} y
\ref{SubSection-EjemploPlykin}).

En la Secci\'on \ref{Section-EjemplosQuasiAtractores} se presentan
diversos ejemplos de quasi-atractores y de difeomorfismos
robustamente transitivos algunos de los cuales son modificaciones
de ejemplos conocidos pero igual consideramos que pueden
representar un aporte al entendimiento de estos fen\'omenos.

Luego, en el Cap\'itulo \ref{Capitulo-FoliacionesEHipParcial}
donde se preparan los resultados sobre foliaciones, que luego
ser\'an utilizados en el Cap\'itulo
\ref{Capitulo-ParcialmenteHiperbolicos}, se obtienen resultados
que pueden tener inter\'es independiente. En particular, vale
mencionar un resultado cuantitativo sobre la existencia de
estructura de producto global para foliaciones presentado en la
Secci\'on \ref{Section-GlobalProductStructure}. Tambi\'en, en la
Secci\'on \ref{Section-DescDomT2}, se da una clasificaci\'on de la
din\'amica de los difeomorfismos globalmente parcialmente
hiperb\'olicos en $\TT^2$ que de alguna manera muestra en un
contexto sencillo lo que se har\'a despu\'es en el Cap\'itulo
\ref{Capitulo-ParcialmenteHiperbolicos}.

Tambi\'en se incluyen en la tesis cuatro ap\'endices donde se
presentan resultados que est\'an desviados del cuerpo central de
la tesis. Se hace particular \'enfasis en el Ap\'endice
\ref{Apendice-PseudoRotaciones}, basado en \cite{PotRecurrence},
donde se prueba un resultado acerca de homeomorfismos del toro con
un \'unico vector de rotaci\'on; y en el Ap\'endice
\ref{Apendice-BCGP}, basado en \cite{BCGP}, donde se presenta un
trabajo conjunto con Bonatti, Crovisier y Gourmelon. En ese
trabajo se estudian bifurcaciones de clases de recurrencia
robustamente aisladas y se dan ejemplos de ese tipo de clases que
no son robustamente transitivas, respondiendo as\'i a una pregunta
de \cite{BC}.

\section{Introduction (Français)}

\subsection{Contexte historique}

On peut\footnote{L'introduction historique que nous pr\'esentons
est subjective et peut ne pas d\'ecrire avec exactitude les faits
historiques. On peut penser que ce qui suit est une possible
explication historique des raisons pour lesquels les sujets ici
pr\'esent\'es ont \'et\'e abord\'es.} dire que l'objectif
principal des syst\`emes dynamiques est de comprendre le
comportement assymptotique d'une loi d'\'evolution avec certaines
conditions initiales dans un espace de configurations. Au d\'ebut,
on \'etudiait les equations diff\'erentielles ordinaires du type

$$ \dot x = X(x)  \qquad  X: \RR^n \to \RR^n $$

\noindent et on cherche \`a r\'esoudre analytiquement
l'\'equation, en essayant de trouver, pour chaque valeur possible
de la condition initiale $x_0\in \R^n$, une solution explicite
$\varphi_t(x_0)$.

On a vite remarqu\'e qu'\`a partir d'\'equations extr\'emement
simples on obtenait des solutions compliqu\'ees qui n'aidaient pas
\`a la compr\'ehension de la loi d'evolution, m\^eme apr\`es
\^etre int\'egr\'ees.

Poincar\'e, int\'er\'ess\'e par l'\'etude du fameux probl\`eme des
3 corps (\cite{Poincare}), a \'et\'e pe\^ut-\^etre le premier \`a
proposer que l'\'etude de l'\'evolution se fasse du point de vue
qualitatif, et d'une certaine façon, il a sugg\'erer d'\'etudier
le comportement de la ``plupart" des orbites pour la ``plupart"
des syst\`emes.

Au début, la recherche se centrait autour de la stabilité du
système. Lyapunov a étudié les orbites stables, qui contiennent un
voisinage de points où le comportement est essentiellement le
même. D'autre part, Andronov et Pontryagin, suivis par Peixoto,
ont étudié les systèmes stables, c'est à dire, ceux dont la
structure dynamique ne change pas sous des perturbations. C'est
peût-être Smale le premier qui a repris la suggérence de Poincaré,
en présentant une formulation plus précise: l'objectif est de
fixer une variété fermée $M$ de dimension $d$ et de comprendre la
dynamique d'un sous-ensemble ``grand" de $Diff^r(M)$, l'ensemble
des difféomorphismes de $M$ muni de la topologie de la convergence
uniforme jusqu'à l'ordre $r$.

Aussi, Smale a proposé qu'un sous-ensemble soit consideré "grand"
s'il est soit ouvert et dense, ou bien s’il est à résidual ou
m\^eme dense (ainsi les ensembles ouverts de systèmes peuvent être
négligés). On ne faira pas allusion à d’autres possibles
définitions de ``grand” et on ne discutera pas les raisons de
notre choix de topologie (voir \cite{B} ou \cite{Crov-Hab} pour
des explications detailles).

Smale a proposé un programme centré autour de l’étude de la
stabilité structurelle. Bien que les dynamiques typiques peuvent
ne pas être stables (elles peuvent même être chaotiques), il
espérait que leurs propriétes dynamiques soient robustes aux
perturbations. Ainsi, Palis et Smale \cite{PalisSmale} ont
conjecturé que les systèmes structurellement stables coïncident
avec les difféomorphismes hyperboliques.

L’hyperbolicité est devenue donc le nouveau paradigme. Robbin et
Robinson (\cite{Robbin,RobinsonEstabilidad}) ont montré que les
systèmes hyperboliques sont stables. Finalement, Mañé
(\cite{ManheIHES}) a characterisé les dynamiques stables en
topologie $C^1$ (on ne sait pas encore ce qui arrive dans d’autres
topologies). Pendant les années 60 et début des 70, la plupart de
l’attention s’est concentrée autour de la description de la
dynamique des difféomorphismes hyperboliques.

Du point de vue de la compréhension de la dynamique des
difféomorphismes hyperboliques, le programme a été très réussi: il
a servi à la compréhension profonde des aspects semi-locaux de la
dynamique (à travers l’étude symbolique et ergodique de ses
propriétés) mais aussi à l’obtention d’importants résultats sur
les aspects globaux et la topologie des pièces basiques.

En ce qui concerne les propriétés semi-locales, nous citons les
contributions de Bowen, Newhouse, Palis, Sinai, Ruelle et Smale.
Pour une vision panoramique de la théorie, voir \cite{KH} Partie
4.

D’un autre coté, les aspects globaux sont associés
fondamentalement aux travaux de Anosov, Bowen, Franks, Shub,
Sullivan et Williams, dont la plupart sont présentés dans le très
beau livre \cite{Franks-Homology}, qui rassemble les travaux
connus sur la dynamique globale des difféomorphismes
hyperboliques. Pour de différentes raisons, ces auteurs ont
abandonné l’étude de ces sujets, ce qui peut expliquer qu’ils
soient peu connus.

Le programme de Smale a échoué avec l’apparition d’exemples de
dynamiques robustement non hyperboliques et non structuralement
stables comme ceux de Abraham-Smale (\cite{AS}) et Newhouse
(\cite{Newhouse}).

Ces exemples ont été à la base de la théorie des bifurcations, de
laquelle Newhouse, Palis et Takens (entre autres) ont été
pionniers. Après plusieurs travails à ce sujet, les premiers
programmes  ont été adaptés afin de, tout en considérant les
bifurcations, maintenir la philosophie initiale de Smale. Malgré
cela, l’approche du programme de Palis \cite{PalisGlobalView} est
surtout semi-locale.

Après la chute du paradigme de l’hyperbolicité, on a commencé à
chercher des notions alternatives comme l’hyperbolicité non
uniforme (due à l’école russe, fondamentalement à Katok et Pesin,
voir \cite{KH} Supplément S) ou l’hyperbolicité partielle
(considérée de façon indépendente par Hirsch, Pugh et Shub
\cite{HPS} et par Brin et Pesin \cite{BrinPesin}). C’est à cette
deuxième géneralisation que s’intéresse surtout cette thèse. Nous
sommes intéréssés par l’approche géométrique de cette
géneralisation (qui contraste avec la condition mesurable de
l’hyperbolicité non uniforme) et par son lien étroit avec les
propriétés robustes de la dynamique. Le lecteur peut se référer à
\cite{BDV} pour une vision générale de la dynamique au delà de
l’hyperbolicité.

En cherchant la preuve de la conjecture de stabilité, Mañe (aussi
Pliss et Liao \cite{pliss,Liao}, de façon indépendente) a
introduit la notion de décomposition dominée et, ce qui est encore
plus important, a montré la relation de la dynamique de
l’application tangente sur les orbites périodiques.

En ce qui concerne l’étude de l’espace des difféomorphismes  sous
la topologie $C^1$, nous savons, grâce aux techniques de
perturbation d’orbites, développées depuis les années 60 par Pugh,
Mañe, Hayashi et plus récemment par Bonatti et Crovisier, que les
orbites périodiques capturent, d’une certaine façon, la dynamique
des difféomorphismes génériques (topologique et statistiquement).
Voir \cite{Crov-Hab} pour une vision générale sur ces sujets.

Récemment Bonatti \cite{B} a proposé un programme réaliste pour la
compréhension d’un grand ensemble de difféomorphismes sous la
topologie $C^1$.  Ce programme prolonge et complète, aussi d’un
point de vue semi-locale, le susmentionné programme général de
Palis.

Du point de vue globale, il y a justes quelques idées sur comment
procéder en quelques cas de dimension 3.

En ce qui suit, nous présentons les contributions de cette thèse
et nous essayons de montrer comment celles-ci s’adaptent à ce
panorama subjectif du développement de la théorie.

\subsection{Attracteurs en dynamique $C^1$-générique.}

Il est toujours possible de décomposer la dynamique d’un
homéomorphisme d’un espace métrique compacte en classes de
récurrence, comme explique la théorie de Conley \cite{Conley}.

Cette décomposition est très utile pour la compréhension des
dynamiques $C^1$-génériques\footnote{Nous utiliserons cette
expression pour faire allusion aux difféomorphismes appartenant à
un ensemble résiduel de $\Diff^1(M)$ sous la topologie $C^1$.},
grâce à un résultat de Bonatti et Crovisier (\cite{BC}) qui assure
que les orbites périodiques sont suffisantes pour détecter les
classes de récurrence d’un difféomorphisme générique. Dans une
bonne mesure, l’attention se place sur la compréhension de la
dynamique autour desdites orbites et on espère qu’elles décrivent
la dynamique des difféomorphismes en question (voir \cite{B}).

Lorsque l’on veut comprendre  la dynamique de la plupart des
points, on trouve des classes de récurrence distinguées: les
\emph{quasi-attracteurs} sont des classes de récurrence qui
admitent une base de voisinages $U_n$ qui vérifient
$f(\overline{U}_n)\en U_n$. Les quasi-attracteurs existent
toujours. En \cite{BC}, les auteurs ont montré que, dans le cas
d’un difféomorphisme générique,  il existe un ensemble résiduel de
points de la variété qui convergent aux dites classes.

Il est possible dans certains cas de prouver que les
quasi-attracteurs sont isolés du reste des classes de récurrence.
Dans ces cas, on les appelle \emph{attracteurs}. Les attracteurs
ont la propriete d'\^etre accumul\'es par l'orbite future des
points proches et d'\^etre dynamiquement ind\'ecomposables. Un des
grands problèmes des systèmes dynamiques consiste à déterminer
leur existence et leur propriétés topologiques et statistiques. En
dimension 2, il est possible de prouver que pour la plupart des
systèmes dynamiques, sous la topologie $C^1$, il en existent. Ceci
a été prouvé d’abord par Araujo \cite{Araujo}, mais la preuve
contenait une erreur et n’a jamais été publiée\footnote{Voir
\cite{BrunoSantiago} pour une correction de la preuve originelle
qui utilise les résultats de Pujals et Sambarino
\cite{PujSamAnnals1}.}. Cette thèse présente une preuve du
résultat suivant, apparue pour la première fois en
\cite{PotExistenceOfAttractors} (voir section
\ref{Section-AtractoresSuperficies}).

\begin{teof*}
Il existe un ensemble ouvert et dense $\cU$ de l’espace de
difféomorphismes d’une surface sous la topologie $C^1$ tel que
tout $f\in \cU$ a un attracteur hyperbolique. Si $f$ ne peut pas
être perturbé de façon à obtenir un nombre infini de points
périodiques attracteurs (puits), alors pour $f$ et ses perturbés
il y a un nombre fini d’attracteurs dont tous sont hyperboliques.
\end{teof*}

Tout quasi-attracteur hyperbolique est un attracteur. Pour montrer
que tout quasi-attracteur est hyerbolique, dans le cas d’un nombre
fini de puits, c’est fondamental de montrer d'abord que les dits
attracteurs possèdent ce qu’on appelle une \emph{décomposition
dominée}: une décomposition du fibré tangent sur le
quasi-attracteur en deux sous-fibrés $Df$-invariants qui vérifient
une condition uniforme de domination de l’un sur l’autre (la
contraction des vecteurs sur un des fibrés est uniformément plus
petit que la contraction sur l’autre).

La décomposition dominée, ainsi que d’autres possibles structures
géométriques $Df$-invariantes, est un outil important pour l’étude
des propriétés dynamiques d’une classe de récurrence et surtout
pour comprendre comment ladite classe est accumulée par d’autres
classes (voir Section \ref{Section-EstructurasInvariantes}).

En quête de comprendre la dynamique proche aux quasi-attracteurs
pour les systèmes $C^1$-génériques, le résultat partiel suivant,
concernant la structure des quasi-attracteurs qui correspondent à
des classes homocliniques, a été obtenu (voir Section
\ref{Section-GenericBiLyapunov} et \cite{PotGenericBiLyapunov}):

\begin{teof*}
Soit $f$ un difféomorphisme $C^1$-générique. Si $H$ est un
quasi-attracteur de $f$ qui contient un point périodique $p$ pour
lequel la différentiel \`a la periode contract le volume, alors
$H$ admet une décomposition dominée non triviale.
\end{teof*}

Il est possible de voir à travers des examples (\cite{BV}) que le
théorème est optimal dans un certain sens. Aussi, l’hypothèse de
l’existence de points périodiques est nécessaire (voir
\cite{BDihes}). Par contre, l’hypothèse sur la dissipativité du
différentiel sous l’orbite périodique semble être une conséquence
des autres hypothèses, bien qu’elle n’a pas pu être éliminée. Pour
le faire, il semble nécessaire d’avoir un résultat du type ergodic
closing lemma sur la classe homoclinique, ce qui constitue un
problème qui n’est pas encore compris correctement (voir \cite{B}
Conjecture 2).

La nouveauté de la preuve de ce théorème se base surtout dans le
fat qu’elle utilise un nouveau résultat dû à Gourmelon
\cite{GouFranksLemma}, qui permet de perturber le différentiel de
l’orbite périodique en gardant un certain contrôle des variétés
invariantes. L’utilisation de ce résultat et la stabilité Lyapunov
du système permettent de garantir que le point reste dans la
classe après une perturbation. Ce resultat répond positivement à
une question posée en \cite{ABD} (Probl\`eme 5.1).

Le r\^eve des dynamiques $C^1$-g\'en\'eriques admettant des
attracteurs\footnote{Voir l’introduction de \cite{BLY}.} a disparu
à cause d’un exemple surprenant dû à \cite{BLY} qui montre que le
dévelopement récent de la dynamique $C^1$-générique a eu une
grande influence sur la façon de comprendre la théorie et a
simplifié des questions qui résultaient à simple vue intraitables.

Les exemples de \cite{BLY} possèdent ce que l’on a appelé des
\emph{attracteurs essentiels} (voir la sous-section
\ref{SubSection-AttractingSets}); il n’est pas encore connu s’il
possèdent ou pas des attracteurs dans les sens Milnor. Dans la
Section \ref{Section-EjemplosQuasiAtractores} nous reverrons ces
exemples et nous présenterons de nouveaux exemples qui
apparaissent sur \cite{PotExistenceOfAttractors} et sur lesquels
nous avons une meilleure compréhension de comment les autres
classes accumulent.

\begin{teof*}
Il existe un ouvert $\cU$ de $\Diff^1(\TT^3)$ tel que:

\begin{itemize}
\item[-] Chaque $f\in \cU$ possède un seul quasi-attracteur et un
attracteur au sens Milnor. Aussi, si l’élément est de classe
$C^2$, il possède un seul mesure SRB dont le bassin est de mesure
totale. \item[-] Pour $f$ dans un ensemble r\'esiduel de $\cU$ ne
possède pas des attracteurs, \item[-] Les classes de récurrence
qui ne sont pas un quasi-attracteurs d’un élément $f$ de $\cU$
sont contenues dans des surfaces périodiques.
\end{itemize}
\end{teof*}

La dernière conclusion du théorème contraste avec les nouveaux
résultats obtenus par Bonatti et Shinohara (\cite{BS}). Dans la
Section \ref{Section-TrappingAttractors}, nous discutons comment
les deux résultats pourraient s’intégrer dans une même théorie.

\subsection{Hyperbolicité partielle sur le tore $\TT^3$.}

La section précédente ayant traité de façon implicite des
problèmes de nature semi-locale (bien qu’il est possible de poser
des questions sur l’existence et la topologie des attracteurs du
point de vue global), nous dédions cette section à des problèmes
des dynamiques globales.

À la suite, nous présentons un théorème connu en dynamique
différentiable, qui rassemble des résultats classiques de Mañe
(\cite{ManheErgodic}) et Franks (\cite{FranksAnosov}).

\begin{teof*}[Ma\~ne-Franks]
Soit $M$ une surface compacte. Pour un difféormorphisme $f\in Diff^1(M)$,  les trois conditions suivantes sont équivalentes:
\begin{itemize}
\item[(i)] $f$  est $C^1$-robustement transitif, \item[(ii)] $f$ f
est Anosov, \item[(iii)] $f$  f est Anosov et conjugué à un
difféomorphisme Anosov linéaire du tore $\TT^2$.
\end{itemize}
\end{teof*}

Essentiellement, Ma\~ne a prouvé  (i) $\Rightarrow$ (ii) et Franks
a prouvé (ii) $\Rightarrow$ (iii). L’autre implication se déduit
des faits que les difféomorphismes Anosov sont robustes et que les
linéaires sont en plus transitifs.

Si on interprète la condition Anosov comme le fait que le
différentiel préserve une certaine structure géométrique, le
résultat (i)$\Rightarrow$ (ii) s’identifie, d’une certaine façon,
à la condition suivante:  ``une propriété dynamique robuste
implique l’existence d’une structure géométrique $Df$-invariante”.
En fait, vu que les perturbations $C^1$ ne peuvent pas casser la
dynamique, il n’est pas surprenant que ladite structure
géométrique ait des propriétés de rigidité aux perturbation.

D’autre part, l’implication (ii)$\Rightarrow$ (i) peut
s’identifier à la condition suivante: ``l’existence d’une
propriété géométrique invariante est liée à l’existence d’une
propriété dynamique robuste'' (dans ce cas, la transitivité).

En dimensions supérieures la relation entre les propriétés
dynamiques robustes et les structures géométriques invariantes est
moins développée, bien qu’il existe des résultats visant à obtenir
une structure géométrique invariante à partir d’une propriété
dynamique robuste. En dimension 3, ceci ce déduit de \cite{DPU} et
a été généralisé à dimensions supérieures en \cite{BDP} (voir
Section \ref{Section-EstructurasInvariantes}):

\begin{teof*}[D\'iaz-Pujals-Ures]
Soit $M$  un variété de dimension 3. Tout difféomorphisme
$C^1$-robustement transitif de $M$ est partiellement hyperbolique.
\end{teof*}

Plusieurs définitions d’hyperbolicité partielle existent dans la
littérature, et celles-ci ont changé le long du temps. Nous
suivons la définition utilisée en \cite{BDV} qui est celle qui
s’adapte le mieux à nôtre approche (voir Section
\ref{Section-EstructurasInvariantes} pour les définitions
précises). Un difféomorphisme  partiellement hyperbolique sera tel
que le fibré tangent se décompose en une somme $Df$-invariante $T
M=E\oplus F$ qui vérifie la propriété de domination et telle que
l’un des fibrés de la décomposition est uniforme. Pour simplifier,
nous éliminons la symétrie et considérons les décompositions de la
forme $T M=E^{cs} \oplus E^u$, ou $E^u$ est uniformément dilate.

Si on s’intéresse à comprendre les propriétés dynamiques robustes
impliquées par l’existence d’une décomposition partiellement
hyperbolique, une première difficulté que l’on trouve concerne le
manque de contrôle de la contraction sur la direction centre
stable $E^{cs}$. Ceci empêche de comprendre à fond la dynamique
sur cette direction, comme on la comprend sur la direction des
fibrés uniformes. Une exception notable dans ce sens est
\cite{PujSamDominated}, travail qui étudie les conséquences
dynamiques de la décomposition dominée pour le cas de dimension 2.

Il mérite de mentionner aussi l’importance de (iii) dans le
Théorème de Ma\~ne et Franks. D’une certaine façon, il est basé
sur l’idée qui suggère que pour obtenir une propriété dynamique
robuste à partir de l’existence d’une structure géométrique
invariante, il peut être important de s’appuyer sur les propriétés
topologiques imposées par ladite structure: la topologie de la
variété et la classe d’isotopie du difféomorphisme. Tenant compte
de ce théorème et des difficultés à trouver des résultats dans le
sens réciproque au Théorème de Diaz-Pujals-Ures , il semble
important de diviser le travail en cherchant des résultats tu type
(ii)$\Rightarrow$ (iii) et du type (ii)+(iii) $\Rightarrow$ (i),
en imitant le Théorème de Ma\~ne et Franks. En particulier, il
pourrait être pertinent d’étudier les propriétés des
difféomorphismes partiellement hyperboliques dans certaines
variétés, ou même dans des classes d’isotopie fixées.

Une autre difficulté en dimension 3 est l’obtention de propriétés
topologiques à partir des structures invariantes. En dimension 2,
le seul fait de préserver un camp de vecteurs impose des énormes
restrictions sur la topologie de la variété. En dimension 3, il
est bien connu que toute variété admet un champ de vecteurs non
nul, et même un feuilletage de codimension 1.

Cette situation pourrait être considérée mauvaise du point de vue
d’obtenir des résultats visant à trouver des propriétés
topologiques d’un difféomorphisme qui préserve une structure
géométrique. En revanche, une belle remarque faite récemment par
Brin-Burago-Ivanov ([BBI, BI]) a redonné vie à l’espoir:

\begin{obsf*}[Brin-Burago-Ivanov]
Dans une variété de dimension 3, si ${\cal F}$ est une feuilletage
transversale à la direction instable $E^u$ d’un difféomorphisme
partiellement hyperbolique, alors ${\cal F}$ n’a pas de
composantes de Reeb. \finobs
\end{obsf*}

Les composantes de Reeb et leur relation avec les difféomorphismes
partiellement hyperboliques avait déjà été étudiée, par exemple
dans \cite{DPU} (Theorem H) et \cite{BonWilk} (Lemma 3.7) en un
contexte plus restrictif (la coherence dynamique et la
transitivite). La puissance de cette remarque est lie au fait
qu'elle d\'epend de la dynamique seulement dans le fait que le
feuilletage instable n'a pas des courbes ferm\'ees.

D’un côté, il est connu que plusieures variétés de dimension 3
n’admettent pas des feuilletage sans composantes de Reeb. D’un
autre côté, il existe beaucoup de théorie sur leur classification
(voir par exemple \cite{Plante, Roussarie}). Il semble donc
possible d'utiliser ces connaisances pour etudier les
diff\'eomorphismes partiellement hyperboliques.

Il apparaît une autre difficulté: on ne sait pas si tout
difféomorphisme partiellement hyperbolique possède une feuilletage
transversale à la direction instable. Dans cette thèse, nous
proposons la notion de \emph{presque-cohérence dynamique} et nous
prouvons qu’il s’agit d’une propriété ouverte et fermée dans
l’espace des difféomorphismes partiellement hyperboliques.

Un des théorèmes principaux (voir Chapitre
\ref{Capitulo-ParcialmenteHiperbolicos} et
\cite{Pot3dimPHisotopic}) garantie, dans certaines classes
d’isotopie de difféomorphismes de $\TT^3$, que cette condition est
suffisante pour que le fibré $E^{cs}$ soit tangent à un
feuilletage invariante:

\begin{teof*}
Soit $f$ un difféomorphisme partiellement hyperbolique et isotope
à un difféomorphisme Anosov linéaire. Si $f$ est
presque-dynamiquement cohérent, alors il est dynamiquement
cohérent.
\end{teof*}

La condition de \emph{cohérence dynamique} signifie justement que
le fibré tangent $E^{cs}$ soit tangent à una feuilletage
$f$-invariante.

On peut dire plus pour le cas \emph{partiellement hyperbolique
fort} (c’est à dire lorsque  il'y a un décomposition $TM=E^s\oplus
E^c\oplus E^u$ qui est $Df$-invariante avec des propriétés de
domination et telle que $E^s$ et $E^u$ sont uniformes).

En \cite{BI} les auteurs ont montré que tout difféomorphisme
partiellement hyperbolique fort est presque-dynamiquement cohérent
et ceci a été repris d’abord par \cite{BI} (suivant \cite{BBI}),
et après par \cite{Parwani}, pour donner des conditions
topologiques nécessaires pour l’hyperbolicité partielle forte.

Avec nos resultats, ca suggère que la suivante conjecture de
Pujals pourrait être abordable dans certains varietes:

\begin{conjf*}[Pujals \cite{BonWilk}]
Soient $M$ une variété de dimension 3 et $f$ un difféomorphisme
partiellement hyperboique fort et transitif de $M$. Une des trois
conditions suivantes doit se vérifier:
\begin{itemize}
\item[-]  $f$ est conjugué par feuilles à un difféomorphisme
Anosov de $\TT^3$. \item[-] $f$ est conjugué par feuilles à un
skew-produc sur un difféomorphisme Anosov en $\TT^2$ (la variété
est donc $\TT^3$ ou une nilvariété). \item[-] $f$ est conjugué par
feuilles au temps 1 d’un flot d'Anosov.
\end{itemize}
\end{conjf*}

Il y a eu des progrès dernièrement en ce qui concerne cette
conjecture. Tout d’abord, le travail de \cite{BonWilk} a fait des
progrès importants sans imposer d’hypothèse sur la topologie de la
variété, en supposant l’existence de courbes fermées tangentes à
la direction centrale. Plus tard, les travaux de Hammerlindl
\cite{Hammerlindl, HNil} donnent une preuve de la conjecture pour
le cas où la variété est $\TT^3$ ou une nilvariété mais pour une
notion plus restrictive de l’hyperbolicité partielle. Bien que de
nombreux exemples vérifient cette définition plus restrictive,
elle est dans un sens artificielle et ne s’adapte pas au résultat
de \cite{DPU}.

Bien sûr, pour obtenir une preuve générale, il faut d’abord
montrer que les difféomorphismes de ce genre sont dynamiquement
cohérents, puis que la définition de ``conjugué par feuilles”
(voir Section \ref{Section-FamiliasDePlacas})  a besoin de
l’existence de feuilletages invariantes. Le travail de Hammerlindl
s’appui sur des travaux précédents de Brin-Burago-Ivanov
(\cite{BBI2}) qui montrent que cette définition restrictive
d’’hyperbolicité partielle, dans le cas de $\TT^3$, implique
cohérence dynamique.

D’un autre côté, un exemple récent de Rodriguez Hertz-Rodriguez
Hertz-Ures (\cite{HHU}) présente un difféomorphisme partiellement
hyperbolique de $\TT^3$ qui n’admet pas de foliations invariantes.
Dans cette thèse nous complétons le paysage avec le suivant
résultat (\cite{Pot3dimPHisotopic} et Chapitre
\ref{Capitulo-ParcialmenteHiperbolicos}).

\begin{teof*}
Tout difféomorphisme partiellement hyperbolique fort de $\TT^3$
qui n’admet pas de tore periodique normallement attracteur ni de
tore periodique normallement répulseur est dynamiquement cohérent.
\end{teof*}

Ce résultat répond à une conjecture posée par Rodriguez
Hertz-Rodriguez Hertz-Ures dans le tore $\TT^3$ et sa preuve
permet de répondre à la conjecture de Pujals sur $\TT^3$ (voir
\cite{HP}).

Nous présentons à la fin du Chapitre
\ref{Capitulo-ParcialmenteHiperbolicos} des résultats en dimension
supérieures qui généralisent les résultats obtenus par Franks,
Newhouse et Manning pour des difféomorphismes Anosov dans le
contexte partiellement hyperbolique.

\subsection{Autres contributions}

Dans cette section nous décrivons des autres contributions de
cette thèse.

D’une part, dans la Section \ref{SubSection-LocalizacionDeClasses}
nous décrivons un mechanisme pour localiser des classes de
récurrence qui a été présenté dans \cite{PotWildMilnor} et qu'on
considéré important par soi-même. Ce mechanisme peut être appliqué
dans des contextes variés (dans cette thèse nous l’utilisons dans
la Section \ref{Section-DescDomT2} et aussi dans les sous-sections
\ref{SubSection-EjemploDA} et \ref{SubSection-EjemploPlykin}).

Dans la Section \ref{Section-EjemplosQuasiAtractores}, nous
présentons divers exemples de quasi-attracteurs et de
difféomorphismes robustement transitifs, dont quelques uns sont
obtenus à partir de modifications d’exemples connus, mais qui
peuvent, d’après nous, représenter une contribution a la
compréhension de ces phénomènes.

Le Chapitre \ref{Capitulo-FoliacionesEHipParcial} est dédié à
présenter les idées sur les feuilletages qui seront utilisées
après dans le Chapitre \ref{Capitulo-ParcialmenteHiperbolicos},
nous obtenons des résultats qui peuvent avoir un intérêt
indépendent. En particulier, un résultat quantitatif sur
l’existence d’une structure de produit global pour les foliations,
présenté dans la Section \ref{Section-GlobalProductStructure}. De
même, dans la Section \ref{Section-DescDomT2}, nous donnons une
classification de la dynamique des difféomorphismes globalement
partiellement hyperbolique en $\TT^2$ qui, d’une certaine façon,
montre ce qu’on faira dans le Chapitre
\ref{Capitulo-ParcialmenteHiperbolicos} dans un contexte plus
simple.

Nous ajoutons aussi quatre annexes où nous présentons des
résultats qui se détachent du corps central de la thèse. Nous
soulignons l’annexe \ref{Apendice-PseudoRotaciones}, basé sur
\cite{PotRecurrence}, où nous prouvons un résultat sur les
homéomorphismes du tore qui possèdent un seul vecteur de
rotations. Aussi, dans l’annexe \ref{Apendice-BCGP}, basé sur
\cite{BCGP} nous présentons un travail en collaboration avec
Bonatti, Crovisier et Gourmelon, où nous étudions les
biffurcations de classes de récurrence robustement isolées et nous
donnons des exemples non robustement transitifs, ce qui répond à
une question posée en \cite{BC}.

\section{Organization of this thesis}

This thesis is organized as follows:

\begin{itemize}
\item[-] In Chapter \ref{Capitulo-Preliminares} we introduce
definitions and known results about differentiable dynamics which
will be used along the thesis. This chapter also presents in a
systematic way the context in which we will work.  \item[-] In
Chapter \ref{Capitulo-Semiconjugaciones} we present background
material on semiconjugacies and we also present a mechanism for
localization of chain-recurrence classes (see Section
\ref{SubSection-LocalizacionDeClasses}). \item[-] In Chapter
\ref{Capitulo-Atractores} we study attractors and quasi-attractors
in $C^1$-generic dynamics. \item[-] In Chapter
\ref{Capitulo-FoliacionesEHipParcial} we give an introduction to
the known results in foliations mainly focused in codimension one
foliations and particularly in foliations of $3$-manifolds. We
prove some new results which we will use later in Chapter
\ref{Capitulo-ParcialmenteHiperbolicos}. Also, this chapter
contains an appendix which shows similar results as the main
results of this thesis in the context of surfaces. \item[-] In
Chapter \ref{Capitulo-ParcialmenteHiperbolicos} we study global
partial hyperbolicity. Most of the chapter is devoted to the study
of partially hyperbolic diffeomorphisms of $\TT^3$ and in the last
section some results in higher dimensions are given. \item[-]
Appendix \ref{Appendix-PerturbacionCociclos} presents some
techniques on perturbations of cocycles over periodic orbits.
\item[-] Appendix \ref{Apendix-PlaneDec} gives an example of a
decomposition of the plane satisfying some pathological
properties. \item[-] Appendix \ref{Apendice-PseudoRotaciones} is
devoted to the study of homeomorphisms of $\TT^2$ with a unique
rotation vector. There we present the results of
\cite{PotRecurrence} and we also give a quite straightforward
extension of the results there to certain homeomorphisms homotopic
to dehn-twists. \item[-] Appendix \ref{Apendice-BCGP} presents the
results of \cite{BCGP}.
\end{itemize}


\section{Reading paths}

Being quite long, it seems reasonable to indicate at this stage
how to get to certain results without having to read the whole
thesis.

First of all, it must be said that Chapter
\ref{Capitulo-Preliminares}, concerning preliminaries, need not be
read for those who are acquainted with the subject. In particular,
those who are familiar with one or more of the excellent surveys
\cite{BDV,Crov-Hab,PujSamHandbook}.

Similar comments go for Chapter
\ref{Capitulo-FoliacionesEHipParcial}, in particular the first
part covers well known results on the theory of foliations which
can be found in \cite{CandelConlon, Calegari} among other nice
books. The final part though may be of interest specially for
those which are not specialist on the theory of foliations such as
the author.

If the reader is interested in the part of this thesis concerned
with attractors for $C^1$-generic dynamics, then, the suggested
path (which can be coupled with the suggestions in the previous
paragraphs) is first reading Chapters \ref{Capitulo-Preliminares}
and Chapter \ref{Capitulo-Semiconjugaciones} and then Chapter
\ref{Capitulo-Atractores}.

If on the other hand, the reader is interested in the part about
global partially hyperbolic dynamics, then some parts of Chapter
\ref{Capitulo-Preliminares} can be skipped, in particular it is
enough with reading Sections \ref{Section-EstructurasInvariantes}
and \ref{Section-HIRSHPUGHSHUB}. Then, Chapter
\ref{Capitulo-FoliacionesEHipParcial} is fundamental in those
results, but the reader which is familiar with the theory of
foliations may skip it in a first read. Finally, the results about
global partial hyperbolicity are contained in Chapter
\ref{Capitulo-ParcialmenteHiperbolicos}. If the reader is
interested in Section \ref{Section-DescDomT2}, then Chapter
\ref{Capitulo-Semiconjugaciones} is suggested.


   \chapter{Preliminaries}\label{Capitulo-Preliminares}


\section{Recurrence and orbit perturbation tools}\label{Seccion-Recurrence-And-Perturbation}

\subsection{Some important dynamically defined sets and transitivity}\label{SubSection-ConjuntosDinamicamenteDefinidos}

\subsubsection{} Let $f: X \to X$ a homeomorphism of $X$ a compact
metric space. For a point $x \in X$ we define the following sets:

\bi \item[-] The \emph{orbit} of $x$ is the set $\cO(x) = \{
f^n(x) \ : \ n \in \ZZ \}$. We can also define the future (resp.
past) orbit of $x$ as $\cO^+(x)= \{ f^n(x) \ : \ n \geq 0\}$
(resp. $\cO^-(x)= \{f^n(x) \ : \ n \leq 0 \}$). \item[-] The
\emph{omega-limit set} (resp.\emph{alpha-limit set}) is the set
$\omega(x,f)= \{ y \in X \ : \ \exists n_j \to + \infty \
\text{such that } f^{n_j}(x) \to y\}$ (resp. $\alpha(x,f)=
\omega(x,f^{-1})$). In general, when $f$ is understood, we shall
omit it from the notation. \ei

We can divide the points depending on how their orbit and the
nearby orbits behave. We define the following sets:

\bi \item[-] $\Fix(f)= \{ x \in X \ : \ f(x)=x \}$ is the set of
\emph{fixed points}. \item[-] $\Per(f)= \{x \in X \ : \ \#\cO(x) <
\infty \}$ is the set of \emph{periodic points}. The period of a
periodic point $x$ is $\#\cO(x)$ which we denote as
$\pi(x)=\#\cO(x)$. \item[-] We say that a point $x$ is
\emph{recurrent} if $x \in \omega(x)\cup \alpha(x)$. \item[-]
$\Lim(f) = \overline{ \bigcup_x \omega(x) \cup \alpha(x)}$ is the
\emph{limit set} of $f$. \item[-] $\Omega(f)= \{ x \in X \ : \
\forall \eps>0 \ , \ \exists n>0 \ ; \ f^n(B_\eps(x)) \cap
B_\eps(x) \neq \emptyset \}$ it the \emph{nonwandering set} of
$f$. \ei

We refer the reader to \cite{KH,Robinson,Shub} for examples
showing the strict inclusions in the following chain of closed
sets which is easy to verify:

$$ \Fix(f) \en \Per(f) \en \Lim(f) \en \Omega(f) $$

In section \ref{SubSection-ChainRecurrence} we shall explore
another type of recurrence which will play a central role in this
text.

\smallskip

 For a point $x\in X$ we will define the following sets (as
before the reference to the homeomorphism $f$ may be omitted when
it is obvious from the context).

\bi \item[-] The \emph{stable set} of $x$ is $W^s(x, f)= \{y \in X
\ : \ d(f^n(x),f^n(y)) \to 0 \text{ as } \ n \to +\infty \}$.
\item[-] The \emph{unstable set} of $x$ is $W^u(x,f)=
W^s(x,f^{-1})$. \ei

It is clear that $f(W^\sigma(x)) = W^\sigma(f(x))$ for
$\sigma=s,u$. The study of these sets and how are they related is
one of the main challenges one faces when trying to understand
dynamical systems.

Sometimes, it is useful to consider instead the following sets
which in some cases (for sufficiently small $\eps$) are related
with the stable and unstable sets:

\bi \item[-] $S_\eps (x,f)= \{ y \in X \ : \ d(f^n(x),f^n(y)) <
\eps \ , \ \forall n \geq 0 \}$. \item[-] $U_\eps(x,f) =
S_\eps(x,f^{-1})$. \ei

Notice the  following two properties which are essentially the
reason for defining these sets:

$$ W^s(x) \en \bigcup_{n \geq 0} f^{-n}(S_\eps(f^n(x))) $$

$$ f(S_\eps(x)) \en S_\eps(f(x)) $$

Similar properties hold for $U_\eps(x)$. It is not hard to make
examples where the inclusions are strict. However, the first
property is an equality for certain special maps (expansive, or
hyperbolic) which will be of importance in this text.

For a set $K \en X$ we define $W^\sigma (K) = \bigcup_{x\in K}
W^\sigma(x)$ with $\sigma=s,u$. Notice that the set is (a priori)
smaller than the set of points whose omega-limit is contained in
$K$.


\smallskip

We say that $f$ is \emph{transitive} if there exists $x\in X$ such
that $\cO(x)$ is dense in $X$. Sometimes, when $f$ is understood
(for example, when $X$ is a compact invariant subset of a
homeomorphism of a larger set), we say that $X$ is
\emph{transitive}.

It is an easy exercise to show the following equivalences (see for
example \cite{KH} Lemma 1.4.2 and its corollaries):

\begin{prop}\label{Proposition-EquivalenciasTransitividad}
The homeomorphism $f: X \to X$ is transitive if and only if for
every $U,V$ open sets there exists $n \in \ZZ$ such that
$f^n(U)\cap V \neq\emptyset$, if and only if there is a residual
subset of points whose orbit is dense.
\end{prop}

We say that $f$ (or as above that $X$) is \emph{minimal} if every
orbit is dense.

Given an open set $U$, we define the following compact
$f$-invariant set $\Lambda=\bigcap_{n\in \ZZ} f^n(\overline U)$
which we call the \emph{maximal invariant set} in $\overline{U}$.

Many of the dynamical properties one obtains are invariant under
what is called \emph{conjugacy}. We say that two homeomorphisms
$f:X\to X$ and $g: Y \to Y$ are (topologically) \emph{conjugated}
if there exists a homeomorphism $h: X \to Y$ such that:

$$ h \circ f = g \circ h $$

When $K$ is an $f$-invariant set and $K'$ a $g$-invariant set we
say that $f$ and $g$ are \emph{locally conjugated at} $K$ if there
exists a neighborhood $U$ of $K$, a neighborhood $V$ of $K'$ and a
homeomorphism $h: U \to V$ such that if a point $x\in U \cap
f^{-1}(U)$ then:

$$ h \circ f(x) = g \circ h(x) $$

\subsection{Hyperbolic periodic points}\label{SubSection-PeriodicosHiperbolicos}

From now on, $f: M^d \to M^d$ will denote a $C^1$ diffeomorphism.

Given $p \in \Per(f)$ we have the following linear map:

$$ D_p f^{\pi(p)}: T_p M \to T_p M$$

We say that $p \in \Per(f)$ is a \emph{hyperbolic periodic point}
if $D_pf^{\pi(p)}$ has no eigenvalues of modulus $1$. We denote
the set of hyperbolic periodic points as $\Per_H(f)$.

It is a direct application of standard linear algebra  to show
that in fact the set $\Per_H(f)$ is $f$-invariant. This implies
that we can also talk about \emph{hyperbolic periodic orbits}.

For a hyperbolic periodic point $p$ we have that $T_pM = E^s(p)
\oplus E^u(p)$ where $E^s(p)$ (resp. $E^u(p)$) corresponds to the
eigenspace of $D_pf^{\pi(p)}$ associated to the eigenvalues of
modulus smaller than $1$ (resp. larger than $1$). We have that
$D_pf (E^\sigma(p)) = E^\sigma(f(p))$ with $\sigma=s,u$.

We define the (\emph{stable}) \emph{index}\footnote{We warn the
reader that some authors define the index of a periodic point as
the unstable dimension.} of a periodic point $p$ as $\dim E^s(p)$.
This also leads to calling \emph{stable eigenvalues} (resp.
\emph{unstable eigenvalues}) to those which are of modulus smaller
(resp. larger) than $1$. We denote the set of index $i$ periodic
points as $\Per_i(f)$.

The importance of hyperbolic periodic points is related to the
fact that their local dynamics is very well understood and it is
persistent under $C^1$-perturbations (for a proof see for example
\cite{KH} chapter 6):

\begin{teo}\label{Teorema-PeriodicosHiperbolicosVariedadEstable}
Let $p$ be a hyperbolic periodic point of a $C^1$-diffeomorphism
$f$, then: \bi \item[(i)] $f^{\pi(p)}$ is locally conjugated to
$D_pf^{\pi(p)}$ at $p$. In particular, there are no periodic
points of period smaller or equal to $\pi(p)$ inside a
neighborhood $U$ of $p$. \item[(ii)] There exists a
$C^1$-neighborhood $\cU$ of $f$ such that for every $g\in \cU$
there is a unique periodic point $p_g$ of period $\pi(p)$ of $g$
inside $U$ which is also hyperbolic. We say that $p_g$ is the
\emph{continuation} of $p$ for $g$. \item[(iii)] There exists
$\eps>0$ such that $S_\eps$ is an embedded $C^1$ manifold tangent
to $E^s(p)$ at $p$ and in particular, one has $S_\eps \en W^s(p)$.
\ei
\end{teo}

As a consequence we have that the set of hyperbolic periodic
points of period smaller than $n$ is finite. We have that the sets
$W^\sigma(p)$ are (injectively) immersed $C^1$-submanifolds of $M$
diffeomorphic to $\RR^{\dim E^\sigma}$. Moreover, one can define
$W^\sigma(\cO(p))$ which will also be an injectively immersed
$C^1$-submanifold with the same number of connected components as
the period of $p$.

When the stable index $s$ of a hyperbolic periodic point is the
same as the dimension of the ambient manifold (resp. $s=0$), we
shall say that it is a \emph{periodic sink} (resp. \emph{periodic
source}). In any other case we shall say that it is a
\emph{periodic saddle} of index $s$.

One of the first perturbation results in dynamics was given by
Kupka and Smale independently showing:

\begin{teo}[\cite{KupkaKupkaSmale,SmaleKupkaSmale}]\label{Teorema-KupkaSmale}
For every $r\geq 1$, there exists a
residual subset $\cG_{KS} \en \Diff^r(M)$ of diffeomorphisms such
that if $f\in \cG_{KS}$: \bi \item[-] All periodic points are
hyperbolic (i.e. $\Per(f)=\Per_H(f)$). \item[-] Given $p,q \in
\Per(f)$ we have that $W^s(p)$ and $W^u(q)$ intersect
transversally (recall that this allows the manifolds not to
intersect at all). \ei
\end{teo}

Transversal intersections between stable and unstable manifolds
yield information on the iterates of those manifolds. A quite
useful tool to treat those intersections is given by the
celebrated $\lambda$-Lemma (or Inclination Lemma) of Palis (see
\cite{PalisLambdaLemma}) which we state as follows (see also
\cite{KH} Proposition 6.2.23). The statement we present is for
fixed points, but considering an iterate one can of course treat
also periodic points as well:

\begin{teo}[$\lambda$-Lemma]\label{Teorema-LambdaLemmaOriginal}
Let $p$ be a hyperbolic fixed point of $f$ a $C^1$-diffeomorphism
of a manifold $M$ and let $D$ be a $C^1$-embedded disk which
intersects $W^s(p)$ transversally. Then, given a compact
submanifold $B$ of $W^u(p)$ and $\eps>0$ there exists $n_0>0$ such
that for every $n>n_0$ there is a compact submanifold $D_n \en D$
such that $f^n(D_n)$ is at $C^1$-distance smaller than $\eps$ of
$B$.
\end{teo}

See also Lemma \ref{LemaLambdaLema}.


\subsection{Homoclinic classes}\label{Subsection-ClasesHomoclinicas}

Given two hyperbolic periodic points $p,q$ we say that they are
\emph{homoclinically related} if $W^s(p) \trans W^u(q) \neq
\emptyset$ and $W^s(q) \trans W^u(p)\neq \emptyset$.

Given a periodic obit $\cO$ we define its \emph{homoclinic class}
$H(\cO)$ as the closure of the set of periodic points
homoclinically related to some point in the orbit $\cO$.

Notice that by the definition of being homoclinically related,
necessarily one has that if two periodic points are homoclinically
related, then they have the same stable index. However, this does
not exclude the possibility of having periodic points of different
index (and even non-hyperbolic periodic points) inside $H(\cO)$
and this will be ``usually'' the case outside the ``hyperbolic
world''.

Homoclinic classes were introduced by Newhouse as an attempt to
generalize for arbitrary diffeomorphisms the basic pieces
previously defined by Smale (\cite{SmaleBulletin}) for Axiom A
diffeomorphisms. The first and probably the main example of
non-trivial homoclinic class is given by the famous
\emph{horseshoe} of Smale (see \cite{SmaleBulletin} or \cite{KH}
2.5). We have the following properties:

\begin{prop}[\cite{Newhouse-homoclinic}]\label{Proposition-PropiedadesClasesHomoclinicas}
For every hyperbolic periodic orbit $\cO$ of a
$C^1$-diffeomorphism $f$ the homoclinic class $H(\cO)$ is a
transitive $f$-invariant set. Moreover, we have that:

$$ H(\cO) = \overline{ W^s(\cO) \trans W^u(\cO)} .$$
\end{prop}

This proposition essentially follows as an application of the
$\lambda$-Lemma (see also \cite{SmaleBulletin}).

For a periodic point $p$ we denote as $H(p)=H(\cO(p))$. Given a
hyperbolic periodic point $p$ of a $C^1$-diffeomorphism, we have
by Theorem \ref{Teorema-PeriodicosHiperbolicosVariedadEstable}
that there exists a continuation of $p$ as well as $W^s(p)$ and
$W^u(p)$ for close diffeomorphisms. We will sometimes make
explicit reference to the diffeomorphism and use the notation
$H(p,f)$.

We have the following fact which follows from the continuous
variation of stable and unstable manifolds with the
diffeomorphism:

\begin{prop}\label{Proposition-ContinuacionClaseHomoclinica}
Given a hyperbolic periodic point $p$ of a $C^1$-diffeomorphism
$f$ and $U$ an open neighborhood of $M$ such that $H(p)\cap U \neq
\emptyset$, there exists $\cU$ a $C^1$-neighborhood such that if
$H(p_g)$ is the homoclinic class for $g\in \cU$ of the
continuation $p_g$ of $p$ we have that $H(p_g)\cap U \neq
\emptyset$.
\end{prop}

\begin{obs}[Semicontinuity]\label{Remark-Semicontinuidad}
The last statement of the proposition can be stated by saying that
homoclinic classes vary \emph{semicontinuously} with respect to
the Hausdorff topology. This means that they cannot \emph{implode}
(if $f_n \to f$ have hyperbolic periodic points $p_n$ which are
the continuation of $p \in \Per_H(f)$ for $f_n$, we have that if
$H(p,f)$ intersects a given open set $U$, then $H(p_n,f_n)$ also
intersects $U$ for large enough $n$). There exists the possibility
that the homoclinic class \emph{explodes} by small perturbations
(see for example \cite{PalisOmega,DiazSantoro}). However, a
classic result in point set topology guaranties that when a map is
semicontinuous then it must be continuous in a residual subset
(see Proposition 3.9 of \cite{Crov-Hab} for a precise statement).
Other compact sets related to the dynamics which have
semicontinuous variation are $\overline{\Per_H(f)}$,
$\overline{\Per_i(f)}$ or for a given hyperbolic periodic point
$p$ of $f$ the set $\overline{W^\sigma(p_g)}$ ($\sigma=s,u$) for
$g$ in a neighborhood of $f$. All these sets cannot implode but
may explode in some situations, from the mentioned result on point
set topology, there is a residual subset $\cG$ of $\Diff^1(M)$
where all these sets vary continuously. In the next subsection we
shall see a set which also varies semicontinuously but in ``the
other sense'', meaning that it can implode but not explode.\finobs
\end{obs}

As it was already mentioned, periodic points in a homoclinic class
may not have the same stable index (even if to be homoclinically
related they must have the same index). The existence of periodic
points of different index in a homoclinic class is one of the main
obstructions for hyperbolicity. Given a homoclinic class $H$ we
say that its \emph{minimal index} (resp. \emph{maximal index}) is
the smallest (resp. largest) stable index of periodic points in
$H$.

\subsection{Chain recurrence and filtrations}\label{SubSection-ChainRecurrence}

In this section we shall review yet another recurrence property,
namely, \emph{chain-recurrence}. From the point of view of
recurrence, it can be regarded as the ``weakest'' form of
recurrence for the dynamics. Indeed, it is so weak that it was
neglected for much time in differentiable dynamics since its
classes seem to have really poor dynamical indecomposability (see
for example Appendix \ref{Apendice-PseudoRotaciones} and
\cite{PotRecurrence}).

On the other hand, it is by far the best notion when one wishes to
decompose the dynamics in pieces, and this is why after being
shown to be quite similar to the rest of the notions for
$C^1$-generic dynamics (in \cite{BC}) it became ``the'' notion of
recurrence used in $C^1$-differentiable dynamics.

We derive the reader to \cite{Crov-Hab} for a more comprehensive
introduction to these concepts (see also \cite{BDV} chapter 10 for
another introduction to these topics which is less up to date but
still a good introduction).

\begin{defi}[Pseudo-orbits]
Given a homeomorphism $f: X \to X$ and points $x, y \in X$ we say
that there exists an $\eps$-\emph{pseudo orbit} from $x$ to $y$
and we denote it as $x \dashv_\eps y$ iff there exists points
$z_0=x, \ldots, z_k =y$ such that $k \geq 1$ and
$$ d(f(z_i),z_{i+1}) \leq \eps  \quad \quad 0\leq i \leq k-1$$
We use the notation $x\dashv y$ to express that for every $\eps>0$
we have that $x\dashv_\eps y$. We also use $x \dashrv y$ to mean
$x\dashv y$ and $y\dashv x$. \finobs
\end{defi}

We define the \emph{chain-recurrent set} of $f: X \to X$ as

$$ \cR(f) = \{x \in X \ : \ x\dashrv x \}  $$

It is easy to show that inside $\cR(f)$ the relation $\dashrv$ is
an equivalence relation so we can decompose $\cR(f)$ in the
equivalence classes which we shall call \emph{chain-recurrence
classes}. For a point $x \in \cR(f)$ we will denote as $\cC(x)$
its chain-recurrence class. Both $\cR(f)$ and the chain recurrence
classes can be easily seen to be closed (and thus compact).

One can regard chain recurrence classes as maximal \emph{chain
transitive sets}. We say that an invariant set $K \en X$ is
chain-transitive if for every $x,y \in K$ we have that $x\dashv
y$. We say that a homeomorphism $f: X \to X$ is \emph{chain
recurrent} if $X$ is a chain transitive set for $f$.

For a chain-transitive set $K$ we define its \emph{chain stable
set} (resp. \emph{chain unstable set}) as $pW^s(K) = \{ y \in X \
: \ \exists z \in K \ \text{such that} \ y\dashv z \}$ (resp. as
$pW^u(K)$, the chain stable set for $f^{-1}$).

\begin{obs}\label{Remark-ClaseHomoclinicaYdeRecurrencia}
Given a hyperbolic periodic point $p$ of a $C^1$-diffeomorphism,
one has that the homoclinic class of $p$ is a chain-transitive
set. In particular, it is always contained in the chain-recurrence
class of $p$. In general, one can have that the inclusion is
strict (see for example \cite{DiazSantoro}). Notice also that if
$p$ is a hyperbolic sink (or source) we have that $\cC(p)=\cO(p)$.
Indeed, since $p$ admits neighborhoods whose closure is sent to
its interior by $f$ (or $f^{-1}$), this prevents small
pseudo-orbits to leave (or enter) any small neighborhood of $p$.
\finobs
\end{obs}

An essential tool for decomposing chain-recurrence classes whose
existence is the content of Conley's theory (see \cite{Conley} and
\cite{Robinson} chapter 9.1) are Lyapunov functions. We remark
that the definition we use of Lyapunov function is not the
standard one in the literature, we have adapted our definition in
order to have the properties of the function given by Conley's
Theorem.

\begin{defi}[Lyapunov Functions]\label{Defi-LyapunovFunc}
Given a homeomorphism $f: X \to X$ we say that $\varphi: X\to
[0,1]$ is a \emph{Lyapunov function} if the following conditions
are satisfied: \bi \item[-] For every $x \in X$ we have that
$\varphi(f(x)) \leq \varphi(x)$ and $\varphi(x)=\varphi(f(x))$ if
and only if $x\in \cR(f)$. \item[-] Given $x,y \in \cR(f)$ then
$\varphi(x)=\varphi(y)$ if and only if $\cC(x)=\cC(y)$. \item[-]
The image of $\cR(f)$ by $\varphi$ has empty interior. \ei \finobs
\end{defi}

It is remarkable that these functions always exist (see
\cite{Robinson} chapter 9.1 for a simple proof of the following
theorem also sometimes called \emph{Fundamental theorem of
dynamical systems}).

\begin{teo}[Conley \cite{Conley}]\label{Teorema-Conley}
For any homeomorphism $f:X\to X$ of a compact metric space $X$
there exists a Lyapunov function $\varphi: X \to \RR$.
\end{teo}

\begin{obs}[Filtrations]\label{Remark-Filtrations}
Lyapunov functions allow to create \emph{filtrations} separating
chain recurrence classes. Indeed, consider $\cC_1$ and $\cC_2$ two
distinct chain recurrence classes, and a Lyapunov function
$\varphi$.

From the definition, we have that without loss of generality, we
can assume $\varphi(\cC_1) < \varphi(\cC_2)$. Since the image by
$\varphi$ of the chain-recurrent set has empty interior, there
exists $a \in [0,1] \setminus \varphi(\cR(f))$ such that
$\varphi(\cC_1) < a < \varphi(\cC_2)$.

Let $U= \varphi^{-1}((-\infty,a))$ an open set. Since every point
such that $\varphi(x)=a$ is not chain-recurrent, we obtain by the
definition of Lyapunov function that $f(\overline{U}) \en U$ and
that $\cC_1 \en U$ and $\cC_2 \en \overline{U}^c$.

Moreover, every chain-recurrence class $\cC$ admits a basis of
neighborhoods $U_n$ such that if $\Lambda_n$ is the maximal
invariant subset of $U_n$, then $\cC= \bigcap_n \Lambda_n$.
Moreover, it verifies that if $\cC_0$ is a chain recurrence class
which intersects $U_n$ then $\cC_0$ is contained in $U_n$ (one can
consider $\varphi^{-1}((a- \eps_n, a+\eps_n)$ with $\eps_n \to
0$). The sets $U_n$ are sometimes called \emph{filtrating
neighborhoods} for $\cC$.

 \finobs
\end{obs}

These filtrating neighborhoods (which persist under $C^0$-small
perturbations) allow one to show that the mapping $f\mapsto
\cR(f)$ is semicontinuous in the sense that it cannot ``explode''
(see Remark \ref{Remark-Semicontinuidad}), so it will vary
continuously in a residual subset of $\Diff^1(M)$ with respect to
the Hausdorff topology on compact sets. See the example of
Appendix \ref{Apendice-BCGP}.

Pseudo-orbits can be thought of real-orbits of $C^0$-perturbations
of the initial system. We have:

\begin{prop}\label{Proposicion-VariacionClasesDeRecurrencia}
Let $f_n:X \to X$ be a sequence of homeomorphisms such that $f_n
\to f$ in $C^0$-topology and let $\Lambda_n$ be chain-transitive
sets for $f_n$. Then, in the Hausdorff topology we have that
$\Lambda = \limsup \Lambda_n$ is a chain-transitive set.
\end{prop}

\dem Consider $x,y \in \Lambda$ and $\eps>0$. We can consider $n$
large enough so that

\bi \item[-] $d_{C^0}(f_n,f) < \frac{\eps}{2}$. \item[-]
$d_H(\Lambda_n,\Lambda) < \frac{\eps}{2}.$ \ei

Since a $\frac{\eps}{2}$-pseudo-orbit for $f_n$ will be an
$\eps$-pseudo-orbit of $f$ and since there are points in
$\Lambda_n$ which are $\frac{\eps}{2}$-close to $x$ and $y$ we
conclude.

\lqqd

\begin{obs}[Trapping regions]\label{Remark-RegionesAtractorasConsequenciaConley}
Conley's Theorem implies in particular that a homeomorphism
$f:X\to X$ is chain-recurrent if and only if there is no proper
(i.e. strictly contained) open set $U \en X$ such that
$f(\overline{U})\en U$. \finobs
\end{obs}

In order to understand the asymptotic behavior of orbits one must
then comprehend the dynamics inside chain recurrence classes as
well as how the classes are related to each other.

An important concept is then that of \emph{isolation}. We say that
a chain-recurrence class $\cC$ is \emph{isolated} iff there exists
a neighborhood $U$ of $\cC$ such that $U \cap \cR(f) = \cC$. This
is equivalent to $\varphi(\cC)$  being isolated in
$\varphi(\cR(f))$ for some Lyapunov function (as defined in
Definition \ref{Defi-LyapunovFunc}) of $f$. Sometimes, non
isolated classes will be referred to as \emph{wild chain
recurrence classes}.

\begin{obs}\label{Remark-ClasesAisladasSonMaximalInvariante}
In particular, one has that a chain-recurrence class is isolated
if and only if it is the maximal invariant set in a neighborhood
of itself. \finobs
\end{obs}

\subsection{Attracting sets}\label{SubSection-AttractingSets}

It seems natural that the goal of understanding the whole orbit
structure for general homeomorphisms should be quite difficult.
This is why, in general, we content ourselves by trying to
understand ``almost every'' orbit of ``almost every'' system. This
informal statement has various ways to be understood, in
particular, it is well known that many different formalizations of
``almost every'' can be quite different (the paradigma of this is
seen in the case of irrational numbers, where Diophantine ones
have total Lebesgue measure while Liouville ones form a disjoint
residual subset of $\RR$).

However, it seems natural in view of the Lyapunov functions to
study certain special chain-recurrence classes which are called
\emph{quasi-attractors}. In this section we shall define plenty
types of attractors which in a certain sense will be the
chain-recurrence classes to which we shall pay more attention in
view of the discussion above.

Given an open set $U \en X$ such that $f(\overline{U}) \en U$ we
can consider the set $\Lambda = \bigcap_{n>0} f^n(U)$ which is
compact and invariant (it is the maximal invariant subset in $U$).
We call $\Lambda$ a \emph{topological attractor}\footnote{This is
the usual definition in the literature related to this subject
(see \cite{BDV} Chapter 10). It seems that it could be better to
call this sets \emph{attracting sets}, since the word attractor
may be sometimes misleading. However, we have chosen to keep this
nomenclature.}.

\begin{prop}\label{ProposicionAtractoresSaturadosPorInestables}
Let $U$ be an open set such that $f(\overline U)\en U$ and let
$\Lambda$ be its maximal invariant set. If $y \in X$ is a point
such that for every $\eps>0$ there exists $z\in \Lambda$ such that
$z \dashv_\eps y$, then, $y\in \Lambda$. In particular, for every
$z\in \Lambda$ we have that $\overline{W^u(z)}\en \Lambda$.
\end{prop}

\dem Let $y \in X$ be such that for every $\eps>0$ there is some
$z\in \Lambda$ such that $z \dashv_\eps y$.

Assume by contradiction that $y \notin \Lambda$. Since $\Lambda$
is invariant, we can assume (by iterating backwards) that $y\notin
\overline U$.

Let $\delta>0$ be such that $d(\partial U, f(\overline{U}))>
\delta$. We will show that there cannot be a $\delta$-pseudo orbit
from $\Lambda$ to $y$.

Indeed, given a point $x \in U$ we have that $f(x)$ is in
$f(\overline U)$ which implies by induction that a
$\delta$-pseudo-orbit starting at $U$ must remain in $U$. This is
a contradiction and proves the proposition.

\lqqd

The problem with topological attractors is they are not
indecomposable in the sense that the dynamics inside $\Lambda$ may
not even be chain-recurrent (and in fact they can admit
topological attractors contained inside themselves). On the other
hand, they have the virtue of always existing (for example, the
whole space is always a topological attractor, and by Remark
\ref{Remark-RegionesAtractorasConsequenciaConley} there always
exist proper topological attractors when the homeomorphism is not
chain-recurrent). To obtain in a sense better suited definitions
we present now the definition of \emph{attractors} and
\emph{quasi-attractors} which will appear throughout this text as
one of the main objects of study.

\begin{defi}[Attractor]\label{Definicion-Attractor}
We say that a compact invariant set $\Lambda$ is an
\emph{attractor} if it is a topological attractor and it is
chain-recurrent. An attractor for $f^{-1}$ is called a
\emph{repeller}. \finobs
\end{defi}

\begin{obs}\label{Remark-Attractor} We remark that it is usual in the
literature also to define attractor by asking the stronger
indecomposability hypothesis of being transitive, we use this
definition since our context is better suited with the use of
chain-recurrence. It is easy to see that if $\Lambda$ is an
attractor, then it is an isolated chain-recurrence class. \finobs
\end{obs}

In general, a homeomorphism may not have any attractors, however,
it will always have what we call \emph{quasi-attractors}.

\begin{defi}\label{Definicion-QuasiAtractor} A compact invariant set
$\cQ$ is a \emph{quasi-attractor} for a homeomorphism $f: X\to X$
if and only if it is a chain-recurrence class and there exists a
nested sequence of open neighborhoods $\{U_n\}$ of $\cQ$ such
that: \bi \item[-] $\cQ = \bigcap_n U_n$, and \item[-]
$f(\overline{U_n}) \en U_n$. \ei A quasi-attractor for $f^{-1}$ is
called a \emph{quasi-repeller}. \finobs
\end{defi}

\begin{obs}\label{Remark-QuasiAttractor}
\bi \item[-]If $\varphi: X \to \RR$ is a Lyapunov function for a
homeomorphism $f:X \to X$, it is clear that $\cQ$, the
chain-recurrence class for which the value of $\varphi$ is the
minimum must be a Lyapunov stable set. Recall that a compact set
$\Lambda$ is \emph{Lyapunov stable} for $f$ if for every
neighborhood $U$ of $\Lambda$ there exists a neighborhood $V$ of
$\Lambda$ such that $f^n(V) \en U$ for every $n\geq 0$. \item[-]It
is almost direct from the definition that a quasi-attractor must
always be a Lyapunov stable set. \item[-] Moreover, although we
shall not use it, it is not hard to see that given a
quasi-attractor of $f$ one can always construct a Lyapunov
function attaining a minimum in the given quasi-attractor. This
follows from the proof of Conley's Theorem (see \cite{Robinson}
Chapter 9.1). \item[-] A quasi-attractor is a topological
attractor if and only if it is an attractor. A quasi-attractor is
an attractor if and only if it is isolated (as a chain-recurrence
class). This implies that if a homeomorphism has no attractors,
then it must have infinitely many distinct chain-recurrence
classes.  \ei \finobs
\end{obs}

We say that a chain-recurrence class $\cC$ is a \emph{bi-Lyapunov
stable class} iff it is Lyapunov stable for both $f$ and $f^{-1}$.

We now state a corollary from Proposition
\ref{ProposicionAtractoresSaturadosPorInestables}:

\begin{cor}\label{Corolario-QuasiAtractoresSaturadosPorInestables}
Let $\cQ$ be a quasi-attractor for a homeomorphism $f:X\to X$. If
$y \in X$ is a point such that for every $\eps>0$ there exists
$z\in \Lambda$ such that $z \dashv_\eps y$, then, $y\in \cQ$. In
particular, for every $z\in \cQ$ we have that
$\overline{W^u(z)}\en \cQ$.
\end{cor}

\dem In Proposition
\ref{ProposicionAtractoresSaturadosPorInestables} it is proved
that if $y$ is as in the statement it must belong to $U_n$ for all
$U_n$ in the definition of quasi-attractor. This concludes.
\lqqd

To finish this section we will define two further notions of
attracting sets which will also appear later in the text.

We say that a quasi-attractor $\cQ$ is an \emph{essential
attractor} (as defined in \cite{BLY}) if it has a neighborhood $U$
which does not intersect any other quasi-attractors. The
importance of these classes is given by a conjecture by Hurley
\cite{Hurley} (known in certain topologies, see Theorem
\ref{Teorema-BonattiCrovisier}): \emph{For typical dynamics,
typical points converge to quasi-attractors}.

Sometimes, the invariant sets which attract important parts of the
dynamics need not be chain-recurrence classes, and even not
Lyapunov stable. Another important kind of ``attracting sets'' are
\emph{Milnor attractors} (see \cite{MilnorAttractor}). To define
them we first define the \emph{topological basin} of a compact
invariant set $K$ as:

$$  \Bas(K) =\{ y\in X \ : \ \omega(y) \en K \}  $$

We say that a compact $f$-invariant set $K$ is a \emph{Milnor
attractor} if $\Leb(\Bas(K)) >0$ and for every $K' \en K$ compact,
invariant different from $K$ one has that $\Leb(\Bas(K')) <
\Leb(\Bas(K))$.

The definition seems a little stronger than demanding that the
basin has positive Lebesgue measure, but a simple Zorn's Lemma
argument gives:

\begin{lema}[Lemma 1 of \cite{MilnorAttractor}]\label{Lemma-MilnorAttractor}
Let $K$ be a compact invariant set such that $\Leb(\Bas(K))>0$
then, there exists $K' \en K$ a compact invariant set which is a
Milnor attractor.
\end{lema}

In some situations, we can have a stronger notion of attractor. We
say that a compact $f$-invariant set $K$ is a \emph{minimal Milnor
attractor} if $\Leb(\Bas(K))>0$ and $\Leb(\Bas(K'))=0$ for every
$K' \en K$ compact invariant subset different from $K$.

\subsection{Connecting lemmas}\label{SubSection-ConnectingLemmas}

In the study of ``typical'' dynamics in the space of
$C^1$-diffeomorphisms, the main tool is the study of periodic
orbits which we hope to describe accurately the recurrent
behavior. It is then important to control perturbations of orbits
in order to create the desired behavior. This section reviews
several orbit perturbation results (in
\ref{Subsection-FranksGourmelon} we shall review perturbations of
the derivative which is the other main tool in the study of
typical $C^1$-behavior) such as the closing and connecting lemmas
and their consequences. All this results were for many time
considered extremely difficult and technical. Nowadays, even if
they remain subtle, many proofs have been considerably improved
(see in particular \cite{Crov-Hab} for simple proofs of some of
the results and sketches of the rest). This section intends to be
a mere presentation of the results, for an introduction see
\cite{Crov-Hab} and \cite{BDV} appendix A.

We first introduce the well known \emph{Closing Lemma} of Pugh
(\cite{PughClosingLemma,PughDenisityTheorem}) and a very important
consequence which together with Kupka Smale's theorem was one of
the first \emph{genericity} results. Its simple and natural
statement may hide its intrinsic difficulties, the references
above explain why it is not that easy to perform such
perturbation.

\begin{teo}[Closing Lemma \cite{PughClosingLemma}]\label{Teorema-ClosingLemma}
Given $f\in \Diff^1(M)$, $\cU$ a neighborhood of $f$ and $x \in
\Omega(f)$, then, there exists $g \in \cU$ such that $x\in
\Per(g)$.
\end{teo}

The extension of the Closing Lemma to the $C^2$-topology is far
beyond reach of the current techniques except in certain cases
where one can control the recurrence and be able perform
perturbations in higher topologies (see for example
\cite{PughASpecialCrClosing,CP}).

As a consequence of the Closing lemma and the fact that
semicontinuous functions are continuous in a residual subset, Pugh
obtained the following consequence from his theorem.  We shall
only sketch the proof to show how to use the techniques (see for
example \cite{Crov-Hab} Corollary 2.8 for a complete proof).

\begin{cor}\label{Corolario-DensityTheorem}
There exists a $C^1$-residual subset $\cG \en \Diff^1(M)$ such
that if $f\in \cG$ one has that $\overline{\Per(f)}=\Omega(f)$.
\end{cor}

\esbozo{}Since hyperbolic periodic points persist under
$C^1$-perturbations (see Theorem
\ref{Teorema-PeriodicosHiperbolicosVariedadEstable}) we get that
the map $f \mapsto \overline{\Per_H(f)}$ which goes from
$\Diff^1(M)$ to $\cK(M)$ is semicontinuous (see Remark
\ref{Remark-Semicontinuidad}).

We obtain that there is a residual subset $\cG \en \Diff^1(M)$
consisting of diffeomorphisms where the map $f \mapsto
\overline{\Per_H(f)}$ varies continuously with respect to the
Hausdorff metric in $\cK(M)$.

We claim that if $f$ is such a continuity point, then
$\overline{\Per(f)} = \Omega(f)$. Indeed, if this was not the
case, using the Closing Lemma we could make the non-wandering set
explode by creating a periodic point far from
$\overline{\Per_H(f)}$ by an arbitrarily small
perturbation\footnote{Notice that it is not hard to perturb a
periodic point in order to make it hyperbolic.} contradicting the
continuity and concluding the proof. \lqqd

It may seem that creating a periodic point out of a recurrent (or
non-wandering) one is equally as difficult as creating a
connection between orbits $x$ and $y$ such that the omega-limit
set of one intersects the alpha-limit set of the other. However,
the difficulties that arise in this context are considerably
larger and it took a long time to handle this case.

The Connecting Lemma was first proven by Hayashi in
\cite{Hayashi}. Then, many improvements appeared (see
\cite{Arnaud-Connecting,WenXia,BC} for example). The statement we
present is taken from \cite{Crov-Haussdorf} Theorem 5 and it is
quite stronger.

\begin{teo}[Connecting Lemma \cite{Hayashi,Crov-Haussdorf}]\label{Teorema-ConnectingLemma}
Let $f\in \Diff^1(M)$ and $\cU$ a neighborhood of $f$. Then, there
exists $N>0$ such that every non periodic point $x \in M$ admits
two neighborhoods $W\en \hat W$ satisfying that: \bi \item[-] The
sets $\hat W, f(\hat W), \ldots, f^{N-1}(\hat W)$ are pairwise
disjoint. \item[-] For every $p,q \in M \backslash (f(\hat W) \cup
\ldots \cup f^{N-1}(\hat W))$ such that $p$ has a forward iterate
$f^{n_p}(p)\in W$ and $q$ has a backward iterate $f^{-n_q}(q)\in
W$, there exists $g \in \mathcal U$ which coincides with $f$ in $M
\backslash (f(\hat W) \cup \ldots \cup f^{N-1}(\hat W))$ and such
that for some $m>0$ we have $g^m(p)=q$. \ei Moreover, $\{p,g(p),
\ldots, g^m(p)\}$ is contained in the union of the orbits $\{p,
\ldots, f^{n_p}(p)\}$, $\{ f^{-n_q}(q), \ldots, q \}$ and the
neighborhoods $\hat W, \ldots, f^N(\hat W)$. Also, the
neighborhoods $\hat W, W$ can be chosen arbitrarily small.
\end{teo}

A much harder problem is to create orbits which realize in some
sense the $\eps$-pseudo-orbits since this will clearly require
making several perturbations. We shall state the consequences of a
connecting lemma for pseudo-orbits obtained in \cite{BC} since we
shall not use the perturbation result itself (previous partial
results can be found in \cite{Abdenur,BDConnexions, CMP, MP}).

\begin{teo}[\cite{BC}]\label{Teorema-BonattiCrovisier}
There exists a $C^1$-residual subset $\cG_{BC}$ of $\Diff^1(M)$
such that for $f \in \cG_{BC}$ one has that: \bi \item[-]
$\overline{\Per(f)} = \cR(f)$. \item[-] For $p \in \Per(f)$ we
have that $\cC(p)=H(p)$. In particular, homoclinic classes of $f$
are disjoint or equal. Moreover, if two periodic points in $H(p)$
have the same index, then they are homoclinically related.
\item[-]\emph{(Hurley's Conjecture)} There exists a residual
subset $R\en M$ such that for every $x\in R$ we have that
$\omega(x)$ is a quasi-attractor. \item[-] If $\cQ$ is an
essential attractor, then there exists $U$ a neighborhood of $\cQ$
such that for a residual set of points in $U$ the $\omega$-limit
is contained in $\cQ$.
 \item[-] If a chain-recurrence
class $\cC$ is isolated, then, there exists $U$ a neighborhood of
$\cC$ and $\cU$ a neighborhood of $f$ such that for every $g\in
\cU$ the maximal invariant subset of $U$ is chain-recurrent.
\item[-] \emph{(\cite{CMP})} The closure of the unstable manifold
of a periodic orbit is a Lyapunov stable set. Moreover, the
homoclinic class of a periodic point $p\in M$ is $H(p)=
\overline{W^s(p)}\cap \overline{W^u(p)}$. \ei
\end{teo}

When a chain recurrence class $\cC$ has no periodic points we say
that $\cC$ is an \emph{aperiodic class}.

As a direct consequence we obtain the following properties:

\begin{cor}[\cite{BC}]\label{Corolario-BonattiCrovisier}
For $f\in \cG_{BC}$ one has that: \bi \item[-] If a
chain-recurrence class $\cC$ is isolated, then $\cC$ is a
homoclinic class. \item[-] If a chain-recurrence class $\cC$ has
non-empty interior,  then it is a bi-Lyapunov stable homoclinic
class. \ei
\end{cor}

\dem The statement of those classes being homoclinic classes is a
direct consequence of the fact that periodic points are dense in
the chain-recurrence set and that chain-recurrence classes
containing periodic points coincide with the homoclinic classes of
the periodic points.

To prove that a homoclinic class with non-empty interior is
bi-Lyapunov stable, notice that since it contains the unstable
manifold of any periodic orbit in its interior, from the last
statement of Theorem \ref{Teorema-BonattiCrovisier} it follows
that it must be a quasi-attractor. The argument is symmetric and
it also shows that it must also be a quasi-repeller. \lqqd

The first statement of this corollary poses the following natural
question (see \cite{BDV} Problems 10.18 and 10.22):

\begin{quest}
Is an isolated homoclinic class of a $C^r$-generic diffeomorphism
$C^r$-robustly transitive\footnote{We say that a chain recurrence
class $C$ of a diffeomorphism $f$ is $C^r$-\emph{robustly
transitive} if there exists a $C^r$-neighborhood $\cU$ of $f$ and
a neighborhood $U$ of $C$ such that the maximal invariant of $U$
for $g\in \cU$ is transitive.}?
\end{quest}

We have given a negative answer to this question with C. Bonatti,
S. Crovisier and N. Gourmelon in \cite{BCGP} (see Appendix
\ref{Apendice-BCGP})

In general, another difficult problem is to control that when one
perturbs a pseudo-orbit in order to create a real orbit then the
new orbit essentially ``shadows'' the pseudo-orbit (this will be
better explained later). In general, this is not possible
(\cite{BDT}), however, one can approach pseudo-orbits for
$C^1$-generic diffeomorphisms with a weak form of shadowing (see
\cite{Arnaud-Haussdorf} for a previous related result).

\begin{teo}[\cite{Crov-Haussdorf} Theorem 1]\label{Teorema-ShadowingCrivisier}
There exists a $C^1$ residual subset $\cG_{H}$ of $\Diff^1(M)$
such that for any $\delta>0$ there exists $\eps>0$ such that for
any $\eps$-pseudo-orbit $\{z_0, \ldots, z_k\}$ there exist a
segment of orbit $\{x, \ldots, f^m(z_k)\}$ which is at Hausdorff
distance smaller than $\delta$ from the pseudo-orbit. Moreover, if
the $\eps$-pseudo-orbit is periodic (i.e. $z_k=z_0$) then one can
choose the orbit to be periodic.
\end{teo}

\subsection{Invariant measures and the ergodic closing lemma}\label{Subsection-MedidasInvariantes}

When studying the orbit structure of diffeomorphisms it is
sometimes important to understand the recurrence from a more
quantitative viewpoint to have better control on how the recurrent
points affect the orbits which pass close to them. A main tool for
measuring recurrence is ergodic theory which treats dynamics of
bi-measurable, measure preserving transformations of measure
spaces. When combined with topological dynamics one can obtain
lots of information some of which will be used in this text. We
will present here some of this theory and refer to the reader to
\cite{ManheLibro} for a more complete review of ergodic theory of
differentiable dynamics.

Let $f\in \Diff^1(M)$ and $\mu$ a regular (probability) measure in
the Borel $\sigma$-algebra of $M$. We shall denote the set of
regular probability measures of $M$ as $\cM(M)$. We say that $\mu
\in \cM(M)$ is \emph{invariant} if $\mu(f^{-1}(A)) = \mu(A)$ for
every measurable set $A$. We denote the set of invariant measures
of $f$ as $\cM_f(M) \en \cM(M)$. With the weak-$\ast$ topology in
the space of  measures of $M$ we know that $\cM(M)$ is convex and
compact. It is easy to see that $\cM_f(M)$ is also convex and
compact.

An invariant measure $\mu$ is called \emph{ergodic} if and only if
$f$-invariant measurable sets have $\mu$-measure $0$ or $1$. We
denote the set of ergodic measures as $\cM_{erg}(M)$ which can be
seen to be the set of extremal points of $\cM_f(M)$. We deduce
from the usual Krein-Milman theorem (see \cite{RudinFA}) the
following consequence which will allow us in general to
concentrate in ergodic measures (see \cite{ManheLibro} II.6 for a
more general statement):

\begin{prop}\label{Proposicion-MedidasErgodicas}
Let $\mu$ be an invariant measure and $\varphi: M \to \RR$ such
that $$ \int_M \varphi d\mu >0$$ \noindent then, there exists an
ergodic measure $\mu'$ whose support is contained in the support
of $\mu$ and such that
$$\int_M \varphi d\mu'
>0.$$ The same result holds\footnote{In fact, for the non-strict inequalities one needs to use the
ergodic decomposition theorem (see \cite{ManheLibro} II.6).} for
$\geq, <, \leq$.
\end{prop}

The importance of measuring the integral of real-valued functions
with invariant measures is given by the well known Birkhoff
ergodic theorem which guaranties that for knowing such integrals
it is enough to average the values obtained in the orbit of a
generic point:

\begin{teo}[Birkhoff Ergodic Theorem]\label{Teorema-Birkhoff}
Let $f\in \Diff^1(M)$ and $\mu \in \cM_f(M)$ an invariant measure.
Given $\varphi: M \to \RR$ a $\mu$-integrable function we have
that for $\mu$-almost every point $x$ there exists
$$ \lim \frac{1}{n} \sum_{i=0}^{n-1} \varphi\circ f^i(x) = \tilde \varphi(x). $$
Moreover, $\tilde \varphi(x)$ is $\mu$-integrable and
$f$-invariant \emph{(} i.e. $\tilde \varphi(f(x))=\tilde
\varphi(x)$\emph{)} and it satisfies that
$$\int_M \varphi d\mu = \int_M \tilde \varphi d\mu
$$
\end{teo}

Since when $\mu$ is ergodic $f$-invariant functions are (almost)
constant we get that for $\mu$-almost every point $x\in M$:

$$ \lim \frac{1}{n} \sum_{i=0}^{n-1} \varphi\circ f^i(x) = \int \varphi d\mu $$

One can define the \emph{statistical basin} of an ergodic measure
$\mu$ as the set of points whose averages with respect to any
continuous function $\varphi$ converge towards $\int \varphi
d\mu$:

$$ \Bas(\mu) = \left\{ y \in M \ : \ \forall \varphi \in C^0(M,\RR) \ , \ \frac{1}{n}\sum_{i=0}^{n-1}\varphi \circ f^i(y) \to \int_M \varphi d\mu \right\} $$

A direct consequence of Birkhoff theorem is that $\mu(\Bas(\mu))
=1$, however, a $C^1$-generic diffeomorphism verifies that the
support of invariant measures is very small from the topological
point of view \cite{ABC} (as well as from the point of view of
Lebesgue measure \cite{AvilaBochi-AbsolutelyContinuous}). This
suggest the following definition of measures which are sometimes
also called \emph{ergodic attractors} (see \cite{BDV} chapter 11
and references therein for an introduction in the context we are
interested in):

\begin{defi}[SRB measures] We say that an ergodic $f$-invariant measures is an \emph{SRB measure} for $f$ iff:

$$ \Leb(\Bas(\mu)) >0 $$\finobs
\end{defi}

Notice that in general one has that $\supp(\mu) \en
\overline{\Bas(\mu)}$ and that $\Bas(\mu) \en \Bas(\supp(\mu))$ so
we obtain that the support of an SRB measure always contains a
Milnor attractor (see Lemma \ref{Lemma-MilnorAttractor}). In some
(quite usual) circumstances, one can in fact guaranty that the
support of an SRB measure is indeed a minimal Milnor attractor.

There is a way of generalizing the notion of eigenvalues of the
derivative for points which are not periodic:

\begin{defi}[Lyapunov Regular Points]
Given a $C^1$-diffeomorphism $f$ of a $d$-dimensional manifold
$M$, we say that a point $x\in M$ is \emph{Lyapunov regular} if
there are numbers $\lambda_1(x) <\lambda_2(x)< \ldots <
\lambda_{m(x)} (x)$ and a decomposition $T_x M = E_1(x) \oplus
\ldots \oplus E_{m(x)}(x)$ such that for every $1\leq j\leq m(x)$
and every $v\in E_j(x)\backslash \{0\}$ we have:

$$ \lim_{n\to \pm \infty} \frac 1 n \log \|D_xf^n v \| = \lambda_j(x) $$

The numbers $\lambda_j(x)$ are called \emph{Lyapunov exponents} of
the point $x$ and the space $E_j(x)$ is called the \emph{Lyapunov
eigenspace} associated to $\lambda_j(x)$. We shall denote as
$\Reg(f)$ to  the set of Lyapunov regular points of $f$.\finobs
\end{defi}

A remarkable fact is that given an invariant measure, typical
points with respect to the measure are Lyapunov regular (see in
contrast Theorem 3.14 of \cite{ABC}).

\begin{teo}[Oseledet's Theorem \cite{Oseledet}]\label{Teorema-Oseldedet}
Given $f\in \Diff^1(M)$ and an $f$-invariant measure $\mu$ we have
that the set of Lyapunov regular points has $\mu$ total measure
\emph{(}$\mu(\Reg(f))=1$\emph{)}.
\end{teo}

\begin{obs}\label{Remark-EspaciosLyapunovYExponentesInvariantes}
Notice that the set of Lyapunov regular points of $f$ is
$f$-invariant. Indeed, $\lambda_i(f(x))=\lambda_i(x)$ and
$E_i(f(x))= D_xf(E_i(x))$. This implies in particular that given
an ergodic measure $\mu$ the Lyapunov exponents of $\mu$-generic
points is constant. We can thus define $\lambda_i(\mu)$ for an
ergodic measure\footnote{We integrate the quantities only as a way
of saying that it equals the value for almost every point.} as
$\int \lambda_i(x) d\mu$ and $m(\mu)= \int m(x) d\mu$. \finobs
\end{obs}

Given a periodic orbit $\cO$ of a diffeomorphism $f$ one can
define the following $f$-invariant ergodic measure:

$$ \mu_\cO = \frac{1}{\# \cO} \sum_{x\in \cO} \delta_x $$

To some extent, the closing lemma and pseudo-orbit connecting
lemma give information on the union of the supports of these
measures for $C^1$-generic diffeomorphisms and Theorem
\ref{Teorema-ShadowingCrivisier} says that in fact every
chain-transitive set can be approached in the Hausdorff topology
by the supports of such measures for a $C^1$-generic
diffeomorphism. A very important tool yet to be presented is the
celebrated Ergodic Closing Lemma of Ma\~ne (\cite{ManheErgodic})
which asserts that in fact one can perturb a generic point of a
measure in order to shadow it. See \cite{Crov-Hab} section 4.1 for
a simple and modern proof of this result. An important consequence
(with some improvement) is the following:

\begin{teo}[Ergodic Closing Lemma \cite{ManheErgodic}, \cite{ABC} Theorem 3.8]\label{Teorema-ErgodicClosingLemma}
There exists a $C^1$-residual subset $\cG_{E} \en \Diff^1(M)$ such
that if $f\in \cG_E$ and $\mu$ is an ergodic measure, then, there
exists $\cO_n$ a sequence of periodic orbits such that:
   \bi \item[-] The measures $\mu_{\cO_n}$ converge towards $\mu$ in the weak-$\ast$ topology.
   \item[-] The supports of the measures $\mu_{\cO_n}$ converge towards the support of $\mu$ in the Hausdorff topology.
   \item[-] $m(\mu_{\cO_n})= m(\mu)$ for every $n$ and the Lyapunov exponents $\lambda_i(\mu_{\cO_n})$ converge towards $\lambda_i(\mu)$.
   \ei
\end{teo}

We say an invariant ergodic measure $\mu$ is \emph{hyperbolic} if
all its Lyapunov exponents are different from $0$.




\section{Invariant structures under the tangent map}\label{Section-EstructurasInvariantes}

\subsection{Cocycles over vector bundles}\label{SubSection-CociclosYFibradosInvariantes}

Consider a homeomorphism $f: X \to X$ of a metric space $X$ and a
vector bundle $p:\cX \to X$. We say that $A: \cX \to \cX$ is a
linear cocycle over $f$ if it is a homeomorphism, we have that $f
\circ p(v) = p \circ A(v)$ and $A_x: p^{-1}(\{x\}) \to
p^{-1}(\{f(x)\})$ is a linear isomorphism.

This general abstract context particularizes to several
applications of which we will only be interested in two:

\bi
 \item[-] The \emph{derivative of a diffeomorphism} $Df: TM \to TM$ is a linear
cocycle over $f$ where the vector bundle is given by the trivial
projection $TM \to M$.
 \item[-] $\cA= (\Sigma, f, E, A, d)$ is a \emph{large period linear cocycle} of dimension $d$ and bounded by $K$ iff
   \bi
   \item[] $f: \Sigma \to \Sigma$ is a bijection such that every point of $\Sigma$ is periodic and such that given $n>0$ there
   are finitely many points in $\Sigma$ of period smaller than $n$ (in particular, $\Sigma$ is at most countable).
   \item[] $E$  is a vector bundle over $\Sigma$, this is, there exists $p: E \to \Sigma$ such that $p^{-1}(\{x\})=E_x$  is a $d$-dimensional vector space endowed with a
   Euclidean metric $\langle \cdot, \cdot \rangle_x$.
   \item[] $A: x \in \Sigma \mapsto A_x \in GL(E_x,E_{f(x)})$ is such that $\|A_x\| < K$ and $\|A_x^{-1}\| <K$.
   \ei
\ei \finobs

In fact, one can think (and it is what it will represent) of large
period linear cocycles as the restriction of the derivative of a
diffeomorphism to a subset of periodic points (which of course can
be a unique periodic orbit). For this reason, in the core of this
text we shall restrict ourselves to the study of the derivative of
diffeomorphisms (and the restriction of that cocycle to invariant
subsets), however, in Appendix \ref{Appendix-PerturbacionCociclos}
we will use the formalism of large period linear cocycles where
some quantitative results can be put in qualitative form.

\subsection{Dominated splitting}\label{SubSection-DescomposicionDominada}

Consider $f \in \Diff^1(M)$ and $\Lambda$ an $f$-invariant subset
of $M$. We say that a subbundle $E \en T_\Lambda M$ is
$Df$-invariant if we have that $$ Df(E(x)) = E(f(x)) \qquad
\forall x \in \Lambda$$

Given $E$ and $F$ two $Df$-invariant subbundles of $T_\Lambda M$
we say that $F$ $\ell$-\emph{dominates} $E$ (and we denote it as
$E\prec_\ell F$) iff for every $x\in \Lambda$ and any pair of unit
vectors $v_E \in E(x)$ and $v_F \in F(x)$ we have that

$$ \| D_x f^\ell v_E \| < \frac 1 2  \|D_x f^\ell v_F \| $$

In general, we say that $F$ \emph{dominates} $E$ (which we denote
$E \prec F$) iff there exists $\ell$ such that $E\prec_\ell F$.

In some examples, there is a stronger form of domination, called
\emph{absolute domination} (in contrast to the previous concept
sometimes called \emph{pointwise domination}). We say that $F$
\emph{absolutely} $\ell$-\emph{dominates} $E$ (and we denote it as
$E\prec^{ab}_\ell F$) iff for any pair of points $x,y \in \Lambda$
and any pair of unit vectors $v_E \in E(x)$ and $v_F \in F(y)$ we
have that:

$$ \| D_x f^\ell v_E \| < \frac 1 2  \|D_y f^\ell v_F \| $$

\noindent In a similar manner as above, we say that $F$
\emph{absolutely dominates} $E$ (and denote it as $E\prec^{ab} F$)
if there exists $\ell>0$ such that $E\prec^{ab}_\ell F$.

It is possible to define domination in other ways. Notice that the
concept of a bundle $\ell$-dominating another one (both absolutely
as pointwisely) is dependent on the metric chosen in $TM$.
However, if a subbundle $E$ is dominated by other subbundle $F$
then this does not depend on the metric chosen. See
\cite{GourmelonAdaptada} for information on possible changes of
metric to get domination.

We say that an $f$-invariant subset $\Lambda$ admits a
\emph{dominated splitting} if $T_\Lambda M =E \oplus F$ where $E$
and $F$ are non-trivial $Df$-invariant subbundles and $E\prec F$.
If one has that $E\prec^{ab} F$ then we say that the dominated
splitting is \emph{absolute}.

More generally, a $Df$-invariant decomposition $T_\Lambda M=E_1
\oplus \ldots \oplus E_k$ over an $f$-invariant subset $\Lambda$
is called a \emph{dominated splitting} if for every $1 < i < k$ we
have that $(E_1 \oplus \ldots \oplus E_{i-1}) \prec (E_i \oplus
\ldots \oplus E_k)$. One can extend to absolute domination in a
trivial manner.

A notational parentheses is that we will make a difference in
$E\oplus F$ and $F \oplus E$ since in the first case we shall
understand that $E\prec F$ and in the second one that $F \prec E$
(this is not always the notation used in the literature).

We shall now give some properties of dominated splittings, the
proofs can be found in \cite{BDV} appendix B.

\begin{prop}[Uniqueness]\label{Proposicion-UnicidadDescDom}
Let $T_\Lambda M=E_1 \oplus \ldots \oplus E_k$ and $T_\Lambda
M=F_1 \oplus \ldots \oplus F_k$ be two dominated splitting for $f$
with $\dim E_i = \dim F_i$ for every $i$. Then, $E_i = F_i$ for
every $i$.
\end{prop}

The uniqueness of the decomposition with fixed dimensions of the
subbundles allows one to consider the maximal possible
decomposition which we shall call the \emph{finest dominated
splitting}. When a set does not admit any dominated splitting we
will say that its finest dominated splitting is in only one
subbundle (we shall explicitly mention the possibility).

\begin{prop}[Finest dominated splitting]\label{Proposition-DescDomMasFina}
If an $f$-invariant set $\Lambda$ admits a (non-trivial) dominated
splitting, then there exists a dominated splitting $T_\Lambda M =
E_1 \oplus \ldots \oplus E_k$ ($k\geq 2$) such that every other
dominated splitting $T_\Lambda M = F_1 \oplus \ldots \oplus F_l$
verifies that $l \leq k$ and that for every $1 \leq j \leq l$ one
has that $F_j = E_i \oplus \ldots \oplus E_{i+t}$ for some $1 \leq
i \leq k$ and $0\leq t \leq k-i$.
\end{prop}

The fact that the invariant set needs not be compact will play an
important role since we shall mainly study the behavior over
periodic orbits and then use the following proposition to extend
that behavior to the closure (it was in fact this property which
made Ma\~ne, and probably also Liao and Pliss, consider this
notion for attacking the stability conjecture, see
\cite{ManheErgodic, ManheIHES})

\begin{prop}[Extension to the closure and to the limit]\label{Proposicion-DominacionPasaALaClausuraYLimite}
Let $f_n \in \Diff^1(M)$ be a sequence of diffeomorphisms
converging to a diffeomorphism $f$ and let $\Lambda_n$ be a
sequence of $f_n$-invariant sets admitting $\ell$-dominated
splittings $T_{\Lambda_n} M = E_1^n \oplus \ldots \oplus E_k^n$
such that $\dim E_i^n$ does not depend on $n$ nor on the point.
Then, if
$$\Lambda = \limsup \Lambda_n= \bigcap_{N>0} \overline{ \bigcup_{n>N} \Lambda_n } $$
\noindent then $\Lambda$ is a compact $f$-invariant set which
admits a dominated splitting $$T_\Lambda M = E_1 \oplus \ldots
\oplus E_k$$ \noindent such that $E_i(x) = \lim E_k^n(x^n_k)$ for
every $x^n_k \in \Lambda_n$ converging to $x$.
\end{prop}

\begin{obs}\label{Remark-ExtensionClausuraAngulosYContinuidad}
As a consequence of the previous proposition we obtain that if
$T_\Lambda M = E_1 \oplus \ldots \oplus E_k$ is a dominated
splitting for an $f$-invariant subset $\Lambda$ the following
properties are verified:
\bi
 \item[-] The closure $\overline \Lambda$ admits a dominated splitting $T_{\overline \Lambda} M = E'_1 \oplus \ldots \oplus E_k'$
 which extends the previous splitting (i.e. restricted to $\Lambda$ it  coincides with $E_1 \oplus \ldots \oplus E_k$).
 This follows by applying the previous proposition with $f_n=f$, $\Lambda_n =\Lambda$ and $E_i^n =E_i$.
 \item[-] The bundles $E_i$ vary continuously, this means, if $x_n$ is a sequence in $\Lambda$ such that $x_n \to x \in \Lambda$ then $E_i(x_n)\to E_i(x)$.
 This follows by applying the previous proposition with $f_n=f$, $\Lambda_n=\Lambda$, $E_i^n=E_i$ and $x_n^k = x_k$ converging to $x$.
 \item[-] There exists $\alpha>0$ such that for $i\neq j$ the angle\footnote{We can define the angle between two subbundles as $\arccos$ of the suppremum of the inner product
 between pairs of unit vectors, one in $E_i$ and the other in $E_j$.} between $E_i$ and $E_j$ is larger than $\alpha$. This follows from the fact that the bundles
 vary continuously and extend to the closure, so, if the angle is not bounded from below then there must be a point where two bundles have non-trivial intersection contradicting the fact that the sum is direct.
  \ei \finobs
\end{obs}

Given a vector space $V$ of dimension $d$ with an inner product
$\langle \cdot, \cdot \rangle$ and a $k$-dimensional subspace $E$
of $V$ we can express every vector of $V$ in a unique way as $v+
v^\perp$ where $v\in E$ and $v^\perp \in E^\perp$. We define
$\cE$, the $\alpha$-\emph{cone of} $E$ as:

$$ \cE = \{ v + v^\perp \in V \ : \ \|v^\perp\| \leq \alpha \|v\| \} $$

The interior of $\cE$ is denoted as $\Int(\cE)$ and is the
topological interior of $\cE$ together with the vector $0$. The
dimension of $\cE$ is the dimension of the largest subspace it
contains.

Given a subset $K$ of a manifold $M$ a $k$-dimensional \emph{cone
field} is a continuous association of cones in $T_xM$ to points
$x$ in $K$. It will be given by a continuous $k$-dimensional
subbundle $E \en T_K M$ together with a continuous function
$\alpha: K \to \RR$: so, the cone field will associate to $x\in K$
the $\alpha(x)$-cone of $E(x)$  which we shall denote as $\cE(x)$.

One can define more general cones and cone fields and this is
useful in other contexts (\cite{BoGou}), but for us it will
suffice to consider this notion.

\begin{prop}[Cone fields]\label{Proposition-ConeFields}
If $\Lambda$ is a $f$-invariant set admitting a dominated
splitting $T_\Lambda M = E\oplus F$. Then, there exists an open
neighborhood $U$ of $\Lambda$ and a $\dim F$-dimensional cone
field $\cE$ defined in $U$ such that for every $x\in U$ such that
$f(x)\in U$ one has that
$$ Df(\cE(x)) \en \Int(\cE(f(x))) $$
Conversely, if there exists a $k$-dimensional cone field $\cE$
defined in an open subset $U$ of $M$ verifying that if $x\in U$
and $f(x) \in U$ it satisfies $Df(\cE(x)) \en \Int(\cE(f(x)))$
then, if $\Lambda$ is the maximal invariant subset of $U$ we have
a dominated splitting $T_\Lambda M = E\oplus F$ with $\dim F=k$.
\end{prop}

We say that a set $\Lambda$ admits a \emph{dominated splitting of
index} $k$ if there is a dominated splitting $T_\Lambda M= E\oplus
F$ with $\dim E= k$.

\begin{obs}[Robustness of dominated splitting]\label{Remark-RobustnessDominatedSplitting}
The previous proposition allows one to show that if $\Lambda$ is
an $f$-invariant set which admits a dominated splitting of index
$k$, then there exists a neighborhood $U$ of $\Lambda$ and a
neighborhood $\cU$ of $f$ such that for every $g\in \cU$ the
maximal invariant set of $U$ for $g$ admits a dominated splitting
of index $k$. Indeed, the first part allows one to construct a
$k$-dimensional cone field in a neighborhood $U$ of $\Lambda$
which is $Df$-invariant in the sense explained in the statement of
the proposition. This invariance is not hard to see is robust in
the $C^1$-topology and it gives an open neighborhood $\cU$ of $f$
such that for every $g\in \cU$ the cone field will verify the
converse part of the proposition.  \finobs
\end{obs}

\subsection{Uniform subbundles}\label{SubSection-UniformBundles}

Given a $f$-invariant subset $\Lambda$ and a $Df$-invariant
subbundles $E\en T_\Lambda M$ we say that $E$ is \emph{uniformly
contracted} (resp. \emph{uniformly expanded}) if there exist $N>0$
(resp. $N<0$) such that for every $x\in \Lambda$:

$$ \|D_xf^N|_{E(x)} \| < \frac 1 2 $$

Given a $C^1$-diffeomorphism $f$ on a Riemannian  manifold $M$ we
define the continuous map $Jf|_E: M \to \RR$ such that $Jf|_E(x)$
is the $\dim E$-dimensional (oriented) volume of the
parallelepiped generated by the $Df$ image of a orthonormal basis
in $E$.

We say that a $Df$-invariant subbundle $E\en T_\Lambda M$ is
\emph{uniformly volume contracted} (resp. \emph{uniformly volume
expanding}) if there exists $N > 0$ (resp. $N<0$) such that for
every $x\in \Lambda$:

$$ |Jf^N|_{E(x)}(x)| < \frac 1 2  $$

Invariant measures may give a criteria for knowing weather
invariant subbundles are uniform (see \cite{Crov-CentralModels}
Claim 1.7):

\begin{prop}\label{Proposicion-ExponentesYFibradosUniformes}
Consider a $Df$-invariant continuous subbundle $E\en T_\Lambda M$
over a compact $f$-invariant subset $\Lambda \en M$. We have that:
\bi
 \item[(i)] $E$ is uniformly contracted if and only if for every invariant ergodic measure $\mu$ such that $\supp(\mu)\en \Lambda$
 we have that the largest Lyapunov exponent of $\mu$ whose Lyapunov eigenspace is contained in $E$ is strictly smaller than $0$.
 \item[(ii)] $E$ is uniformly volume contracted if and only if for every invariant ergodic measure $\mu$ such that $\supp(\mu) \en \Lambda$
 we have that the sum of all the Lyapunov exponents of $\mu$ whose Lyapunov eigenspaces are contained in $E$ is strictly smaller than $0$.
 In particular, if $E$ is one dimensional uniform volume contraction implies uniform contraction.
\ei An analogous statement holds for uniform expansion and uniform
volume expansion.
\end{prop}

\dem We first prove (i). If $E$ is uniformly contracted, it is
clear that for every vector $v \en E$ one has that

$$ \limsup_n \frac 1 n \log \|Df^n v \| < 0 $$

Which gives the direct implication. Now, assuming that $E$ is not
uniformly contracted, one can prove that there must be points $x_n
\in \Lambda$ and $v_n \in E(x_n)$ such that

$$ \|D_{x_n} f^j v_n \| \geq \frac 1 2  \quad 0\leq j \leq n $$

Consider the invariant measures $\mu_n = \frac 1 n
\sum_{i=0}^{n-1} \delta_{f^i (x_n)}$. By compactness of $\cM(M)$
we have that (modulo considering a subsequence)  $\mu_n \to \mu$
which will be an invariant measure as it is easy to check.

We claim that $\mu$ must have a Lyapunov exponent larger or equal
to $0$ inside $E$. Indeed, consider (again modulo considering
subsequences) the limit $F\en E$ of the subspaces generated by
$v_n$ which will be contained in $E$ by continuity of $E$. Let $x$
be a generic point in the support of $\mu$ and $v$ a vector in
$F(x)$. We have that for every $\eps>0$ there are points $x_n$
arbitrarily close to $x$ such that $v_n$ is arbitrarily close to
$v$. This implies that the derivative along $v$ for $x$ will not
be able to contract giving the desired Lyapunov exponent which is
larger or equal to $0$ (see \cite{Crov-CentralModels} Claim 1.7
for more details, in particular, formalizing this idea requires
passing to the unit tangent bundle of the manifold and consider
the measures there).

To prove (ii) one proceeds in a similar way by noticing that the
change of volume is related to the expansions and contractions in
an orthonormal basis and the angles to which they are sent. The
fact that the angles between Lyapunov eigenspaces vary
subexponentially is a consequence of a stronger version of
Oseledet's theorem (see for example \cite{KH} Theorem S.2.9).
\lqqd

With some more work one can prove the following result of Pliss
\cite{pliss} (see \cite{ABC} Lemma 8.4 for a proof):

\begin{lema}[Pliss]\label{Lemma-SiExponentesNegativosPozo}
Let $\mu$ be an ergodic measure such that all of its Lyapunov
exponents are negative, then $\mu$ is supported on a hyperbolic
sink.
\end{lema}

Now we are ready to define some notions which will in some sense
capture robust dynamical behavior as will be reviewed in
subsection \ref{SubSection-PropiedadesRobustas}.

\begin{defi}\label{Definicion-HypPHVH}
Let $f\in \Diff^1(M)$ and $\Lambda$ a compact $f$-invariant set
such that its finest dominated splitting is of the form
$T_{\Lambda} M = E_1 \oplus \ldots \oplus E_k$ (in this case we
allow $k=1$). We will say that:
\begin{itemize}
\item[-] $\Lambda$ is \emph{hyperbolic} if either $k=1$ and $E_1$
is uniformly expanded or contracted or there exists $1<j \leq k$
such that $E_1 \oplus \ldots \oplus E_{j-1}$ is uniformly
contracted and $E_j \oplus \ldots \oplus E_k$ is uniformly
expanded. \item[-] $\Lambda$ is \emph{strongly partially
hyperbolic} if $E_1$ is uniformly contracted and $E_k$ is
uniformly expanded. \item[-] $\Lambda$ is \emph{partially
hyperbolic} if either $E_1$ is uniformly contracted or $E_k$ is
uniformly expanded. \item[-] $\Lambda$ is \emph{volume partially
hyperbolic} if both $E_1$ is uniformly volume contracted and $E_k$
is uniformly volume expanded. \item[-] $\Lambda$ is \emph{volume
hyperbolic} if it is both volume partially hyperbolic and
partially hyperbolic.
\end{itemize}
\finobs
\end{defi}

The definitions of volume partial hyperbolicity and volume
hyperbolicity may vary in the literature as well as those of
partial hyperbolicty and strong partial hyperbolicity. We warn the
reader for that distinction.

\begin{obs}\label{Remark-HyperbolicidadYotras}
Using Proposition \ref{Proposicion-ExponentesYFibradosUniformes}
and Lemma \ref{Lemma-SiExponentesNegativosPozo} we get that if the
finest dominated splitting of a compact invariant set is trivial
(i.e. $k=1$) and the set is either hyperbolic or partially
hyperbolic then it must be a periodic sink or a source. When the
finest dominated splitting is not trivial we have the following
implications:

$$  \text{Hyperbolic} \Rightarrow \text{Strong Partially Hyperbolic} \Rightarrow $$
$$ \Rightarrow \text{Volume hyperbolic} \Rightarrow \text{Partially Hyperbolic} $$

Moreover if one extremal bundle is one-dimensional we have that:

$$ \text{Volume Partially Hyperbolic} \Rightarrow \text{Volume Hyperbolic} \Rightarrow \text{Partially Hyperbolic} $$

\finobs
\end{obs}

\begin{notacion}[Uniform bundles]
Let $\Lambda$  be a compact $f$-invariant set admiting a dominated
splitting of the form  $T_\Lambda M =E_1 \oplus \ldots \oplus E_k$
which is the finest dominated splitting (where $k$ may be equal to
$1$). Assume that $j$ is the largest value such that $E_j$
uniformly contracted and $l$ the smallest such that $E_{j+l}$ is
uniformly expanded. If we denote a $Df$-invariant subbundle of
$T_\Lambda M$ as $E^s$ it will be implicit that $E^s= E_1 \oplus
\ldots \ldots E_t$ with $t\leq j$. In a similar way, if we denote
a $Df$-invariant subbundle as $E^u$ it will be implicit that $E^u
= E_{j+t} \oplus \ldots \oplus E_k$ with $t\geq l$.

In certain situations we may separate $E^s=E^{ss} \oplus E^{ws}$
(or $E^u= E^{wu}\oplus E^{uu}$) which will denote that the
contraction in $E^{ss}$ is stronger than the one in $E^{ws}$.
 \finobs
\end{notacion}

An important part of this thesis will be devoted to study
diffeomorphisms such that the whole manifold is a partially
hyperbolic (or strong partially hyperbolic) set. We shall say that
$f$ is \emph{Anosov} (resp. \emph{partially hyperbolic}, resp.
\emph{strong partially hyperbolic}, resp. \emph{volume
hyperbolic}) if $M$ is a hyperbolic (resp. partially hyperbolic,
resp. strong partially hyperbolic, resp. volume hyperbolic) set
for $f$. We will review these concepts with more detail later.

We remark that there are alternative definitions of these global
concepts, for example, it is usual (see \cite{Crov-Hab}) to name a
diffeomorphism \emph{hyperbolic} if its chain-recurrent set is
hyperbolic, in a similar way, it is defined in \cite{CSY} a
diffeomorphism to be \emph{partially hyperbolic} if its
chain-recurrent set admits a partially hyperbolic splitting.

\begin{notacion}[Absolute Notions] In many examples one gets a stronger version of these
concepts which is given by the fact that the domination provided
by the definitions (between the uniform bundles and the
``central'' or ``neutral'' ones)  may be absolute instead of
pointwise as we have been working with. In those cases we will add
the word absolute before partial hyperbolicity, strong partial
hyperbolicity or volume hyperbolicity depending on the context.
Notice that in the hyperbolic case both notions coincide since
uniform bundles are naturally absolutely dominated (this is
another reason for choosing sometimes the definition of absolute
domination). \finobs
\end{notacion}

We obtain the following robustness property which is quite
straightforward from the definitions and Remark
\ref{Remark-RobustnessDominatedSplitting}.

\begin{prop}[Robustness]\label{Proposition-Robustez}
Assume that $\Lambda$ is a compact $f$-invariant set which is
hyperbolic, then there exists $U$ a neighborhood of $\Lambda$ and
$\cU$ a $C^1$-neighborhood of $f$ such that for every $g\in \cU$
the maximal invariant set of $g$ in $U$ is also hyperbolic. The
same holds for the concepts of partial hyperbolicity, strong
partial hyperbolicity, volume hyperbolicity, volume partial
hyperbolicity and the absolute versions. \end{prop}

\subsection{Franks-Gourmelon's Lemma}\label{Subsection-FranksGourmelon}

In this section we shall review some techniques of perturbation
which allow to change the derivative of the diffeomorphism over a
periodic orbit  by a small $C^1$-perturbation. Notice that this
cannot be done in higher topologies (not even in $C^2$ see
\cite{PujSamDominated}).

The classical Franks' Lemma (\cite{FranksLema}) states the
following:

\begin{teo}[Franks' Lemma \cite{FranksLema}]\label{Teorema-FranksLemma}
Given a $C^1$-neighborhood $\cU$ of a diffeomorphism $f$, there
exists $\eps>0$ such that:
\begin{itemize}
\item[-] given any finite set $\{x_1,\ldots,x_k \}$ in $M$,
\item[-] any neighborhood $U$ of this finite set \item[-] any set
of linear transformations $A_i: T_{x_i}M \to T_{f(x_i)}M$
verifying that $\|A_i - D_{x_i}f \| < \eps$ for every $1\leq i
\leq k$ \end{itemize} \noindent then there exists $g\in \cU$ such
that:
\begin{itemize}
\item[-] $g=f$ outside $U$. \item[-] $g(x_i)=f(x_i)$ for every
$1\leq i\leq k$. \item[-] $D_{x_i}g = A_i$ for every $1\leq i\leq
k$.
\end{itemize}
\end{teo}

In sections \ref{Section-AtractoresSuperficies} and
\ref{Section-GenericBiLyapunov} we will use a stronger version of
this Lemma which allows to have control on invariant manifolds of
periodic points when one perturbs their derivatives. For doing
this it is important to have a better understanding of the way one
perturbs the cocycle of derivatives in order to make the
perturbation step by step and in some sense ``follow'' the
invariant manifolds.

Consider $f$ a $C^1$-diffeomorphism. We denote as $\Per_{j}(f)$
the set of (stable) index $j$ hyperbolic periodic points. Let
$\mathcal{O}$ be a periodic orbit and $E$ a $Df$ invariant
subbundle of $T_{\mathcal O}M$. We denote as $D_{\mathcal
O}f_{/E}$ to the cocycle over the periodic orbit given by its
derivative restricted to the invariant subbundle as defined in
greater generality in subsection
\ref{SubSection-CociclosYFibradosInvariantes}.

Let $\mathcal O$ be a periodic orbit and $\cA_\cO$ be a linear
cocycle \footnote{Recall from subsection
\ref{SubSection-CociclosYFibradosInvariantes} that a linear
cocycle $\cA$ of dimension $n$ over a transformation $f:\Sigma \to
\Sigma$ can in this case be represented by a map $A:\Sigma \to
GL(n,\R)$. When one point $p\in \Sigma$ is $f-$periodic, the
eigenvalues of the cocycle at $p$ are the eigenvalues of the
matrix given by $A_{f^{\pi(p)-1}(p)} \ldots A_p$} over $\mathcal
O$.  We say that $\cA_{\cO}$ has a \emph{strong stable manifold of
dimension} $i$ if the eigenvalues $|\lambda_1|\leq |\lambda_2|
\leq \ldots \leq |\lambda_d|$ of $\cA_\cO$ satisfy that
$|\lambda_i| < \min\{1,|\lambda_{i+1}|\}$.

If the derivative of $\mathcal O$ has strong stable manifold of
dimension $i$ then classical results ensure the existence of a
local, invariant manifold $W^{s,i}_{\eps}(x)$ tangent to the the
subspace generated by the eigenvectors of these $i$ eigenvalues
and imitating the behavior of the derivative (see \cite{KH}
Theorem 6.2.8 for a precise formulation and recall Theorem
\ref{Teorema-PeriodicosHiperbolicosVariedadEstable}). In fact,
$W^{s,i}_\eps$ is characterized for being the set of points in an
$\eps$-neighborhood of $\cO$ such that the distance of future
iterates of those points and $\cO$ goes to zero exponentially at
rate faster than $\lambda_i + \eps$ with small $\eps$.

Let $\Gamma_i(\cO)$ be the set of cocycles over $\mathcal O$ which
have a strong stable manifold of dimension $i$.

We endow $\Gamma_i(\cO)$ with the following distance,
$d(\cA_\cO,\cB_\cO)=\max\{ \|\cA_\cO - \cB_\cO\| , \|\cA_\cO^{-1}
- \cB_\cO^{-1}\|\}$ where the norm is $$\|\cA_\cO\|=\sup_{p\in
\mathcal O } \{ \frac{\|A_p(v)\|}{\|v\|} \ ; \ v \in
T_pM\backslash \{0\}\}.$$

Let $g$ be a perturbation of $f$ such that the cocycles
$D_{\mathcal O}f$ and $D_{\mathcal O} g$ are both in
$\Gamma_i(\cO)$, and let $U$ be a neighborhood of $\mathcal O$. We
shall say that $g$ \emph{preserves locally the $i$-strong stable
manifold of $f$ outside $U$}, if the set of points of the
$i$-strong stable manifold of $\mathcal O$ outside $U$ whose
positive iterates do not leave $U$ once they entered it, are the
same for $f$ and for $g$.

We have the following theorem due to Gourmelon which allows to
perturb the derivative of periodic orbits while controlling the
position of the invariant manifolds of them.

\begin{teo}[\cite{GouFranksLemma}]\label{Teorema-FranksLemmaGourmelon}
Let $f$ be a diffeomorphism, and $\mathcal{O}$ a periodic orbit of
$f$ such that $D_{\mathcal O}f \in \Gamma_i(\cO)$ and let
$\gamma:[0,1]\to \Gamma_i(\cO)$ be a path starting at $D_{\mathcal
O}f$. Then, given a neighborhood $U$ of $\mathcal O$, there is a
perturbation $g$ of $f$ such that $D_{\mathcal O}g=\gamma(1)$, $g$
coincides with $f$ outside $U$ and preserves locally the
$i$-strong stable manifold of $f$ outside $U$. Moreover, given
$\cU$ a $C^1$ neighborhood of $f$, there exists $\eps>0$ such that
if $\diam (\gamma)<\eps$ one can choose $g \in \cU$.
\end{teo}

We observe that the Franks' lemma for periodic orbits (Theorem
\ref{Teorema-FranksLemma}) is the previous theorem with $i=0$.
Also, we remark that Gourmelon's result is more general since it
allows to preserve at the same time more than one strong stable
and more than one strong unstable manifolds (of different
dimensions, see \cite{GouFranksLemma}).

\subsection{Perturbation of periodic cocycles}\label{SubSection-BDPBGVBoBo}

In view of the techniques of perturbation of the derivative over
finite sets of points reviewed in the previous section, it makes
sense to try to understand what type of behavior one can create by
(small) perturbations of the derivative of periodic orbits.

Of course, eigenvalues depend continuously on the matrices, so a
small perturbation has only small effect on the derivative over a
periodic orbit. However, the fact that we can perturb a small
amount but on many points at a time gives that it is sometimes
possible to get large effect by making a small perturbation (of
the diffeomorphism) by accumulation of these effects. It turns out
that the main obstruction for making such perturbations is the
existence of a dominated splitting.

The first results of this kind were obtained by Frank's itself in
his paper \cite{FranksLema}. However, the progress made in
\cite{ManheErgodic} started the systematic study of perturbations
of cocycles over periodic orbits.

Relevant development was obtained in \cite{BDP} where the concept
of \emph{transitions} was introduced. Later, in \cite{BGV} some
results were recovered without the need for transitions, and
recently, in \cite{BoBo}, a kind of optimal result was obtained
which in turn combines in a very nice way with the recent result
of \cite{GouFranksLemma} (see Theorem
\ref{Teorema-FranksLemmaGourmelon}).

In this section we shall present the results we shall use without
proofs.

The first result we shall state is the result from \cite{BDP}
which uses the notion of \emph{transitions}. It gives a dichotomy
between the existence of a dominated splitting and the creation of
homotheties by small perturbations along orbits.

\begin{teo}[\cite{BDP}]\label{Teorema-BDPCociclos}
Let $H$ be the homoclinic class of a hyperbolic periodic point $p$
of a $C^1$-diffeomorphism $f$. Let $\Sigma_p$ be the set of
periodic points homoclinically related to $p$ and assume that $E
\en T_{\Sigma_p}M$ is a $Df$-invariant subbundle. We have the
following dichotomy:
\begin{itemize}
 \item[-] Either $E= E_1 \oplus E_2$ were $E_i$ are $Df$-invariant subbundles and $E_1 \prec E_2$.
 \item[-] Or, for every $\eps>0$ there exists a periodic point $q \en \Sigma_p$ and a periodic
 linear cocycle $A: T_{\cO(q)} M \to T_{\cO(q)} M$ such that $\|A - D_{\cO(q)} f\| < \eps$ and we have that
 $A(f^{\pi(q)-1})\ldots A(q)$ is a linear homothety. Moreover, if $\det(D_pf^{\pi(p)}) \leq 1$ we can consider the homothety to be contracting.
\end{itemize}
\end{teo}

We shall not present a proof of this fact, the reader can consult
\cite{BDV} chapter 7 for  a nice sketch of the proof. We will give
though a proof of the following result to give the reader a taste
on the idea of considering transitions.

\begin{prop}[\cite{BDP}]\label{Proposicion-transitions}
Let $H$ be homoclinic class of a periodic point with
$|det(D_pf^{\pi(p)})|> 1$, then there is a dense subset of
periodic points in $H$ having the same property.
\end{prop}

\dem Let $U$ be an open set in $H$. There is a periodic point
$q\in U$ homoclinically related to $p$. Consider $x\in
W^s(\mathcal O(p))\trans W^u(\mathcal{O} (q))$ and $y \in
W^s(\mathcal O (q))\trans W^u(\mathcal{O}(p))$. The set $\mathcal
O(p) \cup \mathcal O(q) \cup \mathcal O(x) \cup \mathcal O(y)$ is
a hyperbolic set. So, using the shadowing lemma (see \cite{KH}
Theorem 6.4.15 for example) we can obtain a periodic point $r \in
U$, homoclinically related to $p$ such that its orbit spends most
of the time near $\mathcal O(p)$. Thus, it will satisfy that
$|det(D_rf^{\pi(r)})|>1$.

\lqqd

When we wish to use Theorem \ref{Teorema-FranksLemmaGourmelon} we
need to not only make small perturbations but also to make them in
small paths which do not affect the index of the periodic points
during the perturbation. The natural idea of considering the
straight line between the initial cocycle and the homothety falls
short of providing the desired perturbation and it is quite a
difficult problem to really realize the desired perturbation. A
recent result of J.Bochi and C.Bonatti (\cite{BoBo} which extends
previous development in this sense by \cite{BGV}) provides a
solution to this problem as well as it investigates which kind of
paths of perturbations can be realized in relation to the kind of
domination a cocycle admits. We shall state a quite weaker version
of their result and avoid the (very natural) mention to large
period linear cocycles and work instead with the derivative and
paths of perturbations. We refer the interested reader to Appendix
\ref{Appendix-PerturbacionCociclos} in order to get a more
complete account with complete proofs of partial results.

Recall from the previous subsection that given a periodic orbit
$\cO$ we denote $\Gamma_i(\cO)$ to be the set of cocycles over
$\cO$ which have strong stable manifold of dimension $i$ endowed
with the distance considered there.

\begin{teo}[\cite{BGV,BoBo}]\label{Teorema-BonattiBochi}
Let $f: M \to M$ be a $C^1$-diffeomorphism and $p_n$ a sequence of
periodic points whose periods tend to infinity and their orbits
$\cO_n$ converge in the Hausdorff topology to a compact set
$\Lambda$. Let
$$ T_\Lambda M = E_1 \oplus \ldots \oplus E_k $$
be the finest dominated splitting over $\Lambda$ (where $k$ may be
$1$). Then, for every $\eps>0$ there exists $n>0$ such that an
$\eps$-perturbation of the derivative along $\cO_n$ makes all the
eigenvalues of the orbit in the subspace $E_i$ to be equal.
Moreover, if the determinants of $D_{p_n}f^{\pi(p_n)}|_{E_i}$ have
modulus smaller than $1$ for every $n$ and the periodic orbits
have strong stable manifold of dimension $j$ (which must be
strictly larger than the dimension of $E_1 \oplus \ldots \oplus
E_{i-1}$) then, given $\eps>0$ there is $n>0$ and a path
$\gamma:[0,1]\to \Gamma_j(\cO_n)$ such that: \bi \item[-] $\diam
(\gamma) < \eps$. \item[-] $\gamma(0) =D_{\cO_n} f$. \item[-]
$\gamma(1)$ has all its eigenvalues of modulus smaller than $1$ in
$E_i$.\ei
\end{teo}

We recommend reading Lemma 7.7 of \cite{BDV} whose (simple)
argument can be easily adapted to give this result in dimension
$2$.

\subsection{Robust properties and domination}\label{SubSection-PropiedadesRobustas}

In this subsection we shall explain certain results which are
consequence of the perturbation results reviewed in the previous
subsections.

Possibly, one of the departure points of this study was the study
of the \emph{stability conjecture} finally solved in
\cite{ManheIHES}. We say that a diffeomorphism $f$ is
$\cR$-\emph{stable} if and only if there exists a
$C^1$-neighborhood of $f$ such that for every $g\in \cU$  the
diffeomorphism $g$ restricted to $\cR(g)$ is conjugated to $f$
restricted to $\cR(f)$. See \cite{Crov-Hab} section 7.7 for a
modern proof of the following result (which has been an underlying
motivation for the development of differentiable dynamics) and
more references for related results.

\begin{teo}[Stability Conjecture, \cite{PalisSmale,ManheIHES}]\label{Teorema-StabilityConjecture}
A diffeomorphism $f$ is $\cR$-stable if and only if $\cR(f)$ is
hyperbolic.
\end{teo}

One of the robust properties we will be most interested in is in
isolation of classes (and in particular quasi-attractors). For
homoclinic classes of generic diffeomorphisms, isolation is enough
to guaranty the existence of some invariant geometric structure.
The following result of \cite{BDP} was preceded by results in
\cite{ManheErgodic, DPU}.

\begin{teo}\label{TeoremaBonattiDiazPujals}
There exists a residual subset $\cG_{BDP} \en \Diff^1(M)$ such
that if $f\in \cG_{BDP}$ and $\cC$ is an isolated chain-recurrence
class, then, $\cC$ is volume partially hyperbolic. Moreover, if a
homoclinic class $H$ of a diffeomorphism $f\in \cG_{BDP}$ does not
admit any non-trivial dominated splitting then it is contained in
the closure of sinks and sources ($\overline{\Per_0(f) \cup
\Per_d(f)}$).
\end{teo}

\esbozo{} The proof of the first statement goes as follows: First,
by a classical Baire argument, one shows that there is a residual
subset of $\cG_{BDP}$ of $\Diff^1(M)$ such that if $\cC$ is an
isolated chain-recurrence class of a diffeomorphism $f\in
\cG_{BDP}$ then $\cC$ is a homoclinic class and there exists a
neighborhood $\cU$ of $f$ and a neighborhood $U$ of $\cC$ such
that for every $g\in \cG_{BDP} \cap \cU$ we have that the maximal
invariant set of $g$ in $U$ is a homoclinic class (see Theorem
\ref{Teorema-BonattiCrovisier} and the discussions after).

Now, consider $\cC$ an isolated chain-recurrence class of a
diffeomorphism $f\in \cG_{BDP}$. Let $p$ be a hyperbolic periodic
point of $f$ such that $\cC=H(p)$ and consider $\Sigma_p$ the set
of hyperbolic periodic points homoclinically related to $p$.

Assume that $\cC$ does not admit a non-trivial dominated
splitting, then, by Remark
\ref{Remark-ExtensionClausuraAngulosYContinuidad} we know that
$\Sigma_p$ cannot admit a dominated splitting. Now, by Theorem
\ref{Teorema-BDPCociclos} we know that we can make an arbitrarily
small perturbation of the derivative of some periodic point in
order to make it an homothety.

This perturbation can be made dynamically by using Theorem
\ref{Teorema-FranksLemma} creating a sink or a source inside $U$
the neighborhood of $\cC$. The chain-recurrence class of a sink or
a source is the point itself (see Remark
\ref{Remark-ClaseHomoclinicaYdeRecurrencia}), and since sinks and
sources persist under perturbations (Theorem
\ref{Teorema-PeriodicosHiperbolicosVariedadEstable}), we find a
diffeomorphism $g\in \cG_{BDP} \cap \cU$ having more than one
chain-recurrence class in $U$, a contradiction.

Now, let $T_\cC M = E_1 \oplus \ldots \oplus E_k$ be the finest
dominated splitting over $\cC=H(p)$. Assume that $E_1$ is not
uniformly volume contracting (the same argument applied to
$f^{-1}$ will show that $E_k$ is uniformly volume expanding).

By Proposition \ref{Proposicion-ExponentesYFibradosUniformes} we
know that there exists an ergodic measure $\mu$ supported on
$H(p)$ such that the sum of the Lyapunov exponents of $\mu$ in
$E_1$ is larger or equal to $0$. By the Ergodic Closing Lemma
(Theorem \ref{Teorema-ErgodicClosingLemma}) there are periodic
orbits $\cO_n$ converging in the Hausdorff topology towards
$\supp(\mu)$ and such that for $n$ large enough, the sum of the
eigenvalues of $\cO_n$ in the invariant bundle $E_1$ is larger
than or equal to $-\eps$ with small $\eps$. These periodic orbits
belong to $H(p)$ since  we have assumed that it is an isolated
chain-recurrence class.

Now, by the argument of Proposition \ref{Proposicion-transitions}
we get a dense\footnote{The fact that we are able to create a
dense subset of periodic orbits with such behavior is crucial in
the proof, since a priori we do not know if the measure for which
the volume contraction is not satisfied has total support or not,
and whether the finest dominated splitting inside that subbundle
is finer or not than the global one.} subset of $H(p)$ of periodic
orbits such that the sum of the eigenvalues in $E_1$ is larger
than $-\eps$ with arbitrarily small $\eps$. Since $E_1$ does not
admit a subdominated splitting in $H(p)$ we are able by using
again Theorem \ref{Teorema-BDPCociclos} to create sources inside
$U$ by small perturbations of $f$.

For the second statement, the proof is very similar. Consider a
homoclinic class $H$ of a diffeomorphism $f$ and we can assume
that the residual subset $\cG_{BDP}$ verifies that for $g$ in a
neighborhood of $f$ the continuation $H_g$ of $H$ is close to $H$
in the Hausdorff topology.

Using Proposition \ref{Proposicion-transitions} one can make the
periodic points which one can turn into sinks or sources
$\eps$-dense in $H$. A Baire argument then allows to show that if
the homoclinic class admits no-dominated splitting then it
contained in the closure of the set of sinks and sources.

\lqqd

We remark that Theorem \ref{Teorema-BonattiBochi} together with
Franks' Lemma (Theorem \ref{Teorema-FranksLemma}) allows also to
obtain that if a chain-recurrence class is not accumulated by
infinitely many sinks or sources then it admits a non-trivial
dominated splitting (see \cite{ABCProceedings}). This results can
be used in order to re-obtain the examples presented in
\cite{BDConnexions} where homoclinic classes of generic
diffeomorphisms accumulated by infinitely many sinks and sources
were constructed.

An immediate consequence is that we obtain a criterium for
guaranteeing that a homoclinic class is not isolated (this is a
key ingredient in constructing examples of dynamics with no
attractors).

\begin{cor}\label{Corollary-CriterioDeNoAislado} Let $f$ be a $C^1$-generic
diffeomorphism of $M$ and $H$ a homoclinic class of $f$ which is
not volume partially hyperbolic. Then, $H$ is not isolated.
\end{cor}

\begin{obs}\label{Remark-PozosYFuentesDensosEnLaClase}
Indeed, the proof of Theorem \ref{TeoremaBonattiDiazPujals} allows
one to show that if $f \in \cG_{BDP}$ and $H$ is a homoclinic
class  such that: \bi \item[-] The finest dominated splitting in
$H$ is of the form $T_HM = E_1 \oplus \ldots \oplus E_k$. \item[-]
$H$ has a periodic point $q$ verifying that
$\det(Df^{\pi(q)}|_{E_1}) \geq 1$\ei \noindent then $H$ is
contained in the closure of the set of sources of $f$.
\finobs
\end{obs}

Another consequence of Theorem \ref{TeoremaBonattiDiazPujals} is
the global characterization of diffeomorphisms which are robustly
transitive (or even, $C^1$-generic diffeomorphisms which are
transitive). The following result is one of the motivations of
many of the results we will present in this thesis (the optimality
of this result from the point of view of the obtained
$Df$-invariant geometric structure is given by the examples of
\cite{BV}):

\begin{cor}\label{Corollary-RobustamenteTransitivoImplicaVolumenHiperbolico}
If $f\in \cG_{BDP}$ and $\cR(f)=M$, then $f$ is volume partially
hyperbolic. Also, if $f \in \cG_{BDP}$ and there is a
chain-recurrence class $\cC$ of $f$ with non-empty interior, then
$\cC$ admits a non-trivial dominated splitting.
\end{cor}

In dimension 2, this result was shown by Ma\~ne in
\cite{ManheErgodic}: if a $C^1$-generic diffeomorphism of a
surface is transitive, then, $f$ is Anosov (recall Remark
\ref{Remark-HyperbolicidadYotras}). Indeed, together with a result
from Franks (\cite{FranksAnosov}) we get the following
characterization of robust behaviour in terms of the dynamics of
the tangent map:

\begin{teo}[\cite{FranksAnosov,ManheErgodic}]\label{TeoremaFranksManhe}
A diffeomorphism $f$ of a surface has a $C^1$-neighborhood $\cU$
such that every $g\in \cU$ is chain-recurrent if and only if $f$
is an Anosov diffeomorphism of $\TT^2$.
\end{teo}

A remarkable feature of this result is that it leaves in evidence
the fact that robust dynamical behavior is in relation with the
topology of the state space (and even the isotopy class). This
relation is given through the appearance of a geometric structure
which is invariant under the tangent map of the diffeomorphism.
This leads to the following idea whose understanding represents a
main challenge:

$$ \text{Robust dynamical behaviour} \Leftrightarrow \text{Invariant Structures} \Leftrightarrow \text{Topological Properties} $$

Other than in dimension $2$, very few is known in this respect
other than what it was reviewed in this section (which represents
hints on the direction of giving invariant geometric structures by
the existence of robust dynamical behavior). In dimension $3$, the
fact that Corollary
\ref{Corollary-RobustamenteTransitivoImplicaVolumenHiperbolico}
admits a stronger form suggests that it may be possible to search
for results with similar taste as Theorem \ref{TeoremaFranksManhe}
(we shall review some of the known results later):

\begin{teo}[\cite{DPU}]\label{Teorema-DPU} Let $M$ be a $3$-dimensional
manifold and $f \in \cG_{BDP}$ be chain-recurrent. Then, $f$ is
volume hyperbolic (i.e. partially hyperbolic and volume partially
hyperbolic). It can present the following forms of domination:
\begin{itemize}
\item[] $TM = E^{cs} \oplus E^u$, \item[] $TM = E^s \oplus E^c
\oplus E^u$ or \item[] $TM = E^s \oplus E^{cu}$. \end{itemize}
\end{teo}

We recommend reading \cite{BDV} chapter 7 for a review on robust
transitivity and for a survey of examples which show how these
results are optimal from the point of view of the geometric
structures obtained.

In dimension 3 the main examples of transitive strong partially
hyperbolic diffeomorphisms fall in the following classes: Fiber
bundles whose base is Anosov, Time one maps of Anosov flows or
Examples derived from Anosov diffeomorphisms in $\TT^3$. See
\cite{BonWilk}.

\begin{obs}\label{Remark-DominacionPointwise}
It is important to remark that all this results give pointwise
domination and not the absolute one. Indeed, it is not hard to
construct examples which verify all these robust properties and
fail to admit absolute  domination between the invariant
subbundles. \finobs
\end{obs}

\subsection{Homoclinic tangencies and domination}\label{SubSection-TangenciasYDominacion}

In this subsection we shall review certain properties of
diffeomorphisms  which are far from homoclinic tangencies. We
refer the reader to \cite{CSY,LVY} for the latest results on
dynamics of diffeomorphisms $C^1$-far from tangencies.

In this section we shall recall some result whose germ can be
traced to \cite{PujSamAnnals1} where it was proved that a
diffeomorphism of a surface which is far away from homoclinic
tangencies must admit a dominated splitting on the closure of the
saddle hyperbolic periodic points.

First, we define the notion of a homoclinic tangency: Given a
hyperbolic periodic saddle $p$ of a $C^1$-diffeomorphism $f$, we
say that $p$ has a \emph{homoclinic tangency} if there exists a
point of non-transverse intersection between $W^s(\cO(p))$ and
$W^u(\cO(p))$.

We denote as $\Tang \en \Diff^1(M)$ to the set of diffeomorphisms
$f$ having a hyperbolic saddle with a homoclinic tangency. As a
consequence of results which relate the existence of a dominated
splitting with the creation of homoclinic tangencies for periodic
orbits \cite{Wen-Tangencias1,Wen-Tangencias2} as well as some
adaptations of the ergodic closing lemma (see \cite{ABC}), in
\cite{Crov-CentralModels} the following result is proved:

\begin{teo}[\cite{Crov-CentralModels} Corollary 1.3]\label{Teorema-LejosTangenciasSylvain}
Let $f$ be a diffeomorphism in $\Diff^1(M) \backslash
\overline{\Tang}$. Then, there exists a $C^1$-neighborhood $\cU$
of $f$ integers $\ell, N>0$ and constants $\delta,\rho>0$ such
that for every $g\in \cU$ and every ergodic $g$-invariant measure
$\mu$ the following holds:

Let $x$ be a $\mu$-generic point and $T_{\cO(x)}M =E_- \oplus E^c
\oplus E_+$ be the Oseledet's splitting into the Lyapunov
eigenspaces corresponding respectively to Lyapunov exponents in
$(-\infty, -\delta)$, $[-\delta,\delta]$ and $(\delta,+\infty)$,
so:
\begin{itemize}
 \item[-] $\dim E^c \leq 1$.
 \item[-] The splitting $E_-\oplus E^c \oplus E_+$ is $\ell$-dominated (hence it extends to $\supp(\mu)$).
 \item[-] For $\mu$-almost every point we have that:
 $$  \lim_{k\to \infty} \frac 1 k \sum_{i=0}^{k-1} \log \|Df^N|_{E_-(f^{iN}(x))}\| < -\rho \quad  \lim_{k\to \infty} \frac 1 k \sum_{i=0}^{k-1} \log \|Df^{-N}|_{E_+(f^{-iN}(x))}\| > \rho         $$
\end{itemize}
\end{teo}

To finish this subsection, we present the following result of
\cite{ABCDW} which will be used later.

\begin{teo}[\cite{ABCDW}, \cite{GouTangencias}]\label{Teorema-ABCDW}
There exists a residual subset $\cG$ of $\Diff^1(M)$ such that if
$f \in \cG$ and $H$ is a homoclinic class of $f$ having periodic
points of (stable) index $s$ and $s'$. Then, for every $s<j<s'$ we
have that $H$ contains periodic points of index $j$. In
particular, if there is no perturbation of $f$ which creates a
homoclinic tangency for a periodic point in $H$ then $H$ admits a
dominated splitting of the form $T_H M =E \oplus G_1 \oplus \ldots
G_k \oplus F$ with $\dim E=s$, $\dim F=d-s'$ and $\dim G_j=1$ for
every $j$.
\end{teo}

We would like to point out that with the new techniques of
perturbation given by \cite{BoBo} and \cite{GouFranksLemma} one
can give a proof of the fact that a homoclinic class of a
$C^1$-generic diffeomorphism is \emph{index complete} which is
very direct: First, if a homoclinic class has periodic points of
index $i$ and $j$, then, by perturbing the derivative of some
periodic points it is possible to get periodic orbits of index in
between (by the results of \cite{BoBo}). After, the results of
\cite{GouFranksLemma} allows one to make the perturbation in order
to keep the necessary homoclinic relations in order guarantee that
the point remains in the homoclinic class after perturbation.

\subsection{Domination and non-isolation in higher regularity}\label{SubSection-PalisViana}

As well as in the case of $C^1$-topology, we can obtain a similar
criterium to obtain non-isolation of a homoclinic class for
$C^r-$generic diffeomorphisms combining the main results of
\cite{BDAbundanceTangencies} and \cite{PV} (it is worth also
mentioning \cite{Romero}). The only cost will be that we must
consider a new open set and that the accumulation by other classes
is not as well understood.

We state a consequence of the results in those papers in the
following result. We shall only use the result in dimension 3, so
we state it in this dimension, it can be modified in order to hold
in higher dimension but it would imply defining sectionally
dissipative saddles (see \cite{PV}).

\begin{teo}[\cite{BDAbundanceTangencies} and \cite{PV}]\label{Teorema-PalisViana}
Consider $f \in \Diff^r(M)$ with $\dim M =3$, a $C^1-$open set
$\cU$ of $\Diff^r(M)$ ($r\geq 1$), a hyperbolic periodic point $p$
of $f$ such that its continuation $p_g$ is well defined for every
$g\in \cU$ and such that:
   \begin{itemize}
   \item[-] The homoclinic class $H(p_g)$ admits a partially hyperbolic splitting of the form $T_{H(p)}M = E^{cs} \oplus E^u$ for every $g\in \cU$.
   \item[-] The subbundle $E^{cs}$ admits no decomposition in non-trivial $Dg-$invariant subbundles which are dominated.
   \item[-] There is a periodic point $q\in H(p_g)$ such that $\det (Dg^{\pi(q)}_q |_{E^{cs}(q)}) > 1$.
   \end{itemize}
   Then, there exists a $C^1$-open and $C^1$-dense subset $\cU_1 \en \cU$ and a $C^r$-residual subset $\cG_{PV}$ of $\cU_1$ such that for every $g\in \cG_{PV}$ one has that $H(p_g)$ intersects the closure of the set of periodic sources of $g$.
   \end{teo}

The conditions of the Theorem are used in
\cite{BDAbundanceTangencies} in order to create robust tangencies
for a hyperbolic set for diffeomorphisms in an $C^1$-open and
dense subset $\cU_1$ of $\cU$. Then, using similar arguments as in
 section 3.7 of \cite{BLY} one creates tangencies associated with
periodic orbits which are sectionally dissipative for $f^{-1}$
which allows to use the results in \cite{PV} to get the
conclusion.


\section{Plaque families and laminations}\label{Section-FamiliasDePlacas}

\subsection{Stable and unstable lamination}\label{SubSection-VariedadEstable}

As for periodic orbits, when a compact invariant subset admits a
dominated splitting with one uniform extremal subbundle, one can
in a sense integrate the subbundle in order to translate the
uniformity of the bundle in a dynamical property in the manifold.
The following result is classical, the standard proof can be found
in \cite{Shub,HPS} (see also \cite{KH} chapter 6).

\begin{teo}[Strong Unstable Manifold Theorem]\label{Teorema-VariedadEstableFuerte}
Let $\Lambda$ be a compact $f$-invariant set which admits a
dominated splitting of the form $T_\Lambda M = E^{cs} \oplus E^u$
where $E^u$ is uniformly expanded. Then, there exists a
\emph{lamination} $\cF^u$ such that:
\begin{itemize}
 \item[-] For every $x\in \Lambda$ the leaf $\cF^u(x)$ through $x$ is an injectively immersed copy of $\RR^{\dim E^u}$ tangent at $x$ to $E^u(x)$.
 \item[-] The leaves of $\cF^u$ form a partition, this is, for $x,y\in \Lambda$ we have that either $\cF^u(x)$ and $\cF^u(y)$ are disjoint or coincide.
 \item[-] There exists $\rho > 0 $ such that points of $\cF^u(x)$ are characterized in the following way:

   $$ y \in \cF^u(x) \Leftrightarrow \lim_{n\to \infty} \frac 1 n \log( d(f^{-n}(x),f^{-n}(y)))< -\rho $$
 \item[-] Leaves vary continuously in the $C^1$-topology: If $x_n \in \Lambda \to x \in \Lambda$ we have that $\cF^u(x_n)$ tends uniformly in compact subsets to $\cF^u(x)$ in the $C^1$-topology.
 \item[-] There exists a neighborhood $U$ of $\Lambda$ such that the leaves also vary continuously in the $C^1$-topology for points in the maximal invariant subset of $U$ for diffeomorphisms $g$ close to $f$.
\end{itemize}
\end{teo}

By \emph{lamination} on a set $K$ we mean a collection of disjoint
$C^1$ injectively immersed manifolds of the same dimension (called
\emph{leaves}) such that there exists a compact metric space
$\Gamma$ such that for every point $x\in K$ there exists a
neighborhood $U$ and a homeomorphism $\varphi: U \cap K \to \Gamma
\times \RR^d$ such that if $L$ is a leaf of the lamination and
$\tilde L$ a connected component of $L\cap U$ then
$\varphi|_{\tilde L}$ is a $C^1$-diffeomorphism to $\{s\} \times
\RR^d$ for some $s \in \Gamma$ (this implies that in $K$ they are
tangent to a continuous subbundle of $T_KM$).

For a lamination $\cF$ on a compact set $K \en M$ we shall always
denote as $\cF(x)$ to the leaf of $\cF$ through $x$. It is worth
remarking that in Theorem \ref{Teorema-VariedadEstableFuerte} the
set laminated by $\cF^u$ need not coincide with $\Lambda$ as it
may be (and it is in various situations) strictly larger.

For hyperbolic sets, this results gives two transversal
laminations which will admit a \emph{local product structure} and
dynamical properties. This allows to obtain the well known
\emph{shadowing lemma} (see \cite{Shub}). In particular, we obtain
the following corollary in quite a direct way:

\begin{prop}\label{Proposition-HiperbolicoImplicaAislado}
Let $\cC$ be a chain-recurrence class of a diffeomorphism $f$
which is hyperbolic. Then, it is isolated and coincides with the
homoclinic class of any of its periodic points. In particular, it
is transitive.
\end{prop}

\subsection{Locally invariant plaque families}\label{SubSection-FliasDePlacas}

In order to search for dynamical or topological consequences of
having a geometric structure invariant under the tangent map, it
is important to try to ``project'' into the manifold the
information we have on the tangent map.

A model of this kind of projection was given in Theorem
\ref{Teorema-VariedadEstableFuerte} where we saw that the dynamics
of the tangent map on uniform bundles project into similar uniform
behavior in invariant submanifolds of the same dimension.

When the invariant bundles are not uniform, we are not able to
obtain such a description (much in the way it is in general not
possible to understand the local behavior of a real-valued
function when the derivative is $1$) but we are able to obtain
certain plaque-families which sometimes help in reducing the
ambient dimension and transforming problems in high dimensional
dynamics into lower dimensional ones.

\begin{teo}[\cite{HPS} Theorem 5.5]\label{Teorema-HPSFliaPlacas}
Let $\Lambda$ be a compact $f$-invariant set endowed with a
dominated splitting of the form $T_\Lambda M =E \oplus F$. Then,
there exists a continuous map $\cW : x \in \Lambda \mapsto \cW_x
\in \emb^1(E(x), M)$ such that:
\begin{itemize}
 \item[-] For every $x\in \Lambda$ we have that $\cW_x(0) =x$ and the image of $\cW_x$ is tangent to $E(x)$ at $x$.
 \item[-] It is \emph{locally invariant}, i.e. There exists $\rho>0$ such that $f(\cW_x(B_\rho(0))) \en \cW_{f(x)} (E(f(x)))$ for every $x\in\Lambda$.
\end{itemize}
\end{teo}

\begin{obs}\label{Remark-HPSparaelfibradodelmedio} In case one has a dominated splitting of the form $T_\Lambda
M =E \oplus F \oplus G$ one can obtain (by applying the previous
theorem to $E\oplus F$ and to $F\oplus G$ with $f^{-1}$) a locally
invariant plaque family tangent to $F$ as well.\finobs
\end{obs}

We will usually (in case no confusion appears) abuse notation and
denote $\cW_x$ to $\cW_x(E(x))$. Also, $\overline{\cW_x}$ will
denote the closure of $\cW_x(E(x))$ which we can assume is the
image of a closed ball of $\RR^{\dim E}$. The proof of the theorem
allows one to obtain a uniform version of this result, in fact,
one obtains that the locally invariant plaque family can also be
chosen continuous with respect to the diffeomorphism in a
neighborhood of $f$ and defined in the maximal invariant subset by
that diffeomorphism in a neighborhood of $\Lambda$ (see \cite{CP}
Lemma 3.7).

\begin{obs}\label{Remark-NoUnicidadFliasPlacas}
Since these locally invariant manifolds are not dynamically
defined they have no uniqueness properties a priori. They may even
have wild intersections between them (see \cite{BuW} for a
construction which is slightly more ``friendly'' which they call
\emph{fake foliations}).  \finobs
\end{obs}

When an invariant plaque family has dynamical properties one can
often recover certain uniqueness properties (see Chapter 5 of
\cite{HPS} or Lemma 2.4 of \cite{Crov-CentralModels}):

\begin{prop}\label{Proposicion-CoherenciaLocal}
Let $\Lambda$ be a compact set admitting a dominated splitting
$T_\Lambda M = E \oplus F$. There exists $\eps>0$ such that if
there exists a plaque family $\{\cW_x\}_{x\in \Lambda}$ tangent to
$E$ verifying that:
\begin{itemize}
\item[-] Every plaque $\cW_x$ has diameter smaller than $\eps$.
\item[-] The plaques verify the following \emph{trapping
condition}:

$$ \forall x \in \Lambda \quad f(\overline{\cW_x}) \en \cW_{f(x)} $$
\end{itemize}
Then, the following properties are verified:
\begin{itemize}
\item[]{\bf (Uniqueness)} Any locally invariant plaque family
$\{\cW'_x\}_{x\in \Lambda}$ tangent to $E$ verifies that for every
$x \in \Lambda$ the intersection $\cW'_x \cap \cW_x$ is open
relative to both plaques.

\item[]{\bf (Coherence)} Given $x,y\in \Lambda$ such that $\cW_x
\cap \cW_y\neq \emptyset$ then we have that the intersection is
open relative to both plaques.

\item[]{\bf (Robust Trapping)} There exists $\delta>0$ such that
if $y \in \cW_x \cap \Lambda$ is at distance smaller than $\delta$
from $x$ then we have that $f(\overline{\cW_y})\en \cW_{f(x)}$.
Moreover, if $\Lambda$ is the maximal invariant set in a
neighborhood $U$, then, there exists $\cU$ a neighborhood of $f$
such that for every $g \in \cU$ the maximal invariant set in $U$
has a plaque family which verifies the same trapping condition.
\end{itemize}
\end{prop}

\dem One can choose a neighborhood $U$ of $\Lambda$ and a
neighborhood $\cU$ of $f$ such that the maximal invariant set in
$U$ will have a dominated splitting for every $g\in \cU$.
Moreover, there will be a cone field $\cE^E$ around $E$ (resp.
$\cE^F$ of $F$) which is invariant for every point in $U$ whose
backward (resp. forward) iterate is also in $U$: i.e. verifies
that $Df_x^{-1} \cE^E(x) \en \Int \cE^E (f^{-1}(x))$ when $x,
f^{-1}(x) \in U$ (resp. $Df_x\cE^F(x) \en \Int \cE^F(f(x))$ when
$x,f(x)\in U$). See Proposition \ref{Proposition-ConeFields} and
Remark \ref{Remark-RobustnessDominatedSplitting}.

Choose $\eps$-small enough so that the expansion in $F$ dominates
the one in $E$ for every pair of points at distance smaller than
$\eps$ (this is trivially verified for every $\eps$ if the
domination is absolute).

Assuming any of the first two properties of the consequence of the
proposition does not hold (uniqueness or coherence), one can find
points $x,z_1, z_2$ in a ball of radius smaller than $\eps$ such
that $x$ can be connected to both $z_1$ and $z_2$ by curves
$\gamma_1$ and $\gamma_2$ contained in some plaque of the plaque
family and such that $z_1$ and $z_2$ can be joined by a curve
$\eta$ (of positive length) which is tangent to $\cE^F$. Moreover,
for every $n\geq 0$ choose $\eta_n$ the curve tangent to $\cE^F$
joining $z_1$ and $z_2$ such that the length of $f^n(\eta_n)$ is
minimal among curves joining $f^n(z_1)$ and $f^n(z_2)$ and whose
preimage by $f^{-n}$ are tangent to $\cE^F$.

Using the trapping condition, one concludes that for every $n>0$
one has that the points $f^n(x), f^n(z_1), f^n(z_2)$ are contained
in a ball of radius $\eps$ around $f^n(x)$.

This implies that $f^n(\eta_n)$ remains always of length smaller
than $\eps$ and since the initial length was at least $\delta>0$
this implies that the length of $\gamma_1$ and $\gamma_2$
decreases exponentially fast, in particular:

$$ \delta \leq \length(\eta_n) \leq (1+\epsilon)^n \|Df^{-n}|_{F(f^n(x))}\| \length (f^n(\eta_n)) \leq$$
$$\leq (1+\epsilon)^n \|Df^{-n}|_{F(f^n(x))}\| (\length (f^n(\gamma_1)
+\length (f^n(\gamma_2)) \leq $$ $$\leq (1 + \epsilon)^{2n}
\|Df^{-n}|_{F(f^n(x))}\| \|Df^n|_{E(x)} \| (\length (f^n(\gamma_1)
+\length (f^n(\gamma_2)) \leq $$ $$\leq const (\lambda
(1+\epsilon)^2)^n \to 0 $$

\noindent which is a contradiction.

To show that trapping holds after perturbation, it is enough to
use the fact that the plaque families vary continuously and the
trapping condition is $C^0$-open.

\lqqd

This argument uses strongly the fact that plaques are sufficiently
small. There are two reasons for this:

\begin{itemize}
\item[-] It allows to control the domination between the points
involved. \item[-] It allows to control the geometry of the curves
joining the points in different plaques and sharing a point in
their plaques. \end{itemize}

In section \ref{Section-QIyArgumentoBrin} we will review an
argument of Brin (\cite{Brin}) which allows to obtain tangent
foliations under the existence of a partially hyperbolic
splitting. To solve  the first problem, he uses absolute
domination, and for the second one, he introduces the concept of
quasi-isometric foliations which allows him to obtain the desired
geometry for comparing distances and lengths.

\subsection{Holonomy and local manifolds}\label{SubSection-HolonomiaFliasLocales}


When there exists an invariant lamination or foliation $\cF$
tangent to certain bundle $E$ on some invariant set $\Lambda$ we
will denote the \emph{local leaves} through a point $x$ as
$\cF_{loc} (x)$. By this we mean that $\cF_{loc}(x)$ is the
connected component of the leaf $\cF(x)$ containing $x$ in a
neighborhood of $x$. We remark that this notion is of course not
strictly well defined but when we mention local leaves we will
state which are the referred neighborhoods. In some situations, we
will use other notations such as $W^{\sigma}_{loc}(x)$ or
$\cW^{\sigma}_{x,loc}$. This notations hold also for locally
invariant plaque families which we shall sometimes give similar
notation.

When we have two transverse laminations, or even only one and
transverse local leaves we can define the \emph{holonomy} between
the transversals (which in some sense generalizes the holonomy of
foliations, c.f. Chapter \ref{Capitulo-FoliacionesEHipParcial}).

Consider a compact set $\Lambda$ admitting a lamination $\cF$
tangent to a subbundle $F$ of $T_\Lambda M$, denote $\tilde
\Lambda$ to the union of leaves of $\cF$.

Given a plaque family $\{\cW_x\}_{x\in \Lambda}$ tangent to a
bundle $E$ such that $E\oplus F =T_\Lambda M$ we can define the
following set of maps for $x,y \in \Lambda$ in the same leaf of
$\cF$:

$$ \Pi^{\cF}_{x,y} : U \cap \tilde \Lambda \en \cW_x \to \cW_y $$

\noindent given by the intersection between the local leaves of
$\cF$ intersecting $\cW_x$ with $\cW_y$. Notice that the domain
$U\cap \tilde \Lambda$ of $\Pi_{x,y}$ can be chosen in order to be
open in $\cW_x \cap \tilde \Lambda$ and contain a neighborhood of
$x$ there.

When $x,y$ are sufficiently close in the leaf $\cF(x)$, the
domains in the transversal can be chosen arbitrarily large, since
the transversals are very close.

\subsection{Control of uniformity of certain bundles}\label{SubSection-UniformidaddeFibrados}

Sometimes one can deduce that certain extremal bundles are
uniform. In dimension $2$ this follows from a result from
\cite{PujSamAnnals1}:

\begin{teo}[\cite{PujSamAnnals1}, \cite{ABCD} Theorem 2]\label{Teorema-PujSamC1}
There exist a residual subset $\cG_{PS} \en \Diff^1(M^2)$ where
$M^2$ is a surface such that if $f\in \cG_{PS}$ and $\Lambda$ is a
chain recurrence class admitting a dominated splitting, then,
$\Lambda$ is hyperbolic.
\end{teo}

The proof of this result uses approximation by
$C^2$-diffeomorphisms. At the moment, it is not completely
understood the importance of the fact that bundles are extremal
and extending this results to higher dimensions as well as for
non-extremal bundles represents a main challenge (see
\cite{PujSamIHP,CP,CSY} for some progress in that direction).

Another important result we will use relates the existence of
hyperbolic invariant measures for $C^1$-diffeomorphisms whose
support admits a dominated splitting separating positive and
negative Lyapunov exponents. This results extends a well known
result of Katok (see \cite{KH} Supplement S) asserting that
hyperbolic measures of $C^2$-diffeomorphisms are contained in the
support of a homoclinic class. The cost for doing this is
requiring a dominated splitting separating the Lyapunov exponents
of the measures (a necessary hypothesis, see \cite{BCS}):

\begin{teo}[\cite{ABC} and \cite{Crov-CentralModels}, Proposition 1.4]\label{Teorema-PesinC1}
Let $\mu$ be an ergodic hyperbolic measure of a
$C^1$-diffeomorphism $f$ (that is, all the Lyapunov exponents are
different from zero) such that $\supp(\mu)$ admits a dominated
splitting $T_{\supp(\mu)}M =E \oplus F$ such that the Lyapunov
exponents on $E$ are negative and in $F$ positive. Then, the
support of $\mu$ is contained in a homoclinic class containing
periodic orbits of stable index $\dim E$.
\end{teo}

The proof of this theorem follows from careful application of the
existence of locally invariant plaque families as well as ideas in
the vein of Lemma \ref{Lemma-SiExponentesNegativosPozo} (see also
\cite{pliss}).

\subsection{Central models and Lyapunov exponents}\label{SubSection-ModelosCentrales}

We will present the tool of central models first introduced in
\cite{Crov-Birth} and developed in \cite{Crov-CentralModels} which
allows to treat the case where there is no knowledge on the
Lyapunov exponents along a certain $Df$-invariant subbundle of
dimension $1$. The presentation will be incomplete and restricted
to the uses we will make of this tool. We strongly recommend the
reading of \cite{Crov-Hab} Chapter 9 or \cite{Crov-CentralModels}
section 2 if the reader is interested in understanding this tool.

Consider a compact $f$-invariant set $\Lambda$ which is
chain-transitive and we will assume that $\Lambda$ admits a
dominated splitting $T_\Lambda M = E_1 \oplus E^c \oplus E_3$ with
$\dim E^c =1$.

Consider a locally invariant plaque family $\{\cW^c_x\}_{x\in
\Lambda}$ tangent to $E^c$. Recall that each $\cW^c_x$ is an
embedding of $E^c(x)$ in $M$.

By local invariance, there exists $\rho>0$ such that
$f(\cW^c_x([-\rho,\rho])) \en \cW^c_{f(x)} (\RR)$, where we are
identifying $E^c(x)$ with $\RR$. Without loss of generality, we
can take $\rho =1$.

When $Df$-preserves some continuous orientation on $E^c$ (which in
particular implies that $E^c$ is orientable) this allows us to
define two maps:

$$ \hat f_1 : \Lambda \times [0,1] \to \Lambda \times [0, +\infty) $$
$$  \hat f_2 : \Lambda \times [-1,0] \to \Lambda \times (-\infty, 0] $$

\noindent induced by the way $f$ acts on the locally invariant
plaques.

When $Df$-does not preserve any continuous orientation on $E^c$,
(in particular when $E^c$ is not orientable) one can consider the
double covering $\hat \Lambda$ of $\Lambda$ (on which the dynamics
will still be chain-transitive) and in a similar way define one
map (see \cite{Crov-CentralModels} section 2):

$$ \hat f: \hat \Lambda \times [0,1] \to \hat \Lambda \times [0,+\infty) $$

This motivates the study of continuous skew-products of the form
(called \emph{central models}):

$$ \hat f: K \times [0,1] \to K \times [0,+\infty) $$

$$ \hat f (x,t) = (f_1(x), f_2(x,t)) $$

\noindent where $f_1 : K \to K$ is chain-transitive, $f_2(x,0)=0$
and $\hat f$ is a local homeomorphism in a neighborhood of
$K\times \{0\}$.

For this kind of dynamics, in \cite{Crov-CentralModels} the
following classification was proven:

\begin{prop}[Central Models \cite{Crov-CentralModels} Proposition 2.2]\label{Proposition-CentralModels}
For a central model $\hat f: K \times [0,1] \to K \times
[0,+\infty)$ one of the following possibilities holds:
\begin{itemize}
\item[-]The chain-stable and the chain-stable set of $K\times
\{0\}$ are non-trivial. In this case there is a segment $\{x\}
\times [0,\delta]$ which is chain-recurrent. \item[-] The
chain-stable set contains a neighborhood of $K\times \{0\}$ and
the chain-unstable set is trivial. \item[-] The chain-unstable set
contains a neighborhood of $K\times \{0\}$ and the chain-stable
set is trivial. \item[-] Both the chain-stable and the
chain-unstable set of $K\times \{0\}$ are trivial.
\end{itemize}
\end{prop}

As a consequence, for partially hyperbolic dynamics we have (at
least one of) the following types of central dynamics:

\begin{itemize}
\item[-]{\bf Type (R)} For every neighborhood $U$ of $\Lambda$,
there exists a curve $\gamma$ tangent to $E^c$ at a point of
$\Lambda$ such that $\gamma$ is contained in a compact, invariant,
chain-transitive set in $U$. \item[-]{\bf Type (N)} There are
arbitrarily small neighborhoods $U_k$ of the $0$ section of $E^c$
such that $f(\cW^c_x (\overline{U_k})) \en \cW^c_{f(x)}(U_k)$
(which we call \emph{trapping strips} for $f$) and there are
arbitrarily small trapping strips for $f^{-1}$. \item[-]{\bf Type
(H)} There are arbitrarily small trapping strips for $f$ (case
($H_S$)) or for $f^{-1}$ (case ($H_U$)) and the trapping strips
belong to the chain-stable set of $\Lambda$ (case ($H_S$)) or the
chain-unstable set of $\Lambda$ (case ($H_U$)). \item[-]{\bf Type
(P)} This is only possible in the orientable case and corresponds
to the following subtypes: ($P_{SN}$), ($P_{UN}$) and ($P_{SU}$)
and corresponds to the case where there is a mixed behavior with
respect to the types defined above.
\end{itemize}

In \cite{Crov-CentralModels} it is proved that these types are
well defined and more properties are studied.

\subsection{Blenders}\label{SubSection-Blenders}

Blenders represent one of the main tools of differentiable
dynamics, in particular when searching to prove certain robust
properties of diffeomorphisms. They were introduced in
\cite{Diaz-Blenders} and \cite{BDAnnals}. See \cite{BDV} chapter 6
for a nice introduction to these sets, we will only present some
properties which we use later. An explicit construction of these
sets can be found in Appendix \ref{Apendice-BCGP}.

We shall now present $cu$-blenders by its properties: A
$cu$-\emph{blender} $K$ for a diffeomorphism $f: M \to M$ of a
$3$-dimensional manifold\footnote{Of course they can be defined in
more generality, see \cite{BDAnnals,BDAbundanceTangencies}.} is a
compact $f$-invariant hyperbolic set with splitting $T_K M =
E^{ss} \oplus E^s \oplus E^{u}$ such that the following properties
are verified:

\begin{itemize}
\item[-] $K$ is the maximal invariant subset in a neighborhood
$U$. \item[-] There exists a cone-field $\cE^{ss}$ around $E^{ss}$
defined in all $U$ which is invariant under $Df^{-1}$. \item[-]
There exists a compact region $B$ with non-empty interior (which
is called \emph{activating region}) such that every curve
contained in $U$, tangent to $\cE^{ss}$ with length larger than
$\delta$ and intersecting $B$ verifies that it intersects the
stable manifold of a point of $K$. \item[-] There exists an open
neighborhood $\cU$ of $f$ such that for every $g$ in $\cU$ the
properties above are verified for the same cone field, the same
set $B$ and for $K_g$ the maximal invariant set of $U$.
\end{itemize}

For more properties and construction of $cu$-blenders, see
\cite{BDV} chapters 6 and \cite{BDAnnals}. There one can see a
proof of the following:

\begin{prop}[\cite{BDAnnals} Lemma 1.9, \cite{BDV} Lemma 6.2]\label{PropBlender}
If the stable manifold of a periodic point $p\in M$ of stable
index $1$ contains an arc $\gamma$ tangent to $\cE^{ss}$ and
intersecting the activating region of a $cu$-blender $K$, then,
$W^u(p) \en \overline{W^u(q)}$ for every $q$ periodic point in
$K$.
\end{prop}


\subsection{Higher regularity and SRB measures}\label{SubSection-BonattiViana}

We shall briefly review some of the results from \cite{BV} (see
also \cite{VY} for recent advances on this direction) that
guaranty the existence of a unique SRB measure in certain
partially hyperbolic sets whenever there are some properties
verified by the exponents in the center stable direction and $f$
is sufficiently regular (i.e. $C^2$ is enough).

Consider $f: M \to M$ a $C^2$-diffeomorphism of a compact manifold
such that it contains an open set $U$ such that $f(\overline U)
\en U$. We denote:

$$\Lambda= \bigcap_{n\geq 0} f^n(\overline U)$$

\noindent which is a (not-necessarily transitive) topological
attractor.

We shall assume that $\Lambda$ admits a partially hyperbolic
splitting of the form $T_\Lambda M =E^{cs} \oplus E^u$ where $E^u$
is uniformly expanding and $E^{cs}$ is dominated by $E^u$.

Since $\Lambda$ is a topological attractor, we get that $\Lambda$
is saturated by unstable manifolds (see Proposition
\ref{ProposicionAtractoresSaturadosPorInestables}).

To obtain SRB measures for this type of attractors one considers
the push-forward by the iterates of the diffeomorphism of Lebesgue
measure and by controlling the distortion (here the
$C^2$-hypothesis becomes crucial) one can see that the limiting
measure is absolutely continuous with respect to the unstable
foliation (for precise definitions see \cite{BDV} chapter 11).
After this is done, the fact that for $C^2$-diffeomorphisms the
center-stable leaves are also absolutely continuous, one shows
that the limit measures are SRB measures and that their basins
covers a full Lebesgue measure of the basin. To obtain this
results, further hypothesis are required in \cite{BV} which we
pass to review.

We define

$$ \lambda^{cs}(x) = \limsup_{n\to +\infty} \frac 1 n \log \|D_x f^n|_{E^{cs}(x)} \| $$

\noindent which resembles the Lyapunov exponent (only that $x$
needs not be a Lyapunov regular point).

We obtain the following result:

\begin{teo}[Bonatti-Viana \cite{BV} Theorem A]\label{Teorema-BonattiVianaA}
Let $f: M \to M$ be a $C^2$-diffeomorphism such that it admits an
open set $U$ verifying $f(\overline{U}) \en U$ such that $\Lambda
=\bigcap_n f^n(U)$, its maximal invariant set is partially
hyperbolic with splitting $T_\Lambda M= E^{cs} \oplus E^u$. Assume
moreover that for every $D$ disc contained in $\cF^u$ there is a
positive Lebesgue measure of points $x$ such that $\lambda^{cs}(x)
<0$. Then, there exists finitely many SRB measures $\mu_1, \ldots,
\mu_k$ such that $\bigcup_i \Bas(\mu_i)$ has total Lebesgue
measure inside $\Bas(\Lambda)$.
\end{teo}

Under certain assumptions, one can see that there is a unique SRB
measure. We shall state the following theorem which has slightly
more general hypothesis but for which the same proof as in
\cite{BV} works (see also \cite{VY} for a further development of
these results):

\begin{teo}[Bonatti-Viana \cite{BV} Theorem B]\label{Teorema-BonattiVianaB}
Assume that $f$ and $\Lambda$ satisfy the hypothesis of Theorem
\ref{Teorema-BonattiVianaA} and that moreover there is a unique
minimal set of $\cF^u$ inside $\Lambda$, then, $f$ admits a unique
SRB measure in $\Lambda$ whose statistical basin coincides with
the topological one modulo a zero Lebesgue measure set.
\end{teo}

The hypothesis required in \cite{BV} is that $\cF^u$ is minimal
inside $\Lambda$. However, it is not hard to see how the proof of
\cite{BV} works for the hypotheis stated above: See the first
paragraph of section 5 in \cite{BV}, consider the unique minimal
set $\tilde \Lambda$ of the unstable foliation: we get that there
is only one \emph{accessibility class} there as needed for their
Theorem B.


\section{Normal hyperbolicity and dynamical coherence}\label{Section-HIRSHPUGHSHUB}

Consider a lamination $\cF$ in a compact set $\Lambda$ and let
$f:M\to M$ be a $C^1$-diffeomorphism preserving $\cF$. We will say
that $\cF$ is \emph{normally hyperbolic} if there exists a
splitting of $T_\Lambda M= E^s \oplus T\cF \oplus E^u$ as a
$Df$-invariant sum verifying that the decomposition is partially
hyperbolic (in particular, $E^s$ or $E^u$ can be trivial). If the
domination is of absolute nature, we say that $\cF$ is
\emph{absolutely normally hyperbolic}. See \cite{HPS} for the
classical reference and \cite{Berger} for recent results and some
modern proofs of the results (and extensions to general
laminations and endomorphisms).

When the lamination $\cF$ covers the whole manifold, we say that
it is a \emph{foliation} (this corresponds with a $C^0$-foliation
with $C^1$-leaves and tangent to a continuous distribution in the
literature). See Chapter \ref{Capitulo-FoliacionesEHipParcial}.

\subsection{Leaf conjugacy}\label{SubSection-ConjugacionPorHojas}

Given a lamination $\cF$ which is invariant under a diffeomorphism
$f$ one wishes to understand which conditions guaranty the fact
that for perturbations $g$ of $f$ there will still be a foliation
which is $g$-invariant. As hyperbolicity gives a sufficient
condition for structural stability, normal hyperbolicity appears
as a natural requirement when one searches for persistence of
invariant laminations\footnote{Though in this case, the issue of
being a necessary condition is quite more subtle \cite{Berger}.}.

In some situations, one obtains something much stronger than
persistence of an invariant lamination (notice that for a
0-dimensional foliation-by points- the following notion coincides
with the usual conjugacy). For a lamination $\cF$ we denote as
$K_\cF$ to the (compact) set which is the union of the leaves of
$\cF$.

\begin{defi}[Leaf conjugacy]\label{DefinicionLeafConjugacy}
Given $f, g: M \to M$ be $C^1$-diffeomorphisms such that there are
laminations $\cF_f$ and $\cF_g$ invariant under $f$ and $g$
respectively. We say that $(f,\cF_f)$ and $(g,\cF_g)$ are
\emph{leaf conjugate} if there exists a homeomorphism $h:
K_{\cF_f} \to K_{\cF_g}$ such that:
\begin{itemize}
\item[-] For every $x \in M$, $h(\cF_f(x)) = \cF_g(h(x))$.
\item[-] For every $x\in M$ we have that

$$ h (\cF_f(f(x))) = \cF_g (g\circ h(x)) $$
\end{itemize}

If a $C^1$-diffeomorphism $f$ leaves a lamination $\cF$ invariant
we say that the foliation $\cF$ is \emph{structurally stable} if
there exists a neighborhood $\cU$ of $f$ such that for $g \in \cU$
the diffeomorphism $g$ admits a $g$-invariant foliation $\cF_g$
such that the pairs $(f,\cF)$ and $(g,\cF_g)$ are leaf conjugate.
\finobs
\end{defi}

The classical result of \cite{HPS} asserts that normal
hyperbolicity along with a technical condition called \emph{plaque
expansivity} is enough to guarantee structural stability of a
lamination:

\begin{teo}[\cite{HPS} Chapter 7 and \cite{Berger} Remark 2.2]\label{Teorema-HPSEstabilidadEstr}
Let $f$ be a $C^1$-diffeomorphism leaving invariant a foliation
$\cF$ which is normally hyperbolic and plaque expansive we have
that the foliation is structurally stable.
\end{teo}

We shall not give a  definition of plaque-expansivity (we refer
the reader to \cite{HPS,Berger}) but we mention that it is not
known if it is a necessary hypothesis and all known normally
hyperbolic foliations are either known to be structurally stable
or at least suspected.

We do however state the following result which ensures
plaque-expansivity and is useful to treat many important examples:

\begin{prop}\label{Proposicion-C1implicaPlaqueExpansive}
If a normally hyperbolic foliation $\cF$ is of class  $C^1$ (this
means that the change of charts given by Proposition
\ref{Proposicion-FoliacionAtlas} are of class $C^1$) then it is
plaque-expansive.
\end{prop}

This extends also to general laminations where the concept of
being $C^1$ is harder to define. We will use this fact later in
this thesis.

In general, checking plaque-expansiveness is hard and this makes
it an undesirable hypothesis for leaf conjugacy.

\subsection{Dynamical coherence}\label{SubSeccion-CoherenciaDefinicion}

One sometimes wishes to consider the inverse problem. We have seen
in Theorem \ref{Teorema-VariedadEstableFuerte} that if a
diffeomorphism $f: M \to M$ is partially hyperbolic with splitting
$TM = E^{cs} \oplus E^u$ then there exists a (unique) foliation
$\cF^u$ tangent to $E^u$ which we call the \emph{unstable
foliation} of $f$. In general, it may happen that there is no
foliation tangent to $E^{cs}$ (this was remarked by Wilkinson in
\cite{Wilkinson} using an example of Smale \cite{SmaleBulletin},
this is extended in section 3 of \cite{BuW2}).

\begin{defi}[Dynamical coherence]\label{Definicion-Coherencia}
Let $f: M \to M$ be a partially hyperbolic diffeomorphism with
splitting $TM = E^{cs} \oplus E^u$. We say that $f$ is
\emph{dynamically coherent} if there exists a foliation $\cF^{cs}$
everywhere tangent to $E^{cs}$ which is $f$-invariant in the sense
that $f(\cF^{cs}(x)) = \cF^{cs}(f(x))$. When $f$ is strongly
partially hyperbolic with splitting $TM = E^s \oplus E^c \oplus
E^u$ we say that it is dynamically coherent if there exists
$f$-invariant foliations tangent to both $E^s\oplus E^c$ and to
$E^c \oplus E^u$. \finobs
\end{defi}

\begin{obs}[Central Direction]\label{Remark-FoliacionCentral}
When a strong partially hyperbolic diffeomorphism is dynamically
coherent one can intersect the foliations $\cF^{cs}$ and
$\cF^{cu}$ tangent to $E^s \oplus E^c$ and $E^c \oplus E^u$
respectively and obtain a foliation $\cF^c$ tangent to $E^c$.
Moreover, one can show (see Proposition 2.4 of \cite{BuW2}) that
the foliations $\cF^c$ and $\cF^s$ (resp. $\cF^u$) subfoliate the
leaves of $\cF^{cs}$ (resp. $\cF^{cu}$). \finobs
\end{obs}

\begin{obs}[Unique integrability]\label{Remark-UnicidadCoherencia}
We have not made assumptions in the definition of dynamical
coherence about the uniqueness of the $f$-invariant foliation
tangent to $E^{cs}$. There are many ways to require uniqueness:
\begin{itemize}
\item[-] One can ask for $\cF^{cs}$ to be the unique $f$-invariant
foliation tangent to $E^{cs}$. If there exists $n>0$ such that
there exists a unique $f^n$-invariant foliation, then $f$ is
dynamically coherent and with a unique $f$-invariant foliation.
Dynamical coherence in principle does not follow from the
existence of an $f^n$-invariant foliation tangent to $E^{cs}$.
\item[-] One can ask for $\cF^{cs}$ to be the unique foliation
tangent to $E^{cs}$ which is stronger than the previous
requirement. \item[-] One can ask for the following much stronger
statement: Any $C^1$-curve everywhere tangent to $E^{cs}$ is
contained in a leaf of $\cF^{cs}$.
\end{itemize}
These (and more) types of uniqueness properties are discussed
further in section 2 of \cite{BuW2}. See
\cite{PughShubWilkinson,BF} for examples of foliations which
satisfy weak forms of uniqueness.\finobs
\end{obs}

Notice that if $f$, a partially hyperbolic diffeomorphism is
dynamically coherent, then the $f$-invariant foliations
($\cF^{cs}$, $\cF^{cu}$ and $\cF^c$) are automatically normally
hyperbolic.

Notice that again, the lack of knowledge (in general) of whether
the foliations are plaque-expansive does not allow to know if in
general the foliation must be structurally stable, in particular,
the following is an open question:

\begin{quest} Is dynamical coherence an open property?
\end{quest}

\subsection{Classification of transitive $3$-dimensional strong partially hyperbolic diffeomorphisms}
Another main problem in dynamics is to consider a class of systems
and try to classify their possible dynamics. For partially
hyperbolic diffeomorphisms there is the following conjecture posed
by Pujals (see \cite{BonWilk}):

\begin{conj}[Pujals]
Let $f: M^3 \to M^3$ be a transitive strong partially hyperbolic
diffeomorphism of a $3$-dimensional manifold, then $f$ is leaf
conjugate to one of the following models:
\begin{itemize}
\item[-] Finite lifts of a skew product over an Anosov map of
$\TT^2$. \item[-] Finite lifts of time one maps of Anosov flows.
\item[-] Anosov diffeomorphisms on $\TT^3$.
\end{itemize}
\end{conj}

In \cite{BonWilk} this conjecture is treated and some positive
results are obtained without any assumptions on the topology of
the manifold.

Hammerlindl (\cite{Hammerlindl,HNil}) has made some important
partial progress to this conjecture by assuming that the manifold
is a nilmanifold and the partial hyperbolicity admits absolute
domination. The work done in the present thesis allows to
eliminate the need for absolute domination, but it seems that we
still lack of tools to attack the complete conjecture. With A.
Hammerlindl we plan to use several of the techniques used in this
thesis in order to prove Pujals' conjecture for $3$-manifolds with
fundamental group of polynominal growth (see \cite{HP}).

Notice also that the example of \cite{HHU} which is not transitive
does not belong to any of the classes, so that the hypothesis of
transitivity cannot be removed.

\subsection{Accessibility}\label{SubSection-Accessibility}

Let $f: M \to M$ be a strong partially hyperbolic diffeomorphism
with splitting $TM = E^s \oplus E^c \oplus E^u$. As stated in
Theorem \ref{Teorema-VariedadEstableFuerte} there exist foliations
$\cF^s$ and $\cF^u$ tangent to the bundles $E^s$ and $E^u$
respectively.

An important notion introduced by Pugh and Shub in the mid 90's is
the concept of \emph{accessibility} on which their celebrated
conjectures on abundance of ergodicity is based (see
\cite{PughShubConjetura} and also Chapter 8 of \cite{BDV}).

We define the \emph{accessibility class} of a point $x\in M$ as
the set of points $y\in M$ such that there exists an
$su$-\emph{path} from $x$ to $y$. An $su$-path is a concatenation
of finitely many $C^1$-paths alternatively tangent to $E^s$ of
$E^u$.

We say that a diffeomorphism is \emph{accessible} if there is a
unique accessibility class. There has been lots of work devoted to
the understanding of accessibility and its relationship with
ergodicity, but from the point of view of abundance, we know the
following results:

\begin{teo}[\cite{DolgopyatWilkinson}]\label{Teorema-AccesibilidadDensa}
There is a $C^1$-open and dense subset $\cA$ of partially
hyperbolic diffeomorphisms such that for every $f\in \cA$ we have
that $f$ is accessible.
\end{teo}

Moreover, it is proved in \cite{HHUAcc}, that if $\dim E^c =1$ the
set of accessible partially hyperbolic diffeomorphisms forms a
$C^1$-open and $C^\infty$-dense set.

We refer the reader to \cite{BuW} and \cite{HHUInventiones} for
proofs of ergodicity by using accessibility and certain technical
conditions we will not discuss.

\section{Integer $3\times 3$ matrices}

We will denote as $GL(d,\ZZ)$ to the group of invertible $d\times
d$ matrices with integer coefficients. If one interprets this as
having invertibility in the group of integer matrices, it is
immediate that the determinant must be of modulus $1$, but since
there are different uses of this notation in the literature, we
will make explicit mention to this when used.

\subsection{Hyperbolic matrices}\label{SubSection-MatricesHiperbolicas}

Let $A \in GL(3,\ZZ)$ with determinant of modulus $1$ and no
eigenvalues of modulus $1$.


Since $A$ is hyperbolic and the product of eigenvalues is one, we
get that $A$ must have one or two eigenvalues with modulus smaller
than $1$. We say that $A$ has \emph{stable dimension} $1$ or $2$
depending on how many eigenvalues of modulus smaller than one it
has.

We call \emph{stable eigenvalues} (resp. \emph{unstable
eigenvalues}) to the eigenvalues of modulus smaller than one
(resp. larger than one). The subspace $E^s_A=W^s(0,A)$ (resp
$E^u_A=W^u(0,A)$) corresponds to the eigenspace associated to the
stable (resp. unstable) eigenvalues.

We shall review some properties of linear Anosov automorphisms on
$\TT^3$.

We say that a matrix $A \in GL(3,\ZZ)$ (with determinant of
modulus $1$) is \emph{irreducible} if and only if its
characteristic polynomial is irreducible in the field $\QQ$. This
is equivalent to stating that the characteristic polynomial has no
rational roots. It is not hard to prove:

\begin{prop}\label{PropAnosovSonIrreducibles}
Every hyperbolic matrix $A \in GL(3,\ZZ)$ with determinant of
modulus $1$ is irreducible. Moreover, it cannot have an invariant
linear two-dimensional torus.
\end{prop}

\dem Assume that a matrix $A$ is not irreducible, this means that
$A$ has one eigenvalue in $\QQ$.

Notice that the characteristic polynomial of $A$ has the form
$-\lambda^3 + a\lambda^2+ b\lambda \pm 1$. By the rational root
theorem  (see \cite{Hungerford}), if there is a rational root, it
must be $\pm 1$ which is impossible if $A$ is hyperbolic.

Every linear Anosov automorphism is transitive. Let $T$ be a
linear two-dimensional torus which is invariant under $A$. Since
the tangent space of $T$ must also be invariant, we get that it
must be everywhere tangent to an eigenspace of $A$. Since we have
only $3$ eigenvalues, this implies that either $T$ is attracting
or repelling, contradicting transitivity.

\lqqd

We can obtain further properties of hyperbolic matrices acting in
$\TT^3$:

\begin{lema}\label{RemarkValoresPropiosAnosov}
Let $A \in SL(3,\ZZ)$ be a hyperbolic matrix. Then, the
eigenvalues are simple and irrational. Moreover, if there is a
pair of complex conjugate eigenvalues they must be of irrational
angle. \end{lema}

\dem By the previous proposition, we have that the characteristic
polynomial of $A$ is irreducible as a polynomial with rational
coefficients.

It is a classic result in Galois´ theory that in a field of
characteristic zero, irreducible polynomials have simple roots
(see \cite{Hungerford} Definition V.3.10 and the Remark that
follows): In fact, since $\QQ[x]$ is a principal ideals domain, if
a polynomial has double roots then it can be factorized by its
derivative which has strictly smaller degree contradicting
irreducibility.

This also implies that if the roots are complex, they must have
irrational angle since otherwise, by iterating $A$ we would obtain
an irreducible polynomial of degree $3$ and non-simple roots
(namely, the power of the complex conjugate roots which makes them
equal).

\lqqd

When $A$ has two different stable eigenvalues $|\lambda_1| <
|\lambda_2| < 1$ (resp. unstable eigenvalues
$|\lambda_1|>|\lambda_2|>1$) we call \emph{strong stable manifold
of} $A$ (resp. \emph{strong unstable manifold of} $A$) to the
eigenline of $\lambda_1$ which we denote as $E^{ss}_A$ (resp.
$E^{uu}_A$).

\begin{obs}\label{RemarkSubespaciosInvariantesAnosov}
For every $A \in SL(3,\ZZ)$ hyperbolic, we know exactly which are
the invariant planes of $A$. If $A$ has complex eigenvalues, then,
the only invariant plane is the eigenspace associated to that pair
of complex conjugate eigenvalues. If $A$ has $3$ different real
eigenvalues then there are $3$ different invariant planes, one for
each pair of eigenvalues. All these planes are totally irrational
(i.e. their projection to $\TT^3$ is simply connected and dense).
\finobs
\end{obs}

\subsection{Non-hyperbolic partially hyperbolic matrices}

We prove the following result which plays the role of Lemma
\ref{RemarkValoresPropiosAnosov} in the non-Anosov partially
hyperbolic case.

\begin{lema}\label{LemaPlanosInvariantes}
Let $A$ be a matrix in $GL(3,\ZZ)$ with eigenvalues $\lambda^s,
\lambda^c, \lambda^u$ verifying $0<|\lambda^s| < |\lambda^c=1| <
|\lambda^u| = |\lambda^s|^{-1}$. Let $E^s_\ast, E^c_\ast,
E^u_\ast$ be the eigenspaces associated to $\lambda^s, \lambda^c$
and $\lambda^u$ respectively. We have that:
\begin{itemize}
\item[-] $E^c_\ast$ projects by $p$ into a closed circle where
$p:\RR^3 \to \TT^3$ is the covering projection. \item[-] The
eigenlines $E^{s}_\ast$ and $E^u_\ast$ project by $p$ into
immersed lines whose closure coincide with a two dimensional
linear torus.
\end{itemize}
\end{lema}

\dem We can work in the vector field $\QQ^3$ over $\QQ$ where
$f_\ast$ is well defined since it has integer entries.

Since $1$ is an eigenvalue of $f_\ast$ and is rational, we obtain
that there is an eigenvector of $1$ in $\QQ^3$. Thus, the
$\RR$-generated subspace (now in $\RR^3$) projects under $p$ into
a circle.

Since $1$ is a simple eigenvalue for $f_\ast$, there is a rational
canonical form for $f_\ast$ which implies the existence of
two-dimensional $\QQ$-subspace of $\QQ^3$ which is invariant by
$f_\ast$ and corresponds to the other two eigenvalues (see for
example Theorem VII.4.2 of \cite{Hungerford}).

This plane (as a $2$-dimensional $\RR$-subspace of $\RR^3$) must
project by $p$ into a torus since it is generated by two linearly
independent rational vectors. This torus is disjoint from the
circle corresponding to the eigenvalue $1$ and coincides with the
subspace generated by $E^s_\ast$ and $E^u_\ast$.

On the other hand, the lines generated by $E^s_\ast$ and
$E^u_\ast$ cannot project into circles in the torus since that
would imply they have rational eigenvalues which is not possible,
this implies that the closure of their projection is the whole
torus.

\lqqd


\chapter{Semiconjugacies and localization of chain-recurrence classes}\label{Capitulo-Semiconjugaciones}

The purpose of this chapter is to present Proposition
\ref{ProposicionMecanismo} which plays an important role in this
thesis. It gives conditions under which chain-recurrence classes
accumulating a given one to be contained in lower dimensional
normally hyperbolic submanifolds. We profit to introduce some
notions on semiconjugacies and decompositions of spaces in Section
\ref{Section-Descomposiciones} and to state and prove a classical
result on topological stability of hyperbolic sets which will be
useful to then use Proposition \ref{ProposicionMecanismo}.

\section{Fibers, monotone maps and decompositions of
manifolds}\label{Section-Descomposiciones}

Consider two homeomorphisms $f: X \to X$ and $g: Y \to Y$. We say
that $f$ is \emph{semiconjugated} to $g$ (or that $g$ is a
\emph{factor} of $f$) if there exists a continuous map $h: X \to
Y$ such that

$$  h \circ f  = g \circ h $$

Semiconjugacies will play an important role in this text, that is
why we shall make some effort in understanding certain continuous
maps.

\begin{obs}\label{Remark-PropiedadesSemiconjugacion}
Semiconjugacies preserve some dynamical properties. For example,
if $h: X \to Y$ semiconjugates $f:X \to X$ and $g: Y\to Y$ then we
have that: \bi \item[-] If $x \in X$ is a periodic point for $f$,
then $h(x)$ is periodic for $g$. \item[-] If $x\in X$ is recurrent
(resp. non-wandering) for $f$, then $h(x)$ is recurrent (resp.
nonwandering) for $g$. \item[-] If $x_2$ is in the stable set
(resp. unstable set) of $x_1$ for $f$, then $h(x_2)$ is in the
stable set (resp. unstable set) of $h(x_1)$ for $g$. \ei \finobs
\end{obs}

Let $h : X \to Y$ be a continuous map and $y \in Y$, we call
$h^{-1}(\{y\})$ the \emph{fiber} of $y$ by $h$.

Sometimes, the topology of the fiber gives us information about
the map $h$. We say that $h$ is a \emph{monotone map} if all the
fibers are compact and connected. In general we will work with
$X=M$ a topological manifold, in that case we will require a
stronger property and say that $h$ is a \emph{cellular map} if the
fiber of every point is a \emph{cellular set} (i.e. decreasing
intersection of topological balls).

Every time we have a continuous and surjective map $h: X \to Y$ we
can think $Y$ as $X /_\sim$ where the equivalence classes are
given by fibers of $h$.

Special interest is payed to \emph{cellular decompositions} of
manifolds (a partition of a manifold $M$ such that the quotient
map is a cellular map) since these quotient spaces are what is
known as \emph{generalized manifolds} (see \cite{Daverman}).

To be more precise, we say that an equivalence relation $\sim$ in
a manifold $M$ is a \emph{cellular decomposition} if the following
properties are verified:

\bi \item[-] If we denote by $A_x$ to the equivalence class of a
point $x$ we have that the sets $A_x$ are cellular for every $x\in
M$. \item[-] The decomposition is \emph{upper semicontinuous} in
the sense that if $x_n \to x$ then we have that $\limsup A_{x_n}
\en A_x$. \ei

When we have a cellular decomposition of a manifold $M$, we can
define a quotient map $\pi: M \to M/_\sim$ and we give to
$M/_\sim$ the quotient topology. We have that (see \cite{Daverman}
Proposition I.2.2) that:

\begin{prop}\label{Proposition-CocienteCelularEsMetrico}
The topological space $M/\sim$ is metrizable.
\end{prop}

Also, we can define a function $d: M/_\sim \times M/_\sim \to \RR$
by:

$$ d(A_x, A_y) = \min \{ d(z,w) \ : \ z \in A_x \ , \ w \in A_y \} $$

Notice that this function may not be a metric since the triangle
inequality may fail. However, in a certain sense, we have that we
can control the topology of $M/_\sim$ using $d$.

\begin{prop}\label{PropositionCocienteCelularEsMetrico}
The quotient topology on $M/_\sim$ verifies the following: For
every $U$ open set in $M/_\sim$ and $p \in U$ there exists $\eps$
such that $B_\eps^d(p) = \{ \pi(y) \  : \ d(A_y,\pi^{-1}(p)) <
\eps \}$ is contained in $U$. Conversely, for every $\eps>0$ and
$p \in M/_\sim$ we have that $B_\eps^d(p)$ contains a neighborhood
of $p$.
\end{prop}

\dem First consider an open set $U \in M/_\sim$ with the quotient
topology. This means that the preimage $\pi^{-1}(U)$ in $M$ is
open.

We must show that for $\pi(x) \in U$  there exists $\eps$ such
that $B^d_\eps(\pi(x))$ is contained in $U$ (here $B_\eps^d$
denotes the $\eps$-ball for the function defined above, that is,
$B_\eps^d(\pi(x)) = \{ \pi(y) \  : \ d(A_y,A_x) < \eps \}$). For
this, we use that the decomposition is upper semicontinuous, thus,
given an open set $V$ of $A_x = \pi^{-1}(\pi(x))$ there exists
$\eps>0$ such that for every $y \in B_\eps(A_x)$ we have that $A_y
\en V$. Considering $V= \pi^{-1}(U)$ we have that there is an open
set for $d$ contained in $U$ as desired.

Now, let $\eps>0$ and $\pi(x)\in M/_\sim$ we must show that
$B^d_\eps(\pi(x))$ contains an open set for the quotient topology.
This is direct since $\pi^{-1}(B_\eps^d(\pi(x)))$ contains the
$\eps$-neighborhood of $A_x$ and so it contains an open saturated
set again by the upper semicontinuity.

\lqqd

\section{A criterium for localization of chain-recurrence
classes}\label{SubSection-LocalizacionDeClasses}

We give a criteria obtained in \cite{PotWildMilnor} in order to
guaranty that a certain (wild) chain-recurrence class is
accumulated by dynamics of lower dimensions. This goes in
contraposition with other kind of ``wildness'' such as universal
properties or viral ones though it is not clear at the moment how
they are related (see \cite{B}).

Given a homeomorphism $g:\Gamma \to \Gamma$ where $\Gamma$ is a
compact metric space, we say that $g$ is \emph{expansive} if there
exists $\alpha>0$ such that for any pair of distinct points $x
\neq y \in \Gamma$ there exists $n\in \ZZ$ such that
$d(g^n(x),g^n(y))\geq \alpha$.

Given a center-stable plaque family $\{\cW^{cs}_x\}_{x\in
\Lambda}$ for a partially hyperbolic set $\Lambda$ and a set $C$
contained in one of the plaques $\cW^{cs}_x$ we define the
\emph{center-stable frontier} of $C$ which we denote as
$\partial^{cs} C$ to the set of points $z$ in $\cW^{cs}_x$ such
that every ball centered in $z$ intersects both $C$ and
$\cW^{cs}_x \cap C^c$ (i.e. the frontier with respect to the
relative topology).

\begin{prop}\label{ProposicionMecanismo}
Let $f$ be a $C^1$-diffeomorphism and $U$ a filtrating set such
that its maximal invariant set $\Lambda$ admits a partially
hyperbolic structure $T_\Lambda M = E^{cs} \oplus E^u$ such that
it admits a locally invariant plaque family $\{\cW^{cs}_x\}_{x\in
\Lambda}$ tangent to $E^{cs}$ whose plaques are contained in $U$.
Assume that there exists a continuous surjective map $h: \Lambda
\to \Gamma$, a homeomorphism $g : \Gamma \to \Gamma$ and a
chain-recurrence class $Q$ such that:
\begin{itemize}
\item[-] $h \circ f= g \circ h$. \item[-] $h$ is injective in
unstable manifolds. \item[-] For every $x\in \Lambda$ we have that
$h^{-1}(\{h(x)\})$ is contained in $\cW^{cs}_x$ and $\partial^{cs}
h^{-1}(\{h(x)\}) \en Q$. In particular $h(Q)=\Gamma$. \item[-] The
fibers $h^{-1}(\{y\})$ are invariant under unstable holonomy.
\item[-] $g$ is expansive.
\end{itemize}
Then, every chain-recurrence class in $U$ different from $Q$ is
contained in the preimage of a periodic orbit by $h$.
\end{prop}

For simplicity, the reader can follow the proof assuming that $g$
is an Anosov diffeomorphism we shall make some footnotes when some
differences (which are quite small) appear.

We remark that the hypothesis of having the fibers invariant under
unstable holonomy is necessary for proving the result and does not
follow from the others in general.

Before starting with the proof we would like to comment on the
hypothesis which are quite strong. We are asking the fibers to be
invariant under unstable holonomy (which is only defined on
$\Lambda$) and asking that the center-stable frontier of the
fibers to be contained in the chain-recurrence class $Q$. This
implies that in order to have the possibility of $\Lambda$ being
different from $Q$ we must have at least some fibers of $h$ which
have interior in some of the center-stable plaques (and by the
holonomy invariance in many of the center-stable plaques). In
general, the most difficult hypothesis to verify will then be that
the center-stable frontier is always contained in the
chain-recurrence class, to do this we use different arguments
depending on the example.

\dem Let $R \neq Q$ be a chain recurrence class of $f$. Then,
since $\partial h^{-1}(\{y\}) \en Q$ for every $y\in \Gamma$, we
have that $R \cap int(h^{-1}(\{y\})) \neq \emptyset$ for some $y
\in \Gamma$.

Conley's theory gives us an open neighborhood $V$ of $R$ whose
closure is disjoint with $Q$ and such that every two points $x,z
\in R$ are joined by arbitrarily small pseudo-orbits contained in
$V$.

\begin{figure}[ht]
\begin{center}
\input{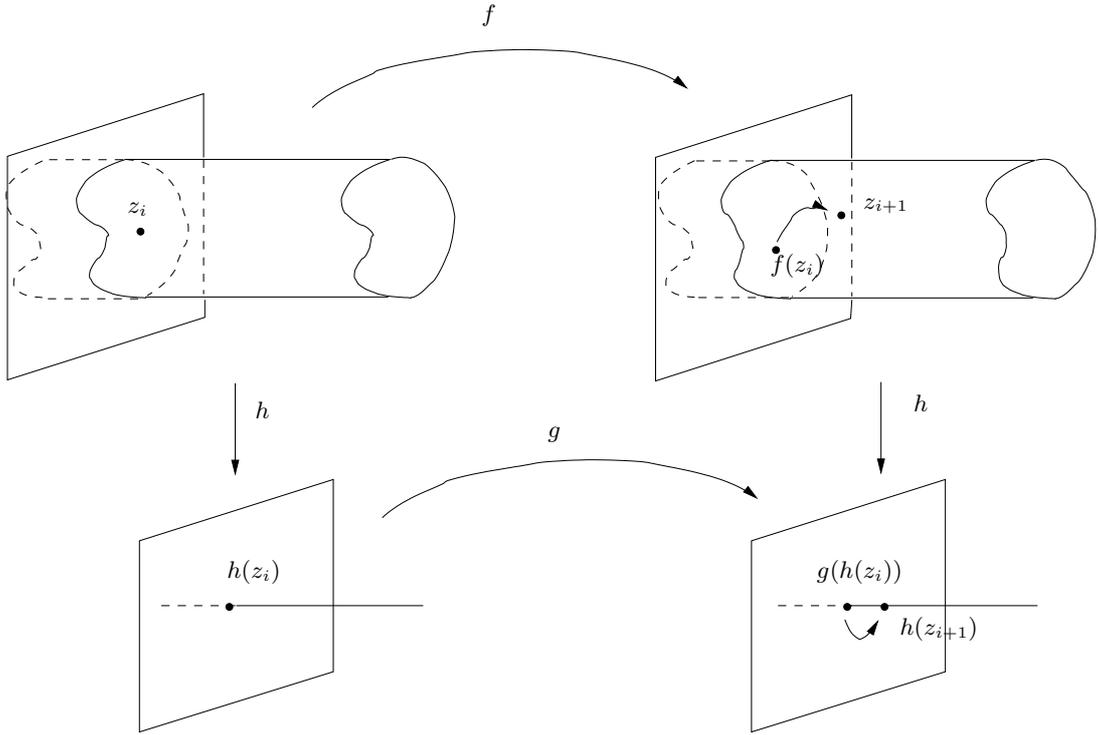}
\caption{\small{Pseudo-orbits for $f$ are sent to pseudo-orbits of
$g$ with jumps in the unstable sets.}} \label{FiguraPseudoOrbitas}
\end{center}
\end{figure}

Since $\overline{V}$ does not intersect $Q$, using the invariance
under unstable holonomy of the fibers, we get that there exists
$\eta_0$ such that if $d(x,z)<\eta_0$ and $x\in V$, then $h(x)$
and $h(z)$ lie in the same local unstable manifold\footnote{The
$\zeta-$local unstable set of a point $x$ for an expansive
homeomorphism $g$ is the set of points whose orbit remains at
distance smaller than $\zeta$ for every past iterate. For an
expansive homeomorphism, this set is contained in the unstable
set.}: In fact, choose $\eta_0< d(\overline V, Q)/2$ and assume
that the image of $h(z)$ is not in the unstable manifold of
$h(x)$, then, we get that if $\gamma:[0,1]\to U$ is the straight
segment joining $x$ and $z$ there is a last point $t_0$ such that
$h(\gamma(t_0)) \en W^u(h(x),g)$, it is not hard to show that
$\gamma(t_0)$ must then belong to
$\partial^{cs}h^{-1}(h(\gamma(t_0))) \en Q$ contradicting that the
straight segment cannot intersect $Q$ from the choice of $\eta_0$.

Given $\zeta>0$ we choose $\eta>0$ such that $d(x,z)<\eta$ implies
$d(h(x), h(z))< \zeta$. The semiconjugacy implies then that if
$z_0, \ldots z_n$ is a $\eta-$pseudo orbit for $f$, then $h(z_0),
\ldots, h(z_n)$ is a $\zeta$-pseudo orbit for $g$ (that is,
$d(g(h(z_i)), h(z_{i+1})) < \zeta$). Also, if $\eta<\eta_0$ and
$z_0, \ldots z_n$ is contained in $U$, then we get that the the
pseudo-orbit $h(z_0), \ldots, h(z_n)$ has jumps inside local
unstable sets (i.e. $h(z_{i+1}) \in W^{u}_{\zeta}(g(h(z_i)))$).

Take $x \in R$. Then, for every $\eta < \eta_0$  we take $x=z_0,
z_1, \ldots, z_n=x$ ($n\geq 1$) a $\eta-$pseudo orbit contained in
$V$ joining $x$ to itself. Thus, we have that

$$g^n(W^{u}(h(x)))= W^u(h(x))$$

\noindent so, $W^u(h(x))$ is the unstable manifold for $g$ of a
periodic orbit $\mathcal O$. Since $R$ is $f$-invariant and since
the semiconjugacy implies that $f^{-n}(x)$ accumulates on
$h^{-1}(\cO)$, we get that $R$ intersects the fiber
$h^{-1}(\mathcal O)$.

We must now prove that $R \en h^{-1}(\cO)$ which concludes.

Given $\eps>0$ there exists $\delta> 0$ such that if $z_0, \ldots
z_n$ is a $\delta-$pseudo orbit for $g$ with jumps in the unstable
manifold, then $z_n \in \cO$ implies that $z_0 \in
W^{u}_{\eps}(\cO)$ (notice that a pseudo orbit with jumps in the
unstable manifold of a periodic orbits can be regarded as a pseudo
orbit for a homothety\footnote{In the general case of $g$ being an
expansive homeomorphism, it is very similar since one has that
restricted to the unstable set of a periodic orbit, one can obtain
a metric inducing the same topology where $g^{-1}$ is an uniform
contraction. This follows from \cite{Fathi} and can also be
deduced using the uniform expansion of $f$ in unstable leaves and
the injectivity of the semi-conjugacy along unstable leaves.} in
$\RR^k$).

Assume that there is a point $z\in R$ such that $h(z) \in
W^u(\mathcal O) \backslash \mathcal O$. Let $\eps$ such that
$d(h(z), \cO) > \eps$. Since $R$ intersects $h^{-1}(\cO)$ there
are $\delta$-pseudoorbits joining $z$ to $h^{-1}(\cO)$ for every
$\delta>0$. This implies that after sending the pseudo orbit by
$h$ we would get $\delta$-pseudo orbits for $g$, with jumps in the
unstable manifold, joining $h(z)$ with $\mathcal O$. This
contradicts the remark made in the last paragraph.

So, we get that $R$ is contained in $h^{-1}(\mathcal O)$ where
$\mathcal O$ is a periodic orbit of $g$. \lqqd

\section{Diffeomorphisms homotopic to Anosov ones, $C^0$
perturbations of hyperbolic
sets}\label{SubSection-SemiconjugacionDeHiperbolicos}

Let $f: \TT^d \to \TT^d$ be a diffeomorphism which is isotopic to
a linear Anosov automorphism $A: \TT^d \to \TT^d$ (i.e. the
diffeomorphism induced by a hyperbolic matrix in $GL(d,\ZZ)$ with
determinant $\pm 1$). Along this text, we assume that $A \in
SL(d,\ZZ)$ which does not represent a loss in generality since the
results are invariant under considering iterates.

In the context we are working, being isotopic to a linear Anosov
automorphism is equivalent to the fact that the induced action
$f_\ast$ of $f$ on the (real) homology (which equals $\RR^d$) is
hyperbolic (see \cite{FranksAnosov}).

We shall denote as $A$ to both the diffeomorphism of $\TT^d$ and
to the hyperbolic matrix $A \in SL(d,\ZZ)$ which acts in $\RR^d$
and is the lift of the torus diffeomorphism $A$ to the universal
cover.

Let $p:\RR^d \to \TT^d$ be the covering projection, and $\tilde f
: \RR^d \to \RR^d$ the lift of $f$ to its universal cover. Notice
that the fact that $f_\ast= A$ implies that there exists $K_0>0$
such that $d(\tilde f(x), A x) < K_0$ for every $x\in \RR^3$.

Classical arguments (see \cite{Walters}) give that there exists
$K_1$ such that for every $x\in \RR^d$, there exists a unique
$y\in \RR^d$ such that $$d(\tilde f^n(x), A^n y) \leq K_1 \qquad
\forall n \in \ZZ$$

We say that the point $y$ \emph{shadows} the point $x$. Notice
that uniqueness implies that the point associated with $x +
\gamma$ is $y+\gamma$ where $\gamma \in \ZZ^d$. We get the
following well known result:

\begin{prop}\label{PropExisteSemiconjugacion}
There exists $H: \RR^d \to \RR^d$ continuous and surjective such
that $H \circ \tilde f = A \circ H$. Also, it is verified that
$H(x + \gamma) = H(x) + \gamma$ for every $x \in \RR^d$ and
$\gamma \in \ZZ^d$ so, there exists also $h: \TT^d \to \TT^d$
homotopic to the identity such that $h \circ f = A \circ h$.
Moreover, we have that $d(H(x), x)< K_1$ for every $x\in \RR^d$.
\end{prop}

\dem Any orbit of $\tilde f$ is a $K_0$-pseudo-orbit of the
hyperbolic matrix $A$. This gives that for every $x$ we can
associate a unique point $y$ such that

 $$d(\tilde f^n(x), A^n y) \leq K_1
\qquad \forall n \in \ZZ$$

We define $H(x)=y$. It is not hard to show that $H$ is continuous.
Since it is at distance smaller than $K_1$ from the identity, we
deduce that $H$ is surjective (this follows from a degree
argument, see \cite{Hatcher} Chapter 2). Periodicity follows from
the fact that all maps here project to the torus. \lqqd

It is well known and easy to show that $H(W^\sigma(x,\tilde f))
\en W^\sigma(H(x),A)$ with $\sigma= s,u$.

The previous result generalizes to general $C^0$-perturbations of
hyperbolic sets. We get the following classical result whose proof
is very similar to the previous one:

\begin{prop}\label{Proposition-PersistenciaHiperbolicos}
Let $\Lambda \en M$ be a hyperbolic set for a diffeomorphism $f$
such that it is maximal invariant in a neighborhood $U$ of
$\Lambda$. Then, there exists $\eps>0$ such that for every
homeomorphism $g$ which is $\eps$-$C^0$-close to $f$ in $U$ we
have that if $\Lambda_g$ is the maximal invariant set for $g$ in
$U$ then there exists a continuous and surjective map $h:
\Lambda_g \to \Lambda$ such that:

$$ f|_{\Lambda} \circ h = h \circ g|_{\Lambda_g} $$
\end{prop}

   \chapter{Attractors and quasi-attractors}\label{Capitulo-Atractores}

This chapter is devoted to the study of attractors and
quasi-attractors for $C^1$-generic dynamics (see subsection
\ref{SubSection-AttractingSets}). We will present the results
obtained in \cite{PotExistenceOfAttractors, PotGenericBiLyapunov,
PotWildMilnor}.

The chapter is organized as follows:

\begin{itemize}
\item[-] In Section \ref{Section-AtractoresSuperficies} we present
a proof of a result by Araujo stating that $C^1$-generic
diffeomorphisms of surfaces have hyperbolic attractors. In
addition, we prove a result from \cite{PotGenericBiLyapunov} in
the context of surfaces which we believe makes its proof more
transparent.

\item[-] In Section \ref{Section-GenericBiLyapunov} we study
quasi-attractors of $C^1$-generic diffeomorphisms in any dimension
and present some results regarding their structure. Then, we give
as an application some results on bi-Lyapunov stable homoclinic
classes.

\item[-] In Section \ref{Section-EjemplosQuasiAtractores} we
present several examples of dynamics without attractors and of
robustly transitive attractors. We present the results from
\cite{BLY} and then the ones of \cite{PotWildMilnor}. We profit to
add some examples which we believe may have some interest in the
general theory.

\item[-] In Section \ref{Section-TrappingAttractors} we present a
definition which covers a certain class of quasi-attractors in
dimension $3$ and explain why we believe this class of examples
should be studied.

\end{itemize}

\section{Existence of hyperbolic attractors in surfaces}\label{Section-AtractoresSuperficies}

In dimension 2, the result of the existence of attractors for
$C^1$-generic diffeomorphisms was announced to be true by Araujo
(\cite{Araujo}) but the result was never published since there was
a gap on its proof. However, it has become a folklore result: with
the techniques of \cite{PujSamAnnals1} the gap in the proof can be
arranged (see \cite{BrunoSantiago}). There has also been an
announcement of this result in \cite{BLY}.

We prove here the following theorem which is similar to the one by
Araujo. The proof we give is quite short but based on the recent
developments of generic dynamics (mainly \cite{ABC},\cite{BC},
\cite{MP} and \cite{PujSamAnnals1}). The proof was made available
in \cite{PotExistenceOfAttractors}.

\begin{teo}\label{Teorema-ExistenciaAtractoresSuperficies}
There is $\cG \en \Diff^1(M^2)$, a residual subset of diffeomorphisms in the surface $M^2$
such that for every $f\in \cG$, there is an hyperbolic attractor.
Moreover, if $f$ has finitely many sinks, then $f$ is essentially
hyperbolic.
\end{teo}

We say that $f$ is \emph{essentially hyperbolic} if it admits
finitely many hyperbolic attractors and such that the union of
their basins cover an open and dense subset of $M$ (Araujo proves
that the basin of atraction has Lebesgue measure one, his
techniques work in this context too, see \cite{BrunoSantiago}).
This definition comes from \cite{PT} and is motivated by a new
result of \cite{CP} which closes a long standing problem posed by
Palis in \cite{PT} (though a stronger formulation remains open and
important).

We will prove the following Theorem in any dimensions in Section
\ref{Section-GenericBiLyapunov}, however, we present it here
before since the proofs are easier to follow in the surface case.

\begin{teo}\label{Teorema-FinitasFuentesImplicaSobreQA}
There is $\cG \en \Diff^1(M^2)$, a residual subset of diffeomorphisms
in the surface $M^2$ such that for every $f\in \cG$ with finitely
many sources satisfies that every homoclinic class which is a
quasi-attractor is an hyperbolic attractor.
\end{teo}

In particular we get the following using results in \cite{MP} and
\cite{BC} (see Theorem \ref{Teorema-BonattiCrovisier}):

\begin{cor} There is $\cG \en \Diff^1(M^2)$, a residual subset of diffeomorphisms in the surface $M^2$
such that for every $f\in \cG$ with finitely many sources
satisfies that generic points converge either to hyperbolic
attractors or to aperiodic classes.
\end{cor}

This last Corollary applies for example in a $C^1$-neighborhood of
the well known Hen\'on attractor (see \cite{BDV} Chapter 4). In
fact, since hyperbolic attractors which are in a disc which is
dissipative are sinks, in the Hen\'on case we get that there are
no non-trivial attractors (aperiodic quasi attractors for generic
diffeomorphisms cannot be attractors).

\subsection{Proof of the Theorem \ref{Teorema-ExistenciaAtractoresSuperficies}}

Let $\mathcal K$ be the set of all compact subsets of $M$ with the
Hausdorff topology.

Let $S: \Diff^1(M^2) \to \mathcal K$ be the map such that $S(f)=
\overline{\Per_0(f)}$ is the clousure of the set of sinks of $f$.

Since $S$ is semicontinuous (see Remark
\ref{Remark-Semicontinuidad}), there exists a residual subset
$\cG_0$ of $\Diff^1(M)$ such that for every $f\in \cG_0$, the
diffeomorphism $f$ is a continuity point of $S$.

This implies that we can write $\cG_0= \cA \cup \cI$ open sets in
$\cG_0$ such that for every $f\in \mathcal A$ the number of sinks
is locally constant and finite (that is, there is a neighborhood
$\U$ of $f$ in $\Diff^1(M)$ such that for every $g\in \U$ the
number of sinks is the same and they vary continuously), and such
that for every $f \in \mathcal I$ there are infinitely many sinks.

To prove the Theorem it is enough to work inside $\tilde{\mathcal
A}$ (an open set in $\Diff^1(M)$ such that $\mathcal A =
\tilde{\mathcal A} \cap \cG_0$) since the theorem is trivially
satisfied in $\mathcal I$.

Let $\cG= \cG_0 \cap \cG_{BC} \cap \cG_{BDP} \cap \cG_{PS}$ (see
Theorems \ref{Teorema-BonattiCrovisier},
\ref{TeoremaBonattiDiazPujals} and \ref{Teorema-PujSamC1}). Let $f
\in \tilde{\mathcal A}\cap \cG$. We must show that $f$ is
essentially hyperbolic.

\medskip

\noindent {\bf Step 1:} We first prove that every quasi-attractor
$\cQ$ is a hyperbolic attractor.

\medskip

We have that $\Lambda$ admits a nested sequence of open
neighborhoods $U_n$ such that $\cQ= \bigcap_{n\geq 0} U_n$ and
such that $f(\overline{U_{n}}) \en U_n$. The following lemma holds
in every dimension:

\begin{lema}\label{Lema-BaseDeEntornosMedidaDisipativa}
Let $\cQ$ be a quasi-attractor. Then, there exist an ergodic
measure $\mu$ supported in $\cQ$ such that $\int
log(|det(Df)|)d\mu \leq 0$.
\end{lema}

\dem Let $m_n$ be the normalized Lebesgue measure in $U_n$.

Consider $\mu_n$ a limit point in the weak-$\ast$ topology of the
sequence of measures given by

$$\nu_k=\frac{1}{k} \sum_{i=1}^k f^i_{\#}(m_n)$$

\noindent which is an invariant measure supported in
$f(\overline{U_n})$. We are here using the following notation:
$f_{\#}(\nu) (A) = \nu(f^{-1}(A))$.

Since $f(\overline{U_n}) \en U_n$, the change of variables theorem
and Jensen's inequality (see \cite{RudinFA}) implies that:

$$\int \log(|\det Df|) dm_n <0. $$

The same argument, using the fact that $f^k(\overline{U_n})\en
f^{k-1}(U_n)$ implies that

$$\int \log(|\det Df|) df^k_{\#}(m_n) <0.$$

We obtain that

$$\int \log(|\det Df|) d\nu_k \leq 0$$

From which we deduce that$\int \log(|\det Df|) d\mu_n \leq 0$.
Now, consider a measure $\mu$ which is a limit point in the
weak-$\ast$ topology of the measures $\mu_n$.

The measure $\mu$ must be an invariant measure, supported on $\cQ$
and satisfying that $$\int \log(|\det Df|) d\mu \leq 0.$$

Using Proposition \ref{Proposicion-MedidasErgodicas} one can
assume that $\mu$ is ergodic. \lqqd

Since the set of sinks varies continuously with $f$ and there are
finitely many of them, we can choose $n$ such that there are no
sinks in $U_n$.

Using the Ergodic closing Lemma (Theorem
\ref{Teorema-ErgodicClosingLemma}) we get that the support of the
measure must admit a dominated splitting: Otherwise we get
periodic points converging in the Hausdorff topology to the
support of the measure and with $\log(|\det Df^{\pi(p)}|)$
arbitrarily close to zero. If they do not admit a dominated
splitting, using a classical argument (see subsection
\ref{SubSection-BDPBGVBoBo}) one can convert them into sinks by
applying Franks' Lemma (Theorem \ref{Teorema-FranksLemma}), a
contradiction.

Also, the measure must be hyperbolic: If no positive exponents
exist, since $\int \log(|\det(Df)|) d\mu \leq 0$, one can approach
the measure by periodic orbits with both exponents smaller than
$\eps$ (arbitrarily small) by using the Ergodic closing Lemma
(Theorem \ref{Teorema-ErgodicClosingLemma}) and one can create a
sink by making a further small perturbation with Franks' Lemma
(Theorem \ref{Teorema-FranksLemma}).

Using Theorem \ref{Teorema-PesinC1}, we deduce that the support of
$\mu$ is contained in a homoclinic class. Since $f \in \cG_{BC}$
we have that $\cQ$ is a homoclinic class.

Also we get periodic points inside the class such that $\log(|\det
Df^{\pi(p)}|)<\eps$ for small $\eps>0$.

Using Proposition \ref{Proposicion-transitions} we get that
periodic points with this property are dense in the homoclinic
class and so we get a dominated splitting $T_{\cQ}M= E\oplus F$ in
the whole class. In fact, since we are far from sinks and $f\in
\cG_{BDP}$, we get that $F$ must be uniformly expanding.

Since we are in $\cG_{PS}$ we get that $\cQ$ is hyperbolic (see
Theorem \ref{Teorema-PujSamC1}), and thus, $\cQ$ is a hyperbolic
attractor (see Proposition
\ref{Proposition-HiperbolicoImplicaAislado}).

This proves the first assertion of the Theorem.

\medskip
\noindent {\bf Step 2:} We now prove that in fact $f$ must be
essentially hyperbolic.

\medskip

Suppose first that there are infinitely many non-trivial
hyperbolic attractors (recall that we are assuming that there are
finitely many sinks). Assume $\cQ_n$ is a sequence of distinct
hyperbolic attractors such that $\cQ_n \to K$ in the Hausdorff
topology. From Proposition
\ref{Proposicion-VariacionClasesDeRecurrencia} we have that $K$
must be a chain-transitive set, this implies that $K \cap
S(f)=\emptyset$ (notice that since $S(f)$ is a finite set it will
be isolated from the chain-recurrent set).

Notice that there are measures $\mu_n$ supported in $\cQ_n$ such
that

$$\int \log(|det (Df)|) d\mu_n \leq 0.$$

Consider a weak-$\ast$ limit $\mu$ of these measures, so we have
that $\mu$ is supported in $K$ and verifies that $$\int \log(|det
(Df)|) d\mu \leq 0.$$

So using the same argument as before we deduce that $K$ is
contained in a hyperbolic homoclinic class, and thus isolated, a
contradiction with the fact that it was accumulated by
quasi-attractors.

Since $f\in \cG_{BC}$ (see Theorem
\ref{Teorema-BonattiCrovisier}), generic points in the manifold
converge to Lyapunov stable chain recurrence classes and we get
that there is an open and dense subset of $M$ in the basin of
hyperbolic attractors. This finishes the proof of the Theorem.

\lqqd

\subsection{Proof of the Theorem
\ref{Teorema-FinitasFuentesImplicaSobreQA}}

Theorem \ref{Teorema-FinitasFuentesImplicaSobreQA} is implied by
the following Theorem from \cite{PotGenericBiLyapunov}. This
theorem will be extended to higher dimensions in Section
\ref{Section-GenericBiLyapunov} but we prefer to present a proof
in this context since it helps to grasp better the idea involved
(which is the use of Theorem \ref{Teorema-FranksLemmaGourmelon}):

\begin{teo}\label{Teorema-QuasiAtractoresEnDimension2}
Let $H$ be a homoclinic class of a $C^1$-generic surface diffeomorphism $f$ which is a quasi-attractor.
Then, if $H$ has a periodic point $p$ such that
$det(D_pf^{\pi(p)})\leq 1$ then $H$ admits a dominated splitting
and thus it is a hyperbolic attractor.
\end{teo}

\dem Consider a generic diffeomorphism $f$. Consider a periodic
point $q\in H$ fixed such that for a neighborhood $\mathcal U$ of
$f$ the class $H(q_g,g)$ is a quasi-attractor for every $g$ in a
residual subset of $\mathcal U$ (see subsection
\ref{SubSection-PersistenciaQA}).

By genericity, we can assume that every periodic point $q \in H$
verifies that $det(D_pf^{\pi(p)}) \neq 1$ and using Proposition
\ref{Proposicion-transitions} we deduce there is a dense set of
points verifying that the determinant is smaller than $1$.

By Theorem \ref{Teorema-BonattiBochi} we know that if $H$ does not
admit a dominated splitting then for every $\eps>0$ we can modify
the derivative of $f$ along a periodic orbit along a curve with
small diameter and which satisfies the hypothesis of Gourmelon's
version of Franks' Lemma (Theorem
\ref{Teorema-FranksLemmaGourmelon}).

We can then apply Theorem \ref{Teorema-FranksLemmaGourmelon} to
perturb the periodic point $p$ preserving its strong stable
manifold locally.

The idea is to make the perturbation in a small neighborhood of
$p$ without breaking the intersection between $W^u(q)$ and
$W^s(p)$. This creates a sink in $p$ but such that $W^u(q)$ still
intersects its basin. Since $H(q_g,g)$ is still a quasi-attractor,
this means that $p_g \in H(q_g,g)$ which is a contradiction (see
Lemma \ref{normadiferencial} for a more detailed proof).

Now by Theorem \ref{Teorema-PujSamC1} we get that $H$ must be
hyperbolic which concludes.

\lqqd

The last theorem has some immediate consequences which may have
some interest on their own.

We say that an embedding $f:\D^2 \to \D^2$ is \emph{dissipative}
if for every $x\in \D^2$ we have that $|det(D_xf)|<b<1$. Recall
that for a dissipative embeddings of the disc, the only hyperbolic
attractors are the sinks (\cite{Plykin}).

\begin{cor}
Let $f:\D^2 \to \D^2$ be a generic dissipative embedding. Then,
every quasi-attractor which is a homoclinic class is a sink.
\end{cor}


\section{Structure of quasi-attractors}\label{Section-GenericBiLyapunov}

In this section we explore the existence of invariant structures
for quasi-attractors of $C^1$-generic dynamics in any dimensions.
We expect that, for a homoclinic class, being a quasi-attractor
imposes sufficiently many structure in the dynamics on the tangent
map. We have obtained partial results in this direction.

The main difficulty is that domination is very much related to
either isolation of the homoclinic class or with being far from
homoclinic tangencies, in this result we manage to deal with
homoclinic classes for which we do not know a priori that either
of these conditions are satisfied. The main idea is to use the
fact that being a quasi-attractor is a somewhat robust property
and the perturbation techniques developed by Gourmelon (Theorem
\ref{Teorema-FranksLemmaGourmelon}). This allows us to prove:

\begin{teo}\label{Teorema-QATieneDescomposicionDominada}
For every $f$ in a residual subset $\cG_1$ of $\Diff^1(M)$, if $H$
is a homoclinic class for $f$ which is a quasi-attractor and there
is a periodic point $p\in H$ such that $det(D_pf^{\pi(p)})\leq 1$,
then, $H$ admits a non-trivial dominated splitting.
\end{teo}

This theorem is proved in subsection
\ref{SubSection-DescDominadaparaQA}. In view of Lemma
\ref{Lema-BaseDeEntornosMedidaDisipativa} one can ask (see
\cite{B} Conjecture 2 for a stronger version of this question):

\begin{quest}
Is it true that every quasi-attractor $\cQ$ of a $C^1$-generic
diffeomorphism which is a homoclinic class has a periodic point
$p$ such that $det(D_pf^{\pi(p)})\leq 1$?
\end{quest}

Another important task is to determine whether some (extremal)
bundles are uniform in order to be able to derive further
dynamical properties. In this direction, we have obtained the
following result motivated by previous results obtained in
\cite{PotSambarino}:

\begin{teo}\label{Teorema-ExtremalDimension1}
There exists a residual set $\cG_2$ of $\Diff^1(M)$ such that if
$f$ is a diffeomorphism in $\cG_2$ and $H$ a homoclinic class
which is a quasi-attractor admitting a codimension one dominated
splitting $T_H M=E \oplus F$  where $\dim(F)=1$. Then, the bundle
$F$ is uniformly expanding for $f$.
\end{teo}

The proof of this theorem is presented in subsection
\ref{SubSection-FibradoDimension1Extremal}.

In dimension $3$ we have the following corollary about the
structure of quasi-attractors:

\begin{cor}\label{Corolario-QuasiAtractoresEnDim3}
There exists a residual subset $\cG_{QA}$ of $\Diff^1(M^3)$ with
$M$ a $3$-dimensional manifold such that if $f\in \cG_{QA}$ and
$\cQ$ a quasi-attractor of $f$ having a periodic point $p$ such
that $det(D_pf^{\pi(p)}) \leq 1$ then, $\cQ$ admits a dominated
splitting of one of the following forms:
\begin{itemize}
\item[-] $T_{\cQ} M = E^{cs} \oplus E^u$ where $E^u$ is one
dimensional and uniformly expanded and $E^{cs}$ may or may not
admit a sub-dominated splitting. \item[-] $T_{\cQ}M = E_1 \oplus
E^{cu}$ where the bundle $E^{cu}$ is two-dimensional and verifies
that periodic points are volume hyperbolic at the period in
$E^{cu}$.
\end{itemize} Moreover, if $\cQ$ is not accumulated by sinks, then
in the second case we have that $E^{cu}$ is uniformly volume
expanding.
\end{cor}

We also explore some properties of quasi-attractors far from
homoclinic tangencies and we deduce some consequences for
homoclinic classes which are both quasi-attractors and
quasi-repellers (bi-Lyapunov stable homoclinic classes) responding
to some questions posed in \cite{ABD}. The results about dynamics
far from tangencies overlap with \cite{Yang-LyapunovStable}.

\subsection{Persistence of quasi-attractors which are homoclinic classes}\label{SubSection-PersistenciaQA}

The following result will be essential for our proofs of Theorems
\ref{Teorema-QATieneDescomposicionDominada} and
\ref{Teorema-ExtremalDimension1}. It states that for $C^1$-generic
diffeomorphisms, the quasi-attractors which are homoclinic classes
have a well defined continuation.

\begin{prop}[Lemma 3.6 of \cite{ABD}]\label{Proposition-PersistenciaQA}
There is a residual subset $\cG_{ABD}$ of $\Diff^1(M)$ such that
if $f \in \cG_{ABD}$ and $H(p,f)$ is a homoclinic class of $f$
then, there exists a neighborhood $\cU$ of $f$ such that the point
$p$ has a continuation in $\cU$ and such that:
\begin{itemize}
\item[-] If $H(p,f)$ is a quasi-attractor, then $H(p_g,g)$ is also
a quasi-attractor for every $g\in \cG_{ABD} \cap \cU$. \item[-] If
$H(p,f)$ is not a quasi-attractor, then, $H(p_g,g)$ is not a
quasi-attractor for any $g\in \cG_{ABD} \cap \cU$.
\end{itemize}
\end{prop}

\dem We reproduce the proof from \cite{ABD}.

Consider $\cG= \cG_{BC} \cap \cG_{cont}$ a residual subset of
$\Diff^1(M)$ such that for every $f\in \cG$ we have that:

\begin{itemize}
\item[-] A homoclinic class $H(p,f)$ of $f$ is a quasi-attractor
if and only if $H(p,f)= \overline{W^u(p,f)}$ (see Theorem
\ref{Teorema-BonattiCrovisier}). \item[-] There exists a
neighborhood $\cU$ of $f$ such that the following maps $g\mapsto
H(p_g,g)$ and $g\mapsto \overline{W^u(p_g,g)}$ are well defined
and continuous on every $g\in \cG \cap \cU$ (see Remark
\ref{Remark-Semicontinuidad}).
\end{itemize}

For a pair $\cU$ and $p$ which has a continuation for every $f\in
\cU$ we let $\cA_{(\cU,p)} \en \cU \cap \cG$ be the set of $g$
such that $H(p_g,g)\neq \overline{W^u(p_g,g)}$. Since both sets
are compact and vary continuously we deduce that $\cA_{(\cU,p)}$
is open in $\cU \cap \cG$.

Let $\cB_{(\cU,p)}$ be the complement of the closure of
$\cA_{(\cU,p)}$ in $\cU \cap \cG$ which is also open and verifies
that if $g\in \cB_{\cU}$ then there is a neighborhood of $g$ in
$\cG$ consisting of diffeomorphisms $\hat g$ such that the
homoclinic class $H(p_{\hat g}, \hat g)$ is a quasi-attractor.

From how we defined $\cA_{(\cU,p)}$ and $\cB_{(\cU,p)}$ we have
that their union is open and dense in $\cG \cap \cU$.

The residual subset $\cG_{ABD}$ is then obtained by considering a
countable collection of pairs $(\cU,p)$ where $\cU$ varies in a
countable basis of the topology of $\Diff^1(M)$ and $p$ is a
hyperbolic periodic point of $f\in \cU$ which has a continuation
to the whole $\cU$ (there are clearly at most countably many of
them by Kupka-Smale's Theorem \ref{Teorema-KupkaSmale}). We
finally define:

$$ \cG_{ABD} = \bigcap_{(\cU,p)}\left( \cA_{(\cU,p)} \cup \cB_{(\cU,p)} \cup \overline{\cU}^c \right)$$

\lqqd

In general, if $X$ is a Baire topological space and $A$ is a Borel
subset of $X$, then there exists a residual subset $G$ of $X$ such
that $A$ is open and closed in $G$. This usually serves as an
heuristic principle, but in general it is not so easy to show that
a certain property is Borel. One must prove this kind of result by
barehanded arguments.

\subsection{One dimensional extremal bundle}\label{SubSection-FibradoDimension1Extremal}

We will first prove Theorem \ref{Teorema-ExtremalDimension1}. The
proof strongly resembles the proof of the main theorem of
\cite{PotSambarino} and it was indeed motivated by it.

The main difference is that the fact that periodic points must be
hyperbolic at the period in the case the homoclinic class has
non-empty interior is quite direct by using transitions and the
fact that the interior has some persistence properties. This is
the content of the following lemma whose proof will serve also as
a model for the proof of Theorem
\ref{Teorema-QATieneDescomposicionDominada}. We will make emphasis
only in the part of the proof which differs from
\cite{PotSambarino} and we refer the reader to that paper for more
details in the rest of the proof.

\begin{lema}\label{normadiferencial} Let $H$ be a homoclinic class
which is a quasi-attractor of a $C^1$-generic diffeomorphism $f$
such that the class has only periodic orbits of stable index
smaller or equal to $\alpha$. So, there exists $K_0>0$, $\lambda
\in (0,1)$ and $m_0 \in \Z$ such that for every $p\in
\Per_{\alpha}(f|_{H})$ of sufficiently large period one has

$$\prod_{i=0}^k \left\|\prod_{j=0}^{m_0-1} Df^{-1}|_{E^u(f^{-im_0 -j}(p))}\right\| < K_0 \lambda^k \qquad k=\left[\frac{\pi(p)}{m_0}\right]. $$
\end{lema}

\dem Let $\Res$ be a residual subset of $\Diff^1(M)$ such that if
$f\in \Res$ and $H$ is homoclinic class of a periodic point $q$ of
index $\alpha$ and a quasi-attractor, there exists a small
neighborhood $\U$ of $f$ where the continuation $q_g$ of $q$ is
well defined and such that for every $g\in \U \cap \Res$ one has
that $H(q_g,g)$ is quasi-attractor and $g$ is a continuity point
of the map $g\mapsto H(q_g,g)$ (see Proposition
\ref{Proposition-PersistenciaQA} and Remark
\ref{Remark-Semicontinuidad}).

Also, being $f$ generic, we can assume that for every $g \in \U
\cap \Res$ and every $p\in \Per_{\alpha}(g)\cap H(q_g,g)$ we have
that $H(q_g,g)=H(p,g)$, so, the orbits of $p$ and $q_g$ are
homoclinically related (see Theorem
\ref{Teorema-BonattiCrovisier}).

We can also assume that $\U$ and $\Res$ were chosen so that for
every $g \in \U \cap \Res$ every periodic point in $H(q_g,g)$ has
index smaller or equal to $\alpha$. We can also assume that $q_g$
has index $\alpha$ for every $g \in \U$.

Lemma II.5 of \cite{ManheErgodic} asserts that to prove the lemma
it is enough to show that there exists $\eps>0$ such that the set
of cocycles

$$\Theta_{\alpha}=\{ D_{\mathcal O(p)}f^{-1}|_{E^u} \ : \ p\in
\Per_{\alpha}(f|_{H}) \}$$

\noindent which all have its eigenvalues of modulus bigger than
one, verify that every $\eps$-perturbation of them preserves this
property. That is, given $p\in \Per_\alpha (f|_{H})$ one has that
every $\eps$-perturbation $\{A_0,\ldots, A_{\pi(p)-1}\}$ of
$D_{\mathcal O(p)}f|_{E^u}$ verifies that $A_{\pi(p)-1}\ldots A_0$
has all its eigenvalues of modulus bigger or equal to one.

Therefore, assuming by contradiction that the Lemma is false, we
get that $\forall \eps>0$ there exists a periodic point $p\in
\Per_\alpha(f|_{H})$ and a linear cocycle $\{A_0, \ldots,
A_{\pi(p)}\}$ over $p$  satisfying that:

\begin{itemize}
\item[-] $\|D_{f^i(p)}f|_{E^u} -A_i\|\leq \eps$, \item[-]
$\|D_{f^i(p)}f^{-1}|_{E^u} -A_i^{-1}\|\leq \eps$ and \item[-]
$\prod_{i=0}^{\pi(p)-1}A_i$ has some eigenvalue of modulus smaller
or equal to $1$.
\end{itemize}

In coordinates $T_{\mathcal O(p)}M =E^u\oplus (E^u)^\perp$, since
$E^u$ is invariant we have that the form of $Df$ is given by

\[ D_{f^i(p)}f= \left(
     \begin{array}{cc}
       D_{f^i(p)}f_{/E^u} & K^1_i(f) \\
       0 & K^2_i(f) \\
     \end{array}
   \right)
 \]

Let $\gamma:[0,1]\to \Gamma_{\alpha}$ given in coordinates
$T_{\mathcal O(p)}M =E^u\oplus (E^u)^\perp$ by

$$\gamma_i(t)= \left(
     \begin{array}{cc}
      (1-t)D_{f^i(p)}f|_{E^u} + t A_i & K^1_i(f) \\
       0 & K^2_i(f) \\
     \end{array}
   \right)$$

\noindent whose diameter is bounded by $\eps$ (see Lemma 4.1 of
\cite{BDP}).

\medskip

Now(\footnote{The following argument will be referred too by the
proof of Theorem \ref{Teorema-QATieneDescomposicionDominada}.}),
choose a point $x$ of intersection between $W^s(p,f)$ with
$W^u(q,f)$ and choose a neighborhood $U$ of the orbit of $p$ such
that:

\begin{itemize}
\item[(i)] It does not intersect the orbit of $q$. \item[(ii)] It
does not intersect the past orbit of $x$. \item[(iii)] It verifies
that once the orbit of $x$ enters $U$ it stays there for all its
future iterates by $f$.
\end{itemize}

It is very easy to choose $U$ satisfying (i) since both the orbit
of $p$ and the one from $q$ are finite. Since the past orbit of
$x$ accumulates in $q$ is not difficult to choose $U$ satisfying
(ii). To satisfy (iii) one has only to use the fact that $x$
belongs to the stable manifold of $p$ so, after a finite number of
iterates, $x$ will stay in the local stable manifold of $p$. It is
then not difficult to choose a neighborhood $U$ which satisfies
(iii) also.

Applying Theorem \ref{Teorema-FranksLemmaGourmelon} we can perturb
$f$ to a new diffeomorphism $\hat g$ so that the orbit of $p$ has
index greater than $\alpha$ and so that it preserves locally its
strong stable manifold. This allows to ensure that the
intersection between $W^u(q_{\hat{g}},\hat g)$ and
$W^s(p,\hat{g})$ is non-empty.

This intersection is transversal so it persist by small
perturbations, the same occurs with the index of $p$ so we can
assume that $\hat g$ is in $\Res \cap \U$.

Using the fact that $H(q_{\hat g},\hat g)$ is a quasi-attractor we
obtain that $p \in H(q_{\hat g},\hat g)$:

This is because quasi-attractors are saturated by unstable sets,
so, since $q_{\hat g} \en H(q_{\hat{g}}, \hat g)$ we have that
$\overline{W^u(q_{\hat{g}},\hat g)} \en H(q_{\hat{g}},\hat g)$ and
since $W^s(p,\hat g)\cap W^u(q_{\hat{g}},,\hat g) \neq \emptyset$,
we get by the $\lambda$-Lemma (Theorem
\ref{Teorema-LambdaLemmaOriginal}) and the fact that the
quasi-attractor is closed that

$$p \in \overline{W^u(p, \hat g)} \en \overline{W^u(q_{\hat{g}},\hat g)} \en H(q_{\hat g},\hat g) $$

This contradicts the choice of $\mathcal U$  since we find a
diffeomorphism in $\U\cap \Res$ with a periodic point with index
bigger than $\alpha$ in the continuation of $H$, and so the lemma
is proved.

\lqqd

\begin{obs}\label{hipeneelperiodocodimensionuno}
One can recover Lemma 2 of \cite{PotSambarino} in this context. In
fact, if there is a codimension one dominated splitting of the
form $T_H M =E\oplus F$ with $\dim F=1$ then (using the adapted
metric given by \cite{GourmelonAdaptada}) for a periodic point of
maximal index one has $\|Df^{-1}|_{F(p)}\|\leq
\|Df^{-1}|_{E^u(p)}\|$ so,
$$\prod_{i=0}^k \left\|\prod_{j=0}^{m_0-1} Df^{-1}|_{F(f^{-im_0 -j}(p))}\right\| < K_0 \lambda^k \qquad k=\left[\frac{\pi(p)}{m_0}\right] $$

And since $F$ is one dimensional one has $\prod_i\|A_i\|
=\|\prod_i A_i\|$ so $\|Df^{-\pi(p)}|_{F(p)}\|\leq K_0
\lambda^{\pi(p)}$ (maybe changing the constants $K_0$ and
$\lambda$).

In fact, there is $\gamma \in (0,1)$ such that for every periodic
point of maximal index and big enough period one has $
\|Df^{-\pi(p)}|_{F(p)}\|\leq \gamma^{\pi(p)}$

Also, it is not hard to see, that if the class admits a dominated
splitting of index bigger or equal than the index of all the
periodic points in the class, then, periodic points should be
hyperbolic in the period along $F$ (for a precise definition and
discussion on this topics one can read \cite{BGY},
\cite{wen-Selecting}). \finobs
\end{obs}

\begin{obs}\label{Remark-Deapedazos}
As a consequence of the proof of the lemma we get that: One can
perturb the eigenvalues along an invariant subspace of a cocycle
without altering the rest of the eigenvalues. The perturbation
will be of similar size to the size of the perturbation in the
invariant subspace. See Lemma 4.1 of \cite{BDP}. Notice also that
we could have perturbed the cocycle $\{K_i^2(f)\}_i$ without
altering the eigenvalues of the cocycle $D_{\cO(p)}f|_{E^u}$.
\finobs
\end{obs}

One can now conclude the proof of Theorem
\ref{Teorema-ExtremalDimension1} with the same techniques as in
the proof of the main Theorem of \cite{PotSambarino}.

\demo{ of Theorem \ref{Teorema-ExtremalDimension1}} We have that
$T_HM=E\oplus F$ with $\dim F=1$. We first prove that the center
unstable curves tangent to $F$ should be unstable and with uniform
size (this is Lemma 3 of \cite{PotSambarino}). To do this, we
first use Lemma \ref{normadiferencial} to get this property in the
periodic points and then use the results from \cite{PujSamBrazil}
and \cite{BC} to show that the property extends to the rest of the
points. This dynamical properties imply also uniqueness of these
central unstable curves.

Assuming the bundle $F$ is not uniformly expanded, one has two
cases: one can apply Liao's selecting lemma or not (see
\cite{Liao, wen-Selecting}).

In the first case one gets weak periodic points inside the class
which contradict the thesis of Lemma \ref{normadiferencial}.

The second case is similar, if Liao's selecting lemma
(\cite{Liao}) does not apply, one gets a minimal set inside $H$
where the expansion along $F$ is very weak and thus $E$ is
uniformly contracting. Using the dynamical properties of the
center unstable curves, classical arguments give that we can
shadow orbits of this minimal sets by periodic points which are
weak in the $F$ direction. Since the stable manifold of this
periodic point will be uniform, it will intersect the unstable
manifold of a point in $H$, and then the fact that $H$ is a
quasi-attractor implies the point is inside the class and again
contradicts Lemma \ref{normadiferencial}.

For more details see \cite{PotSambarino}.

\lqqd

\subsection{Existence of a dominated splitting}\label{SubSection-DescDominadaparaQA}

We  prove here Theorem \ref{Teorema-QATieneDescomposicionDominada}
which state that a homoclinic class which is a quasi-attractor and
has a dissipative periodic orbit admits a dominated splitting.

The idea is the following: in case $H$ does not admit any
dominated splitting we can perturb the derivative of some periodic
point in order to convert it into a sink with the techniques of
\cite{BoBo}.

We pretend to use Theorem \ref{Teorema-FranksLemmaGourmelon} to
ensure that the unstable manifold of a periodic point in the class
intersects the stable set of the sink and reach a contradiction as
we did in the end of the proof of Lemma \ref{normadiferencial}.


\demo{ of Theorem \ref{Teorema-QATieneDescomposicionDominada}} Let
$H$ be a homoclinic class of a $C^1$-generic diffeomorphism $f$
which is a quasi-attractor. Let us assume that $H$ contains
periodic points of index $\alpha$ and we consider
$\Delta^\eta_{\alpha}\en \Per_{\alpha}(f|_{H})$ the set of index
$\alpha$ and $\eta$-disippative periodic points in $H$ for some
$\eta<1$.

It is enough to have one periodic point with determinant smaller
than one to get that for some $\eta<1$, the set
$\Delta^\eta_{\alpha}$ will be dense in $H$ (see Proposition
\ref{Proposicion-transitions}). Notice that from hypothesis, and
the fact that for generic diffeomorphisms the determinant of the
differential at the period is different from one, there is
$\eta<1$ such that $\Delta^\eta_\alpha$ is dense.

Notice that if $H$ admits no dominated splitting, then neither
does the cocycle of the derivatives over $\Delta^\eta_\alpha$.
This implies that we can apply Theorem \ref{Teorema-BonattiBochi}
and there is a periodic point $p \in \Delta^\eta_{\alpha}$ which
can be turned into a sink with a $C^1$-small perturbation done
along a path contained in $\Gamma_{\alpha}$ (which maintains or
increases the index).

Now we are able to use Theorem \ref{Teorema-FranksLemmaGourmelon}
and reach a contradiction. Consider a periodic point $q \in
\Delta^\eta_{\alpha}$ fixed such that for a neighborhood $\mathcal
U$ of $f$ the class $H(q_g,g)$ is a quasi-attractor for every $g
\in \U \cap \cG$.

Suppose the class does not admit any dominated splitting, so, as
we explained above we have a periodic point $p\in
\Delta^\eta_\alpha$ such that $f$ can be perturbed in an
arbitrarily small neighborhood of $p$ to a sink for a
diffeomorphism $g\in \U$ (which we can assume is in $\cG \cap \U$
since sinks are persistent) and preserving locally the strong
stable manifold of $p$. So, we choose a neighborhood of $p$ such
that it does not meet the orbit of $q$ nor the past orbit of some
intersection of its unstable manifold with the local stable
manifold of $p$ with the same argument as in Lemma
\ref{normadiferencial}.

Thus, we get that $W^u(q_g,g) \cap W^s(p,g) \neq \emptyset$ and
using Lyapunov stability we reach a contradiction since it implies
that $p \in H(q_g)$ which is absurd since $p$ is a sink.

\lqqd

One can also deduce some properties on the indices of the possible
dominated splitting depending on the indices of the periodic
points in the class (see \cite{PotGenericBiLyapunov}).

\begin{obs}
\begin{itemize}
\item[-] Also the same ideas give that periodic points in the
class must be volume hyperbolic in the period (not necessarily
uniformly, see \cite{BGY} for a discussion on the difference
between hyperbolicity in the period and uniform hyperbolicity)

\item[-] In fact, we can assume that if a homoclinic class which
is a quasi-attractor admits no dominated splitting, then, there
exists $\eta>1$ such that every periodic point $p$, it has
determinant bigger than $\eta^{\pi(p)}$. Otherwise, there would
exists a subsequence $p_n$ of periodic points with normalized
determinant converging to $1$. After composing with a small
homothety, we are in the hypothesis of Theorem
\ref{Teorema-QATieneDescomposicionDominada}.
\end{itemize}
\finobs
\end{obs}

\subsection{Quasi-attractors far from tangencies}\label{SubSection-QAlejosTangencias}

We present a proof of a result originally proved in
\cite{Yang-LyapunovStable}. We believe that having another
approach to this result is not entirely devoid of interest (see
\cite{CSY} for results that exceed the results here presented on
dynamics far from tangencies).

We remark that in the far from tangencies context, J. Yang
(\cite{Yang-LyapunovStable}) has proved that quasi-attractors of
$C^1$-generic diffeomorphisms are homoclinic classes (see also
\cite{Crov-Hab}) and more recently, C.Bonatti, S.Gan, M.Li and
D.Yang have proved that for $C^1$-generic diffeomorphisms far from
homoclinic tangencies quasi-attractors are in fact essential
attractors (see \cite{BGLY}).

\begin{teo}[\cite{Yang-LyapunovStable} Theorem 3]\label{teoY}
Let $f \in \cG$, where $\cG$ is a residual subset of $\Diff^1(M)
\backslash \overline{\Tang}$, and let $H$ be a Lyapunov stable
homoclinic class for $f$ of minimal index $\alpha$. Let
$T_HM=E\oplus F$ be a dominated splitting for $H$ with $\dim
E=\alpha$, so, one of the following two options holds:

\begin{enumerate}
\item $E$ is uniformly contracting. \item $E$ decomposes as
$E^s\oplus E^c$ where $E^s$ is uniformly contracting and $E^c$ is
one dimensional and $H$ is the Hausdorff limit of periodic orbits
of index $\alpha-1$.
\end{enumerate}
\end{teo}

As in \cite{Yang-LyapunovStable}, the proof has 3 stages, the
first one is to reduce the problem to the central models
introduced by Crovisier, the second one to treat the possible
cases and finally, the introduction of some new generic property
allowing to conclude in the difficult case.

Our proof resembles that of \cite{Yang-LyapunovStable} in the
middle stage (which is the most direct one after the deep results
of Crovisier) and has small differences mainly in the other two.

For the first one, we use a recent result of
\cite{Crov-CentralModels} (Theorem \ref{Teorema-Tricotomia}) and
for the last one, we introduce Lemma \ref{encontrarpuntos} which
can be compared with the main Lemma of \cite{Yang-LyapunovStable}
but the proof and the statement are somewhat different (in
particular, ours is slightly stronger). We believe that this Lemma
can find some applications (see for example \cite{Crov-Hab}).

Before we start the proof of Theorem \ref{teoY} we state the
following theorem due to Crovisier which will be the starting
point for our study:

\begin{teo}[Theorem 1 of \cite{Crov-CentralModels}]\label{Teorema-Tricotomia}
Let $f \in \cG$ where $\cG \en \Diff^1(M) \backslash
\overline{\Tang}$ is residual, and $K_0$ an invariant compact set
with dominated splitting $T_{K_0}M= E \oplus F$. If $E$ is not
uniformly contracted, then, one of the following cases occurs.
\begin{enumerate}
 \item $K_0$ intersects a homoclinic class whose minimal index is
strictly less than $\dim E$. \item $K_0$ intersects a homoclinic
class whose minimal index is $\dim E$ and which contains weak
periodic orbits (for every $\delta$ there is a sequence of
hyperbolic periodic orbits homoclinically related which converge
in the Hausdorff topology to a set $K\en K_0$, whose index is
$\dim E$ but whose maximal exponent in $E$ is in $(-\delta,0)$).
Also, this implies that every homoclinic class $H$ intersecting
$K_0$ verifies that it admits a dominated splitting of the form
$T_HM = E'\oplus E^c \oplus F$ with $\dim E^c =1$. \item There
exists a compact invariant set $K\en K_0$ with minimal dynamics
and which has a partially hyperbolic structure of the form
$T_KM=E^s\oplus E^c\oplus E^u$ where $\dim E^c=1$ and $\dim
E^s<\dim E$. Also, any measure supported on $K$ has zero Lyapunov
exponent along $E^c$.
\end{enumerate}
\end{teo}

Now we are ready to give a proof of Theorem \ref{teoY}.

\demo{ of Theorem \ref{teoY}} Let $\cG \en \Diff^1(M) \backslash
\overline{\Tang}$ be a residual subset such that for every $f \in
\cG$ and every periodic point $p$ of $f$, there exists a
neighborhood $\U$ of $f$, where the continuation $p_g$ of $p$ is
well defined, such that $f$ is a continuity point of the map
$g\mapsto H(p_g,g)$ and such that if $H(p,f)$ is a homoclinic
class which is a quasi-attractor for $f$, then $H_g = H(p_g,g)$ is
also a quasi-attractor for every $g\in \U\cap \cG$. Also, we can
assume that for every $g\in \U\cap \Res$, the minimal index of
$H_g$ is $\alpha$.

The class admits a dominated splitting of the form $T_HM=E \oplus
F$ with $\dim E=\alpha$ (see Theorem
\ref{Teorema-LejosTangenciasSylvain}). We assume that the
subbundle $E$ is not uniformly contracted. This allows us to use
Theorem \ref{Teorema-Tricotomia}.

Since the minimal index of $H$ is $\alpha$, option $1)$ of the
theorem cannot occur.

We shall prove that option $3)$ implies option $2)$. That is, we
shall prove that if $E$ is not uniformly contracted, then we are
in option $2)$ of Theorem \ref{Teorema-Tricotomia}.

This is enough to prove the theorem since if we apply Theorem
\ref{Teorema-Tricotomia} to $E'$ given by option $2)$ we get that
since $\dim E' = \alpha-1$ option $1)$ and $2)$ cannot happen, and
since option $3)$ implies option $2)$ we get that $E'$ must be
uniformly contracted thus proving Theorem \ref{teoY} (observe that
the statement on the Hausdorff convergence of periodic orbits to
the class can be deduced from option $2)$ also by using Frank's
Lemma).

\begin{claim}
To get option $2)$ in Theorem \ref{Teorema-Tricotomia} is enough
to find one periodic orbit of index $\alpha$ in $H$ which is weak
(that is, it has one Lyapunov exponent in $(-\delta,0)$).
\end{claim}

\dem  This follows using the fact that being far from tangencies
there is a dominated splitting in the orbit given by Theorem
\ref{Teorema-LejosTangenciasSylvain} with a one dimensional
central bundle associated with the weak eigenvalue. Let $\mathcal
O$ be the weak periodic orbit, so we have a dominated splitting of
the form $T_{\mathcal O}M =E^s \oplus E^c \oplus E^u$.

Using transitions (see Proposition \ref{Proposicion-transitions}),
we can find a dense subset in the class of periodic orbits that
spend most of the time near the orbit we found, say, for a small
neighborhood $U$ of $\mathcal O$, we find a dense subset of
periodic points $p_n$ such that the cardinal of the set $\{ i \in
\Z \cap [0, \pi(p_n)-1] \ : \ f^i(p_n) \in U \}$ is bigger than
$(1-\eps)\pi(p_n)$.

Since we can choose $U$ to be arbitrarily small, we can choose
$\eps$ so that the orbits of all $p_n$ admit the same dominated
splitting (this can be done using cones for example) and maybe by
taking $\eps$ smaller to show that $p_n$ are also weak periodic
orbits.

\finobs

It rests to prove that option $3)$ implies the existence of weak
periodic orbits in the class. To do this, we shall discuss
depending on the structure of the partially hyperbolic splitting
using the classification given in \cite{Crov-CentralModels}. There
are 3 different cases according to the possibilities given by
Proposition \ref{Proposition-CentralModels}.

We have a compact invariant set $K\en H$ with minimal dynamics and
which has a partially hyperbolic structure of the form
$T_KM=E^s\oplus E^c\oplus E^u$ where $\dim E^c=1$ and $\dim
E^s<\dim E$. Also, any measure supported on $K$ has zero Lyapunov
exponent along $E^c$. We shall assume that the dimension of $E^s$
is minimal in the sense that every other compact invariant
$\tilde{K}$ satisfying the same properties as $K$ satisfies that
$\dim(E^s_{\tilde K}) \geq \dim(E^s_K)$ (this will be used only
for Case C)).

\subsubsection*{Case A): There exists a chain recurrent central segment. $K$ has type $(R)$}

Assume that the set $K\en H$ admits a chain recurrent central
segment. That is, there exists a curve $\gamma$ tangent to $E^c$
in a point of $K$, which is contained in $H$ and such that
$\gamma$ is contained in a compact, invariant, chain transitive
set in $U$, a small neighborhood of $K$.

In this case, the results of \cite{Crov-Birth} (Addendum 3.9)
imply that there are periodic orbits in the same chain recurrence
class as $K$ (i.e. $H$) with index $\dim E^s \leq \alpha-1$, a
contradiction.

\subsubsection*{Case B): $K$ has type $(N)$, $(H)$ or $(P_{SN})$}

If $K$ has type $(H)$, one can apply Proposition 4.4 of
\cite{Crov-CentralModels} which implies that there is a weak
periodic orbit in $H$ giving option $2)$ of Theorem B.1.

Cases $(N)$ and $(P_{SN})$ give a family of central curves
$\gamma_x$ $\forall x \in K$ (tangent to $E^c$, see
\cite{Crov-CentralModels}) which satisfy that $f(\gamma_x)\en
\gamma_{f(x)}$. It is not difficult to see that there is a
neighborhood $U$ of $K$ such that for every invariant set in $U$
the same property will be satisfied (see remark 2.3 of
\cite{Crov-Birth}).

Consider a set $\hat{K}=K\cup \bigcup_n \mathcal O_n$ where
$\mathcal O_n$ are close enough periodic orbits converging in
Hausdorff topology to $K$ (these are given, for instance, by
Theorem \ref{Teorema-ShadowingCrivisier}) which we can suppose are
contained in $U$.

So, since for some $x\in K$, we have that $\cF^{uu}_{loc}(x)$ will
intersect $W^{cs}_{loc}(p_n)$ in a point $z$ (for a point $x$, the
local center stable set, $W^{cs}_{loc}(x)$ is the union of the
local strong stable leaves of the points in $\gamma_x$).

Since the $\omega$-limit set of $z$ must be a periodic point (see
Lemma 3.13 of \cite{Crov-Birth}) and since $H$ is Lyapunov stable
we get that there is a periodic point of index $\alpha$ which is
weak, or a periodic point of smaller index in $H$ which gives a
contradiction.

\subsubsection*{ Case C): $K$ has type ($P_{UN}$) or ($P_{SU}$)}

One has a minimal set $K$ which is contained in a homoclinic class
which is a quasi-attractor and it admits a partially hyperbolic
splitting with one dimensional center with zero exponents and type
$(P_{UN})$ or $(P_{SU})$.

This gives that given a compact neighborhood $U$ of $K$, there
exists a family of $C^1$-curves $\gamma_x:[0,1]\to U$
($\gamma_x(0)=x$) tangent to the central bundle such that
$f^{-1}(\gamma_x([0,1]))\en \gamma_{f^{-1}(x)}([0,1))$. This
implies that the preimages of these curves remain in $U$ for past
iterates and with bounded length.

They also verify that the chain unstable set of $K$ restricted to
$U$ (that is, the set of points that can be reached from $K$ by
arbitrarily small pseudo orbits contained in $U$) contains these
curves. Since $H$ is a quasi-attractor, this implies that these
curves are contained in $H$.

Assume we could extend the partially hyperbolic splitting from $K$
to a dominated splitting $T_{K'}M=E_1\oplus E^c \oplus E_3$ in a
chain transitive set $K'\en H$ containing $\gamma_x([0,t))$ for
some $x\in K$ and for some $t\in (0,1)$.

Since the orbit of $\gamma_x([0,1])$ remains near $K$ for past
iterates, we can assume (by choosing $U$ sufficiently small) that
the bundle $E_3$ is uniformly expanded there. So, there are
uniformly large unstable manifolds for every point in
$\gamma_x([0,1])$ and are contained in $H$.

If we prove that $E_1$ is uniformly contracted in all $K'$, since
we can approach $K'$ by weak periodic orbits, we get weak periodic
in the class since its strong stable manifold (tangent to $E_1$)
will intersect $H$.

To prove this, we use that for $K$ the dimension of $E^s$ is
minimal. So $E_1$ must be stable, otherwise, we would get that,
using Theorem \ref{Teorema-Tricotomia} again, there is a partially
hyperbolic minimal set inside $K'$ with stable bundle of dimension
smaller than the one of $K$, a contradiction.

The fact that we can extend the dominated splitting and approach
the point $y$ in $\gamma_x((0,1))$ by weak periodic points is
given by  Lemma \ref{encontrarpuntos} below.

\begin{lema}\label{encontrarpuntos}
There exists a residual subset $\Res' \en \Diff^1(M) \backslash
\overline{\Tang}$ such that every $f\in \Res'$ verifies the
following. Given a compact invariant set $K$ such that

\begin{itemize}
\item[-]$K$ is a chain transitive set. \item[-]$K$ admits a
partially hyperbolic splitting $T_KM = E^s \oplus E^c \oplus E^u$
where $E^s$ is uniformly contracting, $E^u$ is uniformly expanding
and $\dim E^c= 1$.
 \item[-] Any invariant measure supported in $K$ has zero Lyapunov exponents along $E^c$.
\end{itemize}

Then, for every $\delta>0$, there exists $U$, a neighborhood of
$K$ such that for every $y\in U$ satisfying:
\begin{itemize}
\item[-] $y$ belongs to the local chain unstable set $pW^u(K,U)$
of $K$ (that is, for every $\eps>0$ there exists an $\eps-$pseudo
orbit from $K$ to $y$ contained in $U$)
    \item[-] $y$ belongs to the chain recurrence class of $K$
\end{itemize}

\noindent we have that there exist $p_n \to y$, periodic points,
such that:

\begin{itemize}
\item[-] The orbit $\cO(p_n)$ of the periodic point $p_n$ has its
$\dim E^s+1$ Lyapunov exponent contained in $(-\delta,\delta)$.
\item[-] For large enough $n_0$, if $\tilde K = K \cup
\overline{\bigcup_{n>n_0} \mathcal O(p_n)}$, then we can extend
the partially hyperbolic splitting to a dominated splitting of the
form $T_{\tilde K}M =E_1 \oplus E^c \oplus E_3$.
\end{itemize}
\end{lema}

The following lemma allows to conclude the proof as we mentioned
before. Its proof is postponed to subsection \ref{Lema}.

This concludes the proof of Theorem \ref{teoY} \lqqd

\subsection{Proof of Lemma \ref{encontrarpuntos}}\label{Lema}

We shall first prove a perturbation result and afterwards we shall
deduce Lemma \ref{encontrarpuntos} with a standard Baire argument.
One can compare this lemma with Lemma 3.2 of
\cite{Yang-LyapunovStable} which is a slightly weaker version of
this. See \cite{Crov-Hab} Chapter 9 for possible applications.

\begin{lema}\label{perturbacion}
There exists a residual subset $\Res \en \Diff^1(M)$ such that
every $f\in \Res$ verifies the following. Given:
\begin{itemize}
\item[-] $K$ a compact chain transitive set. \item[-] $U$ a
neighborhood of $K$ and $y\in U$ verifying that $y$ is contained
in the local chain unstable set $pW^u(K,U)$ of $K$ and in the
chain recurrence class of $K$. \item[-] $\U$ a $C^1$-neighborhood
of $f$.
\end{itemize}
Then, there exists $l>0$ such that, for every $\nu>0$ and $L>0$ we
have $g\in \U$ with a periodic orbit $\cO$ with the following
properties:
\begin{itemize}
\item[-] There exists $p_1\in \cO$ such that
$d(f^{-k}(y),g^{-k}(p_1))<\nu$ for every $0\leq k \leq L$.
\item[-] There exists $p_2\in \cO$ such that $\cO \backslash
\{p_2, \ldots, g^l(p_2)\} \en U$.
\end{itemize}
\end{lema}

\dem{\!\!} The argument is similar as the one in section 1.4 of
\cite{Crov-CentralModels}. We must show that after an arbitrarily
small perturbation, we can construct such periodic orbits.

Consider a point $y$ as above.  We can assume that $y$ is not
chain recurrent in $U$, otherwise $y$ we would be accumulated by
periodic orbits contained in $U$ (see Theorem
\ref{Teorema-ShadowingCrivisier}) and that would conclude without
perturbing.

For every $\eps>0$ we consider an $\eps-$pseudo orbit $Y_\eps =
(z_0 , z_1, \ldots, z_n)$ with $z_0 \in K$ and $z_n=y$ contained
in $U$. Using that $y$ is not chain recurrent in $U$, we get that
for $\eps$ small enough we have that $B_\nu(y)\cap Y_\eps = \{y\}$
where $B_\nu(y)$ is the ball of radius $\nu$ and $\nu$ is small.
So if we consider a Hausdorff limit of the sequence $Y_{1/n}$ we
get a compact set $Z^{-}$ for which $y$ is isolated and such that
it is contained in the chain unstable set of $K$ restricted to
$U$. Notice that $Z^-$ is backward invariant.

If we now consider the pair $(\Delta^{-}, y)$ where $\Delta^{-} =
Z^{-}\backslash \{f^{-n}(y)\}_{n\geq 0}$ we get a pair as the one
obtained in Lemma 1.11 of \cite{Crov-CentralModels} where $y$
plays the role of $x^{-}$. Notice that $\Delta^{-}$ is compact and
invariant.

Now, we consider $\tilde U \en U$ a small neighborhood of $K\cup
\Delta^-$ such that $y\notin \tilde U$. Take $x \in H \cap U^c$
where $H$ is the chain recurrence class of $K$.

Consider $X_\eps =(z_0, \ldots z_n)$ an $\eps-$pseudo orbit such
that $z_0=x$ and $z_n\in K$. Take $z_j$ the last point of $X_\eps$
outside $\tilde U$. Since we chose $\tilde U$ small we have $z_j
\in U\backslash \tilde U$. We call $\tilde X_\eps$ to $(z_j,\ldots
z_n)$.

Consider $Z^+$ the Hausdorff limit of the sequence $\tilde
X_{1/n}$ which will be a forward invariant compact set which
intersects $U\backslash \tilde U$. Since $y$ is not chain
recurrent in $U$ we have that $y\notin Z^+$.

We consider a point $x^+ \in Z^+ \cap U\backslash \tilde U$. This
point satisfies that one can reach $K$ from $x^+$ by arbitrarily
small pseudo orbits. We get that the future orbit of $x^+$ does
not intersect the orbit of $y$.

Consider $\U$ a neighborhood of $f$ in the $C^1$ topology.
Hayashi's connecting lemma (Theorem \ref{Teorema-ConnectingLemma})
gives us $N>0$ and neighborhoods $W^+ \en \hat{W}^+$ of $x^+$ and
$W^- \en \hat{W}^-$ of $y$ which we can consider arbitrarily small
so, we can suppose that

\begin{itemize}
\item[-] All the iterates $f^i(\hat{W}^+)$ and $f^{j}(\hat{W}^-)$
for $0\leq i,j \leq N$ are pairwise disjoint. \item[-] The
iterates $f^i(\hat{W}^+)$ for $0\leq i \leq N$ are disjoint from
the past orbit of $y$.
\end{itemize}

Since there are arbitrarily small pseudo orbits going from $y$ to
$x^+$ contained in $H$ and $f$ is generic, Theorem
\ref{Teorema-ShadowingCrivisier} $\{x_0, \ldots, f^l(x_0) \}$ in a
small neighborhood of $H$ and such that $x_0 \in W^-$ and
$f^l(x_0)\in W^+$.

The same argument gives us an orbit $\{x_1, \ldots, f^k(x_1)\}$
contained in $U$ such that $x_1 \in W^+$ and $f^k(x_1) \in W^{-}$.
In fact, we can choose it so that $d(f^{-i}(y),f^{k-i}(x_1))<\nu$
for $0\leq i \leq L$ (this can be done using uniform continuity of
$f^{-1}, \ldots, f^{-L}$ and choosing $f^k(x_1)$ close enough to
$y$).

Using Hayashi's connecting lemma (Theorem
\ref{Teorema-ConnectingLemma}) we can then create a periodic orbit
$\mathcal O$ for a diffeomorphism $g \in \U$ which is contained in
$\{x_0,\ldots, f^l(x_0)\} \cup \hat{W}^+ \cup \ldots \cup
f^N(\hat{W^+}) \cup \{x_1, \ldots , f^k(x_1)\} \cup \hat{W}^- \cup
\ldots \cup f^N(\hat{W}^-)$ (in the proof of Proposition 1.10 of
\cite{Crov-CentralModels}) is explained how one can compose two
perturbations in order to close the orbit).

Notice that from how we choose $\hat{W}^+$ and $\hat{W}^-$ and the
orbit $\{x_1, \ldots, f^k(x_1)\}$ we get that the periodic orbit
we create with the connecting lemma satisfies that $d(f^{-i}(y),
g^{-i}(p_n)) < \nu$ for $0 \leq i \leq L$ and some $p_n$ in the
orbit. This is because $\{x_1, \ldots, f^{k-1}(x_1)\}$ does not
intersect $\hat{W}^- \cup \ldots \cup f^N(\hat{W}^-)$ and
$\{f^{k-L}(x_1), \ldots f^k(x_1)\}$ does not intersect $\hat{W}^+
\cup \ldots \cup f^N(\hat{W^+})$ (in fact, this gives that the
orbit of $p_n$ for $g$ contains $\{f^{k-L}(x_1), \ldots
f^k(x_1)\}$).

Also, since $\hat{W}^+ \cup \ldots \cup f^N(\hat{W^+}) \cup \{x_1,
\ldots, f^k(x_1)\}  \en U$ we get that, except for maybe $l$
consecutive iterates of one point in the resulting orbit, the rest
of the orbit is contained in $U$. This concludes the proof.

\lqqd

\demo{of Lemma \ref{encontrarpuntos}}  Take $x \in M$, an let $m,
e, t \in \N$. We consider $\U(m,e, t, x)$ the set of $C^1$
diffeomorphisms $g\in \Diff^1(M)$ with a periodic orbit $\cO$
satisfying:

\begin{itemize}
\item[-] $\cO$ is hyperbolic. \item[-] Its $e+1$ Lyapunov exponent
of $\cO$ is contained in $(-1/m,1/m)$.
 \item[-] There exists $p\in \cO$ such that $d(p,x) < 1/t$.
\end{itemize}

This set is clearly open in $\Diff^1(M)$. Let $\{x_s\}$ be a
countable dense set of $M$. We define $\Res_{m,e, t, s} =  \U(m,e,
t, x_s) \cup \Diff^1(M)\backslash \overline{\U(m,e,t , x_s)}$
which is open and dense by definition. Consider $\Res_1 =
\bigcap_{m,e, t, s} \Res_{m,e, t, s}$ which is residual. Finally,
taking $\Res$ as in Lemma \ref{perturbacion}, we consider

 $$\Res'=(\Res \cap \Res_1) \backslash \overline{\Tang}$$

\smallskip

Consider $K$ compact chain transitive and with a partially
hyperbolic splitting $T_KM = E^s \oplus E^c \oplus E^u$ with $\dim
E^c=1$. We assume that any invariant measure supported in $K$ has
Lyapunov exponent equal to zero.

Choose $\delta>0$ small enough. Since $f$ is far from tangencies,
Theorem \ref{Teorema-LejosTangenciasSylvain} gives us that every
periodic orbit having its $\dim E^s +1$ Lyapunov exponent in
$(-\delta, \delta)$ admits a dominated splitting $E_1 \oplus E^c
\oplus E_3$ with uniform strength (that is, if there is a set
$\{O_n\}$ of periodic orbits with their $\dim E^s+1$ Lyapunov
exponents in $(-\delta,\delta)$, then the dominated splitting
extends to the closure).

We choose $U_1$, an open neighborhood of $K$ such that every
invariant measure supported in $U_1$ has its Lyapunov exponents in
$(-1/2m_0, 1/2m_0)$ where $1/m_0 <\delta$.

We can assume that $U_1$ verifies that there are $Df$ invariant
cones $\cE^{uu}$ and $\cE^{cu}$ around $E^u$ and $E^c\oplus E^u$
respectively, defined in $U_1$. Similarly, there are in $U_1$,
$Df^{-1}$ invariant cones $\cE^{ss}$,  $\cE^{cs}$ around $E^s$ and
$E^s\oplus E^c$ respectively.

We can assume, by choosing an adapted metric (see \cite{HPS} or
\cite{GourmelonAdaptada}), that for every $v\in \cE^{ss}$ we have
$\|Df v \|< \lambda \|v\|$ and for every $v\in \cE^{uu}$ we have
$\| Df v \|>\lambda^{-1} \|v\|$ for some $\lambda<1$. There exists
$\U_1$, a $C^1$-neighborhood of $f$ such that for every $g\in
\U_1$ the properties above remain true.

Given $U$ neighborhood of $K$ such that $\overline U \en U_1$ we
have that any $g$ invariant set contained in $U$ admits a
partially hyperbolic splitting.

We now consider $y\in pW^u(K,U)$ which is contained in the chain
recurrence class of $K$.

\begin{claim}\label{AfirmacionConos} Given $t$, for any $x_s$ with $d(x_s,y)<1/2t$ we get that  $f\in \U(m_0,\dim E^s, t, s) $.
\end{claim}

\dem Since $f$ is in $\Res_1'$ it is enough to show that every
neighborhood of $f$ intersects $\U(m_0,\dim E^s, t, s)$.

Choose a neighborhood $\U$ of $f$ and consider $\U_0 \en \U$ given
by Franks' Lemma (Theorem \ref{Teorema-FranksLemma}) such that we
can perturb the derivative of some $g\in \U_0$ in a finite set of
points less than $\xi$ and obtain a diffeomorphism in $\U$.

For $\U_0$, Lemma \ref{perturbacion} gives us a value of $l<0$
such that for any $L>0$ and there exists $g_L \in \U_0$ and a
periodic orbit $\cO_L$ of $g_L$ such that there is a point $p_1\in
\cO_L$ satisfying that $d(g^{-i}(p_1), f^{-i}(y))<1/2t$ ($0\leq i
\leq L$) and a point $p_2\in \cO_L$ such that $\cO_L \backslash \{
p_2 , \ldots, g^l(p_2)\}$ is contained in $U$. We can assume that
$\cO$ is hyperbolic.

We must perturb the derivative of $\cO_L$ less than $\xi$ in order
to show that the $\dim E^s +1$ Lyapunov exponent is in
$(-1/m_0,1/m_0)$.

Notice that if we choose $L$ large enough, we can assume that the
angle of the cone $Dg^L(\cC^{\sigma}(g^{-L}(p_2))$ is arbitrarily
small ($\sigma= uu, cu$). In the same way, we can assume that the
angle of the cone $Dg^{-L}(\cC^{\tilde \sigma}(g^{L+l}(p_2)$ is
arbitrarily small ($\tilde \sigma=cs, ss$ respectively).

Since $l$ is fixed, we get that for $p \in \cO_L \cap U^c$ (if
there exists any, we can assume it is $p_2$), it is enough to
perturb less than $\xi$ the derivative in order to get the cones
$Dg^L(\cC^{\sigma}(g^{-L}(p_2))$ and $Dg^{-L-l}(\cC^{\tilde
\sigma}(g^{L+l}(p_2)$ transversal (for $\sigma=uu, cu$ and $\tilde
\sigma=cs, ss$ respectively). This allows us to have a well
defined dominated splitting above $\cO_L$ (which may be of very
small strength) which in turn allow us to define the $\dim E^s+1$
Lyapunov exponent. Since the orbit $\cO_L$ spends most of the time
inside $U$, and any measure supported in $U$ has its center
Lyapunov exponent in $(-1/2m_0,1/2m_0)$ we get the desired
property.

\finobs

Taking $t=n\to \infty$ and using Theorem
\ref{Teorema-LejosTangenciasSylvain}, we get a sequence of
periodic points $p_n \to y$ such that if $\cO(p_n)$ are their
orbits, the set $\tilde K = K \cup \overline{ \bigcup_n \cO(p_n)}$
admits a dominated splitting $T_{\tilde K}M = E_1 \oplus E^c
\oplus E_3$ extending the partially hyperbolic splitting.

This concludes the proof of Lemma \ref{encontrarpuntos}.

\lqqd

\subsection{Application: Bi-Lyapunov stable homoclinic classes}\label{SubSection-BiLyapunovClasses}

In \cite{ABD} the following conjecture was posed (it also appeared
as Problem 1 in \cite{BC})

\begin{conj}[\cite{ABD}]\label{Conjetura-Interior}
There exists a residual set $\cG$ of $\Diff^1(M)$ of
diffeomorphisms such that if $f\in \cG$ admits a homoclinic class
with nonempty interior, then the diffeomorphism is transitive.
\end{conj}

Some progress has been made towards the proof of this conjecture
(see \cite{ABD},\cite{ABCD} and \cite{PotSambarino}), in
particular, it has been proved in \cite{ABD} that isolated
homoclinic classes as well as homoclinic classes admitting a
strong partially hyperbolic splitting
 verify the conjecture. Also, they proved that a homoclinic class with non empty
interior must admit a dominated splitting (see Theorem 8 in
\cite{ABD}).

In \cite{ABCD} the conjecture was proved for surface
diffeomorphisms, other proof for surfaces (which does not use the
approximation by $C^2$ diffeomorphisms) can be found in
\cite{PotSambarino} where the codimension one case is studied.

Also, from the work of Yang (\cite{Yang-LyapunovStable}, see also
subsection \ref{SubSection-QAlejosTangencias} and Proposition
\ref{lejosdetangencias} below) one can deduce the conjecture in
the case $f$ is $C^1$-generic and far from homoclinic tangencies.

When studying some facts about this conjecture, in \cite{ABD} it
was proved that if a homoclinic class of a $C^1$-generic
diffeomorphism has nonempty interior then this class should be
bi-Lyapunov stable. In fact, in \cite{ABD} they proved that
isolated and strongly partially hyperbolic bi-Lyapunov stable
homoclinic classes for generic diffeomorphisms are the whole
manifold.

This concept is a priori weaker than having nonempty interior and
it is natural to ask the following question.

\begin{quest}[Problem 1 of \cite{BC}]\label{Question-BiLyapunov}
Is a bi-Lyapunov stable homoclinic class of a generic diffeomorphism necessarily the whole manifold?
\end{quest}

It is not difficult to deduce from \cite{BC} that, for generic
diffeomorphisms, a chain recurrence class with non empty interior
must be a homoclinic class (see Corollary
\ref{Corolario-BonattiCrovisier}), thus, the answer to Conjecture
\ref{Conjetura-Interior} must be the same for chain recurrence
classes and for homoclinic classes.

However, we know that Question \ref{Question-BiLyapunov} admits a
negative answer if posed for general chain recurrence classes.
Bonatti and Diaz constructed (see \cite{BDihes}) open sets of
diffeomorphisms in every manifold of dimension $\geq 3$ admitting,
for generic diffeomorphisms there, uncountably many bi-Lyapunov
stable chain recurrence classes which in turn have no periodic
points.

Although this may suggest a negative answer for Question
\ref{Question-BiLyapunov} we present here some results suggesting
an affirmative answer. In particular, we prove that the answer is
affirmative for surface diffeomorphisms, and that in three
dimensional manifold diffeomorphisms the answer must be the same
as for Conjecture \ref{Conjetura-Interior}.

The main reason for which the techniques in \cite{ABCD} (or in
\cite{PotSambarino}) are not able to answer Question
\ref{Question-BiLyapunov} for surfaces, is because differently
from the case of homoclinic classes with interior, it is not so
easy to prove that bi-Lyapunov stable classes admit a dominated
splitting (in fact, the bi-Lyapunov stable chain recurrence
classes constructed in \cite{BDihes} do not admit any). However,
as a consequence of Theorem
\ref{Teorema-QATieneDescomposicionDominada} we will have this
property automatically.

\begin{teo}\label{descdom}
For every $f$ in a residual subset $\cG_1$ of $\Diff^1(M)$, if $H$
is a bi-Lyapunov stable homoclinic class for $f$, then, $H$ admits
a dominated splitting. Moreover, it admits at least one dominated
splitting with index equal to the index of some periodic point in
the class.
\end{teo}

\dem This follows directly from Theorem
\ref{Teorema-QATieneDescomposicionDominada} applied either to $f$
or $f^{-1}$.

\lqqd

This theorem solves affirmatively the second part of Problem 5.1
in \cite{ABD}. We remark that Theorem \ref{descdom} does not imply
that the class is not accumulated by sinks or sources. Also, we
must remark that the theorem is optimal in the following sense, in
\cite{BV} an example is constructed of a robustly transitive
diffeomorphism (thus bi-Lyapunov stable) of $\T^4$ admitting only
one dominated splitting (into two two-dimensional bundles) and
with periodic points of all possible indexes for saddles.

We recall now that a compact invariant set $H$ is \emph{strongly
partially hyperbolic} if it admits a three ways dominated
splitting $T_HM= E^s \oplus E^c \oplus E^u$, where $E^s$ is non
trivial and uniformly contracting and $E^u$ is non trivial and
uniformly expanding.

In the context of Question \ref{Question-BiLyapunov} it was shown
in \cite{ABD} that generic bi-Lyapunov stable homoclinic classes
admitting a strongly partially hyperbolic splitting must be the
whole manifold. Thus, it is very important to study whether the
extremal bundles of a dominated splitting must be uniform.

As a consequence of Theorem \ref{Teorema-ExtremalDimension1} we
get the following easy corollaries.

\begin{cor}\label{Cor1} Let $H$ be a bi-Lyapunov stable homoclinic class
for a $C^1$-generic diffeomorphism $f$ such that $T_HM= E^1 \oplus
E^2 \oplus E^3$ is a dominated splitting for $f$ and
$\dim(E^1)=\dim(E^3)=1$. Then, $H$ is strongly partially
hyperbolic and $H=M$.
\end{cor}

\dem The class should be strongly partially hyperbolic by applying
Theorem \ref{Teorema-ExtremalDimension1} applied to both $f$ and
$f^{-1}$.  Corollary 1 of \cite{ABD} (page 185) implies that
$H=M$. \lqqd

We say that a hyperbolic periodic point $p$ is \emph{far from
tangencies} if there is a neighborhood of $f$ such that there are
no homoclinic tangencies associated to the stable and unstable
manifolds of the continuation of $p$. The tangencies are of index
$i$ if they are associated to a periodic point of index $i$, that
is, its stable manifold has dimension $i$.

\begin{cor}\label{corolariolejostang} Let $H$ be a bi-Lyapunov stable homoclinic class
for a $C^1$-generic diffeomorphism $f$ which has a periodic point
$p$ of index $1$ and a periodic point $q$ of index $d-1$ and such
that $p$ and $q$ are far from tangencies . Then, $H=M$.
\end{cor}

\dem Using Theorem \ref{Teorema-ABCDW} we are in the hypothesis of
Corollary \ref{Cor1}
\lqqd

In low dimension, our results have some stronger implications, we
obtain:

\begin{teo}\label{Teorema-BiLyapunovSuperficies}
Let $f$ be a $C^1$-generic surface diffeomorphism having a bi-Lyapunov stable homoclinic class $H$. Then, $H=\TT^2$ and $f$ is Anosov.
\end{teo}

\dem From Theorem \ref{descdom} and Theorem
\ref{Teorema-ExtremalDimension1} we deduce that $H$ must be
hyperbolic. Using the interior and the local product structure, we
obtain that $H=M$ (see \cite{ABD}) and thus $f$ is Anosov.

Now, by Franks' theorem (\cite{FranksAnosov}) $M$ must be $\TT^2$
and $f$ conjugated to a linear Anosov diffeomorphism.

\lqqd

\begin{obs}
Notice that for proving this Theorem we do not need to use the results of \cite{PujSamAnnals1} which involve $C^2$ approximations.
\finobs
\end{obs}

The following proposition gives a complete answer to Problem 5.1
of \cite{ABD} in dimension $3$.

\begin{prop}\label{dimension3}
Let $H$ be a bi-Lyapunov stable homoclinic class for a $C^1$-generic diffeomorphism in dimension 3. Then, $H$ has nonempty interior.
\end{prop}

\dem Applying Theorem \ref{descdom} one can assume that the class
$H$ admits a dominated splitting of the form $E\oplus F$, and
without loss of generality one can assume that $\dim F=1$.

Theorem \ref{Teorema-ExtremalDimension1} thus implies that $F$ is
uniformly expanded so the splitting is $T_HM = E \oplus E^u$.

Assume first that there exist a periodic point $p$ in $H$ of index
$2$. Thus, this periodic point has a local stable manifold of
dimension $2$ which is homeomorphic to a 2 dimensional disc.

Since the class is Lyapunov stable for $f^{-1}$ the stable
manifold of the periodic point is completely contained in the
class.

Now, using Lyapunov stability for $f$ and the lamination by strong
unstable manifolds given by Theorem
\ref{Teorema-VariedadEstableFuerte} one gets (saturating by
unstable sets the local stable manifold of $p$) that the
homoclinic class contains an open set. This implies the thesis
under this assumption.

So, we must show that if all the periodic points in the class have
index $1$ then the class is the whole manifold. As we have been
doing, using the genericity of $f$ we can assume that there is a
residual subset $\Res$ of $\Diff^1(M)$ and an open set $\U$ of $f$
such that for every $g\in \U \cap \Res$ all the periodic points in
the class have index $1$.

We have 2 situations, on the one hand, we consider the case where
$E$ admits two invariant subbundles, $E=E^1\oplus E^2$, with a
dominated splitting and thus, we get that $E^1$ should be
uniformly contracting (using Theorem
\ref{Teorema-ExtremalDimension1}) proving that the homoclinic
class is the whole manifold (Corollary \ref{Cor1}).

If $E$ admits no invariant subbundles then, using Theorem
\ref{Teorema-BonattiBochi}, we can perturb the derivative of a
periodic point in the class, so that the cocycle over the periodic
point restricted to $E$ has all its eigenvalues contracting. So,
we can construct a periodic point of index $2$ inside the class.

\lqqd

\begin{obs}
It is very easy to adapt the proof of this proposition to get
that: If a bi-Lyapunov stable homoclinic class of a generic
diffeomorphism admits a codimension one dominated splitting,
$T_HM=E\oplus F$ with $\dim F=1$, and has a periodic point of
index $d-1$, then, the class has nonempty interior.\finobs
\end{obs}

Using a Theorems \ref{teoY}  and \ref{Teorema-ExtremalDimension1}
we are able to prove a similar result which is stronger than the
previous corollary but which in turn, has hypothesis of a more
global nature. We say that a diffeomorphism $f$ is \emph{far from
tangencies} if it can not be approximated by diffeomorphisms
having homoclinic tangencies for some hyperbolic periodic point.
Notice that in the far from tangencies context, it is proved in
\cite{Yang-LyapunovStable} that a Lyapunov stable chain recurrence
class must be an homoclinic class.

\begin{prop}\label{lejosdetangencias}
There exists a $C^1$-residual subset of the open set of
diffeomorphisms far from tangencies such that if $H$ is a
bi-Lyapunov stable chain recurrence class for such a
diffeomorphism, then, $H=M$.
\end{prop}

\dem First of all, if the class has all its periodic points with
index between $\alpha$ and $\beta$ we know that it admits a 3 ways
dominated splitting of the form $T_HM=E\oplus G \oplus F$ where
$\dim E=\alpha$ and $\dim F= d-\beta$. This is because we can
apply the result of \cite{Wen-Tangencias1} (see Theorem
\ref{Teorema-LejosTangenciasSylvain}) which says that far from
homoclinic tangencies there is an index $i$ dominated splitting
over the closure of the index $i$ periodic points together with
the fact that index $\alpha$ and $\beta$ periodic points should be
dense in the class since the diffeomorphism is generic.

Now, we will show that $H$ admits a strong partially hyperbolic
splitting. If $E$ is one dimensional, then it must be uniformly
hyperbolic because of Theorem \ref{Teorema-ExtremalDimension1}. If
not, suppose $\dim E >1$ then, if it is not uniform, Theorem
\ref{teoY} implies that it can be decomposed as a uniform bundle
together with a one dimensional central bundle, since $\dim E>1$
we get a uniform bundle of positive dimension.

The same argument applies for $F$ using Lyapunov stability for
$f^{-1}$ so we get a strong partially hyperbolic splitting.

Corollary 1 of \cite{ABD} finishes the proof.

\lqqd



\section{Examples}\label{Section-EjemplosQuasiAtractores}

We have seen in section \ref{Section-AtractoresSuperficies} that a
$C^1$-generic diffeomorphism of a compact surface admits a
hyperbolic attractor. Moreover, we have seen that if a
$C^1$-generic diffeomorphism of a manifold has an attractor, then
this attractor must be volume partially hyperbolic (see Theorem
\ref{TeoremaBonattiDiazPujals}, notice that an attractor is an
isolated chain-recurrence class).

It seems natural to ask the following question (see
\cite{PalisPugh} Problem 26, \cite{MilnorAttractor}, \cite{BDV}
Problems 10.1 and 10.30, \cite{BC}):

\begin{quest} Does a $C^r$-generic diffeomorphism of a compact manifold have an attractor?
\end{quest}

The question traces back to R. Thom and S. Smale who believed in a
positive answer to this question. See \cite{BLY} for a more
complete historical account on this problem.

Recently, and surprisingly (notice that even though it was always
posed as a question, it had always follow up questions in case the
answer was positive), it was shown by \cite{BLY} that this is not
the case. We will review their example in subsection
\ref{SubSection-EjemploBLY}.

Their example posses infinitely many sources accumulating on
certain quasi-attractors (which as we have seen always exist by
Theorem \ref{Teorema-BonattiCrovisier}) so it seemed also natural
to ask whether a $C^r$-generic diffeomorphism has either
attractors or repellers. A very subtle modification of the example
of \cite{BLY} allows one to create examples not having either
attractors nor repellers (see \cite{BS}, we shall extend this
comment in the following section).

It seems natural then, to weaken the notion of attractor in order
to continue searching for the chain-recurrence classes which
capture ``most'' of the dynamics of a ``typical'' diffeomorphism.
We have seen in subsection \ref{SubSection-AttractingSets} several
notions of attracting sets, and we have payed special attention
(specially in this chapter) to \emph{quasi-attractors} (which as
we said, always exist by Theorem \ref{Teorema-BonattiCrovisier}).

However, the notion of quasi-attractor is not that satisfying
since it may even have empty basin (see the examples of
\cite{BDihes}). The following natural question was posed in
\cite{BLY} and seems the ``right'' one:

\begin{quest} Does a $C^r$-generic diffeomorphism admit an essential attractor? and a Milnor attractor?
\end{quest}

As we already mentioned, the first question has been answered in
the affirmative for $C^1$-generic diffeomorphisms far away from
tangencies (\cite{BGLY}).

Of course, candidates for such classes will be quasi-attractors,
specially those which are homoclinic classes. In view of our
Corollary \ref{Corolario-QuasiAtractoresEnDim3} it seems that
there are some tools to attack certain partial questions in
dimension $3$, and regarding at partially hyperbolic
quasi-attractors in dimension $3$ which are homoclinic classes
seems a reasonable way to proceed. We still lack of examples, but
in certain cases, we seem to be acquiring the necessary tools to
understand these particular classes and start constructing a
theory. In this section, we will review this bunch of examples and
we will close the chapter by proposing a direction in order to
understand a certain class of quasi-attractors in dimension $3$.

\subsection{The example of Bonatti-Li-Yang}\label{SubSection-EjemploBLY}

We briefly explain the construction of C.Bonatti, M.Li and D.Yang
in \cite{BLY} and the modifications made by C.Bonatti and
K.Shinohara in \cite{BS}.

They prove the following theorem:

\begin{teo}[\cite{BLY,BS}]\label{Teorema-BLYyBS}
Given a $d$-dimensional manifold $M$ ($d\geq 3$) we have that for
every isotopy class of diffeomorphisms of $M$ there exists an open
set $\cU$ such that for a diffeomorphisms $f$ in a $C^r$-residual
subset of $\cU$ we have that:
\begin{itemize}
\item[-] There is no attractor nor repeller for $f$. \item[-]
Every quasi-attractor $\cQ$ of $f$ is an essential attractor.
\item[-] In every neighborhood of $\cQ$ there are aperiodic
classes.
\end{itemize}
\end{teo}

\subsubsection{Sketch of the construction of dynamics without attractors}

We will outline the construction given in \cite{BLY} of generic
dynamics in an open set of diffeomorphisms without attractors. We
will make the construction in dimension $3$ and in an attracting
solid torus, see \cite{BLY} for details on how to extend the
dynamics into an attracting ball and other dimensions. Notice that
in appendix B of \cite{Franks-Homology} it is shown how to
construct an Axiom A diffeomorphism $C^0$-close to a given one
having only finitely many sinks as attractors (a simple surgery
then allows to obtain that in any isotopy class of diffeomorphisms
of a given manifold of dimension $\geq 3$ one can construct the
desired diffeomorphisms).

We assume from now on some acquaintance with the construction of
Smale's solenoid and Plykin's attractor (see for example
\cite{KH,Robinson,Shub}).

One starts with the solid torus $T= \DD^2 \times S^1$. As in the
solenoid, one ``cuts'' $T$ by a disk of the form $\DD^2 \times
\{x\}$ obtaining a solid cylinder of the form $\DD^2 \times
[-1,1]$.

Then, one ``streches'' the resulting filled cylinder in order to
be able to ``wrap'' the torus $T$ more than once. Then one inserts
the resulting cylinder in $T$ such that it does not autointersects
and ``glues'' again in the same place one had cut.

One can do this in order that the following conditions are
satisfied:

\begin{itemize}
\item[-] The resulting map $f: T \to T$ is an injective $C^\infty$
embedding. \item[-] The image of a disk of the form $\DD^2 \times
\{z\}$ is contained in $\DD^2 \times \{z^2 \}$ where one thinks of
$S^1 \en \CC$ (so that the map $z \mapsto z^2$ is the well known
doubling map, see \cite{KH} section 1.7). \item[-]There exists
$\alpha>0$ such that any vector $v$ in the tangent space of a
point of the form $(p,z)$ whose angle with respect to to $\DD^2
\times \{z\}$ is larger than or equal to $\alpha$ verifies that
the image by $Df$ of $v$ makes angle strictly larger than $\alpha$
with $\DD^2 \times \{z^2\}$ and norm larger than twice the one of
$v$.
\end{itemize}

If one also requires that in the $f$-invariant plaque family given
by $\DD^2 \times \{z\}$ one has uniform contraction, one obtains
the well known Smale's solenoid. C.Bonatti, M.Li and D.Yang have
profited from the fact that there is still some freedom to choose
the dynamics in this invariant plaque family in order to construct
their mentioned example.

Notice that the maximal invariant subset $\Lambda$ of $T$ admits a
partially hyperbolic splitting of the form $T_\Lambda T = E^{cs}
\oplus E^u$ where $E^{cs}$ is tangent to $\DD^2 \times \{z\}$ in
any point of the form $(p,z)$, the subbundle $E^u$ is uniformly
expanding and the angle between $E^{cs}$ and $E^u$ is larger than
$\alpha$.

They demand the following two extra properties which are enough to
show that in a $C^1$-neighborhood of $f$, there will be a residual
subset of diffeomorphisms having no attractor at all in $T$
(notice that there must be at least one quasi-attractor since $T$
is a trapping neighborhood, see Theorem \ref{Teorema-Conley}):

\begin{itemize}
\item[-] Let $1 \in S^1$ be the fixed point of $z\mapsto z^2$. We
demand that the dynamics in $\DD^2 \times  \{1\}$ which is
$f$-invariant has a unique fixed point $p$ which will be
hyperbolic and attracting and has complex eigenvalues. \item[-]
Let $z_0$ be a periodic point of $z\mapsto z^2$ of period $k>0$.
We have that $f^k(\DD^2 \times \{z_0\}) \en \DD^2 \times \{z_0\}$.
We will demand that the dynamics of $f^k$ in that disk is the one
of the Plykin attractor on the disk and contains a periodic point
$q$ whose determinant restricted to $E^{cs}$ is larger than $1$.
\end{itemize}

The first property allows one to show the following:

\begin{lema}\label{Lemma-ArgumentoBLYunicoQA}
There exists a $C^1$-neighborhood $\cU$ of $f$ such that every
$g\in \cU$ has a unique quasi-attractor in $T$ which contains the
homoclinic class of the continuation of $p$.
\end{lema}

\dem Notice that the unstable manifold of every point in the
maximal invariant set of $g$ inside $T$ must intersect the stable
manifold of the continuation of $p$. This implies that the closure
of the unstable manifold of $p$ (and therefore its homoclinic
class) must be contained in every quasi-attractor inside $T$. See
Lemma \ref{unicocasiatractor} for more details.

\lqqd

With the second property we can show that generic diffeomorphisms
in a neighborhood of $f$ cannot have attractors (recall that an
attractor is an isolated quasi-attractor) inside $T$:

\begin{teo}[\cite{BLY}]\label{Teorema-EJEMPLOBLY}
There exists a $C^1$-neighborhood $\cU$ of $f$ such that for every
$C^1$-generic diffeomorphism $g\in \cU$ the homoclinic class of
the continuation of $p$ is contained in the closure of the set of
sources of $g$, in particular, $g$ has no isolated
quasi-attractors.
\end{teo}

\esbozo{} The fact that $p$ has complex stable eigenvalues implies
that $E^{cs}$ admits no sub-dominated splitting.

We now show that the continuation of $q$ for diffeomorphisms in a
neighborhood of $p$ belongs to the chain-recurrence class of $p$.
Indeed, the unstable manifold of $p$ intersects the basin of
attraction of the Plykin attractor (which is contained in a
normally hyperbolic disk) and thus, there are arbitrarilly small
pseudo-orbits going from $p$ to the Plykin attractor. Now, since
the chain-recurrence class of $p$ is robustly a quasi-attractor by
the previous Lemma, we get that the Plykin attractor (and thus
$q$) belongs robustly to the unique quasi-attractor.

By Theorem \ref{TeoremaBonattiDiazPujals} we obtain that for
$C^1$-generic diffeomorphisms in a neighborhood of $f$ the
homoclinic class of $p$ (which coincides with the unique
quasi-attractor for generic diffeomorphisms) is contained in the
closure of the set of sources. This concludes.

\lqqd

Indeed, in \cite{BLY} it is proved that the same result holds for
the $C^r$-topology (we refer the reader to the next subsection for
more details).

\subsubsection{Bonatti and Shinohara's result}

In the proof of Theorem \ref{Teorema-EJEMPLOBLY} the non-isolation
of the quasi-attractor follows from the fact that containing a
periodic orbit whose determinant at the period is larger than $1$,
the class cannot be volume hyperbolic, and thus, by Theorem
\ref{TeoremaBonattiDiazPujals} it cannot be isolated.

One could wonder what happens in the event that the
quasi-attractor is indeed volume hyperbolic, so that the creation
of sources is not allowed and the criterium given by Corollary
\ref{Corollary-CriterioDeNoAislado} does not apply.

Bonatti and Shinohara, in \cite{BS} have developed a very subtle
technique which allows them to eject periodic saddles from the
homoclinic class of $p$ even if the class is volume hyperbolic.

Using this technique they are able to construct examples which
have no quasi-attractors in $T$ but do not contain sources either
\footnote{Notice that they must change the Plykin attractor since
it forces the existence of sources even if they are not near the
quasi-attractor, see \cite{Plykin}. In order to change this, the
construction slightly more complicated since it must use Blenders,
see subsection \ref{SubSection-Blenders}.} thus without attractors
nor repellers. A surgery argument allows to prove the second
statement of Theorem \ref{Teorema-BLYyBS}.

The final item of Theorem \ref{Teorema-BLYyBS} hides some deep
consequences of Bonatti and Shinohara's construction. Indeed, they
show that such a quasi-attractor (they in fact work in a more
general framework which applies to this context) has some
\emph{viral} properties as defined in \cite{BCDG} (see also
\cite{B}). This allows them to show the existence of quite
atypical aperiodic classes (for example, aperiodic classes which
are not transitive) as well as to show that there are, for
$C^1$-generic diffeomorphisms in a neighborhood of the constructed
$f$, uncountably many chain-recurrence classes.

The remarkable feature of their construction is that it relies
heavily on the topology of the intersection of the class with the
center stable plaques, and indeed, in the next subsection we will
show some examples from \cite{PotWildMilnor} which have a quite
opposed behavior\footnote{The examples were obtained almost
simultaneously.}.

Both Bonatti-Li-Yang's example and the modifications made by
Bonatti and Shinohara yield essential attractors, it would be nice
to know if indeed:

\begin{quest} Do these examples admit a Milnor attractor?
\end{quest}

\subsection{Derived from Anosov examples}\label{SubSection-EjemploDA}

After the examples of Bonatti,Li and Yang appeared, it became
clear that the use of Theorem \ref{TeoremaBonattiDiazPujals} could
be a tool yielding examples of dynamics without attractors: It
suffices to construct a quasi-attractor which has periodic points
which are sectionally dissipative in some sense. Also, the
question of understanding ergodic properties of attracting sets
and sets whose topological basin is large in some sense becomes an
important question.

On the other hand, Bonatti-Li-Yang's example was in a sense, a new
kind of \emph{wild} homoclinic class, and the understanding of how
the class is accumulated by other classes became a new challenge.
Hoping to answer partially to this, I was able to construct some
examples whose properties are summarized in the following
statement.

\begin{teo}\label{Teorema-EjemploCarvalho} There exists a $C^1$-open set $\cU$ of $\Diff^r(\TT^3)$ such that:
\begin{itemize}
\item[(a)] For every $f\in \cU$ we have that $f$ is partially
hyperbolic with splitting $TM= E^{cs} \oplus E^u$ and $E^{cs}$
integrates into a $f$-invariant foliation $\cF^{cs}$. \item[(b)]
Every $f\in \cU$ has a unique quasi-attractor $\cQ_f$ which
contains a homoclinic class. \item[(c)] Every chain recurrence
class $R\neq \cQ_f$ is contained in the orbit of a periodic disk
in a leaf of the foliation $\cF^{cs}$. \item[(d)] There exists a
residual subset $\cG^r$ of $\cU$ such that for every $f\in \cG^r$
the diffeomorphism $f$ has no attractors. In particular, $f$ has
infinitely many chain-recurrence classes accumulating on $\cQ_f$.
\item[(e)] For every $f\in \cU$ there is a unique Milnor attractor
$\tilde \cQ \en \cQ_f$. \item[(f)] If $r\geq 2$ then every $f \in
\cU$ has a unique SRB measure whose support coincides with a
homoclinic class. Consequently, $\tilde \cQ$ is a minimal
attractor in the sense of Milnor. If $r=1$, then there exists a
residual subset $\cG_M$ of $\cU$ such that for every $f\in \cG_M$
we have that $\tilde \cQ$ coincides with $\cQ_f$ and is a minimal
Milnor attractor.
\end{itemize}
\end{teo}

The goal of this subsection is to prove this theorem.

By inspection in the proofs, one can easily see that in fact the
construction can be made in higher dimensional torus, however, it
can only be done in the isotopy classes of Anosov diffeomorphisms.
Also, it can be seen that condition (d) can be slightly
strengthened in the $C^1$-topology.

It is clear that this example contrasts with the properties
obtained by Bonatti and Shinohara in \cite{BS}. Indeed, for this
example, the following question remains unsolved (see
\cite{PotWildMilnor} for more discussion on this question):

\begin{quest} Does a $C^1$-generic diffeomorphism in $\cU$ have countably many chain-recurrence classes?
\end{quest}

We remark that the answer to this question for the $C^r$ topology
with $r\geq 2$ is false (see \cite{BDV} section 3 and the
discussion in \cite{PotWildMilnor}).

In \cite{PotWildMilnor} it is also proved that the example is in
the hypothesis of the main theorem of \cite{BuFi} and consequently
admits a unique measure of maximal entropy (concept we will not
define in this thesis but which is self-explanatory).


\subsubsection{Construction of the
example}\label{SectionConstruccion1} In this section we shall
construct an open set $\cU$ of $\Diff^r(\TT^3)$ for $r\geq 1$
verifying Theorem \ref{Teorema-EjemploCarvalho}.

The construction is very similar to the one of Carvalho's example
(\cite{Car}) following \cite{BV} with the difference that instead
of creating a source, we create an expanding saddle. We do not
assume acquaintance of the reader with the referred construction
but we will in some stages point the reader to specific parts we
will not reproduce.

We start with a linear Anosov diffeomorphism $A:\T^3 \to \T^3$
admitting a splitting $E^s\oplus E^u$ where $\dim E^s=2$.

We assume that $A$ has complex eigenvalues on the $E^s$ direction
so that $E^s$ cannot split as a dominated sum of other two
subspaces. For example, the matrix

\[      \left(
            \begin{array}{ccc}
              1 & 1 & 0 \\
              0 & 0 & 1 \\
              1 & 0 & 0 \\
            \end{array}
          \right)
\]

\noindent which has characteristic polynomial $1 + \lambda^2
-\lambda^3$ works since it has only one real root, and it is
larger than one.

Considering an iterate, we may assume that there exists
$\lambda<1/3$ satisfying:

 \[ \|(DA)_{/E^s}\|< \lambda \ \ ; \ \ \|(DA)_{/E^u}^{-1}\|<\lambda  \]

Let $q$ and $r$ be different fixed points of $A$.

Consider $\delta$ small enough such that $B_{6\delta}(q)$ and
$B_{6\delta}(r)$ are pairwise disjoint and at distance larger than
$400\delta$ (this implies in particular that the diameter of
$\TT^3$ is larger than $400\delta$).


Let  $\cE^u$ be a family of closed cones around the subspace $E^u$
of $A$ which is preserved by $DA$ (that is $D_xA(\cE^u(x))\en \Int
(\cE^u(Ax))$). We shall consider the cones are narrow enough so
that any curve tangent to $\cE^u$ of length bigger than $L$
intersects any stable disk of radius $\delta$. Let $\cE^{cs}$ be a
family of closed cones around $E^s$ preserved by $DA^{-1}$.

From now on, $\delta$ remains fixed. Given $\eps>0$
such\footnote{If $K$ bounds $\|A\|$ and $\|A^{-1}\|^{-1}$ then
$\frac{\delta}{10K}$ is enough.} that $\eps \ll \delta$, we can
choose $\nu$ sufficiently small such that every diffeomorphism $g$
which is $\nu$-$C^0$-close to $A$ is semiconjugated to $A$ with a
continuous surjection $h$ which is $\eps$-$C^0$-close to the
identity (this is a classical result on topological stability of
Anosov diffeomorphisms, see \cite{Walters} and Proposition
\ref{PropExisteSemiconjugacion}).

We shall modify $A$ inside $B_\delta(q)$ such that we get a new
diffeomorphism $F:\T^3 \to \T^3$ that verifies the following
properties:

\begin{itemize}
\item[-] $F$ coincides with $A$ outside $B_{\delta}(q)$ and lies
at $C^0$-distance smaller than $\nu$ from $A$. \item[-] The point
$q$ is a hyperbolic saddle fixed point of stable index $1$ and
such that the product of its two eigenvalues with smaller modulus
is larger than $1$. We also assume that the length of the stable
manifold of $q$ is larger than $\delta$. \item[-]$D_xF(\cE^u(x))
\en \Int(\cE^{u}(F(x)))$. Also, for every $w\in \cE^u(x)\setminus
\{0\}$ we have $\|DF_x^{-1}w\|< \lambda \|w\|$. \item[-] $F$
preserves the stable foliation of $A$. Notice that the foliation
will no longer be stable. \item[-] For some small $\beta>0$ we
have that $\|D_xF v\| < (1+\beta) \|v\|$ for every $v$ tangent to
the stable foliation of $A$ preserved by $F$ and every $x$.
\end{itemize}

\begin{figure}[ht]
\begin{center}
\input{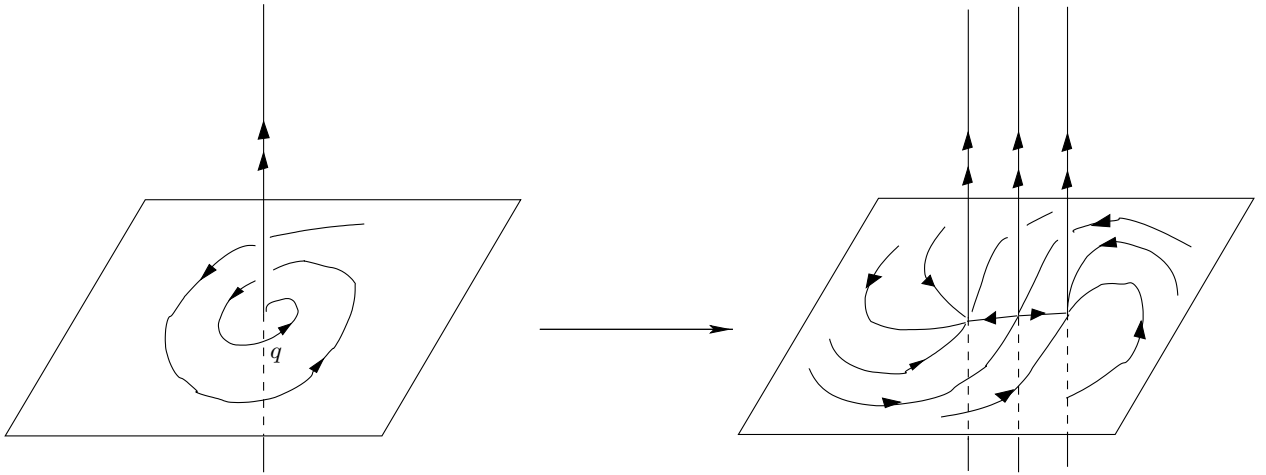}
\caption{\small{Modification of $A$ in a neighborhood of $q$.}}
\label{FiguraPseudoOrbitas}
\end{center}
\end{figure}

This construction can be made using classical methods (see
\cite{BV} section 6). Indeed, consider a small neighborhood $U$ of
$q$ such that $U\en B_{\nu/2}(q)$ such that $U$ admits a chart
$\varphi: U \to \DD^2 \times [-1,1]$ which sends $q$ to $(0,0)$
and sends stable manifolds of $A$ in sets of the form $\DD^2
\times \{t\}$ and unstable ones into sets of the form $\{s\}
\times [-1,1]$. We can modify $A$ by isotopy inside $U$ in such a
way that the sets $\DD^2\times \{t\}$ remain an invariant
foliation but such that the derivative of $q$ becomes the identity
in the tangent space to $\varphi^{-1}(\DD^2\times \{0\})$ which is
invariant and such that the dynamics remains conjugated to the
initial one. At this point, the norm of the images of unit vectors
tangent to the stable foliation of $A$ are not expanded by the
derivative.

Now, one can modify slightly the dynamics in $\varphi^{-1}(\DD^2
\times \{0\})$ in order to obtain the desired conditions on the
eigenvalues of $q$ for $F$. It is not hard to see that for
backward iterates there will be points outside
$\varphi^{-1}(\DD^2\times \{0\})$ which will approach $q$ so one
can obtain the desired length of the stable manifold of $q$ by
maybe performing yet another small modification. All this can be
made in order that the vectors tangent to the stable foliation of
$A$ are expanded by $DF$ by a factor of at most $(1+\beta)$ with
$\beta$ as small as we desire.

The fact that we can keep narrow cones invariant under $DF$ seems
difficult to obtain in view that we made all this modifications.
However, the argument of \cite{BV} (page 190) allows to obtain it:
This is achieved by conjugating the modification with appropriate
homotheties in the stable direction.

The last condition on the norm of $DF$ in the tangent space to the
stable foliation of $A$ seems quite restrictive, more indeed in
view of the condition on the eigenvalues of $q$. This condition
(as well as property (P7) below) shall be only used (and will be
essential) to obtain the ergodic properties of the diffeomorphisms
in the open set we shall construct. Nevertheless, one can
construct such a diffeomorphism as explained above.

There exists a $C^1$-open neighborhood $\cU_1$ of $F$ such that
for every $f \in \cU_1$ we have that:

\begin{itemize}
\item[(P1)] There exists a continuation $q_f$ of $q$ and $r_f$ of
$r$. The point $r_f$ has stable index $2$ and complex eigenvalues.
The point $q_f$ is a saddle fixed point of stable index $1$, such
that the product of its two eigenvalues with smaller modulus is
larger than $1$ and such that the length of the stable manifold is
larger than $\delta$.

\item[(P2)] $D_xf(\cE^u(x)) \en \Int(\cE^{u}(f(x)))$. Also, for
every $w\in \cE^u(x)$ we have $$\|Df_x w\| \geq \lambda \|w\| .$$

\item[(P3)] $f$ preserves a foliation $\cF^{cs}$ which is
$C^0$-close to the stable foliation of $A$. Also, each leaf of
$\cF^{cs}$ is $C^1$-close to a leaf of the stable foliation of
$A$.

\item[(P4)]  For every $x\notin B_\delta(q)$ we have that if $v
\in \cE^{cs}(x)$  then $$\|D_x f v\| \leq \lambda \|v\| .$$ This
is satisfied for $F$ since $F=A$ outside $B_\delta(q)$.

\item[(P5)] There exists a continuous and surjective map $h_f:
\TT^3\to \TT^3$ such that $$h_f \circ f = A\circ h_f$$ \noindent
and $d(h(x), x) < \eps$ for every $x\in \TT^3$.
\end{itemize}

The fact that properties (P1), (P2) and (P4) are $C^1-$robust is
immediate, robustness of (P5) follows from the choice of $\nu$.

Property (P3) holds in a neighborhood of $F$ since $F$ preserves
the stable foliation of $A$ which is a $C^1-$foliation (see
\cite{HPS} chapter 7). The foliation $\cF^{cs}$ will be tangent to
$E^{cs}$ a bidimensional bundle which is $f$-invariant and
contained in $\cE^{cs}$. Other way to proceed in order to obtain
an invariant foliation is to use Theorem 3.1 of \cite{BuFi} of
which all hypothesis are verified here but we shall not state it.

Since the cones $\cE^u$ are narrow and from (P3) one has that:

\begin{itemize}
\item[(P6)] Every curve of length $L$ tangent to $\cE^u$ will
intersect any disc of radius $2\delta$ in $\cF^{cs}$.
\end{itemize}

Finally, there exists an open set $\cU_2 \en \cU_1$ such that for
$f\in \cU_2$ we have:

\begin{itemize}
\item[(P7)]   $\|D_xf v\| \leq (1+\beta) \|v\|$ for every $v\in
\cE^{cs}(x)$ and every $x$.
\end{itemize}

For this examples there exists a unique quasi-attractor for the
dynamics.

\begin{lema}\label{unicocasiatractor}
For every $f\in \U_1$ there exists an unique quasi-attractor
$\cQ_f$. This quasi attractor contains the homoclinic class of
$r_f$, the continuation of $r$.
\end{lema}

\dem We use the same argument as in \cite{BLY}.

There is a center stable disc of radius bigger than $2\delta$
contained in the stable manifold of $r_f$ ((P3) and (P4)). So,
every unstable manifold of length bigger than $L$ will intersect
the stable manifold of $r_f$ ((P6)).

Let $\cQ$ be a quasi attractor, so, there exists a sequence $U_n$,
of neighborhoods of $\cQ$ such that $f(\overline{ U_n}) \en U_n$
and $\cQ=\bigcap_n \overline{U_n}$.

Since $U_n$ is open, there is a small unstable curve $\gamma$
contained in $U_n$. Since $Df$ expands vectors in $\cE^u$ we have
that the length of $f^k(\gamma)$ tends to $+\infty$ as $n\to
+\infty$. So, there exists $k_0$ such that $f^{k_0}(\gamma) \cap
W^s(r_f) \neq \emptyset$. So, since $f(\overline{U_n})\en U_n$ we
get that $U_n \cap W^s(r_f) \neq \emptyset$, using again the
forward invariance of $U_n$ we get that $r_f \in \overline{U_n}$.

This holds for every $n$ so $r_f \in \cQ$. Since the homoclinic
class of $r_f$ is chain transitive, we also get that $H(r_f)\en
\cQ$.

From Conley's theory (subsection
\ref{SubSection-ChainRecurrence}), every homeomorphism of a
compact metric space there is at least one chain recurrent class
which is a quasi attractor. This concludes.
\lqqd

\subsubsection{The example verifies the mechanism of Proposition \ref{ProposicionMecanismo}} We shall
consider $f\in \cU_1$ so that it verifies (P1)-(P6).

Let $\cA^s$ and $\cA^u$ be, respectively, the stable and unstable
foliations of $A$, which are linear foliations. Since $A$ is a
linear Anosov diffeomorphism, the distances inside the leaves of
the foliations and the distances in the manifold are equal in
small neighborhoods of the points if we choose a convenient
metric.

Let $\cA^s_\eta(x)$ denote the ball of radius $\eta$ around $x$
inside the leaf of $x$ of $\cA^s$. For any $\eta>0$, it is
satisfied that $A(\cA^s_\eta(x))\en \cA^s_{\eta/3}(Ax)$ (an
analogous property is satisfied by $\cA^u_\eta(x)$ and backward
iterates).

The distance inside the leaves of $\Fol^{cs}$ is similar to the
ones in the ambient manifold since each leaf of $\cF^{cs}$ is
$C^1$-close to a leaf of $\cA^s$. That is, there exists $\rho
\approx 1$ such that if $x,y$ belong to a connected component of
$\Fol^{cs}(z) \cap B_{10\delta}(z)$ then $\rho^{-1}d_{cs}(x,y)<
d(x,y) <\rho d_{cs}(x,y)$ where $\Fol^{cs}(z)$ denotes the leaf of
the foliation passing through $z$ and $d_{cs}$ the distance
restricted to the leaf.

For $z\in \TT^3$ we define $W^{cs}_{loc}(z)$ (the \emph{local
center stable manifold} of $z$) as the $2\delta$-neighborhood of
$z$ in $\cF^{cs}(z)$ with the distance $d_{cs}$.

Also, we can assume that for some $\gamma < \min
\{\|A\|^{-1},\|A^{-1}\|^{-1}, \delta/10\}$ the plaque
$W^{cs}_{loc}(x)$ is contained in a $\gamma/2$ neighborhood of
$\cA^{s}_{2\delta}(x)$, the disc of radius $2\delta$ of the stable
foliation of $A$ around $x$.

\begin{lema}\label{trapping} We have that $f(\overline{W^{cs}_{loc}(x)}) \en W^{cs}_{loc}(f(x))$.
\end{lema}

\dem Consider around each $x\in \T^3$ a continuous map $b_x:\D^2
\times [-1,1] \to \T^3$ such that $b_x(\{0\} \times
[-1,1])=\cA^u_{3\delta}(x)$ and $b_x(\D^2 \times \{t\})=
\cA^s_{3\delta}(b_x(\{0\}\times\{t\}))$. For example, one can
choose $b_x$ to be affine in each coordinate to the covering of
$\TT^3$.

Thus, it is not hard to see that one can assume also that
$b_x(\frac 1 3 \D^2 \times \{t\} ) =
\cA^s_{\delta}(b_x(\{0\}\times\{t\}))$ and that $b_x(\{y\} \times
[-1/3,1/3])= \cA^u_{\delta}(b_x(\{y\}\times \{0\}))$. Let

  $$B_x = b_x(\D^2 \times [-\gamma/2,\gamma/2]) .$$

We have that $A(B_x)$ is contained in $b_{Ax}( \frac{1}{3}\D^2
\times [-1/2,1/2])$. Since $f$ is $\eps$-$C^0$-near $A$, we get
that $f(B_x) \en b_{f(x)}( \frac{1}{2}\D^2 \times [-1,1])$.

Let $\pi_1: \D^2 \times [-1,1] \to \D^2$ such that $\pi_1(x,t)=x$.
We have that $\pi_1(b_{f(x)}^{-1}(W^{cs}_{loc}(f(x))))$ contains
$\frac{1}{2}\D^2$ from how we chose $\gamma$ and from how we have
defined the local center stable manifolds\footnote{In fact,
$b_{f(x)}^{-1}(W^{cs}_{loc}(g(x))) \cap \frac{1}{2}\D^2\times
[-1,1]$ is the graph of a $C^1$ function from $\frac{1}{2}\D^2$ to
$[-\gamma/2,\gamma/2]$ if $b_x$ is well chosen.}.

Since $f(\Fol^{cs}(x)) \en \Fol^{cs}(f(x))$ and
$f(W^{cs}_{loc}(x)) \en  b_{f(x)}( \frac{1}{2}\D^2 \times [-1,1])$
we get the desired property. \lqqd

The fact that $f\in \U_1$ is semiconjugated with $A$ together with
the fact that the semiconjugacy is $\eps$-$C^0$-close to the
identity gives us the following easy properties about the fibers
(preimages under $h_f$) of the points.

We denote

$$\Pi^{uu}_{x,z}: U\en W^{cs}_{loc}(x) \to W^{cs}_{loc}(z)$$

\noindent the unstable holonomy where $z\in \cF^{u}(x)$ and $U$ is
a neighborhood of $x$ in $W^{cs}_{loc}(x)$ which can be considered
large if $z$ is close to $x$ in $\cF^u(x)$. In particular, let
$\gamma>0$ be such that if $z \in \cF^{u}_{\gamma}(x)$ then the
holonomy is defined in a neighborhood of radius $\eps$ of $x$.

\begin{prop}\label{propiedades} Consider $y= h_f(x)$ for $x\in \TT^3$:
\begin{enumerate}
\item $h_f^{-1}(\{y\})$ is a compact connected set contained in
$W^{cs}_{loc}(x)$. \item If $z\in \cF^{u}_{\gamma}(x)$, then
$h_f(\Pi^{uu}_{x,z}(h^{-1}_f (\{y\})))$ is exactly one point.
\end{enumerate}
\end{prop}
\dem (1)  Since $h_f$ is $\eps$-$C^0$-close the identity, we get
that for every point $y\in \T^3$, $h_f^{-1}(\{y\})$ has diameter
smaller than $\eps$. Since $\eps$ is small compared to $\delta$,
it is enough to prove that $h_f^{-1}(\{y\}) \en W^{cs}_{loc}(x)$
for some $x\in h_f^{-1}(\{y\})$.

Assume that for some $y\in \T^3$, $h_f^{-1}(\{y\})$ intersects two
different center stable leaves of $\Fol^{cs}$ in points $x_1$ and
$x_2$.

Since the points are near, we have that $\cF^{u}_{\gamma}(x_1)
\cap W^{cs}_{loc}(x_2) =\{z\}$. Thus, by forward iteration, we get
that for some $n_0>0$ we have $d(f^{n_0}(x_1), f^{n_0}(z))>
3\delta$.

Lemma \ref{trapping} gives us that $d(f^{n_0}(x_2), f^{n_0}(z)) <
2\delta$. We get that $d(f^{n_0}(x_1),f^{n_0}(x_2)) > \delta$
which is a contradiction since $\{f^{n_0}(x_1),f^{n_0}(x_2)\} \en
h_f^{-1}(\{A^{n_0}(y)\})$ which has diameter smaller than $\eps
\ll \delta$.

Also, since the dynamics is trapped in center stable manifolds, we
get that the fibers must be connected since one can write them as

$$ h^{-1}(\{h(x)\})= \bigcap_{n\geq 0} f^n(W^{cs}_{loc}(f^{-n}(x))) .$$

(2)  Since $f^{-n}(h_f^{-1}(\{y\})) = h_f^{-1}(\{A^{-n}(y)\})$ we
get that $\diametro(f^{-n}(h_f^{-1}(\{y\}))) < \eps$ for every
$n>0$.

This implies that there exists $n_0$ such that if $n>n_0$ then
$f^{-n}(\Pi^{uu}_{x,z}(h^{-1}_f (\{y\})))$ is sufficiently near
$f^{-n}(h_f^{-1}(\{y\}))$. So, we have that
$$\diametro(f^{-n}(\Pi_{x,z}^{uu}(h_f^{-1}(\{y\})))) < 2\eps \ll
\delta .$$

Assume that $h_f(\Pi_{x,z}^{uu}(h_f^{-1}(\{y\})))$ contains more
than one point. These points must differ in the stable coordinate
of $A$, so, after backwards iteration we get that they are at
distance bigger than $3\delta$. Since $h_f$ is $\eps$-$C^0$-close
the identity this represents a contradiction.
 \lqqd

\begin{obs} The second statement of the previous proposition gives that the fibers of $h_f$ are invariant under unstable holonomy.\finobs
\end{obs}

The following simple lemma is essential in order to satisfy the
properties of Proposition \ref{ProposicionMecanismo}.

\begin{lema}\label{conexosgrandes}
For every $f\in \U_1$, given a disc $D$ in $W^{cs}_{loc}(x)$ whose
image by $h_f$ has at least two points, then $D \cap \cF^{u}(r_f)
\neq \emptyset$ and the intersection is transversal.
\end{lema}

\dem Given a subset $K\en \Fol^{cs}(x)$ we define its \emph{center
stable diameter} as the diameter with the metric $d_{cs}$ defined
above induced by the metric in the manifold. We shall first prove
that there exists $n_0$ such that $\diametro_{cs} (f^{-n_0}(D)) >
100 \delta$:

Since $D$ is arc connected so is $h_f(D)$, so, it is enough to
suppose that $\diametro (D)<\delta$. We shall first prove that
$h_f(D)$ is contained in a stable leaf of the stable foliation of
$A$. Otherwise, there would exist points in $h_f(D)$ whose future
iterates separate more than $2\delta$, this contradicts that the
center stable plaques are trapped for $f$ (Lemma \ref{trapping}).

One now has that, since $A$ is Anosov and that $h_f(D)$ is a
connected compact set with more than two points contained in a
stable leaf of the stable foliation, there exists $n_0 > 0$ such
that $A^{-n_0}(h_f(D))$ has stable diameter bigger than
$200\delta$ (recall that $\diametro \TT^3 > 400\delta$). Now,
since $h_f$ is close to the identity, one gets the desired
property.

We conclude by proving the following:

\begin{af} If there exists $n_0$ such that $f^{-n_0}(D)$ has diameter
larger than $100\delta$, then $D$ intersects $\cF^{u}(r_f)$.
\end{af}
\medskip
\dem This is proved in detail in section 6.1 of \cite{BV} so we
shall only sketch it.

If $f^{-n_0}(D)$ has diameter larger than $100\delta$, from how we
choose $\delta$ we have that there is a compact connected subset
of $f^{-n_0}(D)$ of diameter larger than $35\delta$ which is
outside $B_{6\delta}(q)$.

So, $f^{-n_0-1}(D)$ will have diameter larger than $100\delta$ and
the same will happen again. This allows to find a point $x\in D$
such that $\forall n>n_0$ we have that $f^{-n}(x)\notin
B_{6\delta}(q)$.

Now, considering a small disc around $x$ we have that by backward
iterates it will contain discs of radius each time bigger and this
will continue while the disc does not intersect $B_\delta(q)$. If
that happens, since $f^{-n}(x) \notin B_{6\delta}(q)$ the disc
must have radius at least $3\delta$.

This proves that there exists $m$ such that $f^{-m}(D)$ contains a
center stable disc of radius bigger than $2\delta$, so, the
unstable manifold of $r_f$ intersects it. Since the unstable
manifold of $r_f$ is invariant, we deduce that it intersects $D$
and this concludes the proof of the claim.

Transversality of the intersection is immediate from the fact that
$D$ is contained in $\cF^{cs}$ which is transversal to $\cF^u$.
\finobs \lqqd

We obtain the following corollary which puts us in the hypothesis
of Proposition \ref{ProposicionMecanismo}:

\begin{cor}\label{interior}
For every $f\in \U_1$, let $x\in \partial h_f^{-1}(\{y\})$
$($relative to the local center stable manifold of
$h_f^{-1}(\{y\}))$, then, $x$ belongs to the homoclinic class of
$r_f$, and in particular, to $\cQ_f$.
\end{cor}

\dem Notice first that the stable manifold of $r_f$ coincides with
$\cF^{cs}(r_f)$ which is dense in $\TT^3$. This follows from the
fact that when iterating an unstable curve, it will eventually
intersect the stable manifold of $r_f$, since the stable manifold
of $r_f$ is invariant, we obtain the density of $\cF^{cs}(r_f)$.

Now, considering $x\in \partial h_f^{-1}(\{y\})$, and $\eps>0$, we
consider a connected component $\tilde D$ of $\cF^{cs}(r_f) \cap
B_\eps(x)$. Clearly, since the fibers are invariant under holonomy
and $x \in \partial h_f^{-1}(\{y\})$ we get that $\tilde D$
contains a disk $D$ which is sent by $h_f$ to a non trivial
connected set. Using the previous lemma we obtain that there is a
homoclinic point of $r_f$ inside $B_\eps(x)$ which concludes.
\lqqd

The following corollary will allow us to use Theorems
\ref{TeoremaBonattiDiazPujals} and \ref{Teorema-PalisViana}.

\begin{cor}\label{qgenLambdag} For every $f\in \U_1$ we have that $q_f\in H(r_f)$.
\end{cor}

\dem Consider $U$, a neighborhood of $q_f$, and $D$ a center
stable disc contained in $U$.

Since the stable manifold of $q_f$ has length bigger than
$\delta>\eps$, after backward iteration of $D$ one gets that
$f^{-k}(D)$ will eventually have diameter larger than $\eps$, thus
$h_f(D)$ will have at least two points, this means that $q_f \in
\partial h_f^{-1}(\{h(q_f)\})$. Corollary \ref{interior}
concludes. \lqqd

We finish this section by proving the following theorem which is
the topological part of Theorem \ref{Teorema-EjemploCarvalho}.

\begin{teo}\label{TeoParteTopologica}
\begin{itemize}
    \item[(i)] For every $f\in \cU_1$ there exists a unique quasi-attractor $\cQ_f$ which contains the homoclinic class $H(r_g)$ and such that every chain-recurrence class $R\neq \cQ_f$ is contained in a periodic disc of $\cF^{cs}$.
    \item[(ii)] For every $f \in \cG_{BC} \cap \cG_{BDV} \cap \cU_1$ we have that $H(r_f)=\cQ_f$ and is contained in the closure of the sources of $f$.
    \item[(iii)] For every $r\geq 2$, there exists a $C^1$-open dense subset $\cU_3$ of $\cU_1$ and a residual subset $\cG^r \en \cU_3 \cap \Diff^r(\TT^3)$ such that for every $f\in \cG^r$ the homoclinic class $H(r_f)$ intersects the closure of the sources of $f$.
    \item[(iv)] For every $f\in \cU_1$ there exists a unique Milnor attractor contained in $\cQ_f$.
\end{itemize}
\end{teo}

\dem{\!\!} Part (i) follows from Proposition
\ref{ProposicionMecanismo} since $h_f$ is the desired
semiconjugacy: Indeed, Proposition \ref{propiedades} and Corollary
\ref{interior} show that the hypothesis of the mentioned
proposition are verified (notice that $A$ is clearly expansive).

Part (ii) follows from Theorem \ref{TeoremaBonattiDiazPujals}
using Corollary \ref{qgenLambdag}. Notice that $E^{cs}$ cannot be
decomposed in two $Df$-invariant subbundles since $Df$ has complex
eigenvalues in $r_f$.

Similarly, part (iii) follows from Theorem
\ref{Teorema-PalisViana}. The need for considering $\cU_3$ comes
from \cite{BDAbundanceTangencies}.

To prove (iv) notice that every point which does not belong to the
fiber of a periodic orbit belongs to the basin of $\cQ_f$: Since
there are only countably many periodic orbits and their fibers are
contained in two dimensional discs (which have zero Lebesgue
measure) this implies directly that the basin of $\cQ_f$ has total
Lebesgue measure:

Consider a point $x$ whose omega-limit set $\omega(x)$ is
contained in a chain recurrence class $R$ different from $\cQ_f$.
Then, since this chain recurrence class is contained in the fiber
$h_f^{-1}(\cO)$ of a periodic orbit $\cO$ of $A$, which in turn is
contained in the local center stable manifold of some point $z\in
\TT^3$. This implies that some forward iterate of $x$ is contained
in $W^{cs}_{loc}(z)$. The fact that the dynamics in $W^{cs}_{loc}$
is trapping (see Lemma \ref{trapping}) and the fact that $\partial
h_f^{-1}(\cO) \en \cQ_f$ (see Corollary \ref{interior}) gives that
$x$ itself is contained in $h_f^{-1}(\cO)$ as claimed.

Now, Lemma \ref{Lemma-MilnorAttractor} implies that $\cQ_f$
contains an attractor in the sense of Milnor. \lqqd

We have just proved parts (a), (b), (c) and (d) of Theorem
\ref{Teorema-EjemploCarvalho} hold in $\cU_3$. In fact, for the
$C^1$-topology, we have obtain a slightly stronger property than
(d) holds in $\cU_1$. Also, we have proved that (e) is satisfied.

\begin{obs}
The choice of having complex eigenvalues for $A$ was only used to
guaranty that $E^{cs}$ admits no $Df$-invariant subbundles. One
could have started with any linear Anosov map $A$ and modify the
derivative of a given fixed or periodic point $r$ to have complex
eigenvalues and the construction would be the same. \finobs
\end{obs}

\subsubsection{Ergodic properties} In this section we shall work
with $f\in \cU_2$ so that properties (P1)-(P7) are verified.

Consider the open set $U$ defined above such that
$f(\overline{U})\en U$ and consider:

$$ \Lambda_f = \bigcap_n f^n(U) $$

We shall show that the hypothesis of Theorem
\ref{Teorema-BonattiVianaA} are satisfied for $\Lambda_f$, and
thus, we get that there are at most finitely many SRB measures
such that the union of their (statistical) basins has full
Lebesgue measure in the topological basin of $\Lambda_f$. We must
show that in every unstable arc there is a positive Lebesgue
measure set of points such that $\lambda^{cs}(x)<0$.

\begin{prop}\label{Proposicion-EsMostlyContracting}
For every $x\in \TT^3$ and $D\en W^{uu}_{loc}(x)$ an unstable arc,
we have full measure set of points which have negative Lyapunov
exponents in the direction $E^{cs}$.
\end{prop}

\dem The proof is exactly the same as the one in Proposition 6.5
of \cite{BV} so we omit it. Notice that conditions (P2), (P4) and
(P7) in our construction imply conditions $(i)$ and $(ii)$ in
section 6.3 of \cite{BV}. \lqqd

To prove uniqueness of the SRB measure, we must show that there is
a unique minimal set of the unstable foliation inside $\Lambda_f$
to apply Theorem \ref{Teorema-BonattiVianaB}.

However, the fact that the stable manifold of $r_f$ contains
$W^{cs}_{loc}(r_f)$, gives that every unstable manifold intersects
$W^s(r_f)$ and so we get that every compact subset of $\Lambda_f$
saturated by unstable sets must contain $\overline{\cF^u(r_f)}$.
This implies that for every $x\in \overline{\cF^u(r_f)}$ we have
that $\overline{\cF^u(r_f)} = \overline{\cF^u(x)}$ and
$\overline{\cF^u(r_f)}$ is the only compact set with this property
(we say that $\overline{\cF^u(r_f)}$ is the unique minimal set of
the foliation $\cF^u$).

 We get thus that $f$ admits an unique SRB
measure $\mu$ and clearly, the support of this SRB measure is
$\overline{\cF^{u}(r_f)}$.

We claim that $\overline{\cF^{u}(r_f)}=H(r_f)$: this follows from
the fact that the SRB measure $\mu$ is hyperbolic (by Proposition
\ref{Proposicion-EsMostlyContracting}) and that the partially
hyperbolic splitting separates the positive and negative exponents
of $\mu$ and so verifies the hypothesis of Theorem
\ref{Teorema-PesinC1}.

Finally, since the SRB measure has total support and almost every
point converges to the whole support, we get that the attractor is
in fact a minimal attractor in the sense of Milnor. We have
proved:

\begin{prop}\label{Proposicion-UnicaSRBEjemploDA}
If $f\in \cU_2$ is of class $C^2$, then $f$ admits a unique SRB
measure whose support coincides with
$\overline{\cF^{u}(r_f)}=H(r_f)$. In particular,
$\overline{\cF^{u}(r_f)}$ is a minimal attractor in the sense of
Milnor for $f$.
\end{prop}

The importance of considering $f$ of class $C^2$ comes from the
fact that with lower regularity, even if we knew that almost every
point in the unstable manifold of $r_f$ has stable manifolds, we
cannot assure that these cover a positive measure set due to the
lack of absolute continuity in the center stable foliation.

However, the information we gathered for smooth systems in $\cU_2$
allows us to extend the result for $C^1$-generic diffeomorphisms
in $\cU_2$. Recall that for a $C^1$-generic diffeomorphisms $f\in
\cU_2$, the homoclinic class of $r_f$ coincides with $\cQ_f$.

\begin{teo}\label{AtractorMinimalGenerico}
There exists a $C^1$-residual subset $\cG_M \en \cU_2$ such that
for every $f\in \cG_M$ the set $\cQ_f =H(r_f)$ is a minimal Milnor
attractor.
\end{teo}

\dem Notice that since $r_f$ has a well defined continuation in
$\cU_2$, it makes sense to consider the map $f \mapsto
\overline{\cF^{u}(r_f)}$ which is naturally semicontinuous with
respect to the Haussdorff topology. Thus, it is continuous in a
residual subset $\cG_1$ of $\cU_2$. Notice that since the
semicontinuity is also valid in the $C^2$-topology, we have that
$\cG_1 \cap \Diff^2(\TT^3)$ is also residual in $\cU_2 \cap
\Diff^2(\TT^3)$.

It suffices to show that the set of diffeomorphisms in $\cG_1$ for
which $\overline{\cF^{u}(r_f)}$ is a minimal Milnor attractor is a
$G_\delta$ set (countable intersection of open sets) since we have
already shown that $C^2$ diffeomorphisms (which are dense in
$\cG_1$) verify this property.

Given an open set $U$, we define

$$U^+(f)= \bigcap_{n\leq 0} f^n(\overline{U}).$$

Let us define the set $\cO_U(\eps)$ as the set of $f\in \cG_1$
such that they satisfy one of the following (disjoint) conditions

\begin{itemize}
\item[-] $\overline{\cF^{u}(r_f)}$ is contained in $U$ \emph{or}
\item[-] $\overline{\cF^{u}(r_f)} \cap \overline{U}^c \neq
\emptyset$ and $Leb(U^+(f))<\eps$
\end{itemize}

We must show that this sets are open in $\cG_1$ (it is not hard to
show that if we consider an countable basis of the topology and
$\{U_n\}$ are finite unions of open sets in the basis then
$\cG_M=\bigcap_{n,m} \cO_{U_n}(1/m)$).

To prove that these sets are open, we only have to prove the
semicontinuity of the measure of $U^+(f)$ (since the other
conditions are clearly open from how we chose $\cG_1$).

Let us consider the set $\tilde K = \overline U \backslash
U^+(f)$, so, we can write $\tilde K$ as an increasing union
$\tilde K = \bigcup_{n\geq 1} K_n$ where $K_n$ is the set of
points which leave $\overline{U}$ in less than $n$ iterates.

So, if $Leb(U^+(f))<\eps$, we can choose $n_0$ such that
$Leb(\overline{U} \backslash K_{n_0}) < \eps$, and in fact we can
consider $K_{n_0}'$ a compact subset of $K_{n_0}$ such that
$Leb(\overline{U}\backslash K'_{n_0})<\eps$.

In a small neighborhood $\cN$ of $f$, we have that if $f'\in \cN$,
then $K'_{n_0} \en \overline{U} \backslash U^+(f')$. This
concludes. \lqqd

This completes the proof of part (f) of Theorem
\ref{Teorema-EjemploCarvalho}.

\lqqd

\subsection{Example of Plykin type}\label{SubSection-EjemploPlykin}

The examples in subsection \ref{SubSection-EjemploDA} cannot be
embedded in any manifold as the ones of Bonatti-Li-Yang. We were
able to adapt the construction in order to get an example with
similar properties which can be embedded in any isotopy class of
diffeomorphisms of a manifold. However, we were not able to obtain
the same strong ergodic properties (see \cite{PotWildMilnor} for
more discussion and problems).

\begin{teo}\label{Teorema-EjemploPlykin}
For every $d$-dimensional manifold $M$ and every isotopy class of
diffeomorphisms of $M$ there exists a $C^1$-open set $\cU$ of
$\Diff^r(M)$ such that for some open neighborhood $U$ in $M$:
\begin{itemize}
\item[(a)] Every $f\in \cU$ has a unique quasi-attractor $\cQ_f$
in $U$ which contains a homoclinic class and has a partially
hyperbolic splitting $T_{\cQ_f}M = E^{cs} \oplus E^u$ which is
coherent. \item[(b)] Every chain recurrence class $R\neq \cQ_f$ is
contained in the orbit of a periodic leaf of the lamination
$\cF^{cs}$ tangent to $E^{cs}$ at $\cQ_f$. \item[(c)] There exists
a residual subset $\cG^r$ of $\cU$ such that for every $f\in
\cG^r$ the diffeomorphism $f$ has no attractors. In particular,
$f$ has infinitely many chain-recurrence classes. \item[(d)] For
every $f\in \cU$ there is a unique Milnor attractor $\tilde \cQ
\en \cQ_f$.
\end{itemize}
\end{teo}

The examples here are modifications of the product of a Plykin
attractor and the identity on the circle $S^1$. One can also
obtain them in order to provide examples of robustly transitive
attractors in dimension $3$ with splitting $E^{cs} \oplus E^u$.
The author is not aware of other known examples of such attractors
other than Carvalho's example which is only possible to be made in
certain isotopy classes of diffeomorphisms\footnote{N. Gourmelon
communicated me the possibility of constructing examples of this
kind by bifurcating other robustly transitive diffeomorphisms such
as perturbations of time-one maps of Anosov flows.}.

In this section we shall show how to construct an example
verifying Theorem B. We shall see that we can construct a
quasi-attractor with a partially hyperbolic splitting $E^{cs}
\oplus E^u$ such that $E^{cs}$ admits no sub-dominated splitting.
In case $E^{cs}$ is volume contracting, it will turn out that this
quasi-attractor is in fact a robustly transitive attractor (thus
providing examples of robustly transitive attractors with
splitting $E^{cs} \oplus E^u$ in every $3-$dimensional manifold)
and when there is a periodic saddle of stable index $1$ and such
that the product of any two eigenvalues is greater than one and
using Theorems \ref{TeoremaBonattiDiazPujals} and
\ref{Teorema-PalisViana} we shall obtain that the quasi-attractor
will not be isolated for generic diffeomorphisms in a
neighborhood.

We shall work only in dimension $3$. It will be clear that by
multiplying the examples here with a strong contraction, one can
obtain examples in any manifold of any dimension.

A main difference between this construction and the one done in
section \ref{SubSection-EjemploDA} is the use of \emph{blenders}
instead of the argument \`a la Bonatti-Viana. Blenders were
introduced in \cite{BDAnnals} (see subsection
\ref{SubSection-Blenders}) and constitute a very powerful tool in
order to get robust intersections between stable and unstable
manifolds of compact sets. We shall only use the facts presented
in subsection \ref{SubSection-Blenders} and not enter in their
definition or construction for which there are many excellent
references (we recommend chapter 6 of \cite{BDV} in particular).

\subsubsection{Construction of the example}

Let us consider $P: \DD^2 \hookrightarrow \DD^2$ the map given by
the Plykin attractor in the disk $\DD^2$ (see \cite{Robinson}).

We have that $P(\DD^2) \en int(\DD^2)$, there exist a hyperbolic
attractor $\Upsilon \en \DD^2$ and three fixed sources (we can
assume this by considering an iterate).

There is a neighborhood $N$ of $\Upsilon$ which is homeomorphic to
the disc with $3$ holes removed that satisfies that $P(\overline
N) \en N$ and

$$\Upsilon= \bigcap_{n\geq 0} P^n(N) . $$

It is well known that given $\eps>0$, one can choose a finite
number of periodic points $s_1, \ldots, s_N$ and $L>0$ such that
if $A= \bigcup_{i=1}^N W^u_L(s_i)$, then, for every $x\in \Upsilon
\setminus A$ one has that $A$ intersects both connected components
of $W^s_\eps (x) \setminus \{x\}$.


We now consider the map $F_0: \DD^2 \times S^1 \hookrightarrow
\DD^2 \times S^1$ given by $F_0(x,t)= (P(x), t)$ whose chain
recurrence classes consist of the set $\Upsilon \times S^1$ which
is a (non transitive) partially hyperbolic attractor and three
repelling circles.

In \cite{BDAnnals} they make a small $C^\infty$ perturbation $F_1$
of $F_0$, for whom the maximal invariant set in $U= N\times S^1$
becomes a $C^1$-robustly transitive partially hyperbolic attractor
$Q$ which remains homeomorphic to $\Upsilon \times S^1$.

This attractor has a partially hyperbolic structure of the type
$E^s \oplus E^c \oplus E^u$. One can make this example in order
that it fixes the boundary of $\D^2 \times S^1$, this allows to
embed this example (and all the modifications we shall make) in
any isotopy class of diffeomorphisms of any $3$-dimensional
manifold (since every diffeomorphism is isotopic to one which
fixes a ball, then one can introduce this map by a simple
surgery).

In \cite{BDAnnals} the diffeomorphism $F_1$ constructed verifies
the following properties (see \cite{BDAnnals} section 4.a page
391, also one can find the indications in \cite{BDV} section
7.1.3):

\begin{itemize}
\item[(F1)] $F_1$ leaves invariant a $C^1$-lamination $\cF^{cs}$
(see \cite{HPS} chapter 7 for a precise definition) tangent to
$E^s \oplus E^c$ whose leaves are homeomorphic to $\RR\times
\SS^1$. \item[(F2)] There are periodic points $p_1, \ldots, p_N$
of stable index $1$ such that for every $x\in Q$ one has that the
connected component of $\cF^{cs}(x) \setminus
\overline{(W^u_L(p_1) \cup \ldots \cup W^u_L(p_N))}$ containing
$x$ has finite volume for every $x\in Q \setminus \bigcup_{i=1}^N
W^u_L(p_i)$. Here $W^u_L(p_i)$ denotes the $L$-neighborhood of
$p_i$ in its unstable manifold with the metric induced by the
ambient. \item[(F3)] There is a periodic point $q$ periodic point
of stable index $2$ contained in a $cu$-blender $K$ such that for
every $1 \leq i \leq N$ the stable manifold of $p_i$ intersects
the activating region of $K$. By Proposition \ref{PropBlender},
the unstable manifold of $q$ is dense in the union of the unstable
manifolds of $p_i$.

\item[(F4)] The local stable manifold of $q$ intersects every
unstable curve of length larger than $L$.
\end{itemize}

Before we continue, we shall make some remarks on the properties.
The hypothesis (F1) on the differentiability of the lamination
$\cF^{cs}$ will be used in order to apply the results on normal
hyperbolicity of \cite{HPS} (chapter 7, Theorem 7.4, see also
Section \ref{Section-HIRSHPUGHSHUB}).

It can be seen in \cite{BDAnnals} that the construction of $F_1$
is made by changing the dynamics in finitely many periodic circles
and this can be done without altering the lamination $\cF^{cs}$
which is $C^1$ before modification.

This is in fact not necessary; it is possible to apply the
barehanded arguments of the proof of Theorem 3.1 of \cite{BuFi} in
order to obtain that for the modifications we shall make, there
will exist a lamination tangent to the bundle $E^{cs}$.

Hypothesis (F2) is justified by the fact that the Plykin attractor
verifies the same property and the construction of $F_1$ in
\cite{BDAnnals} is made by changing the dynamics in the periodic
points by Morse-Smale diffeomorphisms which give rise property
(F2) (see section 4.a. of \cite{BDAnnals}). Notice that by
continuous variation of stable and unstable sets, this condition
is $C^1$-robust.

Property (F3) is the essence in the construction of
\cite{BDAnnals}, $cs$-blenders are the main tool for proving the
robust transitivity of this examples. As explained in subsection
\ref{SubSection-Blenders} this is a $C^1$-open property.

Property (F4) is given by the fact that the local stable manifold
of $q$ can be assumed to be $W^s_{loc}(s) \times \SS^1$ with a
curve removed, where $s\in \Upsilon$ is a periodic point. This is
also a $C^1$-open property.

\begin{figure}[ht]
\begin{center}
\scalebox{.8}{\input{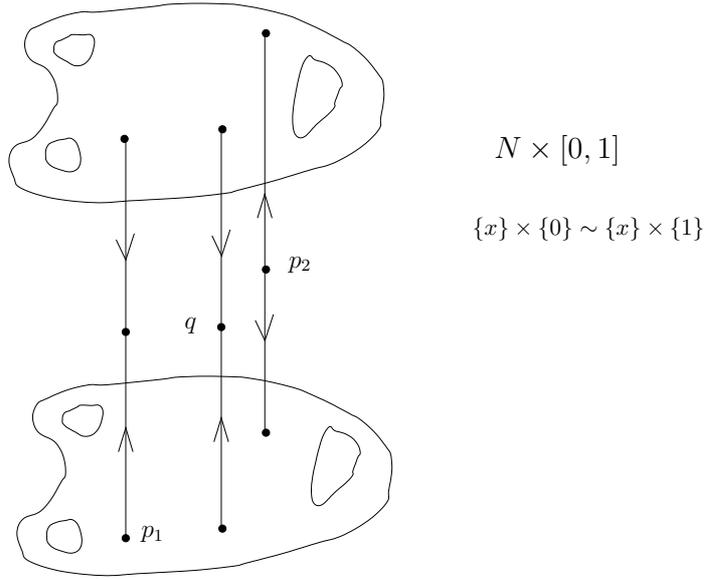}} \caption{\small{How to
construct $F_1$ by small $C^\infty$ perturbations in finitely many
circles.}} \label{FiguraAxRodeado}
\end{center}
\end{figure}

Let us consider a periodic point $r_1 \in Q$ of stable index $1$
and another one $r_2$ of stable index $2$. We can assume they are
fixed (modulo considering an iterate of $F_1$). Consider
$\delta>0$ small enough such that $B_{6\delta}(r_1) \cup
B_{6\delta}(r_2)$ is disjoint from:

\begin{itemize}
\item[-] the periodic points $p_1, \ldots, p_N, q$ defined above,
\item[-] the blender $K$, \item[-] $\overline{(W^u_L(p_1) \cup
\ldots \cup W^u_L(p_N))}$ and \item[-] from $\cF^{u}_{L'}(q)$
(where $L'$ is chosen such that $\cF^{u}_{L'}(q)$ intersects $K$).
\end{itemize}

In the same vein as in subsection \ref{SubSection-EjemploDA} we
shall first construct a diffeomorphism $F_2$ modifying $F_1$ such
that:

\begin{itemize}
\item[-] $F_2$ coincides with $F_1$ outside $B_\delta(r_2)$.
\item[-] $F_2$ preserves the center-stable foliation of $F_1$.
\item[-] $DF_2$ preserves narrow cones $\cE^u$ and $\cE^{cs}$
around the unstable direction $E^u$ and the center stable
direction $E^{s} \oplus E^c$ of $F_1$ respectively. Also, vectors
in $\cE^{u}$ are expanded uniformly by $DF_2$ while every plane
contained in $\cE^{cs}$ verifies that the volume\footnote{This
means with respect to the Riemannian metric which allows to define
a notion of 2-dimensional volume in each plane.} is contracted by
$DF_2$. \item[-] The point $r_2$ remains fixed for $F_2$ but now
has complex eigenvalues in $r_2$.
\end{itemize}

Before we continue with the construction of the example to prove
Theorem \ref{Teorema-EjemploPlykin}, we shall make a small detour
to sketch the following:

\begin{prop} There exists an open $C^1$-neighborhood $\cV$ of $F_2$ such that for every $f\in \cV$ one has that $f$ has a transitive attractor in $U$.
\end{prop}

\esbozo{} Notice that one can choose $\cV$ such that for every
$f\in \cV$ one preserves a center-stable foliation close to the
original one. Also, one can assume that properties (F2) and (F3)
still hold for the continuations $p_i(f)$ and $q(f)$ since $F_2$
coincides with $F_1$ outside $B_\delta(r_2)$ and these are
$C^1$-robust properties.

Also, we demand that for every $f\in \cV$, the derivative of $f$
preserves the cones $\cE^u$ and $\cE^{cs}$, contracts volume in
$E^{cs}\en \cE^{cs}$ (the plane tangent to the center-stable
foliation) and expands vectors in $E^u \en \cE^u$.

Consider now a center stable disk $D$ and an unstable curve
$\gamma$ which intersect the maximal invariant set

$$\cQ_f= \bigcap_{n>0} f^n(U) . $$

Since by future iterations $\gamma$ will intersect the stable
manifold of $q(f)$ (property (F3)) we obtain that by the
$\lambda-$lemma it will accumulate the unstable manifold of
$q(f)$. Since the unstable manifold of $q(f)$ is dense in the
union of the unstable manifolds $W^u(p_1(f))\cup \ldots \cup
W^u(p_N(f))$ we obtain that the union of the future iterates of
$\gamma$ will also be dense there.

Now, iterating backwards the disk $D$ we obtain, using that
$Df^{-1}$ expands volume in the center-stable direction that the
diameter of the disk grows exponentially with these iterates.

Condition (F2) will now imply that eventually the backward
iterates of $D$ will intersect the future iterates of $\gamma$.
This implies transitivity. \lqqd

Now, we shall modify $F_2$ inside $B_\delta(r_1)$ in order to
obtain an open set to satisfy Theorem \ref{Teorema-EjemploPlykin}.
So we shall obtain $F_3$ such that:

\begin{itemize}
\item[-] $F_3$ coincides with $F_2$ outside $B_\delta(r_1)$.
\item[-] $F_3$ preserves the center-stable lamination of $F_2$.
\item[-] $DF_3$ preserves narrow cones $\cE^u$ and $\cE^{cs}$
around the unstable direction $E^u$ and the center stable
direction $E^{cs}$ of $F_2$. Also, vectors in $\cE^{u}$ are
expanded uniformly by $DF_3$. \item[-] $r_1$ is a saddle with
stable index $1$, the product of any pair of eigenvalues is larger
than $1$ and the stable manifold of $r_1$ intersects the
complement of $B_{6\delta}(r_1)$.
\end{itemize}

We obtain a $C^1$ neighborhood $\cU_1$ of $F_3$ where for $f\in
\cU$, if we denote

$$\Lambda_f = \bigcap_{n\geq0} f^n(U): $$

\begin{itemize}
\item[(P1')] There exists a continuation of the points $p_1,
\ldots, p_N, q, r_1, r_2$ which we shall denote as $p_i(f)$,
$q(f)$ and $r_i(f)$. The point $r_1(f)$ is a saddle of stable
index $1$ and its stable manifold intersects the complement of
$B(r_1,6\delta)$. \item[(P2')] There is a $Df$-invariant families
of cones $\cE^u$ in $\cQ_f$ and for every $v\in \cE^u(x)$ we have
that
$$\|D_xf v\|\geq \lambda \|v\|.$$
\item[(P3')] $f$ preserves a lamination $\cF^{cs}$ which is $C^0$
close to the one preserved by $F_3$ and which is trapped in the
sense that there exists a family $W^{cs}_{loc}(x) \en \cF^{cs}(x)$
such that for every point $x\in \cQ_f$ the plaque
$W^{cs}_{loc}(x)$ is homeomorphic to $(0,1)\times \SS^1$ and
verifies that

    $$f(\overline{W^{cs}_{loc}(x)}) \en W^{cs}_{loc}.$$
    Moreover, the stable manifold of $r_1(f)$ intersects the complement of $W^{cs}_{loc}(r_1(f))$.
\item[(P4')] Properties (F2),(F3) and (F4) are satisfied for $f$
and every curve $\gamma$ tangent to $\cE^u$ of length larger than
$L$ intersects the stable manifold of $q(f)$.
\end{itemize}

Notice that (P4') implies that there exists a unique
quasi-attractor $\cQ_f$ in $U$ for every $f\in \cU$ which contains
the homoclinic class $H(q(f))$ of $q(f)$ (the proof is the same as
Lemma \ref{unicocasiatractor}).

\subsubsection{The example verifies the mechanism of Proposition
\ref{ProposicionMecanismo}} We shall show that every $f\in \cU$ is
in the hypothesis of Proposition \ref{ProposicionMecanismo} which
will conclude the proof of Theorem \ref{Teorema-EjemploPlykin} as
in Theorem \ref{TeoParteTopologica}. We shall only sketch the
proof since it has the same ingredients as the proof of Theorem A,
the main difference is that instead of having an a priori
semiconjugacy we must construct one.

To construct the semiconjugacy, one uses property (P3'),
specifically the fact that $f(\overline{W^{cs}_{loc}(x)}) \en
W^{cs}_{loc}(x)$ (compare with Lemma \ref{trapping}) to consider
for each point $x\in \Lambda_f$ the set:

$$ A_x = \bigcap_{n\geq 0} f^n(\overline{W^{cs}_{loc}(f^{-n}(x))}) $$

\noindent  (compare with Proposition \ref{propiedades} (1)). One
easily checks that the sets $A_x$ constitute a partition of
$\Lambda_f$ into compact connected sets contained in local center
stable manifolds and that the partition is upper-semicontinuous.
It is not hard to prove that if $h_f: \Lambda_f \to
\Lambda_f/_{\sim}$ is the quotient map, then, the map $g:
\Lambda_f /_{\sim} \to \Lambda_f/_{\sim}$ defined such that

$$ h_f \circ f = g \circ h_f $$

\noindent is expansive (in fact, $\Lambda_f/_{\sim}$ can be seen
to be homeomorphic to $\Upsilon$ and $g$ conjugated to $P$). See
\cite{Daverman} for more details on this kind of decompositions
and quotients.

\begin{figure}[ht]
\begin{center}
\input{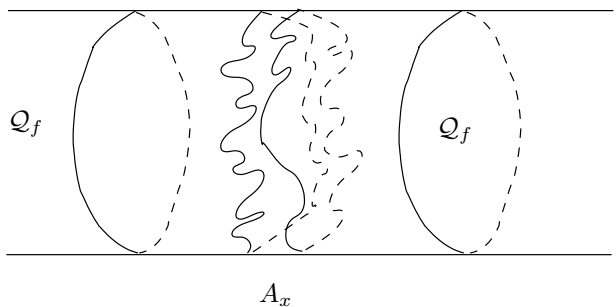}
\caption{\small{The set $A_x$ is surrounded by points in
$\overline{W^u(q)} \en \cQ_f$ }} \label{FiguraAxRodeado}
\end{center}
\end{figure}

Since fibers are contained in center stable sets, we get that
$h_f$ is injective on unstable manifolds and one can check that
the fibers are invariant under unstable holonomy (see the proof of
Proposition \ref{propiedades} (2)). Stable sets of $g$ are dense
in $\Lambda_f/_{\sim}$.

The point $r_1(f)$ will be in the boundary of
$h^{-1}_f(\{h_f(r(f))\})$ since its stable manifold is not
contained in $W^{cs}_{loc}(r_1(f))$.

We claim that the boundary of the fibers restricted to
center-stable manifolds is contained in the unique quasi-attractor
$\cQ_f$. This is proven as follows:

Assume that $x \in \partial h_f^{-1}(\{h_f(x)\})$ and consider a
small neighborhood $V$ of $x$. Consider a disk $D$ in
$W^{cs}_{loc}(x)$, since $x$ is a boundary point, we get that
$h_f(D)$ is a compact connected set containing at least two points
in the stable set of $h_f(x)$ for $g$, so by iterating backwards,
and using (F2) (guaranteed for $f$ by (P4')) we get that there is
a backward iterate of $D$ which intersects $\cF^u(q) \en \cQ_f$
which concludes.

Now, Theorem \ref{Teorema-EjemploPlykin} follows with the same
argument as for Theorem \ref{TeoParteTopologica}, using
Proposition \ref{ProposicionMecanismo} and the fact that $r_1(f)$
is contained in $\cQ_f$.

\lqqd

\subsection{Derived from Anosov revisited}\label{SubSection-EjemploDAalRevesRObustamenteTransitivo}

The first examples of non-hyperbolic $C^1$-robustly transitive
diffeomorphisms were given by Shub (see \cite{HPS} chapter 8) in
$\TT^4$ by considering a skew-product over $\TT^2$. Ma\~ne then
improved the example to obtain a non-hyperbolic $C^1$-robustly
transitive diffeomorphism of $\TT^3$ by deformation of an Anosov
diffeomorphism (\cite{ManheTopology}). These examples were
strongly partially hyperbolic with central dimension $1$. Bonatti
and Diaz (\cite{BDAnnals}) constructed examples with arbitrary
central dimension as well as examples isotopic to the identity but
still strongly partially hyperbolic\footnote{Bonatti also
constructed an (unpublished) example of robustly transitive
partially hyperbolic diffeomorphism of $\TT^3$ which was not
strongly partially hyperbolic using blenders.}.

Finally, Bonatti and Viana (\cite{BV}) constructed examples of
robustly transitive diffeomorphisms without any uniform bundle by
deforming an Anosov diffeomorphism of $\TT^4$ and improving the
ideas from Ma\~ne's example. In all the examples constructed by
deformation of an Anosov diffeomorphism there is an underlying
property which is that the dominated splitting (which must exist
by Theorem \ref{TeoremaBonattiDiazPujals}) has dimensions coherent
with the splitting the Anosov diffeomorphism has. This simplifies
a little the proofs and it is natural to construct the examples
this way.

We present here a very similar construction which does not
introduce new ideas. However, we believe that such an example is
not yet very well understood and may represent an important model
for starting to study partially hyperbolic systems with two
dimensional center and mixed behavior inside it.

\begin{teo}\label{TeoEJEMPLOCONDIMENSIONESCAMBIADAS}
There exists an open set $\cU \en \Diff^1(\TT^3)$ in the isotopy
class of a linear Anosov diffeomorphism $A$ with unstable
dimension 2 such that for every $f\in \cU$ we have that:
\begin{itemize}
\item[-] $f$ is partially hyperbolic with splitting $T\TT^3 =
E^{cs} \oplus E^u$ with $\dim E^u=1$ and such that $E^{cs}$ admits
no subdominated splitting.
\item[-] $f$ is transitive.
\end{itemize}
\end{teo}

\dem Let us consider a linear Anosov automorphism $A \in
SL(3,\ZZ)$ such that the eigenvalues of $A$ verify $0< \lambda_1 <
1/3 < 3 < \lambda_2 < \lambda_3$ and have defined eigenspaces
$E^s$, $E^u$ and $E^{uu}$ respectively.

We denote $p: \RR^3 \to \TT^3$ the covering projection. Notice
that since $\lambda_1 \lambda_2 \lambda_3=1$ we have that
$\lambda_1 \lambda_2 < 1$.

We choose $\delta>0$ to be defined later (in the same way as in
section \ref{SubSection-EjemploDA}) and we make a small
$C^0$-perturbation $f_0$ of $A$ supported on $B_\delta(p(0))$ such
that:

\begin{itemize}
\item[-] The tangent map $Df_0$ preserves narrow cones $\cC^u$ and
$\cC^{cs}$ around $E^{uu}$ and $E^s \oplus E^u$ respectively and
uniformly expands vectors in $\cC^u$. \item[-] The point $p(0)$
becomes a fixed point with stable index $2$ and complex stable
eigenvalues. \item[-] The jacobian of $f_0$ in any $2-$plane
inside $\cC^{cs}$ is smaller than $1$. \item[-] The $C^0$-distance
between $f_0$ and $A$ is smaller than $\eps$. Here, $\eps$ is
chosen in a way that every homeomorphism at $C^0$-distance smaller
than $2\eps$ of $A$ is semiconjugated to $A$ by a continuous map
at distance smaller than $\delta$ from the identity.
\end{itemize}

This modification can be made in the same way as we have done for
the construction of the example for Theorem
\ref{Teorema-EjemploCarvalho}.

We shall first consider a small $C^1-$open neighborhood $\cU$ of
$f_0$ such that for every $f\in \cU$ we have that:

\begin{itemize}
\item[-] The tangent map of $Df$ preserves the cones $\cC^{u}$ and
$\cC^{cs}$, uniformly expands vectors in $\cC^u$ and the jacobian
of $f$ in any $2-$plane inside $\cC^{cs}$ is smaller than $1$.
\item[-] The maximal invariant set of $f$ outside $B_\delta(p(0))$
is a hyperbolic set $\Lambda_f$ with invariant splitting $E^s
\oplus E^u \oplus E^{uu}$ which is close to the invariant
splitting of $A$ at those points and such that the norm of vectors
in an invariant cone around $E^u \oplus E^{uu}$ multiplies by $3$
outside $B_\delta(p(0))$. \item[-] The $C^0$-distance of $f$ and
$A$ is smaller than $2\eps$.
\end{itemize}

We obtain that for every $f\in \cU$ there exists $h_f : \TT^3 \to
\TT^3$ a semiconjugacy with $A$ which is at $C^0-$distance
$\delta$ of the identity, that is, we have that $h_f \circ f  = A
\circ h_f$ and we have that $d(h_f(x), x) < \delta$. Also, we can
assume that $\delta$ is such that the distance between the
connected components of $\tilde B = p^{-1}(B_{3\delta}(p(0)))$ is
larger than $100 \delta$.

We can choose $M$ such that for every disc $D$ of radius $2\delta$
outside $\tilde B$ and every curve of length larger than $M$ and
such that its image by $H_f$ (the lift of $h_f$) to $\RR^3$ is
close to an arc of stable manifold of $A$ we have that it has an
integer translate which intersects $D$.

Now, consider two open sets $U$ and $V$ in $\TT^3$, we must show
that there exists $n>0$ such that $f^n(U) \cap V \neq \emptyset$.
We shall work in the universal cover $\RR^3$ of $\TT^3$. The lift
of $f$ shall be denoted as $\tilde f$. Let us denote $U_0$ and
$V_0$ to connected components of $p^{-1}(U)$ and $p^{-1}(V)$
respectively.

Since there is an invariant cone $\cC^u$ where vectors are
expanded, we get that the diameter of $U_0$ grows exponentially
with future iterates of $\tilde f$.

It is not hard to show that there exists a point $x \in U_0$ and
$n_0>0$ such that for every $n\geq n_0$ we have that $\tilde f^n
(x) \notin \tilde B$. Indeed, once the diameter of $\tilde
f^{n_1}(U_0)$ is larger than $100\delta$, there exists a compact
connected subset $C_1$ of $U_0$ such that $\tilde f^{n_1}(C_1)$
does not intersect $\tilde B$ and $\diam (\tilde f^{-n_1}(C_1))>
40 \delta$. Now, we obtain inductively a decreasing intersection
of sets $C_k$ such that $\tilde f^{n_1+k} (C_k)$ does not
intersect $\tilde B$ and has large enough diameter. This implies
that in the intersection of all those sets one has the desired
point (see the proof the claim inside Lemma \ref{conexosgrandes}).

In a similar fashion, there exists a point $y \in V_0$ such that
its backward iterates after some $n_1<0$ are disjoint from $\tilde
B$ (here it is essential the fact that the jacobian of $f$
contracts uniformly the volume in the cone $\cC^{cs}$).

Now, we consider a small disk $D_1$ tangent to a small cone around
$E^u \oplus E^{uu}$ centered in $\tilde f^{n_0}(x)$ and contained
in $\tilde f^{n_0}(U_0)$. Iterating forward an using the fact that
vectors in that cone are expanded when are outside
$p^{-1}(B(p(0),\delta))$ and the point $\tilde f^{n_0}(x)$ remains
outside $\tilde B$ we get that eventually, $\tilde f^n (U_0)$
contains a disk $D_2 \en \tilde f^n(U_0)$  whose internal radius
is greater than $2\delta$ and whose center is contained in $\tilde
B^{c}$.

For past iterates we consider a disk $D_3$ whose tangent space
belongs to a narrow cone around $E^{cs}$ and we know that by
volume contraction its backward iterates grow exponentially in
diameter\footnote{This could require using dynamical coherence,
but we will ignore this issue. In any case, we can use it because
of the results in Chapter
\ref{Capitulo-ParcialmenteHiperbolicos}.}. We get that if $H_f$ is
the lift of $h_f$ we get that $H_f (\tilde f^{-n_1}(V_0))$
contains a curve transverse to the cone around $E^u \oplus E^{uu}$
and thus when iterating backwards by $A$ we obtain that its length
grows exponentially and it becomes close to the stable direction.
Since $H_f$ is at distance smaller than $\delta$ of the identity,
we get that eventually it has a translate which intersects $D_2$
by the remark made above which concludes. \lqqd


\section{Trapping quasi-attractors and further questions}\label{Section-TrappingAttractors}

To close this chapter, we will introduce a definition of a kind of
quasi-attractors which we believe to be in reach of understanding.
The definition is motivated by the examples of Bonatti-Li-Yang as
well as the examples presented in subsection
\ref{SubSection-EjemploDA}. The rest of the examples presented in
section \ref{Section-EjemplosQuasiAtractores} was essentially
introduced in order to show that understanding this kind of
quasi-attractors is not the end of the story, even in dimension
$3$.

\begin{defi}[Trapping quasi-attractors]\label{Definition-TrappingQA}
Let $\cQ$ be a quasi-attractor of a diffeomorphism $f: M \to M$
admitting a partially hyperbolic splitting of the form $T_\cQ M =
E^{cs} \oplus E^u$. We will say that $\cQ$ is a \emph{trapping
quasi-attractor} if it admits a locally $f$-invariant plaque
family $\{\cW^{cs}_x\}_{x\in \cQ}$ verifying that

$$ f(\overline{\cW^{cs}_x}) \en \cW^{cs}_x $$
\finobs
\end{defi}

This sets must be compared with \emph{chain-hyperbolic}
chain-recurrence classes defined in \cite{CP} which have a similar
yet different definition and have played an important role in the
proof of the $C^1$-Palis' conjecture on dynamics far from
homoclinic bifurcations.

Hyperbolic attractors are of course examples of trapping
quasi-attractors. Both the examples of Bonatti-Li-Yang and the
ones presented in section \ref{SubSection-EjemploDA} are
non-hyperbolic examples of this type.

As was mentioned, the results of Bonatti-Shinohara and the ones
presented in section \ref{SubSection-EjemploDA} present very
different properties and it seems a natural to try to understand
which are the reasons for this different behaviour.

The author's impression is that the different phenomena is related
with the topology of the intersection of the quasi-attractor with
center stable plaques: In one case (Bonatti-Li-Yang's example)
there is room to eject saddle points, and in the other one, the
quasi-attractor surrounds the point which one would like to eject
prohibiting its removal from the class.

What follows is speculation, and moreover, it is restricted to the
$3$-dimensional case. We are also assuming that the splitting is
of the form $T_\cQ M = E^{cs} \oplus E^u$ with $\dim E^u=1$ and
such that $E^{cs}$ does not admit a sub-dominated splitting.

It is in principle not obvious how to define the ``topology of the
intersection of the quasi-attractor with center-stable plaques''.
We propose that the following should be studied, and we hope
pursuing this line in the future:

It should be possible to make a quotient of the dynamics along
center-stable leaves in order to obtain an expansive quotient
which can be embedded as an expansive attractor of a
$3$-dimensional manifold (see Section
\ref{Section-DimensionesMayores} for an explicit construction in a
particular case). Then, the recent classification by A. Brown
\cite{Brown} should allow one to classify these attractors in
terms of the topology of the intersection with center-stable
plaques. We expect to obtain that if the quotient attractor is
one-dimensional then one will be able to perform the techniques of
Bonatti and Shinohara in order to get that these quasi-attractors
share the same properties as the example of Bonatti-Li-Yang.

In the case the dimension of the attractor has topological
dimension $3$, we expect that the properties will be similar to
the ones obtained for the example in subsection
\ref{SubSection-EjemploDA}.

It remains to understand also the case where the quotient has
topological dimension $2$. Although it is not hard to construct
examples of this behavior, it seems that a stronger understanding
of them needs to be acquired in order to have a more clear picture
on the possible dynamics such a quasi-attractor may have.

We finish this section by pointing out that trapping
quasi-attractors of course do not cover all the possibilities,
even in the case where the decomposition is of the form $T_\cQ M =
E^{cs} \oplus E^u$.

On the one hand, the examples of subsection
\ref{SubSection-EjemploPlykin} show that this may not happen,
however, those examples share some property with trapping
quasi-attractors since it is possible to find a family of locally
invariant plaques (whose topology is not of a disk but a cylinder)
which are ``trapped''. Indeed, this property was the key
ingredient in the proof of Theorem \ref{Teorema-EjemploPlykin}.

On the other hand, the example in subsection
\ref{SubSection-EjemploDAalRevesRObustamenteTransitivo} seems the
real challenge if one wishes to completely understand
quasi-attractors in dimension $3$ with splitting of the form
$T_\cQ M = E^{cs} \oplus E^u$. The problem is that the
center-stable direction contains ``unstable'' behavior, and this
is far less understood. C. Bonatti and Y. Shi (\cite{BSh}) have
provided new examples by the study of perturbations the time one
map of the Lorenz attractor (see \cite{BDV} chapter 9) which
essentially share this problems as well as having several other
new interesting properties. It is surely of great interest to have
a big list of examples before one attempts to attack the problem
of understanding general quasi-attractors in dimension $3$.

   \chapter{Foliations}\label{Capitulo-FoliacionesEHipParcial}


This chapter has two purposes. One the one hand, it presents
general results on foliations and gathers well known material
which serves as preliminaries for what we will prove.

On the other hand, we present ``almost new'' results on
foliations: In subsection \ref{SubSection-FoliacionesDeT3} we give
a classification of Reebless foliations on $\TT^3$ (this is
``almost new'' because the proofs resemble quite nearly those of
\cite{BBI2} and similar results exist for $C^2$-foliations).

In Section \ref{Section-GlobalProductStructure} we prove a result
which gives global product structure for certain codimension one
foliations on compact manifolds. Our results are slightly more
general than the ones which appear, for example, in \cite{Hector}
but with more restrictive hypothesis on the topology of the
manifold. Those restrictions on the topology of the manifold have
allowed us to prove this result with weaker hypothesis and a
essentially different proof.

This chapter contains an appendix which presents similar ideas in
the case of surfaces which can be read  independently of the rest
of the chapter and motivates some of the results of the next
chapter.

\section{Generalities on foliations}\label{Seccion-FoliacionesGeneralidades}

\subsection{Definitions}\label{SubSeccion-DefinicionesFoliaciones}

In section \ref{Section-FamiliasDePlacas} we reviewed the concept
of lamination, which consists of a partition of a compact subset
of a manifold by injectively immersed submanifolds which behave
nicely between them. The fact that laminations are only defined in
compact subsets suggests that the information they will give about
the topology of the manifold is not that strong (although there
are many exceptions). In this chapter, we shall work with
foliations, that for us will be laminations of  the whole manifold
and review many results which will give us a lot of information on
the relationship between the dynamics and the topology of the
phase space.

We will give a partial overview of foliations influenced by the
results we use in this thesis. The main sources will be
\cite{Calegari,CamachoLinsNeto,CandelConlon,Hector}.

\begin{defi}[Foliation]\label{Definition-Foliation}
A \emph{foliation} $\cF$ of dimension $k$ on a manifold $M^d$ (or
codimension $d-k$) is a partition of $M$ on injectively immersed
connected $C^1$-submanifolds tangent to a continuous subbundle $E$
of $TM$ satisfying:
\begin{itemize}
\item[-] For every $x\in M$ there exists a neighborhood $U$ and a
continuous homeomorphism $\varphi: U \to \RR^{k} \times \RR^{d-k}$
such that for every $y \in \RR^{d-k}$:

$$ L_y = \varphi^{-1}(\RR^k \times \{y\}) $$

\noindent is a connected component of $L \cap U$ where $L$ is an
element of the partition $\cF$.
\end{itemize}
 \finobs
\end{defi}

In most of the texts about foliations, this notion refers to a
$C^0$-foliation with $C^1$-leaves (or foliations of class
$C^{1,0+}$ in \cite{CandelConlon}).

In dynamical systems, particularly in the theory of Anosov
diffeomorphisms, flows and or partially hyperbolic systems, this
notion is the best suited since it is the one guarantied by these
dynamical properties (see Theorem
\ref{Teorema-VariedadEstableFuerte}, this notion corresponds to a
$C^1$-lamination of the whole manifold).

\begin{notacion}
We will denote $\cF(x)$ to the \emph{leaf} (i.e. element of the
partition) of the foliation $\cF$ containing $x$. Given a
foliation $\cF$ of a manifold $M$, we will always denote as
$\tilde \cF$ to the lift of the foliation $\cF$ to the universal
cover $\tilde M$ of $M$. \finobs
\end{notacion}

We will say that a foliation is \emph{orientable} if there exists
a continuous choice of orientation for the subbundle $E\en TM$
which is tangent to $\cF$. Similarly, we say that the foliation is
\emph{transversally orientable} if there exists a continuous
choice of orientation for the subbundle $E^\perp \en TM$
consisting of the orthogonal bundle to $E$. Notice that if $M$ is
orientable, then the fact that $E$ is orientable implies that
$E^\perp$ is also orientable.

Given a foliation $\cF$ of a manifold $M$, one can always consider
a finite covering of $M$ and $\cF$ in order to get that the lifted
foliation is both orientable and transversally orientable.

We remark that sometimes, the definition of a foliation is given
in terms of atlases on the manifold, we state the following
consequence of our definition:

\begin{prop}\label{Proposicion-FoliacionAtlas}
Let $M$ be a $d$-dimensional manifold and $\cF$ a $k$-dimensional
foliation of $M$. Then, there exists a $C^0$-atlas
$\{(\varphi_i,U_i) \}$ of $M$ such that: \bi \item[-] $\varphi_i:
U_i \to \RR^{k} \times \RR^{d-k}$ is a homeomorphism. \item[-] If
$U_i \cap U_j \neq \emptyset$ one has that $\varphi_i \circ
\varphi_j^{-1}: \varphi_j(U_i \cap U_j) \to \RR^k \times
\RR^{d-k}$ is of the form $\varphi_i \circ \varphi_j^{-1}(x,y) =
(\varphi^1_{ij} (x,y), \varphi_{ij}^2 (y))$. Moreover, the maps
$\varphi^1_{ij}$ are $C^1$. \item[-] The preimage by $\varphi_i$
of a set of the form $\RR^k \times \{y\}$ is contained in a leaf
of $\cF$. \ei
\end{prop}

An important tool in foliation theory is the concept of
\emph{holonomy} (compare with subsection
\ref{SubSection-HolonomiaFliasLocales}). Given two points $x,y$ in
a leaf $\cF(x)$ of a foliation $\cF$ one can consider transverse
disks $\Sigma_x$ and $\Sigma_y$ of dimension $d-k$ and a curve
$\gamma_{x,y}$ joining these two points and contained in $\cF(x)$.
It is possible to ``lift'' this curve to the nearby leaves (by
using the atlas given by Proposition
\ref{Proposicion-FoliacionAtlas}) to define a continuous map from
a neighborhood of $x$ in $\Sigma_x$ to a neighborhood of
$\Sigma_y$. When the curve is understood from the context (for
example, when the foliation $\cF$ is one dimensional) we denote
this map as:

$$ \Pi^{\cF}_{x,y} : U \en \Sigma_x \to \Sigma_y$$

These neighborhoods where one can define the maps may depend on
the curve $\gamma$, however, it can be seen that given two curves
$\gamma_{x,y}$ and $\tilde \gamma_{x,y}$ which are homotopic
inside $\cF(x)$ the maps defined coincide in the intersection of
their domains, thus, a homotopy class of curves defines a germ of
maps from $\Sigma_x$ to $\Sigma_y$.

Considering the curves joining $x$ to itself inside $\cF(x)$ one
can thus define the following map:

$$ \Hol: \pi_1(\cF(x)) \to \Germ(\Sigma_x) $$

Which can be seen to be a group morphism. This is useful in some
cases in order to see that certain leaves are not simply
connected. We call \emph{holonomy group} of a leaf $L$ to the
image of the morphism $\Hol$ restricted to the fundamental group
of $L$.

An important use of the knowledge of holonomy is given by the
following theorem:

\begin{teo}[Reeb's stability theorem]\label{Teorema-EstabilidadChicaDeReeb}
Let $\cF$ be a foliation of $M^d$ of dimension $k$ and let $L$ be
a compact leaf whose holonomy group is trivial. Then, there exists
a neighborhood $U$ of $L$ saturated by $\cF$ such that every leaf
in $U$ is homeomorphic to $L$. Moreover, the neighborhood $U$ can
be chosen arbitrarily small.
\end{teo}

Being an equivalence relation, we can always make a quotient from
the foliation and obtain a topological space (which is typically
non-Hausdorff) called the \emph{leaf space} endowed with the
quotient topology. For a foliation $\cF$ on a manifold $M$ we
denote the leaf space as $M/_\cF$.

\subsection{Generalities on codimension one foliations}\label{SubSection-CodimensionOneFoliations}

Very few is known about foliations in general. However, when the
foliation is of codimension $1$ there is quite a large general
theory (see in particular \cite{Hector}).

The first important property of codimension one foliations, is the
existence of a transverse foliation (which holds in more general
contexts, but for our definition of foliation is quite direct):

\begin{prop}\label{Proposition-ExisteFoliacionTransversal}
Given a codimension $1$ foliation $\cF$ of a compact manifold $M$
there exists a one-dimensional foliation $\cF^\perp$ transverse to
$\cF$. Moreover, the foliations $\cF$ and $\cF^\perp$ admit a
\emph{local product structure}, this means that for every $\eps>0$
there exists $\delta>0$ such that:
\begin{itemize}
\item[-] Given $x,y \in M$ such that $d(x,y)<\delta$ one has that
$\cF_\eps (x) \cap \cF^\perp_\eps (y)$ consists of a unique point.
Here, $\cF_\eps(x)$ and $\cF^\perp_\eps(y)$ denote the local
leaves\footnote{More precisely, $\cF_\eps(x) = cc_x (\cF(x)\cap
B_\eps(x))$ and $\cF^\perp_\eps(y) = cc_y (\cF^\perp(y) \cap
B_\eps(y))$. As defined in the Notation section, $cc_x(A)$ denotes
the connected component of $A$ containing $x$.} of the foliations
in $B_\eps(x)$ and $B_\eps(y)$.
\end{itemize}
\end{prop}

\dem Assume first that $\cF$ is transversally orientable. To prove
the existence of a one dimensional foliation transverse to $\cF$
consider $E$ the continuous subbundle of $TM$ tangent to $\cF$.
Now, there exists an arbitrarily narrow cone $\cE^\perp$
transverse to $E$ around the one dimensional subbundle $E^\perp$
(the orthogonal subbundle to $E$).

In $\cE^\perp$ there exists a $C^1$ subbundle $F$. Since $E^\perp$
is orientable, so is $F$ so we can choose a $C^1$-vector field
without singularities inside $F$ which integrates to a $C^1$
foliation which will be of course transverse to $\cF$.

If $\cF$ is not transversally orientable, one can choose a
$C^1$-line field inside the cone field and taking the double cover
construct a $C^1$-vector field invariant under deck
transformations. This gives rise to an orientable one dimensional
foliation transverse to the lift of $\cF$ which projects to a
non-orientable one transverse to $\cF$.

By compactness of $M$ one checks that the local product structure
holds.
\lqqd

\begin{obs}[Uniform local product structure]\label{Remark-UniformLocalProductStructure}
There exists $\eps>0$ such that for every $x\in M$ there exists
$V_x \en M$ containing $B_\eps(x)$ admitting $C^0$-coordinates
$\psi_x : V_x \to [-1,1]^{d-1} \times [-1,1]$ such that:
\begin{itemize}
\item[-] $\psi_x$ is a homeomorphism and $\psi_x(x)=(0,0)$.
\item[-] $\psi_x$ sends connected components in $V_x$ of leaves of
$\cF$ into sets of the form $[-1,1]^{d-1} \times \{t\}$. \item[-]
$\psi_x$ sends connected components in $V_x$ of leaves of
$\cF^\perp$ into sets of the form $\{s\} \times [-1,1]$.
\end{itemize}
In fact, choose $\eps_0>0$ and consider $\delta$ as in the
statement about local product structure of Proposition
\ref{Proposition-ExisteFoliacionTransversal}, we get that if
$d(x,y)< \delta$ then $\cF_{\eps_0}(x) \cap \cF_{\eps_0}^\perp(y)$
consists of exactly one point. Given a point $z \in M$ we can then
find a continuous map:

$$ \tilde \psi_x : \cF_{\frac \delta 2}(x) \times \cF_{\frac \delta 2}^\perp (x) \to M $$

\noindent such that $\tilde \psi_x(a,b)$ is the unique point of
intersection of $\cF_{\eps_0}(x) \cap \cF_{\eps_0}^\perp(y)$. By
the invariance of domain theorem (see \cite{Hatcher}) we obtain
that $\tilde \psi_x$ is a homeomorphism over its image $\tilde
V_x$. Let $\eps$ be the Lebesgue number of the covering of $M$ by
the open sets $\tilde V_x$. For every point $z \in M$ there exists
$x$ such that $B_\eps(z) \en \tilde V_x$. Consider a homeomorphism
$\nu_x: [-1,1]^{d-1} \times [-1,1] \to \cF_{\frac \delta 2}(x)
\times \cF_{\frac \delta 2}^\perp (x)$ preserving the coordinates,
then it is direct to check that the inverse of $\tilde \psi_x$
composed with $\nu_x$ is the desired $\psi_z$.
 \finobs

\end{obs}

In codimension $1$ the behaviour of the transversal foliation may
detect non-simply connected leafs, this is the content of this
well known result of Haefliger which can be thought of a precursor
of the celebrated Novikov's theorem:

\begin{prop}[Haefliger Argument]\label{Proposicion-ArgumentoHaefliger}
Consider a codimension one foliation $\cF$ of a compact manifold
$M$. Let $\tilde \cF$ and $\tilde \cF^\perp$ be the lifts to the
universal cover of both $\cF$ and the transverse foliation given
by Proposition \ref{Proposition-ExisteFoliacionTransversal}.
Assume that there exists a leaf of $\tilde \cF^\perp$ that
intersects a leaf of $\tilde \cF$ in more than one point, then,
$\tilde \cF$ has a non-simply connected leaf.
\end{prop}

This can be restated in the initial manifold by saying that if
there exists a closed curve in $M$ transverse to $\cF$ which is
nullhomotopic, then there exists a leaf of $\cF$ such that its
fundamental group does not inject in the fundamental group of $M$.

This result was first proven by Haefliger for $C^2$ foliations and
then extended to general $C^0$-foliations by Solodov (see
\cite{Solodov}). The idea is to consider a disk bounding a
transverse curve to the foliation and making general position
arguments (the reason for which Haefliger considered the
$C^2$-case first) in order to have one dimensional foliation with
Morse singularities on the disk, classical Poincare-Bendixon type
of arguments then give the existence of a leaf of $\cF$ with
non-trivial holonomy.

Other reason for considering codimension one foliations is that
leaves with finite fundamental group do not only give a condition
on the local behaviour of the foliation but on the global one:

\begin{teo}[Reeb's global stability theorem]\label{Teorema-EstabilidadCompleta}
Let $\cF$ be a codimension one foliation on a compact manifold $M$
and assume that there is a compact leaf $L$ of $\cF$ with finite
fundamental group. Then, $M$ is finitely covered by a manifold
$\hat M$ admitting a fibration $p: \hat M \to S^1$ whose fibers
are homeomorphic to $\hat L$ which finitely covers $L$ and by
lifting the foliation $\cF$ to $\hat M$ we obtain the foliation
given by the fibers of the fibration $p$.
\end{teo}

\begin{cor}\label{Corolario-EstabilidadReebParaEsferadeHoja}
Let $\cF$ be a codimension one foliation of a $3$-dimensional
manifold $M$ having a leaf with finite fundamental group. Then,
$M$ is finitely covered by $S^2 \times S^1$ and the foliation
lifts to a foliation of $S^2 \times S^1$ by spheres.
\end{cor}

For a codimension one foliation $\cF$ of a manifold $M$, such that
the leafs in the universal cover are properly embedded, there is a
quite nice description of the leaf space $\tilde M /_{\tilde \cF}$
as a (possibly non-Hausdorff) one-dimensional manifold. When the
leaf space is homeomorphic to $\RR$ we say that the foliation is
$\RR$-\emph{covered} (see \cite{Calegari}).


\section{Codimension one foliations in dimension 3}\label{Seccion-Codimension1Dimension3}

\subsection{Reeb components and Novikov's Theorem}\label{SubSeccion-ComponentesDeReebTeoNovikov}

Consider the foliation of the band $[-1,1] \times \RR$ given by
the horizontal lines together with the graphs of the functions
$\displaystyle{x \mapsto \expon \left(\frac{1}{1-x^2}\right) + b}$
with $b\in \RR$.

Clearly, this foliation is invariant by the translation $(x,t)
\mapsto (x,t+1)$ so that it defines a foliation on the annulus
$[-1,1] \times S^1$ which we call \emph{Reeb annulus}.

In a similar way, we can define a two-dimensional foliation on
$\DD^2 \times \RR$ given by the cylinder $\partial \DD^2 \times
\RR$ and the graphs of the maps $\displaystyle{(x,y) \mapsto
\expon \left(\frac{1}{1-x^2-y^2}\right) + b}$.

\begin{defi}[Reeb component]\label{Definition-ReebComponent}
Any foliation of $\DD^2 \times S^1$ homeomorphic to the foliation
obtained by quotienting the foliation defined above by translation
by $1$ is called a \emph{Reeb component}. \finobs
\end{defi}

Another important component of $3$-dimensional foliations are
dead-end components. They consist of foliations of $\TT^2 \times
[-1,1]$ such that any transversal which enters the boundary cannot
leave the manifold again. An example would be the product of a
Reeb annulus with the circle.

\begin{defi}[Dead-end component]\label{Definicion-DeadEndComponent}
A foliation of $\TT^2 \times [-1,1]$ such that no transversal can
intersect both boundary components is called a \emph{dead-end
component}. \finobs
\end{defi}

Novikov's theorem (see \cite{Novikov}) was proved for
$C^2$-foliations by the same reason as with Haefliger's argument
(see Proposition \ref{Proposicion-ArgumentoHaefliger}), with the
techniques of Solodov (the techniques of Solodov can be simplified
with the existence of a transversal one-dimensional foliation) one
can prove it in our context (see \cite{CandelConlon} Theorems
9.1.3 and 9.1.4 and the Remark on page 286):

\begin{teo}[Novikov \cite{Solodov,CandelConlon}]\label{TeoNovikov}
Let $\cF$ be a (transversally oriented) codimension one foliation
on a $3$-dimensional compact manifold $M$ and assume that one of
the following holds:
\begin{itemize}
\item[-] There exist a positively oriented closed loop transverse
to $\cF$ which is nullhomotopic, or, \item[-] there exist a leaf
$S$ of $\cF$ such that the fundamental group of $S$ does not
inject on the fundamental group of $M$. \item[-] $\pi_2(M) \neq
\{0\}$.
\end{itemize}
Then, $\cF$ has a Reeb component.
\end{teo}

\subsection{Reebless and taut foliations}\label{Subseccion-ReeblesTaut}

We will say that a (transversally oriented) codimension one
foliation of a $3$-dimensional manifold is \emph{Reebless} if it
does not contain Reeb components. Similarly, we say that a
Reebless foliation is \emph{taut} if it has no dead-end
components.

As a consequence of Novikov's theorem we obtain the following
corollary on Reebless foliations on $3$-manifolds which we state
without proof. We say that a surface $S$ embedded in a
$3$-manifold $M$ is \emph{incompressible} if the inclusion $\imath
: S \to M$ induces an injective morphism of fundamental groups.

\begin{cor}\label{CorolarioConsecuenciasReeb}
Let $\cF$ be a Reebless foliation on an orientable $3$-manifold
$M$ and $\cF^\perp$ a transversal one-dimensional foliation. Then,
\begin{itemize}
\item[(i)] For every $x\in \tilde M$ we have that $\tilde \cF(x)
\cap \tilde \cF^\perp (x) = \{x\}$. \item[(ii)] The leafs of
$\tilde \cF$ are properly embedded surfaces in $\tilde M$. In fact
there exists $\delta>0$ such that every euclidean ball $U$ of
radius $\delta$ can be covered by a continuous coordinate chart
such that the intersection of every leaf $S$ of $\tilde \cF$ with
$U$ is either empty of represented as the graph of a function
$h_S: \RR^2 \to \RR$ in those coordinates. \item[(iii)] Every leaf
of $\cF$ is incompressible. In particular, $\tilde M$ is either
$S^2\times \RR$ and every leaf is homeomorphic to $S^2$ or $\tilde
M = \RR^3$. \item[(iv)] For every $\delta>0$, there exists a
constant $C_\delta$ such that if $J$ is a segment of $\tilde
\cF^\perp$ then $\Vol(B_\delta(J)) > C_\delta \length(J)$.
\end{itemize}
\end{cor}

Notice that item (iii) implies that every leaf of $\tilde \cF$ is
simply connected, thus, if the manifold $M$ is not finitely
covered by $S^2 \times S^1$ then every leaf is homeomorphic to
$\RR^2$. Also, if $M$ is $\TT^3$ one can see that every closed
leaf of $\cF$ must be a two-dimensional torus (since for every
other surface $S$, the fundamental group $\pi_1(S)$ does not
inject in $\ZZ^3$, see \cite{ClasificacionSuperficies}).

The last statement of (iii) follows by the fact that the leaves of
$\cF$ being incompressible they lift to $\tilde M$ as simply
connected leaves. Applying Reeb's stability Theorem
\ref{Teorema-EstabilidadCompleta} we see that if one leaf is a
sphere, then the first situation occurs, and if there are no
leaves homeomorphic to $S^2$ then all leaves of $\tilde \cF$ must
be planes and by a result of Palmeira (\cite{Palmeira}) we obtain
that $\tilde M$ is homeomorphic to $\RR^3$.

We give a detailed proof of a similar result in the proof of
Corollary \ref{CorolarioConsecuenciasReeb2} (in particular item
(iv)) so we leave this result without proof.

\subsection{Reebles foliations of $\TT^3$}\label{SubSection-FoliacionesDeT3}

This subsection present results whose proofs are essentially
contained in \cite{BBI2}. We mention that there is a paper by
Plante \cite{Plante} which is also based on previous developments
by Novikov (\cite{Novikov}) and Roussarie (\cite{Roussarie}) which
proves essentially the same results for foliations of class $C^2$
and extends it to manifolds with almost solvable fundamental
group. There exists a result of Gabai (\cite{Gabai}) which proves
the result of Roussarie for the foliations of lower regularity.

We consider a codimension one foliation $\cF$ of $\TT^3$ which is
transversally oriented and $\cF^\perp$ a one dimensional
transversal foliation given by Proposition
\ref{Proposition-ExisteFoliacionTransversal} (the only thing we
require to $\cF^\perp$ is to be transversal to $\cF$ in order to
satisfy Remark \ref{Remark-UniformLocalProductStructure}. We shall
assume throughout that $\cF$ has no Reeb components.

Let $p: \RR^3 \to \TT^3$ be the cannonical covering map whose deck
transformations are translations by elements of $\ZZ^3$.

Since $\RR^3$ is simply connected, we can consider an orientation
on $\tilde \cF^\perp$ (since $\cF^\perp$ is oriented, this
orientation is preserved under covering transformations).

Given $x\in \RR^3$ we get that $\tilde \cF^\perp(x) \setminus
\{x\}$ has two connected components which we call $\tilde
\cF^\perp_+(x)$ and $\tilde \cF^\perp_-(x)$ according to the
chosen orientation of $\tilde \cF^\perp$.

By Corollary \ref{CorolarioConsecuenciasReeb} (ii) we have that
for every $x\in \RR^3$ the set $\tilde \cF(x)$ is an embedded
surface in $\RR^3$. It is diffeomorphic to $\RR^2$ by Corollary
\ref{CorolarioConsecuenciasReeb} (iii). It is well known that this
implies that $\tilde \cF(x)$ separates $\RR^3$ into two connected
components\footnote{One can consider the usual one point
compactification of $\RR^3$ and apply the well known
Jordan-Brower's theorem. See for example \cite{Hatcher}
Proposition 2.B.1. This gives that the complement of $\tilde
\cF(x)$ consists of two connected components. The fact that
$\tilde \cF(x)$ is the boundary of both connected components is
proved as follows: First, since $\tilde \cF(x)$ is differentiable,
one can find a normal neighborhood which is an $I$-bundle
(homeomorphic to $\tilde \cF(x) \times [-1,1]$ such that the
homeomorphism maps $\tilde \cF(x)$ to $\tilde \cF(x) \times
\{0\}$). This implies that if a point of $\tilde \cF(x)$ is in the
boundary of a connected component of $\RR^3 \setminus \tilde
\cF(x)$ then the whole $\tilde \cF(x)$ must be in its boundary.
Now, assume that the boundary of one of the connected components
of $\RR^3 \setminus \tilde \cF(x)$ does not coincide with $\tilde
\cF(x)$. This implies that in fact the boundary of the connected
component is empty: there cannot be boundary points in the other
component since is open and contained in the complement and if one
point of $\tilde \cF(x)$ is in the boundary, then from the
argument above, one gets that the boundary coincides with $\tilde
\cF(x)$. This is a contradiction since this connected component
would be open and closed, thus the whole $\RR^3$ which is not  the
case. } whose boundary is $\tilde \cF(x)$. These components will
be denoted as $F_+(x)$ and $F_-(x)$ depending on whether they
contain $\tilde \cF^\perp_+(x)$ or $\tilde \cF^\perp_-(x)$.

Since covering transformations preserve the orientation and send
$\tilde \cF$ into itself, we have that:

$$F_{\pm}(x) + \gamma = F_{\pm}(x+\gamma) \qquad \forall \gamma \in \ZZ^3 $$

\subsubsection{} For every $x\in \RR^3$, we consider the following subsets of $\ZZ^3$ seen as deck transformations:

$$\Gamma_+(x) = \{ \gamma \in \ZZ^3 \ : \ F_+(x) + \gamma \en F_+(x) \} $$

$$\Gamma_-(x) = \{ \gamma \in \ZZ^3 \ : \ F_-(x) + \gamma \en F_-(x) \} $$

We also consider $\Gamma(x) = \Gamma_+ (x) \cup \Gamma_-(x)$.

\begin{obs}\label{RemarkLocalProductStructure}
There exists a uniform local product structure between $\tilde
\cF$ and $\tilde \cF^\perp$ (given by Remark
\ref{Remark-UniformLocalProductStructure}). Since leafs of $\tilde
\cF$ do not intersect there exists $\delta>0$ such that if two
points $x, y$ are at distance smaller than $\delta$, then either
$F_+(x) \en F_+(y)$ or $F_+(y) \en F_+(x)$ (the same for $F_-$).
In particular, if $d(\tilde \cF(x), \tilde \cF(x)+\gamma) <
\delta$, then $\gamma \in \Gamma(x)$. By Corollary
\ref{CorolarioConsecuenciasReeb} (i) we also know that if two
points are in the same leaf of $\tilde \cF^{\perp}$ and are at
distance smaller than $\delta$, then they are connected by a small
arc inside the leaf. \finobs
\end{obs}

\begin{lema}\label{RemarkPosibilidadesdelosFmas}  The following properties hold:
\begin{itemize}
\item[(i)] If both $F_+(x)\cap F_+(y) \neq \emptyset$ and $F_-(x)
\cap F_-(y) \neq \emptyset$ then, either $F_+(x) \en F_+(y)$ and
$F_-(y) \en F_-(x)$ or $F_+(y) \en F_+(x)$ and $F_-(x) \en
F_-(y)$. In both of this cases we shall say that $F_+(x)$ and
$F_+(y)$ are \emph{nested} (similar with $F_-$). \item[(ii)] If
$F_+(x) \cap F_+(y)=\emptyset$ then $F_+(y) \en F_-(x)$ and
$F_+(x) \en F_-(y)$. A similar property holds if $F_-(x) \cap
F_-(y)=\emptyset$. \item[(iii)] In particular, $F_+(x) \en F_+(y)$
if and only if $F_-(y) \en F_-(x)$.
\end{itemize}
\end{lema}

\dem We will only consider the case where $\tilde \cF(x) \neq
\tilde \cF(y)$ since otherwise the Lemma is trivially satisfied
(and case (ii) is not possible).

Assume that both $F_+(x) \cap F_+(y)$ and $F_-(x) \cap F_-(y)$ are
non-empty. Since $\tilde \cF(y)$ is connected and does not
intersect $\tilde \cF(x)$ we have that it is contained in either
$F_+(x)$ or $F_-(x)$. We can further assume that $\tilde \cF(y)
\en F_+(x)$ the other case being symmetric. In this case, we
deduce that $F_+(y) \en F_+(x)$: otherwise, we would have that
$F_-(x) \cap F_-(y) = \emptyset$. But this implies that $\tilde
\cF(x) \en F_-(y)$ and thus that $F_-(x) \en F_-(y)$ which
concludes the proof of (i).

To prove (ii) notice that if $F_+(x) \cap F_+(y) = \emptyset$ then
we have that $\tilde \cF(x) \en F_-(y)$ and $\tilde \cF(y) \en
F_-(x)$. This gives that both $F_+(x) \en F_-(y)$ and $F_+(y) \en
F_-(x)$ as desired.

Finally, if $F_+(x) \en F_+(y)$ we have that $F_-(x) \cap F_-(y)$
contains at least $F_-(y)$ so that (i) applies to give (iii).

See also \cite{BBI2} Lemma 3.8.

\lqqd

\begin{figure}[ht]
\begin{center}
\input{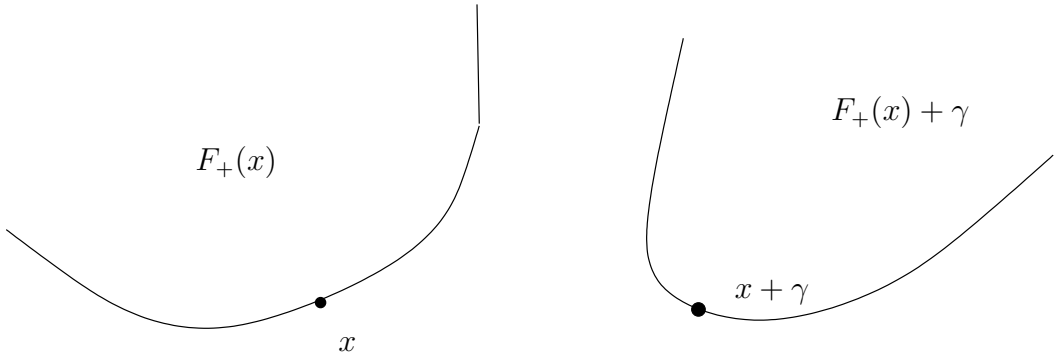}
\caption{\small{When $F_+(x)$ and $F_+(x)+\gamma$ are not
nested.}} \label{FiguraFmasOpciones}
\end{center}
\end{figure}


We can prove (see Lemma 3.9 of \cite{BBI2}):

\begin{lema}\label{LemaGamaEsSubGrupoDeZ3} For every $x\in \RR^3$ we have that $\Gamma(x)$ is a subgroup of $\ZZ^3$.
\end{lema}

\dem Consider $\gamma_1, \gamma_2 \in \Gamma_+(x)$.  Since
$\gamma_1 \in \Gamma_+(x)$ we have that $F_+(x) + \gamma_1 \en
F_+(x)$. By translating by $\gamma_2$ we obtain $F_+(x) + \gamma_1
+ \gamma_2 \en F_+(x) + \gamma_2$, but since $\gamma_2 \in
\Gamma_+(x)$ we have $F_+(x) + \gamma_1 + \gamma_2  \en F_+(x)$,
so $\gamma_1+ \gamma_2 \in \Gamma_+(x)$. This shows that
$\Gamma_+(x)$ is a semigroup.

Notice also that if $\gamma \in \Gamma_+(x)$ then $F_+(x) + \gamma
\en F_+(x)$, by substracting $\gamma$ we obtain that $F_+(x) \en
F_+(x) - \gamma$ which implies that $F_-(x) - \gamma \en F_-(x)$
obtaining that $-\gamma \in \Gamma_-(x)$. We have proved that
$-\Gamma_+(x)= \Gamma_{-}(x)$.

It then remains to prove that if $\gamma_1, \gamma_2 \in
\Gamma_+(x)$, then $\gamma_1- \gamma_2 \in \Gamma(x)$.

Since $F_+(x) + \gamma_1 + \gamma_2$ is contained in both $F_+(x)
+ \gamma_1$ and $F_+(x) + \gamma_2$ we have that
$$(F_+(x)+\gamma_1) \cap (F_+(x)+\gamma_2) \neq \emptyset.$$

By Lemma \ref{RemarkPosibilidadesdelosFmas} (iii) we have that
both $F_-(x)+\gamma_1$ and $F_-(x)+\gamma_2$ contain $F_-(x)$ so
they also have non-empty intersection.

Using Lemma \ref{RemarkPosibilidadesdelosFmas} (i), we get that
$F_+(x)+\gamma_1$ and $F_+(x)+\gamma_2$ are nested and this
implies that either $\gamma_1-\gamma_2$ or $\gamma_2-\gamma_1$ is
in $\Gamma_+(x)$ which concludes. \lqqd

 We close this subsection
by proving the following theorem which provides a kind of
classification of Reebless foliations in $\TT^3$:

\begin{teo}\label{PropObienHayTorosObienTodoEsLindo}
Let $\cF$ be a Reebless foliation of $\TT^3$. Then, there exists a
plane $P \en \RR^3$ and $R>0$ such that every leaf of $\tilde \cF$
lies in an $R$-neighborhood of a translate of $P$. Moreover, one
of the following conditions hold:
\begin{itemize}
\item[(i)] Either for every $x \in \RR^3$ the $R$-neighborhood of
$\tilde \cF(x)$ contains $P+x$, or, \item[(ii)] $P$ projects into
a two-dimensional torus and (if $\cF$ is orientable) there is a
dead-end component of $\cF$ (in particular, $\cF$ has a leaf which
is a two-dimensional torus).
\end{itemize}
\end{teo}

Notice that this theorem is mostly concerned with statements on
the universal cover so that orientability of the foliation is not
necessary. In fact, if one proves the theorem for a finite lift,
one obtains the same result since there cannot be embedded
incompressible Klein-bottles inside $\TT^3$ (see
\cite{ClasificacionSuperficies}). In option (ii), the only thing
we need is the fact that transversals remain at bounded distance
with the plane $P$ (which does not use orientability since it is a
statement on the universal cover).


\begin{figure}[ht]
\begin{center}
\input{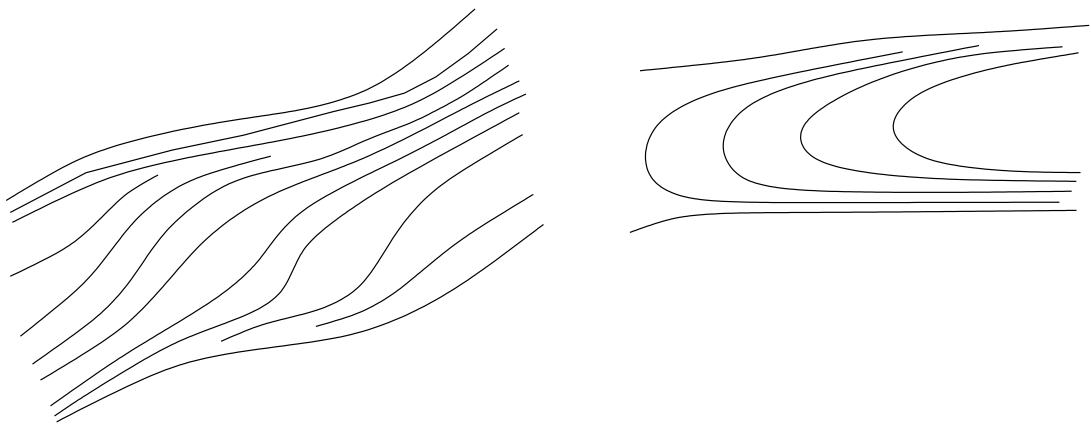}
\caption{\small{How the possibilities on $\tilde \cF$ look like.}}
\label{FiguraFoliOpciones}
\end{center}
\end{figure}


\dem By the remark above, we will assume throughout that the
foliation is orientable and transversally orientable. This allows
us to define as above the sets $F_\pm(x)$ for every $x$.

We define $G_+(x_0)=\bigcap_{\gamma \in \ZZ^3}\overline{F_+(x) +
\gamma}$ and $G_-(x_0)$ in a similar way.

First, assume that there exists $x_0$ such that
$G_+(x_0)=\bigcap_{\gamma \in \ZZ^3} \overline{F_+(x) + \gamma}
\neq \emptyset$ (see Lemma 3.10 of \cite{BBI2}). The case where
$G_-(x_0)\neq \emptyset$ is symmetric. The idea is to prove that
in this case we will get option (ii) of the theorem.

There exists $\delta>0$ such that given a point $z\in G_+(x_0)$ we
can consider a neighborhood $U_z$ containing $B_\delta(z)$ given
by Corollary \ref{CorolarioConsecuenciasReeb} (ii) such that:

\begin{itemize}
\item[-] There is a $C^1$-coordinate neighborhood $\psi_{z}: U_z
\to \RR^2 \to \RR$ such that for every $y\in U_z$ we have that
$\psi(\tilde \cF(y) \cap U_z)$ consists of the graph of a function
$h_y :\RR^2 \to \RR$ (in particular, it is connected).
\end{itemize}

Since $\tilde \cF$ and $\tilde \cF^\perp$ are orientable, we get
that we can choose the coordinates $\psi_z$ in order that for
every $y \in U_z$ we have that $\psi_z(F_+(y) \cap U)$ is the set
of points $(w,t)$ such that $t \geq h_y(w)$.

Notice that for every $\gamma \in \ZZ^3$ we have that $(F_+(x) +
\gamma) \cap U$ is either the whole $U$ or the upper part of the
graph of a function $h_{x+\gamma}:\RR^2 \to \RR$ in some
coordinates in $U$.

This implies that the intersection $G_+(x_0)$ is a $3$-dimensional
submanifold of $\RR^3$ (modeled in the upper half space) with
boundary consisting of leaves of $\tilde \cF$ (since the boundary
components are always locally limits of local leaves).

The boundary is clearly non trivial since $G_+(x_0)\en F_+(x_0)
\neq \RR^3$.

\begin{af}
If $G_+(x_0) \neq \emptyset$ then there exists plane $P$ and $R>0$
such that every leaf of $\tilde \cF(x)$ is contained in an
$R$-neighborhood of a translate of $P$ and whose projection to
$\TT^3$ is a two dimensional torus. Moreover, option (ii) of the
proposition holds.
\end{af}

\dem Since $G_+(x_0)$ is invariant under every integer
translation, we get that the boundary of $G_+(x_0)$ descends to a
closed surface in $\TT^3$ which is union of leaves of $\cF$.

By Corollary \ref{CorolarioConsecuenciasReeb} (iii) we get that
those leaves are two-dimensional torus whose fundamental group is
injected by the inclusion map.

This implies that they are at bounded distance of linear
embeddings of $\TT^2$ in $\TT^3$ and so their lifts lie within
bounded distance from a plane $P$ whose projection is a two
dimensional torus.

Since leafs of $\tilde \cF$ do not cross, the plane $P$ does not
depend on the boundary component. Moreover, every leaf of $\tilde
\cF(x)$ must lie within bounded distance from a translate of $P$
since every leaf of $\cF$ has a lift which lies within two given
lifts of some of the torus leafs.

Consider a point $x$ in the boundary of $G_+(x_0)$. We have that
$\tilde \cF(x)$ lies within bounded distance from $P$ from the
argument above.

Moreover, each boundary component of $G_+(x_0)$ is positively
oriented in the direction which points inward to the interior of
$G_+(x_0)$ (recall that it is a compact $3$-manifold with
boundary).

We claim that $\eta_z$ lies within bounded distance from $P$ for
every $z\in \tilde \cF(x)$ and $\eta_z$ positive transversal to
$\tilde \cF$. Indeed, if this is not the case, then $\tilde
\eta_z$ would intersect other boundary component of $G_+(x_0)$
which is impossible since the boundary leafs of $G_+(x_0)$ point
inward to $G_+(x_0)$ (with the orientation of $\tilde \cF^\perp$).

Now, consider any point $z \in \RR^3$, and $\eta_z$ a positive
transversal which we assume does not remain at bounded distance
from $P$. Then it must intersect some translate of $\tilde
\cF(x)$, and the argument above applies. This is a contradiction.

The same argument works for negative transversals since once a
leaf enters $(G_+(x_0))^c$ it cannot reenter any of its
translates. We have proved that $p(G_+(x_0))$ contains a dead end
component. This concludes the proof of the claim. \finobs

Now, assume that (ii) does not hold, in particular $G_\pm(x)
=\emptyset$ for every $x$. Then, for every point $x$ we have that
$$\bigcup_{\gamma \in \ZZ^3} (F_+(x) + \gamma) =
\bigcup_{\gamma\in \ZZ^3}(F_-(x) + \gamma) = \RR^3$$

\noindent As in Lemma 3.11 of \cite{BBI2} we can prove:

\begin{af} We have that $\Gamma(x)=\ZZ^3$ for every $x\in \RR^3$.
\end{af}
\dem If for some $\gamma_0 \notin \Gamma(x)$ one has that $F_+(x)
\cap (F_+(x) + \gamma_0) = \emptyset$ (the other possibility being
that $F_-(x) \cap (F_-(x)+ \gamma_0) = \emptyset)$) then, we claim
that for every $\gamma \notin \Gamma(x)$ we have that $F_+(x)\cap
( F_+(x)+ \gamma) =\emptyset$.

Indeed, by Lemma \ref{RemarkPosibilidadesdelosFmas} (i) if the
claim does not hold, there would exist $\gamma \notin \Gamma(x)$
such that

$$F_-(x) \cap (F_-(x)+ \gamma) = \emptyset \quad \text{and} \quad F_+(x) \cap (F_+(x) + \gamma) \neq \emptyset.$$

By Lemma \ref{RemarkPosibilidadesdelosFmas} (ii) we have:

\begin{itemize}
\item[-] $F_-(x) \en F_+(x)+\gamma$. \item[-] $F_+(x)+\gamma_0 \en
F_-(x)$.
\end{itemize}

Let $z \in F_+(x) \cap (F_+(x)+\gamma)$ then $z+\gamma_0$ must
belong both to $(F_+(x)+\gamma_0) \en (F_+(x)+ \gamma)$ and to
$(F_+(x)+\gamma+\gamma_0)$. By substracting $\gamma$ we get that
$z+\gamma_0-\gamma$ belongs to $F_+(x) \cap F_+(x) + \gamma_0$
contradicting our initial assumption. So, for every $\gamma \notin
\Gamma(x)$ we have that $F_+(x) \cap (F_+(x)+\gamma) = \emptyset$.

Now, consider the set $$U_+ (x) = \bigcup_{\gamma \in \Gamma(x)}
(F_+(x)+\gamma).$$

From the above claim, the sets $U_+(x) +\gamma_1$ and $U_+(x) +
\gamma_2$ are disjoint (if $\gamma_1-\gamma_2 \notin \Gamma(x)$)
or coincide (if $\gamma_1-\gamma_2 \in \Gamma(x)$).

Since these sets are open, and its translates by $\ZZ^3$ should
cover the whole $\RR^3$ we get by connectedness that there must be
only one. This implies that $\Gamma(x)=\ZZ^3$ and finishes the
proof of the claim. \finobs

Consider $\Gamma_0(x)= \Gamma_+(x) \cap \Gamma_-(x)$, the set of
translates which fix $\tilde \cF(x)$.

If $\rango(\Gamma_0(x))=3$, then $p^{-1}(p(\tilde \cF(x)))$
consists of finitely many translates of $\tilde \cF(x)$ which
implies that $p(\tilde \cF(x))$ is a closed surface of $\cF$. On
the other hand, the fundamental group of this closed surface
should be isomorphic to $\ZZ^3$ which is impossible since there
are no closed surfaces with such fundamental group
(\cite{ClasificacionSuperficies}). This implies that
$\rango(\Gamma_0(x))<3$ for every $x\in \RR^3$.

\begin{af}
For every $x\in \RR^3$ there exists a plane $P(x)$ and translates
$P_+(x)$ and $P_-(x)$ such that $F_+(x)$ lies in a half space
bounded by $P_+(x)$ and $F_-(x)$ lies in a half space bounded by
$P_-(x)$.
\end{af}

\dem Since $\rango(\Gamma_0(x))<3$ we can prove that $\Gamma_+(x)$
and $\Gamma_-(x)$ are half latices (this means that there exists a
plane $P \en \RR^3$ such that each one is contained in a half
space bounded by $P$).

The argument is the same as in Lemma 3.12 of \cite{BBI2} (and the
argument after that lemma).

Consider the convex hulls of $\Gamma_+(x)$ and $\Gamma_-(x)$. If
their interiors intersect one can consider 3 linearly independent
points whose coordinates are rational. These points are both
positive rational convex combinations of vectors in $\Gamma_+(x)$
as well as of vectors in $\Gamma_-(x)$. One obtains that
$\Gamma_0(x)=\Gamma_+(x) \cap \Gamma_-(x)$ has rank $3$
contradicting our assumption.

This implies that there exists a plane $P(x)$ separating these
convex hulls.

Consider $z\in \RR^3$ and let $\cO_+(z)= (z+ \ZZ^3) \cap F_+(x)$.
We have that $\cO_+(z) \neq \emptyset$ (otherwise $z \in G_-(x)$).
Moreover, $\cO_+(z) + \Gamma_+(x) \en \cO_+(z)$ because
$\Gamma_+(x)$ preserves $F_+(x)$. The symmetric statements hold
for $\cO_-(z) = (z+ \ZZ^3) \cap F_-(x)$.

We get that $\cO_+(z)$ and $\cO_-(z)$ are separated by a plane
$P_z$ parallel to $P(x)$. The proof is as follows: we consider the
convex hull $\cC\cO_+(z)$ of $\cO_+(z)$ and the fact that
$\cO_+(z) + \Gamma_+(x) \en \cO_+(z)$ implies that if $v$ is a
vector in the positive half plane bounded by $P(x)$ we have that
$\cC\cO_+(z) + v \en \cC \cO_+(z)$. The same holds for the convex
hull of $\cC\cO_-(z)$ and we get that if the interiors of
$\cC\cO_+(z)$ and $\cC\cO_-(z)$ intersect, then the interiors of
the convex hulls of $\Gamma_+(x)$ and $\Gamma_-(x)$ intersect
contradicting that $\rango(\Gamma_0(x)) < 3$.

Consider $\delta$ given by Corollary
\ref{CorolarioConsecuenciasReeb} (ii) such that every point $z$
has a neighborhood $U_z$ containing $B_\delta(z)$ and such that
$\tilde \cF(y) \cap U_z$ is connected for every $y \in U_z$.

Let $\{z_i\}$ a finite set $\delta/2$-dense in a fundamental
domain $D_0$. We denote as $P_{z_i}^+$ and $P_{z_i}^-$ de half
spaces defined by the plane $P_{z_i}$ parallel to $P(x)$
containing $\cO_+(z_i)$ and $\cO_-(z_i)$ respectively.

We claim that $F_+(x)$ is contained in the $\delta$-neighborhood
of $\bigcup_i P_{z_i}^+$ and the symmetric statement holds for
$F_-(x)$.

Consider a point $y\in F_+(x)$. We get that $\tilde \cF(y)$
intersects the neighborhood $U_y$ containing $B_\delta(y)$ in a
connected component and thus there exists a $\delta/2$-ball in
$U_y$ contained in $F_+(x)$. Thus, there exists $z_i$ and $\gamma
\in \ZZ^3$ such that $z_i + \gamma$ is contained in $F_+(x)$ and
thus $z_i + \gamma \in \cO_+(z_i) \en P_{z_i}^+$. We deduce that
$y$ is contained in the $\delta$-neighborhood of $P_{z_i}^+$ as
desired.

The $\delta$-neighborhood $H^+$ of $\bigcup_i P_{z_i}^+$ is a half
space bounded by a plane parallel to $P(x)$ and the same holds for
$H^-$ defined symmetrically. We have proved that $F_+(x) \en H^+$
and $F_-(x) \en H^-$. This implies that $\tilde \cF(x)$ is
contained in $H^+ \cap H^-$, a strip bounded by planes $P_+(x)$
and $P_-(x)$ parallel to $P(x)$ concluding the claim.

\finobs

We have proved that for every $x \in \RR^3$ there exists a plane
$P(x)$ and translates $P_+(x)$ and $P_-(x)$ such that $F_\pm (x)$
lies in a half space bounded by $P_{\pm}(x)$. Let $R(x)$ be the
distance between $P_+(x)$ and $P_-(x)$, we have that $\tilde
\cF(x)$ lies at distance smaller than $R$ from $P_+(x)$.

Now, we must prove that the $R(x)$-neighborhood of $\tilde \cF(x)$
contains $P_+(x)$. To do this, it is enough to show that the
projection from $\tilde \cF(x)$ to $P_+(x)$ by an orthogonal
vector to $P(x)$ is surjective. If this is not the case, then
there exists a segment joining $P_+(x)$ to $P_-(x)$ which does not
intersect $\tilde \cF(x)$. This contradicts the fact that every
curve from $F_-(x)$ to $F_+(x)$ must intersect $\tilde \cF(x)$.

Since the leaves of $\tilde \cF$ do not intersect, $P(x)$ cannot
depend on $x$. Since the foliation is invariant under integer
translations, we get (by compactness) that $R(x)$ can be chosen
uniformly bounded.

\lqqd

\begin{obs}\label{Remark-UnicidaddelPlanoP}
It is direct to show that for a given Reebless foliation $\cF$ of
$\TT^3$, the plane $P$ given by Theorem
\ref{PropObienHayTorosObienTodoEsLindo} is unique. Indeed, the
intersection of the $R$-neighborhoods of two different planes is
contained in a $2R$-neighborhood of their intersection line $L$.
If two planes would satisfy the thesis of Theorem
\ref{PropObienHayTorosObienTodoEsLindo} then we would obtain that
the complement of every leaf contains a connected component which
is contained in the $2R$-neighborhood of $L$. This is a
contradiction since as a consequence of Theorem
\ref{PropObienHayTorosObienTodoEsLindo} we get that there is
always a leaf of $\tilde \cF$ whose complement contains two
connected components each of which contains a half space of a
plane\footnote{Notice that if case (ii) holds this is direct from
the existence of a torus leaf and in case (i) this follows from
the statement of the last claim in the proof.}. \finobs
\end{obs}

We have used strongly the fact that $\tilde \cF$ is the lift of a
foliation in $\TT^3$ so that the foliation is invariant under
integer translations, this is why there is more rigidity in the
possible foliations of $\RR^3$ which are lifts of foliations on
$\TT^3$. See \cite{Palmeira} for a classification of foliations by
planes of $\RR^3$.

\subsection{Further properties of the foliations}\label{SubSection-FurtherPropertiesReeblesT3}

It is not hard to see that:

\begin{prop}\label{RemarkElCocienteDeLaFoliacionEsR}
Let $\cF$ be a Reebless foliation of $\TT^3$, if option (i) of
Theorem \ref{PropObienHayTorosObienTodoEsLindo} holds, then the
leaf space $\cL= \RR^3 /_{\tilde \cF}$ is homeomorphic to $\RR$.
\end{prop}

\dem The space of leafs $\cL$ with the quotient topology has the
structure of a (possibly non-Hausdorff) one-dimensional manifold
(see \cite{Calegari}). In fact, this follows directly from
Corollary \ref{CorolarioConsecuenciasReeb} as well as the fact
that it is simply connected as a one-dimensional manifold (see
Corollary \ref{CorolarioConsecuenciasReeb} (i)). To prove the
proposition is thus enough to show that it is Hausdorff.

We define an ordering in $\cL$ as follows

$$\tilde \cF(x) \geq \tilde \cF(y) \quad \text{if} \quad  F_+(x) \en F_+(y).$$

If option (i) of Theorem \ref{PropObienHayTorosObienTodoEsLindo}
holds, given $x,y$ we have that $F_+(x)\cap F_+(y) \neq \emptyset$
and $F_-(x) \cap F_-(y) \neq \emptyset$.

Then, Lemma \ref{RemarkPosibilidadesdelosFmas} (i) implies that
$F_+(x)$ and $F_+(y)$ are nested. In conclusion, we obtain that
the relationship we have defined is a total order.

Let $\tilde \cF(x)$ and $\tilde \cF(y)$ two different leaves of
$\tilde \cF$. We must show that they belong to disjoint open sets.

Without loss of generality, since it is a total order, we can
assume that $\tilde \cF(x) < \tilde \cF(y)$. This implies that
$F_+(y)$ is strictly contained in $F_+(y)$. On the other hand,
this implies that $F_-(y) \cap F_+(x) \neq \emptyset$, in
particular, there exists $z$ such that $\tilde \cF(x) < \tilde
\cF(z) < \tilde \cF(y)$.

Since the sets $F_+(z)$ and $F_-(z)$ are open and disjoint and we
have that $\tilde \cF(x) \en F_-(z)$ and $\tilde \cF(y) \in
F_+(z)$ we deduce that $\cL$ is Hausdorff as desired.

\lqqd

Now, since $\tilde \cF$ is invariant under deck transformations,
we obtain that we can consider the quotient action of $\ZZ^3 =
\pi_1(\TT^3)$ in $\cL$. For $[x]=\tilde \cF(x) \in \cL$ we get
that $\gamma \cdot [x] = [x+\gamma]$ for every $\gamma \in \ZZ^3$.

Notice that all leaves of $\cF$ in $\TT^3$ are simply connected if
and only if $\pi_1(\TT^3)$ acts without fixed point in $\cL$. In a
similar fashion, existence of fixed points, or common fixed points
allows one to determine the topology of leaves of $\cF$ in
$\TT^3$.

In fact, we can prove the following result:

\begin{prop}\label{Proposicion-SiPesToroHayHojaToro}
Let $\cF$ be a Reebless foliation of $\TT^3$. If the plane $P$
given by Theorem \ref{PropObienHayTorosObienTodoEsLindo} projects
into a two dimensional torus by $p$, then there is a leaf of $\cF$
homeomorphic to a two-dimensional torus.
\end{prop}

\dem Notice first that if option (ii) of Theorem
\ref{PropObienHayTorosObienTodoEsLindo} holds, the existence of a
torus leaf is contained in the statement of the theorem.

So, we can assume that option (i) holds. By considering a finite
index subgroup, we can further assume that the plane $P$ is
invariant under two of the generators of $\pi_1(\TT^3)\cong \ZZ^3$
which we denote as $\gamma_1$ and $\gamma_2$.

Since leaves of $\tilde \cF$ remain close in the Hausdorff
topology to the plane $P$ we deduce that the orbit of every point
$[x] \in \cL$ by the action of the elements $\gamma_1$ and
$\gamma_2$ is bounded.

Let $\gamma_3$ be the third generator, we get that its orbit
cannot be bounded since otherwise it would fix the plane $P$ since
it is a translation. So, the quotient of $\cL$ by the action of
$\gamma_3$ is a circle. We can make the group generated by
$\gamma_1$ and $\gamma_2$ act on this circle and we obtain two
commuting circle homeomorphisms with zero rotation number. This
implies they have a common fixed point which in turn gives us the
desired two-torus leaf of $\cF$.

\lqqd

Also, depending on the topology of the projection of the plane $P$
given by Theorem \ref{PropObienHayTorosObienTodoEsLindo} we can
obtain some properties on the topology of the leaves of $\cF$.

\begin{lema}\label{RemarkHojasCerradas}
Let $\cF$ be a Reebless foliation of $\TT^3$ and $P$ be the plane
given by Theorem \ref{PropObienHayTorosObienTodoEsLindo}.
\begin{itemize}
\item[(i)] Every closed curve in a leaf of $\cF$ is homotopic in
$\TT^3$ to a closed curve contained in $p(P)$. This implies in
particular that if $p(P)$ is simply connected, then all the leaves
of $\cF$ are also simply connected. \item[(ii)] If a leaf of $\cF$
is homeomorphic to a two dimensional torus, then, it is homotopic
to $p(P)$ (in particular, $p(P)$ is also a two dimensional torus).
\end{itemize}
\end{lema}

\dem To see (i), first notice that leafs are incompressible. Given
a closed curve $\gamma$ in a leaf of $\cF$ which is not
null-homotopic, we know that when lifted to the universal cover it
remains at bounded distance from a linear one-dimensional subspace
$L$. Since $\gamma$ is a circle, we get that $p(L)$ is a circle in
$\TT^3$. If the subspace $L$ is not contained in $P$ then it must
be transverse to it. This contradicts the fact that leaves of
$\cF$ remain at bounded distance from $P$.

To prove (ii), notice that a torus leaf $T$ which is
incompressible must remain close in the universal cover to a plane
$P_T$ which projects to a linear embedding of a $2$-dimensional
torus. From the proof of Theorem
\ref{PropObienHayTorosObienTodoEsLindo} and the fact that $\cF$ is
a foliation we get that $P_T=P$. See also the proof of Lemma 3.10
of \cite{BBI2}.

\lqqd



\section{Global product structure}\label{Section-GlobalProductStructure}

\subsection{Statement of results}

We start by defining global product structure:

\begin{defi}[Global Product Structure]\label{Defi-GPS}
Given two transverse foliations (this in particular implies that
their dimensions are complementary) $\cF_1$ and $\cF_2$ of a
manifold $M$ we say they admit a \emph{global product structure}
if given two points $x,y \in \tilde M$ the universal cover of $M$
we have that $\tilde \cF_1(x)$ and $\tilde \cF_2(y)$ intersect in
a unique point. \finobs
\end{defi}

Notice that by continuity of the foliations and invariance of
domain theorem (\cite{Hatcher}) we have that if a manifold has two
transverse foliations with a global product structure, then, the
universal cover of the manifold must be homeomorphic to the
product of $\tilde \cF_1(x) \times \tilde \cF_2(x)$ for any $x\in
\tilde M$. Indeed, the map

$$ \varphi: \tilde \cF_1(x) \times \tilde \cF_2(y) \to \tilde M \qquad \varphi(z,w) = \tilde \cF_1(z) \cap \tilde \cF_2(w) $$

\noindent is well defined, continuous (by the continuity of
foliations) and bijective (because of the global product
structure), thus a global homeomorphism.

In particular, leaves of $\tilde \cF_i$ must be simply connected
and all homeomorphic between them.

In general, it is a very difficult problem to determine whether
two foliations have a global product structure even if there is a
local one (this is indeed the main obstruction in the
clasification of Anosov diffeomorphisms of manifolds, see
\cite{FranksAnosov}).

However, in the codimension $1$ case we again have much more
information:

\begin{teo}[Theorem VIII.2.2.1 of \cite{Hector}]\label{Teorema-HectorHirsch}
Consider a codimension one foliation $\cF$ of a compact manifold
$M$ such that all the leaves of $\cF$ have trivial holonomy. Then,
for every $\cF^\perp$ foliation transverse to $\cF$ we have that
$\cF$ and $\cF^\perp$ have global product structure.
\end{teo}

This theorem applies for example when every leaf is compact and
without holonomy. The other important case (for this thesis) in
which this result applies is when every leaf of the foliation is
simply connected. Unfortunately, there will be some situations
where we will be needing to obtain global product structure but
not having neither all leaves of $\cF$ simply connected nor that
the foliation lacks of holonomy in all its leaves.

We will instead use the following quantitative version of the
previous result which does not imply it other than it the
situations we will be needing it. The following theorem was proved
in \cite{Pot3dimPHisotopic} and we believe it simplifies certain
parts of the previous theorem (at least for the non-expert in the
theory of foliations and for the more restrictive hypothesis we
include):

\begin{teo}\label{Teorema-EstructuraProductoGlobalMIO}
Let $M$ be a compact manifold and $\delta>0$. Consider a set of
generators of $\pi_1(M)$ and endow $\pi_1(M)$ with the word length
for generators. Then, there exists $K>0$ such that if $\cF$ is a
codimension one foliation and $\cF^\perp$ a transverse foliation
such that:
\begin{itemize}
\item[-] There is a local product structure of size $\delta$
between $\cF$ and $\cF^\perp$ (see Remark
\ref{Remark-UniformLocalProductStructure}). \item[-] The leaves of
$\tilde \cF$ are simply connected and no element of $\pi_1(M)$ of
size less than $K$ fixes a leaf of $\tilde \cF$. \item[-] The leaf
space $\cL = \tilde M /_{\tilde \cF}$ is homeomorphic to $\RR$.
\item[-] The fundamental group of $M$ is abelian.
\end{itemize}
Then, $\cF$ and $\cF^\perp$ admit a global product structure.
\end{teo}

\subsection{Proof of Theorem \ref{Teorema-EstructuraProductoGlobalMIO}}

Notice that the hypothesis of the Theorem are stable by
considering finite lifts and the thesis is in the universal cover
so that we can (and we shall) assume that $\cF$ is both orientable
and transversally orientable.

The first step is to show that leaves of $\tilde \cF$ and $\tilde
\cF^\perp$ intersect in at most one point:

\begin{lema}\label{Lemma-LasHOJAScortanEnALoMasUnPunto}
For every $x\in \tilde M$ one has that $\tilde \cF(x)\cap \tilde
\cF^\perp(x) = \{x \}$.
\end{lema}

\dem Assume otherwise, then, by Proposition
\ref{Proposicion-ArgumentoHaefliger} (Haefliger argument) one
would obtain that there is a non-simply connected leaf of $\tilde
\cF$ a contradiction.
\lqqd

In $\cL = \tilde M /_{\tilde \cF}$ we can consider an ordering of
leafs (by using the ordering from $\RR$). We denote as $[x]$ to
the equivalence class in $\tilde M$ of the point $x$, which
coincides with $\tilde \cF(x)$.

The following condition will be the main ingredient for obtaining
a global product structure:

\begin{itemize}
\item[($\ast$)] For every $z_0\in \tilde M$ there exists $y^-$ and
$y^+ \in \tilde M$ verifying that $[y^-]< [z_0] < [y^+]$ and such
that for every $z_1, z_2\in \tilde M$ satisfying $[y^-]\leq [z_i]
\leq [y^+]$ ($i=1,2$) we have that $\tilde \cF^\perp(z_1) \cap
\tilde \cF(z_2) \neq \emptyset.$
\end{itemize}

We get

\begin{lema}\label{LemaCondicionAsterImplicaGPS} If property $(\ast)$ is satisfied, then $\tilde \cF$ and $\tilde \cF^\perp$ have a global product structure.
\end{lema}

\dem Consider any point $x_0 \in \tilde M$ and consider the set $G
= \{ z \in \tilde M \ : \ \tilde \cF^\perp (x_0) \cap \tilde
\cF(z) \neq \emptyset \}$. We have that $G$ is open from the local
product structure (Remark \ref{RemarkLocalProductStructure}) and
by definition it is saturated by $\tilde \cF$. We must show that
$G$ is closed and since $\tilde M$ is connected this would
conclude.

Now, consider $z_0 \in \overline{G}$, using assumption ($\ast$) we
obtain that there exists $[y^-]< [z] <  [y^+]$ such that every
point $z$ such that $[z^-]<[z]<[z^+]$ verifies that its unstable
leaf intersects both $\tilde \cF(y^-)$ and $\tilde \cF(y^+)$.

Since $z_0\in \overline{G}$ we have that there are points $z_n \in
G$ such that $z_n \to \tilde \cF(z_0)$.

We get that eventually,  $[y^-] < [z_k] < [y^+]$ and thus we
obtain that there is a point $y \in \tilde \cF^\perp(x_0)$
verifying that $[y^-]<[y]<[y^+]$. We get that every leaf between
$\tilde \cF(y^-)$ and $\tilde \cF(y^+)$ is contained in $G$ from
assumption ($\ast$). In particular, $z_0\in G$ as desired.
\lqqd

We must now show that property ($\ast$) is verified. To this end,
we will need the following lemma:

\begin{lema}\label{LemaConLongitudLCortoTodas}
There exists $K>0$ such that if $\ell>0$ is large enough, every
segment of $\cF^\perp(x)$ of length $\ell$ intersects every leaf
of $\cF$.
\end{lema}

We postpone the proof of this lemma to the next subsection
\ref{SubSectionPruebaLema}.

\demo{ of Theorem \ref{Teorema-EstructuraProductoGlobalMIO}} We
must prove that condition ($\ast$) is verified. We consider
$\delta$ given by the size of local product structure boxes (see
Remark \ref{RemarkLocalProductStructure}) and by Lemma
\ref{LemaConLongitudLCortoTodas} we get a value of $\ell>0$ such
that every segment of $\cF^\perp$ of length $\ell$ intersects
every leaf of $\cF$.

There exists $k>0$ such that every curve of length $k\ell$ will
verify that it has a subarc whose endpoints are $\delta$-close and
joined by a curve in $\cF^\perp$ of length larger than $\ell$ (so,
intersecting every leaf of $\cF$).

Consider a point $z_0\in \tilde M$ and a point $z\in \tilde
\cF(z_0)$. Let $\tilde \eta_z$ be the segment in $\tilde
\cF^\perp_+(z)$ of length $k \ell$ with one extreme in $z$. We can
project $\tilde \eta_z$ to $M$ and we obtain a segment $\eta_z$
transverse to $\cF$ which contains two points $z_1$ and $z_2$ at
distance smaller than $\delta$ and such that the segment from
$z_1$ to $z_2$ in $\eta_z$ intersects every leaf of $\cF$. We
denote $\tilde z_1$ and $\tilde z_2$ to the lift of those points
to $\tilde \eta_z$.

\begin{figure}[ht]
\begin{center}
\input{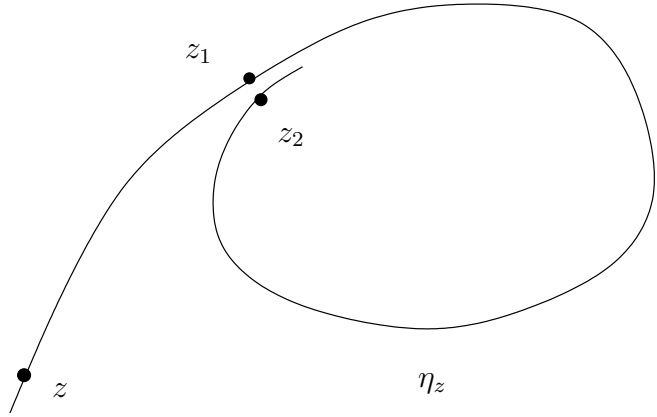}
\caption{\small{The curve $\eta_z$.}} \label{FiguraCurvaEta}
\end{center}
\end{figure}

Using the local product structure, we can modify slightly $\eta_z$
in order to create a closed curve $\eta_z'$ through $z_1$ which is
contained in $\eta_z$ outside $B_\delta(z_1)$, intersects every
leaf of $\cF$ and has length smaller than $k\ell + \delta$.

We can define $\Gamma_+$ as the set of elements in $\pi_1(M)$
which send the half space bounded by $\cF(x)$ in the positive
orientation into itself.

Since $\eta_z$ essentially contains a loop of length smaller than
$k\ell + \delta$ we have that $\tilde \eta_z$ connects $[z_0]$
with $[\tilde z_1 + \gamma]$ where $\gamma$ belongs to $\Gamma_+$
and can be represented by a loop of length smaller than $k\ell +
\delta$. Moreover, since from $z$ to $\tilde z_1$ there is a
positively oriented arc of $\tilde \cF^\perp$ we get that
$[z_0]=[z] \leq [\tilde z_1]$ (notice that it is possible that
$z=\tilde z_1$).

This implies that $[\tilde z_1 + \gamma] \geq [z_0+\gamma] >
[z_0]$, where the last inequality follows from the fact that the
loop is positively oriented and non-trivial (recall that by Lemma
\ref{Lemma-LasHOJAScortanEnALoMasUnPunto} a curve transversal to
$\tilde \cF$ cannot intersect the same leaf twice).

Notice that there are finitely many elements in $\Gamma_+$ which
are represented by loops of length smaller than $k\ell + \delta$.
This is because the fundamental group is abelian so that deck
transformations are in one to one correspondence with free
homotopy classes of loops.

The fact that there are finitely many such elements in $\Gamma_+$
implies the following: There exists $\gamma_0 \in \Gamma_+$ such
that for every $\gamma \in \Gamma_+$ which can be represented by a
positively oriented loop transverse to $\cF$ of length smaller
than $k\ell + \delta$, we have

$$  [z_0] < [z_0+\gamma_0] \leq [z_0+\gamma]$$

We have obtained that for $y^+= z_0 + \gamma_0$ there exists $L=
k\ell>0$ such that for every point $z\in \tilde \cF(z_0)$  the
segment of $\tilde \cF^\perp_+(z)$ of length $L$ intersects
$\tilde \cF(y^+)$.

This defines a continuous injective map from $\tilde \cF(z_0)$ to
$\tilde \cF(y^+)$ (injectivity follows from Lemma
\ref{Lemma-LasHOJAScortanEnALoMasUnPunto}). Since the length of
the curves defining the map is uniformly bounded, this map is
proper and thus, a homeomorphism. The same argument applies to any
leaf $\tilde \cF(z_1)$ such that $[z_0] \leq [z_1] \leq [y^+]$.

For any $z_1$ such that $[z_0]\leq [z_1] \leq [y^+]$ we get that
$\tilde \cF^\perp(z_1)$ intersects $\tilde \cF(z_0)$. Since the
map defined above is a homeomorphism, we get that also $\tilde
\cF^\perp(z_0) \cap \tilde \cF(z_1) \neq \emptyset$.

A symmetric argument allows us to find $y^-$ with similar
characteristics. Using the fact that intersecting with leaves of
$\tilde \cF^\perp$ is a homeomorphism between any pair of leafs of
$\tilde \cF$ between $[y^-]$ and $[y^+]$ we obtain ($\ast$) as
desired.

Lemma \ref{LemaCondicionAsterImplicaGPS} finishes the proof. \lqqd

\subsection{Proof of Lemma \ref{LemaConLongitudLCortoTodas}}\label{SubSectionPruebaLema}

We first prove the following Lemma which allows us to bound the
topology of $M$ in terms of coverings of size $\delta$. Notice
that we are implicitly using that $\pi_1(M)$ as before to be able
to define a correspondence between (free) homotopy classes of
loops with elements of $\pi_1(M)$.

\begin{lema}\label{LemmaCotadePiunoenfunciondedeltacubrimiento}
Given a covering $\{V_1, \ldots , V_n\}$ of $M$ by contractible
open subsets there exists there exists $K>0$ such that if $\eta$
is a loop in $M$ such that it intersects each $V_i$ at most
once\footnote{More precisely, if $\eta$ is $\eta:[0,1] \to M$ with
$\eta(0)=\eta(1)$ this means that $\eta^{-1}(V_i)$ is connected
for every $i$.}, then $[\eta] \in \pi_1(M)$ has norm less than
$K$.
\end{lema}

\dem We can consider the lift $p^{-1}(V_i)$ to the universal cover
of each $V_i$ and we have that each connected component of
$p^{-1}(V_i)$ has bounded diameter since they are simply connected
in $M$. Let $C_V> 0$ be a uniform bound on those diameters.

Let $K$ be such that every loop of length smaller than $2 n C_V$
has norm less than $K$ in $\pi_1(M)$.

Now, consider a loop $\eta$ which intersects each of the open sets
$V_i$ at most once. Consider $\eta$ as a function $\eta: [0,1] \to
M$ such that $\eta(0)=\eta(1)$. Consider a lift $\tilde \eta:
[0,1] \to M$ such that $p(\tilde \eta(t)) = \eta(t)$ for every
$t$.

We claim that the diameter of the image of $\tilde \eta$ cannot
exceed $n C_V$. Otherwise, this would imply that $\eta$ intersects
some $V_i$ more than once. Now, we can homotope $\tilde \eta$
fixing the extremes in order to have length smaller than $2 n
C_V$. This implies the Lemma.

\lqqd

Given $\delta$ of the uniform local product structure (see Remark
\ref{Remark-UniformLocalProductStructure}), we say that a loop
$\eta$ is a $\delta$-\emph{loop} if it is transverse to $\cF$ and
consists of a segment of a leaf of $\cF^\perp$ together with a
curve of length smaller than $\delta$.

\begin{lema}\label{LemaAbiertosSaturadosDeF}
There exists $K\geq 0$ such that if $O \en M$ is an open
$\cF$-saturated set such that $O \neq M$. Then, there is no
$\delta$-loop contained in $O$.
\end{lema}

\dem For every point $x$ consider $N_x = B_\delta(x)$ with
$\delta$ the size of the local product structure boxes. We can
consider a finite subcover $\{N_{x_1}, \ldots, N_{x_n} \}$ for
which Lemma \ref{LemmaCotadePiunoenfunciondedeltacubrimiento}
applies giving $K>0$.

Consider, an open set $O \neq M$ which is $\cF$-saturated. We must
prove that $O$ cannot contain a $\delta$-loop.

Let $\tilde O_0$ a connected component of the lift $\tilde O$ of
$O$ to the universal cover $\tilde M$. We have that the boundary
of $\tilde O_0$ consists of leaves of $\tilde \cF$ and if a
translation $\gamma \in \pi_1(M)$ verifies that

$$ \tilde O_0 \cap \gamma \tilde O_0 \neq \emptyset $$

\noindent then we must have that $\tilde O_0 = \gamma + \tilde
O_0$. This implies that $\gamma$ fixes the boundary leafs of
$\tilde O_0$: This is because the leaf space $\cL= \tilde
M/_{\tilde \cF}$ is homeomorphic to $\RR$ so that $\tilde O_0$
being connected and $\tilde \cF$ saturated is an open interval of
$\cL$. Since deck transformations preserve orientation, if they
fix an open interval then they must fix the boundaries.

The definition of $K$ then guaranties that if an element $\gamma$
of $\pi_1(M)$ makes $\tilde O_0$ intersect with itself, then
$\gamma$ must be larger than $K$. In particular, any $\delta$-loop
contained in $O$ must represent an element of $\pi_1(M)$ of length
larger than $K$.

Now consider a $\delta$-loop $\eta$. Corollary
\ref{CorolarioConsecuenciasReeb} (i) implies that $\eta$ is in the
hypothesis of Lemma
\ref{LemmaCotadePiunoenfunciondedeltacubrimiento}. We deduce that
$\eta$ cannot be entirely contained in $O$ since otherwise its
lift would be contained in $\tilde O_0$ giving a deck
transformation $\gamma$ of norm less than $K$ fixing $\tilde O_0$
a contradiction.



\lqqd

\begin{cor}\label{CorolarioLoopSonTotales}
For the $K\geq 0$ obtained in the previous Lemma, if $\eta$ is a
$\delta$-loop then it intersects every leaf of $\cF$.
\end{cor}

\dem The saturation by $\cF$ of $\eta$ is an open set which is
$\cF$-saturated by definition. Lemma
\ref{LemaAbiertosSaturadosDeF} implies that it must be the whole
$M$ and this implies that every leaf of $\cF$ intersects $\eta$.
\lqqd

\demo{ of Lemma \ref{LemaConLongitudLCortoTodas}}  Choose $K$ as
in Lemma \ref{LemaAbiertosSaturadosDeF}. Considering a covering
$\{V_1, \ldots, V_k\}$ of $M$ by neighborhoods with local product
structure between $\cF$ and $\cF^\perp$ and of diameter less than
$\delta$.

There exists $\ell_0>0$ such that every oriented unstable curve of
length larger than $\ell_0$ traverses at least one of the $V_i's$.
Choose $\ell> (k+1)\ell_0$ and we get that every curve of length
larger than $\ell$ must intersect some $V_i$ twice in points say
$x_1$ and $x_2$. By changing the curve only in $V_i$ we obtain a
$\delta$-loop which will intersect the same leafs as the initial
arc joining $x_1$ and $x_2$.

Corollary \ref{CorolarioLoopSonTotales} implies that the mentioned
arc must intersect all leafs of $\cF$.

\lqqd

\subsection{Consequences of a global product structure}\label{SubSection-ConsecuenciasGPS}

We say that a foliation $\cF$ in a Riemannian manifold $M$ is
\emph{quasi-isometric} if there exists $a,b \in \RR$ such that for
every $x,y$ in a same leaf of $\cF$ we have:

$$ d_{\cF}(x,y) \leq a d(x,y) + b $$

\noindent where $d$ denotes the distance in $M$ induced by the
Riemannian metric and $d_\cF$ the distance induced in the leaves
of $\cF$ by restricting the metric of $M$ to the leaves of $\cF$.
See Section \ref{Section-AlmostDynamicalCoherence} for more
discussion on quasi-isometry.

\begin{prop}\label{PropositionQuasiIsometria}
Let $\cF$ be a codimension one foliation of $\TT^3$ and
$\cF^\perp$ a transverse foliation. Assume the foliations $\tilde
\cF$ and $\tilde \cF^\perp$ lifted to the universal cover have
global product structure. Then, the foliation $\tilde \cF^\perp$
is quasi-isometric. Moreover, if $P$ is the plane given by Theorem
\ref{PropObienHayTorosObienTodoEsLindo}, there exists a cone $\cE$
transverse to $P$ in $\RR^3$ and $K>0$ such that for every $x\in
\RR^3$ and $y\in \tilde \cF^\perp(x)$ at distance larger than $K$
from $x$ we have that $y-x$ is contained in the cone $\cE$.
\end{prop}

\dem Notice that the global product structure implies that $\cF$
is Reebless. Let $P$ be the plane given by Theorem
\ref{PropObienHayTorosObienTodoEsLindo}.

Consider $v$ a unit vector perpendicular to $P$ in $\RR^3$.

Global product structure implies that for every $N>0$ there exists
$L$ such that every segments of $\tilde \cF^\perp$ of length $L$
starting at a point $x$ intersect $P+x+Nv$. Indeed, if this was
not the case, we could find arbitrarily large segments of leaves
of $\tilde \cF^\perp$ not satisfying this property, by taking a
subsequence and translations such that the initial point is in a
bounded region, we obtain a leaf of $\tilde \cF^\perp$ which does
not intersect every leaf of $\tilde \cF$.

This implies quasi-isometry since having length larger than $kL$
implies that the endpoints are at distance at least $kN$.

Moreover, assuming that the last claim of the proposition does not
hold, we get a sequence of points $x_n$, $y_n$ such that the
distance is larger than $n$ and such that the angle between $y_n -
x_n$ with $P$ is smaller than $1/\|x_n - y_n\|$.

In the limit (by translating $x_n$ we can assume that it has a
convergent subsequence), we get a leaf of $\tilde \cF^\perp$ which
cannot intersect every leaf of $\tilde \cF$ contradicting the
global product structure.

\lqqd


\begin{subappendices}

\section{One dimensional foliations of $\TT^2$}\label{Section-DescDomT2}

This appendix is devoted to characterizing foliations by lines in
$\TT^2$ where the ideas of the previous sections can be developed
in an easier way. The goal is to leave this section independent
from the previous ones so that the reader can start by reading
this section (even if it will not be used in the remaining of the
text).


\subsection{Classification of foliations}

Let $\cF$ be a one-dimensional foliation of $\TT^2$ and
$\cF^{\perp}$ any transversal foliation.

We will consider $\tilde \cF$ and $\tilde \cF^{\perp}$ the lifts
of these foliations to $\RR^2$ with the canonical covering map $p:
\RR^2 \to \TT^2$.

Here, \emph{foliation} will mean a partition of $\TT^2$ by
continuous flow tangent to a continuous vector field without
singularities. This definition implies orientability, the proofs
can be easily adapted to cover the non-orientable case.

The first remark is a direct consequence of Poincare-Bendixon's
Theorem (see \cite{KH} 14.1.1):

\begin{prop}\label{Proposicion-PoincareBendixon}
All the leaves of $\tilde \cF$ and $\tilde \cF^{\perp}$ are
properly embedded copies of $\RR$.
\end{prop}

\dem By transversality and compactness, there are local product
structure boxes of uniform size (see Remark
\ref{RemarkEPLUniforme}).

Assume there is a leaf $\tilde \cF(x)$ which intersects a local
product structure box in more than one connected component.

This implies that there exists a leaf of $\tilde \cF$ which is a
circle by the argument of the proof of Poincare-Bendixon's
theorem. This gives a singularity for the foliation $\tilde
\cF^{\perp}$ a contradiction.

\lqqd









This allows us to prove the following:

\begin{prop}\label{Prop-ClasificacionFoliacionesEnT2}
Given a one dimensional orientable foliation $\cF$ of $\TT^2$ we
have that there exists a subspace $L \en \RR^2$ and $R>0$ such
that every leaf of $\tilde \cF$ lies in a $R$-neighborhood of a
translate of $L$. Moreover, one can choose $R$ such that one of
the following properties holds:
\begin{itemize}
\item[(i)] Either the $R$-neighborhood of every leaf of $\tilde
\cF$ contains a translate of $L$ or,

\item[(ii)] The line $L$ projects under $p$ to a circle and there
is no transversal to $\cF$ which intersects every leaf of $\cF$.
\end{itemize}
\end{prop}

See figure \ref{FiguraFoliOpciones}, in fact, in option (ii) it
can be proved that the foliation has a two-dimensional \emph{Reeb
component} (which to avoid confusions we prefer not to define).

\dem Consider a circle $C$ transverse to $\cF$. By Proposition
\ref{Proposicion-PoincareBendixon} we know that $C$ is not
null-homotopic. The existence of $C$ is not hard to show, it
suffices to consider a vector field transverse to $\cF$ and
perturb it in order to have a periodic orbit.

First, assume that $C$ does not intersect every leaf of $\cF$. By
saturating $C$ with the leaves of $\cF$ we construct $O$, an open
$\cF$ saturated set strictly contained in $\TT^2$.

We claim that the boundary of the open set consists of leaves of
$\cF$ homotopic to $C$: Consider $\tilde O_0$ a connected
component of the lift of $O$ to the universal cover. Since $C$ is
contained in $O$ we have that there is a connected component of
the lift of $C$ contained in $O$. This connected component joins a
point $x \in \tilde O_0$ with a point $x + \gamma$ where $\gamma
\in \ZZ^2$ represents $C$ in $\pi_1(\TT^2) \cong \ZZ^2$. This
implies that $\gamma$ fixes $\tilde O_0$ and in particular its
boundary components which must be leafs of $\tilde \cF$ the lift
of $\cF$ to the universal cover.

We have that the one-dimensional subspace $L$ generated by the
vector $\gamma$ in $\RR^2$ verifies that every leaf of $\tilde
\cF$ lies within bounded distance from a translate of $L$. Indeed,
this holds for the boundary leaves of $\tilde O_0$ and by
compactness and the fact that leaves do not cross one extends this
to every leaf.

Now, assume that there is no circle transverse to $\cF$ which
intersects every leaf of $\cF$. We claim that this means that
every transversal to $\cF$ must remain at bounded distance from
$L$ (which is not hard to prove implies (ii) of the Proposition).
Indeed, by the argument above, if this were not the case we would
find two closed leaves of $\cF$ which are not homotopic, a
contradiction with the fact that leaves of $\cF$ do not intersect.

So, we can assume that there exists a circle $C$ which is
transverse to $\cF$ and intersects every leaf of $\cF$. By
composing with a homeomorphism $H: \TT^2 \to \TT^2$ isotopic to
the identity we can assume that $C$ verifies that its lift is a
one-dimensional subspace $C_0$. If we prove (i) for $H(\cF)$ we
get (i) for $\cF$ too since $H$ is at bounded distance from the
identity in the universal cover.

By considering the first return map of the flow generated by $X$
to this circle $C$ we obtain a circle homeomorphism $h:C \to C$.
By the classical rotation number theory, when lifted to the
universal cover $\tilde C \cong \RR$ we have that the orbit of
every point by the lift $\tilde h$ has bounded deviation to the
translation by some number $\rho \in \RR$. We consider the
specific lift given orthogonal projecting into $C_0$ the point of
intersection of the flow line with the first integer translate of
$C_0$ it intersects.

We get that the line $L$ we are looking for is generated by $\rho
\gamma + \gamma^\perp$ where $\gamma \in \ZZ^2$ is a generator of
$C$ and $\gamma^\perp$ is the vector orthogonal to $\gamma$ whose
norm equals the distance of $C_0$ with its closest translate by a
vector of $\ZZ^2$.

\lqqd


\subsection{Global dominated splitting in surfaces}

The goal of this section is to show some of the ideas that will
appear in Chapter \ref{Capitulo-ParcialmenteHiperbolicos}) in a
simpler context.

Consider $f: \TT^2 \to \TT^2$ a $C^1$-diffeomorphism which is
partially hyperbolic. Without loss of generality, we will assume
that the splitting is of the form $T\TT^2 =E \oplus E^u$ where
both bundles are one-dimensional and $E^u$ is uniformly expanded.
By Theorem \ref{Teorema-VariedadEstableFuerte} there exists a
one-dimensional $f$-invariant foliation $\cF^u$ tangent to $E^u$.
Notice that $\cF^u$ cannot have leaves which are circles.

For simplicity, we will assume throughout that the bundles $E$ and
$E^u$ are oriented and their orientation is preserved by $Df$. It
is not hard to adapt the results here to the more general case.

We denote as $\tilde f$ to a lift of $f$ to the universal cover
$\RR^2$ and consider the foliation $\tilde \cF^u$ which is the
lift of $\cF^u$ to $\RR^2$.

Notice that in dimension $2$ being partially hyperbolic is
equivalent to having a global absolute dominated splitting. The
fact that a global dominated splitting implies the existence of a
continuous vector field on the manifold readily implies that in an
orientable surface, the surface must be $\TT^2$.

We will show the following:

\begin{teo}\label{Teorema-ApendiceSuperficies}
Let $f: \TT^2 \to \TT^2$ be a partially hyperbolic diffeomorphism
with splitting $T\TT^2 = E \oplus E^u$. Then: \begin{itemize}
\item[-]  $f$ is semiconjugated to an Anosov diffeomorphism of
$\TT^2$. \item[-] There is a unique quasi-attractor $\cQ$ of $f$.
\item[-] Every chain-recurrence class different from $\cQ$ is
contained in a periodic interval. \end{itemize}
\end{teo}

Since $E$ is a one dimensional bundle uniformly transverse to
$E^u$ we can approximate $E$ by a $C^1$-vector field $X$ which is
still transverse to $E^u$. The vector field $X$ will be integrable
and give rise to a foliation $\cF$ which may not be invariant but
is transverse to $E^u$. As usual, we denote as $\tilde \cF$ to the
lift of $\cF$ to $\RR^2$.

Using this foliation and the things we have proved we will be able
to show:

\begin{lema}\label{Lema-fphenT2esisotopicoaAnosov}
Let $f: \TT^2 \to \TT^2$ be a partially hyperbolic with splitting
$T\TT^2 =E \oplus E^u$, then $f_\ast: \RR^2 \to \RR^2$ is
hyperbolic.
\end{lema}

\dem Assume that $f_\ast$ has only eigenvalues of modulus smaller
or equal to $1$. Then, the diameter of compact sets grows at most
polynomially when iterated forward.

Consider an arc $\gamma$ of $\tilde \cF^u$ and we iterate it
forward. We get that the length of $\gamma$ grows exponentially
while its diameter only polynomially. In $\RR^2$ this implies that
there will be recurrence of $\gamma$ to itself and in particular,
we will obtain a leaf of $\tilde \cF^u$ which intersects a leaf of
$\tilde \cF$ twice, a contradiction with Proposition
\ref{Proposicion-PoincareBendixon}.

Since $f_\ast$ has determinant of modulus $1$ we deduce that
$f_\ast$ must be hyperbolic.

\lqqd

We deduce:

\begin{lema}\label{Lema-EPGparaT2}
There is a global product structure between $\tilde \cF$ and
$\tilde \cF^u$. In particular, $\tilde \cF^u$ is quasi-isometric.
\end{lema}

\dem We apply Proposition \ref{Prop-ClasificacionFoliacionesEnT2}
to $\cF^u$. We obtain a line $L^u$ which will be
$f_\ast$-invariant since $\tilde \cF^u$ is $\tilde f$-invariant.

Since $L^u$ is $f_\ast$-invariant and $f_\ast$ is hyperbolic, we
can deduce that $L^u$ does not project into a circle so that
option (ii) does not hold (recall that hyperbolic matrices have
irrational eigenlines).

Moreover, since $L^u$ must project into a dense line in $\TT^2$,
we get that the foliation $\cF^u$ has no holonomy, and this
implies by Theorem \ref{Teorema-HectorHirsch} that there is a
global product structure between $\cF^u$ and $\cF$.

Quasi-isometry follows exactly as in Proposition
\ref{PropositionQuasiIsometria}.

\lqqd

It is possible to give a proof of Theorem
\ref{Teorema-HectorHirsch} in the lines of the proof of our
Theorem \ref{Teorema-EstructuraProductoGlobalMIO}. In the case
$L\neq L^u$ where $L$ is the line given by Proposition
\ref{Prop-ClasificacionFoliacionesEnT2} for $\cF$ it is almost
direct that there is a global product structure. In the case
$L=L^u$ one must reach a contradiction finding a translation which
fixes the direction contradicting that $L^u$ is totally
irrational.

\begin{obs}\label{Remark-T2FuVaparadireccionInestable}
With the same argument as in Lemma
\ref{Lema-fphenT2esisotopicoaAnosov} we can also deduce that the
line $L^u$ given by Proposition
\ref{Prop-ClasificacionFoliacionesEnT2} for $\cF^u$ must be the
eigenline of $f_\ast$ corresponding to the eigenvalue of modulus
larger than $1$. Indeed, since $\tilde \cF^u$ is $\tilde
f$-invariant, then $L^u$ must be $f_\ast$-invariant. Moreover, if
$L^u$ corresponds to the stable eigenline of $f_\ast$ then the
diameter of forward iterates of an unstable arc cannot grow more
than linearly and the same argument as in Lemma
\ref{Lema-fphenT2esisotopicoaAnosov} applies. \finobs
\end{obs}

This allows us to show (notice that this also follows from
\cite{PujSamBrazil})

\begin{prop}\label{Prop-CoherenciaenT2}
There is a unique $f$-invariant foliation $\cF_E$ tangent to $E$.
\end{prop}

\dem The same argument as in Lemma \ref{Lema-EPGparaT2} gives that
any foliation tangent to $E$ must have a global product structure
with $\cF^u$ when lifted to the universal cover.

We first show there exists one $f$-invariant foliation. To do
this, we consider any foliation $\cF$ transverse to $E^u$ and we
iterate it backwards. Recall that the line $L^u$ close to the
foliation $\tilde \cF^u$ must be the eigenline of the unstable
eigenvalue of $f_\ast$ (see Remark
\ref{Remark-T2FuVaparadireccionInestable}).

Let $L$ be the one dimensional subspace given by Proposition
\ref{Prop-ClasificacionFoliacionesEnT2} for $\cF$. Since there is
a global product structure between $\tilde \cF$ and $\tilde \cF^u$
we get that $L \neq L^u$: Otherwise by considering points $x,y$ at
distance larger than $R$ in the direction orthogonal to $L^u$ we
would get that the leaves of $\tilde F(x)$ and $\tilde \cF^u(y)$
cannot intersect due to Proposition
\ref{Prop-ClasificacionFoliacionesEnT2}.

Iterating backwards, we get that the foliation $\tilde \cF_m=
\tilde f^{-m}(\tilde \cF)$ is close to the line $f_\ast^{-m}(L)$
that as $m\to \infty$ converges to $L^s$, the eigenline of the
stable eigenvalue of $f_\ast$.

Moreover, we can prove that there exists a constant $R$ such that
for every $m$ we have that every leaf of $\tilde \cF_m$ lies at
distance smaller than $R$ from $f_\ast^{-m}(L)$. Indeed, consider
$R \gg \frac{2 \lambda^u K_0}{\cos \alpha}$ where $K_0$ is the
$C^0$-distance from $\tilde f$ and $f_\ast$, $\lambda^u$ the
unstable eigenvalue of $f_\ast$ and $\alpha$ the angle between $L$
and $L^s$. We get that the $R$ neighborhood of any translate of
$L$ is mapped by $f_\ast^{-1}$ into an $\frac{\cos \alpha'}{\cos
\alpha}(\lambda^u)^{-1}R $-neighborhood of $f_\ast^{-1}(L)$ where
$\alpha'< \alpha$ is the angle between $f_\ast^{-1}(L)$ and $L^s$.
Since $R - \frac{\cos \alpha'}{\cos \alpha}(\lambda^u)^{-1}R
 > K_0$ from the choice of $R$ we
get that every leaf of $\tilde \cF_1 = \tilde f^{-1}(\cF)$ lies
within $R$-distance from a translate of $f_\ast^{-1}(L)$.
Inductively, we get that each $\tilde \cF_m$ lies within distance
smaller than $R$ from $f_\ast^{-m}(L)$.

We must show that there exists a unique limit for the backward
iterates of any leaf of $\tilde \cF$. Let us fix $R$ as above.

Let $x\in \RR^2$ and we consider $\tilde \cF_n(x) =\tilde
f^{-n}(\tilde \cF(f^n(x)))$. Notice that $\tilde \cF_n(x)$ is an
embedded line which intersects the unstable leaf of each point of
$\tilde \cF(x)$ in exactly one point. Assume there exists $z \in
\tilde \cF(x)$ such that in $\tilde \cF^u(z)$ there are two
different limit points $z_1$ and $z_2$ of the sequence $\tilde
\cF_n(x) \cap \tilde \cF^u(z)$. We have that forward iterates of
$\tilde f^k(z_i)$ must lie at distance smaller than $R$ from $L^s
+ \tilde f^k(x)$.

Consider $K>0$ such that if two points lie at distance larger than
$K$ inside an unstable leaf then they are at distance larger than
$R$ in the direction transverse to $L$. Then, by choosing $k$
large enough so that the length of the arc of unstable joining
$z_1$ and $z_2$ is larger than $K$ we get that $\tilde f^k(z_1)$
and $\tilde f^k(z_2)$ must be at distance larger than $R$ in the
direction transversal to $L^s$ contradicting the previous claim.

A similar argument implies that there cannot be two different
$f$-invariant foliations tangent to $E$ since both should remain
close to the stable eigenline of $f_\ast$.

\lqqd

Since $f_\ast$ is hyperbolic (Lemma
\ref{Lema-fphenT2esisotopicoaAnosov}) Proposition
\ref{PropExisteSemiconjugacion} gives that there exists a
semiconjugacy $H: \RR^2 \to \RR^2$ which is $C^0$-close to the
identity, is periodic and verifies that:

$$ H \circ \tilde f = f_\ast \circ H$$

We denote as $\tilde \cF_E$ to the lift of $\cF_E$ to the
universal cover, we can prove:

\begin{lema}\label{LemaInjectividadEnInestablesT2}
The preimage by $H$ of every point is contained in a leaf of $\tilde \cF_E$.
\end{lema}

\dem From Proposition \ref{Prop-ClasificacionFoliacionesEnT2} (and
the fact that $\tilde \cF_E$ is $\tilde f$-invariant and has a
global product structure with $\tilde \cF^u$) we get that every
leaf of $\tilde \cF_E$ lies within distance smaller than $R$ from
a translate of $L^s$.

Consider points $x,y$ lying in different leaves of $\tilde \cF_E$.
Now, consider $z= \tilde \cF_E(y) \cap \tilde \cF^u(x)$. We have
that the distance of $z$ and $x$ grows exponentially in the
direction of $L^u$. This implies that by iterating forward, the
distance between $\tilde \cF_E(x)$ and $\tilde \cF_E(y)$ must grow
also exponentially.

We conclude that $d(\tilde f^n(x), \tilde f^n(y)) \to \infty$ with
$n\to +\infty$.

Since $H$ is close to the identity and semiconjugates $f$ with
$f_\ast$ it cannot send $x$ and $y$ to the same point. \lqqd

In order to be able to apply Proposition
\ref{ProposicionMecanismo} we must show the following:

\begin{lema}\label{LemaFronteraCST2} There is a unique quasi-attractor $\cQ$ for $f$.
Moreover, every point $y$ which belongs to the boundary of a fiber
of $H$ relative to its leaf of $\tilde \cF_E$ belongs to $\cQ$.
\end{lema}

\dem By Conley's theorem (Theorem \ref{Teorema-Conley}), there
always exists a quasi-attractor $\cQ$ of $f$. Moreover, we have
seen that such quasi-attractors are saturated by unstable sets
(see \ref{Remark-QuasiAttractor}).

Consider any quasi-attractor $\cQ$. Let $y$ be a point which is in
the boundary of $H^{-1}(\{x\})$ relative to $\tilde \cF_E(y)$.
Given $\eps>0$, since $y$ is in the boundary of $H^{-1}(\{x\})$
relative to $\tilde \cF_E(y)$ we obtain that its image by $H$
cannot be contained in the unstable set of $x$ for $f_\ast$.

Iterating backwards we obtain a connected set of arbitrarily large
diameter in the direction of the stable eigenline of $f_\ast$.
This implies that for large $m$ we have that $\tilde
f^{-m}(B_\eps(y))$ intersects $p^{-1}(\cQ)$. This holds for every
$\eps >0$ so we get that for every $\eps>0$ we can construct an
$\eps$-pseudo-orbit from $y$ to $\cQ$. This implies that $y \in
\cQ$. Since $\cQ$ was arbitrary and quasi-attractors are disjoint
it also implies that there is a unique quasi-attractor.

 \lqqd

We are now able to give a proof of Theorem
\ref{Teorema-ApendiceSuperficies}:

\demo{ of Theorem \ref{Teorema-ApendiceSuperficies}} We have
proved that $f$ is semiconjugated to a linear Anosov
diffeomorphism of $\TT^2$ and that there is a unique
quasi-attractor.

The last claim of the Theorem follows from the fact that we have
proved that the conditions of Proposition
\ref{ProposicionMecanismo} are verified:

\begin{itemize}
\item[-] The partially hyperbolic set is the whole $\TT^2$ (so
that the maximal invariant set in $U$ is also the whole $\TT^2$).
\item[-] The semiconjugacy is the one given by $H$. It is
injective on unstable manifolds by Lemma
\ref{LemaInjectividadEnInestablesT2}. \item[-] Lemma
\ref{LemaFronteraCST2} implies that the frontier of fibers in
center stable leaves are all contained in the unique
quasi-attractor $\cQ$ of $f$. \item[-] Fibers of $H$ are invariant
under unstable holonomy (see the proof of Proposition
\ref{propiedades}).
\end{itemize}

This concludes.

 \lqqd

\end{subappendices}


   \chapter{Global partial hyperbolicity}\label{Capitulo-ParcialmenteHiperbolicos}

This chapter contains the main contributions of this thesis. In
Section \ref{Section-AlmostDynamicalCoherence} we present some
preliminaries and in particular we introduce the concept of
\emph{almost dynamical coherence} which is key in the study we
make in this chapter. In particular, this concept allows us to
prove the following result in Section \ref{Section-PHAnosov}:

\begin{teo*} Dynamical coherence is an open and closed property among partially hyperbolic
diffeomorphisms of $\TT^3$ isotopic to Anosov.
\end{teo*}

We remark that in general it is not known whether dynamical
coherence is neither an open nor a closed property. There are not
known examples where it is not open but in general, to obtain
opennes a technical condition is used (called
\emph{plaque-expansiveness}). See \cite{HPS,Berger}.

Dynamical coherence in the case where the center bundle has
dimension larger than one is a widely open subject. It has been
remarked by Wilkinson (\cite{Wilkinson}) that one can look at some
Anosov diffeomorphisms as partially hyperbolic ones which are not
dynamically coherent (see \cite{BuW2} for an overview of dynamical
coherence). The proof here presented relies heavily both in the
assumption of almost dynamical coherence and in being in the
isotopy class of an Anosov automorphism in $\TT^3$. Several
questions regarding generalizations of these kind of results pop
up even in dimension 3. The one which we believe to be more
important is the following:

\begin{quest} Is it true that every partially hyperbolic diffeomorphism in dimension $3$ is almost dynamically coherent?
\end{quest}

Assuming this question admits a positive answer, one could expect
to make some progress in the direction of classification of both
partially hyperbolic diffeomorphisms of $3$-manifolds, and more
importantly (due to Theorem \ref{Teorema-DPU}) of robustly
transitive diffeomorphisms in $3$-manifolds.

Another quite natural question to be posed, which is related, is
whether some manifolds can admit partially hyperbolic
diffeomorphisms but not strong partially hyperbolic ones. We prove
in Section \ref{Section-PHAnosov} that there are isotopy classes
in $\TT^3$ which admit partially hyperbolic diffeomorphisms but
not strongly partially hyperbolic ones. This poses the following
natural question for which we do not know the answer:

\begin{quest}
Let $M$ be a $3$-manifold different from $\TT^3$. Is every
partially hyperbolic diffeomorphism of $M$ isotopic to a strong
partially hyperbolic one?
\end{quest}

When we treat strong partially hyperbolic systems we are able to
obtain much stronger results concerning integrability. We prove in
Section \ref{SectionCoherenciaPHfuerte} the following:

\begin{teo*} Let $f: \TT^3 \to \TT^3$ be a strong partially hyperbolic diffeomorphism. Then:
\begin{itemize}
\item[-] Either $E^s \oplus E^c$ is tangent to a unique
$f$-invariant foliation, or, \item[-] there exists a $f$-periodic
two-dimensional torus $T$ which is tangent to $E^s\oplus E^c$ and
normally expanding. \end{itemize}
\end{teo*}

This result extends the results of \cite{BBI2} to the pointwise
partially hyperbolic case and answers to a conjecture from
\cite{HHU} where it is shown that the second possibility of the
theorem is non-empty. In the introduction of Section
\ref{SectionCoherenciaPHfuerte} we explain the difference between
our approach and the one of \cite{BBI2}.

In Section \ref{Section-AlmostDynamicalCoherence} we present the
definition of almost dynamical coherence as well as some
properties and we give some preliminaries of results which we will
use afterwards.

Finally, in Section \ref{Section-DimensionesMayores} we comment on
some results in higher dimensions as well as to explore some
results which allow one to characterize the isotopy class of a
partially hyperbolic diffeomorphism.


\section{Almost dynamical coherence and Quasi-Isometry}\label{Section-AlmostDynamicalCoherence}

\subsection{Almost dynamical coherence}\label{SubSectionADCPropiedades}
In general, a partially hyperbolic diffeomorphism may not be
dynamically coherent, and even if it is, it is not known in all
generality if being dynamically coherent is an open property (see
\cite{HPS,Berger}). However, all the known examples in dimension 3
verify the following property which is clearly $C^1$-open:

\begin{defi}[Almost dynamical coherence]\label{DefADC}
We say that $f: M \to M$ partially hyperbolic of the form $TM=
E^{cs} \oplus E^u$ is \emph{almost dynamically coherent} if there
exists a foliation $\cF$ transverse to the direction $E^u$.
\finobs
\end{defi}

The introduction of this definition is motivated by the work of
\cite{BI} where it was remarked that sometimes it is enough to
have a foliation transverse to the unstable direction in order to
obtain conclusions.

Almost dynamical coherence is not a very strong requirement, with
the basic facts on domination we can show:

\begin{prop}\label{PropositionADCCerrado}
Let $\{f_n\}$ a sequence of almost dynamically coherent partially
hyperbolic diffeomorphisms converging in the $C^1$-topology to a
partially hyperbolic diffeomorphism $f$. Then, $f$ is almost
dynamically coherent.
\end{prop}

\dem Let us call $E^{cs}_n \oplus E^u_n$ to the splitting of $f_n$
and $E^{cs}\oplus E^u$ to the splitting of $f$. We use the
following well known facts on domination (see Proposition
\ref{Proposicion-DominacionPasaALaClausuraYLimite} and Remark
\ref{Remark-ExtensionClausuraAngulosYContinuidad}):

\begin{itemize}
\item[-] The subspaces $E^{cs}_n$ and $E^u_n$ converge as $n\to
\infty$ towards $E^{cs}$ and $E^u$. \item[-] The angle between
$E^{cs}$ and $E^u$ is larger than $\alpha>0$.
\end{itemize}

Now, consider $f_n$ such that the angle between $E^{cs}_n$ and
$E^{u}$ is larger than $\alpha/2$. Let $\cF_n$ be the foliation
transverse to $E^u_n$.

By iterating backwards by $f_n$ we obtain that $f^{-m}_n (\cF_n)$
is, when $m$ is large, tangent to a small cone around $E^{cs}_n$.
From our assumptions, we can thus deduce that $f^{-m}_n (\cF_n)$
is transverse also to $E^u$. This implies that $f$ is almost
dynamically coherent as desired. \lqqd

Notice that this proposition implies that if we denote as
$\cP\cH^1(M)$ the set of partially hyperbolic diffeomorphisms of
$M$, and $\cP$ to a connected component: If $\cP$ contains an
almost dynamically coherent diffeomorphisms, every diffeomorphism
in $\cP$ is almost dynamically coherent. In particular, almost
dynamical coherent partially hyperbolic diffeomorphisms contain
the connected component in $\cP\cH^1(\TT^3)$ containing the linear
representatives of the isotopy class when these are partially
hyperbolic.

As a consequence \cite{BI} (Key Lemma 2.1), every strong partially
hyperbolic diffeomorphism of a $3$-dimensional manifold is almost
dynamically coherent. It is important to remark that it is a mayor
problem to determine whether partially hyperbolic diffeomorphisms
in the sphere $S^3$ are almost dynamically coherent (which would
solve the question on the existence of robustly transitive
diffeomorphisms in the sphere\footnote{Indeed, by a result of
\cite{DPU} a robustly transitive diffeomorphism of a
$3$-dimensional manifold should be partially hyperbolic. If it
were almost dynamically coherent, we would get by Novikov's
Theorem that there exists a Reeb component transverse to the
strong unstable direction (if it is partially hyperbolic of type
$TS^3 =E^s \oplus E^{cu}$ one should consider $f^{-1}$). This
contradicts Proposition \ref{Proposition-NoComponentesDeReeb}.}).

The author is not aware of whether the following question is known
or still open:

\begin{quest}\label{QuestionIsotopia}
Are there any examples of partially hyperbolic diffeomorphisms of
$\TT^3$ isotopic to a linear Anosov automorphism which are not
isotopic to the linear Anosov automorphism through a path of
partially hyperbolic diffeomorphisms?
\end{quest}


We end this subsection by stating a property first observed by
Brin,Burago and Ivanov (\cite{BBI,BI}) which makes our definition
a good tool for studying partially hyperbolic diffeomorphisms:

\begin{prop}[Brin-Burago-Ivanov]\label{Proposition-NoComponentesDeReeb}
Let $f: M \to M$ an almost dynamically coherent partially
hyperbolic diffeomorphism with splitting $TM = E^{cs} \oplus E^u$
and let $\cF$ be the foliation transverse to $E^u$. Then, $\cF$
has no Reeb components.
\end{prop}

\dem If a (transversally oriented) foliation $\cF$ on a compact
closed $3$-dimensional manifold $M$ has a Reeb component, then,
every one dimensional foliation transverse to $\cF$ has a closed
leaf (see \cite{BI} Lemma 2.2).

Since $\cF^u$ is one dimensional, transverse to $\cF$ and has no
closed leafs, we obtain that $\cF$ cannot have Reeb components.

\lqqd

This has allowed them to prove (see also \cite{Parwani}):

\begin{teo}[Brin-Burago-Ivanov \cite{BBI,BI,Parwani}]\label{TeoremaLaAccionEnHomologiaEsPH}
If $f$ is an almost dynamically coherent partially hyperbolic
diffeomorphism of a $3$-dimensional manifold $M$ with fundamental
group of polynomial growth, then, the induced map $f_\ast:
H_1(M,\RR) \cong \RR^k \to H_1(M, \RR)$ is partially hyperbolic.
This means, it is represented by an invertible matrix $A\in
GL(k,\ZZ)$ which has an eigenvalue of modulus larger than $1$ and
determinant of modulus $1$ (in particular, it also has an
eigenvalue of modulus smaller than $1$).
\end{teo}

\esbozo{} We prove the result when $\pi_1(M)$ is abelian, so that
it coincides with $H_1(M,\ZZ)$. When the fundamental group is
nilpotent, this follows from the fact that the $3$-manifolds with
this fundamental group are well known (they are circle bundles
over the torus) so that one can make other kind of arguments with
the same spirit (see \cite{Parwani} Theorem 1.12).

Assume that every eigenvalue of $f_\ast$ is smaller or equal to
$1$. Since the universal cover $\tilde M$ is quasi-isometric to
$\pi_1(M)$, it is thus quasi-isometric to $H_1(M,\RR) \cong \RR^k$
(notice that this is trivial if $M=\TT^3$).

Now, we have that $f_\ast$ acting in $H_1(M,\RR)$ has all of its
eigenvalues smaller than one, we obtain that the diameter of a
compact set in $\tilde M$ grows subexponentially by iterating it
with $\tilde f$.

Given $R>0$ the number of fundamental domains needed to cover a
ball of radius $R$ in $\tilde M$ is polynomial in $R$.

Consider an unstable arc $I$. We obtain that $\tilde f^n(I)$ has
subexponential (in $n$) diameter but the length grows
exponentially (in $n$). By the previous observation, we obtain
that given $\eps$ we find points of $\tilde \cF^u$ which are not
in the same local unstable manifold but are at distance smaller
than $\eps$, this implies the existence of a Reeb component for
$\tilde \cF$ (Theorem \ref{TeoNovikov}) and contradicts
Proposition \ref{Proposition-NoComponentesDeReeb}.

\lqqd

Notice that if the growth of the fundamental group is exponential,
one can make partially hyperbolic diffeomorphisms which are
isotopic to the identity (for example, the time-one map of an
Anosov flow). This is because in such a manifold, a sequence $K_n$
of sets with exponentially (in $n$) many points but polynomial (in
$n$) diameter may not have accumulation points.

As a consequence of combining the Proposition
\ref{Proposition-NoComponentesDeReeb} with Novikov's Theorem
\ref{TeoNovikov} we obtain for $\TT^3$ the following consequence
(recall Corollary \ref{CorolarioConsecuenciasReeb}):

\begin{cor}\label{CorolarioConsecuenciasReeb2}
Let $f$ be a partially hyperbolic diffeomorphism of $\TT^3$ of the
form $T\TT^3=E^{cs}\oplus E^u$ ($\dim E^{cs}=2$) which is almost
dynamically coherent with foliation $\cF$. Assume that $\cF$ is
oriented and transversally oriented and let $\tilde \cF$ and
$\tilde \cF^u$ the lifts of the foliations $\cF$ and the unstable
foliation $\cF^u$ to $\RR^3$. Then:
\begin{itemize}
\item[(i)] For every $x\in \RR^3$ we have that $\tilde \cF(x) \cap
\tilde \cF^u(x) = \{x\}$. \item[(ii)] The leafs of $\tilde \cF$
are properly embedded complete surfaces in $\RR^3$. In fact there
exists $\delta>0$ such that every euclidean ball $U$ of radius
$\delta$ can be covered by a continuous coordinate chart such that
the intersection of every leaf $S$ of $\tilde \cF$ with $U$ is
either empty of represented as the graph of a function $h_S: \RR^2
\to \RR$ in those coordinates. \item[(iii)] Each closed leaf of
$\cF$ is a two dimensional torus. \item[(iv)] For every
$\delta>0$, there exists a constant $C_\delta$ such that if $J$ is
a segment of $\tilde \cF^u$ then $\Vol(B_\delta(J)) > C_\delta
\length(J)$.
\end{itemize}
\end{cor}

\dem The proof of (i) is the same as the one of Lemma 2.3 of
\cite{BI}, indeed, if there were two points of intersection, one
can construct a closed loop transverse to $\tilde \cF$ which
descends in $\TT^3$ to a nullhomotopic one. By Novikov's theorem
(Theorem \ref{TeoNovikov}), this implies the existence of a Reeb
component, a contradiction with Proposition
\ref{Proposition-NoComponentesDeReeb}.

Once (i) is proved, (ii) follows from the same argument as in
Lemma 3.2 in \cite{BBI2}. Notice that the fact that the leafs of
$\tilde \cF$ are properly embedded is trivial after (i), with some
work, one can prove the remaining part of (ii) (see also Lemma
\ref{LemaEstructuraProductoLocalUniforme} for a more general
statement).

Part (iii) follows from the fact that if $S$ is an oriented closed
surface in $\TT^3$ which is not a torus, then it is either a
sphere or its fundamental group cannot inject in $\TT^3$ (see
\cite{ClasificacionSuperficies} and notice that a group with
exponential growth cannot inject in $\ZZ^3$).

Since $\cF$ has no Reeb components, we obtain that if $S$ is a
closed leaf of $\cF$ then it must be a sphere or a torus. But $S$
cannot be a sphere since in that case, the Reeb's stability
theorem (Theorem \ref{Teorema-EstabilidadCompleta}) would imply
that all the leafs of $\cF$ are spheres and that the foliated
manifold is finitely covered by $S^2 \times S^1$ which is not the
case.

The proof of (iv) is as Lemma 3.3 of \cite{BBI2}. Since there
cannot be two points in the same leaf of $\tilde \cF^u$ which are
close but in different local unstable leaves, we can find
$\epsilon>0$ and $a>0$ such that in a curve of length $K$ of
$\tilde \cF^u$ there are at least $aK$ points whose balls of
radius $\epsilon$ are disjoint (and all have the same volume).

Now, consider $\delta>0$ and $\tilde \delta= \min\{ \delta,
\epsilon \}$. Let $\{ x_1, \ldots, x_l \}$ with $l > a \length
(J)$ be points such that their $\tilde \delta$-balls are disjoint.
We get that $U= \bigcup_{i=1}^l B_{\tilde \delta}(x_i) \en
B_\delta(J)$ and we have that $\Vol(U) > l \Vol(B_{\tilde
\delta}(x_i))$. We obtain that $C_\delta= \frac{4\pi}{3}a\delta^3$
works.

\lqqd

Notice that most of the previous result can be extended to
arbitrary $3$-dimensional manifolds. In fact, with a similar proof
(see also \cite{Parwani}) one proves that almost dynamically
coherent partially hyperbolic diffeomorphisms can only occur in
certain specific $3$-manifolds:

\begin{cor}\label{Corolario-ConsequenciasReebenMarbitraria}
Let $f$ be an almost dynamically coherent partially hyperbolic
diffeomorphism of a $3$-dimensional manifold $M$ with splitting
$TM = E^{cs} \oplus E^u$. Then:
\begin{itemize}
\item[-] The manifold $M$ is irreducible (i.e. $\pi_2(M)=\{0\}$).
\item[-] The covering space of $M$ is homeomorphic to $\RR^3$.
\item[-] The fundamental group of $M$ is infinite (and different
from $\ZZ$)
\end{itemize}
\end{cor}

\dem Let $\cF$ be the foliation transverse to $E^u$.

The first claim follows from the fact that having $\pi_2(M) \neq
\{0\}$ implies the existence of a Reeb component for $\cF$ by
Novikov's Theorem \ref{TeoNovikov}. The last claim follows by the
same reason. The fact that the fundamental group cannot be $\ZZ$
follows from Proposition \ref{Proposition-NoComponentesDeReeb}
since a manifold with $\ZZ$ as fundamental group (which is of
polynomial growth) has $\ZZ$ as first homology group and admits no
automorphisms with eigenvalues of modulus different from $1$.

To get the second statement, notice that since the fundamental
group of every leaf must inject in the fundamental group of $M$ we
have that every leaf of $\tilde \cF$ must be homeomorphic to
$\RR^2$ or $S^2$. By Reeb's stability theorem (Theorem
\ref{Teorema-EstabilidadCompleta}) leafs must be homeomorphic to
$\RR^2$.

By a result by Palmeira (see \cite{Palmeira}) we obtain that
$\tilde M$ must be homeomorphic to $\RR^3$.

\lqqd

\subsection{Branched Foliations and Burago-Ivanov's result}\label{SubSection-BuragoIvanov}

We follow \cite{BI} section 4.

We define a \emph{surface} in a $3$-manifold $M$ to be a
$C^1$-immersion $\imath: U \to M$ of a connected smooth
2-dimensional manifold (possibly with boundary). The surface is
said to be \emph{complete} if it is complete with the metric
induced in $U$ by the Riemannian metric of $M$ and the immersion
$\imath$. The surface is \emph{open} if it has no boundary.

Given a point $x$ in (the image of) a surface $\imath: U \to M$ we
have that there is a neighborhood $B$ of $x$ such that the
connected component $C$ containing $\imath^{-1}(x)$ of
$\imath^{-1}(B)$ verifies that $\imath(C)$ separates $B$. We say
that two surfaces $\imath_1:U_1 \to M, \imath_2:U_2\to M$
\emph{topologically cross} if there exists a point $x$ in (the
image of) $\imath_1$ and $\imath_2$ and a curve $\gamma$ in $U_2$
such that $\imath_2(\gamma)$ passes through $x$ and intersects
both connected components of a neighborhood of $x$ with the part
of the surface defined above removed. This definition is symmetric
and does not depend on the choice of $B$ (see \cite{BI}) however
we will not use this facts.

\begin{defi}\label{Definicion-BranchingFoliation}
A \emph{branching foliation} on $M$ is a collection of complete
open surfaces tangent to a given continuous 2-dimensional
distribution on $M$ such that every point belongs to at least one
surface and no pair of surfaces of the collection have topological
crossings. \finobs
\end{defi}

We will abuse notation and denote a branching foliation as
$\cF_{bran}$ and by $\cF_{bran}(x)$ to the set of set of surfaces
whose image contains $x$. We call the (image of) the surfaces,
leaves of the branching foliation.

We have the following:

\begin{prop}\label{PropositionBranchingSinBranchingEsFoliacion}
If every point of $M$ belongs to a unique leaf of the branching
foliation, then the branching foliation is a true foliation.
\end{prop}

\dem Let $E$ be the two-dimensional distribution tangent to the
branching foliation and we consider $E^\perp$ a transverse
direction which we can assume is $C^1$ and almost orthogonal to
$E$.

By uniform continuity we find $\eps$ such that for every point $p$
in $M$ the $2\eps$ ball verifies that it admits a $C^1$-chart to
an open set in $\RR^3$ which sends $E$ to an almost horizontal
$xy$-plane and $E^\perp$ to an almost vertical $z$-line.

Let $D$ be a small disk in the (unique) surface through $p$ and
$\gamma$ a small arc tangent to $E^\perp$ thorough $p$. Given a
point $q \in D$ and $t\in \gamma$ we have that inside
$B_{2\eps}(p)$ there is a unique point of intersection between the
curve tangent to $E^{\perp}$ through $q$ and the connected
component of the (unique) surface of $\cF_{bran}$ intersected with
$B_{2\eps}(p)$ containing $t$.

We get a well defined continuous and injective map from $D \times
\gamma \cong \RR^3$ to a neighborhood of $p$ (by the invariance of
domain's theorem, see \cite{Hatcher}) such that it sends sets of
the form $D \times \{t \}$ into surfaces of the branching
foliation. Since we already know that $\cF_{bran}$ is tangent to a
continuous distribution, we get that $\cF_{bran}$ is a true
foliation.

\lqqd

Indeed, the result also follows from the following statement we
will also use:

\begin{prop}[\cite{BonWilk} Proposition 1.6 and Remark 1.10]\label{Proposition-BWProp16}
Let $E$ be a continuous codimension one distribution on a manifold
$M$ and $S$ a (possibly non connected) surface tangent to $E$
which contains a family of disks of fixed radius and whose set of
midpoints is dense in $M$. Then, there exists a foliation $\cF$
tangent to $E$ which contains $S$ in its leaves.
\end{prop}

Invariant branching foliations always exist for strong partially
hyperbolic diffeomorphisms of $3$-dimensional manifolds due to a
remarkable result of Burago and Ivanov:

\begin{teo}[\cite{BI},Theorem 4.1 and Theorem 7.2]\label{TeoBuragoIvanov}Let $f:M^3\to M^3$
be a strong partially hyperbolic diffeomorphism with splitting
$TM= E^s \oplus E^c \oplus E^u$ into one dimensional subbundles.
There exists branching foliations $\cF^{cs}_{bran}$ and
$\cF^{cu}_{bran}$ tangent to $E^{cs}=E^s \oplus E^c $ and
$E^{cu}=E^c \oplus E^u$ which are $f$-invariant\footnote{This
means that for every $\cF_k \in \cF^{\sigma}_{bran}(x)$ there
exists $\cF_{k'} \in \cF^\sigma_{bran}(f(x))$ such that $f(\cF_k)=
\cF_{k'}$.}. Moreover, for every $\eps>0$ there exist foliations
$\cS_\eps$ and $\cU_\eps$ tangent to an $\eps$-cone around
$E^{cs}$ and $E^{cu}$ respectively and continuous maps
$h^{cs}_\eps$ and $h^{cu}_\eps$ at $C^0$-distance smaller than
$\eps$ from the identity which send the leaves of $\cS_\eps$ and
$\cU_\eps$ to leaves of $\cF^{cs}_{bran}$ and $\cF^{cu}_{bran}$
respectively.
\end{teo}

We remark that when there exists an $f$-invariant (branching)
foliation, one can assume that every sequence of leaves through
points $x_k$ such that $x_k \to x$ verifies that it converges to a
leaf through $x$ (see Lemma 7.1 of \cite{BI}).

\begin{convention}
We will assume throughout that every branching foliation is
\emph{completed} in the sense stated above: For every sequence
$L_k$ of leaves in $\cF^{cs}_{bran}(x_k)$ such that $x_k \to x$ we
have that $L_k$ converges in the $C^1$-topology to a leaf $L \in
\cF^{cs}_{bran}(x)$ contained in the branching foliation. \finobs
\end{convention}

Notice that the existence of the maps $h^{cs}_\eps$ and
$h^{cu}_\eps$ implies that when lifted to the universal cover, the
leaves of $\cS_\eps$ (resp. $\cU_\eps$) remain at distance smaller
than $\eps$ from lifted leaves of $\cF^{cs}_{bran}$ (resp.
$\cF^{cu}_{bran}$).

We obtain as a corollary the following result we have already
announced:

\begin{cor}[Key Lemma 2.2 of \cite{BI}]\label{Corolario-PHFuerteESADC}
A strong partially hyperbolic diffeomorphism on a $3$-dimensional
manifold is almost dynamically coherent.
\end{cor}

Using the fact that when $x_k \to x$ the leaves through $x_k$
converge to a leaf through $x$ we obtain:

\begin{prop}\label{Proposition-SucesionDeTorosEnBranchingFol}
Let $\cF_{bran}$ be a branching foliation of $\TT^3$ and consider
a sequence of points $x_k$ such that there are leaves $\cF_k \in
\cF_{bran}(x_k)$ which are compact, incompressible and homotopic
to each other. If $x_k \to x$, then there is a leaf $L \in
\cF_{bran}(x)$ which is incompressible and homotopic to the leaves
$\cF_k$.
\end{prop}

\dem Recall that if $x_k \to x$ and we consider a sequence of
leaves through $x_k$ we get that the leaves converge to a leaf
through $x$.

Consider the lifts of the leaves $\cF_k$ which are homeomorphic to
a plane since they are incompressible. Moreover, the fundamental
group of each of the leaves must be $\ZZ^2$ and the leaves must be
homoeomorphic to 2-torus, since it is the only possibly
incompressible surface in $\TT^3$.

Since all the leaves $\cF_k$ are homotopic, their lifts are
invariant under the same elements of $\pi_1(\TT^3)$. The limit
leaf must thus be also invariant under those elements. Notice that
it cannot be invariant under further elements of $\pi_1(\TT^3)$
since no surface has such fundamental group.

\lqqd

The idea of the proof of the previous proposition can be applied
to other contexts, however, for simplifying the proof we chose to
state it only in this context which is the one of interest for us.


\subsection{Quasi-isometry and dynamical
coherence}\label{Section-QIyArgumentoBrin}

We review in this section a simple criterium given by Brin in
\cite{Brin} which guaranties dynamical coherence for absolutely
dominated partially hyperbolic diffeomorphisms. It involves the
concept of quasi-isometry which we will use after in this thesis.
We present the sketch of the proof by Brin to show the importance
of absolute domination in his argument.

For more information on quasi-isometric foliations we refer the
reader to \cite{HQuasiIsometry}. We recall its definition (which
already appeared in subsection \ref{SubSection-ConsecuenciasGPS}):

\begin{defi}[Quasi-Isometric Foliation]\label{Definicion-QuasiIsom}
Consider a Riemannian manifold $M$ (not necessarily compact) and a
foliation $\cF$ in $M$. We say that the foliation $\cF$ is
\emph{quasi-isometric} if distances inside leaves can be compared
with distances in the manifold. More precisely, for $x,y \in
\cF(x)$ we denote as $d_{\cF}(x,y)$ as the infimum of the lengths
of curves joining $x$ to $y$, we say that $\cF$ is quasi-isometric
if there exists $a,b\in \RR$ such that for every $x,y$ in the same
leaf of $\cF$ one has that:

$$  d_{\cF}(x,y) \leq a d(x,y) + b $$
\finobs
\end{defi}

In general, this notion makes sense in non-compact manifolds, and
it will be used by us mainly in the universal covering space of
the manifolds we work with. Notice that if a foliation of a
compact manifold is quasi-isometric then all leaves must be
compact.

The classic example of a quasi-isometric foliation is a linear
foliation in $\RR^d$ with the euclidean metric. Indeed, it can be
thought that quasi-isometry foliations are in a sense a
generalization of these (notice however that even a one
dimensional foliation of the plane\footnote{Consider for example
the foliation given by $\{(t, t^3 + b)\}_{b}$.} which is
quasi-isometric needs not remain at bounded distance from a
one-dimensional ``direction'').

It is important to remark that the metric in the manifold is quite
important, and as in general we work with the universal cover of a
compact manifold, this metric is also influenced by the topology
of the manifold. See \cite{HQuasiIsometry} for more discussion on
quasi-isometric foliations and topological and restrictions for
their existence.

The argument of the proof of Proposition
\ref{Proposicion-CoherenciaLocal} can be extended to non-local
arguments if one demands that the domination required is absolute
and the geometry of leaves is quite special. In fact, Brin has
proved in \cite{Brin} the following quite useful criterium (see
for example \cite{BBI2}, \cite{Parwani} or \cite{Hammerlindl,HNil}
for applications of this criterium).

\begin{prop}[\cite{Brin}]\label{Proposition-BrinArgument}
Let $f: M \to M$ be an absolutely partially hyperbolic
diffeomorphism with splitting $TM = E^{cs} \oplus E^u$ and such
that the foliation $\tilde \cF^u$ is quasi-isometric in $\tilde M$
the universal cover of $M$. Then, $f$ is dynamically coherent.
\end{prop}

It shows that in fact, the foliation is unique in (almost) the
strongest sense, which is that every $C^1$-embedding of a ball of
dimension $\dim E^{cs}$ which is everywhere tangent to $E^{cs}$ is
in fact contained in a leaf of the foliation $\cF^{cs}$.

We give a sketch of the proof in order to show how the hypothesis
are essential to pursue the argument. See \cite{Brin} for a clear
exposition of the complete argument.

\esbozo{} Assume that there are two embedded balls $B_1$ and $B_2$
through a point $x$ which are everywhere tangent to $E^{cs}$ and
whose intersection is not relatively open in (at least) one of
them.

Then, as in Proposition \ref{Proposicion-CoherenciaLocal} it is
possible to construct a curve $\eta$ which has non-zero length, is
contained in a leaf of $\cF^u$ and joins these two embedded balls.

Let $\gamma_1$ and $\gamma_2$ two curves contained in $B_1$ and
$B_2$r respectively joining $x$ to the extremes of $\eta$.

Since $\eta$ is an unstable curve, its length growths
exponentially, and by quasi-isometry, we know that the extremal
points of the curve are at a distance which grows exponentially
with the same rate as the rate the vectors in $E^u$ expand.

On the other hand, the curves $\gamma_1$ and $\gamma_2$ are forced
to grow with at most an exponential rate which is smaller than the
one in $E^u$ (by using absolute domination) and so we violate the
triangle inequality.

\lqqd

\section{Partially hyperbolic diffeomorphisms isotopic to linear Anosov automorphisms of $\TT^3$}\label{Section-PHAnosov}

In this section we give a proof of the following:

\begin{teo}\label{TEOREMA-COHERENCIAISOTOPICOANOSOV}
Let $f: \TT^3 \to \TT^3$ be an almost dynamically coherent
partially hyperbolic diffeomorphism with splitting of the form
$T\TT^3 = E^{cs} \oplus E^u$. Assume that $f$ is isotopic to
Anosov, then:
\begin{itemize}
\item[-] $f$ is (robustly) dynamically coherent and has a unique
$f$-invariant foliation $\cF^{cs}$ tangent to $E^{cs}$. \item[-]
There exists a global product structure between the lift of
$\cF^{cs}$ to the universal cover and the lift of $\cF^u$ to the
universal cover. \item[-] If $f_\ast$ has two eigenvalues of
modulus larger than $1$ then they must be real and different.
\end{itemize}
\end{teo}

As a consequence of the fact that almost dynamical coherence is an
open and closed property (see Proposition
\ref{PropositionADCCerrado} above) we obtain:

\begin{cor*}
Dynamical coherence is an open and closed property among partially
hyperbolic diffeomorphisms of $\TT^3$ isotopic to Anosov.
\end{cor*}

\dem By Proposition \ref{PropositionADCCerrado} we know that
almost dynamical coherence is an open and closed property. Theorem
\ref{TEOREMA-COHERENCIAISOTOPICOANOSOV} then implies that in the
isotopy class of Anosov dynamical coherence is open and closed too
(since almost dynamical coherence implies dynamical coherence in
this context).

\lqqd

We shall assume that $f: \TT^3 \to \TT^3$ is an almost dynamical
coherent partially hyperbolic diffeomorphism  with splitting of
the form $T \TT^3 = E^{cs} \oplus E^u$ with $\dim E^u=1$ and
isotopic to a linear Anosov automorphism $A: \TT^3 \to \TT^3$.

It is important to remark that we are not assuming that the stable
dimension of $A=f_\ast$ coincides with the one of $E^{cs}$. In
fact, many of the arguments below become much easier in the case
$A$ has stable dimension $2$. The fact that we can treat the case
where $A$  has two eigenvalues of modulus larger than one is in
the authors' opinion, one of the main contributions of this
thesis.

We will denote as $\cF$ the foliation given by the definition of
almost dynamical coherence which we know is Reebless and it thus
verifies the hypothesis of Theorem
\ref{PropObienHayTorosObienTodoEsLindo}.

As before, we denote as $p: \RR^3 \to \TT^3$ the covering
projection and we denote as $\tilde f$, $\tilde \cF$ and $\tilde
\cF^u$ the lifts of $f$, $\cF$ and $\cF^u$ to the universal cover.

We provide an orientation to $\tilde \cF^u$ and denote as $\tilde
\cF^u_+(x)$ and $\tilde \cF^u_-(x)$ to the connected components of
$\tilde \cF^u(x) \setminus \{x\}$. Since $\tilde \cF(x)$ separates
$\RR^3$ (see subsection \ref{SubSection-FoliacionesDeT3}) we
denote $F_+(x)$ and $F_-(x)$ to the connected components of $\RR^3
\setminus \tilde \cF(x)$ containing respectively $\tilde
\cF^u_+(x)$ and $\tilde \cF^u_-(x)$.

Proposition \ref{PropExisteSemiconjugacion} implies the existence
of a continuous and surjective function $H: \RR^3\to \RR^3$ which
verifies
$$H \circ \tilde f = A \circ H$$

\noindent and such that $d(H(x),x)< K_1$ for every $x\in \RR^3$.

\subsection{Consequences of the semiconjugacy}\label{SectionInestableNoAcotada}

We can prove:

\begin{lema}\label{LemaHdemediaFuesNoacotado} For every $x\in \RR^3$ we have that $H(\tilde \cF^u_+(x))$ is unbounded.
\end{lema}
\dem Otherwise, for some $x\in \RR^3$, the unstable leaf $\tilde
\cF^u_+(x)$ would be bounded. Since its length is infinite one can
find two points in $\tilde \cF^u_+(x)$ in different local unstable
leafs at arbitrarily small distance. This contradicts Corollary
\ref{CorolarioConsecuenciasReeb2} (i). \lqqd

\begin{obs}\label{RemarkHAMastieneInestables}
Notice that for every $x\in \RR^3$ the set $F_+(x)$ is unbounded
and contains a half unstable leaf of $\tilde \cF^u$.
\begin{itemize}
\item[-] In the case the automorphism $A$ has stable dimension 2,
this implies that $H(F_+(x))$ contains a half-line of irrational
slope. Indeed, by Lemma \ref{LemaHdemediaFuesNoacotado} we have
that $H(\tilde \cF^u_+(x))$ is non bounded and since we know that
$H(\tilde \cF^u(x))\en W^u(H(x), A)$ we conclude. \item[-] When
$A$ has stable dimension 1, we only obtain that $H(\tilde
\cF^u_+(x))$ contains an unbounded connected set in $W^u(H(x), A)$
which is two dimensional plane parallel to $E^u_A$.
\end{itemize}
\finobs
\end{obs}

 One can push forward Lemma \ref{LemaHdemediaFuesNoacotado} in
order to show that $H$ is almost injective in each unstable leaf
of $\tilde \cF^u$, in particular, a similar argument to the one in
Lemma \ref{LemaHdemediaFuesNoacotado} gives that at most finitely
many points of an unstable leaf can have the same image under $H$.
Later, we shall obtain that in fact, $H$ is injective on unstable
leaves so that we will not give the details of the previous claim
(see Remark \ref{RemarkHesInjectivaEnInestables}).

\subsection{A planar direction for the foliation transverse to $E^u$}

Since $\cF$ is transverse to the unstable direction, we get by
Corollary \ref{CorolarioConsecuenciasReeb} that it is a Reebless
foliation so that we can apply Theorem
\ref{PropObienHayTorosObienTodoEsLindo}. We intend to prove in
this section that option (ii) of Theorem
\ref{PropObienHayTorosObienTodoEsLindo} is not possible when $f$
is isotopic to Anosov (see \cite{HHU} where that possibility
occurs). The following simple remark will be essential in what
follows:

\begin{obs}\label{RemarkSiAplicofElPlanoSeleAplicaA}
Notice that if we apply $\tilde f^{-1}$ to the foliation $\tilde
\cF$, then the new foliation $\tilde f^{-1}(\tilde \cF)$ is still
transverse to $E^u$ so that Theorem
\ref{PropObienHayTorosObienTodoEsLindo} still applies. So, we
obtain a plane $P'$ close to $\tilde f^{-1}(\tilde \cF)$. We claim
that $P' = A^{-1}(P)$ where $P$ is the plane given by Theorem
\ref{PropObienHayTorosObienTodoEsLindo} for $\tilde \cF$. To prove
this, recall that $\tilde f$ and $A$ are at bounded distance so,
leaves of $\tilde f^{-1}(\cF)$ must remain at bounded distance
from $A^{-1}(P)$ and then use the fact that the plane is unique
(Remark \ref{Remark-UnicidaddelPlanoP}). \finobs
\end{obs}

The result that follows can be deduced more easily if one assumes
that $A$ has stable dimension 2.

We say that a subspace $P$ is \emph{almost parallel} to a
foliation $\tilde \cF$ if there exists $R>0$ such that for every
$x\in \RR^3$ we have that $P+x$ lies in an $R$-neighborhood of
$\tilde \cF(x)$ and $\tilde \cF(x)$ lies in a $R$-neighborhood of
$P+x$.

\begin{prop}\label{PropSifIsotopicoAAnosovNoHayTorosEnFoliacion}
Let $f: \TT^3 \to \TT^3$ be a partially hyperbolic diffeomorphism
of the form $T\TT^3 = E^{cs} \oplus E^u$ (with $\dim E^{cs}=2$)
isotopic to a linear Anosov automorphism and $\cF$ a foliation
transverse to $E^u$. Then, there exists a two dimensional subspace
$P \en \RR^3$ which is almost parallel to $\tilde \cF$.
\end{prop}

\dem It is enough to show that option (i) of Proposition
\ref{PropObienHayTorosObienTodoEsLindo} holds, since it implies
the existence of a plane $P$ almost parallel to $\tilde \cF$.

Assume by contradiction that option (ii) of Proposition
\ref{PropObienHayTorosObienTodoEsLindo} holds. Then, there exists
a plane $P \en \RR^3$ whose projection to $\TT^3$ is a two
dimensional torus and such that every leaf of $\tilde \cF^u$,
being transverse to $\tilde \cF$, remains at bounded
distance\footnote{Notice that if $A$ has stable dimension 2, this
already gives us a contradiction since $H(\tilde \cF^u(x)) =
W^u(H(x),A)$ which is totally irrational and cannot acumulate in a
plane which projects into a two-torus.} from $P$. Indeed, when
there is a dead-end component for $\cF$ we get that any transverse
foliation must verify that its leaves remain at bounded distance
from the boundary torus of the dead-end component which in turn
are at bounded distance from the plane $P$.

Since $f$ is isotopic to a linear Anosov automorphism $A$ we know
that $P$ cannot be invariant under $A$ (see Proposition
\ref{PropAnosovSonIrreducibles}). So, we have that $P$ and
$A^{-1}(P)$ intersect in a one dimensional subspace $L$ which
projects into a circle in $\TT^3$ (notice that a linear curve in
$\TT^2$ is either dense or a circle, so, if a line belongs to the
intersection of two linear two dimensional torus in $\TT^3$ which
do not coincide, it must be a circle).


We get that for every point $x$ we have that $\tilde \cF^u(x)$
must lie within bounded distance from $P$ as well as from
$A^{-1}(P)$ (since when we apply $\tilde f^{-1}$ to $\tilde \cF$
the leaf close to $P$ becomes close to $A(P)$, see Remark
\ref{RemarkSiAplicofElPlanoSeleAplicaA}). This implies that in
fact $\tilde \cF^u(x)$ lies within bounded distance from $L$.


On the other hand, we have that $H(\tilde \cF^u_+(x))$ is
contained in $W^u(H(x), A)= E^u_A + H(x)$ for every $x\in \RR^3$.
Since $H$ is at bounded distance from the identity, we get that
$\tilde \cF^u(x)$ lies within bounded distance from $E^u_A$, the
eigenspace corresponding to the unstable eigenvalues of $A$.

Since $E^u_A$ must be totally irrational (see Remark
\ref{RemarkSubespaciosInvariantesAnosov}) and $L$ projects into a
circle $L$, we get that $\cF^u_+(x)$ remains at bounded distance
from $E^u_A \cap L = \{0\}$. This contradicts the fact that
$\tilde \cF^u_+(x)$ is unbounded (Lemma
\ref{LemaHdemediaFuesNoacotado}). \lqqd



\subsection{Global Product Structure in the universal cover}\label{SectionGPS}

When the plane $P$ almost parallel to $\tilde \cF$ is totally
irrational, one can see that the foliation $\cF$ in $\TT^3$ is
without holonomy, and thus there is a global product structure
between $\tilde \cF$ and $\tilde \cF^u$ which follows directly
from Theorem \ref{Teorema-HectorHirsch}.

This would be the case if we knew that the plane $P$ given by
Theorem \ref{PropObienHayTorosObienTodoEsLindo} is
$f_\ast$-invariant (see subsection
\ref{subsection-ApendicedelPaperPHFol}). To obtain the global
product structure in our case we will use the fact that iterating
the plane $P$ backwards by $f_\ast$ it will converge to an
irrational plane and use instead Theorem
\ref{Teorema-EstructuraProductoGlobalMIO}.

Proposition \ref{PropSifIsotopicoAAnosovNoHayTorosEnFoliacion}
implies that the foliation $\tilde \cF$ is quite well behaved. In
this section we shall show that the properties we have showed for
the foliations and the fact that $\tilde \cF^u$ is $\tilde
f$-invariant while the foliation $\tilde \cF$ remains with a
uniform local product structure with $\tilde \cF^u$ when iterated
backwards (see Lemma \ref{LemaEstructuraProductoLocalUniforme})
imply that there is a global product structure. Some of the
arguments become simpler if one assumes that $A$ has stable
dimension $2$.

The main result of this section is thus the following:

\begin{prop}\label{PropositionGlobalProductStructure}
Given $x, y\in \RR^3$ we have that $\tilde \cF(x) \cap \tilde
\cF^u(y) \neq \emptyset$. This intersection consists of exactly
one point.
\end{prop}

Notice that uniqueness of the intersection point follows from
Corollary \ref{CorolarioConsecuenciasReeb2} (i) and will be used
to prove the proposition. We must put ourselves in the conditions
of Theorem \ref{Teorema-EstructuraProductoGlobalMIO}.

We shall proceed with the proof of Proposition
\ref{PropositionGlobalProductStructure}.

We start by proving a result which gives that the size of local
product structure boxes between $f^{-n}(\cF)$ and $\cF^u$ can be
chosen independent of $n$. We shall denote as $\DD^2= \{ z\in \CC
\ : \ |z|\leq 1 \}$.

\begin{lema}\label{LemaEstructuraProductoLocalUniforme}
There exists $\delta>0$ such that for every $x\in \RR^3$ and
$n\geq 0$ there exists a closed neighborhood $V_x^n$ containing
$B_\delta(x)$ such that it admits $C^0$-coordinates $\varphi_x^n:
\DD^2 \times [-1,1] \to \RR^3$ such that:
\begin{itemize}
\item[-] $\varphi_x^n(\DD^2\times [-1,1]) = V_x^n$ and
$\varphi_x^n(0,0)= x$. \item[-] $\varphi_x^n(\DD^2 \times \{t\}) =
\tilde f^{-n}(\tilde \cF(\tilde f^n(\varphi_x^n(0,t)))) \cap
V_x^n$ for every $t \in [-1,1].$ \item[-] $ \varphi_x^n(\{s\}
\times [-1,1]) = \tilde \cF^u(\varphi_x^n(s,0)) \cap V_x^n$ for
every $s\in \DD^2.$
\end{itemize}
\end{lema}


\dem  Notice first that the  tangent space to $f^{-n}(\cF)$
belongs to a cone transverse to $E^u$ and independent of $n$. Let
us call this cone $\cE^{cs}$.

Given $\epsilon>0$ we can choose a neighborhood $V_{\epsilon}$ of
$x$ contained in $B_{\epsilon}(x)$ such that the following is
verified:

\begin{itemize}
\item[-] There exists a two dimensional disk $D$ containing $x$
such that $V_{\epsilon}$ is the union of segments of $\cF^u(x)$ of
length $2\epsilon$ centered at points in $D$. This defines two
boundary disks $D^+$ and $D^-$ contained in the boundary of
$V_\eps$. \item[-] By choosing $D$ small enough, we get that there
exists $\epsilon_1>0$ such that every curve of length $\epsilon_1$
starting at a point $y \in B_{\epsilon_1}(x)$ tangent to
$\cE^{cs}$ must leave $V_{\epsilon}$ and intersects $\partial
V_{\epsilon}$ in $\partial V_{\epsilon} \setminus (D^+ \cup D^-)$.
\end{itemize}

Notice that both $\epsilon$ and $\epsilon_1$ can be chosen
uniformly in $\RR^3$ because of compactness of $\TT^3$ and uniform
transversality of the foliations (see Remark
\ref{RemarkLocalProductStructure}).

This implies that every disk of radius $\epsilon$ tangent to
$\cE^{cs}$ centered at a point $z \in B_{\epsilon_1}(x)$ must
intersect the unstable leaf of every point in $D$, in particular,
there is a local product structure of uniform size around each
point in $\RR^3$.

Now, we can choose a continuous chart (recall that the foliations
are with $C^1$ leaves but only continuous) around each point which
sends horizontal disks into disks transverse to $E^u$ and vertical
lines into leaves of $\tilde \cF^u$ containing a fixed ball around
each point $x$ independent of $n\geq 0$ giving the desired
statement. \lqqd

\begin{obs}\label{RemarkEPLUniforme}
We obtain that there exists $\eps>0$ such that for every $x\in
\RR^3$ there exists $V_x \en \bigcap_{n\geq 0} V_x^n$ containing
$B_\eps(x)$ admitting $C^1$-coordinates $\psi_x : \DD^2 \times
[-1,1] \to \RR^3$ such that:
\begin{itemize}
\item[-] $\psi_x (\DD^2 \times [-1,1]) = V_x$ and $\psi_x(0,0)=x$.
\item[-] If we consider $V_x^\eps= \psi_x^{-1}(B_\eps(x))$ then
one has that for every $y \in B_\eps(x)$ and $n\geq 0$ we have
that:
$$ \psi_x^{-1}(\tilde f^{-n}(\tilde \cF(\tilde f^n(y))) \cap V_x)
$$
    \noindent is the graph of a function $h_y^n: \DD^2 \to [-1,1]$ which has uniformly bounded derivative in $y$ and $n$.
\end{itemize}
Indeed, this is given by considering a $C^1$-chart $\psi_x$ around
every point such that its image covers the $\eps$-neighborhood of
$x$ and sends the $E$-direction to an almost horizontal direction
and the $E^u$-direction to an almost vertical direction (see
Proposition \ref{PropositionBranchingSinBranchingEsFoliacion}).
See for example \cite{BuW} section 3 for more details on this kind
of constructions. \finobs
\end{obs}

This lemma shows that after iterating the foliation backwards, one
gets that it becomes nearly irrational so that we can apply
Theorem \ref{Teorema-EstructuraProductoGlobalMIO}.

\begin{lema}\label{LemaIterandoParaAtrasSeHaceIrracional}
Given $K>0$ there exists $n_0>0$ such that for every $x\in \RR^3$
and for every $\gamma \in \ZZ^3$ with norm less than $K$ we have
that $$\tilde f^{-n_0}(\tilde \cF(x))+ \gamma \neq \tilde
f^{-n_0}(\tilde \cF(x)) \qquad  \forall x \in \RR^3.$$
\end{lema}

\dem Notice that $\tilde f^{-n}(\tilde \cF)$ is almost parallel to
$A^{-n}(P)$. Notice that $A^{-n}(P)$ has a converging subsequence
towards a totally irrational plane $\tilde P$ (see Remarks
\ref{RemarkSubespaciosInvariantesAnosov} and
\ref{RemarkValoresPropiosAnosov}).

We can choose $n_0$ large enough such that no element of $\ZZ^3$
of norm smaller than $K$ fixes $A^{-n_0}(P)$.

Notice first that $\tilde f^{-n_0}(\tilde \cF)$ is almost parallel
to $A^{-n_0}(P)$  (see Remark
\ref{RemarkSiAplicofElPlanoSeleAplicaA}). Now, assuming that there
is a translation $\gamma$ which fixes a leaf of $\tilde
f^{-n_0}(\tilde \cF)$ we get that the leaf contains a loop
homotopic to $\gamma$. This implies that it is at bounded distance
from the line which is the lift of the canonical (linear)
representative of $\gamma$ (see Lemma \ref{RemarkHojasCerradas}).
This implies that $\gamma$ fixes $A^{-n_0}(P)$ and thus has norm
larger than $K$ as desired.

\lqqd

We can now complete the proof of Proposition
\ref{PropositionGlobalProductStructure}.

\demo{ of Proposition \ref{PropositionGlobalProductStructure}} By
Corollary \ref{CorolarioConsecuenciasReeb2} we know that all the
leaves of $\tilde \cF$ are simply connected. Proposition
\ref{RemarkElCocienteDeLaFoliacionEsR} implies that the leaf space
of $\tilde \cF$ is homeomorphic to $\RR$. All this properties
remain true for the foliations $\tilde f^{-n}(\tilde \cF)$ since
they are diffeomorphisms at bounded distance from linear
transformations.

Lemma \ref{LemaEstructuraProductoLocalUniforme} gives that the
size of the local product structure between $\tilde f^{-n}(\tilde
\cF)$ and $\tilde \cF^u$ does not depend on $n$.

Using Lemma \ref{LemaIterandoParaAtrasSeHaceIrracional} we get
that for some sufficiently large $n$ the foliations $\tilde
f^{-n}(\tilde \cF)$ and $\tilde \cF^u$ are in the hypothesis of
Theorem \ref{Teorema-EstructuraProductoGlobalMIO} which gives
global product structure between $\tilde f^{-n}(\tilde \cF)$ and
$\tilde \cF^u$. Since $\tilde \cF^u$ is $\tilde f$-invariant and
$f$ is a diffeomorphism we get that there is a global product
structure between $\tilde \cF$ and $\tilde \cF^u$ as desired.

\lqqd

\begin{figure}[ht]
\begin{center}
\input{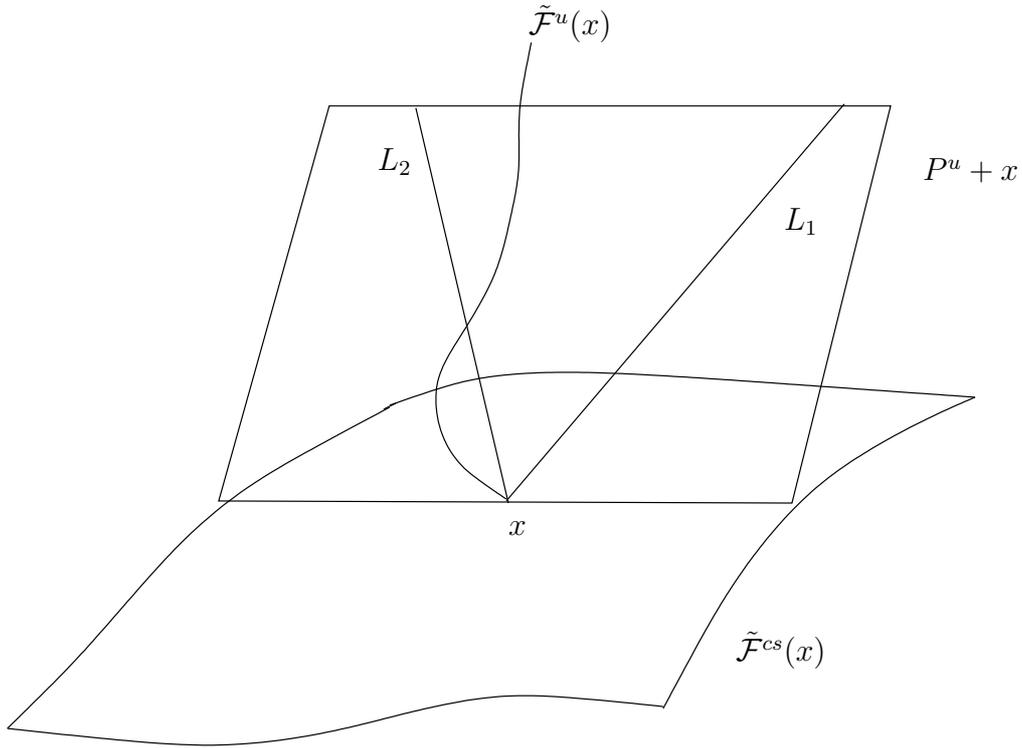}
\caption{\small{The unstable leaf of $x$ remains close to the cone
bounded by $L_1$ and $L_2$.}} \label{FiguraCono}
\end{center}
\end{figure}

Using Proposition \ref{PropositionQuasiIsometria} we deduce the
following (see figure \ref{FiguraCono}) :

\begin{cor}\label{PropositionQuasiIsometriaII}
The foliation $\tilde \cF^u$ is quasi-isometric. Moreover, there
exist one dimensional subspaces $L_1$ and $L_2$ of $E^u_A$
transverse to $P$ and $K>0$ such that for every $x\in \RR^3$ and
$y\in \tilde \cF^u(x)$ at distance larger than $K$ from $x$ we
have that $H(y)-H(y)$ is contained in the cone of $E^u_A$ with
boundaries $L_1$ and $L_2$ and transverse to $P$.
\end{cor}

Notice that if $A$ has stable dimension $2$ then $L_1=L_2 =E^u_A$.

\dem This is a direct consequence of Proposition
\ref{PropositionQuasiIsometria} and the fact that the image of
$\tilde \cF^u(x)$ by $H$ is contained in $E^u_A +H(x)$. \lqqd

\begin{obs}\label{RemarkHesInjectivaEnInestables}
Since points which are sent to the same point by $H$ must have
orbits remaining at bounded distance, the quasi-isometry of
$\tilde \cF^u$ implies that $H$ must be injective on leaves of
$\tilde \cF^u$. \finobs
\end{obs}

\subsection{Complex eigenvalues}\label{SubSectionComplejosYQIsom}

The following proposition has interest only in the case $A$ has
stable dimension 1.

\begin{prop}\label{PropositionAnotieneinestablesComplejos} The matrix $A$ cannot have complex unstable eigenvalues.
\end{prop}

\dem Assume that $A$ has complex unstable eigenvalues, in
particular $E^u_A$ is two-dimensional. Consider a fixed point
$x_0$ of $\tilde f$.

Recall that by Lemma \ref{LemaHdemediaFuesNoacotado} the set $\eta
=H(\tilde \cF^u_+(x_0))$ is an unbounded continuous curve in
$E^u_A$. Since $x_0$ is fixed and since $H$ is a semiconjugacy, we
have that $\eta$ is $A$-invariant.

On the other hand, by Corollary \ref{PropositionQuasiIsometriaII}
we have that $\eta$ is eventually contained in a cone between two
lines $L_1$ and $L_2$.

This implies that $A$ cannot have complex unstable eigenvalues
(recall that they should have irrational angle by Lemma
\ref{RemarkValoresPropiosAnosov}) since a matrix which preserves
an unbounded connected subset of a cone cannot have complex
eigenvalues with irrational angle. \lqqd

\subsection{Dynamical Coherence}\label{SubSectionDynamicalCoherence}

In this section we shall show dynamical coherence of almost
dynamically coherent partially hyperbolic diffeomorphisms isotopic
to linear Anosov automorphisms.

The proof of the following theorem becomes much simpler if one
assumes that the plane $P$ almost parallel to $\tilde \cF$ is
$A$-invariant which as we mentioned before is the most important
case (see also subsection \ref{subsection-ApendicedelPaperPHFol})

\begin{teo}\label{TeoremaCoherencia}
Let $f:\TT^3 \to \TT^3$ be an almost dynamically coherent
partially hyperbolic diffeomorphism of the form $T\TT^3 = E^{cs}
\oplus E^u$ isotopic to a linear Anosov automorphism. Then, there
exists an $f$-invariant foliation $\cF^{cs}$ tangent to $E^{cs}$.
If $\tilde \cF^{cs}$ denotes the lift to $\RR^3$ of this
foliation, then $H(\tilde \cF^{cs}(x)) = P^{cs} + H(x)$ where
$P^{cs}$ is an $A$-invariant subspace and $E^u_A$ is not contained
in $P^{cs}$.
\end{teo}

\dem Consider the foliation $\tilde \cF$, by Proposition
\ref{PropSifIsotopicoAAnosovNoHayTorosEnFoliacion} we have a plane
$P$ which is almost parallel to $\tilde \cF$.

Let $P^{cs}$ be the limit of $A^{-n}(P)$ which is an $A$-invariant
subspace. Since we have proved that $A$ has no complex unstable
eigenvalues (Proposition
\ref{PropositionAnotieneinestablesComplejos}) and since $P$ is
transverse to $E^u_A$ (Proposition
\ref{PropositionQuasiIsometria}), this plane is well defined (see
Remark \ref{RemarkSubespaciosInvariantesAnosov}).

Notice that the transversality of $P$ with $E^u_A$ implies that
$P^{cs}$ contains $E^s_A$, the eigenspace associated with stable
eigenvalues (in the case where $A$ has stable dimension $2$ we
thus have $P^{cs}=E^s_A$).

Since $P^{cs}$ is $A$-invariant, we get that it is totally
irrational so that no deck transformation fixes $P^{cs}$.

Using Remark \ref{RemarkEPLUniforme} we obtain $\eps>0$ such that
for every $x\in \RR^3$ there are neighborhoods $V_x$ containing
$B_\eps(x)$ admitting $C^1$-coordinates $\psi_x : \DD^2 \times
[-1,1] \to V_x$ such that:

\begin{itemize}
\item[-] For every $y\in B_\eps(x)$ we have that if we denote as
$W^x_n(y)$ to the connected component containing $y$ of $V_x \cap
\tilde f^{-n}(\tilde \cF(f^n(y)))$ then the set
$\psi_x^{-1}(W^x_n(y))$ is the graph of a $C^1$-function
$h_{x,y}^n: \DD^2 \to [-1,1]$ with bounded derivatives.
\end{itemize}

By a standard graph transform argument (see \cite{HPS} or
\cite{BuW} section 3) using the fact that these graphs have
bounded derivative we get that $\{h^{x,y}_n\}$ is pre-compact in
the space of functions from $\DD^2$ to $[-1,1]$.

For every $y\in B_\eps(x)$ there exists $\cJ_y^x$ a set of indices
such that for every $\alpha \in \cJ_y^x$ we have a $C^1$-function
$h^{x,y}_{\infty,\alpha}: \DD^2 \to [-1,1]$ and $n_j \to +\infty$
such that:

$$ h^{x,y}_{\infty,\alpha} = \lim_{j \to + \infty} h^{x,y}_{n_j} $$

Every $h^{x,y}_{\infty,\alpha}$ gives rise to a graph whose image
by $\psi_x$ we denote as $W_{\infty,\alpha}^{x}(y)$. This manifold
verifies that it contains $y$ and is everywhere tangent to
$E^{cs}$.

\begin{af}
We have that $H(W^x_{\infty,\alpha}(z)) \en P^{cs}+H(z)$ for every
$z\in B_\eps(x)$ and every $\alpha \in \cJ_z^x$.
\end{af}

\dem Consider $y\in W^x_{\infty,\alpha}(z)$ for  some $\alpha \in
\cJ_z$. One can find $n_j \to \infty$ such that $W^x_{n_j}(z) \to
W^x_{\infty,\alpha}(z)$.

In the coordinates $\psi_x$ of $V_x$, we can find a sequence
$z_{n_j} \in W^x_{n_j}(z) \cap \tilde \cF^u(y)$ such that $z_{n_j}
\to y$. Moreover, we have that $\tilde f^{n_j}(z_{n_j}) \in \tilde
\cF(f^{n_j}(z))$. Assume that $H(y) \neq H(z)$ (otherwise there is
nothing to prove).

We have, by continuity of $H$ that $H(z_{n_j})\to H(y) \neq H(z)$.

We choose a metric in $\RR^3$ so that $(P^{cs})^\perp$ with this
metric is $A$-invariant. We denote as $\lambda$ to the eigenvalue
of $A$ in the direction $(P^{cs})^\perp$.

By Proposition \ref{PropSifIsotopicoAAnosovNoHayTorosEnFoliacion}
and the fact that $H$ is at bounded distance from the identity,
there exists $R>0$ such that for every $n_j\geq 0$ we have that
$A^{n_j}(H(z_{n_j}))$ is at distance smaller than $R$ from $P+
A^{n_j}(H(z))$ since $\tilde f^{n_j}(z_{n_j}) \in \tilde \tilde
\cF(f^{n_j}(z))$.

Suppose that $H(z_{n_j})$ does not converge to $P^{cs}+H(z)$. We
must reach a contradiction.

Consider then $\alpha>0$ such that the angle between $P^{cs}$ and
the vector $H(y)-H(z)$ is larger than $\alpha>0$. This $\alpha$
can be chosen positive under the assumption that $H(z_{n_j})$ does
not converge to $P^{cs} + H(z)$.

Let $n_j>0$ be large enough such that:

 \begin{itemize}
 \item[-] The angle between $A^{-n_j}(P)$ and $P^{cs}$ is smaller than $\alpha/4$,
 \item[-] $\|H(z_{n_j})- H(z)\| > \frac{3}{4} \|H(y)-H(z)\|$,
 \item[-] $\lambda^{n_j} \gg 2R (\sin(\frac{\alpha}{2}) \cos(\beta) \|H(y)-H(z)\|)^{-1}$.
 \end{itemize}

Let $v_{n_j}$ be the vector which realizes $d(H(z_{n_j})-H(z),
A^{-n_j}(P))$  and as $v_{n_j}^\perp$ the projection of $v_{n_j}$
to $(P^{cs})^\perp$. We have that

$$ \|v_{n_j}^\perp \| > \frac{1}{2}\sin\left(\frac{\alpha}{2}\right) \|H(y)-H(z)\|$$

Notice that the distance between $A^{n_j}(H(z_{n_j}))$ and $P +
A^{n_j}(H(z))$ is larger than $\|A^{n_j} v_{n_j}^\perp \|
\cos(\beta)$.

This is a contradiction since this implies that
$A^{n_j}(H(z_{n_j}))$ is at distance larger than

$$ \lambda^{n_j}  \|v_{n_j}^\perp \| \cos(\beta) \gg R $$

\noindent from $P + A^{n_j}(H(z))$. This concludes the claim.

\finobs

Assuming that $P^{cs}$ does not intersect the cone bounded by
$L_1$ and $L_2$ this finishes the proof since one sees that each
leaf of $\tilde \cF^u$ can intersect the pre-image by $H$ of
$P^{cs}+y$ in a unique point, thus showing that the partition of
$\RR^3$ by the pre-images of the translates of $P^{cs}$ defines a
$\tilde f$-invariant foliation (and also invariant under deck
transformations). We leave to the interested reader the task of
filling the details of the proof in this particular case, since we
will continue by giving a proof which works in all cases.

We will prove that $H$ cannot send unstable intervals into the
same plane parallel to $P^{cs}$.

\begin{af}
Given $\gamma:[0,1] \to \RR$ a non-trivial curve contained in
$\tilde \cF^u(x)$ we have that $H(\gamma([0,1]))$ is not contained
in $P^{cs} + H(\gamma(0))$.
\end{af}

\dem Consider $C_\eps$ given by Corollary
\ref{CorolarioConsecuenciasReeb2} (iv) for $\eps$ of the size of
the uniform local product structure. Moreover, consider $L$ large
enough such that $C_\eps L
> \Vol (\TT^3)$.

Since $\tilde \cF^u$ is $\tilde f$-invariant and $P^{cs}$ is
$A$-invariant we deduce that we can assume that the length of
$\gamma$ is arbitrarily large, in particular larger than $2L$.

We will show that $H(B_\eps(\gamma([a,b]))) \en P^{cs} +
H(\gamma(0))$ where $0< a < b < 1$ and the length of
$\gamma([a,b])$ is larger than $L$.

Having volume larger than $\Vol(\TT^3)$ there must be a deck
transformation $\gamma \in \ZZ^3$ such that $\gamma +
B_\eps(\gamma([a,b])) \cap B_\eps(\gamma([a,b]))\neq \emptyset$.
This in turn gives that $\gamma +H(B_\eps(\gamma([a,b]))) \cap
H(B_\eps(\gamma([a,b]))) \neq \emptyset$ and thus $\gamma + P^{cs}
\cap P^{cs} \neq \emptyset$. Since $P^{cs}$ is totally irrational
this is a contradiction.

It remains to show that $H(B_\eps(\gamma([a,b]))) \en P^{cs} +
H(\gamma(0))$. By the previous claim, we know that if $z,w \in
W^x_{\infty,\alpha}(y)$ for some $\alpha \in \cJ_y$, then
$H(z)-H(w) \in P^{cs}$.

Consider $a,b \in [0,1]$ such that $\tilde \cF^u(x) \cap
B_\eps(\gamma([a,b])) \en \gamma([0,1])$. By Corollary
\ref{CorolarioConsecuenciasReeb2} we have that such $a,b$ exist
and we can choose them in order that the length of $\gamma([a,b])$
is larger than $L$.

Let $z\in B_\eps(\gamma([a,b]))$ and choose $w \in \gamma([a,b])$
such that $z \in B_\eps(w)$. We get that for every $\alpha \in
\cJ_z^w$ we have that $W^{w}_{\infty, \alpha}(z) \cap
\gamma([0,1]) \neq \emptyset$. Since $H(\gamma([0,1])) \en P^{cs}
+ H(\gamma(0))$ and by the previous claim, we deduce that $H(w)
\en P^{cs} + H(\gamma(0))$ finishing the proof.

\finobs

Now we are in conditions to show that for every point $x$ and for
every point $y\in B_\eps(x)$ there is a unique manifold
$W^x_{\infty}(y)$ tangent to $E^{cs}$ which is a limit of the
manifolds $W^{x}_n(y)$. Using the same argument as in Proposition
\ref{PropositionBranchingSinBranchingEsFoliacion} we get that the
foliations $\tilde f^{-n}(\tilde \cF)$ converge to a $f$-invariant
foliation $\tilde \cF^{cs}$ tangent to $E^{cs}$ concluding the
proof of the Theorem.

Indeed, assume that the manifolds $W^x_n(y)$ have a unique limit
for every $x \in \RR^3$ and $y\in B_\eps(x)$ and that for any pair
points $y,z\in B_\eps(x)$ these limits are either disjoint or
equal (see the claim below). One has that the set of manifolds
$W^x_\infty(y)$ forms an $f$-invariant plaque family in the
following sense:

\begin{itemize}
\item[-] $\tilde f(W^x_\infty(y)) \cap W^{\tilde f(x)}_\infty
(\tilde f(y))$ is relatively open whenever $\tilde f(y) \in
B_\eps(\tilde f(x))$.
\end{itemize}

We must thus show that these plaque families form a foliation. For
this, we use the same argument as in Proposition
\ref{PropositionBranchingSinBranchingEsFoliacion}. Consider $z,w
\in B_\eps(x)$ we have that $W^x_\infty (z) \cap \tilde \cF^u(w)
\neq \emptyset$ and in fact consists of a unique point (see
Corollary \ref{CorolarioConsecuenciasReeb2} (i)). Since the
intersection point varies continuously and using that plaques are
either disjoint or equal we obtain a continuous map from $\DD^2
\times [-1,1]$ to a neighborhood of $x$ sending horizontal disks
into plaques. This implies that the plaques form an $f$-invariant
foliation as desired.

It thus remains to show the following:

\begin{af}
Given $x \in \RR^3$ and $y,z \in B_\eps(x)$ we have that there is
a unique limit of $W^x_\infty (y)$ and $W^x_\infty(z)$ and they
are either disjoint or coincide. More precisely, for every $\alpha
\in \cJ^x_y$ and $\beta \in \cJ^x_z$ ($z$ could coincide with $y$)
we have that $h^{x,y}_{\infty,\alpha}= h^{x,z}_{\infty,\beta}$ or
the graphs are disjoint.
\end{af}

\dem Assuming the claim does not hold, one obtains $y,z \in
B_\eps(x)$ such that $h^{x,y}_{\infty,\alpha}$ and
$h^{x,z}_{\infty,\beta}$ coincide at some point but whose graphs
are different for some $\alpha \in \cJ_y^x$ and $\beta \in
\cJ_z^x$. In particular, there exists a point $t \in \DD^2$ which
is in the boundary of where both functions coincide. We assume for
simplicity\footnote{If it were not the case we would need to
change the coordinates and perform the same proof, but not to
charge the notation we choose to make this (unnecessary)
assumption.} that $\psi_x (t)$ belongs to $B_\eps(x)$.

Let $\gamma : [0,1] \to B_\eps(x)$ be a non-trivial arc of $\tilde
\cF^u$ joining the graphs of $h^{x,y}_{\infty,\alpha}$ and
$h^{x,z}_{\infty,\beta}$. Since the graphs of both
$h^{x,y}_{\infty,\alpha}$ and $h^{x,z}_{\infty,\beta}$ separate
$V_x$ we have that every point $w \in \gamma((0,1))$ verifies that
for every $\delta \in \cJ_w^x$ one has that $W^x_{\infty,
\delta}(w)$ intersects at least one of $W^x_{\infty,\alpha}(y)$ or
$W^x_{\infty, \beta}(z)$. By the first claim we get that $H(w) \in
P^{cs} + H(y) = P^{cs}+ H(z)$ a contradiction with the second
claim.

\finobs

\lqqd

We can in fact obtain a stronger property since our results allow
us to show that in fact $E^{cs}$ is uniquely integrable into a
foliation. Notice that there are stronger notions of unique
integrability (see \cite{BuW2} and \cite{BF}).

\begin{prop}\label{Prop-UnicidadDeFoliacionTgAEcs}
There is a unique $f$-invariant foliation $\cF^{cs}$  tangent to
$E^{cs}$. Moreover, the plane $P^{cs}$ given by Theorem
\ref{PropObienHayTorosObienTodoEsLindo} for this foliation is
$A$-invariant and contains the stable eigenspace of $A$.
\end{prop}

\dem Assume there are two different $f$-invariant foliations
$\cF^{cs}_1$ and $\cF^{cs}_2$ tangent to $E^{cs}$.

Since they are transverse to $E^u$ they must be Reebless (see
Corollary \ref{CorolarioConsecuenciasReeb}) so that Theorem
\ref{PropObienHayTorosObienTodoEsLindo} applies.

By Remark \ref{RemarkSiAplicofElPlanoSeleAplicaA} we know that
since the foliations are $f$-invariant, the planes $P^{cs}_1$ and
$P^{cs}_2$ given by Theorem
\ref{PropObienHayTorosObienTodoEsLindo} are $A$-invariant. The
fact that $P^{cs}$ contains the stable direction of $A$ is given
by Remark \ref{RemarkSubespaciosInvariantesAnosov} and Corollary
\ref{PropositionQuasiIsometriaII} since it implies that $P^{cs}$
cannot be contained in $E^u_A$.

Assume first that the planes $P^{cs}_1$ and $P^{cs}_2$ coincide.
The foliations remain at distance $R$ from translates of the
planes. By Corollary \ref{PropositionQuasiIsometriaII} we know
that two points in the same unstable leaf must separate in a
direction transverse to $P^{cs}_1=P^{cs}_2$. If $\cF^{cs}_1$ is
different from $\cF^{cs}_2$ we have a point $x$ such that
$\cF^{cs}_1(x) \neq \cF^{cs}_2(x)$. By the global product
structure we get a point $y \in \cF^{cs}_1(x)$ such that $\tilde
\cF^u(y) \cap \cF^{cs}_2(x) \neq \{y\}$. Iterating forward and
using Corollary \ref{PropositionQuasiIsometriaII} we contradict
the fact that leaves of $\cF^{cs}_1$ and $\cF_2^{cs}$ remain at
distance $R$ from translates of $P^{cs}_1=P^{cs}_2$.

Now, if $P^{cs}_1\neq P^{cs}_2$ we know that $A$ has stable
dimension $1$ since we know that $E^s_A$ is contained in both.
Using Corollary \ref{PropositionQuasiIsometriaII} and the fact
that the unstable foliation is $\tilde f$-invariant we see that
this cannot happen.

\lqqd

Notice also that from the proof of Theorem \ref{TeoremaCoherencia}
we deduce that given a foliation $\cF$ transverse to $E^{cs}$ we
have that the backward iterates of this foliation must converge to
this unique $f$-invariant foliation. This implies that:

\begin{cor}\label{Corolario-RobustezDeDireccionTransversal}
Given a dynamically coherent partially hyperbolic diffeomorphism
$f: \TT^3 \to \TT^3$ with splitting $T\TT^3 = E^{cs} \oplus E^u$
isotopic to Anosov we know that it is $C^1$-robustly dynamically
coherent and that the $f_\ast$-invariant plane $P$ given by
Theorem \ref{PropObienHayTorosObienTodoEsLindo} for the unique
$f$-invariant foliation $\cF^{cs}$ tangent to $E^{cs}$ does not
change for diffeomorphisms $C^1$-close to $f$.
\end{cor}

The robustness of dynamical coherence follows from the fact that
being dynamically coherent it is robustly almost dynamically
coherent.

We close this Section with a question we were not able to answer
in full generality:

\begin{quest}
Is it true that $P^{cs}$ corresponds to the eigenspace asociated
to the smallest eigenvalues of $A$?.
\end{quest}

This is true for the case when $A$ has stable index $2$ and we
show in Proposition \ref{PropPlanoCasoIsotopico} that it is the
case in the strong partially hyperbolic case.


\section{Strong partial hyperbolicity and coherence in $\TT^3$}\label{SectionCoherenciaPHfuerte}

 In the strong partially hyperbolic case we are able
to give a stronger result independent of the isotopy class of $f$:

\begin{teo}\label{TEOREMACOHERENCIASPH}
Let $f: \TT^3 \to \TT^3$ be a strong partially hyperbolic
diffeomorphism, then:
\begin{itemize}
\item[-] Either there exists a unique $f$-invariant foliation
$\cF^{cs}$ tangent to $E^s \oplus E^c$ or, \item[-] There exists a
periodic two-dimensional torus $T$ tangent to $E^s \oplus E^c$
which is (normally) repelling.
\end{itemize}
\end{teo}

\begin{obs}\label{RemarkTorosAnosov}
Indeed, it is not hard to show that in the case there is a
repelling torus, it must be an \emph{Anosov tori} as defined in
\cite{HHU3} (see Proposition 2.1 of \cite{BBI} or Lemma
\ref{Lema-fphenT2esisotopicoaAnosov}). In the example of
\cite{HHU} it is shown that the second possibility is not empty.
\finobs
\end{obs}

A diffeomorphism $f$ is \emph{chain-recurrent} if there is no open
set $U$ such that $f(\overline{U})\en U$ (see \cite{Crov-Hab} for
an introduction to this concept in the context of differentiable
dynamics):

\begin{cor*} Let $f: \TT^3 \to \TT^3$ a chain-recurrent strongly partially hyperbolic diffeomorphism.
Then, $f$ is dynamically coherent.
\end{cor*}

In the strong partially hyperbolic case, when no torus tangent to
$E^s\oplus E^c$ nor $E^c\oplus E^u$ exists, we deduce further
properties on the existence of planes close to  the $f$-invariant
foliations. These results are essential to obtain leaf-conjugacy
results (see \cite{Hammerlindl}).

The idea of the proof is to obtain a global product structure
between the foliations involved in order to then get dynamical
coherence. In a certain sense, this is a similar idea to the one
used for the proof of Theorem
\ref{TEOREMA-COHERENCIAISOTOPICOANOSOV}.

However, the fact that global product structure implies dynamical
coherence is much easier in our case due to the existence of
$f$-invariant branching foliations tangent to the center-stable
direction (see subsection \ref{SubSection-GPSimpliesCoherence}).

This approach goes in the inverse direction to the one made in
\cite{BBI2} (and continued in \cite{Hammerlindl}). In \cite{BBI2}
the proof proceeds as follows:

\begin{itemize}
\item[-]First they show that the planes close to the two
foliations are different. To prove this they use absolute
domination. \item[-] They then show (again by using absolute
domination) that leaves of $\tilde \cF^u$ are quasi-isometric.
Here absolute domination is essential (since in the examples of
\cite{HHU} the lift of the unstable foliation is not
quasi-isometric). \item[-] Finally, they use Brin's criterium for
absolutely dominated partially hyperbolic systems (\cite{Brin}) to
obtain coherence. As it was shown in Proposition
\ref{Proposition-BrinArgument} this criterium uses absolute
domination in an essential way.
\end{itemize}

Then, in \cite{Hammerlindl} it is proved that in fact, the planes
$P^{cs}$ and $P^{cu}$ close to the $f$-invariant foliations are
the expected ones in order to obtain global product structure and
then leaf conjugacy to linear models.

Another difference with the proof there is that in our case it
will be important to discuss depending on the isotopy class of $f$
which is not needed in the case of absolute partial hyperbolicity.
In a certain sense, the reason why in each case there is a global
product structure can be regarded as different: In the isotopic to
Anosov case (see subsection
\ref{subsection-ApendicedelPaperPHFol}) we deduce that the
foliations are without holonomy and use Theorem
\ref{Teorema-HectorHirsch} to get global product structure. In the
case which is isotopic to a non-hyperbolic matrix we must first
find out which are the planes close to each foliation in order to
get the global product structure.

\subsection{Preliminary discussions}

Let $f: \TT^3 \to \TT^3$ be a strong partially hyperbolic
diffeomorphism with splitting $T\TT^3= E^s \oplus E^c \oplus E^u$.

We denote as $\cF^s$ and $\cF^u$ to the stable and unstable
foliations given by Theorem \ref{Teorema-VariedadEstableFuerte}
which are one dimensional and $f$-invariant.

As in the previous sections, we will denote as $p:\RR^3 \to \TT^3$
to the covering projection and $\tilde f$ will denote a lift of
$f$ to the universal cover. Recall that $f_\ast :\RR^3 \to \RR^3$
which denotes the linear part of $f$ is at bounded distance
($K_0>0$) from $\tilde f$.

We have already proved:

\begin{teo}\label{TeoremaMio}
Let $f: \TT^3 \to \TT^3$ be a strong partially hyperbolic
diffeomorphism isotopic to Anosov, then $f$ is dynamically
coherent. Moreover, there is a unique $f$-invariant foliation
tangent to $E^{cs}=E^s \oplus E^c$ and a unique $f$-invariant
foliation tangent to $E^{cu}=E^c\oplus E^{u}$.
\end{teo}

This follows from Theorem \ref{TEOREMA-COHERENCIAISOTOPICOANOSOV}
and the fact that strongly partially hyperbolic diffeomorphisms
are almost dynamical coherent (Corollary
\ref{Corolario-PHFuerteESADC}). The uniqueness follows from
Proposition \ref{Prop-UnicidadDeFoliacionTgAEcs}. We will give an
independent proof in subsection
\ref{subsection-ApendicedelPaperPHFol} since in the context of
strong partial hyperbolicity the proof becomes simpler.

The starting point of our proof of Theorem
\ref{TEOREMACOHERENCIASPH} is the existence of $f$-invariant
branching foliations $\cF^{cs}_{bran}$ and $\cF^{cu}_{bran}$
tangent to $E^s\oplus E^c$ and $E^c \oplus E^u$ respectively. By
using Theorem \ref{TeoBuragoIvanov} and Theorem
\ref{PropObienHayTorosObienTodoEsLindo} we can deduce the
following:

\begin{prop}\label{ProposicionDicotomiaEnLaFoliacion}
There exist an $f_\ast$-invariant plane $P^{cs}$ and $R>0$ such
that every leaf of $\tilde \cF^{cs}_{bran}$ lies in the
$R$-neighborhood of a plane parallel to $P^{cs}$.

Moreover, one can choose $R$ such that one of the following
conditions holds:
\begin{itemize}
\item[(i)] The projection of the plane $P^{cs}$ is dense in
$\TT^3$ and the $R$-neighborhood of every leaf of $\tilde
\cF^{cs}_{bran}$ contains a plane parallel to $P^{cs}$, or,
\item[(ii)] The projection of $P^{cs}$ is a linear two-dimensional
torus and there is a leaf of $\cF^{cs}_{bran}$ which is a
two-dimensional torus homotopic to $p(P^{cs})$.
\end{itemize}
An analogous dichotomy holds for $\cF^{cu}_{bran}$.
\end{prop}

\dem We consider sufficiently small $\eps>0$ and the foliation
$\cS_\eps$ given by Theorem \ref{TeoBuragoIvanov}.

Let $h^{cs}_\eps$ be the continuous and surjective map which is
$\eps$-close to the identity sending leaves of $\cS_\eps$ into
leaves of $\cF^{cs}_{bran}$. By taking the lift to the universal
cover, we have that there is $\tilde h^{cs}_\eps : \RR^3 \to
\RR^3$ continuous and surjective which is also at distance smaller
than $\eps$ from the identity such that it sends leaves of $\tilde
\cS_\eps$ homeomorphically into leaves of $\tilde
\cF^{cs}_{bran}$.

This implies that given a leaf $L$ of $\tilde \cF^{cs}_{bran}$
there exists a leaf $S$ of $\tilde \cS_\eps$ such that $L$ is at
distance smaller than $L$ from $S$ and viceversa.

Since the foliation $\cS_\eps$ is transverse to $E^u$ we can apply
Theorem \ref{PropObienHayTorosObienTodoEsLindo} and we obtain that
there exists a plane $P^{cs}$ and $R>0$ such that every leaf of
the lift $\tilde \cS_\eps$ of $\cS_\eps$ to $\RR^3$ lies in an
$R$-neighborhood of a translate of $P^{cs}$. Recall that this
plane is unique (see Remark \ref{Remark-UnicidaddelPlanoP}).

From the previous remark, we get that every leaf of $\tilde
\cF^{cs}_{bran}$ lies in an $R+ \eps$-neighborhood of a translate
of $P^{cs}$ and this is the unique plane with this property.

Since $\tilde \cF^{cs}_{bran}$ is $\tilde f$-invariant, we deduce
that the plane $P^{cs}$ is $f_\ast$-invariant (see also Remark
\ref{RemarkSiAplicofElPlanoSeleAplicaA}).

By Proposition \ref{Proposicion-SiPesToroHayHojaToro} we know that
if $P^{cs}$ projects into a two-dimensional torus, we obtain that
the foliation $\cS_\eps$ must have a torus leaf. The image of this
leaf by $h^{cs}_\eps$ is a torus leaf of $\cF^{cs}_{bran}$. This
gives (ii).

Since a plane whose projection is not a two-dimensional torus must
be dense we get that if option (ii) does not hold, we have that
the image of $P^{cs}$ must be dense. Moreover, option (i) of
Theorem \ref{PropObienHayTorosObienTodoEsLindo} must hold for
$\cS_\eps$ and this concludes the proof of this proposition.

\lqqd

\begin{obs}\label{Remark-PlanosYSubespaciosInvariantes}
Assume that $f: \TT^3 \to \TT^3$ is a strongly partially
hyperbolic diffeomorphism which is not isotopic to Anosov. By
Theorem \ref{TeoremaLaAccionEnHomologiaEsPH} and Corollary
\ref{Corolario-PHFuerteESADC} we have that if $f$ is not  isotopic
to Anosov, then $f_\ast$ is in the hypothesis of Lemma
\ref{LemaPlanosInvariantes}. Let $P$ be an $f_\ast$-invariant
plane, then there are the following $3$ possibilities:
\begin{itemize}
\item[-] $P$ may project into a torus. In this case, $P= E^s_\ast
\oplus E^u_\ast$ (the eigenplane corresponding to the eigenvalues
of modulus different from one). \item[-] If $P=E^s_\ast \oplus
E^c_\ast$ then $P$ projects into an immersed cylinder which is
dense in $\TT^3$. \item[-] If $P=E^c_\ast \oplus E^u_\ast$ then
$P$ projects into an immersed cylinder which is dense in $\TT^3$.
\end{itemize}
\finobs
\end{obs}

\subsection{Global product structure implies dynamical coherence}\label{SubSection-GPSimpliesCoherence}

Assume that $f: \TT^3 \to \TT^3$ is a strong partially hyperbolic
diffeomorphism. Let $\cF^{cs}_{bran}$ be the $f$-invariant
branching foliation tangent to $E^s\oplus E^c$ given by Theorem
\ref{TeoBuragoIvanov} and let $\cS_\eps$ be a foliation tangent to
an $\eps$-cone around $E^s \oplus E^c$ which remains $\eps$-close
to the lift of $\cF^{cs}_{bran}$ to the universal cover for small
$\eps$.

When the lifts of $\cS_\eps$ and $\cF^u$ to the universal cover
have a global product structure, we deduce from Proposition
\ref{PropositionQuasiIsometria} the following:

\begin{cor}\label{PropQuasiIsometria}
If $\cS_\eps$ and $\cF^u$ have global product structure, then, the
foliation $\tilde \cF^u$ is quasi-isometric. Indeed, if $v \in
(P^{cs})^\perp$ is a unit vector, there exists $\ell>0$ such that
for every $n\geq 0$, every unstable curve starting at a point $x$
of length larger than $n \ell$ intersects $P^{cs} + n v + x$ or
$P^{cs} - n v + x$.
\end{cor}

Before we show that global product structure implies coherence, we
will show an equivalence to having global product structure
between $\tilde \cF^u$ and $\tilde \cS_\eps$ which will sometimes
be better adapted to our proofs.

\begin{lema}\label{Lemma-GPSentreSepsybranched}
There exists $\eps>0$ such that $\tilde \cF^u$ and $\tilde
\cS_\eps$ have global product structure if and only if:
\begin{itemize}
\item[-] For every $x,y \in \RR^3$ and for every $L \in \tilde
\cF^{cs}_{bran}(y)$ we have that $\tilde \cF^u(x) \cap L \neq
\emptyset$.
\end{itemize}
\end{lema}

\dem First notice that any of the hypothesis implies that $\tilde
\cS_\eps$ cannot have dead-end components. In particular, there
exists $R>0$ and a plane $P^{cs}$ such every leaf of $\tilde
\cS_\eps$ and every leaf of $\tilde \cF^{cs}_{bran}$ verifies that
it is contained in an $R$-neighborhood of a translate of $P^{cs}$
and the $R$-neighborhood of the leafs contains a translate of
$P^{cs}$ too (see Proposition
\ref{ProposicionDicotomiaEnLaFoliacion}).

We prove the direct implication first. Consider $x,y \in \RR^3$
and $L$ a leaf of $\tilde \cF^{cs}_{bran}(y)$. Now, we know that
$L$ separates in $\RR^3$ the planes $P^{cs}+ y + 2R$ and $P^{cs} +
y -2R$. One of them must be in the connected component of $\RR^3
\setminus L$ which not contains $x$, without loss of generality we
assume that it is $P^{cs} + y+ 2R$. Now, we know that there is a
leaf $S$ of $\tilde \cS_\eps$ which is contained in the half space
bounded by $P^{cs} + y +R$ not containing $L$ (notice that $L$
does not intersect $P^{cs}+ y+R$). Global product product
structure implies that $\tilde \cF^u(x)$ intersects $S$ and thus,
it also intersects $L$.

The converse direction has an analogous proof.

\lqqd

We can prove the following result which does not make use of the
isotopy class of $f$.

\begin{prop}\label{Proposition-GPSIMPLICACOHERENCIA}
Assume that there is a global product structure between the lift
of $\cS_\eps$ and the lift of $\cF^u$ to the universal cover. Then
there exists an $f$-invariant foliation $\cF^{cs}$ everywhere
tangent to $E^s \oplus E^u$.
\end{prop}

\dem We will show that the branched foliation $\tilde
\cF^{cs}_{bran}$ must be a true foliation (it cannot be branched
and use Proposition
\ref{PropositionBranchingSinBranchingEsFoliacion})).

Assume otherwise, i.e. there exists $x\in \RR^3$ such that $\tilde
\cF^{cs}_{bran}(x)$ has more than one complete surface. We call
$L_1$ and $L_2$ different leaves in $\tilde \cF^{cs}_{bran}(x)$.
There exists $y$ such that $y \in L_1 \setminus L_2$. Using global
product structure and Lemma \ref{Lemma-GPSentreSepsybranched} we
get $z \in L_2$ such that:

\begin{itemize}
\item[-] $y \in \tilde \cF^u(z)$.
\end{itemize}

Consider $\gamma$ the arc in $\tilde \cF^u(z)$ whose endpoints are
$y$ and $z$. Let $R$ be the value given by Proposition
\ref{ProposicionDicotomiaEnLaFoliacion} and $\ell>0$ given by
Corollary \ref{PropQuasiIsometria}. We consider $N$ large enough
so that $\tilde f^N(\gamma)$ has length larger than $n \ell$ with
$n \gg R$.

By Corollary \ref{PropQuasiIsometria} we get that the distance
between $P^{cs}+ \tilde f^N(z)$ and $\tilde f^N(y)$ is much larger
than $R$. However, we have that, by $\tilde f$-invariance of
$\tilde \cF^{cs}_{bran}$ there is a leaf of $\tilde
\cF^{cs}_{bran}$ containing both $\tilde f^N(z)$ and $\tilde
f^N(x)$ and another one containing both $\tilde f^N(y)$ and
$\tilde f^N(x)$. This contradicts Proposition
\ref{ProposicionDicotomiaEnLaFoliacion} showing that $\tilde
\cF^{cs}_{bran}$ must be a true foliation.

\lqqd

\subsection{Torus leafs}

This subsection is devoted to the proof of the following:

\begin{lema}\label{LemaNohayToros}
If $\cF^{cs}_{bran}$ contains a leaf which is a two-dimensional
torus, then there exists a leaf of $\cF^{cs}_{bran}$ which is a
torus and it is fixed by $f^k$ for some $k\geq 1$. Moreover, this
leaf is normally repelling.
\end{lema}

\dem Let $T \en \TT^3$ be a leaf of $\cF^{cs}_{bran}$ homeomorphic
to a two-torus. Since $\cF^{cs}_{bran}$ is $f$-invariant and
$P^{cs}$ is invariant under $f_\ast$ we get that the image of $T$
by $f$ is homotopic to $T$ and a leaf of $\cF^{cs}_{bran}$.

Notice that having an $f_\ast$-invariant plane which projects into
a torus already implies that $f_\ast$-cannot be hyperbolic (see
Proposition \ref{PropAnosovSonIrreducibles}).

By Remark \ref{Remark-PlanosYSubespaciosInvariantes} we have that
the plane $P^{cs}$ coincides with $E^s_\ast \oplus E^u_\ast$ (the
eigenspaces corresponding to the eigenvalues of modulus different
from $1$ of $f_\ast$).

Since the eigenvalue of $f_\ast$ in $E^c_\ast$ is of modulus $1$,
this implies that if we consider two different lifts of $T$, then
they remain at bounded distance when iterated by $\tilde f$.
Indeed, if we consider two different lifts $\tilde T_1$ and
$\tilde T_2$ of $T$ we have that $\tilde T_2 = \tilde T_1 +\gamma$
with $\gamma \in E^c_\ast \cap \ZZ^3$. Now, we have that $\tilde f
(\tilde T_2) = \tilde f(\tilde T_1) + f_\ast(\gamma) = \tilde
f(\tilde T_1) \pm \gamma$.

We shall separate the proof depending on how the orbit of $T$ is.
\medskip

{\bf Case 1:} Assume the torus $T$ is fixed by some iterate $f^n$
of $f$ with $n\geq 1$. Then, since it is tangent to the center
stable distribution, we obtain that it must be normally repelling
as desired.
\medskip

{\bf Case 2:} If the orbit of $T$ is dense, we get that
$\cF^{cs}_{bran}$ is a true foliation by two-dimensional torus
which we call $\cF^{cs}$ from now on. This is obtained by the fact
that one can extend the foliation to the closure using the fact
that there are no topological crossings between the torus leaves
(see Proposition \ref{Proposition-BWProp16}).

Since all leaves must be two-dimensional torus  which are
homotopic we get that the foliation $\cF^{cs}$ has no holonomy
(see Theorem \ref{Teorema-EstabilidadCompleta} and Proposition
\ref{Proposition-SucesionDeTorosEnBranchingFol}).

Using Theorem \ref{Teorema-HectorHirsch}, we get that the unstable
direction $\tilde \cF^u$ in the universal cover must have a global
product structure with $\tilde \cF^{cs}$.

Let $S$ be a leaf of $\cF^{cs}$ and consider $\tilde S_1$ and
$\tilde S_2$ two different lifts of $S$ to $\RR^3$.

Consider an arc $J$ of $\tilde \cF^u$ joining $\tilde S_1$ to
$\tilde S_2$. Iterating the arc $J$ by $\tilde f^n$ we get that
its length grows exponentially, while the extremes remain the the
forward iterates of $\tilde S_1$ and $\tilde S_2$ which remain at
bounded distance by the argument above.

By considering translations of one end of $\tilde f^n(J)$ to a
fundamental domain and taking a convergent subsequence we obtain a
leaf of $\tilde \cF^u$ which does not intersect every leaf of
$\tilde \cF^{cs}$. This contradicts global product structure.

\medskip

{\bf Case 3:} Let $T_1, T_2 \in \cF^{cs}_{bran}$ two different
torus leaves. Since there are no topological crossings, we can
regard $T_2$ as embedded in $\TT^2 \times [-1,1]$ where both
boundary components are identified with $T_1$ and such that the
embedding is homotopic to the boundary components (recall that any
pair of torus leaves must be homotopic). In particular, we get
that $\TT^3 \setminus (T_1 \cup T_2)$ has at least two different
connected components and each of the components has its boundary
contained in $T_1 \cup T_2$.

If the orbit of $T$ is not dense, we consider
$\cO=\overline{\bigcup_n f^n(T)}$ the closure of the orbit of $T$
which is an invariant set.

Recall that we can assume completeness of $\cF^{cs}_{bran}$ (i.e.
for every $x_n \to x$ and $L_n \in \cF^{cs}_{bran}(x_n)$ we have
that $L_n$ converges in the $C^1$-topology to $L_\infty \in
\cF^{cs}_{bran}(x)$). We get that $\cO$ is saturated by leaves of
$\cF^{cs}_{bran}$ all of which are homotopic torus leaves (see
Proposition \ref{Proposition-SucesionDeTorosEnBranchingFol}).

Let $U$ be a connected component of the complement of $\cO$. By
the previous remarks we know that its boundary $\partial U$ is
contained in the union of two torus leaves of $\cF^{cs}_{bran}$.

If some component $U$ of $\cO^c$ verifies that there exists $n\geq
1$ such that $f^n(U) \cap U \neq \emptyset$, by invariance of
$\cO^c$ we get that $f^{2n}$ fixes both torus leaves whose union
contains $\partial U$. This implies the existence of a periodic
normally repelling torus as in Case 1.

We claim that if every connected component of $\cO^c$ is
wandering, then we can show that every leaf of $\tilde \cF^u$
intersects every leaf of $\tilde \cF^{cs}_{bran}$ which allows to
conclude exactly as in Case 2.

To prove the claim, consider $\delta$ given by the local product
structure between these two transverse foliations (one of them
branched). This means that given $x,y$ such that $d(x,y)<\delta$
we have that $\tilde \cF^u(x)$ intersects every leaf of $\tilde
\cF^{cs}_{bran}$ passing through $y$.

Assume there is a point $x\in \RR^3$ such that $\tilde \cF^u(x)$
does not intersect every leaf of $\tilde \cF^{cs}_{bran}$. As in
subsection \ref{SubSection-FoliacionesDeT3} we know that each leaf
of $\tilde \cF^{cs}_{bran}$ separates $\RR^3$ into two connected
components so we can choose among the lifts of torus leaves, the
leaf $\tilde T_0$ which is the lowest (or highest depending on the
orientation of the semi-unstable leaf of $x$ not intersecting
every leaf of $\tilde \cF^{cs}_{bran}$) not intersecting $\tilde
\cF^{u}(x)$. We claim that $\tilde T_0$ must project by the
convering projection into a torus leaf which intersects the
boundary of a connected component of $\cO^c$. Indeed, there are
only finitely many connected components $U_1, \ldots, U_N$ of
$\cO^c$ having volume smaller than the volume of a $\delta$-ball,
so if a point is not in $U_i$ for some $i$, we know that it must
be covered by local product structure boxes forcing its unstable
leaf to advance until one of those components.

On the other hand, using $f$-invariance of $\cF^{u}$ and the fact
that every connected component of $\cO^c$ is wandering, we get
that every point in $U_i$ must eventually fall out of $\bigcup_i
U_i$ and then its unstable manifold must advance to other
component. This concludes the claim, and as we explained, allows
to use the same argument as in Case 2 to finish the proof in Case
3.
\lqqd

M.A. Rodriguez Hertz and R. Ures were kind to comunicate an
alternative proof of this lemma by using an adaptation of an
argument due to Haefliger for branched foliations (it should
appear in \cite{HHU}).

\subsection{Obtaining Global Product Structure}

In this section we will prove the following result which will
allow us to conclude in the case where $f_\ast$ is not isotopic to
Anosov.

\begin{prop}\label{PropGPSyPlanos}
Let $f: \TT^3 \to \TT^3$ be a strongly partially hyperbolic
diffeomorphism which is not isotopic to Anosov and does not have a
periodic two-dimensional torus tangent to $E^s\oplus E^c$. Then,
the plane $P^{cs}$ given by Proposition
\ref{ProposicionDicotomiaEnLaFoliacion} corresponds to the
eigenplane corresponding to the eigenvalues of modulus smaller or
equal to $1$. Moreover, there is a global product structure
between $\tilde \cF^{cs}_{bran}$ and $\tilde \cF^u$. A symmetric
statement holds for $\tilde \cF^{cu}_{bran}$ and $\tilde \cF^s$.
\end{prop}

As noted in Remark \ref{Remark-PlanosYSubespaciosInvariantes} we
get that even if a strongly partially hyperbolic diffeomorphism is
not isotopic to Anosov, then, $f_\ast$ still must have one
eigenvalue of modulus larger than one and one smaller than one.

The mentioned remark also gives that there are exactly three
$f_\ast$-invariant lines $E^s_\ast$, $E^c_\ast$ and $E^u_\ast$
corresponding to the eigenvalues of $f_\ast$ of modulus smaller,
equal and larger than one respectively.

\begin{lema}\label{LemaEstableNoSeQuedaCercaDeCentroInestable}
For every $R>0$ and $x\in \RR^3$ we have that $\tilde \cF^u(x)$ is
not contained in an $R$-neighborhood of $(E^s_\ast \oplus
E^c_\ast) + x$. Symmetrically, for every $R>0$ and $x \in \RR^3$
the leaf $\tilde \cF^s(x)$ is not contained in an $R$-neighborhood
of $(E^c_\ast \oplus E^u_\ast) +x$.
\end{lema}

\dem Let $C$ be a connected set contained in an $R$-neighborhood
of a translate of $E^s_\ast \oplus E^c_\ast$, we will estimate the
diameter of $\tilde f(C)$ in terms of the diameter of $C$.

\begin{af}
There exists $K_R$ which depends only on $\tilde f$, $f_\ast$ and
$R$ such that:
$$ \diam (\tilde f(C)) \leq \diam(C) + K_R $$
\end{af}

\dem  Let $K_0$ be the $C^0$-distance between $\tilde f$ and
$f_\ast$ and consider $x,y \in C$ we get that:

$$ d(\tilde f(x),\tilde f(y)) \leq d(f_\ast(x),f_\ast(y)) + d(f_\ast(x),\tilde f(x))+ d(f_\ast(y),\tilde f(y)) \leq $$
$$ \leq d(f_\ast(x), f_\ast(y)) + 2 K_0$$

We have that the difference between $x$ and $y$ in the unstable
direction of $f_\ast$ is bounded by $2R$ given by the distance to
the plane $E^s_\ast \oplus E^u_\ast$ which is transverse to
$E^u_\ast$.

Since the eigenvalues of $f_\ast$ along $E^s_\ast \oplus E^c_\ast$
we have that $f_\ast$ does not increase distances in this
direction: we thus have that $d(f_\ast(x),f_\ast(y)) \leq d(x,y) +
2 |\lambda^u| R$ where $\lambda^u$ is the eigenvalue of modulus
larger than $1$. We have obtained:

$$ d(\tilde f(x),\tilde f(y))  \leq d(x,y) + 2K_0 + 2|\lambda^u| R = d(x,y) + K_R$$

\noindent which concludes the proof of the claim. \finobs

Now, this implies that if we consider an arc $\gamma$ of $\tilde
\cF^u$ of length $1$ and assume that its future iterates remain in
a slice parallel to $E^s_\ast \oplus E^c_\ast$ of width $2R$ we
have that

$$\diam( \tilde f^n(\gamma)) < \diam(\gamma) + n K_R \leq 1 + n K_R$$

So that the diameter grows linearly with $n$.

The volume of balls in the universal cover of $\TT^3$ grows
polynomially with the radius (see Step 2 of \cite{BBI} or page 545
of \cite{BI}, notice that the universal) so that we have that
$B_\delta(\tilde f^{-n}(\gamma))$ has volume which is polynomial
$P(n)$ in $n$.

On the other hand, we know from the partial hyperbolicity that
there exists $C>0$ and $\lambda>1$ such that the length of $\tilde
f^n(\gamma)$ is larger than $C \lambda^n$.

Using Corollary \ref{CorolarioConsecuenciasReeb} (iv), we obtain
that there exists $n_0$ uniform such that every arc of length $1$
verifies that $\tilde f^{n_0}(\gamma)$ is not contained in the
$R$-neighborhood of a translate of $E^s_\ast \oplus E^c_\ast$.
This implies that no unstable leaf can be contained in the
$R$-neighborhood of a translate of $E^s_\ast \oplus E^c_\ast$
concluding the proof of the lemma.

\lqqd

We are now ready to prove Proposition \ref{PropGPSyPlanos}

\demo{ of Proposition \ref{PropGPSyPlanos}} Consider the plane
$P^{cs}$ given by Proposition
\ref{ProposicionDicotomiaEnLaFoliacion} for the branching
foliation $\cF^{cs}_{bran}$.

If option (ii) of Proposition
\ref{ProposicionDicotomiaEnLaFoliacion} holds, we get that there
must be a torus leaf in $\cF^{cs}_{bran}$ which we assume there is
not.

By Lemma \ref{LemaPlanosInvariantes} and Remark
\ref{Remark-PlanosYSubespaciosInvariantes} the plane $P^{cs}$ must
be either $E^s_\ast \oplus E^c_\ast$ or $E^c_\ast \oplus
E^u_\ast$.

Lemma \ref{LemaEstableNoSeQuedaCercaDeCentroInestable} implies
that $P^{cs}$ cannot be $E^c_\ast \oplus E^u_\ast$ since $\tilde
\cF^s$ is contained in $\tilde \cF^{cs}_{bran}$. This implies that
$P^{cs}= E^s_\ast \oplus E^c_\ast$ as desired.

Now, using Lemma \ref{LemaEstableNoSeQuedaCercaDeCentroInestable}
for $\tilde \cF^u$ we see that the unstable foliation cannot
remain close to a translate of $P^{cs}$. This implies that $\tilde
\cF^u$ intersects every translate of $P^{cs}$ and since every leaf
of $\tilde \cS_\eps$ is contained in between two translates of
$P^{cs}$ which are separated by the leaf, we deduce that every
leaf of $\tilde \cF^u$ intersects every leaf of $\tilde \cS_\eps$.
Now, by Lemma \ref{Lemma-GPSentreSepsybranched} gives ``global
prooduct structure'' between $\tilde \cF^u$ and $\tilde
\cF^{cs}_{bran}$ and using Proposition
\ref{Proposition-GPSIMPLICACOHERENCIA}.

\lqqd

\subsection{Proof of Theorem \ref{TEOREMACOHERENCIASPH}}\label{SectionCoherenciaPrueba}

To prove Theorem \ref{TEOREMACOHERENCIASPH}, we first assume that
$f_\ast$ is not isotopic to Anosov.

If there is a torus tangent to $E^s \oplus E^c$, then, by Lemma
\ref{LemaNohayToros} we obtain a periodic normally repelling
torus.

By Proposition \ref{PropGPSyPlanos} we get that if there is no
repelling torus, then there is a global product structure. Now,
Proposition \ref{Proposition-GPSIMPLICACOHERENCIA} gives the
existence of an $f$-invariant foliation $\cF^{cs}$ tangent to $E^s
\oplus E^c$ (see also Lemma \ref{Lemma-GPSentreSepsybranched}).

The proof shows that there must be a unique $f$-invariant
foliation tangent to $E^{cs}$ (and to $E^{cu}$).

Indeed, we get that every foliation tangent to $E^{cs}$ must
verify option (i) of Proposition
\ref{ProposicionDicotomiaEnLaFoliacion} when lifted to the
universal cover and that the plane which is close to the foliation
must correspond to the eigenspace of $f_\ast$ corresponding to the
smallest eigenvalues (Proposition \ref{PropGPSyPlanos}).

Using quasi-isometry of the strong foliations, this implies that
if there is another surface tangent to $E^{cs}$ through a point
$x$, then this surface will not extend to an $f$-invariant
foliation since we get that forward iterates will get arbitrarily
far from this plane (this is proved exactly as Proposition
\ref{Proposition-GPSIMPLICACOHERENCIA}).

This concludes the proof of Theorem \ref{TEOREMACOHERENCIASPH} in
case $f$ is not isotopic to Anosov, Theorem \ref{TeoremaMio}
concludes.

\lqqd

It may be that there are other (non-invariant) foliations tangent
to $E^{cs}$ (see \cite{BF}) or, even if there are no such
foliations there may be complete surfaces tangent to $E^{cs}$
which do not extend to foliations. The techniques here presented
do not seem to be enough to discard such situations.

\subsection{A simpler proof of Theorem \ref{TeoremaMio}. The
isotopy class of Anosov.}\label{subsection-ApendicedelPaperPHFol}

In Section \ref{Section-PHAnosov} Theorem \ref{TeoremaMio} is
obtained as a consequence of a more general result which is harder
to prove. We present here a simpler proof of this result.

\demo{ of Theorem \ref{TeoremaMio}} Let $\cF^{cs}_{bran}$ be the
branched foliation tangent to $E^{cs}$ given by Theorem
\ref{TeoBuragoIvanov}. By Proposition
\ref{ProposicionDicotomiaEnLaFoliacion} we get a
$f_\ast$-invariant plane $P^{cs}$ in $\RR^3$ which we know cannot
project into a two-dimensional torus since $f_\ast$ has no
invariant planes projecting into a torus (see Remark
\ref{RemarkSubespaciosInvariantesAnosov}), this implies that
option (i) of Proposition \ref{ProposicionDicotomiaEnLaFoliacion}
is verified.

Since for every $\eps>0$, Theorem \ref{TeoBuragoIvanov} gives us a
foliation $\cS_\eps$ whose lift is close to $\tilde \cF^{cs}$, we
get that the foliation $\tilde \cS_\eps$ remains close to $P^{cs}$
which must be totally irrational (see Remark
\ref{RemarkSubespaciosInvariantesAnosov}). By Lemma
\ref{RemarkHojasCerradas} (i) we get that all leaves of $\cS_\eps$
are simply connected, thus, we get that the foliation $\cS_\eps$
is without holonomy.

We can apply Theorem \ref{Teorema-HectorHirsch} and we obtain that
for every $\eps>0$ there is a global product structure between
$\tilde \cS_\eps$ and $\tilde \cF^u$ which is transverse to
$\cS_\eps$ if $\eps$ is small enough.

The rest of the proof follows from Proposition.
\ref{Proposition-GPSIMPLICACOHERENCIA}. \lqqd

In fact, using the same argument as in Proposition
\ref{Prop-UnicidadDeFoliacionTgAEcs} we get uniqueness of the
foliation tangent to $E^s\oplus E^c$.

We are also able to prove the following proposition which is
similar to Proposition \ref{PropGPSyPlanos} in the context of
partially hyperbolic diffeomorphisms isotopic to Anosov, this will
be used in \cite{HP} to obtain leaf conjugacy to the linear model.

Notice first that the eigenvalues of $f_\ast$ verify that they are
all different (see Lemma \ref{RemarkValoresPropiosAnosov} and
Proposition \ref{PropositionAnotieneinestablesComplejos}).

We shall name them $\lambda_1, \lambda_2, \lambda_3$ and assume
they verify:

$$ |\lambda_1|< |\lambda_2|< |\lambda_3| \quad ; \quad |\lambda_1|<1 \ , \ |\lambda_2|\neq 1 \ , \ |\lambda_3|>1 $$

\noindent we shall denote as $E^i_\ast$ to the eigenline of
$f_\ast$ corresponding to $\lambda_i$.

\begin{prop}\label{PropPlanoCasoIsotopico}
The plane close to the branched foliation $\tilde \cF^{cs}$
corresponds to the eigenplane corresponding to the eigenvalues of
smaller modulus (i.e. the eigenspace $E^{1}_\ast \oplus E^2_\ast$
corresponding to $\lambda_1$ and $\lambda_2$). Moreover, there is
a global product structure between $\tilde \cF^{cs}$ and $\tilde
\cF^u$. A symmetric statement holds for $\tilde \cF^{cu}$ and
$\tilde \cF^s$.
\end{prop}

\dem This proposition follows from the existence of a
semiconjugacy $H$ between $\tilde f$ and its linear part $f_\ast$
which is at bounded distance from the identity.

The existence of a global product structure was proven above.
Assume first that $|\lambda_2|<1$, in this case, we know that
$\tilde \cF^u$ is sent by the semiconjugacy into lines parallel to
the eigenspace of $\lambda_3$ for $f_\ast$. This readily implies
that $P^{cs}$ must coincide with the eigenspace of $f_\ast$
corresponding to $\lambda_1$ and $\lambda_2$ otherwise we would
contradict the global product structure.

The case were $|\lambda_2|>1$ is more difficult. First, it is not
hard to show that the eigenspace corresponding to $\lambda_1$ must
be contained in $P^{cs}$ (otherwise we can repeat the argument in
Lemma \ref{LemaEstableNoSeQuedaCercaDeCentroInestable} to reach a
contradiction).

Assume by contradiction that $P^{cs}$ is the eigenspace
corresponding to $\lambda_1$ and $\lambda_3$.

First, notice that by the basic properties of the semiconjugacy
$H$, for every $x\in \RR^3$ we have that $\tilde \cF^u(x)$ is sent
by $H$ into $E^u_\ast + H(x)$ (where $E^u_\ast = E^2_\ast \oplus
E^3_\ast$ is the eigenspace corresponding to $\lambda_2$ and
$\lambda_3$ of $f_\ast$).

We claim that this implies that in fact $H(\tilde \cF^u(x)) =
E^2_\ast + H(x)$ for every $x\in \RR^3$. In fact, we know from
Corollary \ref{PropositionQuasiIsometriaII} that points of
$H(\tilde \cF(x))$ which are sufficiently far apart are contained
in a cone of $(E^2_\ast \oplus E^3_\ast) + H(x)$ bounded by two
lines $L_1$ and $L_2$ which are transverse to $P^{cs}$. If
$P^{cs}$ contains $E^3_\ast$ this implies that if one considers
points in the same unstable leaf which are sufficiently far apart,
then their image by $H$ makes an angle with $E^3_\ast$ which is
uniformly bounded from below. If there is a point $y \in \tilde
\cF^u(x)$ such that $H(y)$ not contained in $E^2_\ast$ then we
have that $d(\tilde f^n(y), \tilde f^n(x))$ goes to $\infty$ with
$n$ while the angle of $H(y) - H(x)$ with $E^3_\ast$ converges to
$0$ exponentially contradicting Corollary
\ref{PropositionQuasiIsometriaII}.

Consider now a point $x\in \RR^3$ and let $y$ be a point which can
be joined to $x$ by a finite set of segments $\gamma_1, \ldots,
\gamma_k$ tangent either to $E^s$ or to $E^u$ (an $su$-path, see
subsection \ref{SubSection-Accessibility}). We know that each
$\gamma_i$ verifies that $H(\gamma_i)$ is contained either in a
translate of $E^1_\ast$ (when $\gamma_i$ is tangent to $E^s$, i.e.
it is an arc of the strong stable foliation $\tilde \cF^s$) or in
a translate of $E^2_\ast$ (when $\gamma_i$ is tangent to $E^u$
from what we have shown in the previous paragraph). This implies
that the \emph{accesibility} class of $x$ verifies that its image
by $H$ is contained in $(E^1_\ast \oplus E^2_\ast) + H(x)$. The
projection of $E^1_\ast \oplus E^2_\ast$ to the torus is not the
whole $\TT^3$ so in particular, we get that $f$ cannot be
accesible. From Corollary
\ref{Corolario-RobustezDeDireccionTransversal} this situation
should be robust under $C^1$-perturbations since those
perturbations cannot change the direction of $P^{cs}$.

On the other hand, Theorem \ref{Teorema-AccesibilidadDensa}
implies that by an arbitrarily small ($C^1$ or $C^r$) perturbation
of $f$ one can make it accessible. This gives a contradiction and
shows that $P^{cs}$ must coincide with $E^{1}_\ast \oplus
E^2_\ast$ as desired.

\lqqd




\section{Higher dimensions}\label{Section-DimensionesMayores}

In this section we attempt to find conditions that guarantee a
partially hyperbolic diffeomorphism to be isotopic to an Anosov
diffeomorphism. The progress made so far in this direction is not
as strong though we have obtained some partial results. We present
here part of what will appear in \cite{PotPHtrapping}.

We start by giving the property we will require for a partially
hyperbolic diffeomorphism and hope it implies being isotopic to an
Anosov diffeomorphism.

\begin{defi}[Coherent trapping property]\label{DefiCTP}
Let $f\in \Diff^1(M)$ be a dynamically coherent partially
hyperbolic diffeomorphism of type $TM =E^{cs} \oplus E^u$. We
shall say that it admits the \emph{coherent trapping
property}\footnote{The word coherent (motivated by the existence
of an invariant $cs-$foliation) is included to distinguish it from
the a priori weaker condition of only having a plaque family
trapped by $f$ (it could be that this condition alone implies
coherence, see \cite{BuFi} for progress in that direction).} if
there exists a continuous map $\cD^{cs}: M \to Emb^1(\DD^{cs}, M)$
such that $\cD^{cs}(x)(0)=x$, the image of $\DD^{cs}$ by
$\cD^{cs}(x)$ is always contained in $\cF^{cs}(x)$ and they verify
the following trapping property:
$$f(\cD^{cs}(x)(\DD^{cs})) \en \cD^{cs}(f(x))(int(\DD^{cs}))
\qquad \forall x\in M.$$ \finobs
\end{defi}

For notational purposes, and with the risk of abusing notation, we
shall denote from now on: $\overline{\cD^{cs}_x} =
\cD^{cs}(x)(\DD^{cs})$ and $\cD^{cs}_x =
\cD^{cs}(x)(int(\DD^{cs}))$.

We remark the important point that there is no restriction on the
size of the plaques $\cD^{cs}_x$ so the dynamics can be quite rich
in the center stable plaques.

We will prove the following:

\begin{teo}\label{Teorema-SemiconjALineal}
Let $f : M \to M$ be a partially hyperbolic diffeomorphism with
splitting $TM=E^{cs} \oplus E^u$ having the coherent trapping
property and such that one of the following conditions holds:
\begin{itemize}
\item[-] $\dim E^u=1$ or \item[-] $M = \TT^d$.
\end{itemize}
Then, $M=\TT^d$ and $f$ is isotopic to a linear Anosov
automorphism with stable dimension equal to $\dim E^{cs}$.
\end{teo}

In view of Franks-Manning theory \cite{FranksAnosov,Manning} one
can expect that this result also holds for nilmanifolds.

Notice that in general, obtaining a classification result for
partially hyperbolic diffeomorphisms with the coherent trapping
property should be at least as difficult as having a
classification result for Anosov diffeomorphisms since the latter
are indeed partially hyperbolic with the coherent trapping
property.

The fact that $f$ is dynamically coherent seems to be a strong
hypothesis, more in view of the robustness of the conclusion of
Theorem \ref{Teorema-SemiconjALineal}. However, we have not been
able to remove the hypothesis from our assumptions unless some
strong properties are verified.

Along this section we shall assume that $f\in \Diff^1(M)$ is a
partially hyperbolic diffeomorphism with the coherent trapping
property. Also, we shall call $cs=\dim E^{cs}$ and $u= d- \dim
E^{cs} = \dim E^u$.

\subsection{An expansive quotient of the dynamics}\label{SubSectionCOCIENTE}

We can define for each $x\in M$

$$ A_x = \bigcap_{n\geq 0} f^{n}(\overline{\cD^{cs}_{f^{-n}(x)}})$$

Some obvious properties satisfied by the sets $A_x$ are:

\begin{itemize}
\item[-] $f(A_x) = A_{f(x)}$ for every $x\in M$. \item[-] The set
$A_x$ is a decreasing union of topological balls (it is a cellular
set), so, compact and connected in particular.
\end{itemize}

We would like to prove that the sets $A_x$ constitute a partition
of $M$ and that they vary semicontinuously, so that we can
quotient the dynamics. For this, the following lemma is of great
use:

\begin{lema}\label{LemaRelacionEquivalencia}
For every $y\in \cF^{cs}(x)$, there exists $n_y$ such that
$f^{n_y}(\overline{\cD^{cs}_y}) \en \cD^{cs}_{f^{n_y}(x)}$. The
number $n_y$ varies semicontinuously on the point, that is, there
exists $U$ a small neighborhood of $y$ such that for every $z\in
U$ we have that $n_z\leq n_y$.
\end{lema}

\dem{} Consider in $\cF^{cs}(x)$ the sets

$$ E_n = \{ y\in \cF^{cs}(x) \ : \ f^{n}(\overline{\cD^{cs}_y}) \en \cD^{cs}_{f^{n}(x)} \}$$

Notice that there exists $\delta>0$ (independent of $n$) such that
if $y \in E_n$, then $B_\delta(y) \cap \cF^{cs}(x) \en E_n$. This
is given by continuity of $f$ and of the plaque family (using
compactness of $M$) and by the coherent trapping property.

The sets $E_n$ are thus clearly open and verify that $E_{n} \en
E_{n+1}$ (this is implied by the coherent trapping property).

Now, by the uniform estimate, it is not hard to show that
$\bigcup_{n\geq 0} E_n$ is closed, so, since it is not empty, it
must be the whole $\cF^{cs}(x)$ as claimed.

The fact that the numbers $n_y$ varies semicontinuously is a
consequence of the fact that $E_n$ is open ($n_y$ is the first
integer such that $y \in E_n$). \lqqd

\begin{cor}\label{CorRelacionEquivalencia}
For $x,y\in M$ we have that $A_x = A_y$ or $A_x \cap A_y =
\emptyset$. Moreover, the classes vary semicontinuously, that is,
given $x_n \in M$ such that $\lim x_n = x$:
$$\limsup A_{x_n} = \bigcap_{k>0} \overline{\bigcup_{n>k} A_{x_n}} \en A_x$$.
\end{cor}
\dem{} There exists $n_0$ fixed such that for every $x\in M$ and
$y\in f(\cD^{cs}_{f^{-1}(x)})$ we have that $f^{n_0}(\cD^{cs}_{y})
\en \cD^{cs}_{f^{n_0}(x)}$. This is proved first by showing that
$n_x$ exists for each $x\in M$ (using the Lemma
\ref{LemaRelacionEquivalencia} and compactness of
$f(\overline{\cD^{cs}_{f^{-1}(x)}})$) and then, since the numbers
$n_x$ vary semicontinuously, the uniform bound $n_0$ is found.

We know that for every $z$ such that $z \in A_x$ we have that $A_z
\in \cD^{cs}_x$: Indeed, since $z \in A_x$ we have that
$f^{-n_0}(z) \in \cD^{cs}_{f^{-n_0}(x)}$ and thus
$f^{n_0}(\cD^{cs}_{f^{-n_0}(z)}) \en \cD^{cs}_x$ as desired. In
fact, this shows that if $z \in A_x$ then $A_z \en A_x$. In
particular, by symmetry, we get that if $A_x \cap A_y \neq
\emptyset$, we can find a point $z$ in the intersection and we
have that $A_z = A_x$ and $A_z = A_y$ giving the desired
statement.

To prove semicontinuity: Consider $A_x$ and an $\eps$ neighborhood
$A_x(\eps)$ in $\cF^{cs}(x)$. Now, fix $m$ such that
$f^m(\cD^{cs}_{f^{-m}(x)}) \en A_x(\eps)$. Now, for $n$ large
enough, $x_n$ verifies that $d(f^{-m}(x_n),f^{-m}(x))$ is so small
that $f^m(\cD^{cs}_{f^{-m}(x_n)}) \en A_x(\eps)$ as wanted. The
semicontinuity for points which are not in the same center-stable
manifold follows from Lemma \ref{LemaHolonomia} bellow.

\lqqd

We get thus a continuous projection by considering the equivalence
relation $x\sim y \Leftrightarrow y \in A_x$.

$$\pi : M \to M /_{\sim}$$

We denote as $g:M /_{\sim} \to M /_{\sim}$ the map given by
$g([x]) = [f(x)]$ (that is $g\circ \pi = \pi \circ f$). Since
$\pi$ is continuous and surjective (in fact, it is cellular), it
is a semiconjugacy. Notice that since $f$ is a diffeomorphism, and
$g$ preserves the equivalence classes, one can show that $g$ must
be a homeomorphism of $M/_\sim$.

Notice that a priori, we have no knowledge of the topology of $M
/_{\sim}$ except that it is the image by a cellular map of a
manifold (see Section \ref{Section-Descomposiciones}), for
example, we do not know a priori if the dimension of $M /_{\sim}$
is finite. However, in view of Proposition
\ref{Proposition-CocienteCelularEsMetrico} we know that this
quotient is a metric space.

We will prove that it has finite topological dimension
dynamically after we prove Theorem \ref{TeoremaGesExpansivo}
(combined with \cite{ManheExpansiveTopological}).

We say that a homeomorphism has \emph{local product structure} if
there exists $\delta>0$ such that $d(x,y)<\delta$ implies that
$S_\eps (x)\cap U_\eps (y) \neq \emptyset$ (see Section
\ref{Seccion-Recurrence-And-Perturbation}).

Recall that a homeomorphism $h$ is \emph{expansive} (with
expansivity constant $\alpha$) if for every $x\in X$ we have that
$S_\alpha (x) \cap U_\alpha(x) = \{x\}$.

It is well known that for expansive homeomorphisms we have that
$\diam(h^n(S_\eps(x))) \to 0$ uniformly on $x$ for $\eps<\alpha$
(so this coincides with the usual definitions of stable and
unstable sets). This implies that $S_\eps(x) \en W^s(x)$ for an
expansive homeomorphism ($\eps < \alpha$).

\begin{teo}\label{TeoremaGesExpansivo}
The homeomorphism $g$ is expansive with local product structure.
Moreover, $\pi(\cF^{cs}(x))= W^s(\pi(x))$ and $\pi$ is injective
when restricted to the unstable manifold of any point.
\end{teo}

\dem The last two claims are direct from Lemma
\ref{LemaRelacionEquivalencia} and the definition of the
equivalence classes respectively.

We choose $\eps>0$ such that:

\begin{itemize}
\item[-] $x,y\in M$ and $x\notin f^{-1}(\cD^{cs}_{f(y)})$ then
there exists $n\geq 0$ such that
$d(f^n(x),\overline{\cD^{cs}_{f^n(y)}})>\eps$.
\end{itemize}

Now, let $x,y$ be two points such that $d(f^n(x),f^n(y))\leq \eps$
for every $n\in \ZZ$. From how we choose $\eps$, we have that
$f^{-k}(x)\in \cD^{cs}_{f^{-k}(y)}$ for every $k\geq 0$ so, $x\in
A_y$ as desired.

Since $\pi(\cF^{cs}(x)) = W^s(\pi(x))$, and since $d(x,y)<\delta$
implies that $\cD^{cs}_x \trans W^{uu}_{loc}(y) \neq\emptyset$ we
get that for every two close points, there is a non trivial
intersection between the local stable and unstable sets (here, we
are using the upper semicontinuity of the sets $A_x$ which imply
that if there are two points $\tilde x, \tilde y$ in $M /_{\sim}$
which are near, there are points in $\pi^{-1}(\tilde x)$ and
$\pi^{-1}(\tilde y)$ which are near).

\lqqd

Consider two points $x,y$ such that $y \in W^{uu}(x)$. We denote
$\Pi^{uu}_{x,y}: \cD \en \cD^{cs}_x \to \cD^{cs}_y$ as the
unstable holonomy from a subset of $\cD^{cs}_x $ into a subset of
$\cD^{cs}_y$. An important useful property is the following:

\begin{lema}\label{LemaHolonomia} We have that $\Pi^{uu}_{x,y}(A_x) = A_y$.
\end{lema}

\dem It is enough to show (by the symmetry of the problem) that
$\Pi^{uu}(A_x)\en A_y$. For $n$ large enough we have that
$f^{-n}(\Pi^{uu}(A_x))$ is very close to a compact subset of
$\cD^{cs}_{f^{-n}(x)}$ and thus, by continuity of $\cD^{cs}$ we
have that $f^{-n}(\Pi^{uu}(A_x)) \en \cD^{cs}_{f^{-n}(y)}$ which
concludes. \lqqd

\subsubsection{Some remarks on the topology of the quotient}

We shall cite some results from \cite{Daverman} which help to
understand the topology of $M /_{\sim}$. We refer to the reader to
that book for much more information and precise definitions.

Before, we remark that Ma\~ne proved that a compact metric space
admitting an expansive homeomorphism must have finite topological
dimension (\cite{ManheExpansiveTopological}).

Corollary IV.20.3A of \cite{Daverman} implies that, since $M
/_{\sim}$ is finite dimensional, we have that it is a locally
compact ANR (i.e. absolute neighborhood retract). In particular,
we get that $\dim (M/_{\sim}) \leq \dim M$ (see Theorem III.17.7).
Then, by using Proposition VI.26.1 (or Corollary VI.26.1A) we get
that $M/_{\sim}$ is a $d-$dimensional homology manifold (since it
is an ANR, it is a \emph{generalized manifold}). More properties
of these spaces can be found in section VI.26 of \cite{Daverman}.

Also, in the cited book, one can find a statement of Moore's
theorem (see section IV.25 of \cite{Daverman}) which states that a
cellular decomposition of a surface is approximated by
homeomorphisms (this means that the continuous projection is
approximated by homeomorphisms in the $C^0$-topology). In
particular, in our case, if $\dim E^{cs} = 2$, we get that $M
/_{\sim}$ is a manifold (see also Theorem VI.31.5 of
\cite{Daverman} and its Corolaries).

Some other results are available, in particular, we notice
Edward's cell-like decomposition theorem which asserts that if
$\sim$ is a cellular decomposition of a $d$ dimensional manifold
($d\geq 5$) such that $M /_\sim$ has finite topological dimension
and such that it has the \emph{disjoint disk property} (see
chapter IV.24 of \cite{Daverman}) then the quotient map is
approximated by homeomorphisms. A similar result exists for
dimension $3$ which is even more technical. Notice than in our
case, since we have the decomposition of the center-stable
manifold, we can play with the dimensions in order not to be never
in dimension $4$ by choosing to work with the decomposition on the
center stables or the whole manifold.

Also, we remark that it is known that when multiplying a
decomposition by $\RR^2$ we always get the \emph{disjoint disc
property} and in all known decompositions, after multiplying by
$\RR$ we get a decomposition approximated by homeomorphisms (see
section V.26 of \cite{Daverman}) so in theory, it should be true
that always our space $M  /_\sim$ is a manifold homeomorphic to
$M$. We show this in the case $M=\TT^d$ (the proof should be
adaptable for infranilmanifolds).

\subsection{Transitivity of the expansive homeomorphism}\label{SubSectionTRANSITIVIDADDEG}

In general, it is not yet known if an Anosov diffeomorphism must
be transitive. So, since Anosov diffeomorphisms enter in our
hypothesis, there is no hope of knowing if $f$ or $g$ will be
transitive without solving this long-standing conjecture. We shall
then work with similar hypothesis to the well known facts for
Anosov diffeomorphisms, showing that those hypothesis that we know
guaranty that Anosov diffeomorphisms are transitive imply
transitivity of $g$.

In particular, we shall prove in this section the following
Theorem which implies Theorem \ref{Teorema-SemiconjALineal}:

\begin{teo}\label{TeoremaGtransitivo} The following properties hold:
\begin{itemize}
\item[\emph{(T1)}] If for every $x,y \in M$ we have that
$\cF^{uu}(x) \cap \cD^{cs}(y)\neq \emptyset$, then $g$ is
transitive. \item[\emph{(T2)}] If $\dim E^u =1$, then $g$ is
transitive. Moreover, $M=\TT^d$. \item[\emph{(T3)}] If $M=\TT^d$,
then $g$ is transitive. Moreover, $f$ is homotopic to a linear
Anosov diffeomorphism $A$ the topological space $M /_\sim$ is
homeomorphic to $\TT^d$ and $g$ is conjugated to $A$.
\end{itemize}
\end{teo}

Notice that (T1) is trivial, (T2) can be compared to
Franks-Newhouse theory (\cite{FranksAnosov,NewhouseAnosov}) and
(T3) to Franks-Manning theory (\cite{FranksTori,Manning} see also
\cite{KH} chapter 18.6 which motivated the proof here presented).
It is natural to expect that property (T3) should hold if we
consider $M$ an infranilmanifold.

It is important to notice also that it is natural to extend the
conjecture about transitivity of Anosov diffeomorphisms to
expansive homeomorphisms with local product structure (at least in
manifolds). See the results in \cite{Vieitez, ABP, Hiraide}.

\subsubsection{Proof of (T2)} We shall follow the argument of \cite{NewhouseAnosov}.

From how we defined $g$, we get that for every $y\in M$ we have
that $\pi|_{W^u_{loc}(y)}$ is a homeomorphism over its image and
thus, we get that every point in $M /_\sim$ has a one dimensional
immersed copy of $\RR$ as unstable set.

Also, we have that $g: M/_\sim \to M/_\sim$ is expansive with
local product structure. This implies that there is a spectral
decomposition for $g$:

\begin{lema}\label{LemaGDescEspectral} The homeomorphism $g$ has a spectral decomposition.
\end{lema}

\dem The proof is exactly as the one for Anosov or Axiom A
diffeomorphisms.

It is not hard to show that $Per(g)$ is dense in $\Omega(g)$ with
essentially the same proof as in the Anosov case. Let $x$ be a
nonwandering point of $g$, so, every neighborhood of $\pi^{-1}(x)$
has points which return to the neighborhood in arbitrarily large
backward iterates. The fact that center stable leaves are
invariant and unstable manifolds expand by iterating backwards,
gives the existence of a fixed center stable leaf with a point
returning near itself. Since the center stable disks are trapped,
we obtain a fixed fiber for some iterate, this gives a periodic
point for $g$ which is arbitrarily close to $x$.

The rest of the spectral decomposition, is done by defining
homoclinic classes and that needs no more that the local product
structure of uniform size (see \cite{Newhouse-homoclinic}). \lqqd

By Conley's theory (see Remark \ref{Remark-QuasiAttractor}), we
get a repeller $\Lambda$ for $g$ which will be saturated by stable
sets.

We shall show that $\Lambda = M/_\sim$ which concludes.

To do this, it is enough to show that for every $y\in \Lambda$, we
have that $y$ is accumulated by the intersections of both
connected components of $W^u(y)\setminus \{y\}$ with $\Lambda$.

We can assume that $W^u$ is orientable and $g$ preserves
orientation of $W^u$ (otherwise, we take a double cover and $g^2$,
transitivity at this level is even more general than if we do not
take the cover nor the iterate). So, for every $y \in \Lambda$ we
denote $W^u_+(y)$ and $W^u_-(y)$ the connected components of
$W^u(y)\setminus \{y\}$ depending on the orientation.

We define the set

$$ A^+= \{ y\in \Lambda \ : \ \ W^u_+(y) \cap \Lambda \neq \emptyset \} $$

\noindent which is an invariant set (we define $A^-$ similarly).
It is enough to show that $A^+=\Lambda$: Indeed, this implies, by
compactness that every point intersects $\Lambda$ in a bounded
length of $W^u_+$. This is enough since being invariant, the
length must be zero.

\begin{lema}\label{LemaNonPeriodicInAmas}
Any point which is not periodic by $g$ belongs to $A^+$. In fact,
there are at most finitely many points not in $A^+$.
\end{lema}

\dem The past orbit of every point contains an accumulation point
and they are pairwise not in the same local stable set, thus,
there is one point in the past orbit such that both components of
its unstable set intersect the stable set of the other, and thus
$\Lambda$, invariance concludes.

The fact that there exists $N>0$ such that any set with cardinal
larger than $N$ has $3$ points in a local product structure box,
implies that if a point $x$ does not belong to $A^+$ then its
orbit $\cO(x)$ must have cardinal smaller than $N$. This gives
that there are at most finitely many points outside $A^+$ (which
must be periodic). \lqqd

To prove that periodic points are in $A^+$ too, we assume that it
is not the case and consider $p \in \Lambda$ such that $p\notin
A^+$. We have that $W^{s}(p)\backslash \{p\}$ is connected
\footnote{ We are assuming that $\dim E^{cs} \geq 2$, since if we
remove $\pi^{-1}(p)$ to $\cF^{cs}(\pi^{-1}(p))$ it remains
connected, the claim follows.}.

This implies that if $\phi: W^{s}(p)\backslash \{p\} \to \Lambda$
is the function which sends every point $y\in W^s(p)\backslash
\{p\}$ to the first intersection of $W^u_+(y)$ with $\Lambda$ (the
first point of intersection exists since otherwise we would get
that $p \in A^+$), then we have that the image of
$W^s(p)\backslash \{p\}$ is a unique stable set, say of a point
$z$.

Now, we must show that in fact, we have that $W^u(p)$ must
intersect $W^s(z)$ which will be a contradiction and conclude. So,
consider in $\pi^{-1}(W^s_{loc}(p)) \en \cD^{cs}_{\tilde p}$
(where $\pi(\tilde p)=p$)  a (small) sphere $\Sigma$ around
$\pi^{-1}(p)$ . That is, we assume that $\pi(\Sigma)\en
W^s_{loc}(p)$ (which we can since $\pi$ is a cellular map).

Now, we consider $y\in \Sigma$ and $I \en W^u_+(y)$ the interval
of the unstable manifold of $y$ from $y$ to the only point in
$W^u(y) \cap \pi^{-1}(\phi(\pi(y)))$. This interval can be
parametrized in $[0,1]$. We shall call $y_t$ to the point
corresponding to $t\in [0,1]$.

We consider the set of points $s\in [0,1]$ such that $W^u_+ (p)
\cap W^s_{loc}(\pi(y_t)) \neq \emptyset$ which is open by the
local product structure. We must show that the supremum of this
set, say $t_0$ belongs to the set and we shall conclude.

Let $\Upsilon$ be the set homeomorphic to $\Sigma \times [0,1]$
given by the map $Y: (x,t) \mapsto I_x \cap \pi^{-1}(\pi(I_x) \cap
W^s(\pi(y_t)))$ where $I_x$ is the interval in the unstable
manifold of $x$ such that $\pi(I_x) \en W^u_+(\pi(x))$ and
connects $\pi(x)$ with $\phi(\pi(x))$.

We get that $Y(\Sigma \times \{t_0\})$ is homeomorphic to a
sphere, and thus it separates $\cF^{cs}(y_{t_0})$ in two connected
components, one of which is bounded (and thus compact). The image
$K$ of this compact component by $\pi$ is also compact and
contained in $W^s(\pi(y_{t_0}))$ and thus we get that there is
local product structure well defined around $K$. This implies that
$W^u_+(p) \cap W^s(\pi(y_{t_0})) \neq \emptyset$ and we have shown
that $t_0$ is in the set considered above. This shows that
$A^+=\Lambda$ which as we have already mentioned gives
transitivity of $f$.

Now, following the same proof as in (for example) the appendix of
\cite{ABP}, we get that $M=\TT^d$. \lqqd

\subsubsection{Proof of (T3)} We shall follow the proof given in \cite{KH} chapter 18.6.

Before we start with the proof, we shall recall Theorem 18.5.5 of
\cite{KH} (the statement is modified in order to fit our needs,
notice that for an expansive homeomorphism with local product
structure, we have the shadowing property, and thus, specification
in each basic piece):

\begin{prop}[Theorem 18.5.5 of \cite{KH}]\label{PropConteoDePeriodicos}
Let $X$ a compact metric space and $g:X \to X$ an expansive
homeomorphism with local product structure. Then, there exists
$h,c_1,c_2>0$ such that for $n\in \NN$ we have:
$$ c_1 e^{n h} \leq P_n(g) \leq c_2 e^{n h} $$
\noindent where $P_n(g)$ is the number of fixed points of $g^n$.
\end{prop}

We shall use several time the very well know Lefschetz formula
which relates the homotopy type of a continuous function, with the
index of its fixed points (see \cite{Franks-Homology} Chapter 5).

\begin{defi}
Let $V\en \RR^k$ be an open set, and $F: V\en \RR^k \to \RR^k$ a
continuous map such that $\Gamma\en V$ the set of fixed points of
$F$ is a compact set, then, $I_\Gamma(F) \in \ZZ$ (the index of
$F$) is defined to be the image by $(id-F)_\ast: H_k(V, V-\Gamma)
\to H_k(\RR^k, \RR^k-\{0\})$ of $u_\Gamma$ where $u_\Gamma$ is the
image of $1$ under the composite $H_k(\RR^k , \RR^k - D) \to
H_k(\RR^k, \RR^k- \Gamma) \cong H_k(V, V-\Gamma)$ where $D$ is a
disk containing $\Gamma$. \finobs
\end{defi}

\begin{obs}
In general, if we have a map from a manifold, we can include the
manifold in $\RR^k$ and extend the map in order to be in the
hypothesis of the definition. The value of $I_\Gamma(F)$ does not
depend on how we embed the manifold in $\RR^k$.\finobs
\end{obs}

For hyperbolic fixed points, it is very easy how to compute the
index, it is exactly the sign of $\det(Id - D_pf)$. Since the
definition is topological, any time we have a set which behaves
locally as a hyperbolic fixed point, it is not hard to see that
the index is the same.

Lefshetz fixed point formula for the torus can be stated as
follows:

\begin{teo}[Lefshetz fixed point formula (\cite{Franks-Homology} p.34-38)]\label{teoLefshetz}
Let $h:\TT^d \to \TT^d$ be an homeomorphism, so, the sum of the
Lefshetz index along a covering of $\Fix(h)$ by sets homeomorphic
to balls equals $det(Id- h_\ast)$ where $h_\ast: H_1(\TT^d,\ZZ)
\to H_1(\TT^d, \ZZ)$ is the action of $h$ in homology.
\end{teo}

The first thing we must show, is that the \emph{linear} part of
$f$, that is, the action $A= f_\ast: H_1(\TT^d, \ZZ) \to
H_1(\TT^d, \ZZ) \in SL(d,\ZZ)$ is a hyperbolic matrix.

\begin{lema}\label{LemaParteLinealHiperbolica} The matrix $A$ is hyperbolic.
\end{lema}

\dem We can assume (maybe after considering a double covering and
$f^2$) that $E^{cs}$ and $E^u$ are orientable and its orientations
preserved by $Df$. So, it is not hard to show that for every fixed
point $p$ of $g^n$, the index of $\pi^{-1}(p)$ for $f$ is of
modulus one and always of the same sign.

So, we know from the Lefshetz formula that

$$|\det(Id - A^n)| = \sum_{g^n(p)=p} |I_{\pi^{-1}(p)}(f)| = \# \Fix(g^n).$$

Proposition \ref{PropConteoDePeriodicos} and an easy estimate on
the growth of $|\det(Id -A^n)| = \prod_{i=1}^d |1-\lambda_i^n|$
where $\{\lambda_1, \ldots, \lambda_d\}$ are the eigenvalues of
$A$ gives that $A$ cannot have eigenvalues of modulus $1$ and thus
$A$ must be hyperbolic (see the argument in Lemma 18.6.2 of
\cite{KH}).

\lqqd

Proposition \ref{PropExisteSemiconjugacion} gives the existence of
a semiconjugacy $h:\TT^d \to \TT^d$ isotopic to the identity such
that $h \circ f = A \circ h$. Its lift $H: \RR^d \to \RR^d$ is
given by shadowing, in particular, the iterations of the set
$H^{-1}(x)$ remain of bounded diameter.

\begin{lema}\label{LemaSemiconjIntermedia}
We have that $g$ factors as an intermediate semiconjugacy. More
precisely, there exists $\tilde h : \TT^d /_{ \sim} \to \TT^d$
continuous and surjective such that $\tilde h \circ \pi = h$.
\end{lema}

\dem It is enough to show that for every $x\in \TT^d/_\sim$ there
exists $y\in \TT^d$ such that $\pi^{-1}(x) \subset h^{-1}(y)$.

For this, notice that any lifting  of $\pi^{-1}(x)$ (that is, a
connected component of the preimage under the covering map) to the
universal covering $\RR^d$ verifies that its iterates remain of
bounded size. This concludes by the remark above on $H$.

\lqqd

Now, we shall prove that if $\tilde f : \RR^d \to \RR^d$ is any
lift of $f$, then there is exactly one fixed fiber of $\pi$ for
$\tilde f$.

\begin{lema}\label{LemaUnicaFibraArriba}
Let $\tilde f^n$ be any lift of $f^n$ to $\RR^d$. So, there is
exactly one fixed fiber of $\pi$.
\end{lema}

\dem Since $\tilde f^n$ is homotopic to $A^n$ which has exactly
one fixed point and each fixed fiber of $\pi$ contributes the same
amount to the index of $\tilde f^n$ it must have exactly one fixed
fiber. \lqqd

This allows us to show that $g$ is transitive:

\begin{prop}\label{PropGtransitivo} The homeomorphism $g$ is transitive.
\end{prop}

\dem First, we show that there exists a basic piece of $g$ which
projects by $\tilde h$ to the whole $\TT^d$.

This is easy since otherwise, there would be a periodic point $q$
in $\TT^d \setminus \tilde h (\Omega(g))$ but clearly, the
$g-$orbit of $\tilde h^{-1}(q)$ must contain non-wandering points
(it is compact and invariant).

This concludes, since considering a point $y$  with dense
$A$-orbit and a point in $\Omega(g) \cap \tilde h^{-1}(y)$ we get
the desired basic piece.

Now, let $\Lambda$ be the basic piece of $g$ such that $\tilde h
(\Lambda) = \TT^d$. Assume that there exists $\tilde \Lambda \neq
\Lambda$ a different basic piece and $z$ a periodic point of
$\tilde \Lambda$, naturally, we get that $\tilde h^{-1}(\tilde h
(z))$ contains also a periodic point $z'$ in $\Lambda$. By
considering an iterate, we can assume that $z$ and $z'$ are fixed
by $g$.

So, we get that it is possible to lift $h^{-1}(\tilde h(z))$ and
chose a lift of $f^k$ which fixes $\pi^{-1}(z)$ and $\pi^{-1}(z')$
contradicting the previous lemma.

\lqqd

With this in hand, we will continue to prove that the fibers of
$h$ coincide with those of $\pi$ proving that $g$ is conjugated to
$A$ (in particular, $\TT^d /_\sim \cong \TT^d$).

First, we show a global product structure for the lift of $f$.
Notice that when we lift $f$ to $\RR^d$, we can also lift its
center-stable and unstable foliation. It is clear that both
foliations in $\RR^d$ are composed by leaves homeomorphic to
$\RR^{cs}$ and $\RR^u$ respectively (the unstable one is direct,
the other is an increasing union of balls, so the same holds).

\begin{lema}\label{LemaGlobalProductStructure}
Given $x,y\in \RR^d$, the center stable leaf of $x$ intersects the
unstable leaf of $y$ in exactly one point.
\end{lema}

\dem The fact that they intersect in at most one point is given by
the fact that otherwise, we could find a horseshoe for the lift,
and thus many periodic points contradicting  Lemma
\ref{LemaUnicaFibraArriba} (for more details, see Lemma 18.6.7 in
\cite{KH}).

The proof that any two points have intersecting manifolds, is
quite classical, and almost topological once we know that both
foliations project into minimal foliations (see also Lemma 18.6.7
of \cite{KH}). \lqqd

Now, we can conclude with the proof of Part (T3) of Theorem
\ref{TeoremaGtransitivo}.

To do this, notice that the map $H$ conjugating $\tilde f$ with
$A$ is proper, so the preimage of compact sets is compact. Now,
assume that $A_x, A_y$ are lifts of fibers of $\pi$ such that $H
(A_x) = H (A_y)$ we shall show they coincide.

Consider $K$ such that if two points have an iterate at distance
bigger than $K$ then their image by $H$ is distinct.

We fix $x_0 \in A_x$ and consider a box $D_K^n$ of $\tilde f^n
(x_0)$ consisting of the points $z$ of $\RR^d$ such that $\cF^u(z)
\cap \cF^{cs}_K(x_0) \neq \emptyset$ and $\cF^{cs}(z) \cap
\cF^u_K(x_0) \neq\emptyset$.

It is not hard to show using Lemma
\ref{LemaGlobalProductStructure} that there exists $\tilde K$
independent of $n$ such that every pair of points in $D_K^n$ in
the same unstable leaf of $\cF^u$ have distance along $\cF^u$
smaller than $\tilde K$ (this is a compactness argument). An
analogous property holds for $\cF^{cs}$.

This implies that if $\tilde f^n (A_y) \en D_K^n$ for every $n\in
\ZZ$ then $A_y$ and $A_x$ must be contained in the same leaf of
$\cF^{cs}$. In fact we get that $\tilde{f}^{-n}(A_y) \en
\cF^{cs}_K(\tilde{f}^{-n}(x_0))$  for every $n\geq 0$ and so we
conclude that $A_x=A_y$ using Lemma
\ref{LemaRelacionEquivalencia}.

\lqqd

\subsection{Some manifolds which do not admit this kind of diffeomorphisms}\label{SubSubsectionVARIEDADESQADMITEN}

The arguments used in the previous section also allow to show that
certain manifolds (and even some isotopy classes in some
manifolds) do not admit partially hyperbolic diffeomorphisms
satisfying the coherent trapping property.

To do this, we recall that for general manifolds $M^d$, and a
homeomorphism $h: M\to M$, the Lefschetz number of $h$, which we
denote as $L(h)$ is calculated as $\sum_{i=0}^d trace(h_{\ast,i})$
where $h_{\ast,i}: H_i(M,\QQ) \to H_i(M,\QQ)$ is the induced map
on (rational) homology\footnote{This is just to avoid torsion
elements. Otherwise, one can define the trace with $\ZZ$
coeficients after making a quotient by the torsion.}. We also have
that the sum of index the sum of the Lefshetz index along a
covering of $\Fix(h)$ by sets homeomorphic to balls equals $L(h)$.

A similar argument to the one used in the previous section yields
the following result (see also \cite{NoAnosovEnAlgunasVariedades}
for the analog result for Anosov diffeomorphisms\footnote{Although
this result even refines slightly his result even for Anosov
diffeomorphisms.})

\begin{teo}\label{TeoremaAlgunasVariedadesNoAdmiten}
Let $f$ be a partially hyperbolic diffeomorphism of $M$ with the
coherent trapping property, then, the action $f_\ast :
H_{\ast}(M,\QQ) \to H_{\ast}(M,\QQ)$ is strongly partially
hyperbolic (it has both eigenvalues of modulus $>1$ and $<1$).
\end{teo}

As a consequence, several manifolds cannot admit this kind of
diffeomorphisms (notably $S^d$ and products of spheres of
different dimensions such as $S^1 \times S^2 \times S^3$) and
also, for example, there cannot be diffeomorphisms like this
acting as the identity on homology. This leads to a natural
question: Is every partially hyperbolic diffeomorphism with the
coherent trapping property homotopic to an Anosov diffeomorphism?.

\dem The proof is very similar to the one given in the previous
section, so we shall omit some details.

First, we get by counting the fixed points of $g^n$ (the expansive
quotient of $f^n$ to $M /_\sim$) and we get an exponential growth.

Now, if there are no eigenvalues of modulus greater than one for
$f_\ast$, the trace of the map cannot grow exponentially. The same
argument applies to $f^{-1}$ so we conclude.

\lqqd


   \appendix

\chapter{Perturbation of cocycles}\label{Appendix-PerturbacionCociclos}

\section{Definitions and statement of results}

Before we proceed with the statement and proof of Theorem
\ref{BGV} we shall give some definitions taken from \cite{BGV} and
others which we shall adapt to fit our needs.

Recall that $\mathcal{A} = (\Sigma, f, E, A)$ is a \emph{large
period linear cocycle}\footnote{Sometimes, we shall abuse notation
and call it just \emph{cocycle}.} of dimension $d$ bounded by $K$
over an infinite set $\Sigma$ iff:

\begin{itemize}

\item[-] $f:\Sigma \to \Sigma$ is a bijection such that all points
in $\Sigma$ are periodic and such that given $n>0$ there are only
finitely many with period less than $n$. \item[-]  $E$ is a vector
bundle over $\Sigma$, that is, there is $p: E \to \Sigma$ such
that $E_x=p^{-1}(x)$ is a vector space of dimension $d$ endowed
with an euclidian metric $\langle,\rangle_x$. \item[-] $A:x\in
\Sigma \mapsto A_x\in GL(E_x, E_{f(x)})$ is such that $\|A_x\|\leq
K$ and $\|A_x^{-1}\|\leq K$.

\end{itemize}

In general, we shall denote $A^\ell_x = A_{f^{\ell-1}(x)}  \ldots
A_x$ where juxtaposition denotes the usual composition of linear
transformations.


For every $x\in \Sigma$ we denote by $\pi(x)$ its period and
$M_x^A = A^{\pi(x)}_x$ which is a linear map in $GL(E_x,E_x)$
(which allows to study eigenvalues and eigenvectors).

For $1\leq j \leq d$

$$\sigma^{j}(x,\mathcal A) = \frac{log|\lambda_j|}{\pi(x)}$$

Where $\lambda_1,\ldots \lambda_d$ are the eigenvalues of $M_x^A$
in increasing order of modulus. As usual, we call
$\sigma^j(x,\cA)$ the $j$-th Lyapunov exponent of $\cA$ at $x$.

Given an $f$-invariant subset $\Sigma'\en \Sigma$, we can always
restrict the cocycle to the invariant set defining the cocycle
$\cA|_{\Sigma'} = (f|_{\Sigma'}, \Sigma', E|_{\Sigma'},
A|_{\Sigma'})$.

We shall say that a subbundle $F\en E$ is \emph{invariant} if
$\forall x \in \Sigma$ we have $A_x(F_x) =F_{f(x)}$. When there is
an invariant subbundle, we can write the cocycle in coordinates
$F\oplus F^{\perp}$ (notice that $F^\perp$ may not be invariant)
as

\[    \left(
       \begin{array}{cc}
         \cA_F & C_F \\
         0 & \cA|F \\
       \end{array}
     \right)\]

Where $C_F$ is uniformly bounded. This induces two new cocycles:
$\cA_F =(\Sigma, f, F, A|_{F})$ on $F$ (where $A|_F$ is the
restriction of $A$ to $F$) and $\cA|F = (\Sigma, f, E|F, A|F)$ on
$E|F \simeq F^\perp$ where $(A|F)_x \in GL((F_x)^\perp,
(F_{f(x)})^\perp)$ is given by $p^2_{f(x)} \circ A_x$ where
$p^2_x$ is the projection map from $E$ to $F^\perp$. Notice that
changing only $\cA_F$ affects only the eigenvalues associated to
$F$ and changing only $\cA|F$ affects only the rest of the
eigenvalues, recall Remark \ref{Remark-Deapedazos}. See section
4.1 of \cite{BDP} for more discussions on this decomposition.

As in Section \ref{Section-EstructurasInvariantes}, if $\cA$ has
two invariant subbundles $F$ and $G$, we shall say that $F$ is
$\ell$-\emph{dominated} by $G$ (and denote it as $F\prec_{\ell}
G$) on an invariant subset $\Sigma'\en \Sigma$ if for every $x\in
\Sigma'$ and for every pair of vectors $v\in F_x \backslash
\{0\}$, $w\in G_x \backslash \{0\}$ one has

\[ \frac{\|A^\ell_x(v)\|}{\|v\|} \leq \frac{1}{2} \frac{\|A^\ell_x(w)\|}{\|w\|} \ . \]

We shall denote $F\prec G$ when there exists $\ell>0$ such that
$F\prec_{\ell} G$.

If there exists complementary invariant subbundles $E=F\oplus G$
such that $F \prec G$ on a subset $\Sigma'\en \Sigma$, we shall
say that $\cA$ admits a \emph{dominated splitting} on $\Sigma'$.

As in \cite{BGV}, we shall say that $\cA$ is \emph{strictly
without domination} if it is satisfied that whenever $\cA$ admits
a dominated splitting in a set $\Sigma'$ it is satisfied that
$\Sigma'$ is finite.

Let $\Gamma_x$ be the set of cocycles over the orbit of $x$ with
the distance

$$  d(\cA_x, \cB_x) = {\displaystyle \sup_{0\leq i < \pi(x), v \in E\backslash \{0\}}} \left\{ \frac{\|(A_{f^i(x)}-B_{f^i(x)})v\|}{\|v\|} , \frac{\|(A_{f^i(x)}^{-1}-B_{f^i(x)}^{-1})v\|}{\|v\|} \right\} $$

 and let $\Gamma_{\Sigma}$ (or $\Gamma_{\Sigma,0}$) the set of bounded large period linear cocycles over $\Sigma$. Given $\cA \in \Gamma_\Sigma$ we denote as $\cA_x \in \Gamma_x$ to the cocycle $\{A_x, \ldots, A_{f^{\pi(x)-1}(x)} \}$.

We say that the cocycle $\cA_x$ has \emph{strong stable manifold
of dimension} $i$ if  $\sigma^i(x,\cA_x) < \min
\{0,\sigma^{i+1}(x, \cA_x) \}$.

For $0\leq i \leq d$, let

$$\Gamma_{\Sigma,i} = \{ \cA \in \Gamma_{\Sigma} \ : \ \forall x \in \Sigma \ ; \ \cA_x \text{ has  strong  stable  manifold  of  dimension } i \}$$

\smallskip

Following \cite{BGV} we say that $\mathcal B$ is a
\emph{perturbation} of $\mathcal A$ (denoted by $\mathcal B \sim
\mathcal A$) if for every $\eps>0$ the set of points $x\in \Sigma$
such that $\cB_x$ is not $\eps-$close to the cocycle $\cA_x$ is
finite.

Similarly, we say that $\mathcal B$ is a \emph{path perturbation}
of $\mathcal A$ if for every $\eps>0$ one has that the set of
points $x \in \Sigma$ such that $\cB_x$ is not a perturbation of
$\mathcal A_x$ along a path of diameter $\leq \eps$ is finite.
That is, there is a path $\tilde \gamma: [0,1] \to
\Gamma_{\Sigma}$ such that $\tilde \gamma(0)=\cA$ and $\tilde
\gamma(1) = \cB$ such that $\tilde \gamma_x:[0,1]\to \Gamma_x$ are
continuous paths and given $\eps>0$ the set of $x$ such that
$\tilde \gamma_x([0,1])$ has diameter $\geq \eps$ is finite.

In general, we shall be concerned with path perturbations which
preserve the dimension of the strong stable manifold, so, we shall
say that $\mathcal B$ is a \emph{path perturbation of index $i$}
of a cocycle $\mathcal A \in \Gamma_{\Sigma,i}$ iff: $\cB$ is a
path perturbation of $\cA$ and the whole path is contained in
$\Gamma_{\Sigma,i}$.  This induces a relation in
$\Gamma_{\Sigma,i}$ which we shall denote as $\sim^{\ast}_i$.

We have that $\sim$ and $\sim^\ast_i$ are equivalence relations in
$\Gamma_\Sigma$ and $\Gamma_{\Sigma,i}$ respectively, and clearly
$\sim^\ast_i$ is contained in $\sim$.

\begin{obs}\label{deapedazos} Notice that if there is an invariant subbundle $F\en E$ for a cocycle $\cA$. Then, a perturbation (resp. path perturbation) of $\cA_F$ or $\cA|F$ can be completed to an perturbation (resp. path perturbation) of $\cA$ which does not alter the eigenvalues associated to $\cA|F$ or $\cA_F$. Any of this perturbations also preserves the invariance of $F$, however, one can not control the effect on other invariant subbundles. See section 4.1 of \cite{BDP}. \finobs
\end{obs}

The \emph{Lyapunov diameter} of the cocycle $\mathcal A$ is
defined as

$$\delta(\mathcal A)={\displaystyle \liminf_{\pi(x)\to \infty}} \  [\sigma^d(x,\mathcal A) - \sigma^1(x,\mathcal A)] .$$

If $\mathcal A \in \Gamma_{\Sigma,i}$, we define
$\delta_{min}(\mathcal A) = \inf_{\mathcal B \sim \mathcal A}
\{\delta(\mathcal B)\}$. Similarly, we define
$\delta_{min}^{\ast,i}(\mathcal A)= \inf_{\mathcal B \sim^\ast_i
\mathcal A} \{\delta(\mathcal B)\}$. Notice that
$\delta_{min}^{\ast,i} (\mathcal A) \geq \delta_{min}(\cA)$ and a
priori it could be strictly bigger.

\begin{obs}\label{dimensionextremal}
For any cocycle $\cA$, it is easy to see that
$\delta_{min}^{\ast,0}(\cA) = \delta_{min}(\cA)$. It sufficies to
consider the path $(1-t)\cA + t\cB$ where $\cB$ is a perturbation
of $\cA$ having the same determinant over any periodic orbit and
verifying $\delta(\cB)=\delta_{min}(\cA)$ (see Lemma 4.3 of
\cite{BGV} where it is shown that such a $\cB$ exists).
 \finobs
\end{obs}

The following easy Lemma relates the definitions we have just
introduced. Its proof is contained in \cite{BDP} and \cite{BGV}
except from property (f). We include quick proofs for
completeness.

\begin{lema}\label{BGVenunciado} Let $\cA=(\Sigma,f,E,A)$ and $\cB=(\Sigma,f,E,B)$ be two large period cocycles of dimension $d$ and bounded by $K$.
\begin{itemize}
 \item[(a)] If $\cA$ is strictly without domination and $\cB\sim \cA$ (in particular if $\cB \sim_i^\ast \cA$) then $\cB$ is also strictly without domination.
 \item[(b)] For every $\ell>0$ there exists $\nu>0$ such that if $F$ and $G$ are two invariant subbundles and $F\prec_\ell G$ on $\Sigma'$, then $\delta((\cA_{F\oplus G})|_{\Sigma'}) > \nu$.

 \item[(c)] If $\delta_{min}(\cA)=0$ there exists an infinite subset $\Sigma'\en \Sigma$ such that $\cA|_{\Sigma'}$ is strictly without domination. Conversely, if $\cA$ is strictly without domination, $\delta_{min}(\cA)=0$.
 \item[(d)] If $\cA\sim \cB$, then
 $$\left| \sum_{j=1}^d \sigma^j(x,\cA) - \sum_{j=1}^d \sigma^j(x,\cA) \right| \to 0  \quad as \quad \pi(x)\to \infty$$
  \item[(e)] Let $F,G$ and $H$ be invariant subbundles of $E$. So, $F\prec G\oplus H$ if  $F\prec G$ and $F|G\prec H|G$. Also, if $F\prec G$ and $G\prec H$ then $F\prec G\oplus H$ and $F\oplus G \prec H$.
  \item[(f)] Let $F, G$ and $H$ be invariant subbundles of $E$. So, $F\prec G\oplus H$ implies that $F|G \prec H|G$.
\end{itemize}
\end{lema}

\dem{\!\!} Part $(a)$ is Corollary 2.15 of \cite{BGV}. It uses the
quite standard fact (see Lemma 2.14 of \cite{BGV}) which asserts
that if a cocycle admits certain dominated splitting in a subset
$\Sigma'$, then, there exists $\eps$ such that every
$\eps$-perturbation of the cocycle remains with dominated
splitting in that set.

Assume that $\cB \sim \cA$ admits dominated splitting in a set
$\Sigma'$, then, the previous argument implies that, modulo
removing some finite subset of $\Sigma'$, the cocycle $\cA$ also
admits dominated splitting. This implies that $\Sigma'$ is finite
(otherwise, $\cA$ would admit a dominated splitting in an infinite
set).

Part $(b)$ follows directly from the definition of dominated
splitting ($\nu$ is going to depend on $K$, the dimensions of $F$
and $G$ and $\ell$).

Part $(c)$ is also standard. The existence of a dominated
splitting implies directly the separation of the Lyapunov
exponents in each of the invariant subbundles (see part $(b)$).
That is, given a dominated splitting over a set $\Sigma'$, we get
$\eps>0$ such that for every $x\in \Sigma'$ we have that

  $$\sigma^d(x,\mathcal A) - \sigma^1(x,\mathcal A) >\eps$$

So, if $\delta(\cA)=0$, we get that $\Sigma \backslash \Sigma'$
contains periodic points of arbitrarily large period, so it is
infinite as wanted. Notice that $\cA$ may have infinite sets
admitting a dominated splitting and still verifying
$\delta_{min}(\cA)=0$. The converse part is the main result of
\cite{BGV} (see Theorem 4.1 of \cite{BGV}).

Part $(d)$ is given by the fact that the determinant is
multiplicative, so, an $\eps$-perturbation of a cocycle, can
increase at most $(1+\eps)^{\pi(x)}$ the determinant of $M_x^A$.
So, the sum of the exponents can change at most $\log(1+\eps)$
which converges to zero as $\eps\to 0$.

Part $(e)$ is contained in Lemmas 4.4 and 4.6 of \cite{BDP}.

To prove property (f) we use that the existence of a dominated
splitting admits a change of metric which makes the subbundles
ortogonal (see \cite{BDP} section 4.1). So, we can write the
cocycle restricted to $F\oplus G\oplus H$ in the form (using
coordinates, $G, (G\oplus H)\cap G^\perp  , F$)

$$  \left(
      \begin{array}{ccc}
        \cA_G & \star& 0 \\
        0 & \cA (H) & 0 \\
        0 & 0 & \cA_F \\
      \end{array}
    \right)  $$

Where $\cA(H)= \cA|(G\oplus F)$. So the cocycle $\cA|G$ is written
in coordinates $H|G, F|G$ as

$$ \left(
     \begin{array}{cc}
        \cA (H) & 0 \\
         0 & \cA_F \\
     \end{array}
   \right) $$

Since $F$ is dominated by $G\oplus H$, we get the desired
property. Notice that for any $x\in \Sigma$, we have that $\|
(A_{x}|_{G\oplus H})^{-1} \|^{-1} \leq \| (A(H)_x)^{-1}\|^{-1}$
and this guaranties the domination of $F|G \prec H|G$ as desired.
\lqqd

We are now ready to state the main result of this appendix. A much
stronger version of this result can be found\footnote{In fact, the
result of \cite{BoBo} was announced before I had a proof of this
result but I did not notice the overlap and pursued in this
direction.} in \cite{BoBo}. The proof here presented of this
weaker result has some similarities with their proof but I hope
that its inclusion in the text is not devoid of interest and
introduces some different ways of proving some parts of the
result.

\begin{teo}\label{BGV}
Let $\mathcal{A}=(\Sigma, f, E, A)$ be a bounded large period
linear cocycle of dimension $d$. Assume that
\begin{itemize}
\item[-] $\cA$ is strictly without domination. \item[-]
$\mathcal{A} \in \Gamma_{\Sigma,i}$ \item[-] For every $x\in
\Sigma$, we have $|det(M_x^A)|< 1$ (that is, for all $x\in \Sigma$
we have $\sum_{j=1}^d \sigma^j(x,\cA) <0$).
\end{itemize}
Then, $\delta_{min}^{\ast,i}(\cA)=0$. In particular, given
$\eps>0$ there exists a point $x\in \Sigma$ and a path $\gamma_x$
of diameter smaller than $\eps$ such that $\gamma_x(0)=\cA_x$, the
matrix $M_x^{\gamma_x(1)}$ has all its eigenvalues of modulus
smaller than one, and such that $\gamma(t) \in \Gamma_{i}$ for
every $t \in [0,1]$.
\end{teo}

\section{Proof of Theorem \ref{BGV}}

This section will be devoted to prove this theorem. The proof is
by induction.

The following Lemma allows to find several invariant subbundles in
order to be able to apply induction. It is proved in
\cite{GouFranksLemma} Proposition 6.6.

\begin{lema}\label{valorespropiosdistintos} For every $\mathcal A \in \Gamma_{\Sigma,i}$, there exists $\mathcal B \sim^{\ast}_i \cA$ such that for every $x\in \Sigma$ the eigenvalues of $M_x^B$ have all different modulus and their modulus is arbitrarily near the original one in $M_x^A$, that is, $|\sigma^i(x,\cA)-\sigma^i(x,\cB)| \to 0$ as $\pi(x)\to \infty$ (in particular, $\delta(\cB)=\delta(\cA)$).
\end{lema}

\esbozo{\!\!\!} We proceed by induction. In dimension $2$ the
result is the same as in Proposition 3.7 of \cite{BGV} (the only
perturbations done there can be made along paths without any
difficulty, this was first done in \cite{BC}).  Notice that if the
eigenvalues are both equal for some $x\in \Sigma$, then
necessarily the cocycle belongs to $\Gamma_{\Sigma,2}$ or
$\Gamma_{\Sigma,0}$.

We assume the result holds in dimension $<d$. Since there always
exists an invariant subspace of dimension $2$, you can make
independent perturbations and change the eigenvalues as required
using the induction hypothesis. For this, one should perturb in
the invariant subspace and in the quotient (see Remark
\ref{deapedazos}).

\lqqd

If a cocycle $\cA$ verifies that for every $x\in \Sigma$, the
eigenvalues of $M_x^A$ have all different modulus and different
from $1$, we shall say that $\cA$ is a \emph{diagonal cocycle}.

\begin{obs} For a diagonal cocycle $\cA \in \Gamma_{\Sigma,i}$ one has well defined invariant one-dimensional subspaces $E_1(x,\cA) , \ldots , E_d(x,\cA)$ (we shall in general omit the reference to the point and/or the cocycle) associated to the eigenvalues in increasing order of modulus. Also,  if $F_l=E_1\oplus \ldots \oplus E_l$, one gets that $\cA_{F_l} \in \Gamma_{\Sigma, j}$ for any $0\leq j \leq \min \{i,l\}$.

From now on, for diagonal cocycles we shall name $F_l = E_1 \oplus
\ldots \oplus E_l$ and $G_l = E_l \oplus \ldots \oplus E_d$.
\finobs
\end{obs}

As we said, Theorem \ref{BGV} is easier in dimension 2 (it does
not even need the uniformity hypothesis on the determinant). The
following Lemma is essentially due to Ma\~ne and will be the base
of the induction.

\begin{lema}\label{dimensiondos} Let $\mathcal A \in \Gamma_{\Sigma,i}$ be a bounded large period linear cocycle of dimension $2$ $(0\leq i\leq 2)$ and strictly without domination such that $|det(M_x^A)|\leq 1$ for every $x\in \Sigma$. Then, there exists $\mathcal B \sim^{\ast}_i \mathcal A$  with the following properties:
 \begin{enumerate}
 \item $\delta(\mathcal B)=0$.
 \item $|det(M_x^A)|=|det(M_x^B)|$ for every $x\in \Sigma$.
 \end{enumerate}
\end{lema}

\dem{\!\!\!} This is very standard (see \cite{ManheErgodic} or
section 7.2.1 of \cite{BDV}). With a small perturbation (see
Proposition 6.7 of \cite{GouFranksLemma}) one can make the angle
between the stable and unstable spaces arbitrarily small and not
change the determinant.

After that, one can compose with a rotation, of determinant equal
to $1$ (so, without affecting the product of the modulus of the
eigenvalues), since after rotating a small amount one gets complex
eigenvalues, there is a moment where the eigenvalues are real and
arbitrarily near, there is where we stop.

Notice that if $|det(M^A_x)|>1$ and $i=1$ we would be obliged to
stop the path longtime before the eigenvalues are nearly equal
since the smallest exponent would attain the value $0$ which is
forbidden for path perturbations of index $1$. \lqqd

Before we continue with the induction to prove the general result,
we shall make some general perturbative results which loosely
state that if two bundles are not dominated, then, after a
perturbation which preserves the exponents, we can see the non
domination in the bundles associated to the closest bundles.

First we will state a standard linear algebra result we shall use
in order to perturb two dimensional cocycles.

\begin{lema}\label{amovervectores} Given $\eps>0$  and $K>0$, there exists $\ell>0$ such that if $A_1, \ldots , A_\ell$ is a sequence in $GL(2,\R)$ matrices verifying that $\max_{i} \{\|A_i\|, \|A_i^{-1}\|\} \leq K$ and $v, w\in \R^2$ are vectors with $\|v\|=\|w\|=1$. Suppose that

$$ \| A_\ell  \ldots  A_1 v\| \geq \frac 1 2 \| A_\ell  \ldots  A_1 w \|$$

Then, there exists rotations $R_1, \ldots, R_\ell$ of angle
smaller than $\eps$ verifying that

$$ R_\ell  A_\ell  \ldots  R_1  A_1 \R w = A_\ell  \ldots  A_1 \R  v $$

\end{lema}

\dem{\!\!} For simplicity  we shall assume that $A_i  v = v$ for
every $i$. Since this is made by composing each matrix by a
rotation and a homothety, we should change $K$ by $K^2$ which will
be the new bound for the norm of the matrices.

For $\gamma \in P^1(\R)$, let $\alpha_i(\gamma) = \frac {\| A_i
A_{i-1} \ldots A_1 z \|} { \|A_{i-1} \ldots A_1 z \|}$ where $z$
is a vector in the direction $\gamma$. It is a well known result
in linear algebra that if the function $\alpha_i : P^1(\R) \to
\R^+$ is not constant then it has a maximum and a minimum in
antipodal points and it is monotone in the complement.


We can assume that for every $i\geq 0$, we have that $A_i \ldots
A_1 \R w$ is at distance larger than $\eps$ from $\R v$, otherwise
we can perform the perturbation.

Notice also that given $0<\eta<\nu<1$ there exists $\kappa$ such
that if for some point $\gamma \in P^1(\R)$ we have that
$\prod_{i=k}^j \alpha_i(\gamma) > \kappa$, then (notice that
$\prod_{i=k}^j\alpha_i(\R v)=1$), there is an interval of length
$\nu$ around $\gamma$ which does not contain $\R v$ which is
mapped by $A_j \ldots A_k$ to an interval around $A_j\ldots A_k
\gamma$ of length $\eta$. A similar statement holds for the
inverses in the case the product is smaller than $\kappa^{-1}$.

If we choose $\nu>1-\eps$ and $\eta<\eps$ we get
$\kappa=\kappa(\eps)$ which will verify the following: Assume that
there exists $\gamma \in P^1(\R)$ verifying $\prod_{i=j}^k
\alpha_i(\gamma) < \kappa^{-1}$ for some $j,k$, then, since $A_j
\ldots A_1 \R w$ should be in the interval around $\gamma$, we can
first rotate it to send it to an extreme point of the interval,
and after applying $A_k \ldots A_{j+1}$ the vector will be $\eps$
close to $\R v$ and so we can finish the perturbation.

This implies that for $\kappa=\kappa(\eps)$  we get that for every
$j<k$ we have that $\prod_{i=j}^k \alpha_i(\gamma) > \kappa^{-1}$
for any $\gamma \in P^1(\R)$.

The hypothesis of the Lemma (and the choice of $\alpha_i(\R v)=1$
for every $i$) implies that $\prod_{i=1}^\ell \alpha_i(\R w) \leq
2$.

This implies that also (maybe by rechoosing $\kappa$) that we can
not have a sequence $\alpha_i(\R w)$ verifying $\prod_{i=j}^k
\alpha_i(\R w) > \kappa$. Which in turn implies that this should
also happen for every point $\gamma \in P^1(\R)$ (otherwise an
iterate of $\R w$ would be $\eps$ close to $v$).

This gives us that, there exists $0<\rho<1$ such that the iterates
of an interval of length $\eps$ remain of length bounded from
below by $\rho \eps$. Now, choosing $\ell$ such that $\ell \rho
\eps>1$ we get that we can take any vector to $\R v$ with
rotations of angle less than $\eps$.

\lqqd

The following proposition is the key step of the proof.

Notice that the Proposition does not use the fact that $\cA$ may
be chosen strictly without domination, or even with
$\delta(\cA)=0$.

\begin{prop}\label{dim2tricotomia} Given $K>0$, $k>0$ and $\eps>0$, there exists $N>0$ and $\ell$ such that if
\begin{itemize}
\item[-] There exists a diagonal cocycle $\cA_x$ of dimension $2$
and bounded by $K$ over a periodic orbit of period $\pi(x)>N$.
\item[-] There exists a unit vector $v\in E_x$  such that

$$  \frac 3 2 \| A^\ell_{f^k(x)}|_{E_1}\| \geq \frac{ \|A^\ell_{f^k(x)} A^k_x v\|}{\|A^k_x v\|} $$

\end{itemize}

 Then, there exists $\cB_x$ a path perturbation of $\cA_x$ of diameter smaller than $\eps$ verifying that all along the path the cocycle has the same Lyapunov exponents and

$$    \| B^{k}_{x} |_{E_1(x,\cB)} \| \leq 2 \|A^{k}_{x}v\|$$

\end{prop}

\dem{\!\!}  We use coordinates $E_1 \oplus E_1^\perp$. With this
convention, we have

$$  A_{x} = \left(
             \begin{array}{cc}
                \alpha_x & K_x \\
                 0 & \beta_x \\
             \end{array}
            \right)   \qquad   A^n_x = \left(
              \begin{array}{cc}
                \prod_{j=0}^{n-1}\alpha_{f^j(x)} & \star  \\   
                0 & \prod_{j=0}^{n-1}\beta_{f^j(x)} \\
                \end{array}
             \right)    $$

We also have that $|\alpha_x|$, $|\beta_x|$ and $|K_x|$ are
uniformly bounded from above by $K$ and also $|\alpha_x|$ and
$|\beta_x|$ are bounded from below by $1/K$ for every $x\in
\Sigma$. It is satisfied that $\left| \prod_{j=1}^{\pi(x)}
\alpha_{f^j(x)}\right| < \left|\prod_{j=1}^{\pi(x)}
\beta_{f^j(x)}\right|$.

We fix $k>0$ and $\eps>0$. Let $\ell>0$ given by Lemma
\ref{amovervectores} and let $v\in E_x$ a vector in the hypothesis
of the Proposition.

We consider the set $\Upsilon \en E$ of vectors satisfying that if
$w_0 \in \Upsilon$ is an unitary vector and if we denote $w_j =
\frac{A^j_x w }{\| A^j_x w\|}$ we have

$$  \| A^\ell_{f^k(x)}|_{E_1(x,\cA)}\| \geq \frac 1 2 \|A^\ell_{f^k(x)} w_k\| $$

and

$$  \|A^k_x w_0\| \leq 2 \|A^{k}_{x}v\|$$

We remark that $\Upsilon$ is defined just in terms of $\{A_x
,\ldots A_{f^{k+\ell}(x)}\}$ provided that we maintain the
condition on all $A_{f^j(x)}$ being triangular. Also, it is easy
to see that $\Upsilon$ is a closed under scalar multiplication, so
we shall sometimes consider the unit vectors there and think of it
as a subset of the projective line $P^1(\R)$.

Notice that if $E_1(x,\cA) \in \Upsilon$, the Proposition holds
without need to make any perturbation, so we will assume that it
is not the case. We shall call $\theta$ the distance in $P^1(\R)$
between $E_1(x,\cA)$ and $\Upsilon$ where the distance we consider
is the one given by the inner angle between the generated lines.

We shall consider $k+\ell < L < \pi(x)-\ell$ the largest integer
(if there exists any) verifying that there exists some $w\in
\Upsilon$ satisfying

$$ \|A^\ell_{f^L(x)}  w_L \| \geq \frac 1 2 \| A^\ell_{f^L(x)}|_{E_1(f^L(x),\cA)}\|   \qquad \qquad \qquad (4.1)$$

\begin{claim} If $\pi(x)$ is large enough, there exists some $L$ with the properties above. Moreover, $L \to \infty$ as $\pi(x)\to \infty$.
\end{claim}

\dem{}  Assume that for some $s>0$  and for arbitrarily large
$\pi(x)$ (for simplicity we consider it of the form $\pi(x)=R\ell
+s$ with $R\to \infty$), if $L$ exists is smaller than $s$.

We have that for $w \in \Upsilon$, the nearest vector to
$E_1(x,\cA)$ we have that

$$ \|A^{R\ell}_{f^{s}(x)} w_s \| < \left(\frac 1 2 \right)^R \| A^{R\ell}_{f^{s}(x)}|_{E_1(f^{s}(x), \cA)} \| $$

A simple calculation gives us, choosing $R$ large enough to
satisfy that $K^{s} \left(\frac 1 2 \right)^R \leq \tan
(\theta/2)$, that the cone of angle $\theta/2$ around
$E_1(f^{s}(x), \cA)$ is mapped by $A^{\pi(x)}_{f^{s}(x)}$ inside
itself. This contradicts the fact that $E_1$ is the eigenspace
associated to the smallest eigenvalue (i.e. that  $\left|
\prod_{j=1}^{\pi(x)} \alpha_{f^j(x)}\right| <
\left|\prod_{j=1}^{\pi(x)} \beta_{f^j(x)}\right|$).

\finobs

We shall need a more quantitative version of this growth:

\begin{claim} There exists $N>0$ such that if $\pi(x)> N +2\ell +k$ then

$$\left(\prod_{j=k+l}^{L-1}\alpha_{f^j(x)}(1+\eps)^{-1}\right) K^{2\ell+k}\|A_{f^{L+\ell}(x)}^{\pi(x)-L-\ell}\| \leq \prod_{j=1}^{\pi(x)} \alpha_{f^j(x)} \leq$$
$$\leq \left(\prod_{j=k+l}^{L-1}\alpha_{f^j(x)}\right)
K^{2\ell+k}\|A_{f^{L+\ell}(x)}^{\pi(x)-L-\ell}\|$$.
\end{claim}

\dem The second inequality is direct. We  prove the claim by
contradiction, we assume though that

$$  \prod_{j \notin \{ k+\ell, \ldots , L-1\}} \alpha_{f^j(x)} < \| A^{\pi(x)-L-\ell}_{f^{L+\ell}(x)} \| K^{2\ell+k} (1+\eps)^{-L+k+\ell} $$

Notice that by the previous claim  we have that $L\to \infty$ as
$\pi(x)\to \infty$ so, as $\pi(x)$ grows, the norm of $A^{\pi(x)-L
-\ell}_{f^{L+\ell}(x)}$ grows to infinity compared to the norm of
$A^{\pi(x)-L-\ell}_{f^{L+\ell}(x)}|_{E_1}$.

Also, we get that the distance between $A^{j}_x \Upsilon$ and
$E_1(f^{j}(x),\cA)$ for $j\geq L$ must be bounded from bellow
since the norm of every $A^\ell_{f^i(x)}$ is bounded by $K^\ell$
so being very close implies that $j$ would satisfy $(4.1)$
contradicting the maximality of $L$.

So, if $\pi(x)$ is large enough, we get that the vectors far from
$E_1(f^{L+\ell}(x), \cA)$ must be mapped near the direction of
maximal expansion of $A^{\pi(x)+k +\ell -L}_{f^{L+\ell}(x)}$ and
thus, this allows to find a vector in $\Upsilon$ verifying $(4.1)$
a contradiction with the maximality of $L$. See Lemma
\ref{amovervectores} for a similar argument.

\finobs

We shall now define the perturbation we will make in order to
satisfy the requirements of the Proposition. We shall define a
continuous path $\gamma: [0,1] \to \Gamma_x$ of diameter smaller
than $\eps$ and verifying the required properties. Notice that all
along the path, the determinant is never changed in any point, so,
the product of the eigenvalues of $M_x^{\gamma(t)}$ remains
unchanged for every $t\in [0,1]$.

Using Lemma \ref{amovervectores} we can perturb with small
rotations the transformations $A_{f^L(x)}, \ldots,
A_{f^{L+\ell-1}(x)}$ in order to send $E_1(f^L(x), \cA)$ into the
one dimensional subspace $\R \tilde w$ such that
$A^{\pi(x)-(L+\ell)}_{f^{L+\ell-1}(x)} \R \tilde w=\R w$ where
$w\in \Upsilon$ is the vector defining $L$.

Since this perturbations are made by composing with small
rotations, they can be made along small paths. Let $A_{f^j(x)}^t:
[0,1] \to GL(E_{f^j(x)}, E_{f^{j+1}(x)})$ where $j \in \{L,\ldots,
L+\ell -1\}$ such that

$$ \left(A^{\pi(x)-L}_{f^L(x)} A^1_{f^{L+\ell-1}(x)} \ldots A^1_{f^L(x)}  A^{L-(k+\ell)}_{f^{k+\ell}(x)}\right) E_1(f^{k+\ell}(x), \cA) = \R \tilde w$$

The hypothesis we made and the previous Lemma allow us to send $\R
w$ to $E_1(f^{k+\ell}(x), \cA)$ by composing the matrices
$A_{f^j(x)}$ with $k\leq j <k+\ell$ with small rotations. Clearly,
by choosing properly the rotations, for $k \leq j \leq k+\ell-1$
we can find also paths $A_{f^j(x)}^t: [0,1] \to GL(E_{f^j(x)},
E_{f^{j+1}(x)})$ verifying that

$$ \left(A^t_{f^{k+\ell-1}(x)} \ldots A^t_{f^{k}(x)} A^k_{x}  A^{\pi(x)-L}_{f^L(x)}  A^t_{f^{L+\ell-1}(x)}  \ldots  A^t_{f^L(x)}  A^{L-(k+\ell)}_{f^{k+\ell}(x)}\right) E_1(f^{k+\ell}(x), \cA) = E_1(f^{k+\ell}(x), \cA) $$

We shall also perturb the linear transformations  $A_{f^j(x)}$
with $j\in \{k+\ell,\ldots,L-1\}$ by multiplying them by matrices
of the form

$$  \left(
      \begin{array}{cc}
        \frac{1}{\alpha(t)} & 0 \\
        0 & \alpha(t) \\
      \end{array}
    \right) $$

Where $\alpha(t)$ is conveniently chosen in $[1, 1+\eps]$ in order
to get that for every $t \in [0,1]$ the two exponents of
$M^{\gamma(t)}_{f^{k+\ell}(x)}$ coincide.

To show that the latter can be made, notice that the perturbations
we made imply that in our coordinates, the matrix
$M^{\gamma(t)}_{f^{k+\ell}(x)}$ is of the same triangular form,
so, after applying the first $L-k-\ell$ transformations, the $E_1$
direction will remain horizontal, so, Claim 2 implies that for
every $t\in [0,1]$, there exists $\alpha(t)$ such that the
exponents are equal, the fact that this $\alpha(t)$ varies
continuously is given by the fact that the eigenvalues of a path
of matrices vary continuously. This concludes.

\lqqd

We shall extend this two dimensional result to a more general
context using this kind of two dimensional perturbations. This
will allow us to reduce all the problems to a two dimensional
context that we know well how to treat (see Lemma
\ref{dimensiondos}).

\begin{prop}\label{dimntricotomia} Let $\cA \in \Gamma_{\Sigma,i}$ a bounded diagonal cocycle.
Assume that for $0<j<d$ we have that $F_j$ is not dominated by $G_{j+1}$.
Then, there exists $\cB$, a path perturbation of $\cA$ along a path which does not change any
of the Lyapunov exponents, which verifies that $E_j(\cB)$ is not dominated by $E_{j+1}(\cB)$.
\end{prop}

\dem We shall prove that if $F_j$ is not dominated by $G_{j+1}$ we
can perturb as above in order to break the domination between
$F_j$ and $E_{j+1}$. A symmetric argument gives the desired
property.

We prove this by induction. So, we will fix $j$ and assume that
the proposition holds for every $d<d_0>j+2$ and prove it in the
case $d=d_0$. We assume that $F_j$ is not dominated by $G_{j+1}$
but that $F_j \prec E_{j+1}$ (otherwise there is nothing to
prove).

This implies (by property (e) of Lemma \ref{BGVenunciado}) that
$F_j | E_{j+1}$ is not dominated by $G_{j+2} | E_{j+1}$, so, we
can by induction, find a perturbation respecting all the
eigenvalues such that $F_j | E_{j+1}$ is not dominated by $E_{j+2}
| E_{j+1}$. Now, property (f) of Lemma \ref{BGVenunciado} implies
that $F_j$ is not dominated by $E_{j+1} \oplus  E_{j+2}$. Using
induction again we obtain a perturbation which respects the
eigenvalues and such that that $F_j$ is not dominated by $E_{j+1}$
as wanted.

Finally we must prove the Proposition in the case $d=j+2$ in order
to conclude. Assume then that $F_j$ is dominated by $E_{j+1}$
(otherwise there is nothing to prove).

We have that there exists $s$ such that $F_j \prec_{s} E_{j+1}$.
For simplicity, we take $s=1$, that is, for every $x\in \Sigma$
(maybe by considering an infinite subset), and unitary vectors
$v_j \in F_j $ and $v_{j+1} \in E_{j+1}$ we have that

$$  \| A_x v_j \| \leq \frac 1 2   \| A_x v_{j+1} \| $$

Since we have assumed that $F_j$ is not dominated by
$E_{j+1}\oplus E_{j+2}$ we have that for every $n>0$ there exists
$N$ such that if $\pi(x)>N$ we have that for some point of the
orbit of $x$ (which without loss of generality we suppose is $x$)
one has that for some unitary vectors $v_j\in E_j $ and  $v\in
E_{j+1}\oplus E_{j+2} $ we have that

$$    2 \| A^n_x v_j \| >  \| A^n_x v \| $$

Using standard arguments (see for example Pliss' Lemma
\cite{pliss,wen-Selecting}) we get that for every $k$ and $\ell$
there exists $N$ such that if $\pi(x)>N$ we have that (again
choosing $x$ conveniently and maybe by changing the vectors by
their normalized iterates)

$$  \frac {\| A^\ell_{f^k(x)} A^k_x v \|}{\|A^k_x v\|} < \frac 3 2 \frac{ \| A^\ell_{f^k(x)} A^k_x v_j \|}{\|A^k_x v_j\|} < \frac 3 2 \| A^\ell_{f^k(x)}|_{E_{j+1}}\|$$

and

$$ \| A^k_x v \| < \frac 3 2 \| A^k_{x} v_j \|  $$

This puts us in the hypothesis of Proposition \ref{dim2tricotomia}
which allow us to make a perturbation of $\cA_{E_{j+1} \oplus
E_{j+2}}$ without changing the Lyapunov exponents and that breaks
the domination between $F_j$ and $E_{j+1}$.

Notice that given $k$ and $\eps$ we get that for periodic points
with large period we can perform these perturbations with size
smaller than $\eps$, so, we get that we can perturb a sequence of
periodic orbits with period going to infinity with arbitrarily
small perturbations and break the domination.

\lqqd

This proposition allows us to complete the proof of Theorem
\ref{BGV}. Before, we shall make some reductions.

\begin{lema}\label{realizarelminimo} Let $\cA\in \Gamma_{\Sigma,i}$. Then, there exists $\cB \sim^{\ast}_i \cA$ such that $\delta(\cB)=\delta^{\ast,i}_{min}(\cA)$. Moreover, we can assume that $\cB$ is a diagonal cocycle and $\delta(\cB_{F_l}) = \delta_{min}^{\ast,j}(\cB_{F_l})$ for every $1\leq l \leq d$ and $j = \min\{i,l\}$.
\end{lema}

\dem{\!\!} The existence of $\cB$ is proved by following verbatim
the proof of Lemma 4.3 in \cite{BGV}, since the proof does not
introduce new perturbations. The idea is to take a sequence
$\cB_n$ of path perturbations of index $i$ with $\delta(\cB_n)$
converging to $\delta_{min}^{\ast,i}(\cA)$ with different
eigenvalues (see Lemma \ref{valorespropiosdistintos}), and then
considering the cocycle $\cB$ defined as coinciding with $\cB_n$
over the periodic points of period $n$.

To prove that we can choose
$\delta(\cB_{F_l})=\delta_{min}^{\ast,j}(\cB_{F_l})$ we make
another diagonal process to first take $\delta(\cB_{F_{d-1}})$ to
$\delta_{min}^{\ast,i}(\cB)$, then $\delta(\cB_{F_{d-2}})$ and so
on. \lqqd

From now on, we shall use the following notation

$$ \overline{\sigma}^j(\cA) = {\displaystyle \limsup_{\pi(x)\to \infty} \sigma^j(x,\cA)} \qquad  \underline{\sigma}^j(\cA) = {\displaystyle \liminf_{\pi(x)\to \infty} \sigma^j(x,\cA)} $$

\begin{obs}\label{restriccioninfinita} Notice that if $\Sigma'\en \Sigma$ is an invariant infinite subset, we get that $\delta_{min}^{\ast,i}(\cA|_{\Sigma'}) \geq \delta_{min}^{\ast,i}(\cA)$. So, we can always restrict to an infinite invariant subset to prove the Theorem.

This implies that to prove the Theorem, we can assume that the
cocycle $\cA$ satisfies that for every $1\leq j \leq d$ we have
$\underline \sigma^j (\cA)=\overline \sigma^j(\cA) =
\sigma^j(\cA)$. To do this it is enough to make a diagonal process
showing that there is an infinite subset $\Sigma_1\en \Sigma$
where $\underline \sigma^1(\cA|_{\Sigma_1})= \overline
\sigma^1(\cA|_{\Sigma_1})$. Then, inductively, we can construct
$\Sigma_k \en \ldots \en \Sigma_1 \en \Sigma$ an infinite subset
such that for every $1\leq j \leq k$ we have $\underline
\sigma^j(\cA|_{\Sigma_k}) = \overline \sigma^j(\cA|_{\Sigma_k})$.
Finally, we restrict $\cA$ to $\Sigma_d$ and use the previous
remark.

In this context, we get also that $\delta(\cA) = \sigma^d(\cA) -
\sigma^1(\cA)$. \finobs
\end{obs}

We shall say that a cocycle $\cA\in \Gamma_{\Sigma,i}$ of
dimension $d$ is $i$-\emph{incompressible} if it satisfies the
following properties (notice that they are quite more restrictive
than the ones used in \cite{BGV}):

\begin{itemize}
\item[-] $\cA$ is a diagonal cocycle. \item[-] $\delta(\cA_{F_l})
= \delta_{min}^{\ast, j}(\cA_{F_l})$ where $j=\min\{i, l\}$ for
every $1\leq l \leq d$. \item[-] For every $1\leq j \leq d$ we
have that $\underline{\sigma}^j(\cA)= \overline \sigma^j(\cA)=
\sigma^j(\cA)$.
\end{itemize}

The previous Lemma and Remark \ref{restriccioninfinita} show that
to prove Theorem \ref{BGV} it is enough to work with
$i-$incompressible cocycles strictly without domination. Notice
that trivially, if $\cB \sim_i^\ast \cA$ is a diagonal path
perturbation such that for every $j$ we have
$\sigma^j(\cB)=\sigma^j(\cA)$, then, $\cB$ is also
$i-$incompressible.

\dem{of Theorem \ref{BGV}} We shall prove this Theorem by
induction. As we said we can work with $i-$incompressible
cocycles.

Now we make the standing induction hypothesis, which holds for two
dimensional cocycles after Lemma \ref{dimensiondos}. It is easy to
see that  proving this implies directly the Theorem.

\begin{itemize}
\item[(H)] Let $\cD\in \Gamma_{\Sigma,i}$ be an $i-$incompressible
cocycle of dimension $k<d$ ($0\leq i\leq k$) verifying that for
every $x\in \Sigma$ we have that $|det(M_x^D)|<1$. Then, if $j<k$
is the first number such that $\sigma^j(\cD) < \sigma^{j+1}(\cD)$
then, it holds that $F_j \prec G_{j+1}$.
 \end{itemize}

Now, let us consider an $i-$incompressible cocycle $\cA\in
\Gamma_{\Sigma,i}$ of dimension $d$. Let $j$ be the smallest
number such that $\sigma^j(\cA) < \sigma^{j+1}(\cA)$ (if no such
$j$ exists there is nothing to prove).

If $j=d-1$, we shall show that $F_{d-1} \prec E_d$.

We notice first that since the sum of all exponents is $\leq 0$,
we have that $\sigma^1(\cA) =\sigma^{d-1}(\cA) < 0$.

Assume that $F_{d-1}$ is not dominated by $E_d$. So, by
Proposition \ref{dimntricotomia} (used for the inverses) we get
that for some $\cB\sim_i^\ast \cA$ that remains $i-$incompressible
(since it does not change the eigenvalues) we have that $E_{d-1}$
is not dominated by $E_d$. But this is a contradiction since Lemma
\ref{dimensiondos} allows us to decrease the Lyapunov diameter of
$\cB_{E_{d-1} \oplus E_d}$ contradicting the $i-$incompressibility
(notice that this would make the last exponent to decrease).

So it rest to prove the theorem in the case $j<d-1$.

First of all, by induction we get that $F_j \prec E_{j+1} \oplus
\ldots \oplus E_{d-1}$.

Now, if $F_j$ is not dominated by $G_{j+1}$ we get that a
perturbation which preserves the $i-$incompressibility (given by
Proposition \ref{dimntricotomia}), allows us to break the
domination $F_j \prec E_{j+1} \oplus \ldots \oplus E_{d-1}$ a
contradiction.

\lqqd

   \chapter{Plane decompositions}\label{Apendix-PlaneDec}

We present a construction of a plane diffeomorphism $f$ of bounded
$C^\infty$ norm which is semiconjugated to the homothety $x\mapsto
x/2$ by a continuous map $h:\RR^2\to \RR^2$ whose fibers are all
cellular sets of diameter smaller than $K$ and such that it has
two disjoint attracting neighborhoods which project by $h$ to the
whole plane. We derive some unexpected consequences of the
existence of such an example.

We shall denote as $d_2:\RR^2 \to \RR^2$ to the map

$$ d_2 (x) = \frac x 2 .$$

The goal of this note is to prove the following theorem (recall
subsection \ref{Section-Descomposiciones}):

\begin{teo}\label{teoPPal} There exists a $C^\infty$-diffeomorphism $f:\RR^2 \to \RR^2$ and a constants $K>0$ and $a_K>0$ such that the following properties are verified:
\bi \item[-] There exists a (H\"{o}lder) continuous cellular map
$h:\RR^2\to \RR^2$ such that $d_{C^0}(h,id) < K$ and $d_2 \circ h
= h \circ f$. \item[-] There exist open sets $V_1$ and $V_2$ such
that \bi \item $\overline{V_1}\cap \overline{V_2}=\emptyset$
    \item $h(V_i)= \RR^2$ for $i=1,2$.
    \item $f(\overline{V_i}) \en V_i$ for $i=1,2$.
\ei \item[-] The $C^{\infty}$ norm of $f$ and $f^{-1}$ is smaller
than $a_K$. \ei
\end{teo}

A direct consequence of this Theorem is the existence of $h:\RR^2
\to \RR^2$ whose fibers are all non trivial and cellular
(decreasing intersection of topological disks), the existence of
these decompositions of the plane had been shown by Roberts
\cite{Roberts}.

\section{Construction of $f$.}\label{sectionConstruccion}

For simplicity, we start by considering a curve $\gamma= \{0\}
\times [-\frac 1 4, \frac 1 4]$ (the construction can be made
changing $\gamma$ for any cellular set\footnote{A cellular set is
a decreasing intersection of compact topological disks.}).

Clearly, $\gamma \subset B_0 = B_1(0)$ the ball of radius one on
the origin. We shall also consider the sets $B_n = B_{2^n}(0)$ for
every $n\geq 0$. We have that

$$ \RR^2= \bigcup_{n\geq 0} B_n $$

We shall define $f:\RR^2 \to \RR^2$ with the desired properties in
an inductive manner, starting by defining it in $B_0$ and then in
the annulus $B_n \setminus B_{n-1}$.

Let us define $f_0 : \overline{B_0} \to B_0$ a $C^\infty$
embedding and disjoint open sets $V_1^0$ and $V_2^0$ such that:

\bi \item[(a)] $f_0$ coincides with $d_2$ in a small neighborhood
of $\partial B_0$. \item[(b)] $\bigcap_{n\geq 0} f_0^n (B_0) =
\gamma$. \item[(c)] $f_0 (\overline{V_i^0})\en V_i^0$ for $i=1,2$.
\item[(d)] The sets $V_i^0$ are diffeomorphic to $[-1,1] \times
\RR$, separate $B_0$ in two connected components and intersect
$\{0\} \times [-1/4,1/4]$ in disjoint closed intervals.
\ei

Now, we assume that we have defined a $C^\infty$-diffeomorphism
$f_n: B_n \to B_{n-1}$ and disjoint open connected sets $V_1^n$
and $V_2^n$ (homeomorphic to a band $\RR \times (0,1)$) such that:

\bi \item[(I1)] $f_n|_{B_{n-1}}=f_{n-1}$ and $V_i^{n-1} \en V_i^n$
for $i=1,2$. \item[(I2)] The $C^\infty$-distance between $f_n$ and
$d_2$ in $B_n$ is smaller than $a_K$. \item[(I3)]
$(f_n(\overline{V_i^n}) \setminus \partial B_{n-1}) \en V_i^{n-1}$
for $i=1,2$ and $f_n^n(\overline{V_i^n})$ disconnects $B_0$.
\item[(I4)] $V_i^n$ is $K/2$-dense in $B_n$. \item[(I5)] $f_n$
coincides with $d_2$ in a $K/10$-neighborhood of $\partial B_n$.
\ei

We must now construct $f_{n+1}$ assuming we had constructed $f_n$
and this will define a diffeomorphism $f:\RR^2\to \RR^2$ which we
shall show has the desired properties.

To construct $f_{n+1}$ and $V_i^n$ we notice that in order to
verify (I1), it is enough to define $f_{n+1}$ in $B_n \backslash
B_{n-1}$ as well as to add to $V_i^n$ an open set in
$B_{n+1}\backslash B_{n}$ in order to verify the hypothesis.

We consider $d_2^{-1}(V_i^n) \cap B_{n+1}\backslash B_n$ which
since $V_i^n$ was $K/2$-dense in $B_n$ becomes $K$-dense in
$B_{n+1}\backslash B_{n}$ for $i=1,2$.

We shall use the following lemma whose proof we delay to the end
of this section.

\begin{lema}\label{Lema-Cinfinitochicp}
There exists $a_K$ which only depends on $K$ such that:
\begin{itemize}
\item[-] Given two open sets $A_1, A_2$ which are $K$-dense inside
a set of the form $B_n \setminus B_{n-1}$ with sufficiently large
$n$. \item[-] The sets $A_i$ verify that for every point in $A_i$
there is a curve going from $\partial B_n$ to $\partial B_{n-1}$
and contained in $A_i$.
\end{itemize}
Then there exists a $C^\infty$ diffeomorphism $g$ of
$C^\infty$-norm less than $a_K$ such that coincides with the
identity in $K/10$-neighborhood of the boundaries and such that
the image by $g$ of the open sets is $K/2$-dense in $B_n \setminus
B_{n-1}$.
\end{lema}

We consider a diffeomorphism $g$ given by the previous lemma
which is $a_K-C^\infty$-close to the identity, coincides with the
identity in the $K/10$-neighborhoods of $\partial B_{n+1}$ and
$\partial B_n$ and such that $g(V_i^n)$ is $K/2$-dense for
$i=1,2$.

We define then $f_{n+1}$ in $B_{n+1}\backslash B_n$ as $d_2 \circ
g^{-1}$ which clearly glues together with $f_n$ and satisfies
properties (I2) and (I5).

To define $V_i^{n+1}$ we consider a very small $\eps>0$ (in order
that $g(V_i^n)$ is also $K/2-\eps$-dense) and for each boundary
component $C$ of $g(V_i^n)$ (which is a curve) we consider a curve
$C'$ which is at distance less than $\eps$ of $C$ inside
$g(V_i^n)$ and such that each when it approaches $C\cap \partial
B_n$ the distance goes to zero and when it approaches $C \cap
\partial B_{n+1}$ the distance goes to $\eps$. This allows to
define new $V_i^{n+1}$ as the open set delimited by these curves
united with the initial $V_i^n$. It is not hard to see that it
will satisfy (I3) and (I4).

We have then constructed a $C^\infty$-diffeomorphism $f:\RR^2 \to
\RR^2$ which is at $C^\infty$ distance $a_K$ of $d_2$ and such
that there are two disjoint open connected sets $V_1$ and $V_2$
such that $f(\overline{V_i})\en V_i$. and such that both of them
are $K/2$-dense in $\RR^2$.

We now indicate the proof of the Lemma we have used:

\demo{ of Lemma \ref{Lema-Cinfinitochicp}} The proof follows from
the following simple bound:

\begin{af}
There exists $A>0$ such that for every pair of curves
$\gamma_1,\gamma_2$ in the square $[-2,2]^2$ which touch both
boundaries and intersect $[-1,1]^2$ we have that there exists a
$C^\infty$-diffeomorphism $h$ of $C^\infty$-norm less than $A$ and
which coincides with the identity in the boundaries of the cube
such that the image of both curves is $1/4$-dense.
\end{af}

We can assume that $n$ is large enough since we can get a bound by
hand on the rest of $B_n$'s.

Now, to prove the Lemma it is enough to subdivide the complement
of the $K/10$-neighborhoods of the boundaries of $B_n \setminus
B_{n-1}$ into sets $S_k$ such that they contain balls of radius
$4K$ and are contained in balls of radius $5K$. Moreover, we can
choose these sets $S_k$ in order to verify that for some positive
constant $B$ we have:

\begin{itemize}
\item[-] $S_k$ is diffeomorphic to $[-2,2]^2$ via a $C^\infty$
diffeomorphism $l_k: [-2,2]^2 \to S_k $ of norm less than $B$.
\item[-] The diffeomorphism $l_k$ sends $1/4$-dense subsets of
$[-2,2]^2$ into $K/2$-dense subsets of $S_k$.
\end{itemize}

Now, using the claim it is not hard to see that we can construct
the desired diffeomorphism $g$ of $C^\infty$-norm less than $AB$.

\lqqd

\section{Proof of the Theorem}

We first show the existence of a continuous function $h:\RR^2 \to
\RR^2$ conjugating $f$ to $d_2$ which is close to the identity.

This is quite classical, consider a point $x \in \RR^2$, so, since
$d_{C^0}(f,d_2)<K$ we get that the orbit $\{f^n(x)\}$ is in fact a
$K-$pseudo-orbit of $d_2$. Since $d_2$ is infinitely expansive,
there exists only one orbit $\{d_2^n(y)\}$ which
$\alpha(K)$-shadows $\{f^n(x)\}$ and we define $h(x)=y$ (in fact,
in this case, it suffices with the past pseudo-orbit to find the
shadowing).

We get that $h$ is continuous since when $x_n \to x$ then the
pseudo-orbit which shadows must rest near for more and more time,
and then, again by expansivity, one concludes. This implies also
that $h$ is onto since it is at bounded distance of the identity.

Now, consider any ball $B$ of radius $100 \alpha(K)$ in $\RR^2$,
it is easy to see that $f(B)$ is contained in a ball of radius
$50\alpha(K)$ and then, we get a way to identify the preimage of
points by $h$. Consider a point $x \in \RR^2$, we get that

$$h^{-1}(h(x)) = \bigcap_{n>0} f^n( B_{100 \alpha(K)} (f^{-n}(x))) $$

So, $h$ is also cellular.

It only remains to show that the image under $h$ of both $V_1$ and
$V_2$ is the whole plane. Since they share equal properties, it
will be enough to prove it for one of them.

\begin{lema} $h(V_i)= \RR^2$ for $i=1,2$.
\end{lema}

\dem We shall show that $h(\overline{V_i})$ is dense. Since $h$ is
proper, it is closed: this will imply that it is in fact the whole
plane. Using the semiconjugacy and the fact that
$f(\overline{V_i})\en V_i$ this would prove the lemma.

To prove that $h(\overline{V_i})$ is dense, we consider an
arbitrary open set $U\en \RR^2$. Now, choose $n_0$ such that
$d_2^{-n_0}(U)$ contains a ball of radius $10 \alpha(K)$. We get
that $h^{-1}(d_2^{-n_0}(U))$ contains a ball of radius $9
\alpha(K)$ and thus, since $\alpha(K)>K$, we know that since $V_i$
is $K/2$-dense, we get that $V_i \cap h^{-1}(d_2^{-n_0}(U)) \neq
\emptyset$. So, since $f(\overline{V_i}) \en V_i$ we get that $V_i
\cap f^{n_0} \circ h^{-1}(d_2^{-n_0}(U)) \neq \emptyset$ which
using the semiconjugacy gives us that $h(V_i) \cap U \neq
\emptyset$.

This concludes.\lqqd

\section{H\"{o}lder continuity}

In this section, we shall prove that in fact $h$ is
$\alpha$-Holder continuous (see also Theorem 19.2.1 of \cite{KH}).
Since the boundary of $\partial V_i$ for each $i$ is a space
filling curve in arbitrarily small domains and by some easy
estimates on the change of Hausdorff dimension by H\"{o}lder maps,
we see easily that $\alpha \leq \frac{1}2$.

To prove the existence of $\alpha>0$ such that $h$ is
$\alpha$-H\"{o}lder, consider $C$ to be a (uniform) bound on
$\|Df^{-1}\|$ (recall that $f$ can be choosen ``close'' to $d_2$).
We choose also $\alpha>0$ such that $C^\alpha <2 $.

Also, from how we constructed the semiconjugacy $h$, we see that
there exists $A_1$ and $A_2$ such that $d(x,y)<A_1$ implies that
$d(h(x),h(y))<A_2$. Now, consider a pair of points $x,y \in \RR^2$
such that $\delta= d(x,y)$ is sufficiently small (say, smaller
than $A_1$).

We consider $n_0$ such that $C^{-n_0}\delta < A_1 \leq C^{-n_0-1}
\delta$. We have that

$$ 2^{-n_0} A_1^\alpha \leq 2^{-n_0} C^{-\alpha n_0} C^{-\alpha} \delta^\alpha \leq C^{-\alpha} \delta^\alpha $$

Now, since $d_2^n \circ h \circ f^{-n}(x) = h(x)$ for every $n$
and $x$ we get

$$ d(h(x),h(y)) = d(d_2^{n_0} \circ h \circ f^{-n_0}(x),d_2^{n_0} \circ h \circ f^{-n_0}(y)) = 2^{-n_0} d(h(f^{-n_0}(x)), h(f^{-n_0}(y)))$$

But, from how we choose $C$, we get that $d(f^{-n_0}(x),
f^{-n_0}(y))<A_1$, so

$$ d(h(x),h(y)) \leq 2^{-n_0} A_2 $$

From the above, we obtain thus

$$ d(h(x),h(y)) \leq 2^{-n_0} A_1^\alpha \frac{A_2}{A_1^\alpha} \leq C^{-\alpha} \frac{A_2}{A_1^\alpha} \delta^\alpha = \left(C^{-\alpha} \frac{A_2}{A_1^\alpha}\right) d(x,y)^\alpha$$

   \chapter{Irrational pseudo-rotations of the torus}\label{Apendice-PseudoRotaciones}

This appendix was taken from \cite{PotRecurrence}, but we have
added Section \ref{Section-DehnTwist}.

We consider $\Homeo_0(\TT^2)$ to be the set of homeomorphisms
homotopic to the identity. We shall say that $f\in
\Homeo_0(\TT^2)$ is \emph{non-resonant} if the rotation set of $f$
is a unique vector $(\alpha,\beta)$ and the values $1, \alpha,
\beta$ are irrationally independent (i.e. $\alpha, \beta$ and
$\alpha/\beta$ are not rational). This ammounts to say that given
any lift $F$ of $f$ to $\RR^2$, for every $z\in \RR^2$ we have
that:

\begin{equation}\label{equationRotacion}  \lim_{n\to \infty} \frac{F^n(z)-z}{n} = (\alpha, \beta) (\modulo \ZZ^2) \end{equation}

In general, one can define the rotation set of a homeomorphism
homotopic to the identity (see \cite{MZ}). In fact, although we
shall not make it explicit, our constructions work in the same way
for homeomorphisms of the torus whose rotation set is contained in
a segment of slope $(\alpha,\beta)$ with $\alpha, \beta$ and
$\alpha/\beta$ irrational and not containing zero.

Non-resonant torus homeomorphisms\footnote{These are called
\emph{irrational pseudo-rotations} by several authors, but since
some of them use the term exclusively for conservative ones, we
adopt the definition used in \cite{Kwakkel}.} have been
intensively studied in the last years looking for resemblance
between them and homeomorphisms of the circle with irrational
rotation number (see \cite{Kwap}, \cite{LeCalvez}, \cite{Jager})
and also constructing examples showing some difference between
them (see \cite{Fayad}, \cite{BCL}, \cite{BCJL}, \cite{Jager2}).

In \cite{Kwakkel} the possible topologies of minimal sets these
homeomorphisms admit are classified and it is shown that under
some conditions, these minimal sets are unique and coincide with
the non-wandering set\footnote{A point $x$ is \emph{wandering} for
a homeomorphism $f$ if there exists a neighborhood $U$ of $x$ such
that $f^n(U) \cap U = \emptyset$ for every $n\neq 0$. The
\emph{non-wandering} set is the closed set of points which are not
wandering.}. However, there is one kind of topology of minimal
sets where the question of the uniqueness of minimal sets remains
unknown. When the topology of a minimal set is of this last kind,
\cite{BCJL} constructed an example where the non wandering set
does not coincide with the unique minimal set, in fact, they
construct a transitive non-resonant torus homeomorphism containing
a proper minimal set as a skew product over an irrational
rotation.

A natural example of non-resonant torus homeomorphism is the one
given by a homeomorphism semiconjugated to an irrational rotation
by a continuous map homotopic to the identity. In \cite{Jager} it
is proved that a non-resonant torus homeomorphism is
semiconjugated to an irrational rotation under some quite mild
hypothesis.

Under the hypothesis of being semiconjugated by a monotone
map\footnote{A monotone map is a map whose preimages are all
compact and connected.} which has points whose preimage is a
singleton, it is not hard to show the uniqueness of a minimal set
(see for example \cite{Kwakkel} Lemma 14). However, as shown in
Appendix \ref{Apendix-PlaneDec} (see also \cite{Roberts}), a
continuous monotone map may be very degenerate and thus even if
there exist such a semiconjugation, it is not clear whether there
should exist a unique minimal set nor the kind of recurrence the
homeomorphisms should have. Moreover, for general non-resonant
torus homeomorphisms, there does not exist a semiconjugacy to the
irrational rotation (even when there is ``bounded mean motion'',
see \cite{Jager2}).

Here, we give a simple and self-contained proof (based on some
ideas of \cite{Kwakkel} but not on the classification of the
topologies of the minimal sets) of a result which shows that even
if there may be more than one minimal set, the dynamics is in some
sense irreducible. Clearly, transitivity of $f$ may not hold for a
general non-resonant torus homeomorphism (it may even have
wandering points, as in the product of two Denjoy counterexamples;
some more elaborate examples may be found in \cite{Kwakkel}), but
we shall show that, in fact, these homeomorphisms are weakly
transitive. For a homeomorphism $f$ we shall denote $\Omega(f)$ to
the non-wandering set of $f$ (i.e. the set of points $x$ such that
for every neighborhood $U$ of $x$ there exists $n>0$ with
$f^n(U)\cap U \neq \emptyset$).

\begin{teo}\label{MainTeo} Let $f\in \Homeo_0(\TT^2)$ be a non-resonant torus homeomorphism, then, $f|_{\Omega(f)}$ is weakly transitive.
\end{teo}

Recall that for $h:M\to M$ a homeomorphism, and $K$ an
$h-$invariant compact set, we say that $h|_K$ is \emph{weakly
transitive} if given two open sets $U$ and $V$ of $M$ intersecting
$K$, there exists $n>0$ such  that $h^n(U)\cap V \neq \emptyset$
(the difference with being transitive is that for transitivity one
requires the open sets to be considered relative to $K$).

This allows to re-obtain Corollary E of \cite{Jager}:

\begin{cor}\label{CorolarioTransitividad} Let $f\in \Homeo_0(\TT^2)$ be a non-resonant torus homeomorphism such that $\Omega(f)=\TT^2$. Then, $f$ is transitive.
\end{cor}

In fact, as a consequence of weak-transitivity, we can obtain also
the more well known concept of chain-transitivity for non-resonant
torus homeomorphisms.

\begin{cor}\label{Addendum} Let $f\in \Homeo_0(\TT^2)$ be a non-resonant torus homeomorphism, then, $f$ is chain-transitive.
\end{cor}

Recall that a homeomorphism $h$ of a compact metric space $M$ is
\emph{chain-transitive} if for every pair of points $x,y\in M$ and
every $\eps>0$ there exists an $\eps-$pseudo-orbit $x= z_0,
\ldots, z_n=y$ with $n\geq 1$ (i.e. $d(z_{i+1}, h(z_i))<\eps$).

\dem Consider two points $x,y\in M$ and $\eps>0$.

 We first assume that $x \neq y$ are both nonwandering points which shows the idea in a simpler way. From Theorem A we know that there exists a point $z$ and $n>0$ such that $d(z,f(x))<\eps$ and $d(f^{n+1}(z), y)<\eps$. We can then consider the $\eps-$pseudo-orbit: $\{x, z, \ldots, f^n(z), y\}$.

Now, for general $x, y \in \TT^2$ we consider $n_0\geq 1$ such
that $d(f^{n_0+1}(x), \Omega(f))<\eps/2$ and $d(f^{-n_0}(y),
\Omega(f))<\eps/2$. Now, by Theorem A there exists $z\in \TT^2$
and $n>0$ such that $d(z,f^{n_0}(x))<\eps$ and
$d(f^{n+1}(z),f^{-n_0}(y))<\eps$. Considering the following
$\eps-$psudo-orbit $\{x, \ldots, f^{n_0-1}(x), z, \ldots, f^n(z),
f^{-n_0}(y), \ldots, y\}$ we obtain a pseudo-orbit from $x$ to $y$
and thus proving chain-transitivity. \lqqd

\begin{obs} We have proved that in fact, for every $\eps>0$ the pseudo-orbit can be made with only two ``jumps''.
\end{obs}

As a consequence of our study, we obtain the following result
which may be of independent interest:

\begin{prop}\label{PropConexosCpctsSeCortan} Let $f \in \Homeo_0(\TT^2)$ be a non-resonant torus homeomorphism and $\Lambda_1$ a compact connected set  such that $f(\Lambda_1) \en \Lambda_1$. Then, for every $U$ connected neighborhood of $\Lambda_1$, there exists $K>0$ such that:
\bi \item[-] If $\Lambda_2$ is a compact set which has a connected
component in the universal cover of diameter larger than $K$
then\footnote{This holds  if $\Lambda_2$ is a connected set such
that $f^i(\Lambda_2)\en \Lambda_2$ for some $i\in\ZZ$ for
example.}, \ei \noindent $U \cap \Lambda_2 \neq \emptyset$.
\end{prop}

One could wonder if the stronger property of $\Omega(f)$ being
transitive may hold. However, in section \ref{SectionExample} we
present an example  where $\Omega(f)$ is a Cantor set times
$\SS^1$, but for which the nonwandering set is not transitive.

\section{Reduction of the proofs of Theorem \ref{MainTeo} and Proposition \ref{PropConexosCpctsSeCortan}}

In this section we shall reduce the proofs of Theorem
\ref{MainTeo} and Proposition \ref{PropConexosCpctsSeCortan} to
Proposition \ref{MainProp} and its Addendum \ref{AdeStronger}.

We shall use the word \emph{domain} to refer to an open and
connected set. We shall say a domain $U \in \TT^2$ is
\emph{inessential, simply essential or doubly essential} depending
on whether the inclusion of $\pi_1(U)$ in $\pi_1(\TT^2)$ is
isomorphic to $0, \ZZ$ or $\ZZ^2$ respectively\footnote{In
\cite{Kwakkel} these concepts are called \emph{trivial, essential}
and \emph{doubly-essential}.}. If $U$ is simply essential or
doubly essential, we shall say it is \emph{essential}.

\begin{obs}\label{ObsDoublyEssentialsecortan} Notice that if $U$ and $V$ are two doubly essential domains, then $U\cap V\neq \emptyset$. This is because the intersection number of two closed curves is a homotopy invariant and given two non-homotopic curves in $\TT^2$, they have non-zero intersection number, thus, they must intersect. Since clearly, being doubly essential, $U$ and $V$ contain non homotopic curves, we get the desired result. \finobs
\end{obs}

We claim that Theorem \ref{MainTeo} can be reduced to the
following proposition.

\begin{prop}\label{MainProp}
Given $f\in \Homeo_0(\TT^2)$ a non-resonant torus homeomorphism
and $U$ an open set such that $f(U) \en U$ and $U$ intersects
$\Omega(f)$, then we have that $U$ has a connected component which
is doubly essential.
\end{prop}

Almost the same proof also yields the following statement which
will imply Proposition B:

\begin{ad}\label{AdeStronger} For $f$ as in Proposition \ref{MainProp}, if $\Lambda$ is a compact connected set such that $f(\Lambda) \en \Lambda$, then, for every connected open neighborhood $U$ of $\Lambda$, we have that $U$ is doubly-essential.
\end{ad}

Notice that the fact that $f(\Lambda)\en \Lambda$ for $\Lambda$
compact implies that it contains recurrent points, and in
particular, $\Lambda\cap \Omega(f)\neq \emptyset$.

\demo{of Theorem A and Proposition B} Let us consider two open
sets $U_1$ and $V_1$ intersecting $\Omega(f)$, and let $U =
\bigcup_{n > 0} f^n(U_1)$ and $V = \bigcup_{n < 0} f^n(V_1)$.
These sets verify that $f(U) \en U$ and $f^{-1}(V)\en V$ and both
intersect the nonwandering set.

Proposition \ref{MainProp} (applied to $f$ and $f^{-1}$) implies
that both $U$ and $V$ are doubly essential, so, they must
intersect. This implies that for some $n > 0$ and $m < 0$ we have
that $f^n(U_1) \cap f^m(V_1) \neq \emptyset$, so, we have that
$f^{n-m}(U_1)\cap V_1 \neq \emptyset$ and thus $\Omega(f)$ is
weakly transitive.

Proposition B follows directly from Addendum \ref{AdeStronger}
since given a doubly-essential domain $U$ in $\TT^2$, there exists
$K>0$ such that its lift $p^{-1}(U)$ intersects every connected
set of diameter larger than $K$. \lqqd

\begin{obs}\label{ObsDimensionmayor} Notice that in higher dimensions, Remark \ref{ObsDoublyEssentialsecortan} does not hold. In fact, it is easy to construct two open connected sets containing closed curves in every homotopy class which do not intersect. So, even if we could show a result similar to Proposition \ref{MainProp}, it would not imply the same result. \finobs
\end{obs}

\section{Proof of Proposition \ref{MainProp}}

Consider a non-resonant torus homeomorphism $f \in \Homeo_0(\TT^2)$, and let us assume that $U$ is an open set which verifies $f(U)\en U$ and $U \cap \Omega(f)\neq \emptyset$. 

Since $U\cap \Omega(f)\neq \emptyset$, for some $N>0$ we have that
there is a connected component of $U$ which is $f^N$-invariant. We
may thus assume from the start that $U$ is a domain such that
$f(U)\en U$ and $U\cap \Omega(f) \neq \emptyset$.

Let $p:\RR^2 \to \TT^2$ be the canonical projection. Consider $U_0
\en p^{-1}(U)$ a connected component. We can choose $F$ a lift of
$f$ such that $F(U_0) \en U_0$.

We shall denote $T_{p,q}$ to the translation by vector $(p,q)$,
that is, the map from the plane such that $T_{p,q}(x)= x+(p,q)$
for every $x\in \RR^2$.

\begin{lema}\label{LemaNoTrivial} The domain $U$ is essential.
\end{lema}

\dem \ Consider $x\in U_0$ such that $p(x)\in \Omega(f)$. And
consider a neighborhood $V\en U_0$ of $x$. Assume that there
exists $n_0>0$ and $(p,q)\in \ZZ^2 \setminus \{(0,0)\}$ such that
$F^{n_0}(V)\cap (V +(p,q)) \neq \emptyset$. Since $U_0$ is
$F$-invariant, we obtain two points in $U_0$ which differ by an
integer translation, and since $U_0$ is connected, this implies
that $U$ contains a non-trivial curve in $\pi_1(\TT^2)$ and thus,
it is essential.

To see that there exists such $n_0$ and $(p,q)$, notice that
otherwise, since $x$ is not periodic (because $f$ is a
non-resonant torus homeomorphism) we could consider a basis $V_n$
of neighborhoods of $p(x)$ such that $f^k(V_n) \cap V_n =
\emptyset$ for every $0<  k \leq n$. Since $x$ is non-wandering,
there exists some $k_n > n$ such that $f^{k_n}(V_n)\cap V_n \neq
\emptyset$, but since we have that $F^{k_n}(V_n)\cap (V_n +(p,q))
= \emptyset$ for every $(p,q)\in \ZZ^2 \setminus \{(0,0)\}$, we
should have that $F^{k_n}(V_n)\cap V_n \neq \emptyset$ for every
$n$. Since $k_n \to \infty$, we get that $f$ has zero as rotation
vector, a contradiction. \lqqd


We conclude the proof of by showing the following lemma which has
some resemblance with Lemma 11 in \cite{Kwakkel}.

\begin{lema}\label{LemaNoEsencial} The domain $U$ is doubly-essential.
\end{lema}

\dem Assume by contradiction that $U$ is simply-essential.

Since the inclusion of $\pi_1(U)$ in $\pi_1(\TT^2)$ is non-trivial
by the previous lemma, there exists a closed curve $\eta$ in $U$
such that when lifted to $\RR^2$ joins a point $x \in U_0$ with
$x+(p,q)$ (which will also belong to $U_0$ because $\eta$ is
contained in $U$ and $U_0$ is a connected component of
$p^{-1}(U)$).

We claim that in fact, we can assume that $\eta$ is a simple
closed curve and such that $g.c.d(p,q)=1$ (the greatest common
divisor). In fact, since $U$ is open, we can assume that the curve
we first considered is in general position, and by considering a
subcurve, we get a simple one (maybe the point $x$ and the vector
$(p,q)$ changed, but we shall consider the curve $\eta$ is the
simple and closed curve from the start). Since it is simple, the
fact that $g.c.d(p,q)=1$ is trivial.

If $\eta_0$ is the lift of $\eta$ which joins $x\in U_0$ with
$x+(p,q)$, we have that it is compact, so, we get that

$$ \tilde \eta = \bigcup_{n\in \ZZ} T_{np,nq}\eta_0$$

\noindent is a proper embedding of $\RR$ in $\RR^2$. Notice that
$\tilde \eta \en U_0$.

By extending to the one point compactification of $\RR^2$ we get
by using Jordan's Theorem (see \cite{Moise} chapter 4) that
$\tilde \eta$ separates $\RR^2$ in two disjoint unbounded
connected components which we shall call $L$ and $R$ and such that
their closures $L\cup \tilde \eta$ and $R\cup \tilde \eta$ are
topologically a half plane (this holds by Sch\"onflies Theorem,
see \cite{Moise} chapter 9).

Consider any pair $a,b$ such that\footnote{We accept division by
$0$ as being infinity.} $\frac{a}{b} \neq \frac{p}{q}$, we claim
that $T_{a, b} (\tilde \eta) \cap U_0 = \emptyset$. Otherwise, the
union $T_{a,b}(\tilde \eta) \cup U_0$ would be a connected set
contained in $p^{-1}(U)$ thus in $U_0$ and we could find a curve
in $U_0$ joining $x$ to $x+(a,b)$ proving that $U$ is doubly
essential (notice that the hypothesis on $(a,b)$ implies that
$(a,b)$ and $(p,q)$ generate a subgroup isomorphic to $\ZZ^2$), a
contradiction.

Translations are order preserving, this means that $T_{a,b}(R)
\cap R$ and $T_{a,b}(L)\cap L$ are both non-empty and either
$T_{a,b}(R)\en R$ or $T_{a,b}(L) \en L$ (both can only hold in the
case $\frac{a}{b} = \frac{p}{q}$). Also, one can easily see that
$T_{a,b}(R)\en R$ implies that $T_{-a,-b}(L)\en L$.

Now, we choose $(a,b)$ such that there exists a curve $\gamma$
from $x$ to $x+(a,b)$ satisfying:

\bi \item[-] $T_{a,b}(\tilde \eta) \en L$. \item[-] $\gamma$ is
disjoint from $T_{p,q}(\gamma)$. \item[-] $\gamma$ is disjoint
from $T_{a,b}(\tilde \eta)$ and $\tilde \eta$ except at its
boundary points. \ei

We consider $\tilde \eta_1 = T_{a,b}(\tilde \eta)$ and $\tilde
\eta_2 = T_{-a,-b}(\tilde \eta)$. Also, we shall denote $\tilde
\gamma= \gamma \cup T_{-a,-b}(\gamma)$ which joins $x-(a,b)$ with
$x+(a,b)$.

We obtain that $U_0$ is contained in $\Gamma= T_{a,b}(R)\cap
T_{-a,-b}(L)$ a band whose boundary is $\tilde \eta_1 \cup \tilde
\eta_2$.

Since $U_0$ is contained in $\Gamma$ and is $F$-invariant, for
every point $x\in U_0$ we have that $F^n(x)$ is a sequence in
$\Gamma$, and since $f$ is a non-resonant torus homeomorphism, we
have that $\lim \frac{F^n(x)}{n} = \lim \frac{F^n(x)-x}{n}=
(\alpha,\beta)$ is totally irrational.

However, we notice that $\Gamma$ can be written as:

$$ \Gamma = \bigcup_{n \in \ZZ} T_{np, nq} (\Gamma_0)$$

\noindent where $\Gamma_0$ is a compact set in $\RR^2$. Indeed, if
we consider the curve $\tilde \gamma \cup T_{a,b}(\eta_0) \cup
T_{p,q}(\tilde \gamma) \cup T_{-a,-b}(\eta_0)$ we have a Jordan
curve. Considering $\Gamma_0$ as the closure of the bounded
component we have the desired fundamental domain.

So, if we consider a sequence of points $x_n \in \Gamma$ such that
$\lim \frac{x_n}{n}$ exists and is equal to $v$ it will verify
that the coordinates of $v$ have the same proportion as $p/q$,
thus cannot be totally irrational. This is a contradiction and
concludes the proof of the Lemma.

\lqqd


We conclude this section by showing how the proof adapts to the
case stated in Addendum \ref{AdeStronger}. Consider a compact
connected set $\Lambda$ such that $f(\Lambda) \en \Lambda$, then,
we have that $\Lambda$ contains points which are
recurrent\footnote{Since it is a compact invariant set, it
contains a minimal set whose points will be all recurrent.}.

Let $\tilde \Lambda$ be a connected component of $p^{-1}(\Lambda)$
which is $F$-invariant. Now, if $U$ is an open connected
neighborhood of $\Lambda$ and $U_0$ is a connected component of
$p^{-1}(U)$ containing $\tilde \Lambda$. Notice that $d(\partial
U, \Lambda) > \delta >0$ so $d(\partial U_0, \tilde \Lambda) >
\delta$ also.

Now, the same argument in Lemma \ref{LemaNoTrivial} can be used in
order to show that $U$ must be essential: We can choose a point
$x\in \tilde \Lambda \en U_0$ such that $p(x)$ is recurrent and
the same argument shows that there will exist $(p,q) \in \ZZ^2
\setminus \{(0,0)\}$ such that $F^{n_0}(x)$ is $\delta$-close to
$x +(p,q)$ and since $F^{n_0}(x)$ must be contained in $\tilde
\Lambda$ we get that $x + (p,q)$ is contained in $U_0$ showing
that $U$ is essential.

The proof that in fact $U$ is doubly-essential is now the same as
in Lemma \ref{LemaNoEsencial} since one can see that invariance of
$U$ was not used in the proof, one only needs that there are
points in $U_0$ such that the orbits by $F$ remain in $U_0$ and
this holds for every point in $\tilde \Lambda$.

\section{An example where $f|_{\Omega(f)}$ is not transitive}\label{SectionExample}

The example is similar to the one in section 2 of \cite{Jager2},
however, we do not know a priori if our specific examples admit or
not a semiconjugacy.

Consider $g_1:S^1 \to S^1$ and $g_2: S^1 \to S^1$ Denjoy
counterexamples with rotation numbers $\rho_1$ and $\rho_2$ which
are irrationally independent and have minimal invariant sets $M_1$
and $M_2$ properly contained in $S^1$. We shall consider the
following skew-product map $f_\beta:\TT^2 \to \TT^2$ given by:

$$ f_\beta(s,t) = (g_1(s), \beta(s)(t)) $$

\noindent where $\beta: S^1 \to \Homeo_+(S^1)$ is continuous and
such that $\beta(s)(t)=g_2(t)$ for every $(s,t)\in M_1 \times
S^1$.

The same proof as in Lemma 2.1 of \cite{Jager2} yields:

\begin{lema}\label{LemaUnicoMinimal} The map $f_\beta$ is a non-resonant torus homeomorphism and $M_1 \times M_2$ is the unique minimal set.
\end{lema}

\dem The proof is the same as the one in Lemma 2.1 of
\cite{Jager2}. Indeed any invariant measure for $f$ must be
supported in $M_1 \times M_2$ and the dynamics there is the
product of two Denjoy counterexamples and thus uniquely-ergodic.
Since rotation vectors can be computed with ergodic measures, we
also get that $f_\beta$ has a unique rotation vector
$(\rho_1,\rho_2)$ which is totally irrational by hypothesis. \lqqd

Clearly, if we restrict the dynamics of $f_\beta$ to $M_1 \times
S^1$ it is not hard to see that the nonwandering set will be
$M_1\times M_2$ (it is a product system there). So, we shall prove
that if $\beta$ is properly chosen, we get that $\Omega(f_\beta)=
M_1 \times S^1$. In fact, instead of constructing a specific
example, we shall show that for ``generic'' $\beta$ in certain
space, this is satisfied, this will give the existence of such a
$\beta$.

First, we define $\cB$ to be the set of continuous maps $\beta:
S^1 \to \Homeo_+(S^1)$ such that $\beta(s)=g_2$ for every $s\in
M_1$. We endow $\cB$ with the topology given by restriction from
the set of every continuous map from $S^1$ to $\Homeo_+(S^1)$.
With this topology, $\cB$ is a closed subset of the set of
continuous maps from $S^1 \to \Homeo_+(S^1)$ which is a Baire
space, thus, $\cB$ is a Baire space.

So, the existence of the desired $\beta$ is a consequence of:

\begin{lema} There exists a dense $G_\delta$ (residual) subset of $\cB$ of maps such that the induced map $f_\beta$ verifies that $\Omega(f_\beta)= M_1 \times \SS^1$.
\end{lema}

\dem \ First, we will prove the lemma assuming the following
claim:

\begin{af} Given $\beta \in \cB$, $x\in M_1 \times S^1$, $\eps>0$ and $\delta>0$ there exists $\beta' \in \cB$ which is $\delta-$close to $\beta$ such that there exists $k>0$ with $f_{\beta'}^k (B(x,\eps)) \cap B(x,\eps) \neq \emptyset$.
\end{af}

Assuming this claim, the proof of the Lemma is a standard Baire
argument: Consider $\{x_n\} \en M_1\times S^1$ a countable dense
set. Using the claim, we get that the sets $\cB_{n,N}$ consisting
of the functions $\beta \in \cB$ such that there exists a point
$y$ and a value $k>0$ such that $y$ and $f^k_{\beta}(y)$ belong to
$B(x_n,1/N)$ is a dense set. Also, the set $\cB_{n,N}$ is open,
since the property is clearly robust for $C^0$ perturbations of
$f_\beta$. This implies that the set $\cR= \bigcap_{n,N}
\cB_{n,N}$ is a residual set, which implies, by Baire's theorem
that it is in fact dense.

For $\beta \in \cR$ we get that given a point $x\in M_1\times S^1$
and $\eps>0$, we can choose $x_n \in B(x,\eps/2)$ and $N$ such
that $1/N<\eps/2$. Since $\beta \in \cB_{n,N}$ we have that there
exists $k>0$ such that $f_\beta^k(B(x,\eps)) \cap B(x,\eps) \neq
\emptyset$ proving that $M_1 \times S^1$ is nonwandering for
$f_\beta$ as desired.

\demo{of the Claim} The point $x \in M_1 \times S^1$ can be
written as $(s,t)$ in the canonical coordinates.

Choose an interval $(a,b) \en (s-\eps, s+\eps)$ contained in a
wandering interval of $g_1$. Then, there exists a sequence of
integers $k_n\to +\infty$ such that $g_1^{k_n}((a,b)) \en
(s-\eps,s+\eps)$ for all $n\geq 0$. Further, the orbits of $a$ and
$b$ are disjoint and do not belong to $M_1$. Let $\gamma= (a,b)
\times \{t\}$.

We can assume that $f_\beta^{k_n}(\gamma) \cap B(x,\eps) =
\emptyset$ for every $n>0$, otherwise, there is nothing to prove.

We shall thus consider a $\delta-$perturbation of $\beta$ such
that it does not modify the orbit of $(a,t)$ but moves the orbit
of $(b,t)$ in one direction making it give a complete turn around
$\SS^1$ and thus an iterate of $\gamma$ will intersect
$B(x,\eps)$.

Let $s_n = g_1^n(b)$ and $\beta^n(s_0) = \beta(s_{n-1})\circ
\ldots \circ \beta(s_0)$. Note that $\frac{1}{n}\beta^n(s_0)(t)
\to \rho(g_2)$ as $n\to \infty$ since $\beta(s_k) \to g_2$ as
$k\to \infty$. At the same time, if we let

$$\beta_\theta^n(s_0)= R_\theta \circ \beta(s_{n-1}) \circ R_\theta \circ \beta_{s_{n-2}} \circ \ldots \circ R_\theta \circ \beta(s_0) $$

Then, $\beta_\theta^n(s_0)(t) \to \rho' > \rho(g_2)$ since
$R_\theta \circ \beta_{s_k}$ converges to $R_\theta \circ g_2$
which has rotation number strictly greater than $g_2$ (see for
example \cite{KH} Proposition 11.1.9). If we denote by $\tilde
\beta^n$, respectively $\tilde \beta^n_\theta$ the lifts of
$\beta^n$ and $\beta_\theta^n$ to $\RR$, then, this implies that
there exists $n_0$ such that for $n>n_0$ one has

$$|\tilde \beta_\theta^n(t) - \tilde \beta^n(t) | > 1$$

So, if we consider $k_n>n_0$ and we choose $\beta'$ such that:

\bi \item[-] it coincides with $\beta$ in the $g_1$-orbit of $a$,
\item[-] it coincides with $R_\theta \circ \beta$ in the points
$\{b, g_1(b), \ldots, g_1^{k_n}(b)\}$, \item[-] is at distance
smaller than $\delta$ from $\beta$, \ei

we have that $f^{k_n}_{\beta'}(\gamma) \cap B(x, \eps) \neq
\emptyset$ as desired.

\finobs \lqqd

\section{The homotopy class of the dehn-twist}\label{Section-DehnTwist}

We provide here a simple extension of the results of this appendix
to homeomorphisms which are not homotopic to the identity. We call
\emph{dehn twist} to the torus homeomorphism whose lift to the
universal cover can be written in the form $(x,y) \mapsto (x,
x+y)$.

We say that a homeomorphism $f : \TT^2 \to \TT^2$ homotopic to the
dehn-twist is \emph{non-resonant} if there exists $\alpha \in \RR
\setminus \QQ$ such that for every $x \in \RR^2$ and for some lift
$F: \RR^2 \to \RR^2$ one has that:

$$ \lim_{n\to \infty} \frac{p_1(F^n(x)) - p_1(x)}{n} = \alpha $$

\noindent where $p_1: \RR^2 \to \RR$ is the projection in the
first coordinate. A classical example is given by the projection
of the torus of the following plane map:

$$ (x,y) \mapsto (x + \alpha, x+y) $$

With essentially the same proof as Theorem \ref{MainTeo} we can
prove:

\begin{teo}
Let $f : \TT^2 \to \TT^2$ be a non-resonant torus homeomorphism homotopic to a dehn-twist. Then $\Omega(f)$ is weakly transitive.
\end{teo}

The proof is essentially the same as the one of Theorem
\ref{MainTeo} so we will only give a sketch.

In fact, the result can be reduced to the following statement
analogous to Proposition \ref{MainProp}.

\begin{prop}
Given $f:\TT^2 \to \TT^2$ a non-resonant torus homeomorphism
homotopic to a dehn-twist and $U$ an open set such that $f(U) \en
U$ and $U$ intersects $\Omega(f)$, then we have that $U$ has a
connected component which is doubly essential.
\end{prop}

The proof of the reduction is exactly the same as for Theorem
\ref{MainTeo}.

Let us now prove the proposition.

\esbozo{} As in the proof of Proposition \ref{MainProp} we can
assume that $U$ is connected and such that $f(U) \en U$.

Let us call $U_0$ to a connected component of $p^{-1}(U)$ and we
will choose a lift $F:  \RR^2 \to \RR^2$ which fixes $U_0$.

We consider $x \in U_0$ such that $p(x) \in \Omega(f)$. Let $V \en
U_0$ be a connected neighborhood of $x$.

Assume that there exists $n_0>0$ such that $F^{n_0}(V) \cap
(V+(p,q))\neq \emptyset$. Then, as in Lemma \ref{LemaNoTrivial} we
can prove that there is a loop $\gamma$ contained in $U$ which is
homotopic to the loop joining $(0,0)$ with $(p,q)$.

Now, assume that $p\neq 0$, then, we have that $f(\gamma)$ which
is homotopic to $(p,q+1)$ is also contained in $U$ by invariance
of $U$ and the fact that $f$ is homotopic to a dehn-twist. Since
$(p,q)$ and $(p,q+1)$ are linearly independent when $p\neq 0$ we
get that $U$ should be doubly essential.

To prove that there exists $n_0$ such that $F^{n_0}(V) \cap
(V+(p,q))\neq \emptyset$ for some $(p,q)$ with $p\neq 0$ we use
the fact that the rotation number defined above is irrational and
an argument very similar to that of Lemma \ref{LemaNoTrivial}.

\lqqd

   \chapter{Tame non-robustly transitive diffeomorphisms}\label{Apendice-BCGP}

In this appendix, we review the results from \cite{BCGP}. Recall
the discussion after Corollary \ref{Corolario-BonattiCrovisier}.

Given an $n-$dimensional manifold $M$ (with $n\geq 3$), we
consider $\Diff^r(M)$ the set of diffeomorphisms of $M$ endowed
with the $C^r$ topology ($r\geq 1$).

\begin{teo*} There exists a $C^1$-open set $\cU \en \Diff^r(M)$ ($1\leq r \leq \infty$),

a $C^r$-dense subset $\cD$ of $\cU$ and an open set $U\en M$ with
the following properties:

\begin{itemize}
\item[(I)] \emph{Isolation}: For every $f\in \cU$, the set
$\cC_f:=U\cap \cR(f)$ is a chain-recurrence class.

\item[(II)] \emph{Non-robust transitivity}: For every $f\in \cD$,
the class $\cC_f$ is not transitive.
\end{itemize}

More precisely:

\begin{itemize}
\item[(1)] For any $f\in \cU$ there exist a subset $H_f \subset
\cC_f$ which coincides with the homoclinic class of any hyperbolic
periodic $x\in \cC_f$. Moreover, each pair of hyperbolic periodic
points in $\cC_f$ with the same stable dimension is homoclinically
related.

\item[(2)] For any $f \in \cU$ there exist two hyperbolic periodic
points $p,q \in \cC_f$ satisfying $\dim E^s_p>\dim E^s_q$ and
$\cC_f$ is the disjoint union of $H_f$ with $W^u(p)\cap W^s(q)$.
Moreover the points of $W^u(p)\cap W^s(q)$ are isolated in
$\cC_f$.

In particular, if $W^u(p)\cap W^s(q)\neq \emptyset$, the class
$\cC_f$ is not transitive.

\item[(3)] One has $\cD:= \{ f\in \cU \ : \ W^u(p)\cap W^s(q)\neq
\emptyset\}$.

Moreover, this set is a countable union of one-codimensional
submanifolds of $\cU$.

\item[(4)] The chain-recurrent set of any $f\in \cU$ is the union
of $\cC_f$ with a finite number of hyperbolic periodic points
(which depend continuously on $f$).
\end{itemize}

\end{teo*}

\begin{obs}
In the construction, the chain recurrence class $\cC_f$ is
partially hyperbolic with a one-dimensional central bundle. Thus,
it is also far from homoclinic tangencies. \finobs
\end{obs}

\begin{obs}The isolated points $\cC_f\setminus H_f$ are nonwandering for $f$.
However, they do not belong to $\Omega(f|_{\Omega(f)})$ (since
they are isolated in $\Omega(f)$ and non-periodic). \finobs
\end{obs}

\section{A mechanism for having isolated points in a chain recurrence class}\label{sectionMecanisme}

\subsection{Preliminaries on invariant bundles}

Consider $f \in \Diff^1(M)$ preserving a set $\Lambda$.

A $Df$-invariant subbundle $E \en T_\Lambda M$ is {\em uniformly
contracted} (resp. {\em uniformly expanded}) if there exists $N>0$
such that for every unit vector $v\in E$, we have

$$ \|Df^N v \| < \frac 1 2  \qquad (\mbox{resp. } > 2 ).$$

A $Df$-invariant splitting $T_\Lambda M = E^{ss} \oplus E^c \oplus
E^{uu}$ is \emph{partially hyperbolic} if $E^{ss}$ is uniformly
contracted, $E^{uu}$ is uniformly expanded, both are non trivial,
and if there exists $N>0$ such that for any $x\in \Lambda$ and any
unit vectors $v_s\in E^{ss}_{x}, v_c \in E^c_{x}$ and $v_u \in
E^{uu}_{x}$ we have:

\[  \|Df^N v_s \| < \frac 1 2 \|Df^N v_c \| < \frac 1 4 \|Df^N v_u \|. \]

$E^{ss}$, $E^c$ and $E^{uu}$ are called the {\em strong stable},
{\em center}, and {\em strong unstable} bundles.

\begin{obs}\label{RmqOrientationCentral}
We will sometimes consider a $Df$-invariant continuous orientation
of $E^c$. When $\Lambda$ is the union of two different periodic
orbits $O_{p},O_{q}$ and of a heteroclinic orbit
$\{f^n(x)\}\subset W^u(O_{p})\cap W^s(O_{q})$, such an orientation
exists if and only if above each orbit $O_{p}, O_{q}$, the tangent
map $Df$ preserves an orientation of the central bundle.

On a one-dimensional bundle, an orientation corresponds to a unit
vector field tangent. \finobs
\end{obs}

\subsection{Cuspidal periodic points}

Let $p$ be a hyperbolic periodic point whose orbit is partially
hyperbolic with a one-dimensional central bundle. When the central
space is stable, there exists a strong stable manifold $W^{ss}(p)$
tangent to $E^{ss}_{p}$ that is invariant by the iterates $f^\tau$
that fix $p$. It is contained in and separates the stable manifold
$W^s(p)$ in two {\em half stable manifolds} which contain
$W^{ss}(p)$ as a boundary.

Let us consider an orientation of $E^c_{p}$. The unit vector
defining the orientation goes inward on one half stable manifold
of $p$, that we call the {\em right half stable manifold}
$R^s(p)$. The other one is called the {\em left half stable
manifold} $L^s(p)$.

These half stable manifolds are invariant by an iterate $f^\tau$
which fixes $p$ if and only if the orientation of $E^c_p$ is
preserved by $Df^\tau_p$.

When the central space is unstable, one defines similarly the
right and left half unstable manifolds $R^u(p),L^u(p)$.

\begin{defi}
A hyperbolic periodic point $p$ is \emph{stable-cuspidal} if:
\begin{itemize}

\item[-] its orbit is partially hyperbolic, the central bundle is
one-dimensional and stable;

\item[-] one half stable manifold of $p$ intersects the
chain-recurrence class of $p$ only at $p$.

\end{itemize}
\finobs
\end{defi}

When the chain-recurrence class $\cC$ containing $p$ is not
reduced to the orbit $O_{p}$ of $p$, this forces the existence of
a $Df$-invariant orientation on the central bundle of $O_p$.

\begin{figure}[ht]\begin{center}
\input{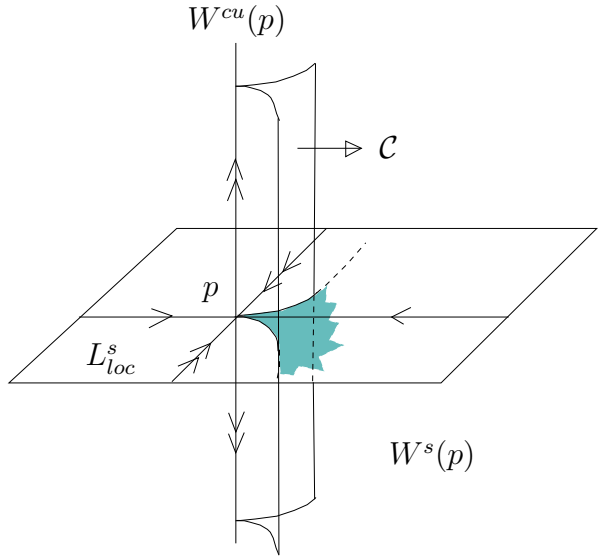}
\caption{\small{Geometry of a chain recurrence class $\cC$ near a
stable-cuspidal fixed point.}} \label{FiguraCostado}
\end{center}\end{figure}

In this case, the other half stable manifold intersects $\cC$ at
points different from $p$. The choice of the name has to do with
the geometry it imposes on $\cC\cap W^s(p)$ in a neighborhood of
$p$, see Figure \ref{FiguraCostado}.

This notion appears in~\cite{fragile}. It is stronger than the
notion of \emph{stable-boundary points} in \cite{CP}.

We can define in a similar way the \emph{unstable-cuspidal
points}.

\begin{obs}\label{RmqCuspidalRobust} If $p$ is a stable-cuspidal point,
then the hyperbolic continuation $p_g$ is still stable-cuspidal
for every $g$ that is $C^1$-close to $f$.

Indeed, there exists a compact set $\Delta\subset L^s(p)$ which
meets every orbit of $L^s(p)\setminus \{p\}$ and which is disjoint
from $\cR(f)$.

By semi-continuity of the chain-recurrent set, a small
neighborhood $V$ of $\Delta$ is disjoint from $\cR(g)$ for any $g$
close to $f$ and meets every orbit of the continuation of
$L^s(p)\setminus \{p\}$. \finobs
\end{obs}

\subsection{Description of the mechanism} Let $x$ be a point in a chain-recurrence class $\cC$. We introduce the following assumptions
(see figure \ref{FiguraMecanismo}).

\begin{itemize}
\item[(H1)] The class $\cC$ contains two periodic points $p,q$
such that $\dim(E^s_{p})=\dim(E^s_{q})+1$.

\item[(H2)] \begin{itemize}\item[(i)] The point $p$ is a
stable-cuspidal point.
 \item[(ii)] The point $q$ is an unstable-cuspidal point.\end{itemize}

\item[(H3)] The point $x$ belongs to $W^{u}(p)\cap W^{s}(q)$. The
union $\Lambda$ of the orbits of $x,p,q$ has a partially
hyperbolic decomposition with a one-dimensional central bundle.

Moreover there exists a $Df$-invariant continuous orientation of
the central bundle over $\Lambda$ such that $\cC$ is disjoint from
the half manifolds $L^s(p)$ and $R^u(q)$.

\end{itemize}

Note that from remark \ref{RmqOrientationCentral} and the fact
that a central orientation is preserved for cuspidal points, a
$Df$-invariant continuous orientation of the central bundle over
$\Lambda$ always exists.

\begin{figure}[ht]\begin{center}
\input{mecanismo.pstex_t}
\caption{\small{Hypothesis (H1)-(H3).}} \label{FiguraMecanismo}
\end{center}\end{figure}

\begin{prop}\label{propMecanisme}
Under (H1)-(H3), the point $x$ is isolated in the chain-recurrence
class $\cC$. In particular, $\cC$ is not transitive.
\end{prop}

\subsection{Proof of proposition~\ref{propMecanisme}}

Let $q$ be a periodic point whose orbit is partially hyperbolic
and whose central bundle is one-dimensional and unstable.  We
shall assume that there is an orientation in $E^c_q$ which is
preserved by $Df$.

We fix such an orientation of the central bundle $E^c_q$, so that
the left and right half unstable manifolds of $q$ are defined.

We denote by $d^u+1$ the unstable dimension of $q$.

Any $x \in W^s(q)$ has uniquely defined stable $E^s_x$ and center
stable $E^{cs}_x$ directions: the first one is the tangent space
$T_xW^s(q)$; a vector $v\in T_xM\setminus \{0\}$ belongs to the
second if the direction of its positive iterates $Df^{n}(v)$ stays
away from the directions of $E^{uu}_q$.

If $E'\subset E$ are two vector subspaces of $T_xM$ such that $E$
is transverse to $E^s_x$ and $E'$ is transverse to $E^{cs}_x$
(hence $E'$ is one-codimensional in $E$), then $F=E^{cs}_x\cap E$
is a one-dimensional space whose forward iterates converge to the
unstable bundle over the orbit of $q$. As a consequence, there
exists an orientation of $F$ which converges to the orientation of
the central bundle by forward iterations.

There is thus a connected component of $E\setminus E'$, such that
it intersects $F$ in the orientation of $F$ which converges
towards the central orientation, its closure is the \emph{right
half plane} of $E\setminus E'$.

The closure of the other component is the \emph{left half plane}
of $E\setminus E'$.










Consider a $C^1$-embedding $\varphi\colon [-1,1]^{d^u} \to M$ such
that $x:=\varphi(0)$ belongs to $W^s(q)$.

\begin{defi}
The embedding $\varphi$ is {\em coherent with the central
orientation at $q$} if
\begin{itemize}
\item[-] $E:=D_0\varphi(\RR^{d^u+1})$ and
$E':=D_0\varphi(\{0\}\times \RR^{d^u})$ are transverse to
$E^s_z,E^{cs}_z$ respectively;

\item[-] the half-plaque $\varphi([0,1]\times [-1,1]^{d^u})$ is
tangent to the right half-plane of $E\setminus E'$.
\end{itemize}
\finobs
\end{defi}

\medskip

Let $\Delta^u$ be a compact set contained in $R^u(q)\setminus
\{q\}$ which meets each orbit of $R^u(q)\setminus \{q\}$.

\begin{lema}\label{LemaLambdaLema} Let $\{\varphi_a\}_{a\in \cA}$ be a continuous family
of $C^1$-embeddings that are coherent with the central orientation
at $q$. Consider some $a_0\in \cA$ and a neighborhood $V^u$ of
$\Delta^u$. Then, there exist $\delta>0$ and some neighborhood $A$
of $a_0$ such that any point $z\in \varphi_a([0,\delta]\times
[-\delta,\delta]^{d^u})$ different from $\varphi_a(0)$ has a
forward iterate in $V^u$.
\end{lema}

\begin{proof} Let $\tau\geq 1$ be the period of $q$ and
$\chi\colon [-1,1]^d\to M$ be some coordinates such that

\begin{itemize}
\item[-] $\chi(0)=q$;

\item[-] the image $D^u:=\chi((-1,1)\times\{0\}^{d-d^u-1}\times
(-1,1)^{d^u})$ is contained in $W^u_{loc}(q)$;

\item[-] the image $D^{uu}:=\chi(\{0\}^{d-d^u}\times
(-1,1)^{d^u})$ is contained in $W^{uu}_{loc}(q)$;

\item[-] the image $D^{u,+}:=\chi([0,1)\times\{0\}^{d-d^u-1}\times
(-1,1)^{d^u})$ is contained in $R^u(q)$;

\item[-] $f^{-\tau}(\overline{D^u})$ is contained in $D^{u}$.

\end{itemize}

One deduces that there exists $n_0\geq 0$ such that:

\begin{itemize}

\item[(i)] Any point $z$ close to $\overline{D^{u,+}}\setminus
f^{-\tau}(D^u)$ has an iterate $f^k(z)$, $|k|\leq n_0$, in $V^u$.

\end{itemize}

The graph transform argument (see for instance~\cite[section
6.2]{KH}) gives the following generalization of the
$\lambda$-lemma.

\begin{af} There exists $N\geq 0$ and, for all $a$ in a neighborhood $A$ of
$a_0$, there exist some decreasing sequences of disks $(D_{a,n})$
of $[-1,1]^{d^u+1}$ and $(D'_{a,n})$ of $\{0\}\times [-1,1]^{d^u}$
which contain $0$ and such that for any $n\geq N$ one has, in the
coordinates of $\chi$:
\begin{itemize}
\item[-] $f^{n\tau}(D_{a,n})$ is the graph of a function $D^u\to
\RR^{d-d^u-1}$ that is $C^1$-close to $0$; \item[-]
$f^{n\tau}(D'_{a,n})$ is the graph of a function $D^{uu}\to
\RR^{d-d^u}$ that is $C^1$-close to $0$.
\end{itemize}
\end{af}

Let us consider $a\in A$.

The image by $f^{n\tau}$ of each component of $D_{a,n}\setminus
D'_{a,n}$ is contained in a small neighborhood of a component of
$D^{u}\setminus D^{uu}$.

The graph $f^{n\tau}(D'_{a,n})$ which is transverse to a constant
cone field around the central direction at $q$.

Since $\varphi$ is coherent with the central orientation at $q$,
one deduces that

\begin{itemize}
\item[(ii)] $f^{n\tau}\circ\varphi_a\left(([0,1]\times
[-1,1]^{d^u})\cap D_{a,n}\right)$ is contained in a small
neighborhood of $D^{u,+}$.
\end{itemize}

For $\delta>0$ small, any point $z\in
\varphi_a([-\delta,\delta]\times [-\delta,\delta]^{d^u})$
different from $\varphi_a(0)$ belongs to some $D_{a,n}\setminus
D_{a,n+1}$, with $n\geq N$. Consequently:

\begin{itemize}

\item[(iii)] Any $z\in \varphi_a([-\delta,\delta]\times
[-\delta,\delta]^{d^u})\setminus \{\varphi_a(0)\}$ has a forward
iterate in $D^u\setminus f^{-1}(D^u)$.

\end{itemize}

Putting the properties (i-iii) together, one deduces the announced
property.

\end{proof}

\begin{proof}[Proof of proposition~\ref{propMecanisme}]

We denote by $d^s+1$ (resp. $d^u+1$) the stable dimension of $p$
(resp. the unstable dimension of $q$)

so that the dimension of $M$ satisfies $d=d^s+d^u+1$.

Consider a point $x \in W^s(q) \cap W^u(p)$ satisfying (H3) and a
$C^1$-embeddeding $\varphi: [-1,1]^{d} \to M$ with $\varphi(0)=x$
such that:

\begin{itemize}

\item[-] $\varphi(\{0\}\times [-1,1]^{d^s} \times \{0\}^{d^u})$ is
contained in $W^s(q)$;

\item[-] $\varphi(\{0\}\times \{0\}^{d^s} \times [-1,1]^{d^u})$ is
contained in $W^u(p)$;

\item[-] $D_0\varphi.(1,0^{d^s},0^{d^u})$ is tangent to $E^c_x$
and has positive orientation.

\end{itemize}

Note that all the restrictions of $\varphi$ to $[-1,1]\times
\{a^s\}\times [-1,1]^{d^u}$ for $a^s\in \RR^{d^s}$ close to $0$,
are coherent with the central orientation at $q$.

Consider a compact set $\Delta^u\subset R^u(q)\backslash \{q\}$
that meets each orbit of $R^u(q)\backslash \{q\}$. Since $\cC$ is
closed and $q$ is unstable-cuspidal, there is a neighborhood $V^u$
of $\Delta^u$ in $M$ that is disjoint from $\cC$.

The lemma~\ref{LemaLambdaLema} can be applied: the points in
$\varphi([0,\delta]\times \{a^s\}\times [-\delta,\delta]^{d^u})$
distinct from $\varphi(0,a^s,0^{d^u})$ have an iterate in $V^u$,
hence do not belong to $\cC$.

This shows that

$$\cC\cap \varphi\left([0,\delta]\times [-\delta,\delta]^{d-1}\right) \subset \varphi\left(\{0\}\times [-\delta,\delta]^{d^s}\times \{0\}^{d^u}\right).$$

From (H3), if one reverses the central orientation and if one
considers the dynamics of $f^{-1}$, then all the restrictions of
$\varphi$ to $[-1,1]\times [-1,1]^{d^s}\times \{a^u\}$ for $a^u\in
\RR^{d^u}$ close to $0$, are coherent with the central orientation
at $p$.

One can thus argues analogously and gets:

$$\cC\cap \varphi\left([-\delta,0]\times [-\delta,\delta]^{d-1}\right) \subset \varphi\left(\{0\}\times \{0\}^{d^s}\times [-\delta,\delta]^{d^u}\right).$$

Both inclusions give that

$$\cC\cap \varphi\left([-\delta,\delta]^{d}\right)=\{\varphi(0)\},$$

which says that $x=\varphi(0)$ is isolated in $\cC$.

\end{proof}

\section{Construction of the example}\label{sectionExample}

In this part we build a collection of diffeomorphisms satisfying
the properties (I) and (II) stated in the theorem.

The construction will be made only in dimension $3$ for notational
purposes.

The generalization to higher dimensions is straightforward.

\subsection{Construction of a diffeomorphism}

Let us consider an orientation-preserving $C^\infty$
diffeomorphism $H$ of the plane $\RR^2$ and a closed subset
$D=D^-\cup C\cup D^+$ such that:

\begin{itemize}
\item[-] $H(\overline D)\subset \interior(D)$ and
$H(\overline{D^-\cup D^+})\subset \interior(D^-)$; \item[-] the
forward orbit of any point in $D^-$ converges towards a sink $S\in
D^-$;

\item[-] $C$ is the cube $[0,5]^2$ whose maximal invariant set is
a hyperbolic horseshoe.

\end{itemize}

On $C\cap H^{-1}(C)$ the map $H$ is piecewise linear, it preserves
and contracts by $1/5$ the horizontal direction and it preserves
and expands by $5$ the vertical direction (see
figure~\ref{herradura}):

\begin{itemize}
\item[-] The set $C\cap H(C)$ is the union of $4$ disjoint
vertical bands $I_1,I_2,I_3,I_4$ of width $1$. We will assume that
$I_1\cup I_2 \en (0,2+\frac 1 3)\times [0,5]$ and $I_3\cup I_4 \en
(2+\frac 2 3,5)\times [0,5]$.

\item[-] The preimage $H^{-1}(C)\cap C$ is the union of $4$
horizontal bands $H^{-1}(I_i)$. We will assume that
$H^{-1}(I_1\cup I_2) \en [0,5]\times (0,2+\frac 1 3)$ and
$H^{-1}(I_3\cup I_4) \en [0,5]\times (2+\frac 2 3,5)$.

\end{itemize}

\begin{figure}[ht]\begin{center}
\input{herradurade4V2.pstex_t}
\caption{\small{The map $H$.}} \label{herradura}
\end{center}\end{figure}

\begin{figure}[ht]\begin{center}
\input{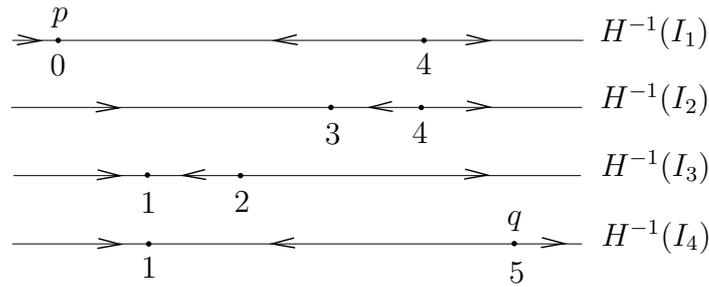}
\caption{\small{The map $g_x$ above each rectangle
$H^{-1}(I_i)$.}} \label{dimension1}
\end{center}\end{figure}

We define a $C^\infty$ diffeomorphism $F$ of $\RR^3$ whose
restriction to a neighborhood of $D\times [-1,6]$ it is a skew
product of the form

$$F\colon (x,t)\mapsto (H(x), g_x(t)),$$

where the diffeomorphisms $g_x$ are orientation-preserving and
satisfy (see figure \ref{dimension1}):

\begin{itemize}
\item[(P1)] $g_x$ does not depend on $x$ in the sets $H^{-1}(I_i)$
for every $i=1,2,3,4$.

\smallskip
\item[(P2)] For every $(x,t)\in D\times [-1,6]$ one has $4/5 <
g_x'(t) <6/5$.

\smallskip

\item[(P3)] $g_x$ has exactly two fixed points inside $[-1,6]$,
which are $\{0,4\}$, $\{3,4\}$, $\{1,2\}$ and $\{1,5\}$, when $x$
belongs to $H^{-1}(I_i)$ for $i$ respectively equal to $1,2,3$ and
$4$. All fixed points are hyperbolic, moreover,

\begin{itemize}

\item[-] $g'_x(t) < 1$ for $t \in [-1, 3+1/2]$ and $x\in
H^{-1}(I_1) \cup H^{-1}(I_2)$.

\item[-] $g'_x(t) > 1$ for $t\in [1+1/2,6]$ and $x\in  H^{-1}(I_3)
\cup H^{-1}(I_4)$.

\end{itemize}

\smallskip

\item[(P4)] For every $(x,t)\in (D^-\cup D^+)\times[-1,6]$ one has
$g_x(t)>t$.

\end{itemize}

We assume furthermore that the following properties are satisfied:

\begin{itemize}

\item[(P5)] $F(D\times [6,8])\subset \interior(D\times [6,8])$;

\item[(P6)] there exists a sink which attracts the orbit of any
point of $D\times [6,8]$;

\item[(P7)] $F$ coincides with a linear homothety outside a
compact domain;

\item[(P8)] any forward orbit meets $D\times [-1,8]$.

\end{itemize}

One can build a diffeomorphism which coincides with the identity
on a neighborhood of the boundary of $D_0\times (-2,9)$ and
coincides with $F$ in $D\times (-1,8)$ ($D_0$ denotes a small
neighborhood of $D$ in $\RR^2$).

This implies that, on any $3$-dimensional manifold, every isotopy
class of diffeomorphisms contains an element whose restriction to
an invariant set is $C^\infty$-conjugated to $F$.

\medskip

On any $3$-dimensional manifold, one can consider an
orientation-preserving Morse-Smale diffeomorphism and by surgery
replace the dynamics on a neighborhood of a sink by the dynamics
of $F$. We denote by $f_0$ the obtained diffeomorphism.

\subsection{First robust properties}\label{ss.robust}

We list some properties satisfied by $f_0$, which are also
satisfied by any diffeomorphism $f$ in a small $C^1$-neighborhood
$\cU_0$ of $f_0$.

\begin{description}

\item[Fixed points] By (P3), in each rectangle
$\interior(I_i)\times (-1,6)$, there exists two hyperbolic fixed
points $p_i,q_i$. Their stable dimensions are respectively equal
to $2$ and $1$. Since $p_1$ and $q_4$ will play special roles, we
shall denote them as $p=p_1$ and $q=q_4$.

\smallskip

\item[Isolation] The two open sets $V_0=\interior(D)\times (-1,8)$
and $V_1= V_0\setminus (C\times [-1,6])$ are isolating blocks,
i.e. satisfy $f(\overline{V_0}) \en V_0$ and $f(\overline{V_1})
\en V_1$.

For $V_0$, the property follows immediatly from the construction.

The closure of the second set $V_1$ can be decomposed as the union
of:

\begin{itemize}

\item[-] $D^+\times [-1,6]$, which is mapped into $(D^-\times
[-1,6])\cup (D\times [6,8])$,

\item[-] $D^-\times [-1,6]$ which is also mapped into $(D^-\times
[-1,6])\cup (D\times [6,8])$ and moreover has a foward iterate in
$D\times [6,8]$ by (P4),

\item[-] $D\times [6,8]$ which is mapped into itself and whose
limit set is a sink.

\end{itemize}

Hence, any chain-recurrence class which meets the rectangle $C
\times [-1,6]$ is contained inside. The maximal invariant set in
$C\times [-1, 6]$ will be denoted by $\cC$.

Any chain-recurrence class which meets $V_1$ coincides with the
sink of $D\times [6,9]$.

\smallskip

\item[Partial hyperbolicity] On $C\times [-1,6]\subset \RR^3$,
there exists some narrow cone fields $\cE^s,\cE^{cs}$ around the
coordinate direction $(1,0,0)$ and the plane $(x,0,z)$ which are
invariant by $Df^{-1}$. The vectors tangent to $\cE^s$ are
uniformly expanded by $Df^{-1}$.

Similarly there exists some forward invariant cone fields
$\cE^u,\cE^{cu}$ close to the direction $(0,1,0)$ and the plane
$(0,y,z)$.

In particular $\cC$ is partialy hyperbolic.

Moreover the tangent map $Df$ preserves the orientation of the
central direction such that any positive unitary central vector is
close to the vector $(0,0,1)$.

\smallskip

\item[Central expansion] Property (P2) holds for $f$ when one
replaces the derivative $g'_x(t)$ by the tangent map
$\|Df|_{E^c}(x,t)\|$ along the central bundle.

\smallskip

\item[Properties (H2) and (H3)] The point $p$ is stable-cuspidal
and the point $q$ is unstable-cuspidal. More precisely the left
half plaque of $W^s(p)$ and the right half plaque of $W^u(q)$ are
disjoint from $\cC$: since the chain-recurrence classes of $p$ and
$q$ are contained in $\cC$ this implies property (H2).

Moreover if there exists an intersection point $x\in W^u(p)\cap
W^s(q)$ for $f$, then by the isolating property it is contained in
$\cC$. By preservation of the central orientation, (H3) holds
also.

Let us explain how to prove these properties: it is enough to
discuss the case of the left half-plaque of $W^{s}(p)$ and
(arguing as in remark~\ref{RmqCuspidalRobust}) to assume that
$f=f_0$.

From (P2) and (P3), we have:

\begin{itemize}

\item[-] every point in $C\times [-1,0)$ has a backward iterate
outside $C\times [-1,6]$;

\item[-] the same holds for every point in $(C\setminus I_1)\times
\{0\}$;

\item[-] any point in $I_1\times \{0\}$ has some backward image in
$(C\setminus I_1)\times \{0\}$, unless it belongs to $W^{u}(p)$.

\end{itemize}

Combining these properties, one deduces that the connected
component of $W^{s}(p)\cap (C\times [-1,0])$ containing $p$
intersects $\cC$ only at $p$.

Note that this is a left half plaque of $W^s(p)$, giving the
required property.

\smallskip

\item[Hyperbolic regions] By (P3), the maximal invariant set in
$Q_p:=[0,5]\times [0,2+\frac 1 3] \times [-1,3+\frac 1 2]$ and
$Q_q:=[0,5]\times [2+\frac 2 3, 5] \times [1+\frac 1 2,6]$ are two
locally maximal transitive hyperbolic sets, denoted by $K_p$ and
$K_q$.

Their stable dimensions are $2$ and $1$ respectively. The first
one contains $p,p_2$, the second one contains $q,q_3$.

\smallskip

\item[Tameness (property (4) of the theorem)] since $f_0$ has been
obtained by surgery of a Morse-Smale diffeomorphism, the
chain-recurrent set in $M\setminus \cC$ is a finite union of
hyperbolic periodic orbits.

\end{description}

\medskip

Any $x\in \cC$ has a strong stable manifold $W^{ss}(x)$. Its
\emph{local} strong stable manifold $W^{ss}_{loc}(x)$ is the
connected component containing $x$ of the intersection $W^{ss}(x)
\cap C\times [-1,6]$.

It is a curve bounded by $\{0,5\} \times [0,5] \times [-1,6]$.
Symmetrically, we define $W^{uu}(x)$ and $W^{uu}_{loc}(x)$.

%




\subsection{Central behaviours of the dynamics}

We analyze the local strong stable and strong unstable manifolds
of points of $\cC$ depending on their central position.

\begin{lema}\label{lemaRegions}
There exists an open set $\cU_1 \en \cU_0$ such that for every $f\in \cU_1$ and $x\in\cC$:
\begin{itemize}
\item[(R1)] If $x\in R_1 := C\times [-1,4+\frac 1 2]$, then
$W^{uu}_{loc}(x) \cap W^s(p) \neq \emptyset$.

\item[(R2)] If $x\in R_2 := C\times [\frac 1 2,6]$, then
$W^{ss}_{loc}(x) \cap W^u(q) \neq \emptyset$.

\item[(R3)] If $x\in R_3 := C\times [\frac 1 2, 2+ \frac 1 2]$,
then $W^{ss}_{loc}(x) \cap W^{uu}_{loc}(y) \neq \emptyset$ for
some $y\in K_p$.

\item[(R4)] If $x\in R_4 := C\times [2+\frac 1 2, 4+\frac 1 2]$,
then $W^{uu}_{loc}(x) \cap W^{ss}_{loc}(y) \neq \emptyset$ for
some $y\in K_q$.

\end{itemize}

Moreover $p_2$ belongs to $R_2$ and $q_3$ belongs to $R_1$.

\end{lema}

\begin{proof}

Properties (R1) and (R2) follow directly from the continuous
variation of the stable and unstable manifolds. Similarly $p_2\in
R_2$ and $q_3\in R_1$ by continuity.

We prove (R3) with classical blender arguments (see
\cite{BDAnnals} and \cite[chapter 6]{BDV} for more details). The
set $K_p$ is a called \emph{blender-horseshoes} in \cite[section
3.2]{BDAbundanceTangencies}.

A \emph{cs-strip} $\cS$ is the image by a diffeomorphism $\phi:
[-1,1]^2 \to Q_p=[0,5]\times [0,2+\frac 1 3]\times [-1,3+\frac 1
2]$ such that:

\begin{itemize}
\item[-] The surface $\cS$ is tangent to the center-stable cone
field and meets $C\times [\frac 1 2, 2+\frac 1 2]$. \item[-] The
curves $\phi(t,[-1,1])$, $t\in [-1,1]$, are tangent to the strong
stable cone field and crosses $Q_p$, i.e. $\phi(t,\{-1,1\})
\subset \{0, 5\}\times [0,2+\frac 1 3] \times [-1,3+\frac 1 2]$.

\item[-] $\cS$ does not intersect $W^{u}_{loc}(p)\cup
W^u_{loc}(p_2)$.

\end{itemize}

The \emph{width} of $\cS$ is the minimal length of the curves
contained in $\cS$, tangent to the center cone, and that joins
$\phi(-1,[-1,1])$ and $\phi(1,[-1,1])$.

Condition (P2) is important to get the following (see \cite[lemma
6.6]{BDV} for more details):

\begin{af}
There exists $\lambda>1$ such that if $\cS$ is a $cs$-strip of
width $\eps$, then, either $f^{-1}(\cS)$ intersects
$W^{u}_{loc}(p)\cup W^u_{loc}(p_2)$ or it contains at least one
$cs$-strip with width $\lambda \eps$.
\end{af}

\dem
 Using (P2), the set $f^{-1}(\cS) \cap C\times [-1,6]$ is the
union of two bands crossing $C\times [-1,6]$: the first has its
two first coordinates near $H^{-1}(I_1)$, the second near
$H^{-1}(I_2)$.

Their width is larger than $\lambda \eps$ where $\lambda>1$ is a
lower bound of the expansion of $Df^{-1}$ in the central direction
inside $Q$. We assume by contradiction that none of them
intersects $W^{u}_{loc}(p)\cup W^{u}_{loc}(p_2)$, nor $C\times
[\frac 1 2, 2+\frac 1 2]$.

Since $\cS$ intersects $C\times[\frac 1 2, 2+\frac 1 2]$, from
conditions (P2) and (P3) the first band intersects $C\times[\frac
1 2, 4]$. By our assumption it is thus contained in
$C\times(2+\frac 1 2, 4]$. Using (P2) and (P3) again, this shows
that $\cS$ is contained in $C\times (2, 4]$.

The same argument with the second band shows that $\cS$ is
contained in $C\times [-1, 2)$, a contradiction.

\finobs

Repeating this procedure, we get an intersection point between
$W^{u}_{loc}(p)\cup W^u_{loc}(p_2)$ and a backward iterate of the
$cs-$strip. It gives in turn a transverse intersection point $z$
between the initial $cs-$strip and $W^{u}(p)\cup W^u(p_2)$. By
construction, all the past iterates of $z$ belong to $Q_p$. Hence
$z$ has a well defined local strong unstable manifold. In
particular, the intersection $y$ between $W^{uu}_{loc}(z)$ and
$W^s_{loc}(p)$ (which exists by (R1)) remains in $Q_p$ both for
future and past iterates, thus, it belongs to $K_p$.

For any point  $x\in \cC \cap R_3$, one builds a $cs$-strip by
thickening in the central direction the local strong stable
manifold. We have proved that this $cs-$strip intersects
$W^{uu}_{loc}(y)$ for some $y\in K_p$.

One con consider a sequence of thiner strips. Since $K_p$ is
closed and the local strong unstable manifolds vary continuously,
we get at the limit an intersection between $W^{ss}_{loc}(x)$ and
$W^{uu}_{loc}(y')$ for some $y'\in K_p$ as desired.

This gives (R3). Property (R4) can be obtained similarly.

\end{proof}

\bigskip

We have controled the local strong unstable manifold of points in
$R_1\cup R_4$ and the local strong stable manifold of points in
$R_2\cup R_3$.

Since neither $R_1 \cup R_4$ nor $R_2 \cup R_3$ cover completely
$C\times [-1,6]$ we shall also make use of the following result:

\begin{lema}\label{lemaUniquePoints} For every diffeomorphism
in a small $C^1$-neighborhood $\cU_2 \en \cU_0$ of $f_0$, the only
point whose complete orbit is contained in $C\times [-1,\frac 1
2]$ is $p$; symmetrically, the only point whose complete orbit is
contained in $C\times [4+\frac 1 2, 6]$ is $q$.
\end{lema}

\begin{proof}
We argue as for property (H2) in section~\ref{ss.robust}: the set
of points whose past iterates stay in $C\times [-1,\frac 1 2]$ is
the local strong unstable manifold of $p$. Since $p$ is the only
point in its local unstable manifold whose future iterates stay in
$C\times [-1,\frac 1 2 ]$ is $p$ we conclude.
\end{proof}

\subsection{Properties (I) and (II) of the theorem}

We now check that (I) and (II) hold for the region
$U=\interior(C\times [-1,6])$ and the neighborhood
$\cU:=\cU_1\cap \cU_2$.

\begin{prop}\label{p.approx-strong}
For any $f\in \cU$, $x\in \cC$, there are arbitrarily large
$n_q,n_p\geq 0$ such that $W^{uu}_{loc}(f^{n_q}(x))\cap
W^{ss}(y_q)\neq \emptyset$ and $W^{ss}_{loc}(f^{-n_p}(x))\cap
W^{uu}(y_p)\neq \emptyset$ for some $y_q\in K_q$, $y_p\in K_p$.
\end{prop}

\begin{proof}
If $\{f^n(x), n\geq n_0\}\subset C\times [4+\frac 1 2, 6]$, for
some $n_0\geq 0$, then $x\in W^{ss}(q)$ by
lemma~\ref{lemaUniquePoints}.

In the remaining case, there exist some arbitrarily large forward
iterates $f^n(x)$ in $R_1$, so that $W^{uu}_{loc}(f^n(x))$ meets
$W^s(p)$ by lemma~\ref{lemaRegions}.

Since $p$ is homoclinically related with $p_2$, by the
$\lambda$-lemma there exists $k\geq 0$ such that
$f^k(W^{uu}_{loc}(f^n(x)))$ contains $W^{uu}_{loc}(x')$ for some
$x'\in W^{s}(p_2)\cap R_4$ because $p_2\in R_4$.

By lemma~\ref{lemaRegions}, $f^k(W^{uu}_{loc}(f^n(x)))$ intersects
$W^{ss}_{loc}(y'_q)$ for some $y'_q\in K_q$ showing that
$W^{uu}_{loc}(f^n(x))\cap W^{ss}(y_q)\neq \emptyset$ with
$y_q=f^{-k}(y'_q)$ in $K_q$.

We have obtained the first property in all the cases. The second
property is similar.

\end{proof}

\medskip

The following corollary (together with the isolation property of
section~\ref{ss.robust}) implies that for every $f\in \cU$, the
properties (I) and (H1) are verified.

\begin{cor}\label{PropH1} For every $f\in \cU$ the set $\cC$ is
contained in a chain-transitive class.
\end{cor}

\begin{proof}
For any $\varepsilon>0$ and $x\in \cC$, there exists a
$\varepsilon$-pseudo-orbit $p=x_0,x_1,\dots,x_n=p$, $n\geq 1$,
which contains $x$.

Indeed by proposition~\ref{p.approx-strong}, and using that
$K_p,K_q$ are transitive and contain respectively $p$ and $q_3$,
there exists a $\varepsilon$-pseudo-orbit from $p$ to $q_3$ which
contains $x$.

By lemma~\ref{lemaRegions}, the unstable manifold of $q_3$
intersects the stable manifold of $p$, hence there exists a
$\varepsilon$-pseudo-orbit from $q_3$ to $p$.

We take the concatenation of these pseudo-orbits.

\end{proof}

\bigskip

Now, we show that (H3) holds for a $C^r$ dense set $\cD$ of $\cU$.

Since (H1) and (H2) are satisfied, proposition~\ref{propMecanisme}
implies that the property (II) of the theorem holds with the set
$\cD\subset \cU$.

In fact, as we noticed in section~\ref{ss.robust} it is enough to
get the following.

\begin{cor}\label{corH3densamente} For every $r\geq 1$,
the set
$$\cD=\{f\in \cU, \; W^{u}(p)\cap W^{s}(q)\neq \emptyset\}$$

is dense in $\cU\cap \Diff^r(M)$. It is a countable union of
one-codimensional submanifolds.
\end{cor}

In the $C^1$ topology, this result is direct consequence of the
connecting lemma (together with
proposition~\ref{p.approx-strong}).

The additional structure of our specific example allows to make
these perturbations in any $C^r$-topology.

\begin{proof} Fix any $f\in \cU$.

By proposition~\ref{p.approx-strong}, there exists $x\in K_q$ such
that $W^{u}(p)$ intersects $W^{ss}(x)$ at a point $y$ (notice that
$y\not \in K_q\cup \{p\}$).

Let $U$ be a neighborhood of $y$ such that:

\begin{itemize}
\item[-] $U$ is disjoint from the iterates of $y$, i.e. $\{f^n(y)
\ : \ n \in \ZZ \} \cap U = \{y\}$;

\item[-] $U$ is disjoint from $K_q\cup \{p\}$.

\end{itemize}

Given a $C^r$ neighborhood $\cV$ of the identity, there exists a
neighborhood $V\en U$ of $y$ such that, for every $z\in V$, the
set $\cV$ contains a diffeomorphism $g_z$ which coincides with the
identity in the complement of $U$ and maps $y$ at $z$.

Since $K_q$ is locally maximal, there exists $\bar x\in K_q\cap
W^s(q)$ near $x$. In particular $W^{ss}_{loc}(\overline{x})$
intersects $V$ in a point $z$ whose backward orbit is disjoint
from $U$.

For the diffeomorphism $h= g_z\circ f$ (which is $C^r$-close to
$f$) the manifolds $W^{s}(q)$ and $W^{u}(p)$ intersect.

Indeed both $f$ and $h$ satisfy $f^{-1}(y)\in W^{u}(p)$ and $z\in
W^{ss}_{loc}(\overline{x})$.

Since  $W^{ss}_{loc}(\overline{x})\en W^{ss}(q)$ and
$h(f^{-1}(y))= z$ we get the conclusion.

For each integer $n\geq 1$, the manifolds $f^n(W^{uu}_{loc}(p))$
and $W^{ss}_{loc}(q)$ have disjoint boundary and intersect in at
most finitely many points.

One deduces that the set $\cD_n$ of diffeomorphisms such that they
intersect is a finite union of one-codimensional submanifold of
$\cU$.

The set $\cD$ is the countable union of the $\cD_n$.

\end{proof}

\subsection{Other properties}\label{sectionAutresProprietes}

We here show properties (1), (2) and (3) of the theorem.

\begin{prop}\label{propPerteneceHomoclinica} For every $f\in \cU$ and $x\in \cC$ we have:
\begin{itemize}
\item[-] If $x\not\in W^{s}(q)$, there exist large $n\geq 0$ such
that $W^{uu}_{loc}(f^n(x))\cap W^s(p)\neq \emptyset$.

\item[-] If $x \not\in W^{u}(p)$, there exists large $n\geq 0$
such that $W^{ss}_{loc}(f^{-n}(x))\cap W^u(q)\neq\emptyset$.

\end{itemize}

Moreover, in the first case $x$ belongs to the homoclinic class of
$p$ and in the second it belongs to the homoclinic class of $q$.
\end{prop}

\begin{proof}

By lemma~\ref{lemaUniquePoints}, any point $x\in \cC\setminus
W^{s}(q)$ has arbitrarily large iterates $f^n(x)$ in $R_1$,
proving that $W^{uu}_{loc}(f^n(x))\cap W^s(p)\neq \emptyset$.

In particular, $W^s(p)$ intersects transversaly $W^{uu}_{loc}(x)$
at points arbitrarily close to $x$. On the other hand by
proposition~\ref{p.approx-strong}, there exists a sequence $z_n$
converging to $x$ and points $y_n\in K_p$ such that $z_n\in
W^{u}(y_n)$ for each $n$, proving that $W^{uu}_{loc}(z_n)$
intersects $W^u(p)$ transversaly at a point close to $x$ when $n$
is large.

By the $\lambda$-lemma, $W^{uu}_{loc}(y_n)$ is the $C^1$-limit of
a sequence of discs contained in $W^u(p)$. This proves that
$W^u(p)$ and $W^s(p)$ have a transverse intersection point close
to $x$, hence $x$ belongs to the homoclinic class of $p$.

The other properties are obtained analogously.
\end{proof}

\medskip

Let $H_f$ denotes the homoclinic class of $p$. The next gives
property (1) of the theorem.

\begin{cor}\label{corHomoclinicCoincide} For every $f\in \cU$,
the homoclinic class of any hyperbolic periodic point of $\cC$
coincides with $H_f$. Moreover, the periodic points in $\cC$ of
the same stable index are homoclinically related.
\end{cor}

\begin{proof}

Let $z\in \cC$ be a hyperbolic periodic point whose stable index
is $2$.

By proposition~\ref{p.approx-strong} $W^{ss}(z)$ intersects
$W^{uu}_{loc}(y)$ for some $y\in K_p$, this implies that $W^s(z)$
intersects $W^{uu}_{loc}(y)$ and since $W^{uu}_{loc}(y)$ is
accumulated by $W^{u}(p)$ we get that $W^s(z)$ intersects
$W^u(p)$. Now, by proposition~\ref{propPerteneceHomoclinica},
$W^{u}(z)$ intersects $W^{s}(p)$. Moreover the partial
hyperbolicity implies that the intersections are transversal,
proving that $z$ and $p$ are homoclinically related.

One shows in the same way that any hyperbolic periodic point whose
stable index is $1$ is homoclinically related to $q$.

It remains to prove that the homoclinic classes of $p$ and $q$
coincide.

The homoclinic class of $q$ contains a dense set of points $x$
that are homoclinic to $q_3$. In particular, $x$ does not belong
to $W^{u}(q)$, hence belongs to the homoclinic class of $p$ by
proposition~\ref{propPerteneceHomoclinica}.

This gives one inclusion. The other one is similar.

\end{proof}

\medskip

Properties (2) and (3) of the theorem follow from
corollary~\ref{corH3densamente} and the following.

\begin{cor}\label{PropUnicasIntersecciones} For every $f\in \cU$
we have $\cC \backslash H_f= W^{s}(q)\cap W^{u}(p)$.
\end{cor}

\begin{proof}
By corollary~\ref{corHomoclinicCoincide}, a point $x\in \cC
\backslash H_f$ does not belong to the homoclinic class of $q$
(nor to the homoclinic class of $p$ by definition of $H_f$).

Proposition~\ref{propPerteneceHomoclinica} gives $\cC \backslash
H_f\subset W^{s}(q)\cap W^{u}(p)$.

Proposition~\ref{propMecanisme} proves that the points of
$W^{s}(q)\cap W^{u}(p)$ are isolated in $\cC$. Since any point in
a non-trivial homoclinic class is limit of a sequence of distinct
periodic points of the class we conclude that $W^{s}(q)\cap
W^{u}(p)$ and $H_f$ are disjoint.

\end{proof}

The proof of the theorem is now complete.

\newpage

\end{document}